%% file: main.tex
\documentclass[11pt]{book}

\usepackage{makeidx}          
\usepackage{amsmath}          
\usepackage{amssymb}          
\usepackage{latexsym}         
\usepackage{amsxtra}          
\usepackage[latin1]{inputenc} 
\usepackage[T1]{fontenc}       
\usepackage{eucal}            
\usepackage{bbm}              
\usepackage{longtable}        
\usepackage{exscale}          
\usepackage{theorem}          
\usepackage{graphicx}         
\usepackage[sort]{cite}       
\usepackage{mathrsfs}         
\usepackage{array}            
\usepackage{stmaryrd}         
\usepackage{paralist}         
\usepackage{color}
\usepackage[a4paper]{geometry} 
\usepackage{wasysym}          
\usepackage{mathtools}        
\usepackage{fancyhdr}         
\usepackage{ifthen}           
\usepackage{gitinfo}          
\usepackage[expansion=false]{microtype}  
\usepackage[final]{hyperref}  

\usepackage[curve]{xypic}           
\newcommand{\pictype}{pdf_t}
\newcommand{\rawpictype}{pdf}

\usepackage{diagxy}                 

\makeatletter
\def\cleardoublepage{\clearpage\if@twoside \ifodd\c@page\else
    \hbox{}
    \vspace*{\fill}
\begin{center}

\end{center}                                                %
    \vspace{\fill}                                              %
    \thispagestyle{empty}                                       %
    \newpage                                                    %
    \if@twocolumn\hbox{}\newpage\fi\fi\fi}                      %
\makeatother

\makeindex                                       
\newcommand{\Index}[1]   {#1\index{#1}}          
\newcommand{\emIndex}[1] {\emph{#1}\index{#1}}   

\renewcommand{\mathbb}[1]{\mathbbm{#1}}          

\allowdisplaybreaks                              

\newcommand{\refitem}[1] {~\textit{\ref{#1}.)}}

\theoremheaderfont{\normalfont\bfseries}
\theorembodyfont{\itshape}
\newtheorem{lemma}{Lemma}[section]
\newtheorem{proposition}[lemma]{Proposition}
\newtheorem{theorem}[lemma]{Theorem}
\newtheorem{corollary}[lemma]{Corollary}
\newtheorem{definition}[lemma]{Definition}

\theorembodyfont{\rmfamily}

\newtheorem{example}[lemma]{Example}
\newtheorem{remark}[lemma]{Remark}

\newtheorem{Blemma}{Lemma}[chapter]
\newtheorem{Btheorem}[Blemma]{Theorem}
\newtheorem{Bdefinition}[Blemma]{Definition}

\theorembodyfont{\rmfamily}

\newtheorem{Bremark}[Blemma]{Remark}

\makeatletter
\newcommand\qedsymbol{\hbox{$\boxempty$}}
\newcommand\qed{\relax\ifmmode\boxempty\else
  {\unskip\nobreak\hfil\penalty50\hskip1em\null\nobreak\hfil\qedsymbol
  \parfillskip=\z@\finalhyphendemerits=0\endgraf}\fi}

\newcommand\subqedsymbol{\hbox{$\triangledown$}}
\newcommand\subqed{\relax\ifmmode\triangledown\else
  {\unskip\nobreak\hfil\penalty50\hskip1em\null\nobreak\hfil\subqedsymbol
  \parfillskip=\z@\finalhyphendemerits=0\endgraf}\fi}
\makeatother

\newenvironment{proof}[1][{}]{\par\noindent Proof{#1}. }{\qed}
\newenvironment{subproof}[1][{}]{\par\noindent Proof{#1}. }{\subqed}
\newenvironment{altproof}[1][{}]{\par\noindent Alternative
  Proof{#1}. }{\qed}

\newenvironment{remarklist}{\begin{compactenum}[\itshape i.)]}{\end{compactenum}}
\newenvironment{examplelist}{\begin{compactenum}[\itshape i.)]}{\end{compactenum}}
\newenvironment{lemmalist}{\begin{compactenum}[\itshape i.)]}{\end{compactenum}}
\newenvironment{theoremlist}{\begin{compactenum}[\itshape i.)]}{\end{compactenum}}
\newenvironment{propositionlist}{\begin{compactenum}[\itshape i.)]}{\end{compactenum}}
\newenvironment{definitionlist}{\begin{compactenum}[\itshape i.)]}{\end{compactenum}}

\input{macros}

\title{{\normalsize Lecture Notes}\\[0.5cm]
  Geometric Wave Equations}

\author{\textbf{Stefan Waldmann}
  \\[0.5cm]
  Department Mathematik\\
  Friedrich-Alexander Universität Erlangen-Nürnberg\\
  Cauerstrasse 11\\
  91058 Erlangen\\
  Germany \\[0.5cm]
  {\small Contact: \texttt{Stefan.Waldmann@math.fau.de}}
}

\date{\input{abstract}}

\begin{document}

%
%

\thispagestyle{empty}
\pagestyle{empty}

%
%

\maketitle

\newpage

%
%

\input{preface}

%
%

\pagenumbering{roman}
\pagestyle{fancy}
\setcounter{page}{3}

%
%

\tableofcontents

%
%

%
%

\input{intro}

%
%

\makeatletter
\renewcommand{\theequation}{\thesection.\arabic{equation}}
\@addtoreset{equation}{section}
\makeatother

%
%

\input{chap1}
\input{chap2}
\input{chap3}
\input{chap4}

%
%

\begin{appendix}

%
%

\input{appendixA}

%
%

\input{appendixB}

\end{appendix}

%
%

\addcontentsline{toc}{chapter}{References}

\input{refs.tex}
%
%

\cleardoublepage
\addcontentsline{toc}{chapter}{Index}
\printindex

%
%

\end{document}

%% file: macros.tex
\newcommand{\I}          {\mathrm{i}}
\newcommand{\E}          {\mathrm{e}}
\newcommand{\D}          {\operatorname{\mathrm{d}}}
\newcommand{\cc}[1]      {\overline{{#1}}}

\newcommand{\RE}         {\mathsf{Re}}

\newcommand{\cl}         {\mathrm{cl}}

\newcommand{\at}[1]      {\big|_{#1}}
\newcommand{\At}[1]      {\Big|_{#1}}
\newcommand{\argument}   {\,\cdot\,}

\newcommand{\id}         {\operatorname{\mathsf{id}}}

\newcommand{\pr}         {\mathrm{pr}}

\newcommand{\image}      {\operatorname{{\mathrm{im}}}}

\newcommand{\tr}         {\operatorname{\mathsf{tr}}}
\newcommand{\rank}       {\operatorname{\mathrm{rank}}}

\newcommand{\diag}       {\operatorname{\mathrm{diag}}}
\newcommand{\Trans}      {{\mathrm{\scriptscriptstyle{T}}}}
\newcommand{\Anti}       {\Lambda}
\newcommand{\Sym}        {\mathrm{S}}
\newcommand{\SP}[1]      {\left\langle{#1}\right\rangle}
\newcommand{\ann}        {\mathrm{ann}}

\newcommand{\supp}       {\operatorname{\mathrm{supp}}}
\newcommand{\Lie}        {\operatorname{\mathscr{L}\!}}
\newcommand{\Laplace}    {\operatorname{\Delta}}
\DeclareMathSymbol\dAlembert  {\mathord}{AMSa}{"03}

\newcommand{\divergenz}  {\operatorname{\mathrm{div}}}
\newcommand{\gradient}   {\operatorname{\mathrm{grad}}}

\newcommand{\Fun}[1][{k}] {\mathcal{C}^{#1}}
\newcommand{\Cinfty}     {\Fun[\infty]}
\newcommand{\Sec}[1][{k}]   {\Gamma^{#1}}
\newcommand{\Secinfty}   {\Sec[\infty]}

\newcommand{\ins}        {\operatorname{\mathrm{i}}}

\newcommand{\inss}       {\operatorname{\mathrm{i}_{\mathrm{s}}}}
\newcommand{\HdR}        {\mathrm{H}_{\scriptscriptstyle\mathrm{dR}}}

\newcommand{\SymD}       {\operatorname{\mathsf{D}}}
\newcommand{\Tor}        {\operatorname{\mathrm{Tor}}}
\newcommand{\Ric}        {\operatorname{\mathrm{Ric}}}
\newcommand{\scal}       {\operatorname{\mathrm{scal}}}

\newcommand{\Diffop}     {\operatorname{\mathrm{DiffOp}}}

\newcommand{\loc}        {\mathrm{loc}}
\newcommand{\Dichten}[1] {|\Anti^{\mathrm{top}}|#1}
\newcommand{\vol}        {\operatorname{\mathrm{vol}}}
\newcommand{\stdrep}     {\operatorname{\varrho_{\scriptscriptstyle\mathrm{Std}}}}
\newcommand{\prol}       {\mathrm{prol}}
\newcommand{\PB}[1] {\left\{{#1}\right\}}

\newcommand{\Pol}        {\operatorname{\mathrm{Pol}}}

\newcommand{\group}[1]   {\mathrm{#1}}

\newcommand{\Hom}        {\operatorname{\mathsf{Hom}}}
\newcommand{\End}        {\operatorname{\mathsf{End}}}

\newcommand{\norm}[1]    {\left\|{#1}\right\|}

\newcommand{\singsupp}   {\operatorname{{\mathrm{sing\,supp}}}}
\newcommand{\seminorm}[1][{}] {\operatorname{\mathrm{p}}_{#1}}
\newcommand{\ord}        {\operatorname{\mathrm{ord}}}

\newcommand{\tensor}[1][{}]   {\mathbin{\otimes_{\scriptscriptstyle{#1}}}}
\newcommand{\extensor}[1][{}] {\mathbin{\boxtimes_{\scriptscriptstyle{#1}}}}


%% file: abstract.tex
%
%

\begin{minipage}{0.8\linewidth}
    \begin{small}
        In these lecture notes we discuss the solution theory of
        geometric wave equations as they arise in Lorentzian geometry:
        for a normally hyperbolic differential operator the existence
        and uniqueness properties of Green functions and Green
        operators is discussed including a detailed treatment of the
        Cauchy problem on a globally hyperbolic manifold both for the
        smooth and finite order setting. As application, the classical
        Poisson algebra of polynomial functions on the initial values
        and the dynamical Poisson algebra coming from the wave
        equation are related. The text contains an introduction to the
        theory of distributions on manifolds as well as detailed
        proofs.
    \end{small}
\end{minipage}


%% file: preface.tex
%
%

\chapter*{Preface}
\thispagestyle{empty}

These lecture notes grew out of a two-semester course on wave
equations on Lorentz manifolds which I gave in Freiburg at the physics
department in the winter term 2008/2009 and the following summer term
2009. This lecture originated from a long term project on the
deformation quantization of classical field theories started some nine
years before: the aim was to understand recent developments on
quantization following the works of Dütsch and Fredenhagen
\cite{duetsch.fredenhagen:2001a, duetsch.fredenhagen:2001b,
  duetsch.fredenhagen:1999a}. As time passed, the beautiful book of
Bär, Ginoux, and Pfäffle \cite{baer.ginoux.pfaeffle:2007a} on the
global theory of wave equations appeared and provided the basis for a
revival of that old project. So the idea of presenting the results of
\cite{baer.ginoux.pfaeffle:2007a} to a larger audience of students was
born.  The resulting lectures aimed at master and PhD students in
mathematics and mathematical physics with some background in
differential geometry and a lively interest in the analysis of
hyperbolic partial differential equations. Though both, the lecture
and these notes, followed essentially the presentation of
\cite{baer.ginoux.pfaeffle:2007a}, I added more detailed proofs and
some background material which hopefully make this material easily
accessible already for students.

During the preparation of these lecture notes many colleagues and
friends gave me their help and support. To all of them I am very
grateful: First of all, I would like to thank Frank Pfäffle for his
continuous willingness to explain many details of
\cite{baer.ginoux.pfaeffle:2007a} to me. Without his help, neither the
lecture nor these lecture notes would have been possible in the
present form. Also, I would like to thank Michael Dütsch and Klaus
Fredenhagen for continuing discussions concerning their works as well
as on related questions on deformation quantization of classical field
theories, thereby constantly raising my interest in the whole subject.
Florian Becher helped not only with the exercise and discussion group
for the students but is ultimately responsible for this project by
pushing me to ``give a lecture on the book of Bär, Ginoux, and
Pfäffle''. I am also very grateful to Domenico Giulini who helped me
out in many questions on general relativity and gave me access and
guidance to various references. Moreover, I am indebted to Stefan Suhr
for helping me in many questions on Lorentz geometry and improving
various arguments during the lecture. I would like to thank also all
the participant of the course who brought the lecture to success by
their constant interest, their questions, and their remarks on the
manuscript of these notes, in particular Jan Paki. Finally, I am very
much obliged to Jan-Hendrik Treude for taking care of the
\LaTeX-files, the Xfig-pictures, and all the typing as well as for his
numerous comments and remarks. Without his help, the manuscript would
have never been finished.

Most grateful I am for my children Silvia, Richard, Sonja, and Robert,
as a source of unlimited inspiration and for Viola, for her continuous
support, both morally and scientifically.

\medskip
\medskip
\medskip
\medskip
\medskip

\begin{flushright}
    \noindent
    Erlangen, August 2012 \hfill \textit{Stefan Waldmann}\\
\end{flushright}


%% file: intro.tex
%
%

\makeatletter
\chapter*{Introduction and Overview\@mkboth{INTRODUCTION AND OVERVIEW}{INTRODUCTION AND OVERVIEW}}
\addcontentsline{toc}{chapter}{Introduction and Overview}
\makeatother
\pagenumbering{arabic}
\setcounter{page}{1}

The theory of linear partial differential equations can be divided
into three principal parts: the first is the elliptic theory of
equations like the Laplace equation, the next is the parabolic theory
being the habitat of the heat equation, and the third is the
hyperbolic theory. All three differ in their behaviour, concepts, and
applications\index{Caesar}.

It will be the hyperbolic theory where the wave equation
\begin{equation}
    \label{eq:TheWaveEquation}
    \frac{1}{c^2} \frac{\partial^2 u}{\partial t^2}
    -
    \sum_{i=1}^{n-1} \frac{\partial^2 u}{\partial x_i^2}
    =
    0
\end{equation}
on $\mathbb{R}^n$ provides the first and most important example. While
for the elliptic theory the boundary problem is characteristic, for
the hyperbolic situation the main task is to understand an initial
value problem: for time $t=0$ one specifies the solution $u(0, x)$ and
its first time derivative $\frac{\partial u}{\partial t}(0, x)$ for
all $x \in \mathbb{R}^{n-1}$ and seeks a solution of the wave equation
with these prescribed initial values. Of course, also for the wave
equation one can pose boundary condition on top of the initial value
problem. Together with the question of how (continuous) the solution
depends on the initial conditions this becomes the Cauchy problem for
hyperbolic equations.

The relevance of the wave equation as coming from the science and in
particular from physics is overwhelming; we indicate just two major
occurrences: on a phenomenological level it describes propagating
waves in elastic media in a linearized approximation. This
approximation is typically well justified as long as the displacements
are not too big. Then the wave equation provides a good model for many
everyday situations like water waves, elastic vibrations of solids, or
propagation of sound. The constant $c$ in the wave equation is then
the speed of propagation and a characteristic quantity of the
material. On a more fundamental level, and more important for our
motivation, is the appearance of the wave equation in various physical
theories of fundamental interactions. Most notable here is Maxwell's
theory of electromagnetic fields. In this context, the wave equation
appears as an exact and fundamental equation describing the
propagation of electromagnetic waves (light, radio waves, etc.) in the
vacuum. Remarkably, it is a field equation not relying on any sort of
carrier material like the hypothetical ether. The constant $c$ becomes
the speed of light, one of the few truly fundamental constants in
physics. But even beyond Maxwell's theory the wave equation and its
generalizations like the Klein-Gordon equation provide the linear part
of all known fundamental field theories.

Needless to say, it is worth studying such wave equations. But which
framework should be taken to formulate the problem in a mathematically
meaningful and yet still interesting way?

A short look at the wave equation shows that it is invariant under the
affine pseudo-orthogonal group $\group{O}(1, n-1) \ltimes
\mathbb{R}^n$ in the sense that the natural affine action of
$\group{O}(1, n-1) \ltimes \mathbb{R}^n$ on $\mathbb{R}^n$ pulls back
solutions of the wave equation to solutions again. In more physical
terms we have the invariance group of special relativity, the
Poincar\'e group. This already indicates to take a \emph{geometric}
point of view and interpret the wave equation as coming from the
d'Alembert operator $\dAlembert$ corresponding to the Minkowski metric
$\eta = \diag(+1, -1, \ldots, -1)$. Indeed, this point of view opens
the door for various generalizations if we replace $\eta$ and
$\mathbb{R}^n$ by an arbitrary Lorentz metric $g$ on an arbitrary
manifold $M$: we still have a d'Alembert operator (coming from $g$)
and hence a wave equation. In more physical terms we pass from special
to general relativity. But even if one is not interested in geometry a
priori, generalizations of the wave equation like
\begin{equation}
    \label{eq:NonConstantCoefficients}
    \sum_{i, j} A^{ij} \frac{\partial^2 u}{\partial x^i \partial x^j}
    +
    \sum_i B^i \frac{\partial u}{\partial x^i}
    +
    C u
    =
    0,
\end{equation}
with coefficient \emph{functions} $A^{ij}$, $B^i$, and $C$ on
$\mathbb{R}^n$ such that the matrix $(A^{ij}(x))$ has signature $(+,
-, \ldots, -)$ at every point $x \in \mathbb{R}^n$, can be treated
best only after a geometric interpretation of the functions
$A^{ij}$. Otherwise, it will be almost impossible to get hands on the
Cauchy problem of such a wave equation with non-constant
coefficients. In fact, the first naive idea would be to find adapted
coordinates on order to bring \eqref{eq:NonConstantCoefficients} to
the form \eqref{eq:TheWaveEquation}, at least concerning the second
order derivatives. However, generically this has to fail since the
typically non-zero curvature of the metric corresponding to the
coefficients $A^{ij}$ is precisely the obstruction to get constant
coefficients in front of the leading orders of differentiation by a
change of coordinates.  This brings us back to a geometric point of
view which we will take in the following.

%
%

\section*{The Geometric Framework}

The wave equations we will discuss are located on a Lorentz manifold,
i.e. on a smooth $n$-dimensional manifold $M$ equipped with a smooth
Lorentz metric $g$. We choose the signature $(+, -, \ldots, -)$ as
common in (quantum) field theory but probably less common in general
relativity. The notions of light-, time-, and spacelike vectors,
future and past, causality, etc. which we will develop in the sequel,
have their origin in the theory of general relativity which is the
main source of inspiration in Lorentz geometry. In particular, the
notion of a spacetime will be used synonymously for a Lorentz
manifold.

The metric allows to speak of the d'Alembert operator $\dAlembert$
acting on the smooth functions on $M$. While this gives already many
interesting wave equations there are still two directions of
generalization: first, we would like to incorporate also lower order
terms of differentiation as in
\eqref{eq:NonConstantCoefficients}. Second, many application like e.g.
Maxwell's theory require to go beyond the \emph{scalar} wave equations
and need ``multicomponent'' functions $u^\alpha$ instead of a single,
scalar one. These components may even be coupled in a non-trivial way.

Both situations can be combined into the following framework. We take
a vector bundle $E \longrightarrow M$ over $M$ and consider a linear
second order differential operator $D$ acting on the sections of $E$
with leading symbol being the same as for the scalar d'Alembert
operator. Such a \emph{normally hyperbolic} differential operator will
have the local form
\begin{equation}
    \label{eq:LocalFormForNormallyHyperbolic}
    Du =
    \left(
    g^{ij} \frac{\partial^2 u^\alpha}{\partial x^i \partial x^j}
    +
    B^{i\alpha}_\beta \frac{\partial u^\beta}{\partial x^i}
    +
    C^\alpha_\beta u^\beta\right) e_\alpha,
\end{equation}
where the section $u = u^\alpha e_\alpha$ is expressed locally in
terms of a local frame $\{e_\alpha\}$ of $E$ and we use local
coordinates $\{x^i\}$ on $M$. Here $g^{ij}$ are the coefficients of
the (inverse) metric tensor $g$ while $B^{i\alpha}_\beta$ and
$C^\alpha_\beta$ are coefficient functions determined by $D$. In this
expression and from now on we shall use Einstein's summation
convention that pairs of matching coordinate or frame indexes are
automatically summed over their range.

A differential operator $D$ like in
\eqref{eq:LocalFormForNormallyHyperbolic} makes sense even on any
semi-Riemannian manifold. For the formulation of the Cauchy problem we
need the Lorentz signature and two extra structures beside the
metric. The first is a time orientation which separates future from
past. This will allow for notions of causality and thus for the
notions of advanced and retarded solutions of the wave equation. From
a physics point of view such a time orientation is absolutely
necessary to have a true interpretation of $(M, g)$ as a
spacetime. The second ingredient is that of a hypersurface $\Sigma$ in
$M$ on which we can specify the initial values. Thus $\Sigma$
corresponds to ``$t=0$'' in this geometric context. At first sight any
spacelike hypersurface might be suitable. However, already in
$\mathbb{R}^n$ the $t=0$ hypersurface has additional properties: it
divides $\mathbb{R}^n$ into two disjoint pieces, the future and the
past of $t=0$. Moreover, \emph{every} inextensible causal curve has to
pass through this $t=0$ hypersurface in precisely one
point. Physically speaking, this means that knowing things on $\Sigma$
allows to compute the entire time evolution in a deterministic
way. This is the main idea behind an initial value problem. Thus we
can already anticipate that this feature will turn out to be crucial
for a good Cauchy problem. In general, a spacelike hypersurface
$\Sigma$ will be called a \emph{Cauchy hypersurface} if it satisfies
this condition: every inextensible causal curve passes in exactly one
point through $\Sigma$. It is a non-trivial and in fact quite recent
theorem that the existence of such a \emph{smooth} Cauchy hypersurface
is equivalent to the notion of a globally hyperbolic Lorentz
manifold. Moreover, having one such Cauchy hypersurface allows already
to split $M$ into a time axis and spacelike directions, i.e. $M \cong
\mathbb{R} \times \Sigma$, in such a way that also the metric becomes
block-diagonal. We will have to explain all these notions in more
detail.

%
%

\section*{The Analytic Framework}

After setting the geometric stage we also have to specify the analytic
aspects properly in order to obtain a complete formulation of the
Cauchy problem. Handling linear partial differential equations allows
for various approaches. Most notably, one can use Sobolev space
techniques or distribution theory. In the sequel, we will exclusively
use the distributional approach for reasons which are not even that
easy to explain. Nevertheless, let us try to motivate our choice:

At first, physicists are usually more adapted to the notions of
distributions, at least on a heuristic level, than to Sobolev spaces
and their usage. Moreover, and more important, the solution to the
Cauchy problem using distribution theory relies on the notion of
\emph{Green functions} also called \emph{fundamental solutions}. These
are particular distributional solutions of the wave equation with a
$\delta$-distribution as inhomogeneity. The collection of all these
Green functions can be combined into a single operator, the Green
operator. Very informally, this will be an ``inverse'' of the
differential operator $D$. Now these Green operators allow for a very
efficient description of the solutions to the Cauchy problem and are
hence worth to be studied. Finally, and this might be the most
important reason to choose the distributional approach, these Green
operators appear as fundamental ingredients, the \emph{propagators},
for every quantum field theory build on top of the classical field
theory described by the wave equation. Even though we do not enter the
discussion of quantizing the classical field theory we at least
provide the starting point by constructing the Poisson algebra of the
classical theory. The Poisson bracket is then defined by means of the
Green operators and will allow us to view the time evolution of the
initial values as a ``Hamiltonian system'' with infinitely many
degrees of freedom. The interest in this Hamiltonian picture is the
ultimate reason for us to favour the distributional approach over the
Sobolev one.  Even though we do not discuss this here, there is yet
another reason why the distributional approach is interesting: it is
within this framework where one can discuss the propagation of
singularities most naturally by means of wavefront analysis.

Within the distributional approach we will have an interplay of very
singular objects, the distributional sections of vector bundles, and
very regular ones, the smooth sections of the corresponding dual
bundles. Here smooth stands for $\Cinfty$, i.e. infinitely often
differentiable. However, at many places we will pay attention to the
number of derivatives which are actually needed. This will result in
certain ``finite order'' statements. Even though there is also a
well-developed theory of real analytic wave equations and their
solutions we will exclusively stick to the $\Cinfty$- and $\Fun$-case.

Throughout this work, we will avoid techniques from Fourier analysis
and stay exclusively in ``coordinate space''. It is clear that in a
geometric framework there is no intrinsic definition of a
\emph{global} Fourier transform. In principle, one can pass to a
\emph{microlocal} version of Fourier transform between tangent and
cotangent spaces. However, we shall not need this more sophisticated
approach here, even though this will lead to some deeper insights in
the nature of the singularities of the Green operators by means of a
wavefront analysis. As this text should serve as a first reading in
this area we decided to concentrate on the more basic formulations.

%
%

\section*{A User's Guide for Reading}

This text addresses mainly master and PhD students who want to get a
fast but yet detailed access to an important research topic in global
analysis and partial differential equations on manifolds of great
recent interest.  The reader should have some background knowledge in
differential geometry. We use the language of manifolds, vector
bundles, and tensor calculus without further explanations. Some
previous exposure to locally convex analysis and distribution theory
on $\mathbb{R}^n$ might be useful but will not be required: all
relevant notions will either be explained in detail or accompanied
with explicit references to other textbooks for detailed
proofs. Knowledge in Lorentz geometry is of course useful as well, but
we will develop those parts of the theory which are relevant for our
purposes, essentially the notions of causality. We assume that the
reader has at least some vague interest in the physical applications
of the theory as we will take this often as motivation.

The presented material is entirely standard and can also be found in
various other sources. We mainly follow the beautiful exposition of
Bär, Ginoux, and Pfäffle \cite{baer.ginoux.pfaeffle:2007a} but rely
also on the textbooks \cite{guillemin.sternberg:1990a,
  hoermander:2003a, friedlander:1975a} for certain details and further
aspects on distributions on manifolds and geometric wave equations not
discussed in \cite{baer.ginoux.pfaeffle:2007a}. Concerning Lorentz
geometry we refer to the textbook of O'Neill \cite{oneill:1983a} and
the recent review article of Minguzzi and Sanchez
\cite{minguzzi.sanchez:2006a:pre} on the causal structure. Other
resources on Lorentz geometry and general relativity are the classical
texts \cite{hawking.ellis:1973a, wald:1984a, straumann:1988a,
  beem.ehrlich.easley:1996a}. More details on distribution theory and
locally convex analysis can be found in the standard textbooks
\cite{treves:1967a, rudin:1991a, jarchow:1981a}.  For further reading
one should consult the recent booklet \cite{baer.fredenhagen:2009a} as
well as the articles \cite{brunetti.fredenhagen.verch:2003a,
  brunetti.fredenhagen:2000a} for approaches to (quantum) field
theories on curved spacetimes based on the construction of Green
functions for geometric wave equations.  Though we do not touch this
subject, background information on axiomatic approaches to quantum
field theory might be helpful and can be found in the classical
textbooks \cite{streater.wightman:1964a, haag:1993a}.  Beside these
general references we will provide more detailed ones throughout the
text.

\medskip

\noindent
The material is divided into four chapters and two supplementary
appendices:

\medskip

In the first chapter we set the stage for the relevant analysis on
manifolds. In Section~\ref{sec:TestFunctions} we introduce test
function and test section spaces and investigate their locally convex
topologies. The central result will be
Theorem~\ref{theorem:inductive-limit-topology} establishing the LF
topology for compactly supported smooth sections as well as important
properties like completeness of this topology. Moreover, we study
continuous linear maps between test section spaces: on one hand
pull-backs with respect to bundle maps and on the other hand various
multilinear pairings between sections. Finally, we show that the
smooth sections with compact support are sequentially dense in all
other $\Fun$- and $\Fun_0$-sections. Then in
Section~\ref{sec:DifferentialOperators} we discuss differential
operators and their symbols. In particular, we introduce a global
symbol calculus based on the usage of covariant
derivatives. Differential operators are then shown to be continuous
linear maps for the test section spaces. We show that differential
operators have adjoints for various natural pairings and compute the
adjoints explicitly by using the global symbol calculus in
Theorem~\ref{theorem:Neumaier-Theorem}. We arrive in
Section~\ref{sec:Distributions} at the definition of distributions or,
more precisely, of generalised sections. Here we first present the
intrinsic definition. Later on, we interpret distributions always with
respect to a fixed reference density: this way, one can avoid carrying
around the additional density bundle everywhere. We define the
weak$^*$ topology and explain the support and singular support of
generalized sections. Important for later use will be the
characterization of generalized sections with compact support in
Theorem~\ref{theorem:gensec-with-compact-support}. We introduce the
push-forward, the action of differential operators as well as the
external tensor product of distributional sections. Parallel to the
smooth case we develop the $\Fun$-case, both for test sections and
distributions of finite order.

Chapter~\ref{cha:LorentzGeometry} contains a rough overview on Lorentz
geometry where we focus on particular topics rather than on a general
presentation. In Section~\ref{sec:SemiRiemannianManifold} we recall
some basic concepts from semi-Riemannian geometry like parallel
transport and the exponential map of a connection, the Levi-Civita
connection and the d'Alembert operator. Still for general
semi-Riemannian manifolds we introduce the notion of a connection
d'Alembertian and provide a definition and characterization of
normally hyperbolic differential operators. We pass to true Lorentz
geometry in Section~\ref{sec:CausalStructure} where we mainly focus on
aspects related to the causal structure. As motivation, also for the
wave equations, we recall some features of general relativity. This
gives us the notions of time orientability, causality, and ultimately,
of Cauchy hypersurfaces. Here we discuss the characterization of
globally hyperbolic spacetimes by the existence of smooth Cauchy
hypersurfaces in Theorem~\ref{theorem:globalhyp-and-cauchy-surface}
and present some important consequences of this ``splitting
theorem''. Throughout this section our proofs are rather sketchy but
illustrated by simple geometric (counter-) examples. Even without
explicit proofs this should help to develop the right intuition.  We
conclude this chapter with some general remarks on wave equations, the
Cauchy problem, and advanced and retarded Green functions in
Section~\ref{sec:CauchyProblemGreenFunctions}.

Even though Chapter~\ref{cha:LocalTheory} deals with the local
construction of Green functions we need already here geometric
concepts like parallel transport and the exponential map. As warming
up we start in Section~\ref{sec:dAlembertOnMinkowskiSpaceTime} with
the wave equation \eqref{eq:TheWaveEquation} on flat Minkowski
spacetime and obtain the advanced and retarded Green functions by
constructing an entirely holomorphic family $\{R^\pm(\alpha)\}_{\alpha
  \in \mathbb{C}}$ of distributions, the \emph{Riesz
  distributions}. For $\alpha = 2$ one obtains the Green functions of
$\dAlembert$. We examine these Riesz distributions in great detail as
they will be the crucial tool to construct local Green functions in
general. The case of spacetime dimensions $n = 1$ (only time) and $n =
2$ is discussed explicitly as one obtains a drastically simpler
approach here.  In Section~\ref{sec:RieszOnConvexDomain} we use the
exponential map to transfer the Riesz distributions also to the curved
situation, at least in a small normal neighborhood of a given
point. However, the curvature will now cause slightly different
features of the Riesz distributions which results in the failure of
$R^\pm(p, 2)$ being a Green function of the scalar d'Alembert
operator. Nevertheless, the defect can be computed explicitly enough
to use the Riesz distributions in
Section~\ref{sec:HadamardCoefficients} to formulate an heuristic
Ansatz for the true Green function, now for a general normally
hyperbolic differential operator, as a series expansion in the
``degree of singularity''. This Ansatz leads to transport equations
similar to the WKB approximation whose solutions will be the
\emph{Hadamard coefficients}. Even though working on a small
coordinate patch the construction of the Hadamard coefficients in
Theorem~\ref{theorem:hadamard-coefficients} requires the full
machinery of differential geometry and would be hard to understand
without the usage of covariant derivatives and their parallel
transports. As an application of this general approach we compute the
Hadamard coefficients for the Klein-Gordon equation in flat spacetime
explicitly and obtain an explicit formula for the advanced and
retarded Green functions in
Theorem~\ref{theorem:green-function-for-klein-gordon}. Back in the
general situation we show in the rather technical
Section~\ref{satz:LocalFundamentalSolution} how a true Green function
with good causal properties can be obtained from the Hadamard
coefficients. Here one first enforces the convergence of the above
Ansatz thereby destroying the property of a Green function. The result
is a parametrix which can be modified in a second step to obtain the
Green functions in Theorem~\ref{theorem:local-green-functions}. As a
first application we use the local Green functions to construct
particular solutions of the inhomogeneous wave equation for
distributional and smooth inhomogeneities in
Section~\ref{satz:SolvingWaveEqLocally} in
Theorem~\ref{theorem:dual-of-fund-solution-gives-inhom-solution}.

Chapter~\ref{satz:global-theory} is now devoted to the global
situation. First we have to recall the notion of the time separation
on a Lorentz manifold in Section~\ref{satz:uniqueness-properties}
which is then used to prove uniqueness of solutions in
Theorem~\ref{theorem:future-past-comp--hom-solution-is-zero} with
either future or past compact support provided the global causal
structure is well-behaved enough. Section~\ref{satz:cauchy-problem}
contains the precise formulation of the global Cauchy problem as well
as its solution for globally hyperbolic spacetimes. We discuss both
the smooth situation as well as certain finite differentiability
versions of the Cauchy problem in
Theorem~\ref{theorem:global-solutions}. The continuous dependence on
the initial values in the Cauchy problem follows from general
arguments using the open mapping theorem. This feature is then used in
Section~\ref{satz:global-fund-and-green-operators} to obtain global
Green functions and the corresponding global Green operators. The
difference of the advanced and retarded Green operator provides an
``inverse'' to the wave operator in the sense of a specific exact
sequence discussed in Theorem~\ref{theorem:propagator-complex}.
Moreover, it constitutes the core ingredient for the classical Poisson
algebra of the field theory corresponding to the wave equation as
discussed then in Section~\ref{satz:poisson-algebra}. We give two
alternative definitions of the Poisson algebra: one as polynomial
algebra on the initial conditions depending on the choice of the
Cauchy hypersurface with the canonical ``symplectic'' Poisson
bracket. The other version is obtained as quotient of the polynomial
algebra on all field configurations with Poisson bracket coming from
the Green operators. The equivalence of both is shown in
Theorem~\ref{theorem:the-poisson-isomorphisms} and gives an easy proof
of the ``time-slice'' axiom of the classical field theory in
Theorem~\ref{theorem:time-slice-axiom}, analogously to the quantum
field theoretic formulation. Also a classical analog of the
``locality'' axiom is proved in
Theorem~\ref{theorem:local-net-of-observables}.

Appendix~\ref{cha:parall-transp-jacobi} contains background
information on parallel transports and the Taylor expansion of various
geometric objects like the exponential map and the volume density. In
Appendix~\ref{satz:stokes-theorem} we recall some basic applications
of Stokes' theorem.

The text does not contain exercises. However, it is understood that
students who really want to learn these topics in a profound way have
to delve deep into the text. Some of the proofs are sketched and
require some extra thoughts, others contain rather long computations
which can and should be repeated.


%% file: chap1.tex
%
%

\chapter{Distributions and Differential Operators on Manifolds}
\label{cha:DistributionDiffOps}

In this chapter we discuss the basic ingredients for analysis on
smooth manifolds: first we introduce the canonical locally convex
topologies for the smooth functions (with compact support) on $M$ as
well as for smooth sections of vector bundles. These spaces will
constitute the spaces of test functions and test sections,
respectively. We have to discuss convergence of test functions as well
as the completeness of the test function spaces. In a second step we
consider differential operators acting on test functions and test
sections. After discussing elementary algebraic and topological
properties we compute the adjoint of a differential operator with
respect to a given positive density explicitly: here a symbol calculus
is introduced and basic properties are shown. Finally, we introduce
distributions as the continuous linear functionals on the various test
function spaces. This allows to dualize all operations on test
functions in an appropriate way. In particular, differential operators
will act on distributions as well. We discuss the module structure of
distributions, give first basic examples and define the support, and
singular support of distributions.

%
%

\section{Test Functions and Test Sections}
\label{sec:TestFunctions}

\input{testfunction}

%
%

\section{Differential Operators}
\label{sec:DifferentialOperators}

\input{diffops}

%
%

\section{Distributions on Manifolds}
\label{sec:Distributions}

\input{distribution}


%% file: testfunction.tex
%
%

A good understanding of the topological properties of test sections of
vector bundles is crucial. The manifold $M$ will be $n$-dimensional.
In the following, we shall use Einstein's summation convention: the
summation over dual pairs of indexes in multilinear expressions is
automatic.

%
%

\subsection{The Locally Convex Topologies of Test Functions and Test
  Sections} 
\label{subsec:locally-convex-topologies}

In this subsection, we give several different but equivalent
descriptions of the locally convex topology of test functions and test
sections.  Let $E \longrightarrow M$ be a vector bundle of rank $N$.
The first collection of seminorms is obtained as follows. For a chart
$(U, x)$ we consider a compact subset $K \subseteq U$ together with a
collection $\{e_\alpha\}_{\alpha=1, \ldots, N}$ of local sections
$e_\alpha \in \Secinfty(E\at{U})$ such that $\{e_\alpha(p)\}_{\alpha =
  1, \ldots, N}$ is a basis of the fiber $E_p$. We always assume that
$U$ is sufficiently small or e.g. contractible such that local base
sections exist. The collection $\{e_\alpha\}_{\alpha = 1, \ldots, N}$
will also be called a \emph{local frame}\index{Frame}.  The dual frame
will then be denoted by $\{e^\alpha\}_{\alpha = 1, \ldots, N}$ where
$e^\alpha \in \Secinfty(E^*\at{U})$ are the local sections with
$e^\alpha(e_\beta) = \delta^\alpha_\beta$. For $s \in \Secinfty(E)$ we
have unique functions $s^\alpha = e^\alpha(s) \in \Cinfty(U)$ such
that
\begin{equation}
    \label{eq:local-form-sections}
    s\at{U} = s^\alpha e_\alpha.
\end{equation}
We define the seminorms
\begin{equation}
    \label{eq:p-K-l-Seminorm}
    \seminorm[U, x, K, \ell, \{e_\alpha\}](s)
    = \sup_{\substack{p\in K \\
        |I| \leq \ell \\
        \alpha = 1, \ldots, N}}
    \left|
        \frac{\partial^{|I|} s^\alpha}{\partial x^I}(p)
    \right|,
\end{equation}
where $I = (i_1, \ldots, i_n) \in \mathbb{N}_0^n$ denotes a multiindex
of total length $|I| = i_1 + \cdots + i_n$. Clearly, the seminorm
depends on the choice of the chart, the compactum, the integer $\ell
\in \mathbb{N}_0$ as well as on the choice of the local base
sections. In case we have just functions, i.e. sections of the trivial
vector bundle $E = M \times \mathbb{C}$, we can use the canonical
trivialization which results in the simpler form
\begin{equation}
    \label{eq:p-K-l-Seminorm-for-functions}
    \seminorm[U, x, K, \ell](f)
    =
    \sup_{\substack{p \in K \\
        |I|\leq \ell}}
    \left|
        \frac{\partial^{|I|} f}{\partial x^I}(p) 
    \right|
\end{equation}
of the seminorm for $f \in \Cinfty(M)$.
\begin{lemma}
    \label{lemma:pKl-is-Seminorm}
    \index{Seminorm!for sections}%
    For all choices of a chart $(U, x)$, a compact subset $K \subseteq
    U$, an integer $\ell \in \mathbb{N}_0$ and local base sections
    $\{e_\alpha\}$ of $E$ on $U$, the map
    \begin{equation}
        \label{eq:pKl-as-Map}
        \seminorm[U, x, K, \ell, \{e_\alpha\}]:
        \Secinfty(E\at{U}) \longrightarrow \mathbb{R}_0^+
    \end{equation}
    is a well-defined seminorm.
\end{lemma}
\begin{proof}
    Clearly, the supremum over $K$ is finite as all partial
    derivatives are continuous. The remaining properties of a seminorm
    are checked easily.
\end{proof}

An alternative construction of seminorms is as follows. On $E
\longrightarrow M$ we choose a covariant derivative $\nabla^E$ and on
$TM \longrightarrow M$ a torsion-free covariant derivative $\nabla$,
e.g. the Levi-Civita connection for some (semi-) Riemannian metric.
Moreover, on $E$ we choose a Riemannian fiber metric if $E$ is a real
vector bundle or a Hermitian fiber metric if $E$ is complex,
respectively. Finally, we shall use a Riemannian metric on $M$. Then
the two metric structures give rise to fiber metrics on all bundles
constructed out of $TM$ and $E$ via tensor products etc.  Moreover, we
have the following operator of symmetrized covariant differentiation:
\begin{definition}[Symmetrized covariant differentiation]
    \label{definition:Symmetric-Covariant-Differentiation}
    \index{Covariant derivative}%
    \index{Symmetrized covariant derivative}%
    Let $\nabla^E$ be a covariant derivative for a vector bundle $E
    \longrightarrow M$ and let $\nabla$ a torsion-free covariant
    derivative on $M$. Then
    \begin{equation}
        \label{eq:Symmetric-Covariant-Derivative}
        \SymD^E:
        \Secinfty(\Sym^k T^*M \tensor E)
        \longrightarrow 
        \Secinfty(\Sym^{k+1} T^*M \tensor E)
    \end{equation}
    is defined by
    \begin{equation}
        \label{eq:Symm-Cov--Der-Formula}
        \SymD^E(\alpha \tensor s) (X_1, \ldots, X_{k+1})
        =
        \sum_{\ell=1}^{k+1}
        \left(
            \nabla_{X_\ell}\alpha \tensor s 
            + \alpha \tensor \nabla_{X_\ell}^E s
        \right)
        (X_1, \ldots, \stackrel{\ell}{\wedge}, \ldots, X_{k+1}),
    \end{equation}
    where $\alpha \in \Secinfty(\Sym^k T^*M)$, $s \in \Secinfty(E)$,
    and $X_1, \ldots, X_{k+1} \in \Secinfty(TM)$.
\end{definition}
\begin{proposition}
    \label{proposition:Symmetric-Covariant-Derivative}
    The operator $\SymD^E$ is linear, well-defined, and satisfies the
    following properties:
    \begin{propositionlist}
    \item \label{item:SymDisdeRhamd} For $E = M \times \mathbb{C}$
        with the canonical flat covariant derivative and $f \in
        \Cinfty(M)$ we have
        \begin{equation}
            \label{eq:flat-symm-cov-derivative}
            \SymD f=\D f.
        \end{equation}
    \item \label{item:SymDProductRule} For $\alpha \in
        \Secinfty(\Sym^k T^*M)$ and $\beta \tensor s \in
        \Secinfty(\Sym^\ell T^*M \tensor E)$ we have
        \begin{equation}
            \label{eq:Symm-Cov-Der-Product-Rule}
            \SymD^E\left((\alpha \vee \beta) \tensor s\right)
            =
            \left(\SymD \alpha \vee \beta\right) \tensor s
            +
            \alpha \vee \SymD^E(\beta \tensor s).
        \end{equation}
    \item \label{item:SymDLocalForm} Locally in a chart $(U,x)$ we have
        \begin{equation}
            \label{eq:Symm-Cov-Der-Local-Form}
            \SymD^E(\alpha\otimes s)\At{U}
            =
            \left(
                \D x^i \vee \nabla_{\frac{\partial}{\partial x^i}}
                \alpha
            \right)
            \tensor s
            +
            \D x^i \vee \alpha \tensor 
            \nabla^E_{\frac{\partial}{\partial x^i}}s.
        \end{equation}
    \end{propositionlist}
\end{proposition}
\begin{proof}   
    Clearly, \eqref{eq:Symm-Cov--Der-Formula} gives a well-defined
    $E$-valued symmetric $(k+1)$-form. On the trivial line bundle the
    flat connection is $\nabla_X f = \Lie_X f = (\D f)(X)$ from
    which \eqref{eq:flat-symm-cov-derivative} is obvious. The Leibniz
    rule \eqref{eq:Symm-Cov-Der-Product-Rule} is a direct consequence
    from \eqref{eq:Symm-Cov-Der-Local-Form} but can also be obtained
    in a coordinate free way. We prove
    \eqref{eq:Symm-Cov-Der-Local-Form} by an explicit computation.
    \begin{align*}
        &\left(\SymD^E(\alpha \tensor s)\right)(X_1, \ldots, X_{k+1}) 
        \\
        &\qquad=
        \sum_{\ell=1}^{k+1}
        \left(
            \nabla_{X_\ell}\alpha \tensor s 
            + \alpha \tensor \nabla_{X_\ell}^E s
        \right)
        (X_1, \ldots, \stackrel{\ell}{\wedge}, \ldots, X_{k+1}) \\
        &\qquad=
        \sum_{\ell=1}^{k+1}
        \left(
            \D x^i(X_\ell)\nabla_{\frac{\partial}{\partial x^i}}\alpha 
            \tensor s
            +
            \D x^i(X_\ell) \alpha \tensor
            \nabla_{\frac{\partial}{\partial x^i}}^E s 
        \right)
        (X_1, \ldots, \stackrel{\ell}{\wedge}, \ldots, X_{k+1})\\
        &\qquad=
        \left(
            \D x^i\vee 
            \left(
                \nabla_{\frac{\partial}{\partial x^i}} \alpha 
                \tensor s 
                +
                \alpha \tensor
                \nabla^E_{\frac{\partial}{\partial x^i}} s
            \right)
        \right)
        (X_1, \ldots, X_{k+1}).
    \end{align*}
\end{proof}
Using this symmetrized covariant differentiation we can construct a
seminorm for $s \in \Secinfty(E)$ as follows. First we consider
$(\SymD^E)^\ell s \in \Secinfty(\Sym^\ell T^*M \tensor E)$. Then we
can use the fiber metric $h$ on $\Sym^\ell T^*M \tensor E$ to get a
fiberwise norm $\norm{\argument}_h$. Then for every compact subset $K
\subseteq M$ we consider
\begin{equation}
    \label{eq:PkL-with-Covariant}
    \seminorm[K, \ell](s)
    =
    \sup_{p \in K}\norm{(\SymD^E)^\ell s\at{p}}_h,
\end{equation}
where we suppress the dependence of $\seminorm[K, \ell]$ on the
choices of $\nabla$, $\nabla^E$ and $h$ to simplify our notation.
\begin{lemma}
    \label{lemma:PKl-with-Covariant-is-Norm}
    \index{Seminorm}%
    For all choices of a compactum $K \subseteq M$ and $\ell \in
    \mathbb{N}_0$ the map
    \begin{equation}
        \label{eq:PKl-with-Covariant-as-Map}
        \seminorm[K, \ell]: \Secinfty(E) \longrightarrow \mathbb{R}_0^+
    \end{equation}
    is a well-defined seminorm.
\end{lemma}
\begin{proof}
    Thanks to the continuity of $\norm{(D^E)^\ell s}$ the supremum is
    actually a maximum over the compact subset $K$. Thus $\seminorm[K,
    \ell](s) \in \mathbb{R}_0^+$ is finite. The remaining properties
    of a seminorm follow at once.
\end{proof}

We can now use both types of seminorms to construct locally convex
topologies for $\Secinfty(E)$. Since neither the system of the
$\seminorm[U, x, K, \ell, \{e_\alpha\}]$ nor the $\seminorm[K, \ell]$
are filtrating, we have to take maximums over finitely many of them in
each of the following cases:
\begin{compactitem}
    \index{Seminorm!filtrating system}
\item[\textbf{A}] Choose an atlas ${(U,x)}$ with local base sections
    $\{e_\alpha\}$ on each chart and consider \emph{all} seminorms
    $\seminorm[U, x, K, \ell, \{e_\alpha\}]$ arising from the charts
    of this atlas, all $\ell \in \mathbb{N}_0$, and all compact
    subsets $K \subseteq U$.
\item[\textbf{B}] Choose $\nabla$, $\nabla^E$ and fiber metrics and
    consider \emph{all} seminorms $\seminorm[K, \ell]$ arising from
    all compact subsets $K\subseteq M$ and all $\ell \in\mathbb{N}_0$.
\end{compactitem}
As a slight variation of \textbf{A} we can also consider the locally
convex topology where we only take countably many compacta:
\begin{compactitem}
\item[\textbf{A'}] Take only at most countably many charts and in each
    chart $(U, x)$ only an exhausting sequence $\ldots \subseteq K_n
    \subseteq \mathring{K}_{n+1}\subseteq K_{n+1}\subseteq \ldots
    \subseteq U$ of compacta.
\end{compactitem}
Analogously we can use only an exhausting sequence of compacta in the
second version:
\begin{compactitem}
\item[\textbf{B'}] Take the $\seminorm[K, \ell]$ seminorms for an
    exhausting sequence $\ldots \subseteq K_n \subseteq
    \mathring{K}_{n+1} \subseteq K_{n+1}\subseteq \ldots \subseteq M$
    of $M$ by compacta.
\end{compactitem}

Note that for second countable manifolds we can indeed find a
countable atlas together with a choice of countably many compacta,
each contained in a chart, which cover the whole manifold $M$.
\begin{theorem}
    \label{theorem:Ck-Topologie}
    \index{Locally convex topology!for sections}%
    \index{Frechet space@Fr{\'e}chet space}%
    Let $E \longrightarrow M$ be a vector bundle over $M$.
    \begin{theoremlist}
    \item The four locally convex topologies induced by the choices
        \textbf{A}, \textbf{B}, \textbf{A'}, and \textbf{B'} of
        seminorms coincide. Thus $\Secinfty(E)$ has an intrinsic
        locally convex topology not depending on any of the above
        choices.
    \item $\Secinfty(E)$ is a Fr\'echet space with respect to the
        above natural topology.
    \item When restricting to those seminorms with $\ell \leq k$ for a
        fixed $k \in \mathbb{N}_0$, then we obtain natural Fr\'echet
        topologies for $\Sec(E)$.
    \end{theoremlist}
\end{theorem}
\begin{proof}
    First we note that the topologies induced by \textbf{B} and
    \textbf{B'} are the same: indeed \textbf{B} is clearly finer than
    \textbf{B'} as it contains all the seminorms of \textbf{B'}.
    Conversely, we have $\seminorm[K, \ell](s) \leq \seminorm[K',
    \ell](s)$ for $K \subseteq K'$. Now if $K_n$ is an exhausting
    sequence of compacta, then $K \subseteq K_n$ for sufficiently
    large $n$, hence the seminorm $\seminorm[K, \ell]$ can be
    dominated by $\seminorm[K_n, \ell]$.  Thus the induced topologies
    are equivalent.

    For the first version it is clear that \textbf{A} induces a finer
    topology than \textbf{A'} as \textbf{A} contains all seminorms
    from \textbf{A'}. Now let $(U,x)$ be a chart of the chosen atlas
    and $U_n$ the sequence of charts which already cover $M$ which
    works since $M$ is assumed to be second countable.  Moreover, let
    $K_{n, m} \subseteq U_n$ be the exhausting sequence of compacta
    and let $K \subseteq U$ be given. Since $K$ is compact, finitely
    many $U_{n_1}, \ldots, U_{n_k}$ already cover $K$.  Furthermore,
    since the $\mathring{K}_{n, m}$ cover $U_n$, already finitely many
    $\mathring{K}_{n, m}$ cover $K$.  Thus the compactum $K$ is
    covered by finitely many of the $K_{n, m}$'s.  From the chain rule
    it is clear that there are smooth functions $\Phi_{IJ}\in\Cinfty
    (U\cap\widetilde{U})$ such that for $|I|\leq \ell$
    \[
    \frac{\partial^{|I|}f}{\partial x^I} \At{U \cap \widetilde{U}} =
    \sum_{|J|\leq \ell} \Phi_{IJ}
    \frac{\partial^{|J|}f}{\partial\widetilde{x}^J} \At{U \cap
      \widetilde{U}}
    \]
    on the overlap of two charts $(U, x)$ and $(\widetilde{U},
    \widetilde{x})$.  In fact, the $\Phi_{IJ}$ are certain polynomials
    in the partial derivatives of the Jacobian of the coordinate
    change.  It follows that there is a constant $c$ with
    \[
    \seminorm[U, x, K, \ell](f) \leq c
    \seminorm[\widetilde{U},\widetilde{x},K,\ell](f)
    \]
    for all $f \in \Cinfty(U\cap\widetilde{U})$ and $K \subseteq U
    \cap \widetilde{U}$ compact. The constant depends on $U$, $x$,
    $K$, $\ell$, $\widetilde{U}$, and $\widetilde{x}$ but not on $f$.
    The precise form of $c$ is irrelevant, it can be obtained from the
    maximum of the $\Phi_{IJ}$ over $K$ where the $\Phi_{IJ}$ can be
    obtained recursively from the chain rule. From this we see that
    \[
    \seminorm[U,x,K,\ell] 
    \leq \max_{n,m}  c_{n, m} \seminorm[U_n,x_n,K_{n,m},\ell],
    \]
    where the maximum is taken over the finitely many $n, m$ such that
    the $K_{n,m}$ cover $K$. This shows that the topology induced by
    \textbf{A'} is finer than the one obtained by \textbf{A}. Thus all
    together, they coincide.

    Finally, let $\ell \in \mathbb{N}_0$ be given. By induction and
    the local expressions
    \[
    \nabla_{\frac{\partial}{\partial x^i}}^E e_\alpha 
    =
    A_{i\alpha}^\beta e_\beta
    \quad
    \textrm{and}
    \quad
    \nabla_{\frac{\partial}{\partial x^i}} \D x^j
    =
    -\Gamma_{ik}^j \D x^k
    \]
    with the connection one-forms\index{Connection one-form} and
    Christoffel symbols\index{Christoffel symbol} of $\nabla^E$ and
    $\nabla$, respectively, we see that there exist smooth functions
    $a^J_{i_1 \ldots i_\ell} {}_\alpha^\gamma \in \Cinfty(U)$ such
    that
    \[
    (\SymD^E)^\ell s\At{U}
    =
    \sum_{|J|\leq \ell}
    \frac{1}{\ell!}
    a^J_{i_1 \ldots i_\ell} {}_\alpha^\gamma 
    \D x^{i_1} \vee \cdots \vee \D x^{i_l} \tensor
    e_\gamma
    \frac{\partial^{|J|}s^\alpha}{\partial x^J}.
    \tag{$*$}
    \]
    The precise form of the $a^J_{i_1 \ldots i_\ell} {}_\alpha^\gamma$
    is irrelevant, they can be obtained recursively as polynomials in
    the partial derivatives of the $A_{i\alpha}^\beta$ and
    $\Gamma_{ij}^k$. Moreover, for the term with highest derivatives,
    i.e. where $|J| = \ell$, we have the following explicit expression
    \[
    (\SymD^E)^\ell s \At{U}
    =
    \D x^{i_1}\vee \cdots \vee\D x^{i_\ell} \tensor e_\alpha
    \frac{\partial^\ell s^\alpha}
    {\partial x^{i_1} \cdots \partial x^{x_\ell}}
    +
    (\textrm{lower order terms}).
    \tag{$**$}
    \]
    This can easily be obtained by induction since the difference
    between partial derivatives and covariant derivatives is given by
    additional terms involving the $A_{i\alpha}^\beta$ and
    $\Gamma_{ij}^k$. But these terms do not involve derivatives of the
    functions $s^\alpha$.  Now let $K$ be a compactum. Then we find
    finitely many compacta $K_n\subseteq U_n$ contained in charts
    $(U_n,x_n)$ such that the $K_n$ cover $K$. In each chart $(U,x)$
    we see that there are constants $c_U > 0$ with
    \[
    \norm{(\SymD^E)^\ell s\at{p}}
    \leq c_U \max_{\substack{|J| \leq \ell \\ \alpha}} 
    \left|
        \frac{\partial^{|J|}s^\alpha}{\partial x^J}(p)
    \right| 
    \qquad
    \textrm{for}
    \quad
    p \in U,
    \]
    where $c_U$ is obtained from the maximum of the $a^J_{i_1 \cdots
      i_\ell} {}_\alpha^\gamma$ and the norms of the $\D x^{i_1}\vee
    \cdots \vee\D x^{i_\ell} \tensor e_\gamma$ with respect to the
    chosen fiber metrics according to ($*$). But this shows that
    \[
    \seminorm[K, \ell] 
    \leq c \max_{n} \seminorm[U_n,x_n,K_n, \ell,\{e_\alpha\}],
    \]
    where the maximum is taken over the finitely many $n$ such that
    $K_n$ cover $K$ and $c = \max_{n} c_{U_n}$. This shows that the
    topology induced by \textbf{A} is finer than the one induced by
    \textbf{B}. Conversely, given a $\seminorm[U, x, K, \ell,
    \{e_\alpha\}]$ we see from ($**$) that we can estimate the partial
    derivatives $\frac{\partial^{|J|}s^\alpha}{\partial x^J}$ with
    $|J| \leq \ell$ by norms of $(\SymD^E)^\ell s$ and norms of
    partial derivatives $\frac{\partial^{|J'|}s^\alpha}{\partial
      x^{J'}}$ for $|J'|< \ell$.  By induction on $\ell$ we conclude
    that we can estimate the partial derivatives
    $\frac{\partial^{|J|}s^\alpha}{\partial x^J}$ with $|J|= \ell$ by
    norms of $(\SymD^E)^{\ell'} s$ with $\ell' \leq \ell$. Since the
    relative coefficient functions are all smooth this gives a
    constant $c>0$ such that
    \[
    \seminorm[U, x, K, \ell, \{e_\alpha\}] (s) 
    \leq c \max_{\ell' \leq \ell} \seminorm[K, \ell'] (s).
    \]
    This shows that the topology induced by \textbf{B} is finer than
    the one induced by \textbf{A}. Thus, we have shown that all four
    topologies coincide. Since the version \textbf{A} does not depend
    on the choices of $\nabla^E$, $\nabla$ and the fiber metrics and
    since the version \textbf{B} does not depend on an atlas and local
    trivializations we see that the topology itself does not depend on
    any of the chosen data. Note however, that the particular systems
    of seminorms certainly do depend on these choices, only the
    resulting topology is independent.

    For the second part, we first notice that the topology is
    certainly Hausdorff: the seminorms $\seminorm[\{p\}, 0]$ with $p
    \in M$ are already separating. Moreover, the versions \textbf{A'}
    and \textbf{B'} consist of countably many seminorms which define
    the topology. Here it is crucial to have second countable
    manifolds. Thus we only have to show completeness and thanks to
    the countably many seminorms we only have to consider Cauchy
    sequences and not general Cauchy nets. Thus let $s_n \in
    \Secinfty(E)$ be a Cauchy sequence with respect to e.g.
    \textbf{A}.  Taking $K = \{p\}$ a point and $\ell = 0$ we see that
    the sequence $s_n^\alpha(p) \in \mathbb{C}$ (or $\mathbb{R}$) is a
    Cauchy sequence and hence a convergent sequence. Thus $s_n(p)
    \longrightarrow s(p) \in E_p$ for some unique vector $s(p)$. This
    shows that there is a section $s: M \longrightarrow E$ of which we
    have to show smoothness. However, smoothness is a local concept
    which we can check in a local chart. But then the seminorms
    $\seminorm[U, x, K, \ell, \{e_\alpha\}]$ just define the usual
    $\Cinfty$-topology of functions on $U$ with values in
    $\mathbb{R}^N$ or $\mathbb{C}^N$, respectively, via the
    trivialization $\{e_\alpha\}$. Hence we conclude that all
    functions $s^\alpha = e^\alpha (s)$ are smooth and thus $s \in
    \Secinfty(E)$ is a smooth section everywhere since by \textbf{A}
    we can cover the whole manifold with charts $(U, x)$.  Again we
    can argue locally to show that $s_n \longrightarrow s$ in the
    sense of \textbf{A}. This shows that $\Secinfty(E)$ is
    (sequentially) complete which gives the second part. The third
    part is clear, we have shown the most difficult part $k = +\infty$
    already.
\end{proof}

In the following, we shall always endow $\Secinfty(E)$ as well as
$\Sec(E)$ with these naturally defined topologies.
\begin{definition}[$\Cinfty$-Topology]
    \label{definition:CinftyTopology}
    \index{Locally convex topology!Cinfinity-topology@$\Cinfty$-topology}%
    \index{Locally convex topology!Ck-topology@$\Fun$-topology}%
    The natural Fr\'echet topology of $\Secinfty(E)$ is called the
    $\Cinfty$-to\-po\-lo\-gy. Analogously, we call the natural
    Fr\'echet topology of $\Sec(E)$ the $\Fun$-topology.
\end{definition}
\begin{remark}[$\Cinfty$-Topology]
    \label{remark:CinftyTopology}
    ~
    \begin{remarklist}
    \item A sequence $s_n\in\Secinfty(E)$ converges to $s$ with
        respect to the $\Cinfty$-topology if and only if $s_n$
        converges uniformly on all compact subsets of $M$ with all
        derivatives to $s$. Similar, the convergence in the
        $\Fun$-topology is the locally uniform convergence in the first
        $k$ derivatives.
    \item If $M$ is compact, we can use $K=M$ in the seminorms of
        \textbf{A} and \textbf{B}. This shows that the $\Fun$-topology
        is even a \emIndex{Banach topology} since we can also take the
        maximum $0 \leq \ell \leq k$. Thus for this particular case,
        techniques from Banach space analysis become available.
        However, the $\Cinfty$-topology is not Banach, even if $M$ is
        compact. In the non compact situation, none of the
        $\Fun$-topologies is Banach.
    \item The case of smooth functions instead of smooth sections is
        somewhat easier. Here we do not need the additional local base
        sections $\{e_\alpha\}$, hence from \textbf{A} we obtain
        seminorms $\seminorm[U, x, K, \ell]$. In the second version,
        we do not need the additional covariant derivative $\nabla^E$
        nor the fiber metric on $E$ but only $\nabla$ and a Riemannian
        metric on $M$.
    \end{remarklist}
\end{remark}
\begin{remark}
    \label{remark:SymbolicSeminorms}
    In the following we can use either types of seminorms to
    characterize the $\Cinfty$-topology. Since the main importance of
    the seminorms is to control derivatives of order up to $\ell$ on a
    compactum $K$ we shall sometimes symbolically write $\seminorm[K,
    \ell]$ for the seminorms obtained from \emph{either} the maximum
    of some finitely many $\seminorm[U_n, x_n, K_n, \ell,
    \{e_{n\alpha}\}]$ where the $K_n$ are such that they cover $K$
    from the seminorms of type \textbf{A} \emph{or} the maximum of the
    $\seminorm[K, \ell']$ with $\ell' \le \ell$ from the seminorms of
    type \textbf{B}. Clearly, the seminorms $\seminorm[K, \ell]$
    obtained this way specify the topology already completely and are
    filtrating and Hausdorff. It should become clear from the context
    whether we apply these symbolic seminorms or the more concrete
    ones as in \textbf{A} or \textbf{B}.
\end{remark}

On a non compact manifold the space $\Cinfty_0(M)$ is a proper
subspace of all smooth functions $\Cinfty(M)$. Analogously,
$\Secinfty_0(M)$ is a proper subspace of $\Secinfty(M)$ for every
vector bundle of positive rank. The following proposition shows that
we can use sections with compact support to approximate arbitrary
ones.
\begin{proposition}
    \label{proposition:compact-sections-dense-in-smooth-sections}
    \index{Approximation}%
    \index{Dense subspace}%
    For a vector bundle $E\longrightarrow M$ the subspace
    $\Secinfty_0(M)$ of compactly supported sections is dense in
    $\Secinfty(M)$ with respect to the $\Cinfty$-topology.
    Analogously, $\Sec_0(E)$ is dense in $\Sec(E)$ in the
    $\Fun$-topology for all $k \in \mathbb{N}_0$.
\end{proposition}
\begin{proof}
    We choose an exhausting sequence $\ldots K_n \subseteq
    \mathring{K}_{n+1} \subseteq K_{n+1} \subseteq \ldots \subset M$
    of compacta and appropriate functions $\chi_n \in \Cinfty_0(M)$
    with the property
    \[
    \chi_n \at{K_n} = 1
    \quad 
    \textrm{and}
    \quad
    \supp(\chi_n)\subseteq K_{n+1}.
    \]
    Clearly, such $\chi_n$ exists thanks to the $\Cinfty$-version of
    the \Index{Urysohn Lemma}, see
    e.g. \cite[Kor.~A.1.5]{waldmann:2007a}.  Then for $s \in
    \Secinfty(E)$ we define $s_n = \chi_n s \in \Secinfty_0(E)$ and
    have for all $\ell \in \mathbb{N}_0$
    \[
    \seminorm[K_n, \ell] (s - s_m) = 0,
    \]
    for $m \geq n$. This shows that $s_n \longrightarrow s$ in the
    $\Cinfty$-topology. For the $\Fun$-topology the argument is the
    same.
\end{proof}

While on one hand, the above statement will be very useful to
approximate sections by compactly supported sections, it shows on the
other hand that the $\Cinfty$-topology is not appropriate for
$\Secinfty_0(M)$ as this subspace is \emph{not complete} in the
$\Cinfty$-topology. Thus we are looking for a \emph{finer} locally
convex topology which makes $\Secinfty_0(M)$ complete. The
construction is based on the following observation:
\begin{lemma}
    \label{lemma:sections-with-compact-support-in-fixed-compacta}
    \index{Locally convex topology!Ck-topology@$\Fun$-topology!closed subspace}%
    Let $A \subseteq M$ be a closed subset and let
    \begin{equation}
        \label{eq:SecADef}
        \Sec_A(E) =
        \left\{
            s \in \Sec(E) \big| \supp(s) \subseteq A
        \right\}.
    \end{equation}
    Then $\Sec_A(E) \subseteq \Sec(E)$ is a closed subspace with
    respect to the $\Fun$-topology for all $k \in \mathbb{N}_0 \cup
    \{+\infty\}$.
\end{lemma}
\begin{proof}
    Since we are in a Fr\'echet situation, it is sufficient to
    consider sequences in order to approach the closure. Thus let $s_n
    \in \Sec_A(E)$ with $s_n \longrightarrow s \in \Sec(E)$ be given.
    Since $\Fun$-convergence implies pointwise convergence we see that
    for $p \in M \backslash A$
    \[
    0 = s_n(p) \longrightarrow s(p),
    \]
    whence $s(p) = 0$. Thus $\supp(s) \subseteq A$ as desired and $s
    \in \Sec_A(E)$ follows.
\end{proof}

This way, the $\Sec_A(E)$ become Fr\'echet spaces themselves being
closed subspaces of the Fr\'echet space $\Sec(E)$. We call the
resulting topology the $\Fun_A$-topology%
\index{Locally convex topology!CkA-topology@$\Fun_A$-topology}. With
respect to their induced topology, the inclusion maps
\begin{equation}
    \label{eq:inclusion-of-sections-with-fixed-support}
    \Sec_A(E) \hookrightarrow \Sec_{A'}(E)
\end{equation}
for $A \subseteq A'$ are continuous and have closed image. This is
clear as the seminorms $\seminorm[K, \ell]$ needed for $\Sec_{A}(E)$
are also continuous seminorms on $\Sec_{A'}(E)$. Moreover, the induced
topology on $\Sec_A(E)$ by
\eqref{eq:inclusion-of-sections-with-fixed-support} is again the
$\Fun_A$-topology. Thus
\eqref{eq:inclusion-of-sections-with-fixed-support} is an
embedding\index{Embedding} and not just an injective continuous map.
We shall now focus on compact subsets $K \subseteq M$ and choose an
exhausting sequence $K_n$ as before. Then the corresponding sequence
\begin{equation}
    \label{eq:inclusion-sequence-for-inductive-limit-topology}
    \Sec_{K_0}(E) 
    \hookrightarrow \Sec_{K_1}(E)
    \hookrightarrow \cdots \hookrightarrow \Sec_{K_n}(E)
    \hookrightarrow \Sec_{K_{n+1}}(E)
    \hookrightarrow \cdots \hookrightarrow \Sec_0(E)
\end{equation}
allows to endow the ``limit'' $\Sec_0(E)$ with the \emph{inductive
  limit topology}. Since all the inclusions are embeddings and since
we only need countably many compacta, we have a \emph{countable strict
  inductive limit topology} (or \emph{LF topology}) for
$\Sec_0(E)$. By general nonsense on such limit topologies, see
e.g. \cite[Sect.~4.6]{jarchow:1981a}, we obtain the following
characterization of a locally convex topology on $\Sec_0(E)$, which we
call the $\Cinfty_0$-topology:
\begin{theorem}[$\Cinfty_0$-topology]
    \label{theorem:inductive-limit-topology}
    \index{Locally convex topology!inductive limit}%
    \index{Locally convex topology!LF topology}%
    \index{Locally convex topology!CinfinityNull-topology@$\Cinfty_0$-topology}%
    \index{Locally convex topology!CkNull-topology@$\Fun_0$-topology}%
    \index{Support}%
    \index{Locally convex topology!Hausdorff}%
    \index{Locally convex topology!complete}%
    \index{Locally convex topology!sequentially complete}%
    \index{Locally convex topology!first countable}%
    \index{Locally convex topology!metrizable}%
    Let $k \in \mathbb{N}_0 \cup \{+\infty\}$. The inductive limit
    topology on $\Sec_0(E)$ obtained from
    \eqref{eq:inclusion-sequence-for-inductive-limit-topology} enjoys
    the following properties:
    \begin{theoremlist}
    \item \label{item:LFisHausdorffComplete} $\Sec_0(E)$ is a
        Hausdorff locally convex complete and sequentially complete
        topological vector space. The topology does not depend on the
        chosen sequence of exhausting compacta.
    \item \label{item:LFinclusionContinuous} All the inclusion maps
        \begin{equation}
            \label{eq:inclusion-fixed-support-in-compact-support}
            \Sec_K(E) \hookrightarrow \Sec_0(E)
        \end{equation}
        are continuous and the $\Fun_0$-topology is the finest locally
        convex topology on $\Sec_0(E)$ with this property. Every
        $\Sec_K(E)$ is closed in $\Sec_0(E)$ and the induced topology
        on $\Sec_K(E)$ is the $\Fun_K$-topology.
    \item \label{item:LFConvergentSequence} A sequence $s_n \in
        \Sec_0(E)$ is a $\Fun_0$-Cauchy sequence if and only if there
        exists a compact subset $K \subseteq M$ with $s_n \in
        \Sec_K(E)$ for all $n$ and $s_n$ is a $\Fun_K$-Cauchy
        sequence. An analogous statement holds for convergent
        sequences.
    \item \label{item:LFcontinuity} If $V$ is a locally convex vector
        space, then a linear map $\Phi: \Sec_0(E) \longrightarrow V$
        is $\Fun_0$-continuous if and only if each restriction $\Phi
        \at{\Sec_K(E)}: \Sec_K(E) \longrightarrow V$ is $\Fun_K$
        continuous. It suffices to consider an exhausting sequence of
        compacta.
    \item \label{item:LFNotMetrizable} If $M$ is non compact
        $\Sec_0(E)$ is not first countable and hence not metrizable.
    \end{theoremlist}
\end{theorem}
\begin{proof}
    We shall only sketch the arguments and refer to
    \cite[Sect.~4.6]{jarchow:1981a} for details on strict inductive
    limit topologies.  The first part follows from general nonsense on
    countable strict inductive limit topologies since all the
    constituents $\Sec_K(E)$ are Fr\'echet spaces. The second part is
    an alternative characterization of inductive limit topologies.
    Part \refitem{item:LFConvergentSequence} and
    \refitem{item:LFcontinuity} are also general facts on inductive
    limit topologies. The last part follows essentially from
    \Index{Baire's theorem}.
\end{proof}
\begin{remark}
    \label{remark:Ck0-topology}
    In the sequel, we only need the properties
    \refitem{item:LFisHausdorffComplete} --
    \refitem{item:LFcontinuity} of the $\Fun_0$-topology, not its
    precise definition. In fact, it will turn out that the actual
    handling of this rather complicated LF topology is fairly easy. We
    refer to the literature for more background information on LF
    topologies, see e.g.~\cite[Sect.~4.6]{jarchow:1981a} or
    \cite{treves:1967a, koethe:1979a, koethe:1969a}. Of course, we are
    mainly interested in the case $k = \infty$.
\end{remark}
\begin{remark}
    \label{remark:FunNullToFunStetig}
    We also remark that the inclusion maps $\Sec_0(E) \hookrightarrow
    \Sec(E)$ are continuous for all $k \in \mathbb{N}_0 \cup
    \{+\infty\}$.
\end{remark}

%
%

\subsection{Continuous Maps between Test Section Spaces}
\label{subsec:cont-maps-between-test-section-spaces}

In this subsection we shall collect some basic examples of maps
between test function and test section spaces which on one hand have a
geometric origin, and which on the other hand are continuous in the
$\Fun$- and $\Fun_0$-topologies, respectively. We start with the
following situation:
\begin{proposition}
    \label{proposition:pullback-of-functions}
    \index{Pull-back}%
    \index{Continuity!pull-back}%
    Let $\phi: M \longrightarrow N$ be a smooth map. Then the
    pull-back $\phi^*: \Cinfty(N) \longrightarrow \Cinfty(M)$ is a
    continuous linear map with respect to the $\Cinfty$-topology.
\end{proposition}
\begin{proof}
    Let $K \subseteq M$ be a compact subset and $\ell \in
    \mathbb{N}_0$ be given. Moreover, let $(U, x)$ be a chart with $K
    \subseteq U$. Then we consider the compact subset $\phi(K)
    \subseteq N$. This will be covered by finitely many charts $(V,
    y)$ of $N$ and we can assume that already one chart will do the
    job. Then we compute by the chain rule
    \[
        \seminorm[U, x, K, \ell] (\phi^* f)
        = \sup_{\substack{p \in K \\ |I| \leq \ell}}
        \left|
            \frac{\partial^{|I|} \phi^* f}{\partial x^I}(p)
        \right|
        = \sup_{\substack{p \in K \\ |I| \leq \ell}}
        \left|
            \sum_{J \leq I} \Phi_{IJ}(p)
            \frac{\partial^{|J| f}}{\partial y^J}(\phi(p))
        \right|,
    \]
    where again the $\Phi_{IJ}$ are smooth functions on $U$ obtained
    from polynomials in the derivatives of the Jacobi matrix of the
    map $\phi$ with respect to the charts $(V, y)$ and $(U, x)$. Since
    $\phi$ is smooth the maps $\Phi_{IJ}$ turn out to be smooth, too,
    hence on $K$ they are bounded. Moreover, the partial derivatives
    of $f$ on $\phi(K)$ are bounded as well so we finally obtain an
    estimate
    \[
    \seminorm[U, x, K, \ell] (\phi ^* f)
    \leq c \seminorm[V, y, \phi(K), \ell] (f),
    \]
    where the constant $c$ depends on the maxima of the functions
    $\Phi_{IJ}$ over $K$ and thus on $\phi$ but not on $f$. But this
    is the desired continuity.
\end{proof}
\begin{remark}
    \label{remark:pullback-of-Ck}      
    Since in the proof we estimated a seminorm with order of
    differentiation $\ell$ again by a seminorm with order of
    differentiation $\ell$, the statement remains true for a
    $\Fun$-map $\phi: M \longrightarrow N$: the pullback $\phi^*:
    \Fun(N) \longrightarrow \Fun(M)$ is $\Fun$-continuous.
\end{remark}

For functions with compact support the pull-back $\phi^* f$ with an
arbitrary map $\phi: M \longrightarrow N$ will no longer have compact
support in general. Take e.g. any smooth map $\phi: M \longrightarrow
N$ from a non compact manifold $M$ into a compact one, then $\phi^*
1_N = 1_M$ but $1_N \in \Cinfty(N)=\Cinfty_0(N)$ and $1_M \notin
\Cinfty_0(M)$. Thus we need an extra condition to assure that $\phi^*$
maps $\Cinfty_0(N)$ into $\Cinfty_0(M)$:
\begin{definition}[Proper map]
    \label{definition:proper-map}
    \index{Proper map}%
    A smooth map $\phi: M \longrightarrow N$ is called proper if
    $\phi^{-1}(K) \subseteq M$ is compact for all compact $K \subseteq
    N$.
\end{definition}
The above definition makes perfect sense in a general topological
context, the smoothness and the manifold structure of $M, N$ are not
needed. Note that a continuous map maps compact subsets to compact
subsets, but inverse images of compact subsets need not be compact as
the above example shows.
\begin{proposition}
    \label{proposition:pullback-with-proper-map}
    \index{Continuity!pull-back}%
    Let $\phi: M \longrightarrow N$ be a smooth proper map. Then
    \begin{equation}
        \label{eq:pullback-with-proper-map}
        \phi^*: \Cinfty_0(N) \longrightarrow \Cinfty_0(M)
    \end{equation}
    is continuous in the $\Cinfty_0$-topology.
\end{proposition}
\begin{proof}
    Let $K \subseteq N$ be compact. By
    Theorem~\ref{theorem:inductive-limit-topology},
    \refitem{item:LFcontinuity} we have to show that the restriction
    $\phi^* \at{\Cinfty_K(N)}: \Cinfty_K(N) \longrightarrow
    \Cinfty_0(M)$ is continuous. Now $\phi^{-1}(K)$ is compact since
    $\phi$ is proper and thus we know
    \[
    \phi^*: \Cinfty_K(N) \longrightarrow \Cinfty_{\phi^{-1}(K)}(M),
    \]
    since in general $\supp(\phi^* f) = \phi^{-1}(\supp(f))$. The
    proof of Proposition~\ref{proposition:pullback-of-functions} shows
    that the $\seminorm[\phi^{-1}(K), \ell]$-seminorms of the images
    of $\phi^*$ can be estimated by the $\seminorm[K,
    \ell]$-seminorms. Thus $\phi^*$ is continuous. Finally, we know
    that
    \[
    \Cinfty_{\phi^{-1}(K)}(M) \hookrightarrow \Cinfty_0(M)
    \]
    is continuous by Theorem~\ref{theorem:inductive-limit-topology},
    \refitem{item:LFinclusionContinuous}. Thus the criterion for the
    continuity of $\phi^*$ is fulfilled.
\end{proof}
\begin{remark}
    \label{remark:proper-pullback-of-Ck-functions}
    Again, there is a $\Fun_0$-version of this statement since we only
    used the same $\ell$ for the estimation in the proof of
    Proposition~\ref{proposition:pullback-of-functions}.
\end{remark}

In a last step, we shall treat test sections of vector bundles. Let $E
\longrightarrow M$ and $F \longrightarrow M$ be vector bundles. Since
a smooth map $\phi: M \longrightarrow N$ alone does not yield any map
between $\Secinfty(E)$ and $\Secinfty(F)$ by itself, we need a vector
bundle morphism. Recall that a \emph{vector bundle morphism}
\index{Vector bundle morphism}%
$\Phi: E \longrightarrow F$ is a smooth map such that $\Phi$ maps
fibers of $E$ into fibers of $F$ and $\Phi$ is linear on each
fiber. Thus $\Phi$ induces a smooth map $\phi$ such that
\begin{equation}
    \label{eq:VectorBundleMorphism}
    \bfig
    \square|arrb|<800,500>[E`F`M`N;\Phi`\pi_E`\pi_F`\phi]
    \morphism(0,0)|l|/{@{>}@/^1em/}/<0,500>[M`E;\iota_E]
    \efig
\end{equation}
commutes. Indeed, $\phi = \pi_F \circ \Phi \circ \iota_E$, where
$\iota_E: M \longrightarrow E$ denotes the zero section.
\begin{lemma}
    \label{lemma:pullback-of-sections}
    Let $\Phi:E \longrightarrow F$ be a vector bundle morphism and
    $\omega \in \Secinfty(F^*)$. Then
    \begin{equation}
        \label{eq:pullback-of-sections}
        (\Phi^* \omega)\at{p} (s_p) = \omega \at{\phi(p)} (\Phi(s_p))
    \end{equation}
    for $s_p \in E_p$ and $p \in M$ defines a smooth section $\Phi^*
    \omega \in \Secinfty(E^*)$ called the pull-back of $\omega$ by
    $\Phi$.
\end{lemma}
\begin{proof}
    It is easy to check that for $s \in \Secinfty(E)$ the function $p
    \mapsto \Phi^*\omega\at{p} (s(p))$ is smooth whence $\Phi^*\omega$
    is smooth itself. Moreover, $\Phi^*\omega \at{p}: E_p
    \longrightarrow \mathbb{R}$ (or $\mathbb{C}$) is clearly linear
    hence the statement follows.
\end{proof}

The pull-back indeed obeys the usual properties of a pull-back, i.e.
for vector bundle morphisms $E \stackrel{\Phi}{\longrightarrow} F
\stackrel{\Psi}{\longrightarrow} G$ we have
\begin{equation}
    \label{eq:pullback-is-functor}
    (\Psi \circ \Phi)^* = \Phi^* \circ \Psi^* \qquad \textrm{and}
      \qquad (\id_E)^* = \id_{\Secinfty(E^*)}.
\end{equation}
We claim that this gives again continuous maps.
\begin{proposition}
    \label{proposition:pullback-of-sections-is-continous}
    \index{Continuity!pull-back}%
    Let $\Phi: E \longrightarrow F$ be a vector bundle morphism. Then
    $\Phi^*: \Secinfty(F^*) \longrightarrow \Secinfty(E^*)$ is
    continuous with respect to the $\Cinfty$-topology.
\end{proposition}
\begin{proof}
    We first need some local expressions. Let $e_\alpha \in
    \Secinfty(E\at{U})$ and $f_\beta \in \Secinfty(F\at{V})$ be local
    base sections defined over open subsets $U \subseteq M$ and $V
    \subseteq N$. We assume that on $V$ we have local coordinates $y$
    and $x$ on $U \subseteq \phi^{-1} (V)$. By choosing $V$ and $U$
    sufficiently small this is possible. Then $\Phi \at{E|_{U}}$ can
    be written as follows. For $s_p = s^\alpha_p e_\alpha(p) \in E_p$
    there exist coefficients $\Phi_\alpha^\beta (p)$ such that
    \[
    \Phi (s_p) = s^\alpha_p \Phi_\alpha^\beta (p) f_\beta (\phi(p)),
    \]
    since $\Phi(s_p) \in F_{\phi(p)}$ for all $s_p \in E_p$ and $p \in
    M$. The smoothness of $\Phi$ gives the smoothness of the locally
    defined functions $\Phi_\alpha^\beta \in \Cinfty(U)$. Now let
    $\omega \in \Secinfty(F^*)$ be given as
    \[
    \omega \at{V} = \omega_\beta f^\beta,
    \]
    where $f^\beta \in \Secinfty(F^*\at{V})$ are the dual base
    sections of the $f_\beta$ as usual. Then
    \[
    (\Phi^* \omega) (s) \at{U}
    = (\omega \circ \phi) (\Phi(s))
    = \left((\omega_\beta \circ \phi) (f^\beta \circ \phi)\right)
    (\Phi^\gamma_\alpha (f_\gamma \circ \phi) s^\alpha)
    = \phi^* (\omega_\beta) \Phi^\beta_\alpha s^\alpha.
    \]
    Hence we have $\Phi^* \omega \at{U} = \phi^* (\omega_\beta)
    \Phi^\beta_\alpha e^\alpha$.  Now we can estimate
    \begin{align*}
        \seminorm[U,x,K,\ell,\{e_\alpha\}] (\Phi^* \omega)
        &= \sup_{
          \substack{|I|\leq \ell \\
            p \in K \\
            \alpha = 1, \ldots, \rank(E)}
        }
        \left|
            \frac{\partial^{|I|}}{\partial x^I} (\Phi^* \omega)_\alpha
            (p)
        \right| \\
        &=
        \sup_{\substack{|I|\leq \ell \\
            p \in K \\
            \alpha = 1, \ldots, \rank(E)}
        }
        \left|
            \frac{\partial^{|I|}}{\partial x^I}
            \left(\phi^*(\omega_\gamma) \Phi^\gamma_\alpha\right) (p)
        \right| \\
        &\leq
        c \seminorm[V, y, \phi(K), \ell, \{f^\beta\}] (\omega),
    \end{align*}
    by the same kind of computation as for the proof of
    Proposition~\ref{proposition:pullback-of-functions}. The constant
    $c$ involves the maxima of polynomials in the partial derivatives
    of the Jacobi matrix of $\phi$ as well as of $\Phi^\beta_\alpha$,
    again by the chain rule and the Leibniz rule. But then the
    continuity is clear.
\end{proof}
\begin{remark}[Pull-back of sections]
    \label{remark:PullBackSections}
    \index{Pull-back!of sections}%
    ~
    \begin{remarklist}
    \item For the support of $\Phi^*\omega$ we obtain
        \begin{equation}
            \label{eq:supp-of-pullback}
            \supp \Phi^*\omega \subseteq \phi^{-1} (\supp \omega),
        \end{equation}
        which is immediate from the definition. Note that due to
        possible degeneration in the fiberwise maps $\Phi\at{E_p}$ the
        support may be strictly smaller than the right hand side.
    \item Again, for a vector bundle morphism $\Phi: E \longrightarrow
        F$ of class $\Fun$ we obtain a continuous map
        \begin{equation}
            \label{eq:Ck-pullback-sections}
            \Phi^*: \Sec(F^*) \longrightarrow \Sec(E^*)
        \end{equation}
        with respect to the $\Fun$-topology.
    \end{remarklist}
\end{remark}
\begin{example}[Tangent map]
    \label{example:TangentMapIsVBMorph}
    \index{Tangent map}%
    Let $\phi: M \longrightarrow N$ be a smooth map. Then $T\phi: TM
    \longrightarrow TN$ is a smooth vector bundle morphism over
    $\phi$. Thus the pull-back gives $(T\phi)^*: \Secinfty(T^*N)
    \longrightarrow \Secinfty(T^*M)$. Clearly, the pull-back
    $(T\phi)^*$ in the sense of Lemma~\ref{lemma:pullback-of-sections}
    coincides with the usual pull-back $\phi^*$ of one-forms in this
    case. Note that if $\phi$ is $\Fun$ then $T\phi$ is only of class
    $\Fun[k-1]$. Moreover, $T\phi$ extends to vector bundle morphisms
    $(T\phi)^{\tensor r}: \tensor^r TM \longrightarrow \tensor^r TN$
    hence we also obtain pull-backs $\phi^*: \Secinfty(\tensor^r T^*N)
    \longrightarrow \Secinfty(\tensor^r T^*M)$ being continuous linear
    maps with respect to the $\Cinfty$-topology.
\end{example}

The case of compactly supported sections is treated analogously to the
case of $\Cinfty_0(N)$. Using \eqref{eq:supp-of-pullback} we can copy
the proof of Proposition~\ref{proposition:pullback-with-proper-map}
and obtain the following result:
\begin{proposition}
    \label{proposition:pulback-of-compact-sections}
    \index{Continuity!pull-back}%
    Let $\Phi: E \longrightarrow F$ be a vector bundle morphism such
    that the induced map $\phi: M \longrightarrow N$ is proper. Then
    the pull-back
    \begin{equation}
        \label{eq:pullback-of-compact-sections}
        \Phi^*: \Secinfty_0(F^*) \longrightarrow \Secinfty_0(E^*)
    \end{equation}
    is continuous with respect to the $\Cinfty_0$-topology. Analogous
    statements hold for the $\Fun$ case.
\end{proposition}

We conclude this section with yet another type of maps, namely the
module structures and various tensor products.
\begin{proposition}
   \label{proposition:modul-tensor-are-continuous}
   \index{Frechet algebra@Fr{\'e}chet algebra}%
   \index{Frechet module@Fr{\'e}chet module}%
   \index{Continuity!tensor product}%
   \index{Continuity!natural pairing}%
   Let $E \longrightarrow M$ and $F \longrightarrow M$ be vector
   bundles.
   \begin{propositionlist}
   \item \label{item:MultiplicationContinuous} The pointwise
       multiplication
       \begin{equation}
           \label{eq:multiplication-of-functions}
           \Cinfty(M) \times \Cinfty(M) \ni (f, g)
           \; \mapsto fg \; 
           \in \Cinfty(M)
       \end{equation}
       is continuous with respect to the $\Cinfty$-topology, hence
       $\Cinfty(M)$ becomes a Fr\'echet algebra.
   \item \label{item:ModuleStructureContinuous} The module structure
       \begin{equation}
           \label{eq:modul-struct-of-sections}
           \Cinfty(M) \times \Secinfty(E) \ni (f, s)
           \; \mapsto \;
           f\cdot s \in \Secinfty(E)
       \end{equation}
       is continuous with respect to the $\Cinfty$-topology, hence
       $\Secinfty(E)$ becomes a Fr\'echet module over the Fr\'echet
       algebra $\Cinfty(M)$.
   \item \label{item:TensorProductContinuous} The tensor product
       \begin{equation}
           \label{eq:tensor-product-of-sections}
           \Secinfty(E) \times \Secinfty(F) \ni (s, t)
           \; \mapsto \; 
           s \tensor t \in \Secinfty(E \tensor F)
       \end{equation}
       is continuous with respect the $\Cinfty$-topology.
   \item \label{item:NaturalPairingContinuous} The natural pairing
       \begin{equation}
           \label{eq:natural-pairing-of-section-with-dual}
           \Secinfty(E^*) \times \Secinfty(E) \ni (\omega, s)
           \; \mapsto \;
           \omega (s) \in \Cinfty(M)
       \end{equation}
       is continuous with respect to the $\Cinfty$-topology.
   \end{propositionlist}
   Analogous statements hold for the $\Fun$ case.
\end{proposition}
\begin{proof}
    All the above statements rely only on the Leibniz rule for
    differentiation of products. Let $K \subseteq U$ be compact and
    let $x$ be local coordinates on $U$, then
    \[
    \seminorm[U, x, K, \ell] (fg)
    = \sup_{\substack{p \in K \\ |I| \leq \ell}}
    \left|
        \frac{\partial^{|I|}}{\partial x^I} (fg)\at{p}
    \right|
    = \sup_{\substack{p \in K \\ |I| \leq \ell}}
    \left|
        \sum_{J \leq I} \binom{I}{J}
        \frac{\partial^{|J|} f}{\partial x^J}(p) 
        \frac{\partial^{|I-J|} g}{\partial x^{I-J}}(p)
    \right|
    \leq 
    c \seminorm[U, x, K, \ell](f) \seminorm[U, x, K, \ell](g),
    \]
    with a constant only depending on the combinatorics of the
    multinomial coefficients $\binom{I}{J}$ and hence only on $\ell$.
    This shows the first part. Writing out the local expressions for
    all the other parts in terms of coefficient functions and local
    base sections shows that all other parts can be reduced to
    part~\refitem{item:MultiplicationContinuous} and hence the above
    computation.
\end{proof}

\begin{remark}
    \label{remark:modul-tensor-are-cont-in-Ck}
    As usual there are $\Fun_0$-versions of this statement. Moreover,
    we have analogous statements for various multilinear pairings and
    applications of endomorphisms to sections etc.
\end{remark}

%
%

\subsection{Approximations}
\label{Approximations}

In this subsection we shall sketch some approximation results of how
less differentiable functions can be approximated by smooth ones. This
rather technical section will turn out to be useful in many places.
\begin{theorem}
    \label{theorem:compact-and-smooth-sections-dense-in-others}
    \index{Approximation}%
    \index{Dense subspace}%
    Let $E \longrightarrow M$ be a vector bundle. Then
    $\Secinfty_0(E)$ is (sequentially) dense in $\Sec(E)$ for all $k
    \in \mathbb{N}_0$ with respect to the $\Fun$-topology.
\end{theorem}
\begin{proof}
    First we know from
    Proposition~\ref{proposition:compact-sections-dense-in-smooth-sections}
    that $\Sec_0(E)$ is dense in $\Sec(E)$. Thus we only have to show
    that $\Secinfty_0(E)$ is dense in $\Sec_0(E)$ with respect to the
    $\Fun$-topology thanks to the continuous embedding of $\Fun_0(E)$
    into $\Fun(E)$ according to
    Remark~\ref{remark:FunNullToFunStetig}. Let $s \in \Sec_0(E)$ be
    given. Then we choose charts $(U_i, x_i)$ of $M$ together with
    local base sections $e_{\alpha i} \in
    \Secinfty(E\at{U_i})$. Moreover, we choose a partition of unity
    $\varphi_i \in \Cinfty_0(M)$ with $\supp \varphi_i \subseteq U_i$
    being compact and $\sum \varphi_i = 1$.  The compactness of $\supp
    s$ guarantees that finitely many $U_i$ already cover $\supp s$ and
    hence
    \[
    s = \sum\nolimits_i \varphi_i s,
    \]
    with a finite sum. Thus we only have to approximate a $\varphi_i s
    \in \Sec_0(E)$ where $\supp (\varphi_i s) \subseteq U_i$ is in the
    domain of a chart. Now
    \[
    \varphi_i s = s^\alpha_i e_{\alpha i}
    \]
    with $s^\alpha_i \in \Fun_0(U_i)$. From the local theory we know
    that we find functions $s^\alpha_{i m} \in \Cinfty_0(U_i)$ with
    $s^\alpha_{i m} \rightarrow s^\alpha_i$ in the $\Fun$-topology:
    e.g. one can use a convolution\index{Convolution} of the
    $s^\alpha_i$ with a function $\chi_m (x) = m^n \chi (mx)$, where
    $\chi \in \Cinfty_0(\mathbb{R}^n)$ is a function with $\int
    \chi(x) \D^n x = 1$. Then for sufficiently large $m$
    \[
    (s^\alpha_i \circ x_i) * \chi_m \in \Cinfty_0(x(U)),
    \]
    hence
    \[
    s^\alpha_{i m} = ((s^\alpha_i \circ x_i) * \chi_m) \circ x^{-1}
    \in \Cinfty_0(U)
    \]
    is smooth and fulfills $s^\alpha_{i m} \longrightarrow s^\alpha_i$
    in the $\Fun$-topology. For details see
    e.g.~\cite[Thm.~1.3.2]{hoermander:2003a}. Since we can approximate
    each $s^\alpha_i$ we also can approximate $s_i = s^\alpha_i
    e_{\alpha i}$ and thus $s = \sum_i \varphi_i s = \sum_i s_i$ since
    the sums are always finite.
\end{proof}


%% file: diffops.tex
%
%

In this section we introduce differential operators on sections of
vector bundles and discuss their continuity properties with respect to
the various $\Fun$- and $\Fun_0$-topologies.

%
%

\subsection{Differential Operators and their Symbols}
\label{sec:diff-ops-and-symbols}

There are several equivalent definitions of differential operators on
manifolds. We present here the most pragmatic one. Let $E
\longrightarrow M$ and $F \longrightarrow M$ be vector bundles over
$M$.
\begin{definition}[Differential operators]
    \label{definition:Diffops}
    \index{Differential operator}%
    \index{Differential operator!order}%
    Let $D:\Secinfty(E) \longrightarrow \Secinfty(F)$ be a linear
    map. Then $D$ is called differential operator of order $k \in
    \mathbb{N}_0$ if the following conditions are fulfilled.
    \begin{definitionlist}
    \item $D$ can be restricted to open subsets $U \subseteq M$,
        i.e. for any open subset $U \subseteq M$ there exists a linear
        map $D_U: \Secinfty(E\at{U}) \longrightarrow
        \Secinfty(F\at{U})$ such that
        \begin{equation}
            \label{eq:restriction-of-diffop}
            D_U (s\at{U}) = (D s)\at{U}
        \end{equation}
        for all sections $s \in \Secinfty(E)$.
    \item In any chart $(U, x)$ of $M$ and for every local base
        sections $e_\alpha \in \Secinfty(E\at{U})$ and $f_\beta \in
        \Secinfty(F\at{U})$ we have
        \begin{equation}
            \label{eq:local-form-diffop}
            Ds \at{U} = \sum_{r=0}^{k} \frac{1}{r!}
            D_U^{i_1 \ldots i_r} {}_\alpha^\beta
            f_\beta
            \frac{\partial^r s^\alpha}
            {\partial x^{i_1}
              \cdots \partial x^{i_r}}
        \end{equation}
        with locally defined functions $ D_U^{i_1 \ldots i_r}
        {}_\alpha^\beta \in \Cinfty(U)$, totally symmetric in $i_1,
        \ldots, i_r$.
    \end{definitionlist}
\end{definition}
The set of differential operators $D: \Secinfty(E) \longrightarrow
\Secinfty(F)$ of order $k \in \mathbb{N}_0$ is denoted by
$\Diffop^k(E; F)$ and we define
\begin{equation}
    \label{eq:Diffops-all-orders}
    \Diffop^{\bullet}(E; F) = \bigcup_{k=0}^{\infty} \Diffop^k(E; F).
\end{equation}
\begin{remark}[Differential operators]
    \label{remark:Diffops}
    ~
    \begin{remarklist}
    \item Clearly, $\Diffop^k(E; F)$ is a vector space and we have
        \begin{equation}
            \label{eq:diffops-sitting-inside-higher}
            \Diffop^k(E; F) \subseteq \Diffop^{k+1}(E; F)
        \end{equation}
        for all $k \in \mathbb{N}_0$. Thus $\Diffop(E; F)$ is a
        filtered vector space. Note however that
        \eqref{eq:Diffops-all-orders} does not yield a \emph{graded}
        vector space.
    \item The restriction of a differential operator $D$ is important
        since we also want to apply $D$ to sections which are only
        locally defined.
    \item If we are given an atlas of charts and local bases and
        locally defined functions $ D_U^{i_1 \ldots i_r}
        {}_\alpha^\beta$, then we can define a global differential
        operator $D$ by specifying its local form as in
        \eqref{eq:local-form-diffop}, \emph{provided} the functions
        $D_U^{i_1 \ldots i_r} {}_\alpha^\beta$ transform in such a way
        that two definitions agree on the overlap of any two charts in
        that atlas. In fact, the precise transformation law of the $
        D_U^{i_1 \ldots i_r} {}_\alpha^\beta$ is rather complicated
        thanks to the complicated form of the chain rule for multiple
        partial derivatives.
    \item Differential operators are \emph{local}, i.e. $\supp (Ds)
        \subseteq \supp (s)$. \index{Differential operator!locality}
    \end{remarklist}
\end{remark}
\begin{lemma}[Leading symbol]
    \label{lemma:Leading-Symbol}
    \index{Differential operator!leading symbol}%
    If $D: \Secinfty(E) \longrightarrow \Secinfty(F)$ is a
    differential operator of order $k \in \mathbb{N}_0$, locally given
    by \eqref{eq:local-form-diffop}, then the definition
    \begin{equation}
        \label{eq:leading-symbol-definition}
        \sigma_k(D) \at{U} = \frac{1}{k!}
        D_U^{i_1 \ldots i_k} {}_\alpha^\beta
        \frac{\partial}{\partial x^{i_1}} 
        \vee \cdots \vee
        \frac{\partial}{\partial x^{i_k}}
        \tensor f_\beta \tensor e^\alpha
    \end{equation}
    yields a globally well-defined tensor field, called the leading
    symbol of $D$
    \begin{equation}
        \label{eq:leading-symbol-tensor}
        \sigma_k(D) \in \Secinfty(\Sym^k TM \tensor F \tensor E^*).
    \end{equation}
\end{lemma}
\begin{proof}
    This is a straightforward computation since the terms with maximal
    number of derivatives of $s^\alpha$ in
    \eqref{eq:local-form-diffop} transform nicely.
\end{proof}

Note that there is no intrinsic way to define ``sub-leading'' symbols
of a differential operator of order $k \geq 2$. The functions $
D_U^{i_1 \ldots i_r} {}_\alpha^\beta$ do not have a tensorial
transformation law. In fact, terms with different $r$ even mix. This
is also the reason that we can only speak of the maximal number of
partial derivatives appearing in \eqref{eq:local-form-diffop} as
``order''. There is no intrinsic way to characterize differential
operators ``with exactly $k$ partial derivatives'': this would be a
chart dependent statement.

Since canonically $F \tensor E^*\simeq \Hom(E, F)$, we can interpret
the leading symbol $\sigma_k(D)$ also as a section
\begin{equation}
    \label{eq:leading-symbols-as-Hom-section}
    \sigma_k(D) \in \Secinfty(\Sym^k TM \tensor \Hom(E, F)).
\end{equation}

We shall sketch now another, more conceptual approach to differential
operators, see \cite[Def.~16.8.1]{grothendieck:1967a}: it is
essentially based on the observation that for a differential operator
$D$ the commutator $[D,f]$ with a left multiplication by $f \in
\Cinfty(M)$ is a differential operator of at least one order less than
$D$ because of the \Index{Leibniz rule}. We consider an associative,
commutative algebra $\mathcal{A}$ over some ground field
$\mathbb{k}$. Of course, we are mainly interested in $\mathcal{A} =
\Cinfty(M)$ and $\mathbb{k} = \mathbb{C}$. Next we consider two
$\mathcal{A}$-modules $\mathcal{E}, \mathcal{F}$ and set for $k<0$
\begin{equation}
    \label{eq:diffops-of-negative-order}
    \Diffop^k(\mathcal{E}; \mathcal{F}) = \{0\}
\end{equation}
and for $k \geq 0$ inductively
\begin{equation}
    \label{eq:Algebraic-Diffops}
    \Diffop^k(\mathcal{E}; \mathcal{F}) = 
    \left\{
        D \in \Hom_{\mathbb{k}}(\mathcal{E}, \mathcal{F})
        \; \Big| \;
        [D,L_a] \in
        \Diffop^{k-1}(\mathcal{E}; \mathcal{F}) 
        \; 
        \forall a \in \mathcal{A}
    \right\},
\end{equation}
where $L_a$ denotes the left multiplication of elements in the module
with $a$. As before we set
\begin{equation}
    \label{eq:all-algebraic-diffops}
    \Diffop^{\bullet}(\mathcal{E}; \mathcal{F}) 
    = 
    \bigcup_{k \in \mathbb{Z}}
    \Diffop^k(\mathcal{E}; \mathcal{F}).
\end{equation}
By general considerations it is rather easy to show that
$\Diffop^k(\mathcal{E}, \mathcal{F}) \subseteq
\Diffop^{k+1}(\mathcal{E}, \mathcal{F})$ whence
\eqref{eq:all-algebraic-diffops} is again filtered. Moreover,
$\Diffop^k(\mathcal{E}; \mathcal{F})$ is a $\mathbb{k}$-vector space
and a left $\mathcal{A}$-module via
\begin{equation}
    \label{eq:diffop-left-module}
    (a\cdot D) (e) = a\cdot D(e),
\end{equation}
where $a \in \mathcal{A}$, $D \in \Diffop^k(\mathcal{E};
\mathcal{F})$, and $e \in \mathcal{E}$.  If $\mathcal{G}$ is yet
another $\mathcal{A}$-module then the composition of differential
operators is defined and yields again differential operators. In fact,
\begin{equation}
    \label{eq:composition-of-diffops}
    \Diffop^k(\mathcal{F}; \mathcal{G})
    \circ
    \Diffop^\ell(\mathcal{E}; \mathcal{F})
    \subseteq 
    \Diffop^{k+\ell}(\mathcal{E}; \mathcal{G})
\end{equation}
holds for all $k, \ell \in \mathbb{Z}$. It follows that
\begin{equation}
    \label{eq:diffops-E-E}
    \Diffop^{\bullet}(\mathcal{E})
    =
    \Diffop^{\bullet}(\mathcal{E}; \mathcal{E})
\end{equation}
is a filtered subalgebra of all $\mathbb{k}$-linear endomorphisms
$\End_{\mathbb{k}}(\mathcal{E})$ of $\mathcal{E}$. Moreover, by
definition we have
\begin{equation}
    \label{eq:diffops-zero-order-are-Hom}
    \Diffop^0(\mathcal{E}; \mathcal{F}) =
    \Hom_{\mathcal{A}}(\mathcal{E}, \mathcal{F}).
\end{equation}
\begin{theorem}
    \label{theorem:abstract-diffops-are-usual-diffops}
    \index{Differential operator!algebraic characterization}%
    For $\mathcal{A} = \Cinfty(M)$ and $\mathcal{E} = \Secinfty(E)$,
    $\mathcal{F} = \Secinfty(F)$ the algebraic definition of
    $\Diffop^{\bullet}(\mathcal{E}; \mathcal{F})$ yields the usual
    differential operators $\Diffop^{\bullet}(E; F)$.
\end{theorem}
The proof is contained e.g. in \cite[App.~A.5]{waldmann:2007a}. We
omit it here as we shall mainly work with the local description of
differential operators.

%
%

\subsection{A Global Symbol Calculus for Differential Operators}
\label{subsec:glob-symb-calc}

The leading symbol of a differential operator is in many aspects a
much nicer object as it is a tensor field. The problem of having no
canonical definition of sub-leading symbols can be cured at the price
of a covariant derivative. We choose a torsion-free covariant
derivative $\nabla$ for the tangent bundle as well as a covariant
derivative $\nabla^E$ for $E$. Then for the operator of symmetrized
covariant differentiation $\SymD^E$ as in
Definition~\ref{definition:Symmetric-Covariant-Differentiation} we
have in any chart $(U, x)$ and with respect to any local base sections
$e_\alpha$
\begin{equation}
    \label{eq:symmcov-highest-order}
    \left(\SymD^E\right)^\ell s \At{U} = 
    \frac{\partial^\ell s^\alpha}
    {\partial x^{i_1} \ldots \partial x^{i_\ell}}
    \D x^{i_1} \vee \cdots \vee \D x^{i_\ell} \tensor e_\alpha
    + (\textrm{lower order terms}),
\end{equation}
for every section $s \in \Secinfty(E)$. This was used in the proof of
Theorem~\ref{theorem:Ck-Topologie} and is an easy consequence of the
local expression $D^E \at{U} = \D x^i \vee
\nabla_{\frac{\partial}{\partial x^i}}$ together with a simple
induction on $\ell$.

Now let $X \in \Secinfty(\Sym^k TM \tensor \Hom(E, F))$ be given. Then
locally we can write
\begin{equation}
    \label{eq:local-hom(ef)-valued-symm-section}
    X\at{U}
    =
    \frac{1}{k!} X^{i_1 \ldots i_k} {}^\beta_\alpha
    \frac{\partial}{\partial x^{i_1}}
    \vee \cdots \vee 
    \frac{\partial}{\partial x^{i_k}}
    \tensor f_\beta \tensor e^\alpha.
\end{equation}
This indicates how we can define a differential operator out of $X$
and $D^E$. We use the natural pairing of the $\Sym^k TM$-part of $X$
with the $\Sym^k T^*M$-part of $\left(\SymD^E\right)^k s$ and apply
the $\Hom(E, F)$-part of $X$ to the $E$-part of
$\left(\SymD^E\right)^k s$. This gives a well-defined section of $F$.
In the literature, different conventions concerning the pairing of
symmetric tensor fields are used. We adopt the following convention,
best expressed locally as
\begin{equation}
    \label{eq:pairing-convention}
    \SP{X, \left(\SymD^E\right)^k s}
    = k! X^{i_1 \ldots i_k} {}^\beta_\alpha
    \frac{\partial^k s^\alpha} 
    {\partial x^{i_1} \ldots \partial x^{i_k}} f_\beta
    + (\textrm{lower order terms}).
\end{equation}
With other words, this is the natural pairing of $\underbrace{V
  \tensor \cdots \tensor V}_{k-\textrm{times}}$ with $\underbrace{V^*
  \tensor \cdots \tensor V^*}_{k-\textrm{times}}$ restricted to
symmetric tensors \emph{without} additional pre-factors. Indeed, note
that the tensor indexes of $\left(\SymD^E\right)^ks$ are given by
\begin{equation}
    \label{eq:TensorIndicesOfDks}
    \left(\SymD^E\right)^k s \At{U} =     
    k!
    \frac{\partial^\ell s^\alpha}
    {\partial x^{i_1} \ldots \partial x^{i_\ell}}
    \D x^{i_1} \tensor \cdots \tensor \D x^{i_\ell} \tensor e_\alpha
    + (\textrm{lower order terms})
\end{equation}
according to our convention for the symmetrized tensor product $\vee$.
\begin{definition}[Standard ordered quantization]
    \label{definition:standard-ordered-quantization}
    \index{Standard ordered quantization}%
    Let $X \in \Secinfty(\Sym^{\bullet} TM \tensor \Hom(E, F))$ be a
    not necessarily homogeneous section and let $\hbar > 0$. Then the
    standard ordered quantization $\stdrep(X): \Secinfty(E)
    \longrightarrow \Secinfty(F)$ of $X$ is defined by
    \begin{equation}
        \label{eq:standard-ordered-quantization-formula}
        \stdrep(X) s = \sum_{r=0}^{\infty} \frac{1}{r!}
        \left(
            \frac{\hbar}{\I}
        \right)^r
        \left<
            X^{(r)}, \frac{1}{r!}\left(\SymD^E\right)^r s
        \right>,
    \end{equation}
    for $s \in \Secinfty(E)$, where $X = \sum_r X^{(r)}$ with $X^{(r)}
    \in \Secinfty(\Sym^r TM \tensor \Hom(E, F))$ are the homogeneous
    parts of $X$.
\end{definition}
Note that by definition of the direct sum there are only finitely many
$X^{(r)}$ different from zero whence the sum in
  \eqref{eq:standard-ordered-quantization-formula} is always
  \emph{finite}.
\begin{theorem}[Global symbol calculus]
    \label{theorem:Global-Symbol-Calculus}
    \index{Global symbol calculus}%
    The standard ordered quantization provides a filtration preserving
    $\Cinfty(M)$-linear isomorphism
    \begin{equation}
        \label{eq:standard-ordered-quantization}
        \stdrep:
        \bigoplus_{k=0}^{\infty} 
        \Secinfty(\Sym^k TM \tensor \Hom(E, F)) 
        \longrightarrow \Diffop^{\bullet}(E; F),
    \end{equation}
    such that for $X \in \Secinfty(\Sym^k TM \tensor \Hom(E, F))$ we
    have
    \begin{equation}
        \label{eq:leading-symbol-and-rhostd}
        \sigma_k (\stdrep(X)) = \left(\frac{\hbar}{\I}\right)^k X.
    \end{equation}
\end{theorem}
\begin{proof}
    From the local expression of $\left(\SymD^E\right)^\ell s$ as in
    the proof of Theorem~\ref{theorem:Ck-Topologie} it is clear that
    $\stdrep(X)$ is indeed a differential operator. Note that the sum
    is finite and for $X = X^{(k)} \in \Secinfty(\Sym^k TM \tensor
    \Hom(E, F))$ the differential operator $\stdrep(X)$ has order $k$.
    For $f \in \Cinfty(M)$ we clearly have $\stdrep(f X) = f
    \stdrep(X)$ since the natural pairing is $\Cinfty(M)$-bilinear.
    This shows that $\stdrep$ is a filtration preserving
    $\Cinfty(M)$-linear map. Let $X \in \Secinfty(\Sym^k TM \tensor
    \Hom(E, F))$ be homogeneous of degree $k \in \mathbb{N}_0$. Then
    locally
    \begin{align*}
        \stdrep(X) s \at{U}
        &= \frac{1}{k!} \left( \frac{\hbar}{\I} \right)^k
        \SP{
          X,
          \frac{1}{k!}\left(\SymD^E\right)^k s
        } \At{U} \\
        &= \frac{1}{k!k!} \left( \frac{\hbar}{\I} \right)^k
        X^{i_1 \ldots i_k} {}^\beta_\alpha 
        k!
        \frac{\partial^k s^\alpha}
        {\partial x^{i_1} \ldots \partial x^{i_k}} f_\beta
        + (\textrm{lower order terms}),
    \end{align*}
    hence \eqref{eq:leading-symbol-and-rhostd} is clear by the
    definition of $\sigma_k$ as in
    \eqref{eq:leading-symbol-definition}. Now let $D \in \Diffop^k(E;
    F)$ be given. Then
    \[
    \sigma_k
    \left(
        D - \left( \frac{\I}{\hbar} \right)^k \stdrep(\sigma_k(D))
    \right)
    = 0,
    \]
    hence $D - \left( \frac{\I}{\hbar} \right)^k \stdrep(\sigma_k(D))$
    is a differential operator of order at most $k-1$. By induction we
    can find $D_k = \sigma_k(D), D_{k-1}, \ldots, D_0$ with $D_\ell
    \in \Secinfty(\Sym^\ell TM \tensor \Hom(E, F))$ such that
    \begin{equation}
        \label{eq:Diffop-in-symbol-partition}
        D = \stdrep
        \left(
            \sum_{r=0}^k \left(\frac{\I}{\hbar}\right)^r D_r
        \right),
    \end{equation}
    which proves surjectivity. The injectivity is also clear, as
    $\sigma_k(D)$ is uniquely determined by $D$ and by induction the
    above $D_{k-1}, \ldots, D_0$ are unique as well.
\end{proof}
\begin{remark}[Global symbol calculus]
    \label{remark:standard-ordered-quantization}
    \index{Quantization}%
    \index{Phase space}%
    \index{Configuration space}%
    ~
    \begin{remarklist}
    \item The standard ordered quantization and its inverse map
        $\stdrep^{-1}$, i.e. the \emph{global symbol calculus}, come
        indeed from quantization theory, where $E = F = M \times
        \mathbb{C}$ is the trivial line bundle and
        $\Secinfty(\Sym^{\bullet} TM)$ is identified in the usual,
        canonical way with functions on $T^*M$ being polynomial in the
        fibers. Indeed, there is a unique algebra isomorphism
        \begin{equation}
            \label{eq:symmetric-forms-and-polynomials}
            \mathcal{J}:
            \bigoplus_{k=0}^{\infty} \Secinfty(\Sym^k TM) \ni X 
            \; \mapsto \;
            \mathcal{J}(X) \in \Pol^{\bullet}(T^*M)
        \end{equation}
        with $\mathcal{J}(f) = \pi^*f$ and $\mathcal{J}(X)(\alpha_p) =
        \alpha_p (X(p))$ for $f \in \Cinfty(M) = \Secinfty(\Sym^0 TM)$
        and $X \in \Secinfty(TM)$, where $\alpha_p \in T_p^*M$. The
        pre-factor $\frac{\hbar}{\I}$ in
        \eqref{eq:standard-ordered-quantization-formula} is due to the
        physical conventions since we can interpret functions in
        $\Pol^1 (T^*M)$ to be linear in the momenta on the phase space
        $T^*M$ corresponding to the configuration space $M$. In the
        case $M = \mathbb{R}^n$ with the flat covariant derivative
        $\nabla$, the map $\stdrep$ is indeed the standard ordered
        quantization on $T^* M = \mathbb{R}^{2n}$, i.e. first all
        ``momenta to the right''. A more detailed discussion can be
        found in \cite[Sect.~5.4]{waldmann:2007a}.
    \item For $X \tensor A \in \Secinfty(TM \tensor \Hom(E, F))$ with $X
        \in \Secinfty(TM)$ and $A \in \Secinfty(\Hom(E, F))$ we simply
        have
        \begin{equation}
            \label{eq:quantization-degree-one}
            \stdrep (X \tensor A) s = \frac{\hbar}{\I} A (\nabla^E_X s).
        \end{equation}
        In particular, the choice of $\nabla$ does not yet enter. This
        is of course no longer the case for higher symmetric degrees.
        Also
        \begin{equation}
            \label{eq:quantization-degree-zero}
            \stdrep (A) = A
        \end{equation}
        is just a $\Cinfty(M)$-linear operator, not yet
        differentiating.
    \end{remarklist}
\end{remark}

%
%

\subsection{Continuity Properties of Differential Operators}
\label{subsec:continuity-of-diffops}

From the local form of differential operators we immediately obtain
the following continuity statement:
\begin{theorem}[Continuity of differential operators]
    \label{theorem:Continuity-of-Diffops}
    \index{Differential operator!continuity}%
    \index{Continuity!differential operator}%
    Let $D \in \Diffop^k(E;F)$ be a differential operator of order
    $k$. Then for all $\ell \in \mathbb{N}_0$ the map
    \begin{equation}
        \label{eq:Diffop-on-Seck}
        D: \Sec[k+\ell](E) \longrightarrow \Sec[\ell](E)
    \end{equation}
    is well-defined and continuous with respect to the $\Fun[k+\ell]$-
    and $\Fun[\ell]$-topology.
\end{theorem}
\begin{proof}
    Clearly, if $s \in \Sec[k+\ell](E)$ then $\left(\SymD^E\right)^k s
    \in \Sec[\ell](E)$ is still $\ell$ times continuously
    differentiable. Since the natural pairing does not lower the
    degree of differentiability, we can define $\stdrep(X) s$ in the
    obvious way. Since furthermore every differential operator $D$ of
    order $k$ is of the form $\stdrep(X)$ with $X$ having at most
    tensorial degree $k$, the extension \eqref{eq:Diffop-on-Seck} is
    defined in a unique way.  If $(U, x)$ is a chart and $e_\alpha \in
    \Secinfty(E \at{U})$ and $f_\beta \in \Secinfty(F \at{U})$ are
    local base sections then
    \begin{align*}
        \seminorm[U, x, K, \ell, \{f_\beta\}] (Ds)
        & = \sup_{\substack{p \in K \\ |I| \leq \ell \\ \beta}}
        \left|
            \frac{\partial^{|I|}}{\partial x^I}
            \sum_{r=0}^{\ell} \frac{1}{r!}
            D_U^{i_1 \ldots i_r} {}_\alpha^\beta (p)
            \frac{\partial^r s^\alpha}
            {\partial x^{i_1} \cdots \partial x^{i_r}} (p)
        \right| \\
        & \leq c \sum_{\substack{i_1, \ldots, i_r \\ r}}
        \sup_{\substack{p \in K \\ |I| \leq \ell \\ \beta, \alpha}}
        \left|
            \frac{\partial^{|I|}}{\partial x^I}
            D_U^{i_1 \ldots i_r} {}_\alpha^\beta (p)
        \right|
        \sup_{\substack{p \in K \\ |J| \leq \ell \\ \alpha}}
        \left|
            \frac{\partial^{|J|}}{\partial x^J}
            \frac{\partial^r s^\alpha}
            {\partial x^{i_1} \cdots \partial x^{i_r}} (p)
        \right| \\
        & \leq c' 
        \max_{\substack{i_1, \ldots, i_r \\ \beta, \alpha}}
        \seminorm[U, x, K, \ell] (D_U^{i_1 \ldots i_r} {}_\alpha^\beta)
        \max_r \seminorm[U, x, K, \ell+r, \{e_\alpha\}] (s) \\
        & \leq c'
        \seminorm[U, x, K, \ell, \{e_\alpha\}, \{f_\beta\}] (D)
        \seminorm[U, x, K, \ell+k, \{e_\alpha\}] (s),
    \end{align*}
    where $c'$ is a combinatorial factor depending only on $\ell$ and
    $k$ and
    \[
    \seminorm[U,x,K,\ell,\{e_\alpha\},\{f_\beta\}](D)
    = \sup_{
      \substack{p \in K \\
        \alpha, \beta \\
        |I| \leq \ell \\
        i_1, \ldots, i_r}
    }
    \left|
        \frac{\partial^{|I|}  D_U^{i_1 \ldots i_r} {}_\alpha^\beta}
        {\partial x^I} (p)
    \right|.
    \]
    But this is the desired estimate to conclude the continuity with
    respect to the $\Fun[k+\ell]$- and $\Fun[\ell]$-topology.
\end{proof}
\begin{corollary}
    \label{corollary:Continuity-of-Diffops}
    A differential operator $D \in \Diffop^\bullet (E; F)$ is
    continuous with respect to the $\Cinfty$-topology.
\end{corollary}

In the proof of Theorem~\ref{theorem:Continuity-of-Diffops} we have
made use of the quantities
\begin{equation}
    \label{eq:Seminorms-on-Diffops}
    \seminorm[U, x, K, \ell, \{e_\alpha\}, \{f_\beta\}](D)
    = \sup_{\substack{p \in K \\ \alpha, \beta \\ |I| \leq \ell \\
        i_1, \ldots, i_r}}
    \left|
        \frac{\partial^{|I|}  D_U^{i_1 \ldots i_r} {}_\alpha^\beta}
        {\partial x^I} (p)
    \right|,
\end{equation}
which are easily shown to be seminorms on $\Diffop^\bullet (E;
F)$. For a fixed $k \in \mathbb{N}_0$, these make $\Diffop^k(E; F)$
again a Fr\'echet space, a simple fact which we shall not prove
here. Moreover, the standard ordered quantization is then a
\emph{continuous} isomorphism with continuous inverse
\begin{equation}
    \label{eq:standrep-is-continuous}
    \index{Differential operator!Frechet space@Fr{\'e}chet space}%
    \stdrep:
    \bigoplus_{\ell = 0}^k \Secinfty(\Sym^\ell TM \tensor \Hom(E, F))
    \longrightarrow \Diffop^k(E; F).
\end{equation}
However, \emph{all} differential operators $\Diffop^\bullet (E; F)$
will have to be equipped with an inductive limit topology similar to
the construction of the $\Cinfty_0$-topology. In any case, we shall
not need these aspects here.

Instead, we consider now the restriction of $D \in \Diffop^k(E; F)$ to
compactly supported sections $\Sec[k+\ell]_K(E)$. Since $\supp(Ds)
\subseteq \supp s$ we have
\begin{equation}
    \label{eq:Diffop-on-compact-supps}
    D: \Sec[k+\ell]_A(E) \longrightarrow \Sec[\ell]_A(F)
\end{equation}
for all closed subsets $A \subseteq M$. Since in the estimate
\begin{equation}
    \label{eq:continuity-estimate-of-diffop}
     \seminorm[U,x,K,\ell,\{f_\beta\}] (Ds)
     \leq c
     \seminorm[U,x,K,\ell,\{e_\alpha\},\{f_\beta\}] (D)
     \seminorm[U,x,K,\ell+k,\{e_\alpha\}] (s)
\end{equation}
we have the same compactum on both sides, we find that
\begin{equation}
    \label{eq:diffop-cont-in-CK-Top}
    D: \Sec[k+\ell]_K(E) \longrightarrow \Sec[\ell]_K(F)
\end{equation}
is continuous in the $\Fun[k+\ell]_K$- and
$\Fun[\ell]_K$-topology. From this we immediately obtain the following
continuity statement:
\begin{theorem}
    \label{theorem:Diffop-is-cont-on-compact-supps}
    \index{Continuity!differential operator}%
    \index{Differential operator!continuity}%
    Let $D \in \Diffop^k(E; F)$ be a differential operator of order $k
    \in \mathbb{N}_0$. Then for all $\ell \in \mathbb{N}_0$ the
    restriction
    \begin{equation}
        \label{eq:diffop-on-Ck-compacts}
        D: \Sec[k+\ell]_0(E) \longrightarrow \Sec[\ell]_0(F)
    \end{equation}
    is continuous in the $\Fun[k+\ell]_0$- and the
    $\Fun[\ell]_0$-topology. Moreover
    \begin{equation}
        \label{eq:diffop-on-cinfty-compacts}
        D: \Secinfty_0(E) \longrightarrow \Secinfty_0(F)
    \end{equation}
    is continuous in the $\Cinfty_0$-topology.
\end{theorem}
\begin{proof}
    This follow immediately from \eqref{eq:diffop-cont-in-CK-Top} and
    the characterization of continuous maps as in
    Theorem~\ref{theorem:inductive-limit-topology},
    \refitem{item:LFcontinuity}.
\end{proof}

%
%

\subsection{Adjoints of Differential Operators}
\label{subsec:adjoints-of-diffops}

For a section $s \in \Secinfty(E)$ and $\mu \in \Secinfty(E^* \tensor
\Dichten T^*M)$ the natural pairing of $E$ and $E^*$ gives a
density\index{Density} $\mu(s) \in \Secinfty(\Dichten T^*M)$ which we
can integrate, provided the support is compact. Therefore we define
\begin{equation}
    \label{eq:pairing-on-sections-and-densities}
    \SP{s, \mu} = \int_M \mu(s) = \int_M s \cdot \mu,
\end{equation}
whenever the support of at least one of $s$ or $\mu$ is compact.
\begin{lemma}
    \label{lemma:pairing-on-sections-and-densities}
    The pairing \eqref{eq:pairing-on-sections-and-densities} is
    bilinear and non-degenerate. Moreover $\SP{s, f \mu} = \SP{f s,
      \mu}$ for $f \in \Cinfty(M)$.
\end{lemma}
\begin{proof}
    Let $s \in \Secinfty(E)$ be not the zero section and let $p \in M$
    be such that $s(p)\neq 0$. Then we find an open neighborhood $U$
    of $p$ and a section $\varphi \in \Secinfty_0(E^*)$ with compact
    support $\supp \varphi \subseteq U$ such that
    \[
    \varphi (s) \geq 0 
    \quad
    \textrm{and}
    \quad \varphi(s) \at{p} > 0.
    \]
    Using local base sections this is obvious. Now choose a positive
    density\index{Density!positive} $\nu \in \Secinfty(\Dichten
    T^*M)$, then $\varphi \tensor \nu \in \Secinfty_0(E^* \tensor
    \Dichten T^*M)$ will satisfy $\SP{s, \varphi \tensor \nu} \neq
    0$. This shows that \eqref{eq:pairing-on-sections-and-densities}
    is non-degenerate in the first argument. The other non-degeneracy
    is shown analogously.  The second statement is clear.
\end{proof}

In particular, $\SP{\argument, \argument}$ restricts to a non-degenerate
pairing
\begin{equation}
    \label{eq:restriction-of-pairing-to-both-compact}
    \SP{\argument, \argument}:
    \Secinfty_0(E) \times \Secinfty_0(E^* \tensor \Dichten T^*M)
    \longrightarrow \mathbb{C}.
\end{equation}
As an immediate consequence we obtain the following statement. First
recall that an operator
\begin{equation}
    \label{eq:adjoinabel-operator}
    D: V \longrightarrow W
\end{equation}
is \emph{adjointable}\index{Adjointable operator} with respect to
bilinear pairings
\begin{equation}
    \label{eq:3}
    \SP{\argument, \argument}_{V,\widetilde{V}}:
    V \times \widetilde{V} \longrightarrow \mathbb{C}
    \quad \textrm{and} \quad
    \SP{\argument, \argument}_{W,\widetilde{W}}:
    W \times \widetilde{W} \longrightarrow \mathbb{C},
\end{equation}
if there is a map $D^\Trans: \widetilde{W} \longrightarrow \widetilde{V}$
such that
\begin{equation}
    \label{eq:adjoint-operator-definition}
    \SP{Dv, \widetilde{w}}_{W, \widetilde{W}}
    = \SP{v, D^\Trans \widetilde{w} }_{V, \widetilde{V}}.
\end{equation}
If the pairings are non-degenerate then an adjoint $D^\Trans$ is
necessarily unique (if it exists at all) and both maps $D$, $D^\Trans$
are linear maps. Clearly, $D^\Trans$ is adjointable, too, with
$(D^\Trans)^\Trans = D$. Thus in our situation, adjointable maps with
respect to the pairing \eqref{eq:pairing-on-sections-and-densities} or
\eqref{eq:restriction-of-pairing-to-both-compact} have unique adjoints
and are necessarily linear.
\begin{proposition}
    \label{proposition:existence-of-adjoints}
    \index{Differential operator!adjoint}%
    Let $D \in \Diffop^k(E; F)$ be a differential operator of order
    $k$. Then $D: \Secinfty_0(E) \longrightarrow \Secinfty_0(F)$ is
    adjointable with respect to
    \eqref{eq:pairing-on-sections-and-densities} and the (unique)
    adjoint
    \begin{equation}
        \label{eq:adjoint-diffop}
        D^\Trans: \Secinfty(F^* \tensor \Dichten T^*M) \longrightarrow
        \Secinfty(E^* \tensor \Dichten T^*M)
    \end{equation}
    is again a differential operator of order $k$.
\end{proposition}
\begin{proof}
    Let $\{(U_i, x_i)\}_{i \in I}$ be a locally finite atlas and let
    $e_{i \alpha} \in \Secinfty(E \at{U_i})$ and $f_{i \beta} \in
    \Secinfty(F \at{U_i})$ be local base sections. Moreover let
    $\{\chi_i\}_{i \in I}$ be a locally finite partition of unity
    subordinate to the atlas with $\supp \chi_i$ being compact. As
    usual, we write
    \[
    Ds \at{U_i} = \sum_{r=0}^k \frac{1}{r!}
    D_{U_i}^{i_1 \ldots i_r} {}_\alpha^\beta f_\beta 
    \frac{\partial^r s_i^\alpha}
    {\partial x^{i_1}_i \cdots \partial x^{i_r}_i},
    \]
    where $s\at{U_i} = s_i^\alpha e_{i \alpha}$ with $s_i^\alpha =
    e^\alpha_i(s) \in \Cinfty(U_i)$. For $\mu \in \Secinfty(F^*
    \tensor \Dichten T^*M)$ we write
    \[
    \mu \at{U_i}
    = \mu_{i \beta} f^\beta |\D x^1_i \wedge \cdots \wedge \D x^n_i|,
    \]
    with $\mu_{i \beta} \in \Cinfty(U_i)$. Here $|\D x^1 \wedge \cdots
    \wedge \D x^n|$ denotes the unique locally defined density with
    value $1$ when evaluated on the coordinate base vector fields
    $\frac{\partial}{\partial x^1}, \ldots, \frac{\partial}{\partial
      x^n}$. Then we compute
    \begin{align*}
        \SP{Ds, \mu} & = \int_M \mu(Ds)
        = \sum_i \int_{x_i(U_i)} \left( \chi_i \mu(Ds) \right)
        \circ x_i^{-1} \D^n x_i \\
       & = \sum_i \int_{x_i(U_i)}
        \left(
            \chi_i \mu_{i\beta} \sum_{r=0}^{k} \frac{1}{r!}
            D_{U_i}^{i_1 \ldots i_r} {}_\alpha^\beta 
            \frac{\partial^r s_i^\alpha}
            {\partial x^{i_1}_i \cdots \partial x^{i_r}_i}
        \right)
        \circ x_i^{-1} \D^n x_i.
    \end{align*}
    Note that the integrand consists of compactly supported functions
    only. Thus we can integrate by parts and obtain
    \[
    \SP{Ds, \mu} = \sum_i \int_{x_i(U_i)}
    \left(
        \sum_{r=0}^{k} \frac{(-1)^r}{r!}
        \frac{\partial^r}
        {\partial x^{i_1} \cdots \partial x^{i_r}}
        \left(
            \chi_i \mu_{i\beta}
            D_{U_i}^{i_1 \ldots i_r} {}_\alpha^\beta
        \right)
        s_i^\alpha
    \right)
    \circ x_i^{-1} \D^n x_i.
    \]
    Now the function $\chi_i \mu_{i\beta} D_{U_i}^{i_1 \ldots i_r}
    {}_\alpha^\beta$ has compact support in $U_i$ thanks to the choice
    of the $\chi_i$. Thus it defines a global function in
    $\Cinfty_0(M)$. It follows that
    \[
    \mu_i =  \sum_{r=0}^{k} \frac{(-1)^r}{r!}
    \frac{\partial^r}
    {\partial x^{i_1} \cdots \partial x^{i_r}}
    \left(
        \chi_i \mu_{i\beta}
        D_{U_i}^{i_1 \ldots i_r} {}_\alpha^\beta
    \right)
    e_i^\alpha \tensor |\D x_i^1 \wedge \cdots \wedge \D x_i^n|
    \]
    is a global section in $\Secinfty_0(E^* \tensor \Dichten T^*M)$
    with compact support in $U_i$. Since the $\chi_i$ are locally
    finite, the sum
    \[
    D^\Trans \mu = \sum_i \mu_i
    \]
    is well-defined and yields a global section $D^\Trans \mu \in
    \Secinfty(E^* \tensor \Dichten T^*M)$ such that
    \[
    \SP{Ds, \mu} = \SP{s, D^\Trans \mu}.
    \]
    This shows that $D$ is adjointable. From the actual computation
    above it is clear that $D^\Trans$ differentiates again $k$ times.
    Thus $D^\Trans \in \Diffop^k(F^* \tensor \Dichten T^*M, E^*
    \tensor \Dichten T^*M)$ follows. However, there is also another
    nice argument based on the algebraic definition of differential
    operators: Let $D: \Secinfty(E) \longrightarrow \Secinfty(F)$ be a
    differential operator of order zero. Thus $D$ can be viewed as a
    section of $\Hom(E, F)$, i.e. $D \in \Secinfty(\Hom(E, F))$. Then
    in $\mu (Ds)$ we can simply apply the pointwise transpose of $D$
    to the $F^*$-part of $\mu$. This defines $D^\Trans \mu$ pointwise
    in such a way that $(D^\Trans \mu)(s) = \mu (Ds)$. Clearly
    $\SP{Ds, \mu} = \SP{s, D^\Trans \mu}$
    follows. Now we proceed by induction. We assume that the adjoint
    always exists (what we have shown already) and for differential
    operators of order $\ell \leq k-1$ the adjoint has order $\ell$,
    too. Thus let $D \in \Diffop^k(E; F)$ and $f \in \Cinfty(M)$. Then
    we have
    \[
    \SP{fDs, \mu} = \SP{Ds, f\mu}
    = \SP{s, D^\Trans f\mu},
    \]
    and on the other hand
    \begin{align*}
        \SP{fDs, \mu}
        & = \SP{[f, D] s, \mu}
        + \SP{D(fs), \mu} \\
        & = \SP{s, [f, D]^\Trans \mu}
        + \SP{fs, D^\Trans \mu} \\
        & = \SP{s, [f, D]^\Trans \mu}
        +\SP{s, f D^\Trans \mu}.
    \end{align*}
    Hence by the non-degeneracy of $\SP{\argument, \argument}$ we
    conclude that
    \[
    [f, D^\Trans] = [f, D]^\Trans \in \Diffop^{k-1}
    (F^* \tensor \Dichten T^*M, E^* \tensor \Dichten T^*M)
    \]
    by induction. But this shows $D^\Trans \in \Diffop^k (F^* \tensor
    \Dichten T^*M, E^* \tensor \Dichten T^*M) $ as wanted.
\end{proof}
\begin{corollary}
    \label{corollary:symbol-of-adjoint}
    \index{Differential operator!leading symbol!adjoint}%
    Let $D \in \Diffop^k(E; F)$. Then for the leading symbol $\sigma_k
    (D^\Trans) \in \ \Secinfty(\Sym^k TM \tensor \Hom(F^* \tensor
    \Dichten T^*M, E^* \tensor \Dichten T^*M))$ we have
    \begin{equation}
        \label{eq:symbol-of-adjoint}
        \sigma_k (D^\Trans)
        = (-1)^k \sigma_k (D)^\Trans \tensor \id_{\Dichten T^*M},
    \end{equation}
    where $\sigma_k(D)^\Trans$ denotes the pointwise transpose from
    $\Hom(E, F)$ to $\Hom(F^*, E^*)$.
\end{corollary}
\begin{proof}
    From the local computations in the proof of
    Proposition~\ref{proposition:existence-of-adjoints} we obtained
    \begin{align*}
        \mu_i
        &= \sum_{r=0}^{k} \frac{(-1)^r}{r!}
        \frac{\partial^r}
        {\partial x^{i_1} \cdots \partial x^{i_r}}
        \left(
            \chi_i \mu_{i\beta}
            D_{U_i}^{i_1 \ldots i_r} {}_\alpha^\beta
        \right)
        e_i^\alpha \tensor |\D x_i^1 \wedge \cdots \wedge \D x_i^n| \\
        &=
        \frac{(-1)^k}{k!} \chi_i
        D_{U_i}^{i_1 \ldots i_k} {}_\alpha^\beta
        \frac{\partial^k \mu_{i\beta}}
        {\partial x^{i_1} \cdots \partial x^{i_k}}
        e_i^\alpha \tensor |\D x_i^1 \wedge \cdots \wedge \D x_i^n|
        + (\textrm{lower order terms}).
    \end{align*}
    Since $D^\Trans \mu = \sum_i \mu_i$ and $\sum_i \chi_i = 1$, we
    conclude that
    \begin{align*}
        D^\Trans \mu \at{U_i}
        &=
        \frac{(-1)^k}{k!} D_{U_i}^{i_1 \ldots i_k} {}_\alpha^\beta
        \frac{\partial^k \mu_{i\beta}}
        {\partial x^{i_1}_i \cdots \partial x^{i_k}_i}
        e_i^\alpha \tensor |\D x_i^1 \wedge \cdots \wedge \D x_i^n|
        + (\textrm{lower order terms}) \\
        & = (-1)^k \sigma_k(D)^\Trans \tensor \id_{\Dichten T^*M} (\mu)
        + (\textrm{lower order terms}).
    \end{align*}
\end{proof}
\begin{remark}[Other pairings]
    \label{remark:other-pairings}
    \index{Density!pairings}%
    ~
    \begin{remarklist}
    \item There are several variations of the above proposition. On
        one hand one can consider the natural pairing of $\alpha$- and
        $(1-\alpha)$-densities for any $\alpha \in \mathbb{C}$ to
        obtain
        \begin{equation}
            \label{eq:pairing-of-alpha-densities}
            \SP{\argument, \argument} :
            \Secinfty_0(E \tensor \Dichten^\alpha T^*M)
            \times
            \Secinfty_0(E^* \tensor \Dichten^{1-\alpha} T^*M)
            \longrightarrow \mathbb{C}
        \end{equation}
        via pointwise natural pairing and integration of the remaining
        $1$-density. This is again non-degenerate. Thus we can also
        compute the adjoints of differential operators
        \begin{equation}
            \label{eq:diffops-for-alpha-sections}
            D: \Secinfty_0(E \tensor \Dichten^{\alpha} T^*M)
            \longrightarrow
            \Secinfty_0(F \tensor \Dichten^{\beta} T^*M)
        \end{equation}
        and obtain differential operators
        \begin{equation}
            \label{eq:4}
            D^\Trans:
            \Secinfty(F^* \tensor \Dichten^{1-\beta} T^*M)
            \longrightarrow
            \Secinfty(E^* \tensor \Dichten^{1-\alpha} T^*M)
        \end{equation}
        by the same kind of computation as in
        Proposition~\ref{proposition:existence-of-adjoints}. There, we
        considered the case $\alpha = 0 = \beta$.
    \item Another important case is for complex bundles $E$ with a
        (pseudo-) Hermitian fiber metric $h_E$. Then we can use the
        $\mathbb{C}$-sesquilinear pairings
        \begin{equation}
            \label{eq:pairing-with-metric-on-factorizing-sections}
            \SP{s, t \tensor \mu} = \int_M h(s,t) \: \mu,
        \end{equation}
        where $s, t \in \Secinfty(E)$ and $\mu \in \Secinfty(\Dichten
        T^*M)$ and at least one has compact support. Clearly, this
        extends to
        \begin{equation}
            \label{eq:pairing-with-metric}
            \SP{\argument, \argument}: 
            \Secinfty(E) \times
            \Secinfty_0(E \tensor \Dichten T^*M)
            \longrightarrow \mathbb{C}
        \end{equation}
        in a $\mathbb{C}$-sesquilinear way. While $D \mapsto D^\Trans$
        is $\mathbb{C}$-linear, now the adjoint $D^*$ depends on $D$
        in an \emph{antilinear} way.
    \item A very important situation is obtained by merging the above
        possibilities. For a Hermitian vector bundle $E
        \longrightarrow M$ with Hermitian fiber metric $h$ we consider
        the sections $\Secinfty_0(E \tensor \Dichten^{\frac{1}{2}}
        T^*M)$. On factorizing sections we can define
        \begin{equation}
            \label{eq:pairing-of-half-densities-on-factorizing-sections}
            \SP{s \tensor \mu, t \tensor \nu}
            = \int_M h(s,t) \: \cc{\mu} \nu,
        \end{equation}
        since $\cc{\mu} \nu$ is a $1$-density. Then the pairing
        extends to a $\mathbb{C}$-sesquilinear pairing
        \begin{equation}
            \label{eq:pairing-of-half-densities}
            \SP{\argument, \argument}:
            \Secinfty_0(E \tensor \Dichten^{\frac{1}{2}} T^*M) \times
            \Secinfty_0(E \tensor \Dichten^{\frac{1}{2}} T^*M)
            \longrightarrow \mathbb{C},
        \end{equation}
        which is not only non-degenerate but \emph{positive definite}.
        Thus $\Secinfty_0(E \tensor \Dichten^{\frac{1}{2}} T^*M)$
        becomes a \emph{pre-Hilbert space}. Moreover, taking $E$ to be
        the trivial line bundle with the canonical fiber metric gives
        a pre-Hilbert space $\Secinfty_0(\Dichten^{\frac{1}{2}}T^*M)$
        of half densities. Its completion to a Hilbert space is the
        so-called \emph{intrinsic Hilbert space} on $M$.
        \index{Intrinsic Hilbert space}%
        \index{Pairing!positive definite}%
        \index{Pairing!half-densities}%
    \end{remarklist}
\end{remark}

While the above constructions are always slightly asymmetric unless we
take half-densities, we obtain a more symmetric situation if we
integrate with respect to a given positive density. Thus we choose
once and for all a positive density $\mu > 0$ on $M$. Later on, this
will be the (pseudo-) Riemannian volume density, but for now we do not
need this additional property. For a vector bundle $E \longrightarrow
M$ we then have the pairing
\begin{equation}
    \label{eq:pairing-with-fixed-densitiy}
    \SP{s, \varphi}_\mu = \int_M \varphi(s) \: \mu,
\end{equation}
for $s \in \Secinfty(E)$ and $\varphi \in \Secinfty(E^*)$, at least
one having compact support. Clearly,
\begin{equation}
    \label{eq:relation-pairing-with-fixed-densitiy-pairing-without-fixed-density}
    \SP{s, \varphi}_\mu
    = \SP{s, \varphi \tensor \mu} 
\end{equation}
with the original version \eqref{eq:pairing-on-sections-and-densities}
of the pairing $\SP{\argument, \argument}$. Since $\mu > 0$ it
easily follows that
\eqref{eq:relation-pairing-with-fixed-densitiy-pairing-without-fixed-density}
is non-degenerate and satisfies
\begin{equation}
    \label{eq:fixpairing-functions-left-right}
    \SP{fs, \varphi}_\mu
    = \SP{s, f \varphi}_\mu
\end{equation}
for all $f \in \Cinfty(M)$. For the action of differential operators
we again have adjoints:
\begin{theorem}
    \label{theorem:agjoint-for-fixed-density}
    \index{Differential operator!adjoint}%
    Let $D \in \Diffop^k(E; F)$ be a differential operator of order $k
    \in \mathbb{N}_0$. Then there exists a differential operator
    $D^\Trans \in \Diffop^k(F^*; E^*)$ such that
    \begin{equation}
        \label{eq:adjoint-for-fixed-density}
        \SP{Ds, \varphi}_\mu
        =
        \SP{s, D^\Trans \varphi}_\mu
    \end{equation}
    for all $s \in \Secinfty(E)$ and $\varphi \in \Secinfty(F^*)$, at
    least one having compact support.
\end{theorem}
\begin{proof}
    The proof is now fairly simple. Since $D$ has an adjoint, denoted
    by $\widetilde{D}$ for a moment, with respect to
    \eqref{eq:pairing-on-sections-and-densities} we have
    \[
    \SP{Ds, \varphi}_\mu
    = \SP{Ds, \varphi \tensor \mu}
    = \SP{s, \widetilde{D} (\varphi \tensor \mu)},
    \]
    and locally
    \begin{align*}
        \widetilde{D} (\varphi \tensor \mu) \At{U}
        &=
        \sum_{r=0}^{k} \frac{1}{r!}
        \widetilde{D}^{i_1 \ldots i_r}_U {}^\beta_\alpha
        \frac{\partial^r}{\partial x^{i_1} \cdots \partial x^{i_r}}
        \left(
            \varphi_\beta \mu_U
        \right)
        e^\alpha \tensor |\D x^1 \wedge \cdots \wedge \D x^n| \\
        &=
        \sum_{r=0}^k \frac{1}{r!}
        \widetilde{D}^I_U {}^\beta_\alpha
        \sum_{J \leq I} \binom{I}{J}
        \frac{\partial^{|J|} \varphi_\beta}{\partial x^J}
        \frac{\partial^{|I-J|} \mu_U}{\partial x^{I-J}}
        e^\alpha \tensor |\D x^1 \wedge \cdots \wedge \D x^n| \\
        &=
        \sum_{\substack{r=0 \\ |I| = r \\ J \leq I}}^k
        \frac{1}{r!} \binom{I}{J}
        \widetilde{D}^I_U {}^\beta_\alpha
        \frac{\partial^{|J|} \varphi_\beta}{\partial x^J}
        \frac{1}{\mu_U}
        \frac{\partial^{|I-J|} \mu_U}{\partial x^{I-J}}
        e^\alpha \tensor \mu_U |\D x^1 \wedge \cdots \wedge \D x^n| \\
        &=
        \Big(
        \sum_{\substack{r=0 \\ |I| = r \\ J \leq I}}^k
        \frac{1}{r!} \binom{I}{J}
        \widetilde{D}^I_U {}^\beta_\alpha
        \frac{\partial^{|J|} \varphi_\beta}{\partial x^J}
        \frac{1}{\mu_U}
        \frac{\partial^{|I-J|} \mu_U}{\partial x^{I-J}}
        e^\alpha
        \Big)
        \tensor \mu \at{U},
    \end{align*}
    since $\mu_U > 0$ thanks to $\mu > 0$. This shows that with
    \[
    D^\Trans \varphi \At{U}
    = \sum_{\substack{r=0 \\ |I| = r \\ J \leq I}}^k
    \frac{1}{r!} \binom{I}{J}
    \widetilde{D}^I_U {}^\beta_\alpha
    \frac{\partial^{|J|} \varphi_\beta}{\partial x^J}
    \frac{1}{\mu_U}
    \frac{\partial^{|I-J|} \mu_U}{\partial x^{I-J}}
    e^\alpha
    \tag{$*$}
    \]
    we obtain a locally defined differential operator $D^\Trans$ such
    that
    \[
    \widetilde{D} (\varphi \tensor \mu) \At{U}
    = (D^\Trans \varphi) \tensor \mu \At{U}.
    \]
    Now the left hand side is globally well-defined and hence the
    right hand side is chart independent as well. This shows that
    $D^\Trans$ is indeed a global object, locally given by ($*$).
    Obviously, it is a differential operator of order $k$.
\end{proof}

\begin{remark}
    \label{remark:adjoint-for-fixed-density}
    ~
    \begin{remarklist}
    \item Note that $D^\Trans$ as in
        Theorem~\ref{theorem:agjoint-for-fixed-density} depends on the
        choice of $\mu > 0$ while the adjoint as in
        Proposition~\ref{proposition:existence-of-adjoints} is
        intrinsically defined, though of course between different
        vector bundles.  However, we shall not emphasize the
        dependence of $D^\Trans$ on $\mu$ in our notation. It should
        become clear from the context which version of adjoint we use.
    \item Analogously to Corollary~\ref{corollary:symbol-of-adjoint}
        we see that the leading symbol of $D^\Trans$ is given by
        \begin{equation}
            \label{eq:symbol-of-adjoint-for-fixed-density}
            \index{Differential operator!leading symbol!adjoint}%
            \sigma_k(D^\Trans) = (-1)^k \sigma_k(D)^\Trans,
        \end{equation}
        where again $\sigma_k(D)^\Trans \in \Secinfty(\Hom(F^*, E^*))$
        is the pointwise adjoint of $\sigma_k(D) \in \Secinfty(\Hom(E,
        F))$. This is obvious from the local computations in the proof
        as we have to collect those terms with all $k$ derivatives
        hitting the $\varphi_\beta$ instead of the $\mu_U$.
    \end{remarklist}
\end{remark}

Sometimes it will be important to compute the adjoint of $D^\Trans$
more explicitly. Here we can use our global symbol calculus developed
in Section~\ref{subsec:glob-symb-calc}. To this end, we introduce the
following divergence operators. If $X \in \Secinfty(TM)$ is a vector
field then its \emph{covariant divergence} is defined by
\begin{equation}
    \label{eq:covariant-divergence}
    \index{Covariant divergence}%
    \divergenz_\nabla (X) = \tr (Y \; \mapsto \; \nabla_Y X),
\end{equation}
where the trace is understood to be the pointwise trace: indeed $Y
\mapsto \nabla_Y X$ is a $\Cinfty(M)$-linear map $\Secinfty(TM)
\longrightarrow \Secinfty(TM)$ which therefor can be identified with a
section in $\Secinfty(\End(TM))$. Thus the trace is well-defined. More
explicitly, in local coordinates $(U, x)$ we have
\begin{equation}
    \label{eq:divergence-local-expression}
    \divergenz_\nabla (X) \at{U}
    = 
    \D x^i \left(\nabla_{\frac{\partial}{\partial x^i}} X\right).
\end{equation}
Clearly, we have for $f \in \Cinfty(M)$ and $X \in \Secinfty(TM)$ the
relation
\begin{equation}
    \label{eq:divergence-product rule}
    \divergenz_\nabla (fX) = X(f) + f \divergenz_\nabla (X)
\end{equation}
This \Index{Leibniz rule} suggests to extend the covariant divergence
to higher symmetric multivector fields as follows.
\begin{definition}[Covariant divergence]
    \label{definition:CovariantDivergence}
    \index{Covariant divergence}%
    Let $\nabla$ be a torsion-free covariant derivative for $M$ and
    let $\nabla^E$ be a covariant derivative for $E$. For $X \in
    \Secinfty(\Sym^\bullet TM \tensor E)$ we define
    \begin{equation}
        \label{eq:covariant-divergence-of-multivector-fields}
        \divergenz^E_\nabla (X) = \inss (\D x^i)
        \nabla_{\frac{\partial}{\partial x^i}} X.
    \end{equation}
\end{definition}
\begin{lemma}
    \label{lemma:covariant-divergence-of-multivector-fields}
    By \eqref{eq:covariant-divergence-of-multivector-fields} we obtain
    a globally well-defined operator
    \begin{equation}
        \label{eq:covariant-divergence-of-multivectors-as-a-map}
        \divergenz_\nabla^E: \Secinfty(\Sym^\bullet TM \tensor E)
        \longrightarrow \Secinfty(\Sym^{\bullet -1} TM \tensor E),
    \end{equation}
    which is given on factorizing sections by
    \begin{align}
        \label{eq:covariant-divergence-multivectors-product-rule}
        \divergenz^E_\nabla (X_1 \vee \cdots \vee X_k \tensor s)
        = & \sum_{\ell=1}^{k}
        X_1 \vee \cdots \stackrel{\ell}{\wedge} \cdots \vee X_k
        \tensor
        \left(
            \divergenz_\nabla (X_\ell) s + \nabla^E_{X_\ell} s
        \right) \\
        & + \sum_{\substack{\ell,m=1 \\ \ell \neq m}}^k
        (\nabla_{X_\ell} X_m) \vee X_1 \vee \cdots \stackrel{\ell}{\wedge}
        \cdots \vee X_k \tensor s,
    \end{align}
    where $X_1, \ldots, X_k \in \Secinfty(TM)$ and $s \in
    \Secinfty(E)$.
\end{lemma}
\begin{proof}
    First it is clear that the transformation properties of
    $\frac{\partial}{\partial x^i}$ and $\D x^i$ under a change of
    local coordinates guarantee that $\divergenz^E_\nabla$ is indeed
    well-defined and independent of the chart. Thus
    $\divergenz^E_\nabla$ is a globally defined operator lowering the
    symmetric degree by one. Now let $X_1, \ldots, X_k \in
    \Secinfty(TM)$ and $s \in \Secinfty(E)$ be given. Then we compute
    \begin{align*}
        & \divergenz^E_\nabla (X_1 \vee \cdots \vee X_k \tensor s) \\
        &\quad= \inss (\D x^i) \nabla_{\frac{\partial}{\partial x^i}}
        (X_1 \vee \cdots \vee X_k \tensor s) \\
        &\quad= \inss (\D x^i)
        \left(
            \sum_{\ell=1}^k X_1 \vee \cdots \vee
            \nabla_{\frac{\partial}{\partial x^i}} X_\ell \vee \cdots X_k
            \tensor s
            + X_1 \vee \cdots \vee X_k \tensor
            \nabla^E_{\frac{\partial}{\partial x^i}} s
        \right) \\
        &\quad=
        \sum_{\substack{\ell,m=1 \\ \ell \neq m}}
        X_1 \vee \cdots \vee 
        \nabla_{\frac{\partial}{\partial x^i}} X_\ell
        \vee \cdots \vee \D x^i (X_m) \vee \cdots \vee X_k \tensor s \\
        &\quad\quad+
        \sum_{\ell=1}^{k} X_1 \vee \cdots \vee
        \D x^i \left(
            \nabla_{\frac{\partial}{\partial x^i}} X_\ell
        \right) 
        \vee \cdots \vee X_k \tensor s
        + \sum_{\ell=1}^k 
        X_1 \vee \cdots \vee \D x^i (X_\ell) \vee \cdots
        \vee X_k \tensor \nabla_{\frac{\partial}{\partial x^i}}^E s \\
        &\quad=
        \sum_{\substack{\ell,m=1 \\ \ell \neq m}}
        X_1 \vee \cdots \vee \nabla_{X_m} X_\ell
        \vee \cdots \stackrel{m}{\wedge} \cdots \vee X_k \tensor s \\
        &\quad\quad+
        \sum_{\ell=1}^{k} X_1 \vee \cdots \vee
        \divergenz_\nabla (X_\ell) \vee \cdots \vee X_k \tensor s
        + \sum_{\ell=1}^k 
        X_1 \vee \cdots \stackrel{\ell}{\wedge}  \cdots
        \vee X_k \tensor \nabla_{X_\ell}^E s.
  \end{align*}
\end{proof}

The covariant derivative $\nabla$ also acts on densities hence we can
compute the derivative $\nabla_X \mu$ of the positive density
$\mu$. This defines a function
\begin{equation}
    \label{eq:alpha-of-vector-field}
    \alpha (X) = \frac{\nabla_X \mu}{\mu},
\end{equation}
depending $\Cinfty(M)$-linearly on $X$. Thus we obtain a one-form
$\alpha \in \Secinfty(T^*M)$ which measures how much $\mu$ is
\emph{not} covariantly constant. Similarly, we can define the
$\mu$-divergence of a vector field by
\begin{equation}
    \label{eq:mu-divergence}
    \index{Density!divergence}%
    \index{Divergence}%
    \divergenz_\mu (X) = \frac{\Lie_X \mu}{\mu}.
\end{equation}
\begin{lemma}
    \label{lemma:mu-divergence-and-covariant-divergence}
    For $X \in \Secinfty(TM)$ we have
    \begin{equation}
        \label{eq:mu-divergence-covariant-divergence}
        \divergenz_\mu (X) = \divergenz_\nabla (X) +\alpha(X).
    \end{equation}
\end{lemma}
\begin{proof}
    This can be obtained from a simple computation in local
    coordinates which we omit here, see e.g.
    \cite[Sect.~2.3.4]{waldmann:2007a}.
\end{proof}

Writing this as
\begin{equation}
    \label{eq:mu-divergence-extension-idea}
    \divergenz_\mu(X) = \divergenz_\nabla(X) + \inss(\alpha) X,
\end{equation}
we can motivate the following definition. For $X \in
\Secinfty(\Sym^\bullet TM \tensor E)$ we set
\begin{equation}
    \label{eq:mu-divergence-of-multivector-fields}
    \divergenz^E_\mu (X) = \divergenz^E_\nabla (X) + \inss(\alpha) X,
\end{equation}
where $\inss(\alpha)$ acts on the $\Sym^\bullet TM$-part as usual.
\begin{lemma}
    \label{lemma:mu-divergence-of-multis}
    On factorizing section we have
    \begin{align}
        \label{eq:mu-divergence-on-factorizing-multis}
        \divergenz^E_\mu (X_1 \vee \cdots \vee X_k \tensor s)
        & = \sum_{\ell=1}^k
        X_1 \vee \cdots \stackrel{\ell}{\wedge} \cdots \vee X_k
        \tensor
        \left(
            \divergenz_\mu (X_\ell) s + \nabla^E_{X_\ell} s
        \right) \\
        & +\sum_{\substack{\ell, m=1 \\ \ell \neq m}}^k
        \nabla_{X_\ell} X_m \vee X_1 \vee \cdots
        \stackrel{\ell}{\wedge} \cdots \stackrel{m}{\wedge}
        \cdots \vee X_k \tensor s.
    \end{align}
\end{lemma}
\begin{proof}
    The proof of \eqref{eq:mu-divergence-on-factorizing-multis} is
    completely analogous to the proof of
    Lemma~\ref{lemma:covariant-divergence-of-multivector-fields}.
\end{proof}

We can now use the divergence operator to compute the adjoint of a
differential operator in a symbol calculus explicitly:
\begin{theorem}[Neumaier]
    \label{theorem:Neumaier-Theorem}
    Let $X \in \Secinfty(\Sym^k TM \tensor \Hom(E, F))$ and let
    $\nabla$ and $\nabla^E$, $\nabla^F$ be given. Then the adjoint
    operator to $\stdrep (X)$ with respect to $\SP{\argument,
      \argument}_\mu$ is explicitly given by
    \begin{equation}
        \label{eq:adjoint-of-stdrep(X)}
        \stdrep (X)^\Trans = (-1)^k \stdrep(N^2 X^\Trans),
    \end{equation}
    where
    \begin{equation}
        \label{eq:Neumaier-Operator}
        N = \exp
        \left(
            \frac{\hbar}{2\I} \divergenz_\mu^{\Hom(E, F)}
        \right)
    \end{equation}
    and where we use the induced covariant derivative on $\Hom(E, F)$
    and $\Hom(F^*, E^*)$.
\end{theorem}
\begin{proof}
    By a partition of unity argument we can reduce the problem to the
    case where the involved tensor fields have compact support in a
    chart $(U, x)$. In this chart we first note that from the
    definition of the covariant derivative of a density we obtain the
    local expression
    \[
    \alpha
    = \left(
        \frac{\Lie_{\frac{\partial}{\partial x^i}} \mu}{\mu}
        - \Gamma_{i \ell}^\ell
    \right)
    \D x^i
    \]
    for the one-form $\alpha$. Now let $\omega \in \Secinfty(\Sym^\ell
    T^*M \tensor F^*)$ and $s \in \Secinfty(E)$. In the following of
    this proof, we shall simply write $\divergenz_\mu$ for all
    divergences instead of specifying the vector bundle explicitly,
    just to simplify our notation.  For $X \in \Secinfty(\Sym^{k+\ell}
    TM \tensor \Hom(E, F))$ we compute
    \begin{align*}
        &\Lie_{\frac{\partial}{\partial x^i}}
        \left(
            \SP{
              \inss (\D x^i) X,
              \omega \tensor \left(\SymD^E\right)^{k-1} s
            } \mu
        \right) \\
        &=
        \frac{\partial}{\partial x^i}
        \left(
            \SP{
              \inss (\D x^i) X,
              \omega \tensor \left(\SymD^E\right)^{k-1} s
            }
        \right) \mu
        + \SP{
          \inss (\D x^i) X,
          \omega \tensor \left(\SymD^E\right)^{k-1} s
        }
        \Lie_{\frac{\partial}{\partial x^i}} \mu \\
        &=
        \SP{
          \nabla_{\frac{\partial}{\partial x^i}}
          \left(\inss (\D x^i) X\right),
          \omega \tensor \left(\SymD^E\right)^{k-1} s
        } \mu
        +
        \SP{
          \inss(\D x^i) X,
          \nabla_{\frac{\partial}{\partial x^i}}           
          \omega \tensor \left(\SymD^E\right)^{k-1} s
        } \mu \\
        &\quad+
        \SP{
          \inss\left(
              \frac{\Lie_{\frac{\partial}{\partial x^i}} \mu}{\mu}
              \D x^i
          \right) X,
          \omega \tensor \left(\SymD^E\right)^{k-1} s
        } \mu \\
        &=
        \SP{
          \inss(\D x^i) \nabla_{\frac{\partial}{\partial x^i}} X, 
          \omega \tensor \left(\SymD^E\right)^{k-1} s
        } \mu
        + 
        \SP{
          \inss\left(-\Gamma_{i \ell}^i \D x^\ell\right) X,
          \omega \tensor \left(\SymD^E\right)^{k-1} s
        } \mu \\
        &\quad+ 
        \SP{
          \inss(\D x^i) X,
          \nabla_{\frac{\partial}{\partial x^i}} \omega \tensor
          \left(\SymD^E\right)^{k-1} s
          + \omega \tensor \nabla_{\frac{\partial}{\partial x^i}}
          \left(\SymD^E\right)^{k-1} s
        } \mu
        \\
        &\quad+ 
        \SP{
          \inss\left(
              \frac{\Lie_{\frac{\partial}{\partial x^i}} \mu}{\mu}
              \D x^i
          \right) X,
          \omega \tensor \left(\SymD^E\right)^{k-1} s
        } \mu \\
        &=
        \SP{
          \divergenz_\nabla (X),
          \omega \tensor \left(\SymD^E\right)^{k-1} s
        } \mu
        + 
        \SP{
          \inss(\alpha) X,
          \omega \tensor \left(\SymD^E\right)^{k-1} s
        } \mu \\
        &\quad+
        \SP{
          X,
          \SymD^{F^*} \omega
          \tensor \left(\SymD^E\right)^{k-1} s
        } \mu
        + 
        \SP{
          X,
          \omega \tensor \left(\SymD^E\right)^k s
        } \mu \\
        &=
        \SP{
          \divergenz_\mu (X),
          \omega \tensor \left(\SymD^E\right)^{k-1} s
        } \mu
        +
        \SP{
          X,
          \SymD^{F^*} \omega 
          \tensor \left(\SymD^E\right)^{k-1} s
        } \mu 
        +
        \SP{
          X,
          \omega \tensor \left(\SymD^E\right)^k s
        } \mu.
    \end{align*}
    Integrating this equality over $M$ gives immediately
    \[
    \int_M
    \SP{
      X,
      \omega \tensor \left(\SymD^E\right)^k s
    } \mu
    = - \int_M
    \SP{
      X,
      \SymD^{F^*} \omega \tensor \left(\SymD^E\right)^{k-1} s
    } \mu
    -
    \int_M
    \SP{
      \divergenz_\mu (X),
      \omega \tensor \left(\SymD^E\right)^{k-1} s
    } \mu.
    \tag{$*$}
    \]
    This result is now again true for general compactly supported
    sections by the above partition of unity argument. We claim now
    that for all $\ell \leq k$ we have
    \[
    \int_M
    \SP{
      X,
      \omega \tensor \left(\SymD^E\right)^k s
    } \mu
    = (-1)^\ell \sum_{r=0}^\ell \binom{\ell}{r} \int_M
    \SP{
      \divergenz_\mu^r (X),
      (\SymD^{F^*})^{\ell-r} \omega 
      \tensor \left(\SymD^E\right)^{k-\ell} s
    } \mu.
    \]
    Indeed, a simple induction gives this formula as we can
    successively apply ($*$)
    \begin{align*}
        &(-1)^\ell \sum_{r=0}^\ell \binom{\ell}{r} \int_M
        \SP{
          \divergenz_\mu^r (X),
          (\SymD^{F^*})^{\ell-r} \omega 
          \tensor \left(\SymD^E\right)^{k-\ell} s
        } \mu \\
        & =
        (-1)^\ell \sum_{r=0}^\ell \binom{\ell}{r}
        \Bigg(
        \int_M
        \SP{
          \divergenz_\mu^{r+1} (X),
          (\SymD^{F^*})^{\ell-r} \omega
          \tensor \left(\SymD^E\right)^{k-\ell-1} s
        } \mu \\
        &\quad\quad\quad+ 
        \int_M
        \SP{
          \divergenz_\mu^r (X),
          (\SymD^{F^*})^{\ell-r+1} \omega
          \tensor \left(\SymD^E\right)^{k-\ell+1} s
        } \mu
        \Bigg) \\
        & =
        \sum_{r=1}^{\ell+1} (-1)^{\ell+1} \binom{\ell}{r-1} \int_M
        \SP{
          \divergenz_\mu^r (X),
          (\SymD^{F^*})^{\ell-r+1} \omega
          \tensor \left(\SymD^E\right)^{k-\ell-1} s
        } \mu \\
        & \quad + \sum_{r=0}^\ell \binom{\ell}{r} \int_M
        \SP{
          \divergenz_\mu^r (X),
          (\SymD^{F^*})^{\ell+1-r} \omega
          \tensor \left(\SymD^E\right)^{k-(\ell+1)} s
        } \mu \\
        & =
        (-1)^{\ell+1} \int_M
        \SP{
          \divergenz_\mu^{\ell+1} (X),
          \omega \tensor \left(\SymD^E\right)^{k-\ell+1} s
        } \mu \\
        & \quad +\sum_{r=1}^\ell (-1)^{\ell+1}
        \left(
            \binom{\ell}{r-1} + \binom{\ell}{r}
        \right)
        \int_M
        \SP{
          \divergenz_\mu^{r} (X),
          (\SymD^{F^*})^{\ell+1-r} \omega 
          \tensor \left(\SymD^E\right)^{k-(\ell+1)} s
        } \mu \\
        & \quad + (-1)^{\ell+1} \int_M
        \SP{
          \divergenz_\mu (X),
          (\SymD^{F^*})^{\ell+1} \omega
          \tensor \left(\SymD^E\right)^{k-\ell+1} s
        } \mu \\
        & =
        (-1)^{\ell+1} \sum_{r=0}^{\ell+1} \binom{\ell+1}{r} \int_M
        \SP{
          \divergenz_\mu^{r} (X),
          (\SymD^{F^*})^{\ell+1-r} \omega
          \tensor
          \left(\SymD^E\right)^{k-(\ell+1)} s
        } \mu.
    \end{align*}
    In particular, for $k = \ell$ we obtain the formula
    \[
    \int_M
    \SP{
      X,
      \omega \tensor \left(\SymD^E\right)^k s
    } \mu
    =
    (-1)^k \sum_{r=0}^k \binom{k}{r} \int_M
    \SP{
      \divergenz_\mu^r (X),
      (\SymD^{F^*})^{k-r} \omega \tensor s
    } \mu
    \]
    with no derivatives acting on $s$ anymore. Thus we have computed
    the adjoint of $\stdrep (X)$. Indeed, collecting the pre-factors
    gives
    \begin{align*}
        \int_M (\stdrep(X) s) \; \mu
        &=
        \frac{1}{k!} \left(\frac{\hbar}{\I}\right)^k \int_M
        \SP{
          X,
          \omega \tensor \left(\SymD^E\right)^k s
        } \mu \\
        &=
        \frac{(-1)^k}{k!} \left(\frac{\hbar}{\I}\right)^k
        \sum_{r=0}^k \binom{k}{r} \int_M
        \SP{
          \divergenz_\mu^r (X),
          (\SymD^{F^*})^{k-r} \omega \tensor s
        } \mu \\
        &=
        \frac{(-1)^k}{k!} \left(\frac{\hbar}{\I}\right)^k
        \sum_{r=0}^k \binom{k}{r} \int_M
        (k-r)! \left(\frac{\hbar}{\I}\right)^{k-r}
        \left(\stdrep(\divergenz_\mu^r(X^\Trans)) \omega\right) (s)
        \; \mu \\
        &=
        (-1)^k \int_M
        \left(
            \stdrep
            \left(
                \sum_{r=0}^k \frac{1}{r!}
                \left(\frac{\hbar}{\I}\right)^r 
                \divergenz_\mu^r(X^\Trans)
            \right)
            \omega
        \right)
        (s)
        \; \mu \\
        &=
        (-1)^k \int_M
        \left(
            \stdrep(N^2 X^\Trans) \omega
        \right)
        (s)
        \; \mu,
    \end{align*}
    with $N$ as in \eqref{eq:Neumaier-Operator}.
\end{proof}
\begin{remark}
    \label{remark:Neumaier-Theorem}
    The reason for the unpleasant prefactor $(-1)^k$ is that we have
    not used a sesquilinear pairing. Indeed, if we have the situation
    as in Remark~\ref{remark:other-pairings} then we would have the
    following result: For simplicity we consider the scalar case only,
    i.e. $E = F = M \times \mathbb{C}$ are both the trivial line
    bundle hence $\Secinfty(E) = \Cinfty(M)$. Then consider
    \begin{equation}
        \label{eq:pairing-of-functions}
        \SP{\varphi, \psi}_\mu
        = \int_M \cc{\varphi} \psi \; \mu
    \end{equation}
    for $\varphi, \psi \in \Cinfty_0(M)$ instead of
    \eqref{eq:pairing-with-fixed-densitiy}. The additional complex
    conjugation uses the sign $(-1)^k$ to obtain
    \begin{equation}
        \label{eq:2}
        \stdrep(X)^\Trans = \stdrep(N^2 (\cc{X})^\Trans)
    \end{equation}
    for $X \in \Secinfty(\Sym^k TM)$ in this case. This also
    generalizes to the case of Hermitian vector bundles, see
    \cite{bordemann.neumaier.pflaum.waldmann:2003a} for an additional
    discussion.
\end{remark}


%% file: distribution.tex
%
%

In this section we introduce distributions as continuous linear
functionals and discuss several of their basic properties. In
particular, the behaviour under smooth maps and differential operators
will be discussed.

%
%

\subsection{Distributions and Generalized Sections}
\label{subsec:distribution-and-generalized-sections}

As in the well-known case of $M = \mathbb{R}^n$ we define
distributions as continuous linear functionals on the test function
spaces:
\begin{definition}[Distribution]
    \label{definition:distribution}
    \index{Distribution}%
    A distribution $u$ on $M$ is a continuous linear functional
    \begin{equation}
        \label{eq:distribution}
        u: \Cinfty_0(M) \longrightarrow \mathbb{C}.
    \end{equation}
    The space of all distributions is denoted by $\Cinfty_0(M)'$ or
    $\mathcal{D}'(M)$.
\end{definition}
\begin{remark}[Distributions]
    \label{remark:distributions}
    ~
    \begin{remarklist}
    \item \label{item:continuity-of-distributions} The continuity of
        course refers to the LF topology of $\Cinfty_0(M)$ as
        introduced in Theorem~\ref{theorem:inductive-limit-topology}.
        In particular, a linear functional is continuous if and only
        if for all compacta $K \subseteq M$ the restriction
        \index{Locally convex topology!CinfinityNull-topology@$\Cinfty_0$-topology}%
        \index{Continuity!distribution}%
        \begin{equation}
            \label{eq:restriction-of-distribution}
            u \at{\Cinfty_K(M)}: \Cinfty_K(M) \longrightarrow \mathbb{C}
        \end{equation}
        is continuous in the $\Cinfty_K$-topology. This is the case if
        and only if for all $\varphi \in \Cinfty_K(M)$ we have a
        constant $c > 0$ and $\ell \in \mathbb{N}_0$ such that
        \begin{equation}
            \label{eq:continuity-estimate-of-distribution}
            |u (\varphi)| \leq
            c \max_{\ell' \le \ell} \seminorm[K,\ell'] (\varphi).
        \end{equation}
        Analogously, we could have used the seminorms $p_{U,x,K,\ell}$
        avoiding the usage of a covariant derivative but taking a
        maximum over finitely many compacta in the domain of a
        chart. With the symbolic seminorms of
        Remark~\ref{remark:SymbolicSeminorms} we can combine this to
        \begin{equation}
            \label{eq:ContinuitySymbolic}
            |u (\varphi)| \leq
            c \seminorm[K,\ell] (\varphi).
        \end{equation}
        In the following, we shall mainly use this version of the
        continuity.  Since each $\Cinfty_K(M)$ is a Fr\'echet space,
        $u$ restricted to $\Cinfty_K(M)$ is continuous iff it is
        sequentially continuous. This gives yet another criterion: A
        linear functional $u: \Cinfty_0(M) \longrightarrow \mathbb{C}$
        is continuous iff for all $\varphi_n \in \Cinfty_0(M)$ with
        $\varphi_n \longrightarrow \varphi$ in the
        $\Cinfty_0$-topology we have
        \begin{equation}
            \label{eq:sequential-continuity-of-distributions}
            \index{Continuity!sequential}%
            u (\varphi_n) \longrightarrow u (\varphi).
        \end{equation}
    \item \label{item:order-of-distributions} The minimal $\ell \in
        \mathbb{N}_0$ such that
        \eqref{eq:continuity-estimate-of-distribution} is valid is
        called the \emph{local order} $\ord_K (u)$ of $u$ on $K$.
        Clearly, this is a quantity independent of the connection used
        for $\seminorm[K, \ell]$ and can analogously be obtained from
        the seminorms $\seminorm[U, x, K, \ell]$ as well. The
        independence follows at once from the various estimates
        between the seminorms as in the proof of
        Theorem~\ref{theorem:Ck-Topologie}. The \emph{total order} of
        $u$ is defined as
        \begin{equation}
            \label{eq:total-order}
            \index{Distribution!total order}%
            \index{Distribution!local order}%
            \ord (u) = \sup_K \ord_K(u) \in \mathbb{N}_0 \cup \{+\infty\},
        \end{equation}
        and the distributions of total order $\leq k$ are sometimes
        denoted by $\mathcal{D}'^k(M)$. Their union is denoted by
        $\mathcal{D}'_F(M)$ and called \emph{distributions of finite
          order}.  \index{Distribution!finite order}
    \item \label{item:finite-order-distributions} The distributions
        $\mathcal{D}'(M)$ as well as $\mathcal{D}'^k(M)$ and
        $\mathcal{D}'_F(M)$ are vector spaces. We have
        $\mathcal{D}'^k(M) \subseteq \mathcal{D}'^\ell(M)$ for $k \leq
        \ell$. It can be shown that already for $M = \mathbb{R}^n$ all
        the inclusions $\mathcal{D}'^k(M) \subseteq
        \mathcal{D}'^\ell(M) \subseteq \mathcal{D}'_F(M) \subseteq
        \mathcal{D}'(M)$ are proper.
    \item \label{item:finite-order-extend-to-Ck} If $u$ has order
        $\leq k$ it can be shown that $u$ extends uniquely to a
        continuous linear function
        \begin{equation}
            \label{eq:distributions-on-Ck}
            u: \Fun[\ell]_0(M) \longrightarrow \mathbb{C}
        \end{equation}
        with respect to the $\Fun[\ell]_0$-topology provided $\ell \ge
        k$. This follows essentially from the approximation
        Theorem~\ref{theorem:compact-and-smooth-sections-dense-in-others},
        see e.g \cite[Thm~2.16]{hoermander:2003a}.
    \end{remarklist}
\end{remark}
\begin{example}[$\delta$-functional]
    \label{example:delta-functional}
    \index{Delta-Functional@$\delta$-Functional}%
    For $p \in M$ the evaluation functional
    \begin{equation}
        \label{eq:delta-functional}
        \delta_p:
        \Cinfty_0(M) \ni \varphi \; \mapsto \; \varphi(p) \in \mathbb{C}
    \end{equation}
    is clearly continuous and has order zero. More generally, if $v_p
    \in T_p M$ is a tangent vector then
    \begin{equation}
        \label{eq:delta-functional-with-tangent-vector}
        v_p: \varphi \; \mapsto \; v_p(\varphi)
    \end{equation}
    is again continuous and has order one.
\end{example}
\begin{example}[Locally integrable densities]
    \label{example:locally-integrable-densities}
    \index{Density!locally integrable}%
    Let $\mu: M \longrightarrow \Dichten T^*M$ be a not necessarily
    continuous section. Then $\mu$ is called \emph{locally integrable}
    if for all charts $(U, x)$ and all $K \subseteq U$ the function
    $\mu_U$ in $\mu \at{U} = \mu_U |\D x^1 \wedge \cdots \wedge \D
    x^n|$ is integrable over $K$ with respect to the Lebesgue measure
    on $x(U)$. Since the $|\D x^1 \wedge \cdots \wedge \D x^n|$
    transform with the \emph{smooth} absolute value of the Jacobian of
    the change of coordinates, it follows at once that local
    integrability is intrinsically defined and it is sufficient to
    check it for an atlas and an exhausting sequence of compacta. It
    is then easy to see that
    \begin{equation}
        \label{eq:locally-integrable-density-as-distribution}
        \mu: \Cinfty_0(M) \ni \varphi
        \; \mapsto \; \int_M \varphi \; \mu
        \in \mathbb{C}
    \end{equation}
    is continuous. Indeed, if $K \subseteq M$ is compact then
    $\vol_\mu(K) = \int_K |\mu| < \infty$ is well-defined and we have
    \begin{equation}
        \label{eq:loc-int-densities-estimate}
        \left|
            \int_M \varphi \; \mu
        \right|
        \leq \vol_\mu(K) \seminorm[K,0] (\varphi).
    \end{equation}
    Note that $|\mu| = \sqrt{\cc{\mu} \mu}$ is well-defined as
    $1$-density and still locally integrable. In particular,
    \eqref{eq:locally-integrable-density-as-distribution} is a
    distribution of order zero.
\end{example}
\begin{remark}[Generalized densities]
    \label{remark:generalized-densities}
    \index{Density!generalized}%
    \index{Generalized density|see{Distribution}}%
    \index{Generalized function}%
    The last example shows that we can identify densities of quite
    general type (locally integrable) with certain distributions. For
    this reason, we call distributions also ``\emph{generalized
      densities}'', following e.g.  \cite{guillemin.sternberg:1990a,
      friedlander:1975a}. Note however that e.g. Hörmander takes a
    different point of view and treats distributions as
    ``\emph{generalized functions}''. In \cite{hoermander:2003a} a
    distribution is \emph{not} a continuous linear functional on
    $\Cinfty_0(M)$ but has a slightly different transformation
    behaviour under local diffeomorphisms. In fact, his generalized
    functions can be viewed as continuous linear functionals on
    $\Secinfty_0(\Dichten T^*M)$. To emphasize the generalized density
    aspect from now on we adopt the notation of
    \cite{guillemin.sternberg:1990a} and write
    \begin{equation}
        \label{eq:generalized-densities}
        \Sec[- \infty] (\Dichten T^*M)
        = \left\{
            u: \Cinfty_0(M) \longrightarrow \mathbb{C}
            \; \big| \;
            u
            \;
            \textrm{is linear and continuous}
        \right\}.
    \end{equation}
    This point of view will be very useful when we discuss the
    transformation properties of distributions. Later on, both
    versions will be combined anyway since we consider distributional
    sections of arbitrary vector bundles. Thus speaking of generalized
    functions will be non ambiguous.
\end{remark}
We can now generalize the notion of distributions to test sections
instead of test functions.
\begin{definition}[Generalized section]
    \label{definition:generalizedSections}
    \index{Generalized section}%
    \index{Distributional section|see{Generalized section}}%
    Let $E \longrightarrow M$ be a smooth vector bundle. Then a
    generalized section (or: distributional section) of $E$ is a
    continuous linear functional
    \begin{equation}
        \label{eq:generalized-section}
        s: \Secinfty_0(E^* \tensor \Dichten T^*M)
        \longrightarrow \mathbb{C}.
    \end{equation}
    The generalized sections will be denoted by $\Sec[-\infty](E)$.
\end{definition}
\begin{remark}
    \label{remark:generalized-density-notation-clash}
    Note that here we have some mild clash of notations since we
    defined a distribution already as a generalized density $u \in
    \Sec[-\infty](\Dichten T^*M)$ while a generalized density
    according to Definition~\ref{definition:generalizedSections} is a
    continuous linear functional
    \begin{equation}
        \label{eq:generalized-section-of-density-bundle}
        u:
        \Secinfty_0
        \left(
            \left(\Dichten T^*M\right)^* \tensor \Dichten T^*M
        \right)
        \longrightarrow \mathbb{C},
    \end{equation}
    and not $u: \Cinfty_0(M) \longrightarrow \mathbb{C}$. However, for
    any line bundle $L$ we have canonically $L^* \tensor L \simeq M
    \times \mathbb{C}$ hence we can (and will) canonically identify
    $\Secinfty_0((\Dichten T^*M)^* \tensor \Dichten T^*M)$ with
    $\Cinfty_0(M)$. Thus
    Definition~\ref{definition:generalizedSections} and
    Definition~\ref{definition:distribution} are consistent.

    Moreover, a section of $E$ is always a generalized section of $E$
    since for $s \in \Secinfty(E)$ we can integrate $\omega(s)$ with
    $\omega \in \Secinfty_0(E^* \tensor \Dichten T^*M)$ over $M$ and
    obtain a continuous linear functional which we can identify with
    an element in $\Sec[-\infty](E)$. In fact, the section $s$ is
    uniquely determined be the values $\int_M \omega(s)$ for all
    $\omega \in \Secinfty_0(E^* \tensor \Dichten T^*M)$ hence this is
    indeed an injection. Therefor we have
    \begin{equation}
        \label{eq:inclusion-of-sections-in-generalized-sections}
        \Secinfty(E) \subseteq \Sec[-\infty](E).
    \end{equation}
    More generally, we also have
    \begin{equation}
        \label{eq:inclusion-of-Cksections-in-gen-sections}
        \Sec(E) \subseteq \Sec[-\infty](E)
    \end{equation}
    for all $k \in \mathbb{N}_0$ by the same argument.
\end{remark}
\begin{remark}
    \label{remark:generalized-sections-for-fixed-density}
    If we choose a smooth positive density $\mu > 0$ then we can also
    identify $\Sec[-\infty](E)$ with the topological dual of
    $\Secinfty_0(E^*)$. Indeed, if $s \in \Sec[-\infty](E)$ then we
    can define
    \begin{equation}
        \label{eq:generalized-secion-for-fixed-density}
        I_\mu (s): \Secinfty_0(E^*) \ni \omega
        \; \mapsto \; s(\omega \tensor \mu) \in \mathbb{C},
    \end{equation}
    and clearly obtain an element $I_\mu(s) \in \Secinfty_0(E^*)'$ in
    the topological dual. The reason is that the map
    \begin{equation}
        \label{eq:tensor-with-density-is-continuous}
        \Secinfty_0(E^*) \ni \omega
        \; \mapsto \;
        \omega \tensor \mu \in \Secinfty_0(E^* \tensor \Dichten T^*M)
    \end{equation}
    is continuous in the $\Cinfty_0$-topology according to
    Proposition~\ref{proposition:modul-tensor-are-continuous} and
    Remark~\ref{remark:modul-tensor-are-cont-in-Ck}. Moreover, since
    \eqref{eq:tensor-with-density-is-continuous} is even a bijection
    with continuous inverse, we obtain an isomorphism
    \begin{equation}
        \label{eq:gensec-topdualofdualsections-isomorphism}
        I_\mu: \Sec[-\infty](E) \longrightarrow \Secinfty_0(E^*)'.
    \end{equation}
    In case of $M = \mathbb{R}^n$ one uses the Lebesgue measure $\D^n
    x \in \Secinfty(\Dichten T^* \mathbb{R}^n)$ to provide such an
    identification. Note however that
    \eqref{eq:gensec-topdualofdualsections-isomorphism} does not
    behave well under vector bundle morphisms as we shall see later
    since $\mu$ needs not to be invariant. Finally, if the choice of
    $\mu$ is clear from the context, we shall omit the symbol $I_\mu$
    and identify $\Sec[-\infty](E)$ directly with the dual space
    $\Secinfty(E^*)'$ to simplify our notation. This will frequently
    happen starting from Chapter~\ref{cha:LocalTheory}.
\end{remark}
\begin{remark}[Module structure]
    \label{remark:module-structure-on-generalized-section}
    \index{Generalized section!module structure}%
    The generalized sections $\Sec[-\infty](E)$ become a
    $\Cinfty(M)$-module via the definition
    \begin{equation}
        \label{eq:module-structure-on-generalized-sections}
        (f \cdot s)(\omega) = s (f\omega).
    \end{equation}
    Indeed $\omega \mapsto f \omega$ is $\Cinfty_0$-continuous and
    hence \eqref{eq:module-structure-on-generalized-sections} is
    indeed a continuous linear functional $f \cdot s \in
    \Sec[-\infty](E)$. The module property is clear.
\end{remark}
\begin{remark}[Order of generalized sections]
    \label{remark:order-of-generalized-section}
    \index{Generalized section!order}%
    The continuity of $s \in \Sec[-\infty](E)$ is again expressed
    using the seminorms of $\Secinfty(E^* \tensor \Dichten T^*M)$ in
    the following way. For every compactum $K \subseteq M$ there are
    constants $c > 0$ and $\ell \in \mathbb{N}_0$ such that
    \begin{equation}
        \label{eq:cont-estimate-for-gen-sections}
        |s(\omega)|
        \leq c \max_{\ell' \leq \ell} \seminorm[K, \ell'] (\omega),
    \end{equation}
    for all $\omega \in \Secinfty_K(E^* \tensor \Dichten
    T^*M)$. Again, the \emph{local order} of $s$ on $K$ is defined to
    be the smallest $\ell$ such that
    \eqref{eq:cont-estimate-for-gen-sections} holds. This also defines
    the \emph{global order}
    \begin{equation}
        \label{eq:global-order-of-gen-section}
        \ord(s) = \sup_K \ord_K (s)
    \end{equation}
    as before. As in the scalar case, a generalized section $s \in
    \Sec[-\infty](E)$ with global order $\ord(s) \leq k$ extends
    uniquely to a $\Fun[\ell]_0$-continuous functional
    \begin{equation}
        \label{eq:extension-of-finite-order-sections}
        s: \Sec[\ell]_0(E^* \tensor \Dichten T^*M) \longrightarrow \mathbb{C}
    \end{equation}
    for all $\ell \geq k$.  We shall denote the distributional
    sections of order $\leq \ell$ by $\Sec[-\ell](E)$.  Note that
    $\Sec[-0](E)$ are \emph{not} just the continuous sections.
\end{remark}

We also want to topologize the distributions. Here we use the most
simple locally convex topology: the weak$^*$ topology:
\begin{definition}[Weak$^*$ topology]
    \label{definition:weak-star-topology}
    \index{Weakstar topology@Weak$^*$ topology}%
    \index{Locally convex topology!weakstar@weak$^*$}%
    \index{Seminorm!weakstar@weak$^*$}%
    The weak$^*$ topology for $\Sec[-\infty](E)$ is the locally convex
    to\-po\-lo\-gy obtained from all the seminorms
    \begin{equation}
        \label{eq:weak-top-seminorms}
        \seminorm[\omega] (s) = |s(\omega)|,
    \end{equation}
    where $\omega \in \Secinfty_0(E^* \tensor \Dichten T^*M)$.
\end{definition}
In the following we always use the $\mathrm{weak^*}$ topology for
$\Sec[-\infty](E)$. We have the following properties:
\begin{theorem}[Weak$^*$ topology of \protect{$\Sec[-\infty](E)$}]
    \label{theorem:weak-topology-on-gensecs}
    \index{Weakstar topology@Weak$^*$ topology!sequential completeness}%
    \index{Continuity!weakstar@weak$^*$}%
    ~
    \begin{theoremlist}
    \item \label{item:convergence-of-gensecs} A sequence $s_n
        \in \Sec[-\infty](E)$ converges to $s \in \Sec[-\infty](E)$ if
        and only if for all $\omega \in \Secinfty_0(E^* \tensor
        \Dichten T^*M)$
        \begin{equation}
            \label{eq:weak-star-convergence}
            s_n (\omega) \longrightarrow s(\omega).
        \end{equation}
    \item \label{item:gensecs-are-sequentially-comlete}
        $\Sec[-\infty](E)$ is sequentially complete, i.e. every
        weak$^*$ Cauchy sequence converges.
    \item \label{item:inclusion-sec-into-gensec-is-cont} The
        inclusions $\Sec(E) \subseteq \Sec[-\infty](E)$ are continuous
        in the $\Fun[k]$- and weak$^*$ topology for all $k \in
        \mathbb{N}_0 \cup \{+\infty\}$.
    \item \label{item:gensec-modulstructure-is-cont} The map
        $\Sec[-\infty](E) \ni s \mapsto fs \in \Sec[-\infty](E)$ is
        weak$^*$ continuous for all $f \in \Cinfty(M)$.
    \item \label{item:sec-dense-in-gensec} The sections
        $\Secinfty_0(E)$ are sequentially weak$^*$ dense in
        $\Sec[-\infty](E)$.
    \end{theoremlist}
\end{theorem}
\begin{proof}
    The first part is clear since $s_n \longrightarrow s$ means for
    every seminorm $\seminorm[\omega]$ we have
    \[
    \seminorm[\omega] (s_n -s) \longrightarrow 0,
    \]
    which is \eqref{eq:weak-star-convergence}. Thus the notion of
    convergence in $\Sec[-\infty](E)$ is pointwise convergence on the
    test sections $\Secinfty_0(E^* \tensor \Dichten T^*M)$.  The
    second part is non-trivial but follows from general arguments:
    first one shows that the topological dual $V'$ of a Fr\'echet
    space $V$ is sequentially complete by a Banach-Steinhaus argument.
    Here Fr\'echet is crucial.  Second, one extends this result to LF
    spaces like our $\Secinfty_0(E^* \tensor \Dichten T^*M)$, see e.g.
    \cite[Thm.~2.1.8]{hoermander:2003a} or
    \cite[Thm.~6.17]{rudin:1991a} for details.  Note however that
    $\Sec[-\infty](E)$ is \emph{not} complete; in fact, the completion
    is the full algebraic dual \cite[p.147]{jarchow:1981a}.  The third
    part is easy since for a $\Fun$-section $s \in \Sec(E)$ we have
    for all $\omega \in \Secinfty_0(E^* \tensor \Dichten T^*M)$
    \[
    \seminorm[\omega] (s) = |s(\omega)|
    = \left| \int_M \omega(s) \right| \leq c \seminorm[K, 0](s),
    \]
    with some constant $c > 0$ depending on $\omega$ but not on $s$
    and a compactum $\supp \omega \subseteq K$. Essentially, $c$ is
    the volume of $K$ times the maximum of $\omega$ with respect to
    the metrics used to define $\seminorm[K, 0]$. From this the
    continuity is obvious. For the fourth part we compute
    \[
    \seminorm[\omega] (fs)
    = |fs(\omega)|
    = |s(f \omega)|
    = \seminorm[f\omega] (s),
    \]
    which already shows the continuity. The last part is slightly more
    tricky. We have to construct a sequence $s_n \in \Secinfty_0(E)$
    with $s_n \longrightarrow s$ in the weak$^*$ topology using of
    course the identification of $s_n$ with an element of
    $\Sec[-\infty](E)$. We choose a countable atlas of charts $(U_n,
    x_n)$ and a partition of unity $\chi_n$ subordinate to this atlas.
    Then we consider the distributions $\chi_n s \in
    \Sec[-\infty](E)$. We claim that
    \[
    \sum_{n=0}^\infty \chi_n s = s
    \]
    in the weak$^*$ topology. To prove this, let $\omega \in
    \Secinfty_0(E^* \tensor \Dichten T^*M)$ be given and let $K =
    \supp \omega$. Then only finitely many $\chi_n$ are nonzero on
    $K$, hence
    \[
    \sum_n (\chi_n s) (\omega)
    = \sum_n^{\mathrm{finite}} s(\chi_n \omega)
    = s \left( \sum_n^{\mathrm{finite}} \chi_n \omega \right)
    = s(\omega).
    \]
    This proves convergence. Since the $\chi_n s $ are countable, it
    is sufficient to prove that each $\chi_n s$ can be approximated by
    a sequence of sections in $\Secinfty_0(E)$. Since $\supp \chi_n
    \subseteq U_n$ we also conclude that $(\chi_n s) (\omega) = 0$ if
    $\supp \omega \cap U_n = \emptyset$. Thus we are left with the
    problem to approximate a distribution on a chart which can be done
    by some appropriate convolution, see
    e.g. \cite[Thm.~6.32]{rudin:1991a}.
\end{proof}

\begin{remark}[Weak$^*$ topology of \protect{$\Sec[-\infty](E)$}]
    \label{remark:weakstar-topology}
    ~
    \begin{remarklist}
    \item \label{item:gensec-not-frechet} It should be noted that
        $\Sec[-\infty](E)$ is not Fr\'echet, in fact it is not
        metrizable. Thus sequential completeness is weaker than
        completeness: $\Sec[-\infty](E)$ is not complete and its
        completion is the full \emph{algebraic} dual of
        $\Secinfty_0(E^* \tensor \Dichten T^*M)$.
    \item \label{item:convergence-of-sec-in-weakstar} The importance
        of continuity of the inclusion is that for sections $s_n \in
        \Sec(E)$ with $s_n \longrightarrow s$ in the $\Fun$-topology
        we also have $s_n \longrightarrow s$ in the weak$^*$ topology
        of $\Sec[-\infty](E)$ for all $k \in \mathbb{N}_0 \cup
        \{+\infty\}$.
    \item \label{relation-genset-to-sec} The last part shows that
        $\Sec[-\infty](E)$ is, on one hand, a large extension of
        $\Secinfty_0(E)$ and also $\Secinfty(E)$ which, on the other
        hand, is still not ``too large'': continuous operations with
        distributions are already determined by their restrictions to
        $\Secinfty_0(E)$. This justifies the name ``generalized
        section''.
    \end{remarklist}
\end{remark}

%
%

\subsection{Calculus with Distributions}
\label{subsec:calc-with-distributions}

In this subsection we shall extend various constructions with sections
to generalized sections. The main idea is to ``dualize'' continuous
linear operations on test sections in an appropriate way.

We begin with the definition of the support of a distribution and its
restriction to open subsets.
\begin{definition}[Restriction and support]
    \label{definition:support-of-gensec}
    \index{Generalized section!restriction}%
    \index{Generalized section!support}%
    Let $U \subseteq M$ be open and $s \in \Sec[-\infty](E)$.
    \begin{definitionlist}
    \item \label{item:restriction-of-gensec} The restriction $s
        \at{U}$ is defined by
        \begin{equation}
            \label{eq:restriction-of-gensec}
            s\at{U} (\omega) = s(\omega)
        \end{equation}
        for $\omega \in \Secinfty_0(E^* \tensor \Dichten T^*M
        \at{U})$, i.e. for $\omega \in \Secinfty_0(E^* \tensor
        \Dichten T^*M)$ with $\supp \omega \subseteq U$.
    \item \label{item:support-of-gensec} The support of $s$ is defined by
        \begin{equation}
            \label{eq:support-of-gensec}
            \supp s = \bigcap_{\substack{A \subseteq M \textrm{closed} \\
                s|_{M \backslash A} = 0}} A.
        \end{equation}
    \end{definitionlist}
\end{definition}
\begin{remark}[Restriction and Support]
    \label{remark:restriction-and-support-of-gensec}
    ~
    \begin{remarklist}
    \item \label{item:gensec-are-presheaf} It is easy to show that
        $s\at{U} \in \Sec[-\infty](E\at{U})$. Moreover, we clearly have
        \begin{equation}
            \label{eq:restriction-of-restricion}
            (s\at{U})\at{V} = s\at{V},
        \end{equation}
        for $V \subseteq U$. In more sophisticated terms this means
        that $\Sec[-\infty](E)$ has the structure of a
        \emph{presheaf}\index{Presheaf} over $M$ with values in
        locally convex vector spaces.
    \item \label{item:gensec-are-sheaf} If $U_\alpha \subseteq M$ is an
        open cover of $M$ and if we have $s_\alpha \in
        \Sec[-\infty](E\at{U_\alpha})$ given such that
        \begin{equation}
            \label{eq:patching-condition-for-local-gensecs}
            s_\alpha \at{U_\alpha \cap U_\beta}
            = s_\beta \at{U_\alpha \cap U_\beta},
        \end{equation}
        whenever $U_\alpha \cap U_\beta \neq \emptyset$ then there
        exists a unique $s \in \Sec[-\infty](E)$ with $s\at{U_\alpha}
        = s_\alpha$. The proof of this fact uses a partition of unity
        argument to glue together the locally defined $s_\alpha$. In
        fact, if $\chi_\alpha$ is a subordinate partition of unity one
        checks that the definition
        \begin{equation}
            \label{eq:glueing-definition}
            s(\omega) = \sum_\alpha s_\alpha (\chi_\alpha \omega)
        \end{equation}
        indeed gives the desired $s$, independent of the choice of the
        partition of unity. Moreover, if $s, t \in \Sec[-\infty](E)$
        are given then
        \begin{equation}
            \label{eq:local-agreement}
            s\at{U_\alpha} = t\at{U_\alpha}
        \end{equation}
        for all $\alpha$ implies $s=t$. This is obvious. Again, with
        more high-tech language this means that $\Sec[-\infty](E)$ is
        in fact a \emph{sheaf}\index{Sheaf} and not only a presheaf.
    \item \label{item:vanish-out-of-supp} The support $\supp s$ of $s
        \in \Sec[-\infty](E)$ is the smallest closed subset with
        $s\at{M \backslash \supp s} = 0$ and we have $p \in \supp s$
        if and only if for every open neighborhood $U$ of $p$ we find
        $\omega \in \Secinfty_0(E^* \tensor \Dichten T^*M)$ with
        $\supp \omega \subseteq U$ and $s(\omega) \neq 0$.
    \item \label{item:supp-of-product} For $s \in \Sec[-\infty](E)$,
        $f \in \Cinfty(M)$, $t \in \Sec[0](E)$ and $\omega \in
        \Secinfty_0(E^* \tensor \Dichten T^*M)$ we have
        \begin{equation}
            \label{eq:support-of-product}
            \supp (fs) \subseteq \supp f \cap \supp s
        \end{equation}
        \begin{equation}
            \label{eq:gensec-zero-if-supps-dont-overlap}
            s(\omega) = 0 \quad \textrm{if} \quad
            \supp s \cap \supp \omega = \emptyset,
        \end{equation}
        and the support of $t$ as a \emph{continuous section} in
        $\Sec[0](E)$ coincides with the support of $t$ viewed as
        distribution. Thus the notion of support has the usual
        properties as known from continuous or smooth sections.
    \end{remarklist}
\end{remark}

After the support we also have a more refined notion, namely the
\emph{singular support}. It characterizes where a generalized section
is not just a smooth section but actually ``singular''.
\begin{definition}[Singular support]
    \label{definition:singsupp}
    \index{Generalized section!regular point}%
    \index{Generalized section!singular support}%
    Let $s \in \Sec[-\infty](E)$.
    \begin{definitionlist}
    \item \label{item:regular-gensec} $s$ is called regular in $p \in
        M$ if there is an open neighborhood $U \subseteq M$ of $p$
        such that
        \begin{equation}
            \label{eq:regular-condition}
            s\at{U} \in \Secinfty(E\at{U}).
        \end{equation}
    \item \label{item:singsupp} The singular support of $s$ is
        \begin{equation}
            \label{eq:singsupp}
            \singsupp s =
            \left\{
                p \in M \; | \; s \; \textrm{is not regular in} \; p
            \right\}.
        \end{equation}
    \end{definitionlist}
\end{definition}
The singular support of $s$ indeed behaves similar to the support.
\begin{remark}[Singular support]
    \label{remark:singsupp}
    Let $s \in \Sec[-\infty](E)$, $t \in \Secinfty(E)$ and $f \in
    \Cinfty(M)$.
    \begin{remarklist}
    \item \label{regular-outside-singsupp} The singular support
        $\singsupp s$ is the smallest closed subset of $M$ with
        \begin{equation}
            \label{eq:gensec-on-complement-of-singsupp-is-smooth}
            s\at{M \backslash \singsupp s} \in \Secinfty(E).
        \end{equation}
        This follows easily from the fact that smooth sections are
        determined by their restrictions to open subsets and by
        \eqref{eq:patching-condition-for-local-gensecs} in
        Remark~\ref{remark:restriction-and-support-of-gensec}.
    \item \label{item:singsupp-properties} We have
        \begin{equation}
            \label{eq:singsupp-in-supp}
            \singsupp s \subseteq \supp s,
        \end{equation}
        \begin{equation}
            \label{eq:singsupp-of-product}
            \singsupp (fs) \subseteq \singsupp s,
        \end{equation}
        and
        \begin{equation}
            \label{eq:singsupp-of-smooth-section}
            \singsupp t = \emptyset.
        \end{equation}
        Again these properties follow in a rather straightforward way
        from the very definition.
    \end{remarklist}
\end{remark}
Having a notion of support of distributions it is interesting to
consider those elements of $\Sec[-\infty](E)$ with \emph{compact}
support. The following theorem gives a full description:
\begin{theorem}[Generalized sections with compact support]
    \label{theorem:gensec-with-compact-support}
    \index{Generalized section!compact support}%
    Let $s \in \Sec[-\infty](E)$ have compact support. Then we have:
    \begin{theoremlist}
    \item \label{item:compact-supp-finite-order} $s$ has finite global
        order $\ord(s) < \infty$.
    \item \label{item:compact-supp-extension-to-Cinfty} $s$ has a
        unique extension to a linear functional
        \begin{equation}
            \label{eq:extension-by-compact-supp}
            s: \Secinfty(E^* \tensor \Dichten T^*M)
            \longrightarrow \mathbb{C},
        \end{equation}
        which is continuous in the $\Cinfty$-topology.
    \end{theoremlist}
    Conversely, if $s:\Secinfty(E^* \tensor \Dichten T^*M)
    \longrightarrow \mathbb{C}$ is a continuous linear functional then
    its restriction to $\Secinfty_0(E^* \tensor \Dichten T^*M)$ is a
    generalized section of $E$ with compact support.
\end{theorem}
\begin{proof}
    Thanks to the compactness of $\supp s$ we can find an open
    neighborhood $U$ of $\supp s$ such that $U^\cl \subseteq M$ is
    still compact. Hence there is a $\chi \in \Cinfty_0(M)$ with $\chi
    \at{U^\cl} = 1$. It follows from \eqref{eq:support-of-product}
    that
    \[
    \chi s = s.
    \]
    For $K = \supp \chi$ we find some $\ell \in \mathbb{N}_0$ and $c >
    0$ such that for all $\omega \in \Secinfty_K(E^* \tensor \Dichten
    T^*M)$ we have
    \[
    |s(\omega)| \leq c \seminorm[K, \ell] (\omega),
    \]
    since $s$ is continuous with the seminorms of
    Remark~\ref{remark:SymbolicSeminorms}. If $\omega \in
    \Secinfty_0(E^* \tensor \Dichten T^*M)$ is arbitrary we have $\chi
    \omega \in \Secinfty_K(E^* \tensor \Dichten T^*M)$, hence
    \[
    |s(\omega)| = |s(\chi \omega)|
    \leq c \seminorm[K, \ell](\chi \omega)
    \leq c' \seminorm[K, \ell](\omega)
    \leq c' \seminorm[M, \ell](\omega)
    \]
    by the Leibniz rule and the compactness of $\supp \omega$. From
    this we immediately see that $s$ has global order $\ord(s) \leq
    \ell$. For the second part consider $\omega \in \Secinfty(E^*
    \tensor \Dichten T^*M)$ then $\chi \omega$ has compact support and
    we can set
    \[
    s(\omega) = s(\chi \omega).
    \]
    This clearly provides a linear extension of $s$ and since $\supp
    (\chi \omega) \subseteq K$ we have
    \[
    |s(\omega)| = |s(\chi \omega)|
    \leq c \seminorm[K,\ell](\chi \omega)
    \leq c' \seminorm[K,\ell](\omega),
    \]
    which is the continuity in the $\Cinfty$-topology. Thus $s$ is a
    continuous extension. Since
    \[
    \Secinfty_0(E^* \tensor \Dichten T^*M)
    \subseteq \Secinfty(E^* \tensor \Dichten T^*M)
    \]
    is dense by
    Proposition~\ref{proposition:compact-sections-dense-in-smooth-sections},
    such an extension is necessarily unique.  Now let $s:
    \Secinfty(E^* \tensor \Dichten T^*M) \longrightarrow \mathbb{C}$
    be linear and continuous in the $\Cinfty$-topology. Then there
    exists a compactum $K \subseteq M$ and $\ell \in \mathbb{N}_0$, $c
    > 0$ with
    \[
    |s(\omega)| \leq c \seminorm[K,\ell](\omega)
    \]
    for all $\omega \in \Secinfty(E^* \tensor \Dichten T^*M)$. From
    this it follows easily that $s \at{\Secinfty_{K'}(E^* \tensor
      \Dichten T^*M)}$ is continuous in the $\Cinfty_{K'}$-topology
    for all compacta $K'$.  Moreover, for $\supp \omega \cap K =
    \emptyset$ we have $s(\omega) = 0$, hence $\supp s \subseteq K$
    follows.
\end{proof}

\begin{definition}
    \label{definition:gensec-with-compact-support}
    The generalized sections of $E$ with compact support are denoted
    by $\Sec[-\infty]_0(E)$.
\end{definition}

After having identified the distributions with compact support we can
extend this construction under slightly milder assumptions: if only
the overlap $\supp s \cap \supp \omega$ is compact then the pairing
$s(\omega)$ is already well-defined:

\begin{proposition}
    \label{proposition:supp-overlap-is-compact}
    Let $s \in \Sec[-\infty](E)$ be a generalized section. Then there
    exists a unique extension $\widetilde{s}$ of $s$ to a linear
    functional
    \begin{equation}
        \label{eq:extension-for-supp-overlap-is-compact}
        \widetilde{s}:
        \left\{
            \omega \in \Secinfty(E^* \tensor \Dichten T^*M)
            \; \big| \;
            \supp \omega \cap \supp s \mathrm{\;is\; compact}
        \right\}
        \longrightarrow \mathbb{C},
    \end{equation}
    such that
    \begin{propositionlist}
    \item \label{item:extension-compact-overlap-is-extension}
        $\widetilde{s}$ coincides with $s$ on $\Secinfty_0(E^* \tensor
        \Dichten T^*M)$,
    \item \label{item:supp-overlap-zero-then-sw-zero}
        $\widetilde{s}(\omega) = 0$ if $\supp s \cap \supp \omega =
        \emptyset$.
    \end{propositionlist}
\end{proposition}
\begin{proof}
    Assume first that $\widetilde{s}'$ is another such extension and
    let $\omega$ be a test sections as in
    \eqref{eq:extension-for-supp-overlap-is-compact}. Then we choose a
    cut-off function $\chi \in \Cinfty(M)$ with $\chi=1$ on an open
    neighborhood $U$ of $K = \supp s \cap \supp \omega$. Thus $\omega
    = \chi \omega + (1-\chi)\omega$ with $\chi \omega$ having compact
    support and $\supp (1-\chi)\omega \cap \supp s = \emptyset$. Hence
    for the extension $\widetilde{s}$ we get by linearity and
    \refitem{item:extension-compact-overlap-is-extension} and
    \refitem{item:supp-overlap-zero-then-sw-zero}
    \begin{align*}
        \widetilde{s}(\omega)
        = \widetilde{s}( \chi \omega + (1-\chi)\omega)
        = \widetilde{s}(\chi \omega) + \widetilde{s}((1-\chi) \omega)
        = s(\chi \omega).
    \end{align*}
    The same arguments hold for $\widetilde{s}'$ whence
    $\widetilde{s}'(\omega) = s(\chi \omega) = \widetilde{s}(\omega)$
    follows. This shows that such an extension is necessarily unique.
    To show existence we simply define $\widetilde{s}(\omega) = s(\chi
    \omega)$ where $\chi$ is chosen as above. Clearly, two different
    choices of $\chi$ lead to the same extension by the above
    uniqueness argument. Since for $\omega, \omega'$ we can find a
    common $\chi$, satisfying the requirements with respect to both
    $\omega$ and $\omega'$, we see that the above definition is
    linear. For $\supp \omega$ compact we find a $\chi$ with $\chi
    \omega = \omega$ whence
    \refitem{item:extension-compact-overlap-is-extension} follows.
    Finally, if $\supp s \cap \supp \omega = \emptyset$ then $\chi =
    0$ will do the job and so
    \refitem{item:supp-overlap-zero-then-sw-zero} holds.
\end{proof}

One can also put a certain locally convex topology on the vector space
of such test functions such that the extension is actually
continuous. In the following we will denote this extension simply by
$s$.
\begin{remark}
    \label{remark:extension-for-finite-order}
    A slight variation of this proposition is the following. If
    $\ord_K s \leq \ell$ for some compact subset $K$ then $s$ extends
    uniquely to a linear functional
    \begin{equation}
        \label{eq:extension-for-finite-order}
        s:
        \left\{
            \omega \in \Sec[\ell](E^* \tensor \Dichten T^*M)
            \; \big| \;
            \supp \omega \cap \supp s \subseteq K
        \right\}
        \longrightarrow \mathbb{C},
    \end{equation}
    such that
    \begin{remarklist}
    \item $\widetilde{s}$ coincides with the continuous extension of
        $s$ to $\Sec[\ell]_K(E^* \tensor \Dichten T^*M)$ on those
        $\omega$ with $\supp \omega \subseteq K$.
    \item $s(\omega)=0$ if $\supp \omega \cap \supp s = \emptyset$.
    \end{remarklist}
\end{remark}

After the discussion of supports we can now move distributions around
by using smooth maps between manifolds and vector bundle
morphisms. The latter one clearly includes the case of smooth maps by
viewing smooth functions as sections of the trivial line bundle and
extending a smooth map in the unique way to a vector bundle morphism
of the trivial line bundles.

Thus let $E \longrightarrow M$ and $F \longrightarrow M$ by vector
bundles and let $\Phi: E \longrightarrow F$ be a smooth vector bundle
morphism over the smooth map $\phi: M \longrightarrow N$. We can now
obtain pull-backs and push-forwards of distributions by dualizing the
statements of the
Propositions~\ref{proposition:pullback-of-sections-is-continous} and
\ref{proposition:pulback-of-compact-sections} appropriately. We start
with the scalar case:
\begin{definition}[Push-forward of distributions]
    \label{definition:pushforward-of-distributions}
    \index{Distribution!push-forward}%
    Let $\phi: M \longrightarrow N$ be a smooth map. The push-forward
    of compactly supported generalized densities
    \begin{equation}
        \label{eq:pushforward-of-gendensities}
        \phi_*: \Sec[-\infty]_0(\Dichten T^*M)
        \longrightarrow \Sec[-\infty]_0(\Dichten T^*N)
    \end{equation}
    is defined on $f \in \Cinfty(M)$ by
    \begin{equation}
        \label{eq:formula-for-pushforward-of-gendensities}
        (\phi_* \mu) (f) = \mu (\phi^* f).
    \end{equation}
\end{definition}
\begin{proposition}[Push-forward of distributions]
    \label{proposition:push-forward-of-distributions}
    \index{Push-forward!weakstar continuity@weak$^*$ continuity}%
    \index{Proper map}%
    Let $\phi: M \longrightarrow N$ be a smooth map.
    \begin{propositionlist}
    \item \label{item:pushforward-of-distribution} The push-forward
        $\phi_* \mu$ of $\mu \in \Sec[-\infty]_0(\Dichten T^*M)$ is a
        well-defined generalized density with compact support
        \begin{equation}
            \label{eq:push-forward-of-gengensitiy-is-gendensity}
            \phi_* \mu \in \Sec[-\infty]_0(\Dichten T^*N).
        \end{equation}
        The map $\phi_*$ is linear and continuous with respect to the
        weak$^*$ topologies.
    \item \label{item:proper-pushforward} Assume $\phi$ is in addition
        proper. Then the push-forward extends uniquely to
        $\Sec[-\infty](\Dichten T^*M)$ and gives a linear continuous
        map
        \begin{equation}
            \label{eq:proper-push-forward-of-gendens}
            \phi_*: \Sec[-\infty](\Dichten T^*M)
            \longrightarrow \Sec[-\infty](\Dichten T^*N)
        \end{equation}
        with respect to the $\textrm{weak}^*$ topologies. Explicitly,
        for all $\varphi \in \Cinfty_0(N)$ the push-forward
        $\phi_*\mu$ of $\mu$ is given by
        \begin{equation}
            \label{eq:formula-for-proper-push-forward-of-gendens}
            (\phi_* \mu) (\varphi) = \mu (\phi^* \varphi).
        \end{equation}
    \item \label{item:pushforward-functorial} We have
        \begin{equation}
            \label{eq:functorial-properties-of-push-forward-of-gendens}
            \left( \id_M \right)_*
            = \id_{\Sec[-\infty](\Dichten T^*M)}
            \quad
            \textrm{and}
            \quad
            \left( \phi \circ \psi \right)_* = \phi_* \circ \psi_*.
        \end{equation}
    \end{propositionlist}
\end{proposition}
\begin{proof}
    Since by Proposition~\ref{proposition:pullback-of-functions} the
    pull-back $\phi^*: \Cinfty(N) \longrightarrow \Cinfty(M)$ is
    $\Cinfty$-continuous, by
    \eqref{eq:formula-for-pushforward-of-gendensities} one obtains a
    well-defined transpose map of $\phi^*$ which ---consequently--- is
    denoted by $\phi_*$. Clearly, $\phi_*$ is linear and
    \[
    \seminorm[f](\phi_* \mu) = |\phi_* \mu (f)| = |\mu (\phi^* f)| =
    \seminorm[\phi^* f](\mu)
    \]
    shows immediately that $\phi_*$ is weak$^*$ continuous.  The
    second part follows analogously, now using
    Proposition~\ref{proposition:pullback-with-proper-map} instead.
    The uniqueness of this extension follows since $\phi_* \mu$ is
    continuous and since the compactly supported distributions
    $\Sec[-\infty]_0(\Dichten T^*M)$ are sequentially dense in
    $\Sec[-\infty](\Dichten T^*M)$. The later follows from
    Theorem~\ref{theorem:weak-topology-on-gensecs},
    \refitem{item:sec-dense-in-gensec} since already
    $\Secinfty_0(\Dichten T^*M) \subseteq \Sec[-\infty]_0(\Dichten
    T^*M) \subseteq \Sec[-\infty](\Dichten T^*M)$ is sequentially
    dense. The last part is obvious and follows immediately from the
    corresponding properties of the pull-back of functions.
\end{proof}

\begin{remark}[Push-forward of smooth densities]
    \label{remark:pushforward-of-smooth-densities}
    \index{Density!push-forward}%
    Since by Remark~\ref{remark:restriction-and-support-of-gensec},
    \refitem{item:supp-of-product} we have $\Secinfty_0(\Dichten T^*M)
    \subseteq \Sec[-\infty]_0(\Dichten T^*M)$, we can always
    push-forward compactly supported smooth densities in the sense of
    generalized densities by
    \eqref{eq:formula-for-pushforward-of-gendensities}.  However, even
    though $\mu$ is smooth, $\phi_* \mu$ needs not to be smooth at
    all. A simple example is obtained as follows: Let $\iota: C
    \longrightarrow M$ be a submanifold of positive codimension and
    let $\mu \in \Secinfty(\Dichten T^*C)$ be a smooth density on $C$.
    Then for $f \in \Cinfty_0(M)$ we have
    \begin{equation}
        \label{eq:pushforward-of-smooth-densitiy-on-submanifold}
        \iota_* \mu (f) = \int_C \iota^*f \: \mu,
    \end{equation}
    which can not be written as $\int_M f \: \nu$ with some \emph{smooth}
    $\nu \in \Secinfty(\Dichten T^*M)$.
    \begin{figure}
        \centering
        \input{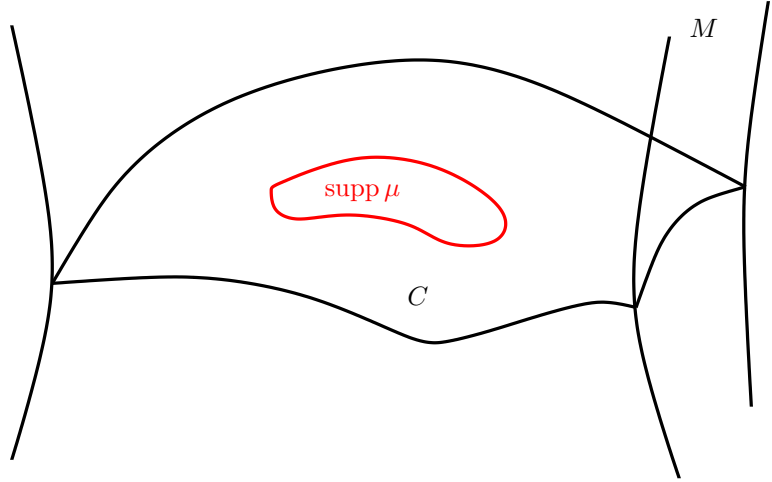}
        \caption{\label{fig:singsupp}%
          The push-forward has now singular support.
        }
    \end{figure}
    In fact, one can show rather easily that
    \begin{equation}
        \label{eq:pushforward-from-proper-submanifold-not-nice}
        \supp \iota_* \mu = \singsupp \iota_* \mu = \iota (\supp \mu)
    \end{equation}
    in this case, see also Figure~\ref{fig:singsupp}. The simplest
    case of this class of examples is given by $C = \{\mathrm{pt}\}$
    and $\mu = \delta_{\mathrm{pt}}$ the evaluation functional on
    $\Cinfty(\mathrm{pt})=\mathbb{C}$. On $C$, the $\delta$-functional
    is actually a \emph{smooth} density but on any higher dimensional
    manifold this is of course no longer the case.
\end{remark}
\begin{remark}
    \label{remark:vector-valued-pushforward}
    \index{Generalized section!push-forward}%
    \index{Vector bundle morphism}%
    There is also a vector-valued version of push-forward. Since for a
    vector bundle morphism $\Phi: E \longrightarrow F$ over $\phi: M
    \longrightarrow N$ we have a continuous pull-back
    \begin{equation}
        \label{eq:pull-back-of-smooth-sections}
        \Phi^*: \Secinfty(F^*) \longrightarrow \Secinfty(E^*),
    \end{equation}
    this dualizes to a push-forward
    \begin{equation}
        \label{eq:vectorvalued-pushforward}
        \Phi_*: \Sec[-\infty]_0(E \tensor \Dichten T^*M)
        \longrightarrow \Sec[-\infty]_0(F \tensor \Dichten T^*N)
    \end{equation}
    being again linear and weak$^*$ continuous. In case $\phi$ is
    proper we get an extension
    \begin{equation}
        \label{eq:proper-vectorvalued-pushforward}
        \Phi_*: \Sec[-\infty](E \tensor \Dichten T^*M)
        \longrightarrow \Sec[-\infty](F \tensor \Dichten T^*N),
    \end{equation}
    which is again linear, unique and weak$^*$ continuous. In general,
    a smooth section of $E \tensor \Dichten T^*M$ is pushed forward to
    a singular section of $F \tensor \Dichten T^*N$. Note however,
    that there are conditions on $\Phi$ and $\phi$ such that $\Phi_*
    s$ is again smooth for a smooth $s$, see e.g. the discussion in
    \cite[p.~307]{guillemin.sternberg:1990a}.
\end{remark}

Analogously to the pull-backs we shall now dualize the action of
differential operators to find an extension to distributional
sections. As we had (at least) two versions of dualizing differential
operators, we again obtain several possibilities for distributions.

We start with the ``intrinsic'' version.  Thus let $D: \Secinfty(E)
\longrightarrow \Secinfty(F)$ be a differential operator. Then its
adjoint is a differential operator
\begin{equation}
    \label{eq:transpose-of-diffop}
    D^\Trans: \Secinfty_0(F^* \tensor \Dichten T^*M)
    \longrightarrow \Secinfty_0(E^* \tensor \Dichten T^*M)
\end{equation}
of the same order as $D$. This motivates the following definition:
\begin{definition}[Differentiation of generalized sections]
    \label{definition:differentiation-of-gensec}
    \index{Generalized section!differentiation}%
    Let $D \in \Diffop^\bullet (E; F)$ then
    \begin{equation}
        \label{eq:diff-of-gensec}
        D: \Sec[-\infty](E) \longrightarrow \Sec[-\infty](F)
    \end{equation}
    is defined by
    \begin{equation}
        \label{eq:diff-of-gensec-formula}
        (Ds) (\mu) = s (D^\Trans \mu)
    \end{equation}
    for all $s \in \Sec[-\infty](E)$ and $\mu \in \Secinfty_0(F^*
    \tensor \Dichten T^*M)$.
\end{definition}
This definition indeed gives a reasonable notion of differentiation of
generalized sections as the following theorem shows:
\begin{theorem}
    \label{theorem:differentiation-of-gensec}
    Let $D \in \Diffop^k(E; F)$.
    \begin{theoremlist}
    \item \label{item:diffops-are-continuous} For all $s \in
        \Sec[-\infty](E)$ the definition
        \eqref{eq:diff-of-gensec-formula} gives a well-defined
        generalized section $Ds \in \Sec[-\infty](F)$ and the map
        \begin{equation}
            \label{eq:diff-of-gen-sec-is-cont}
            D: \Sec[-\infty](E) \longrightarrow \Sec[-\infty](F)
        \end{equation}
        is linear and weak$^*$ continuous. Moreover, we have for all
        $\ell \in \mathbb{N}_0$
        \begin{equation}
            \label{eq:DiffopOrderkDistributions}
            D: \Sec[-\ell](E) \longrightarrow \Sec[-\ell-k](F).
        \end{equation}
    \item \label{item:diffops-on-gensec-are-weakstar-extensions} The
        map $D$ is the unique extension of $D: \Secinfty(E)
        \longrightarrow \Secinfty(F)$ which is linear and weak$^*$
        continuous.
    \item \label{item:diffop-on-gensec-is-diffop} With respect to the
        $\Cinfty(M)$-module structure of $\Sec[-\infty](E)$ and
        $\Sec[-\infty](F)$, the map $D$ as in
        \eqref{eq:diff-of-gen-sec-is-cont} is a differential operator
        of order $k$ in the sense of the algebraic definition of
        differential operators, i.e.
        \begin{equation}
            \label{eq:diffop-of-gensec-is-alg-diffop}
            D \in \Diffop^k (\Sec[-\infty](E), \Sec[-\infty](F)).
        \end{equation}
    \item \label{item:supp-of-Ds} We have
        \begin{equation}
            \label{eq:supp-of-Ds}
            \supp (Ds) \subseteq \supp s
        \end{equation}
        and
        \begin{equation}
            \label{eq:singsupp-of-Ds}
            \singsupp (Ds) \subseteq \singsupp s.
        \end{equation}
    \item \label{item:diffop-for-gensec-is-local} For every open
        subset $U \subseteq M$ we have
        \begin{equation}
            \label{eq:diffop-for-gensec-is-local}
            Ds \at{U} = D\at{U} (s\at{U}).
        \end{equation}
    \end{theoremlist}
\end{theorem}
\begin{proof}
    Since $D^\Trans: \Secinfty_0(F^* \tensor \Dichten T^*M)
    \longrightarrow \Secinfty_0(E^* \tensor \Dichten T^*M)$ is again a
    differential operator of order $k$ by
    Proposition~\ref{proposition:existence-of-adjoints} and since
    differential operators are $\Cinfty_0$-continuous by
    Theorem~\ref{theorem:Diffop-is-cont-on-compact-supps}, the
    definition \eqref{eq:diff-of-gensec-formula} yields indeed a
    continuous linear functional $Ds \in \Secinfty_0(F^* \tensor
    \Dichten T^*M)' = \Sec[-\infty](F)$. Clearly, $D$ is linear and we
    have
    \[
    \seminorm[\mu](Ds) = |Ds(\mu)| = |s (D^\Trans \mu)|
    = \seminorm[D^\Trans \mu] (s),
    \]
    from which we obtain the $\textrm{weak}^*$-continuity at once. The
    claim \eqref{eq:DiffopOrderkDistributions} is clear by counting.
    The second part follows easily since by
    Theorem~\ref{theorem:weak-topology-on-gensecs},
    \refitem{item:sec-dense-in-gensec} the space $\Secinfty_0(E)
    \subseteq \Sec[-\infty](E)$ is weak$^*$ dense hence any weak$^*$
    continuous extensions is necessarily unique. For $s \in
    \Secinfty(E)$ the definition \eqref{eq:diff-of-gensec-formula}
    coincides with the usual application of $D$ by
    Proposition~\ref{proposition:existence-of-adjoints}: the
    definition \eqref{eq:diff-of-gensec-formula} was made precisely
    that way to have an \emph{extension} of $D: \Secinfty(E)
    \longrightarrow \Secinfty(F)$. For the third part, we first
    consider a differential operator $D \in \Diffop^0(E; F) =
    \Hom_{\Cinfty(M)}(\Secinfty(E), \Secinfty(F))= \Secinfty(\Hom(E,
    F))$ of order zero. For $s \in \Sec[-\infty](E)$ we have then for
    all $\mu \in \Secinfty_0(F^* \tensor \Dichten T^*M)$ the relation
    \[
    D(f \cdot s) (\mu) = (f \cdot s) (D^\Trans \mu)
    = s(f D^\Trans(\mu))
    = s(D^\Trans (f\mu)) = (f \cdot Ds)(\mu),
    \]
    hence $D(f \cdot s) = f \cdot D(s)$ follows. Thus $D$ as in
    \eqref{eq:diff-of-gensec-formula} is a $\Cinfty(M)$-linear map and
    hence a differential operator of order zero in the sense of
    definition \eqref{eq:Algebraic-Diffops}. Now we can proceed by
    induction on the order: assume that $D \in \Diffop^k(E; F)$ yields
    a differential operator $D \in \Diffop^k(\Sec[-\infty](E),
    \Sec[-\infty](F))$ of the same order $k$ for all $k \leq
    \ell$. Then for $D \in \Diffop^{\ell+1}(E,F)$ we have
    \begin{align*}
        \left(
            D(f \cdot s) - f \cdot D(s)
        \right) (\mu)
        = (f \cdot s) (D^\Trans \mu) -D(s) (f\mu)
        = s \left( f D^\Trans \mu - D^\Trans (f\mu) \right)
        = s ( [f,D^\Trans] \mu).
    \end{align*}
    Since for $A = [f,D] \in \Diffop^\ell(E,F)$ we have $A^\Trans =
    [f,D^\Trans]$, we see that $[f, D]: \Sec[-\infty](E)
    \longrightarrow \Sec[-\infty](F)$ is a differential operator of
    order $\ell$ by induction. Thus $D$ is again a differential
    operator of order $\ell +1$, since $f$ was arbitrary.  This shows
    the third part. Now let $\mu \in \Secinfty_0(F^* \tensor T^*M)$
    with $\supp \mu \subseteq M \backslash \supp s$ then $\supp
    D^\Trans \mu \subseteq M \backslash \supp s$ as well hence
    $(Ds)(\mu) = s(D^\Trans \mu) =0$ by
    Remark~\ref{remark:restriction-and-support-of-gensec},
    \refitem{item:supp-of-product}. Thus \eqref{eq:supp-of-Ds}
    follows. Let $t \in \Secinfty(E\at{M \backslash \singsupp s})$ be
    the smooth section such that for all $\mu \in \Secinfty_0(F^*
    \tensor \Dichten T^*M)$ with $\supp \mu \subseteq M \backslash
    \singsupp s $ we have $s(\mu) = \int_M t \: \mu$. Then for those
    $\mu$ we have
    \[
    (Ds)(\mu)
    = s(D^\Trans \mu)
    = \int_M t \; D^\Trans \mu
    = \int_M (Dt) \mu,
    \]
    since $\supp D^\Trans \mu \subseteq \supp \mu$. Thus $Ds$ is
    regular on $M \backslash \singsupp s$, too, hence for the singular
    support we get $\singsupp(Ds) \subseteq M \backslash (M \backslash
    \singsupp s) = \singsupp s$.  For the last part let $\mu \in
    \Secinfty_0(F^* \tensor \Dichten T^*M)$ be a test section with
    $\supp \mu \subseteq U$. Then
    \begin{align*}
        Ds\at{U} (\mu) = Ds(\mu) = s(D^\Trans \mu)
        = s(D^\Trans\at{U} \mu) = s\at{U} (D^\Trans\at{U} \mu)
        = \left(
            D\at{U} (s\at{U})
        \right)
        (\mu),
    \end{align*}
    since $D^\Trans \at{U}(\mu)$ has still support in $U$ by the
    locality of differential operators.
\end{proof}
\begin{remark}
    \label{remark:differentiation-with-fixed-density}
    In Theorem~\ref{theorem:agjoint-for-fixed-density} we have defined
    a different adjoint $D^\Trans \in \Diffop(E^*; F^*)$ of $D \in
    \Diffop(E; F)$ with respect to an a priori chosen positive density
    $\mu > 0$. We can use this adjoint to extend $D$ to distributional
    sections as well. To this end we first observe that every section
    in $\Secinfty_0(F^* \tensor \Dichten T^*M)$ is a tensor product
    $\omega \tensor \mu$ of a uniquely determined section $\omega \in
    \Secinfty_0(F^*)$ and the positive density $\mu$, since $\mu$
    provides a trivialization of $\Dichten T^*M$. Thus it is
    sufficient to consider $\omega \tensor \mu \in \Secinfty_0(F^*
    \tensor \Dichten T^*M)$ in the following. For $s \in
    \Sec[-\infty](E)$ we define $Ds: \Secinfty_0(F^* \tensor \Dichten
    T^*M) \longrightarrow \mathbb{C}$ by
    \begin{equation}
        \label{eq:differentiation-with-fixed-density}
        (Ds) (\omega \tensor \mu) = s( (D^\Trans \omega) \tensor \mu),
    \end{equation}
    which gives a well-defined linear map. Since $D^\Trans$ is
    continuous and since the tensor product is continuous too, $Ds \in
    \Sec[-\infty](F)$. Moreover,
    \begin{equation}
        \label{eq:cont-estimate-for-diff-with-density}
        \seminorm[\omega \tensor \mu](Ds)
        = |Ds(\omega \tensor \mu)|
        = \seminorm[D^\Trans \omega \tensor \mu] (s)
    \end{equation}
    shows that $D: \Sec[-\infty](E) \longrightarrow \Sec[-\infty](F)$
    is weak$^*$ continuous. Since by construction $D$ coincides with
    $D: \Secinfty(E) \longrightarrow \Secinfty(F)$ on the smooth
    sections $\Secinfty(E) \subseteq \Sec[-\infty](E)$, we conclude
    that the definition \eqref{eq:differentiation-with-fixed-density}
    and the intrinsic definition from
    Definition~\ref{definition:differentiation-of-gensec} actually
    \emph{coincide}. In particular, even though $D^\Trans$ in
    \eqref{eq:differentiation-with-fixed-density} depends on $\mu$
    explicitly, the combination $s (D^\Trans \omega \tensor \mu)$ only
    depends on the combination $\omega \tensor \mu$. In
    \cite[Sect.~1.1.2]{baer.ginoux.pfaeffle:2007a} the approach
    \eqref{eq:differentiation-with-fixed-density} was used to define
    the extension of $D$ to generalized sections.
\end{remark}

%
%

\subsection{Tensor Products}
\label{subsec:tensor-products}

In this section we consider various tensor product constructions for
distributions. The first one is about the values of a distribution and
provides a rather trivial extension of our previous considerations.
\begin{definition}[Vector-valued generalized sections]
    \label{definition:vector-valued-gensecs}
    \index{Generalized section!vector-valued}%
    Let $E \longrightarrow M$ be a vector bundle and $V$ a
    finite-dimensional vector space. Then a $V$-valued generalized
    section of $E$ is a continuous linear map
    \begin{equation}
        \label{eq:V-valued-gensec}
        s: \Secinfty_0(E^* \tensor \Dichten T^*M) \longrightarrow V.
    \end{equation}
    The set of all $V$-valued generalized sections of $E$ is denoted
    by $\Sec[-\infty](E; V)$.
\end{definition}
Since we always assume that the target vector space $V$ is
finite-dimensional, all Hausdorff locally convex topologies on $V$
coincide. Thus the notion of continuity of \eqref{eq:V-valued-gensec}
is non-ambiguous. It is clear that all the previous operations on
distributions can be carried over to the vector-valued case since they
were constructed from operations on the \emph{arguments} of $s$.
\begin{proposition}
    \label{proposition:splitting-of-vector-valued-gensec}
    For a finite-dimensional vector space $V$ and a vector bundle $E
    \longrightarrow M$ we have the canonical isomorphism
    \begin{equation}
        \label{eq:splitting-of-vector-valued-gensec}
        \Sec[-\infty](E) \tensor V \ni s \tensor v
        \; \mapsto \;
        (\omega \; \mapsto \; s(\omega) v) \in \Sec[-\infty](E; V).
    \end{equation}
\end{proposition}
\begin{proof}
    First we note that the map $\omega \mapsto s(\omega) v$ is linear
    and continuous with respect to the $\Cinfty_0$-topology of
    $\Secinfty_0(E^* \tensor \Dichten T^*M)$. Indeed, if $|s(\omega)|
    \leq c \seminorm[K,\ell] (\omega)$ for $K \subseteq M$ compact, $c
    > 0$ and $\ell \in \mathbb{N}_0$ and all $\omega \in
    \Secinfty_K(E^* \tensor \Dichten T^*M)$ then
    \[
    \norm{s(\omega) v} \leq c \seminorm[K,\ell] (\omega) \norm{v},
    \]
    where $\norm{\argument}$ is any norm on $V$. Thus the right hand side
    of \eqref{eq:splitting-of-vector-valued-gensec} is a vector-valued
    distribution. Clearly, the map is bilinear in $s$ and $v$ hence it
    indeed defines a linear map
    \[
    \Sec[-\infty](E) \tensor V \longrightarrow \Sec[-\infty](E; V).
    \]
    Let $e_1, \ldots, e_k \in V$ be a vector space basis. For a
    $V$-valued distribution $s \in \Sec[-\infty](E; V)$ we have scalar
    distributions $s^\alpha = e^\alpha \circ s$ since for
    finite-dimensional vector spaces the algebraic and topological
    duals coincide. Thus $s = s^\alpha e_\alpha$ in the sense that
    $s(\omega) = s^\alpha (\omega) e_\alpha$. Moreover, the $s^\alpha
    (\omega)$ are uniquely determined hence the $s^\alpha$ are unique.
    It follows that $s^\alpha \tensor e_\alpha$ is a pre-image of $s$
    under \eqref{eq:splitting-of-vector-valued-gensec}, hence
    \eqref{eq:splitting-of-vector-valued-gensec} is surjective.
    Injectivity is clear since the $s^\alpha$ are unique.
\end{proof}

In the following we shall use this isomorphism to identify
$\Sec[-\infty](E) \tensor V$ with $\Sec[-\infty](E; V)$. In
particular, the weak$^*$ topology of $\Sec[-\infty](E; V)$ is just the
component-wise weak$^*$ topology of $\Sec[-\infty](E)$. One can endow
$\Sec[-\infty](E) \tensor V$ with a tensor product topology such that
\eqref{eq:splitting-of-vector-valued-gensec} is even an isomorphism of
locally convex vector spaces. However, we shall not need this here.
Note also that for arbitrary locally convex $V$ the map
\eqref{eq:splitting-of-vector-valued-gensec} is still defined and
injective, but usually no longer surjective.

The next tensor product is based on the tensor product of the
arguments. We consider a product manifold $M \times N$ with the
canonical projections
\begin{equation}
    \label{eq:CanonicalProjections}
    M
    \stackrel{\pr_M}{\longleftarrow}
    M \times N
    \stackrel{\pr_N}{\longrightarrow}
    N.
\end{equation}
For this situation, we first prove the following statement which is of
independent interest:
\begin{theorem}
    \label{theorem:tensor-product-of-Ck-functions}
    \index{Tensor product!functions}%
    Let $M$, $N$ be manifolds. Then for all $k \in \mathbb{N}_0 \cup
    \{+\infty\}$ the map
    \begin{equation}
        \label{eq:tensor-product-of-Ck-functions}
        \Fun_0(M) \tensor \Fun_0(N) \ni f \tensor g
        \; \mapsto \;
        \pr_M^*f \pr_N^* g
        \in \Fun_0(M \times N)
    \end{equation}
    is a continuous injective algebra homomorphism with sequentially
    dense image with respect to the $\Fun_0$-topologies. In more
    detail, we have estimates
    \begin{equation}
        \label{eq:cont-estimate-for-tensor-product-of-Ck-functions}
        \seminorm[K \times L, k] (\pr_M^*f \pr_N^*g)
        \leq
        c \max_{\ell \leq k} \seminorm[K,\ell](f)
        \max_{\ell \leq k} \seminorm[L, \ell](g),
    \end{equation}
    if we use factorizing data to define the seminorms $\seminorm[K
      \times L, k]$ on $M \times N$.
\end{theorem}
\begin{proof}
    First we discuss the linear algebra aspects. Since the algebraic
    tensor product of two associative algebras is canonically an
    associative algebra, we can indeed speak of an algebra
    homomorphism.  It follows immediately that
    \eqref{eq:tensor-product-of-Ck-functions} is bilinear in $f$ and
    $g$ and thus well-defined on the tensor product. Then the
    homomorphism property is clear. The injectivity is clear as for
    linear independent $f_\alpha$ and linear independent $g_\beta$ the
    images of $f_\alpha \tensor g_\beta$ are still linear independent.
    This can be seen by evaluating at appropriate points $(x,y) \in M
    \times N$. Thus we can identify $f \tensor g$ with $\pr_M^*f
    \pr_N^*g$ and avoid the latter, more clumsy notation.  We come now
    to the continuity property. Thus let $\nabla^M$ and $\nabla^N$ be
    torsion-free covariant derivatives and let $\nabla^{M \times N}$
    be the corresponding covariant derivative on $M \times N$. By
    $\SymD_M$, $\SymD_N$, and $\SymD_{M \times N}$ we denote the
    corresponding symmetrized covariant derivatives. Now let $K
    \subseteq M$ and $L \subseteq N$ be compact. Then $K \times L
    \subseteq M \times N$ is compact, too, and every compact subset of
    $M \times N$ is contained in such a compactum for appropriate $K$
    and $L$. Thus it suffices to consider $K \times L \subseteq M
    \times N$. For $f \in \Fun_K(M)$ and $g \in \Fun_L(N)$ we compute
    \begin{align*}
        \SymD_{M \times N}^k (\pr_M^*(f) \pr_N^*(g))
        &=
        \sum_{\ell = 0}^k \binom{k}{\ell}
        \SymD_{M \times N}^\ell (\pr_M^* f) \vee
        \SymD_{M \times N}^{k-\ell} (\pr_N^* g) \\
        &=
        \sum_{\ell = 0}^k \binom{k}{\ell}
        \pr_M^*\left(\SymD_M^\ell f\right) \vee
        \pr_N^*\left(\SymD_N^{k-\ell} g\right),
    \end{align*}
    since $\SymD_{M \times N}$ is a derivation and since $\SymD_{M
      \times N} \pr_M^* = \pr_M^* \SymD_M$ as well as $\SymD_{M \times
      N} \pr_N^* = \pr_N^* \SymD_N$. If we also choose the Riemannian
    metric on $M \times N$ to be the product metric of $g_M$ on $M$
    and $g_N$ on $N$ we obtain for the $\seminorm[K \times L, k]$
    seminorm
    \begin{align*}
    \seminorm[K \times L, k](\pr_M^*(f) \pr_N^*(g))
    & = \sup_{(x,y) \in K \times L}
    \norm{\SymD_{M \times N}^k (\pr_M^*(f) \pr_N^*(g)) \At{(x,y)}}_{M
      \times N} \\
    & \leq \sup_{\substack{x \in K \\ 0 \leq \ell \leq k}} \sup_{y \in L}
    c \norm{\SymD_M^\ell f \at{x}}_M
    \norm{\SymD_N^{k-\ell} g \at{y}}_N \\
    & = c \max_{\ell \leq k} \seminorm[K,\ell](f) \seminorm[L,\ell](g),
    \end{align*}
    which shows the continuity property of
    \eqref{eq:tensor-product-of-Ck-functions}. We are left with the
    task to show that finite sums of factorizing functions are
    sequentially dense. Thus let $F \in \Fun_0(M \times N)$ be given.
    We choose atlases $\{(U_\alpha, x_\alpha)\}$ of $M$ and
    $\{(V_\beta, y_\beta)\}$ of $N$ together with subordinate
    partitions of unity $\{ \chi_\alpha \}$ and $\{ \psi_\beta \}$,
    respectively. Then the $\{ (U_\alpha \times V_\beta, x_\alpha
    \times y_\beta, \chi_\alpha \tensor \psi_\beta) \}$ provides an
    atlas of $M \times N$ with a corresponding partition of unity.
    Since $\supp F$ is compact, it follows that
    \[
    F = \sum_{\alpha, \beta} \chi_\alpha \tensor \psi_\beta \cdot F
    \]
    is a finite sum and each term $\chi_\alpha \tensor \psi_\beta
    \cdot F$ has compact support in $U_\alpha \tensor V_\beta$. Thus
    it will be sufficient to find a sequence in $\Fun_0(U_\alpha)
    \tensor \Fun_0(V_\beta)$ which approximates a function in
    $\Fun_0(U_\alpha \times V_\beta)$. This reduces the problem to the
    following \emph{local} problem: We have to show that
    $\Fun_0(\mathbb{R}^n) \tensor \Fun_0(\mathbb{R}^m)$ is
    sequentially dense in $\Fun_0(\mathbb{R}^{n+m})$. We will need the
    following technical lemma:
    \begin{lemma}
        \label{lemma:approx-by-polynomial}
        Let $f \in \Fun_0(\mathbb{R}^n)$ and $K \subseteq
        \mathbb{R}^n$ compact. For every $\epsilon > 0$ there exists a
        polynomial $p_\epsilon \in \Pol(\mathbb{R}^n)$ such that for
        all $\ell \le k$
        \begin{equation}
            \label{eq:approx-by-polynomial}
            \seminorm[K,\ell](p_\epsilon - f) < \epsilon.
        \end{equation}
    \end{lemma}
    \begin{subproof}
        We only sketch the proof which uses some convolution tricks.
        We consider the normalized Gaussian
        \[
        G_\delta(x) = \frac{1}{\sqrt{\pi \delta}^n}
        \E^{-\frac{1}{\delta} x²}
        \]
        for $\delta > 0$. Then the integral of $G_\delta$ equals one
        for all $\delta$. For $f \in \Fun_0(\mathbb{R}^n)$ the
        convolution\index{Convolution}
        \[
        (G_\delta * f)(x)
        = \int_{\mathbb{R}^n} G_\delta(x-y) f(y) \D^n y
        \]
        is a smooth function $G_\delta * f \in \Cinfty(\mathbb{R}^n)$
        and we have $\frac{\partial^{|I|}}{\partial x^I} (G_\delta *
        f) = G_\delta * \frac{\partial^{|I|} f}{\partial x^I}$ for all
        multiindexes $I$ with $|I| \leq k$. It is now a well-known
        fact that $G_\delta * f$ approximates $f$ uniformly on
        $\mathbb{R}^n$, i.e.
        \[
        \norm{G_\delta * f - f}_{\infty} \longrightarrow 0
        \quad \textrm{for} \quad \delta \longrightarrow 0
        \]
        in the sup-norm $\norm{\argument}_\infty$ over
        $\mathbb{R}^n$. If $k \geq 1$ we can repeat the argument and
        obtain that
        \[
        \sup_{x \in \mathbb{R}^n}
        \left|
            \frac{\partial^{|I|}}{\partial x^I}
            (G_\delta * f) -  \frac{\partial^{|I|} f}{\partial x^I}
        \right|
        \longrightarrow 0
        \]
        for $\delta \longrightarrow 0$ and for all multiindexes with
        $|I| \leq k$. In a second step we approximate the Gaussian by
        its Taylor series. Since on every compact subset $K \subseteq
        \mathbb{R}^n$ the Taylor series converges to $G_\delta$ in the
        $\Cinfty_K$-topology we find a polynomial
        \[
        p_{\epsilon, K, k, \delta} (x)
        = \frac{1}{\sqrt{\pi \delta}^n}
        \sum_{ r\geq 0}^{\textrm{finite}} \frac{1}{r!}
        \left(
            - \frac{x^2}{\delta}
        \right)^r
       \]
       such that for $\ell \le k$
       \[
       \seminorm[K,\ell] (G_\delta - p_{\epsilon, K, k, \delta})
       < \epsilon.
       \]
       The convolution
       \[
       f_{\epsilon,K,k,\delta} (x)
       = \int_{\mathbb{R}^n} p_{\epsilon,K,k,\delta}(x-y) f(y) \D^n y
       \]
       is again a polynomial of $x$ of the same order as
       $p_{\epsilon,K,k,\delta}$ and we use this to approximate $f$ on
       a compact subset. Thus let $K \subseteq \mathbb{R}^n$ be fixed
       and consider $x \in K$. Then
       \begin{align*}
           &\left|
               \frac{\partial^{|I|}}{\partial x^I}
               f_{\epsilon,B_R(0),k,\delta}(x)
               -
               \frac{\partial^{|I|}}{\partial x^I}
               (G_\delta *f)(x)
           \right| \\
           &\quad=
           \left|
               \int
               \left(
                   p_{\epsilon,B_R(0),k,\delta}(x-y)
                   \frac{\partial^{|I|} f}{\partial x^I}(y)
                   -
                   G_\delta(x-y)
                   \frac{\partial^{|I|} f}{\partial x^I}(y)
               \right) \D^n y
           \right| \\
           &\quad\leq
           \int_{\supp f}
           \left|
               p_{\epsilon,B_R(0),k,\delta}(x-y) - G_\delta(x-y)
           \right|
           \left|
               \frac{\partial^{|I|} f}{\partial x^I}(y)
           \right| \\
           &\quad\leq
           \epsilon \int_{\supp f}
           \left|
               \frac{\partial^{|I|} f}{\partial x^I}(y)
           \right| \D^n y,
       \end{align*}
       if we choose $B_R(0)$ large enough such that $K - \supp f
       \subseteq B_R(0)$. This is clearly possible since both $K$ and
       $\supp f$ are compact. It follows that on a compact subset $K
       \subseteq \mathbb{R}^n$ the polynomial
       $f_{\epsilon,B_R(0),k,\delta}$ approximates $G_\delta * f$ in
       the $\Fun_K$-topology. Thus we obtain that
       $f_{\epsilon,B_R(0),k,\delta}$ also approximates $f$ in the
       $\Fun_K$-topology as well. Rescaling $\epsilon$ appropriately
       gives the polynomials $p_\epsilon$ as desired.
   \end{subproof}

   Using this lemma we can proceed as follows: Let $F \in
   \Fun_0(\mathbb{R}^{n+m})$ be given and choose $\chi \in
   \Cinfty_0(\mathbb{R}^n)$ and $\psi \in \Cinfty_0(\mathbb{R}^m)$
   such that their tensor product $\chi \tensor \psi$ is equal to one
   on $\supp F$. This is clearly possible. Then $F = \chi \tensor \psi
   \cdot F$ can be approximated by polynomials on every compact
   subset. Let $K \times L \subseteq \mathbb{R}^{n+m}$ be a compactum
   such that $\chi \tensor \psi \in \Cinfty_{K \times
     L}(\mathbb{R}^{n+m})$ and choose $p_r \in \Pol(\mathbb{R}^{n+m})$
   such that on $K \times L$ the polynomials $p_r$ converge to $F$ in
   the $\Fun_{K \times L}$-topology by
   Lemma~\ref{lemma:approx-by-polynomial}. Since for polynomials we
   have
   \[
   \Pol(\mathbb{R}^{n+m}) =
   \Pol(\mathbb{R}^{n}) \tensor \Pol(\mathbb{R}^{m}),
   \]
   we find that $\chi \tensor \psi \cdot p_r \in \Cinfty_{K \times
     L}(\mathbb{R}^{n+m})$ is actually in $\Cinfty_K(\mathbb{R}^n)
   \tensor \Cinfty_L(\mathbb{R}^m)$. Now
   \[
   \seminorm[K \times L, k](F - \chi \tensor \psi \cdot p_r)
   = \seminorm[K \times L, k]
   (\chi \tensor \psi \cdot F - \chi \tensor \psi \cdot p_r)
   \leq c \seminorm[K \times L, k](F - p_r)
   \longrightarrow 0
   \]
   for $r \longrightarrow \infty$, which shows the density with
   respect to the $\Fun_0$-topology since all members of the sequence
   are in one fixed compactum $K \times L$.
\end{proof}
\begin{remark}
    \label{remark:tensor-product-of-functions}
    In fact, the proof even shows that
    \begin{equation}
        \label{eq:dense-inclusion-of-tensor-product}
        \Cinfty_0(M) \tensor \Cinfty_0(N) \subseteq \Fun_0(M \times N)
    \end{equation}
    is sequentially dense in the $\Fun_0$-topology for all $k \in
    \mathbb{N}_0 \cup \{+\infty\}$. Note that this gives an
    independent proof of
    Theorem~\ref{theorem:compact-and-smooth-sections-dense-in-others}
    at least for the scalar case as we can choose $N =
    \{\mathrm{pt}\}$ hence $\Cinfty_0(N) = \mathbb{C}$ and $\Fun_0(M
    \times N) \simeq \Fun_0(M)$. Thus we recover that
    \begin{equation}
        \label{eq:smooth-function-dense-in-Ck-functions}
        \Cinfty_0(M) \subseteq \Fun_0(M)
    \end{equation}
    is dense in the $\Fun_0$-topology.
\end{remark}
\begin{corollary}
    \label{corollary:tensor-product-without-compact-supp}
    For all $k \in \mathbb{N}_0 \cup \{+\infty\}$ the map
    \begin{equation}
        \label{eq:tensor-product-without-compact-supp}
        \Fun(M) \tensor \Fun(N) \ni f \tensor g
        \; \mapsto \;
        ((x,y) \; \mapsto \; f(x)g(y)) \in \Fun(M \times N)
    \end{equation}
    extends to a linear injective continuous algebra homomorphism with
    dense image with respect to the $\Fun$-topology.
\end{corollary}
\begin{proof}
    The estimates
    \eqref{eq:cont-estimate-for-tensor-product-of-Ck-functions} also
    show that \eqref{eq:tensor-product-without-compact-supp} is
    continuous. The fact that the image is dense follows from
    Theorem~\ref{theorem:tensor-product-of-Ck-functions} since it
    contains the images of $\Fun_0(M) \tensor \Fun_0(N)$ which is
    dense in $\Fun_0(M \times N)$ in the $\Fun_0$-topology. By
    Proposition~\ref{proposition:compact-sections-dense-in-smooth-sections}
    the subspace $\Fun_0(M \times N)$ is dense in $\Fun(M \times N)$
    in the $\Fun$-topology. Since $\Fun_0$-convergence implies
    $\Fun$-convergence, the statement follows. The remaining
    statements are clear.
\end{proof}

We can also extend the above statements to vector bundles. To this end
we recall the following construction of the \emph{external tensor
  product} of two vector bundles $E \longrightarrow M$ and $F
\longrightarrow N$. Over the Cartesian product $M \times N$ we
consider the vector bundle
\begin{equation}
    \label{eq:external-tensor-product}
    \index{Tensor product!external}%
    E \extensor F = \pr_M^\#(E) \tensor \pr_N^\#(F),
\end{equation}
where $\pr_M$ and $\pr_N$ are the usual projections and $\pr_M^\#(E)
\longrightarrow M \times N$ as well as $\pr_N^\#(F) \longrightarrow M
\times N$ denote the pull-backs of the vector bundles $E$ and $F$,
respectively. More informally, $E \extensor F$ is the vector bundle
with fiber $E_x \tensor F_y$ over $(x,y) \in M \times N$ and vector
bundle structure coming from \eqref{eq:external-tensor-product}. If
$e_\alpha \in \Secinfty(E\at{U})$ and $f_\beta \in \Secinfty(F\at{V})$
are local base sections then $\pr_M^\#(e_\alpha) \tensor
\pr_N^\#(f_\beta) \in \Secinfty\left(E \extensor F \at{U \times
      V}\right)$ are local base sections, too. To simplify our
notation we shall write
\begin{equation}
    \label{eq:external-tensor-product-of-sections}
    s \extensor t = \pr_M^\#(s) \tensor \pr_N^\#(t)
\end{equation}
for $s \in \Secinfty(E)$ and $t \in \Secinfty(F)$ in the sequel.
Without going into the details, the local trivializations of $E$ and
$F$ allow to use Theorem~\ref{theorem:tensor-product-of-Ck-functions}
and Corollary~\ref{corollary:tensor-product-without-compact-supp} to
obtain the following analogue for vector bundles:
\begin{theorem}
    \label{theorem:sections-on-external-tensor-product}
    \index{Tensor product!external!sections}%
    Let $k \in \mathbb{N}_0 \cup \{+\infty\}$ and let $E
    \longrightarrow M$ and $F \longrightarrow N$ be vector
    bundles. Then
    \begin{equation}
        \label{eq:sections-on-external-tensor-product-compact-supp}
        \Sec_0(E) \tensor \Sec_0(F) \ni s \tensor t \; \mapsto \;
        s \extensor t \in \Sec_0(E \extensor F)
    \end{equation}
    is an injective continuous $\Fun_0(M) \tensor \Fun_0(N)$-module
    morphism with sequentially dense image in the $\Fun_0$-topology.
    Analogously,
    \begin{equation}
        \label{eq:sections-on-external-tensor-product}
         \Sec(E) \tensor \Sec(F) \ni s \tensor t \; \mapsto \;
         s \extensor t \in \Sec(E \extensor F)
    \end{equation}
    is an injective continuous $\Fun(M) \tensor \Fun(N)$-module
    morphism with dense image in the $\Fun$-topology.
\end{theorem}

Note that on the left hand side the tensor product is taken over
$\mathbb{R}$ or $\mathbb{C}$, depending on the type of the vector
bundles. The module structures on both sides are the canonical ones.
\begin{remark}
    \label{remark:support-of-external-tensor-product}
    It should be noted that for $s \in \Secinfty(E)$ and $t \in
    \Secinfty(F)$ we have
    \begin{equation}
        \label{eq:support-of-external-tensor-product}
        \supp(s \extensor t) = \supp s \times \supp t.
    \end{equation}
\end{remark}
\begin{remark}
    \label{remark:external-tensor-product-of-densities}
    \index{Density!external tensor product}%
    \index{Tensor product!external!densities}%
    For the density bundles we have canonically
    \begin{equation}
        \label{eq:densities-on-product-factorize}
        \Dichten T^*M \extensor \Dichten T^*N
        \cong \Dichten T^*(M \times N),
    \end{equation}
    where the isomorphism is defined by
    \begin{equation}
        \label{eq:map-for-densitiy-factorizing}
        \Dichten T^*_x M \tensor \Dichten T^*_y N
        \ni \mu_x \tensor \nu_y
        \; \mapsto \;
        \mu_x \extensor \nu_y \in \Dichten T^*_{(x,y)} (M \times N),
    \end{equation}
    with
    \begin{equation}
        \label{eq:how-densities-on-product-factorize}
        (\mu_x \extensor \nu_y) (v_1, \ldots, v_m, w_1, \ldots, w_n)
        = \mu_x (v_1, \ldots, v_m) \nu_y (w_1, \ldots, w_n),
    \end{equation}
    where $v_1, \ldots, v_m \in T_x M$ and $w_1, \ldots, w_n \in T_y
    N$.  Moreover, for $E_i \longrightarrow M$ and $F_i
    \longrightarrow N$ with $i = 1, 2$ we have the compatibility
    \begin{equation}
        \label{eq:extensor-is-component-wise}
        (E_1 \tensor E_2) \extensor (F_1 \tensor F_2)
        \simeq (E_1 \extensor F_1) \tensor (E_2 \extensor F_2),
    \end{equation}
    and in particular
    \begin{equation}
        \label{eq:extensor-of-vecvalued-densities}
        (E^* \tensor \Dichten T^*M) \extensor (F^* \tensor \Dichten
        T^*N)
        \simeq (E^* \extensor F^*) \tensor \Dichten T^*(M \times N),
    \end{equation}
    which we shall frequently use in the following.
\end{remark}
In order to define the tensor product of distributions we need the
following technical lemma:
\begin{lemma}
    \label{lemma:parameter-differentiation-under-integral}
    Let $X \subseteq \mathbb{R}^n$ and $Y \subseteq \mathbb{R}^m$ be
    open and let $\phi \in \Cinfty(X \times Y)$ be smooth. Assume that
    there is a compact subset $K \subseteq X$ such that $\supp \phi
    \subseteq K \times Y$. Let $u \in \Cinfty_0(X)'$ be a scalar
    distribution. Then the function
    \begin{equation}
        \label{eq:integral-with-parameter}
        y \; \mapsto \; u(\phi(\argument, y))
    \end{equation}
    is smooth on $Y$. Moreover, for all multiindexes $I \in
    \mathbb{N}_0^m$ we have
    \begin{equation}
        \label{eq:differentiation-rule-for-parameter-under-integral}
        \frac{\partial^{|I|}}{\partial y^I} u(\phi(\argument, y))
        =
        u \left(
            \frac{\partial^{|I|}}{\partial y^I} \phi(\argument, y)
        \right).
    \end{equation}
    Finally, for $f \in \Cinfty(Y)$ we have
    \begin{equation}
        \label{eq:function-linearity-of-integral-in-parameter}
        u( f(y) \phi(\argument, y) ) = f(y) u(\phi(\argument, y)),
    \end{equation}
    i.e. the map $\phi \mapsto (y \mapsto u(\phi(\argument,y)))$ is
    $\Cinfty(Y)$-linear.
\end{lemma}
\begin{proof}
    In \eqref{eq:integral-with-parameter} we apply $u$ to the function
    $x \mapsto \phi(x,y)$ for fixed $y \in Y$. By assumption, this
    function has compact support in $K \subseteq X$ with respect to
    the $x$-variables for every fixed $y \in Y$, hence
    \eqref{eq:integral-with-parameter} is a well-defined function.  We
    shall now consider a slightly more detailed statement. On the
    compact subset $K$ the distribution $u$ has some finite order
    $\ell = \ord_K (u)$. Thus for all test functions $\varphi \in
    \Cinfty_K(X)$
    \[
    |u(\varphi)| \leq c \seminorm[K,\ell](\varphi),
    \tag{$*$}
    \]
    and we can extend $u$ to a continuous linear functional on
    $\Fun[\ell]_K(X)$ such that ($*$) still holds for $\varphi \in
    \Fun[\ell]_K(X)$. We refine the claim as follows: for $\phi \in
    \Fun(X \times Y)$ with $\supp \phi \subseteq K \times Y$ and $k
    \geq \ell$ the function $y \mapsto u(\phi(\argument,y))$ is in
    $\Fun[k-\ell](Y)$ and
    \eqref{eq:differentiation-rule-for-parameter-under-integral} holds
    for all $|I| \leq k-\ell$. Clearly, this statement includes
    \eqref{eq:differentiation-rule-for-parameter-under-integral} and
    \eqref{eq:integral-with-parameter} for the smooth case $k =
    \infty$.  Let $y_0 \in Y$ be fixed and consider some $B_r(y_0)^\cl
    \subseteq Y$. Then on the compact subset $K \times B_r(y_0)^\cl
    \subseteq X \times Y$ the function $\frac{\partial^{|J|}
      \phi}{\partial x^J}$ is uniformly continuous as long as $|J|
    \leq k$. Thus for $\epsilon > 0$ there is a $\delta > 0$ such that
    for $(x, y_0)$, $(x, y_0 + h) \in K \times B_r(y_0)^\cl$ with $|h|
    \leq \delta$ we have
    \[
    \left|
        \frac{\partial^{|J|} \phi}{\partial x^J} (x,y_0)
        - \frac{\partial^{|J|} \phi}{\partial x^J} (x, y_0 + h)
    \right|
    < \epsilon.
    \]
    It follows that
    \[
    \seminorm[K, \ell]( \phi(\argument, y_0) - \phi(\argument, y_0 + h))
    < \epsilon
    \]
    for those $h$ since $\ell \leq k$. Thus the continuity of $u$ with
    respect to the norm $\seminorm[K, \ell]$ on $\Fun[\ell]_K(X)$ as
    in ($*$) yields
    \[
    u(\phi(\argument, y_0 + h)) \longrightarrow u(\phi(\argument, y_0))
    \]
    for $h \longrightarrow 0$ for all $y_0 \in Y$. Thus
    \eqref{eq:integral-with-parameter} is continuous, This proves the
    case $k = \ell$. Now assume $k \geq \ell + 1$ hence we have some
    orders of differentiation for ``free''. Thus let $e \in
    \mathbb{R}^n$ be a unit vector and $y_0 \in Y$ together with a
    sufficiently small ball $B_r(y_0)^\cl \subseteq Y$ as before. Then
    for $|J| \leq \ell$ the partial derivatives $\frac{\partial^{|J|}
      \phi}{\partial x^J}$ are at least once continuously
    differentiable. Hence for $0 < |t| < r$
    \[
    \frac{1}{t}
    \left(
        \frac{\partial^{|J|} \phi}{\partial x^J} (x, y_0 + te)
        - \frac{\partial^{|J|} \phi}{\partial x^J} (x, y_0)
    \right)
    =
    e^i \frac{\partial^{|J|+1} \partial \phi}{\partial x^J \partial y^i}
    (x, y_0 + t_0 e)
    \]
    with some appropriate $t_0 \in [0,t]$. Since the $(|J| + 1)$-st
    derivatives are still continuous, on $K \times B_r(y_0)^\cl$ they
    are uniformly continuous. Thus for all $0 < |t| \leq \delta$
    \[
    \left|
        \frac{1}{t}
        \left(
            \frac{\partial^{|J|} \phi}{\partial x^J}(x, y_0 + te)
            -  \frac{\partial^{|J|} \phi}{\partial x^J}(x, y_0)
            -  \frac{\partial^{|J|+1} \phi}{\partial x^J \partial y^i}
            (x, y_0) e^i
        \right)
    \right|
    < \epsilon
    \]
    with some appropriately chosen $\delta > 0$. This means that
    \[
    \seminorm[K,\ell]
    \left(
        \frac{1}{t}
        \left(
            \phi(\argument, y_0 + te) - \phi(\argument, y_0)
        \right)
        - \partial_e \phi(\argument, y_0)
    \right)
    < \epsilon.
    \]
    Hence again by the continuity of $u$ we get for the directional
    derivative in direction $e$
    \[
    \partial_e u(\phi(\argument, y_0)) = u( \partial_e \phi(\argument, y_0)).
    \]
    Since $y_0$ was arbitrary and since $\partial_e \phi$ is
    $\Fun[k-1]$ with $k-1 \geq \ell$ we see that all directional
    derivatives at all $y \in Y$ exist and are continuous. This proves
    that \eqref{eq:integral-with-parameter} is in $\Fun[1](Y)$ and
    \eqref{eq:differentiation-rule-for-parameter-under-integral} is
    valid for $|I| = 1$. By induction we can proceed as long as $k
    \geq \ell$. The last statement is clear since $u$ acts only on the
    $x$-variables and not on the $y$-variables.
\end{proof}

In a geometric context the above lemma, in its refined version,
becomes the following statement:
\begin{proposition}
    \label{proposition:geometric-parameter-differentiation}
    Let $E \longrightarrow M$ and $F \longrightarrow N$ be vector
    bundles and let $\mu \in \Sec((E^* \extensor F^*) \tensor \Dichten
    T^*(M \times N))$ be a density such that there exists a compact
    subset $K \subseteq M$ with $\supp\mu \subseteq K \times N$. Let
    $s \in \Sec[-\infty](E)$ be a generalized section such that
    $\ord_K(s) \leq \ell$. Then the map
    \begin{equation}
        \label{eq:reduced-section}
        (s \tensor \id)(\mu): y \; \mapsto \; s(\mu(\argument, y))
    \end{equation}
    defines a $\Fun[k-\ell]$-section of $F^* \tensor \Dichten T^*N$.
    If $F' \longrightarrow N$ is another vector bundle and $D \in
    \Diffop^m (F \tensor \Dichten T^*N; F' \tensor \Dichten T^*N)$ a
    differential operator of order $m \leq k-\ell$ then $D$ applied to
    \eqref{eq:reduced-section} coincides with the section
    \begin{equation}
        \label{eq:geometric-parameter-differentiation}
        y \; \mapsto \; s( (\id \extensor D)(\mu)(\argument, y)),
    \end{equation}
    where $\id \extensor D$ means that $D$ acts only on the
    $y$-variables. For the support of $(s \tensor \id) (\mu)$ we have
    \begin{equation}
        \label{eq:support-of-reduced-section}
        \supp (s \tensor \id)(\mu) \subseteq \pr_N (\supp \mu).
    \end{equation}
\end{proposition}
\begin{proof}
    By the usual partition of unity argument with the usual local
    trivialization of the involved bundles we can reduce the above
    statements to the local and scalar case. Thus
    Lemma~\ref{lemma:parameter-differentiation-under-integral} yields
    that \eqref{eq:reduced-section} is a well-defined
    $\Fun[k-\ell]$-section and the combination of
    \eqref{eq:differentiation-rule-for-parameter-under-integral} and
    \eqref{eq:function-linearity-of-integral-in-parameter} gives
    \eqref{eq:geometric-parameter-differentiation}. It remains to show
    \eqref{eq:support-of-reduced-section}. Thus let $y \in N
    \backslash \pr_N (\supp \mu)$. Thus for all $x \in M$ we have
    $\mu(x, y) = 0$. This gives immediately $s(\mu(\argument, y)) =
    0$.  Since $N \backslash \pr_N(\supp \mu)$ is open,
    \eqref{eq:support-of-reduced-section} follows.
\end{proof}

\begin{remark}
    \label{remark:smooth-reduced-section}
    In particular, for all $s\in \Sec[-\infty](E)$ and $\mu \in
    \Secinfty_0((E^* \extensor F^*) \tensor \Dichten T^*(M \times N))$
    we have $(s \tensor \id) (\mu) \in \Secinfty_0(F^* \tensor
    \Dichten T^*N)$.
\end{remark}
We use this proposition now to prove the following statement on the
(external) tensor product of distributions.
\begin{theorem}[Tensor product of generalized sections]
    \label{theorem:tensor-product-of-generalized-sections}
    \index{Generalized section!external tensor product}%
    Let $E \longrightarrow M$ and $F \longrightarrow N$ be vector
    bundles and let $s \in \Sec[-\infty](E)$ and $t \in
    \Sec[-\infty](F)$ be generalized sections. Then there exists a
    unique generalized section $s \extensor t \in \Sec[-\infty](E
    \extensor F)$ such that
    \begin{equation}
        \label{eq:defining-property-for-gensec-tensor-product-on-fact-sections}
        (s \extensor t) (\mu \extensor \nu) = s(\mu) t(\nu)
    \end{equation}
    for $\mu \in \Secinfty_0(E^* \tensor \Dichten T^*M)$ and $\nu \in
    \Secinfty_0(F^* \tensor \Dichten T^*N)$. Moreover, for $\omega \in
    \Secinfty_0((E^* \extensor F^*) \tensor \Dichten T^*(M \times N))$
    we have
    \begin{equation}
        \label{eq:fubini-for-tensorproduct}
        (s \extensor t) (\omega)
        = t( (s \tensor \id)(\omega))
        = s( (\id \tensor t) (\omega)).
    \end{equation}
\end{theorem}
\begin{proof}
    Since $\Secinfty_0(E^* \tensor \Dichten T^*M) \tensor
    \Secinfty_0(F^* \tensor \Dichten T^*N)$ is dense in $\Secinfty_0(
    (E^* \extensor F^*) \tensor \Dichten T^*(M \times N))$ by
    Theorem~\ref{theorem:sections-on-external-tensor-product} and the
    identification \eqref{eq:extensor-of-vecvalued-densities} of
    Remark~\ref{remark:external-tensor-product-of-densities}, the
    uniqueness of $s \extensor t$ with the property
    \eqref{eq:defining-property-for-gensec-tensor-product-on-fact-sections}
    is clear. The idea is now to use the feature
    \eqref{eq:fubini-for-tensorproduct} to actually construct $s
    \extensor t$: Thus let $\omega \in \Secinfty_0(E^* \extensor F^*
    \tensor \Dichten T^*(M \times N))$ be given. We can assume that
    $\supp \omega \subseteq K \times L$ with compact subsets $K
    \subseteq M$ and $L \subseteq N$, respectively. For $s$ and $t$ we
    have estimates of the form
    \[
    |s(\mu)| \leq c \seminorm[K,k](\mu)
    \tag{$*$}
    \]
    \[
    |t(\nu)| \leq c' \seminorm[L, \ell](\nu),
    \tag{$**$}
    \]
    for the seminorms of Remark~\ref{remark:SymbolicSeminorms}
    whenever $\supp \mu \subseteq K$ and $\supp \nu \subseteq L$. By
    Proposition~\ref{proposition:geometric-parameter-differentiation}
    we know that
    \[
    (s \tensor \id) (\omega): y \; \mapsto \; s(\omega(\argument, y))
    \]
    is a smooth section of $F^* \tensor \Dichten T^*N$. Moreover,
    since the application of $s$ is $\Cinfty(N)$-linear and commutes
    with differentiation in $N$-direction we immediately conclude that
    \[
    \seminorm[L, \ell] ((s \tensor \id) (\omega))
    \leq c'' \seminorm[K \times L, \ell](\omega).
    \tag{$*$$**$}
    \]
    Finally, by Remark~\ref{remark:smooth-reduced-section} we have
    \[
    \supp ((s \tensor \id) (\omega)) \subseteq \pr_N (\supp \omega)
    \subseteq L,
    \]
    hence $(s \tensor \id) (\omega)$ has compact support. Thus we can
    apply $t$ and obtain by ($**$)
    \[
    | t( (s \tensor \id)(\omega)) |
    \leq c' \seminorm[L,\ell]( (s \tensor \id) (\omega))
    \leq c' c'' \seminorm[K \times L,\ell](\omega).
    \]
    Thus $\omega \mapsto t((s \tensor \id)(\omega) )$ is a continuous
    linear functional on $\Secinfty_{K \times L}( (E^* \extensor F^*)
    \tensor \Dichten T^*(M \times N))$ for all $K \times L$ with
    respect to the $\Cinfty_{K \times L}$-topology. Hence it defines a
    generalized section in $\Sec[-\infty](E \extensor F)$ by the
    characterization of
    Theorem~\ref{theorem:inductive-limit-topology},
    \refitem{item:LFcontinuity}.  If $\omega = \mu \extensor \nu$ is
    an external tensor product itself, we obtain
    \begin{align*}
        t( (s \tensor \id) (\mu \extensor \nu))
        & = t( y \; \mapsto \; s( (\mu \extensor \nu)(\argument, y)) ) \\
        & = t( y \; \mapsto \; s( \mu(\argument) \nu(y) ) ) \\
        & = t( y \; \mapsto \; \nu(y) s(\mu) ) \\
        & = t(\nu) s(\mu).
    \end{align*}
    This shows that the distribution $t \circ (s \tensor \id)$
    satisfies
    \eqref{eq:defining-property-for-gensec-tensor-product-on-fact-sections}.
    Hence it is the \emph{unique} solution $s \extensor t$ we are
    looking for. This proves existence of $s \extensor t$ and the
    first half of \eqref{eq:fubini-for-tensorproduct}. However, we
    could have constructed $s \extensor t$ by taking $s \circ (\id
    \tensor t)$ as well which gives, by uniqueness, the same $s
    \extensor t$. Thereby we have \eqref{eq:fubini-for-tensorproduct}.
\end{proof}
\begin{remark}
    \label{remark:external-tensor-product-of-generalized-sections}
    For the external tensor product
    \begin{equation}
        \label{eq:external-tensor-product-of-generalized-sections}
        \extensor: \Sec[-\infty](E) \tensor \Sec[-\infty](F)
        \longrightarrow \Sec[-\infty](E \extensor F)
    \end{equation}
    one immediately obtains
    \begin{equation}
        \label{eq:support-of-extensor-product-of-gensecs}
        \supp (s \extensor t) = \supp s \times \supp t,
    \end{equation}
    whence we also have
    \begin{equation}
        \label{eq:extensor-of-compact-supp-gensecs}
        \extensor: \Sec[-\infty]_0(E) \tensor \Sec[-\infty]_0(F)
        \longrightarrow \Sec[-\infty]_0(E \extensor F).
    \end{equation}
    It can be shown that for compactly supported $s$ and $t$ the
    conclusions of
    Theorem~\ref{theorem:tensor-product-of-generalized-sections}
    remain valid for $\mu$, $\nu$, $\omega$ not necessarily compactly
    supported.
\end{remark}
\begin{remark}[``Internal'' tensor product of distributions]
    \label{remark:internal-tensor-product-of-distributions}
    \index{Generalized section!internal tensor product}%
    For vector bundles $E \longrightarrow M$ and $F \longrightarrow M$
    over the same manifold, one may wonder whether there is an
    ``internal'' tensor product of generalized sections, i.e. a map
    \begin{equation}
        \label{eq:internal-tensor-product-of-distributions}
        \tensor: \Sec[-\infty](E) \tensor \Sec[-\infty](F)
        \longrightarrow \Sec[-\infty](E \tensor F),
    \end{equation}
    extending the tensor product of smooth sections, which is now
    $\Cinfty(M)$-bilinear with respect to the $\Cinfty(M)$-module
    structures of generalized sections. If such an extension of the
    usual tensor product of smooth section would exist in general,
    this would result in an algebra structure on $\Cinfty_0(M)'$ if we
    take $E = F$ to be the trivial line bundles. Here on meets serious
    problems: such a multiplication (obeying the usual properties) can
    be shown to be impossible. A ``definition'' of $s \tensor t$ like
    \begin{equation}
        \label{eq:wannabe-internal-tensor-product-definition}
        \textrm{``}
        (s \tensor t) (\mu \tensor \nu) = s(\mu) t(\nu)
        \textrm{''}
    \end{equation}
    is \emph{not} well-defined since the tensor product $\mu \tensor
    \nu$ of sections is $\Cinfty(M)$-bilinear while the right hand
    side of \eqref{eq:wannabe-internal-tensor-product-definition} is
    certainly not $\Cinfty(M)$-bilinear.

    Note however, that under certain circumstances the tensor product
    $s \tensor t$ can indeed be defined in a reasonable way. However,
    a much more sophisticated analysis of the singularities of $s$ and
    $t$ is needed.
\end{remark}


%% file: chap2.tex
%
%

\chapter{Elements of Lorentz Geometry and Causality}
\label{cha:LorentzGeometry}

In this second chapter we set the stage for wave equations on
spacetime manifolds. First we recall some basic properties and notions
for manifolds with covariant derivative, positive densities and
semi-Riemannian metrics. We shall discuss their relations and
introduce concepts like parallel transport as well as certain
canonical differential operators arising from the choice of a
semi-Riemannian metric. In particular, the d'Alembert operator will
provide the prototype of a wave operator. We generalize this to
arbitrary vector bundles and discuss several physical examples of wave
equations resulting from these differential operators.

After discussing the basics of semi-Riemannian and Lorentz metrics we
introduce the notions of causality on Lorentz manifolds. To this end
we first have to endow the Lorentz manifold with a time orientation
which then gives rise to the notions of future and past. The most
important notion in this context for us will be that of Cauchy
hypersurfaces. On one hand, the existence of a Cauchy hypersurface
will yield a particularly nice causal structure of the Lorentz
manifold. On the other hand, they will serve as the natural starting
point where we can pose initial value problems for a wave equation.

Such initial value problems for wave equations will then be the
subject of the last part of this chapter. Closely related will be the
notion of Green functions of advanced and retarded type. They are
particular elementary solutions of the wave equations subject to
``boundary conditions'' referring to the causal structure of the
spacetime.

For several theorems we will not provide proofs in this chapter as
this would lead us too far into the realm of Lorentz geometry. Instead
we refer to the literature, in particular to the textbooks
\cite{friedlander:1975a, oneill:1983a, beem.ehrlich.easley:1996a} as
well as to the review article \cite{minguzzi.sanchez:2006a:pre}.

%
%

\section{Preliminaries on Semi-Riemannian Manifolds}
\label{sec:SemiRiemannianManifold}

\input{semiriemann}

%
%

\section{Causal Structure on Lorentz Manifolds}
\label{sec:CausalStructure}

\input{causal}

%
%

\section{The Cauchy Problem and Green Functions}
\label{sec:CauchyProblemGreenFunctions}

\input{cauchy}


%% file: semiriemann.tex
%
%

In this section we collect some further properties of covariant
derivatives on vector bundles and their curvature, specializing to the
Levi-Civita connection of a semi-Riemannian metric. All of the
material is very much standard and can be found in textbooks like
\cite{oneill:1983a, beem.ehrlich.easley:1996a, lang:1999a}.

%
%

\subsection{Parallel Transport and Curvature}
\label{subsec:paralltransport-curvature}

Let $\nabla^E$ be a covariant derivative for a vector bundle $E
\longrightarrow M$ as before. Recall that the curvature tensor $R$ of
$\nabla^E$ is defined by
\begin{equation}
    \label{eq:curvature-tensor}
    \index{Covariant derivative!curvature}%
    R(X,Y)s
    = \nabla_X \nabla_Y s - \nabla_Y \nabla_X s - \nabla_{[X,Y]} s
\end{equation}
for $X, Y \in \Secinfty(TM)$ and $s \in \Secinfty(E)$. A simple
computation shows that $R$ is $\Cinfty(M)$-linear in each argument and
thus defines a tensor field
\begin{equation}
    \label{eq:curvature-tensor-field}
    R \in \Secinfty(\End(E) \tensor \Anti^2 T^*M).
\end{equation}
There are certain contractions we can build out of $R$. The most
important one is the pointwise trace of the $\End(E)$-part of
$R$. This gives a two-from
\begin{equation}
    \label{eq:trace-of-curvature}
    \tr R (X,Y) = \tr (s \mapsto R(X,Y)s),
\end{equation}
i.e. a section $\tr R \in \Secinfty(\Anti^2 T^*M)$. The following
lemma gives an interpretation of $\tr R$:
\begin{lemma}
    \label{lemma:trace-of-curvature}
    Let $\nabla^E$ be a covariant derivative for a vector bundle $E
    \longrightarrow M$.
    \begin{lemmalist}
    \item \label{item:trR-is-closed} The two-form $\tr T \in
        \Secinfty(\Anti^2 T^*M)$ is closed, $\D \tr R = 0$.
    \item \label{item:trR-is-exact} The two-form $\tr R$ is exact. In
        fact,
        \begin{equation}
            \label{eq:trR-is-exact}
            \tr R = -\D \alpha,
        \end{equation}
        where $\alpha \in \Secinfty(T^*M)$ is defined by
        \begin{equation}
            \label{eq:potential-for-trR}
            \alpha(X) = \frac{\nabla^E_X \mu}{\mu},
        \end{equation}
        with respect to any chosen positive density $\mu \in
        \Secinfty(\Dichten E^*)$.
    \end{lemmalist}
\end{lemma}
\begin{proof}
    Clearly, we only have to show \refitem{item:trR-is-exact}. Note
    that \refitem{item:trR-is-closed} would also follow rather easily
    from the Bianchi identity. Let $\mu \in \Secinfty(\Dichten E^*)$
    be a positive density. Then the covariant derivative $\nabla^E$ is
    extended as usual to $\Dichten E^*$ and $\alpha$ is a well-defined
    one-form. A simple computation shows that the curvature of
    $\nabla^{\Dichten E^*}$ is given by $\D \alpha$. On the other
    hand, the curvature of $\nabla^{\Dichten E^*}$ is given by $-\tr
    R$, see e.g.  \cite[Prop.~2.2.43]{waldmann:2007a}.
\end{proof}

With other words, $\tr R = 0$ is a necessary condition for the
existence of a covariantly constant density $\mu \in \Secinfty(\Dichten
E^*)$. In fact, the condition is locally also sufficient and globally
the deRham class $[\alpha] \in \HdR^1(M)$ might be an obstruction.
\begin{definition}[Unimodular covariant derivative]
    \label{definition:unimodular-covariant-derivative}
    \index{Covariant derivative!unimodular}%
    A covariant derivative $\nabla^E$ is called unimodular if $\tr R^E
    = 0$.
\end{definition}

Let $\gamma: I \subseteq \mathbb{R} \longrightarrow M$ be a smooth
curve defined on an open interval $I$ and let $a, b \in I$. In
general, the fibers of $E$ at $\gamma(a)$ and $\gamma(b)$ are not
related in a canonical way. Using the covariant derivative, this can
be done as follows. We are looking at a section $s$ along $\gamma$
such that $s$ is covariantly constant in the direction $\dot{\gamma}$.
More precisely, we consider the pull-back bundle $\gamma^\#E
\longrightarrow I$ together with the pull-back $\nabla^\#$ of
$\nabla^E$. Then we want to find a section $s \in
\Secinfty(\gamma^\#E)$ with
\begin{equation}
    \label{eq:covariant-constant-section}
    \nabla^\#_{\frac{\partial}{\partial t}} s = 0.
\end{equation}
If $\{e_\alpha\}$ are local base sections of $E$ over some open subset
$U \subseteq M$ and $\gamma(I) \subseteq U$ then
\eqref{eq:covariant-constant-section} is equivalent to
\begin{equation}
    \label{eq:ComputenablasLocally}
    0 = \nabla^\#_{\dot{\gamma}} (s^\alpha (t) e_\alpha (\gamma(t)))
    = \dot{s}^\alpha(t) e_\alpha (\gamma(t))
    + s^\alpha(t) A^\beta_\alpha (\dot{\gamma}(t)) e_\beta(\gamma(t)),
\end{equation}
i.e.
\begin{equation}
    \label{eq:covariant-constant-local}
    \dot{s}^\beta(t) + A^\beta_\alpha (\dot{\gamma}(t)) s^\alpha(t) = 0.
\end{equation}
Since \eqref{eq:covariant-constant-local} is an ordinary linear
differential equation for the coefficient functions $s^\alpha: I
\longrightarrow \mathbb{R}$, they have unique solutions $s^\alpha(t)$
for all $t$ and all initial conditions $s^\alpha(a)$. Moreover, the
resulting time evolution $s^\alpha(a) \mapsto s^\alpha(t)$ is a
\emph{linear} map and by uniqueness even an isomorphism. If the image
of $\gamma$ is not within the domain of a single bundle chart we can
cover it with several ones (finitely many for compact time intervals)
and use the uniqueness statement to glue the local solutions together
in the usual way. The uniqueness will then guarantee that the result
will not depend on the choice how we covered the curve with bundle
charts. Finally, this gives the following result:
\begin{proposition}
    \label{proposition:parallel-transport}
    Let $\nabla^E$ be a covariant derivative for $E \longrightarrow M$
    and let $\gamma: I \subseteq \mathbb{R} \longrightarrow M$ be a
    smooth curve. Let $a, b \in I$.
    \begin{propositionlist}
    \item \label{item:uniquenss-of-parallel-transport} For every
        initial condition $s_{\gamma(a)} \in E_{\gamma(a)}$ there
        exists a unique solution $s(t) \in E_{\gamma(t)}$ of
        \eqref{eq:covariant-constant-section}.
    \item \label{item:parallel-transport-is-isomorphism} The map
        $s_{\gamma(a)} \mapsto s(b)$ is a linear isomorphism
        $E_{\gamma(a)} \longrightarrow E_{\gamma(b)}$ which is denoted
        by
        \begin{equation}
            \label{eq:parallel-transport}
            P_{\gamma, a \rightarrow b}:
            E_{\gamma(a)} \longrightarrow E_{\gamma(b)}.
       \end{equation}
    \end{propositionlist}
\end{proposition}
\begin{definition}[Parallel transport]
    \label{definition:parallel-transport}
    \index{Parallel transport}%
    The linear isomorphism $ P_{\gamma, a \rightarrow b}:
    E_{\gamma(a)} \longrightarrow E_{\gamma(b)}$ is called the
    parallel transport along $\gamma$ with respect to $\nabla^E$.
\end{definition}
\begin{remark}[Parallel transport]
    \label{remark:parallel-transport}
    ~
    \begin{remarklist}
    \item \label{item:partrans-is-curve-dependent} In general,
        $P_{\gamma, a \rightarrow b}$ depends very much on the choice
        of the curve $\gamma$ connecting $\gamma(a)$ and $\gamma(b)$.
    \item \label{item:partrans-for-piecewise-smooth-curve} We can
        define $P_{\gamma, a \rightarrow b}$ also for piecewise smooth
        curves by composing the parallel transports of the smooth
        pieces appropriately.
    \item \label{item:partrans-and-curvature} If the curvature $R^E$
        is zero then the parallel transport $P_{\gamma, a \rightarrow
          b}$ is \emph{independent} of the curve $\gamma$ but only
        depends on $\gamma(a)$ and $\gamma(b)$, \emph{provided} the
        points are close enough. More precisely, if $\gamma$ and
        $\widetilde{\gamma}$ are two curves with $\gamma(a) =
        \widetilde{\gamma}(a)$ and $\gamma(b) = \widetilde{\gamma}(b)$
        such that there is a smooth homotopy between $\gamma$ and
        $\widetilde{\gamma}$ then $P_{\gamma, a \rightarrow b} =
        P_{\widetilde{\gamma}, a \rightarrow b}$. Note however that
        $R^E = 0$ is a rather strong condition which implies certain
        strong topological properties of the vector bundle $E
        \longrightarrow M$.
    \item \label{item:ParallelTransportReparametrization} If
        $\gamma: I \longrightarrow M$ is a smooth curve and $\sigma: J
        \longrightarrow I$ is a smooth reparametrization then the
        parallel transports along $\gamma$ and $\tilde{\gamma} =
        \gamma \circ \sigma$ coincide. More precisely, for $a', b' \in
        J$ we have
        \begin{equation}
            \label{eq:ReparametrizationInvariance}
            P_{\gamma, \sigma(a') \rightarrow \sigma(b')}
            =
            P_{\gamma \circ \sigma, a' \rightarrow b'}.
        \end{equation}
    \end{remarklist}
\end{remark}
Since the parallel transport ``connects'' the fibers of $E$ at
different points, a covariant derivative is also called
\emph{connection}\index{Connection|see{Covariant derivative}}. Some
further properties of the parallel transport are collected in the
Appendix~\ref{sec:taylor-parallel-transport}.

%
%

\subsection{The Exponential Map}
\label{subsec:exponential-map}

In the case $E = TM$ a covariant derivative has additional features we
shall discuss now. First, we have another contraction of the curvature
tensor $R$ given by
\begin{equation}
    \label{eq:ricci-curvature}
    \Ric(X,Y) = \tr (Z \mapsto R(Z, X)Y)
\end{equation}
for $X, Y \in \Secinfty(TM)$. The resulting tensor field
\begin{equation}
    \label{eq:ricci-tensor-field}
    \Ric \in \Secinfty(T^*M \tensor T^*M)
\end{equation}
is called the \emIndex{Ricci tensor} of $\nabla$. Note that the trace
in \eqref{eq:ricci-curvature} only can be defined for $E = TM$. The
third contraction $\tr (Z \mapsto R(X, Z)Y)$ would give again the
Ricci tensor up to a sign. Thus \eqref{eq:ricci-curvature} is the only
additional interesting contraction.

For a covariant derivative $\nabla$ on $TM$ we have yet another tensor
field, the \emph{torsion}\index{Torsion}
\begin{equation}
    \label{eq:torsion}
    \Tor(X,Y) = \nabla_X Y - \nabla_Y X - [X,Y],
\end{equation}
which gives a tensor field
\begin{equation}
    \label{eq:torsion-tensor}
    \Tor \in \Secinfty(\Anti^2 T^*M \tensor TM).
\end{equation}
Then $\nabla$ is called \emph{torsion-free} if $\Tor = 0$. The
relation between $R$ and $\Tor$ is encoded in the first Bianchi
identity, see e.g. \cite[Chap.~III]{kobayashi.nomizu:1963a}:
\begin{lemma}[First Bianchi identity]
    \label{lemma:first-bianchi}
    \index{First Bianchi identity}%
    \index{Covariant derivative!torsion-free}%
    For any covariant derivative for $TM$ we have
    \begin{equation}
        \label{eq:first-bianchi-identity}
        R(X, Y)Z + \textrm{cycl.}(X, Y, Z)
        = (\nabla_X Tor)(Y, Z) + \Tor(\Tor(X, Y), Z)
        + \textrm{cycl.}(X, Y, Z).
    \end{equation}
    In particular, for a torsion-free $\nabla$ we have
    \begin{equation}
        \label{eq:first-bianchi-without-torsion}
        R(X, Y)Z + R(Y, Z)X + R(Z, X)Y = 0,
    \end{equation}
    for all $X, Y, Z \in \Secinfty(TM)$.
\end{lemma}
\begin{proof}
    The proof consists in a straightforward algebraic manipulation
    using only the definitions.
\end{proof}
\begin{corollary}
    \label{corollary:ric-symmetry}
    \index{Covariant derivative!unimodular}%
    Let $\nabla$ be torsion-free. Then
    \begin{equation}
        \label{eq:ric-symmetry}
        \Ric(X,Y) - \Ric(Y,X) + (\tr R) (X,Y) = 0,
    \end{equation}
    whence $\Ric$ is symmetric if in addition $\nabla$ is unimodular.
\end{corollary}

\begin{figure}
    \centering
    \input{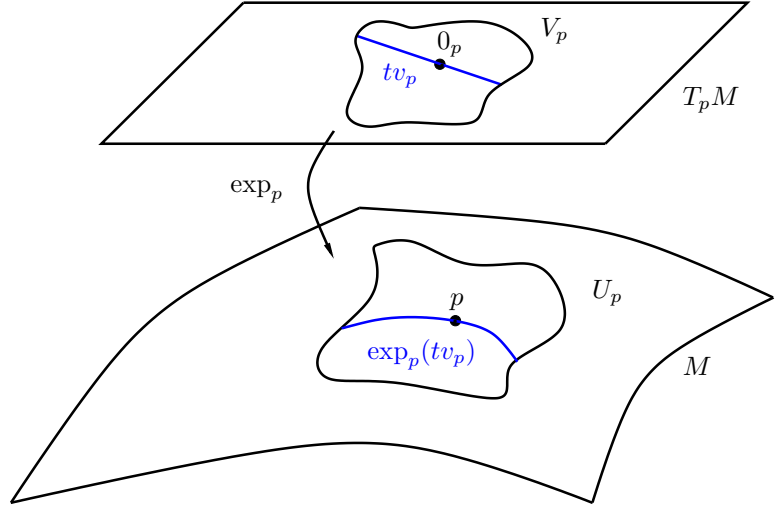}
    \caption{\label{fig:ExponentialMap}%
      The exponential map gives a normal chart.
    }
\end{figure}

In case of the tangent bundle the parallel transport can be used to
motivate the following question. For a starting point $p \in M$ and a
starting velocity $v_p \in T_pM$, is there a curve $\gamma$ with
$\dot{\gamma}(0) = v_p$ such that $\dot{\gamma}$ is parallel along
$\gamma$? Such an
\emph{auto-parallel}\index{Auto-parallel|see{Geodesic}} curve will be
called a \emph{geodesic}\index{Geodesic}.  To get an idea we consider
this condition, which globally reads
\begin{equation}
    \label{eq:geodesic-condition}
    \nabla^\#_{\frac{\partial}{\partial t}} \dot{\gamma} = 0,
\end{equation}
in a local chart $(U, x)$. We denote by
\begin{equation}
    \label{eq:christoffel-symbols}
    \index{Christoffel symbol}%
    \nabla_{\frac{\partial}{\partial x^i}}
    \frac{\partial}{\partial x^j}
    = \Gamma^k_{ij} \frac{\partial}{\partial x^k}
\end{equation}
the locally defined \emph{Christoffel symbols} $\Gamma^k_{ij} \in
\Cinfty(U)$. Then \eqref{eq:geodesic-condition} means for the curve
$\gamma: I \subseteq \mathbb{R} \longrightarrow M$ with
$\dot{\gamma}(t) = \dot{\gamma}^i(t) \frac{\partial}{\partial x^i}$
and $\gamma^i = x^i \circ \gamma \in \Cinfty(I)$ explicitly
\begin{equation}
    \label{eq:local-geodesic-condition}
    \ddot{\gamma}^i(t)
    + \Gamma^i_{k\ell} (\gamma(t))
    \dot{\gamma}^k(t) \dot{\gamma}^\ell(t)
    = 0.
\end{equation}
This is a (highly nonlinear) ordinary second order differential
equation. Hence we have unique solutions for every initial condition
$\dot{\gamma}^i(0) \frac{\partial}{\partial x^i} \At{\gamma(0)} = v_p$,
where $p = \gamma(0)$, at least for small times. Since locally
\begin{equation}
    \label{eq:torsion-locally-with-christoffels}
    \Tor^k_{ij} = \Gamma^k_{ij} - \Gamma^k_{ji},
\end{equation}
we see that the torsion $\Tor$ of $\nabla$ does not enter the geodesic
equation \eqref{eq:local-geodesic-condition}. We collect a few
well-known facts about the solution theory of the geodesic equation:
\begin{theorem}[Geodesics]
    \label{theorem:geodesics}
    Let $\nabla$ be a covariant derivative for $TM \longrightarrow M$.
    \begin{theoremlist}
    \item \label{item:existence-uniqueness-geodesics} For every $v_p
        \in T_p M$ there exists a unique solution $\gamma: I_{v_p}
        \subseteq \mathbb{R} \longrightarrow M$ of
        \eqref{eq:local-geodesic-condition} with $\dot{\gamma}(0) =
        v_p$ and maximal open interval $I_{v_p} \subseteq \mathbb{R}$
        around $0$.
    \item \label{item:rescaling-of-geodesics} Let $\lambda \in
        \mathbb{R}$ and $v_p \in T_p M$. If $\gamma$ denotes the
        geodesic with $\dot{\gamma}(0) = v_p$ then $\gamma_\lambda(t)
        = \gamma(\lambda t)$ is the geodesic with
        $\dot{\gamma}_\lambda(0) = \lambda v_p$.
    \item \label{item:exponential-map} There exists an open
        neighborhood $\mathcal{V} \subseteq TM$ of the zero section
        such that for all $v_p \in \mathcal{V}$ the geodesic with
        $\dot{\gamma}(0) = v_p$ is defined for all $t \in [0,1]$. We
        set $\exp_p (v_p) = \gamma(1)$ for this geodesic.
    \item \label{item:geodesic-through-exp} For $v_p \in \mathcal{V}
        \subseteq TM$ the curve $t \mapsto \exp_p (t v_p)$ is the
        geodesic $\gamma$ with $\dot{\gamma}(0) = v_p$.
    \item \label{item:exp-is-smooth} The map $\exp: \mathcal{V}
        \subset TM \longrightarrow M$ is smooth.
    \item \label{item:exp-x-pi-diffeomorphism} The map
        \begin{equation}
            \label{eq:exp-x-pi}
            \pi \times \exp: \mathcal{V} \subseteq TM \ni v_p
            \; \mapsto \;
            (p, \exp_p(v_p)) \in M \times M
        \end{equation}
        is a local diffeomorphism around the zero-section. It maps the
        zero section onto the diagonal and for all $p \in M$
        \begin{equation}
            \label{eq:exp-differential-at-zero}
            T_{0_p} \exp_p = \id_{T_p M}.
        \end{equation}
    \end{theoremlist}
\end{theorem}
\begin{proof}
    The proof can be found e.g. in \cite[Chap.~VIII, §5]{lang:1999a}
    or \cite[§11 and §12]{broecker.jaenich:1990a}.
\end{proof}
\begin{definition}[Exponential map]
    \label{definition:exponential-map}
    \index{Exponential map}%
    \index{Covariant derivative!exponential map}%
    For a given covariant derivative $\nabla$, the map $\exp:
    \mathcal{V} \subseteq TM \longrightarrow M$ given by
    \refitem{item:exp-is-smooth} of Theorem~\ref{theorem:geodesics} is
    called the exponential map of $\nabla$.
\end{definition}
\begin{remark}[Exponential map]
    \label{remark:expnonential-map}
    Let $\nabla$ be a covariant derivative on $M$.
    \begin{remarklist}
    \item \label{item:can-use-torsion-free-cov} Since the geodesic
        equation does not depend on the antisymmetric part of the
        $\Gamma^k_{ij}$ we can safely pass from $\nabla$ to a
        torsion-free covariant derivative by adding the appropriate
        multiple of the torsion tensor. The geodesics do not change
        and neither does the exponential map.
    \item \label{item:spray} The exponential map is best understood in
        terms of \emph{spray vector fields}\index{Spray} on $TM$, see
        e.g.~\cite[Chap.~VIII, §5]{lang:1999a} or \cite[§11 and
        §12]{broecker.jaenich:1990a}. In fact, $\exp$ is just the
        projection of the time-one-flow of the spray vector field
        associated to $\nabla$ by the bundle projection.
    \item \label{item:exp-induces-local-diffeo} It follows from
        \eqref{eq:exp-differential-at-zero} that the exponential map
        $\exp_p$ at a given point $p \in M$ induces a diffeomorphism
        \begin{equation}
            \label{eq:exp-local-diffeo}
            \index{Normal chart}%
            \index{Geodesic chart}%
            \index{Normal coordinates}%
            \exp_p: V_p \subseteq T_p M
            \longrightarrow U_p \subseteq M
        \end{equation}
        between a sufficiently small open neighborhood $V_p \subseteq
        T_p M$ of $0_p$ and its image $U_p \subseteq M$ which becomes
        an open neighborhood of
        \begin{equation}
            \label{eq:5}
            p = \exp_p (0_p) \in U_p \subseteq M.
        \end{equation}
        Thus the map $(\exp_p \at{V_p})^{-1}: U_p \longrightarrow V_p$
        yields a chart of $M$ centered around $p$ which is called a
        \emph{normal} or \emph{geodesic chart} with respect to
        $\nabla$. The choice of linear coordinates on $V_p \subset T_p
        M$ induces then \emph{normal coordinates} on $ U_p \subseteq
        M$, see also Figure~\ref{fig:ExponentialMap}.
      \end{remarklist}
\end{remark}
More details on properties of the exponential map can be found in the
Appendix~\ref{sec:jacobi-and-tangent-of-exp} where we compute, among
other things, the Taylor expansions of various objects in normal
coordinates.
The following definition is motivated by the flat situation where the
notions of ``star-shaped'' and ``convex'' have an immediate meaning.
\begin{definition}
    \label{definition:geodesic-starshaped-convex}
    \index{Geodesically star-shaped}%
    \index{Geodesically convex}%
    \index{Convex|see{Geodesically convex}}%
    \index{Star-shaped|see{Geodesically star-shaped}}%
    An open subset $U \subseteq M$ is called
    \begin{definitionlist}
    \item \label{item:geod-starshaped} geodesically star-shaped with
        respect to $p \in M$ if there is a star-shaped $V \subseteq
        V_p$ with $\exp_p \at{V}: V \stackrel{\simeq}{\longrightarrow}
        U = \exp_p(V)$.
    \item \label{item:geod-convex} geodesically convex if it is
        geodesically star-shaped with respect to any point $p \in U$.
    \end{definitionlist}
\end{definition}
Usually, we simply speak of star-shaped and convex open subsets of $M$
if the reference to $\nabla$ is clear.  Note that the properties
described in Definition~\ref{definition:geodesic-starshaped-convex}
depend on the choice of $\nabla$ and are \emph{not} invariant under an
arbitrary change of coordinates.

For a general covariant derivative it might well be that the domain of
definition of $\exp$ is a proper open subset: geodesics need not be
defined for all times but can ``fall of the manifold''. The simplest
example is obtained from $\mathbb{R}^2 \backslash \{0\}$ with the flat
connection. Geodesics are straight lines. Thus the geodesic starting
at $(-1,0)$ with tangent vector $(1,0)$ stops being defined at $t=1$
since it would reach $0$ which is not a part of $\mathbb{R}^2
\backslash \{0\}$. While this example looks rather artificial there
are more difficult situations where one can not just ``add a few
points''. These considerations motivate the following definition:
\begin{definition}[Geodesic completeness]
    \label{definition:geodesic-completeness}
    \index{Geodesic completeness}%
    The covariant derivative $\nabla$ is called geodesically complete
    if all geodesics are defined for all times.
\end{definition}

%
%

\subsection{Levi-Civita-Connection and the d'Alembertian}
\label{subsec:levi-civita-and-dAlembertian}

We shall now specialize the connection $\nabla$ further and add one
more structure, namely a semi-Riemannian metric:
\begin{definition}[Semi-Riemannian metric]
    \label{definition:semiRiem-metric}
    \index{Metric!semi-Riemannian}%
    \index{Metric!Riemannian}%
    \index{Metric!Lorentzian}%
    \index{Signature}%
    A section $g \in \Secinfty(\Sym^2 T^*M)$ is called
    semi-Rie\-mann\-ian metric if the bilinear form $g_p \in \Sym^2
    T^*_pM$ on $T_pM$ is non-degenerate for all $p \in M$. If in
    addition $g_p$ is positive definite for all $p \in M$ then $g$ is
    called Riemannian metric. If $g_p$ has signature $(+, -, \ldots,
    -)$ then g is called Lorentz metric.
\end{definition}
\begin{remark}[Semi-Riemannian metrics]
    \label{remark:semiRiem-metrics}
    ~
    \begin{remarklist}
    \item \label{item:signature-locally-constant} The signature of a
        semi-Riemannian metric is locally constant and hence constant
        on a connected manifold, since it depends continuously on $p$
        and has only discrete values.
    \item \label{item:signature-lorentz} For Lorentz metrics also the
        opposite signature $(-, +, \ldots, +)$ is used in the
        literature. This causes some confusions and funny signs. So be
        careful here! Our convention is the more common one in quantum
        field theory, while the other one is preferred in general
        relativity.
    \end{remarklist}
\end{remark}
A semi-Riemannian metric specifies a unique covariant derivative and a
unique positive density:
\begin{proposition}
    \label{proposition:levi-civita-connection-etc}
    \index{Covariant derivative!Levi-Civita}%
    \index{Density!metric}%
    \index{Covariant derivative!unimodular}%
    Let $g$ be a semi-Riemannian metric on $M$.
    \begin{propositionlist}
    \item \label{iten:levi-civita-connection} There exists a unique
        torsion-free covariant derivative $\nabla$, the Levi-Civita
        connection, such that
        \begin{equation}
            \label{eq:levi-civita-is-metric}
            \nabla g = 0.
        \end{equation}
    \item \label{item:Riem-density} There exists a unique positive
        density $\mu_g \in \Secinfty(\Dichten T^*M)$ such that
        \begin{equation}
            \label{eq:Riem-density-normalization}
            \mu_g \at{p} (v_1, \ldots, v_n) = 1,
        \end{equation}
        whenever $v_1, \ldots, v_n$ form a basis of $T_pM$ with
        $|g_p(v_i, v_j)| = \delta_{ij}$. In a chart $(U, x)$ we have
        \begin{equation}
            \label{eq:Riem-density-coordinate-expression}
            \mu_g \at{U}
            =
            \sqrt{|\det (g_{ij})|} \;
            |\D x^1 \wedge \cdots \wedge \D x^n|,
        \end{equation}
        with $g_{ij} = g( \frac{\partial}{\partial x^i},
        \frac{\partial}{\partial x^j} )$.
    \item \label{item:Riem-densitiy-is-cov-const} The density $\mu_g$
        is covariantly constant with respect to the Levi-Civita
        connection,
        \begin{equation}
            \label{eq:Riem-density-is-cov-const}
            \nabla \mu_g = 0.
        \end{equation}
        Thus $\nabla$ is unimodular.
    \end{propositionlist}
\end{proposition}
\begin{proof}
    The proof is very much standard and will be omitted here, see e.g.
    \cite[Aufgabe~3.7 and 5.10]{waldmann:2007a}.
\end{proof}
\begin{remark}[Semi-Riemannian metrics]
    \label{remark:semiRiem-metrics-2}
    Let $g$ be a semi-Riemannian metric on $M$.
    \begin{remarklist}
    \item \label{item:semiRiem-geodesics} For a semi-Riemannian metric
        we have a notion of geodesics, namely those with respect to
        the corresponding Levi-Civita connection.
    \item \index{Covariant divergence}%
        \index{Divergence}%
        \index{Density!divergence}%
        \label{item:cov-div-equals-density-div} The covariant
        divergence $\divergenz_\nabla (X)$ of a vector field $X \in
        \Secinfty(TM)$ and the divergence with respect to the density
        $\mu_g$, i.e.
        \begin{equation}
            \label{eq:semiRiem-density-covariance}
            \divergenz_{\mu_g} (X) = \frac{\Lie_X \mu_g}{\mu_g}
        \end{equation}
        coincide: We have
        \begin{equation}
            \label{eq:cov-div-equals-density-div}
            \divergenz_\nabla (X) = \divergenz_{\mu_g} (X),
        \end{equation}
        which follows immediately from
        Lemma~\ref{lemma:mu-divergence-and-covariant-divergence}, see
        also \cite[Sect.~2.3.4]{waldmann:2007a}, since $\nabla \mu_g =
        0$.  Thus we shall speak of \emph{the} divergence and simply
        write
        \begin{equation}
            \label{eq:divergence}
            \divergenz(X) = \divergenz_\nabla (X)
            = \divergenz_{\mu_g} (X)
        \end{equation}
        on a semi-Riemannian manifold.
    \item \index{Musical isomorphism}%
        \label{item:sharp-and-flat}
        Since $g \in \Secinfty(\Sym^2 T^*M)$ is non-degenerate it
        induces a \emph{musical isomorphism}
        \begin{equation}
            \label{eq:flat}
            \flat: T_pM \ni v_p \mapsto v_p^\flat
            = g(v_p, \argument) \in T_p^*M,
        \end{equation}
        which gives a vector bundle isomorphism
        \begin{equation}
            \label{eq:flat-isomorphism}
            \flat: TM \longrightarrow T^*M.
        \end{equation}
        The inverse of $\flat$ is usually denoted by
        \begin{equation}
            \label{eq:sharp-isomoprhism}
            \sharp: T^*M \longrightarrow TM.
        \end{equation}
        Extending $\flat $ and $\sharp$ to higher tensor powers we get
        musical isomorphisms also between all corresponding
        contravariant and covariant tensor bundles. If locally in a
        chart $(U, x)$
        \begin{equation}
            \label{eq:g-local}
            g\at{U} = \frac{1}{2} g_{ij} \D x^i \vee \D x^j,
        \end{equation}
        then $v^\flat = g_{ij} v^i \D x^j$, where $v = v^i
        \frac{\partial}{\partial x^i}$. If $g^{ij}$ denotes the
        inverse matrix to the $g_{ij}$ from \eqref{eq:g-local}, i.e.
        $g^{ij} g_{jk} = \delta_{ik}$, then
        \begin{equation}
            \label{eq:sharp-local}
            \alpha^\sharp
            = g^{ij} \alpha_i \frac{\partial}{\partial x^j}
        \end{equation}
        for a one-form $\alpha = \alpha_i \D x^i$. This motivates the
        notion ``musical'' as $\flat$ lowers the indexes while
        $\sharp$ raises them. Finally, we have the dual metric locally
        given by
        \begin{equation}
            \label{eq:dual-metric-local}
            g^{-1} \at{U} =
            \frac{1}{2} g^{ij} \frac{\partial}{\partial x^i}
            \vee \frac{\partial}{\partial x^j},
        \end{equation}
        which is a global section $g^{-1} \in \Secinfty(\Sym^2 TM)$.
    \item \index{Kinetic energy}%
        \index{Hamiltonian}%
        \index{Geometric mechanics}%
        \index{Euler-Lagrange equation}%
        \index{Hamilton equation}%
        \index{Lagrangian}%
        \index{Hamiltonian flow}%
        \label{item:kinetic-energy} The metric $g \in \Secinfty(\Sym^2
        T^*M)$ can equivalently be interpreted as a homogeneous
        quadratic function on $TM$ via the usual canonical isomorphism
        from Remark~\ref{remark:standard-ordered-quantization}. The
        function
        \begin{equation}
            \label{eq:kinetic-energy-Lagrange}
            T = \mathcal{J}(g) \in \Pol^2 (TM)
        \end{equation}
        is then usually called the \emph{kinetic energy function} in
        the Lagrangian picture of mechanics. Analogously, $g^{-1} \in
        \Secinfty(\Sym^2 TM)$ gives a homogeneous quadratic function
        \begin{equation}
            \label{eq:kinetic-energy-Hamilton}
            T = \mathcal{J}(g^{-1}) \in \Pol^2 (T^*M)
        \end{equation}
        on $T^*M$, the kinetic energy in the Hamiltonian picture of
        mechanics. It turns out that all notions of geodesics etc. can
        be understood in this geometric mechanical framework. For
        example, geodesics are just the base point curves of solutions
        of the Euler-Lagrange equations and Hamilton's equations with
        respect to the Lagrangian $L = T$ and Hamiltonian $H = T$,
        respectively. Thus geodesic motion is motion \emph{without}
        additional forces induced by some addition potentials. The
        exponential map $\exp$ is then just the Hamiltonian flow of
        $T$ at time $t = 1$ projected back to $M$. For more on this
        mechanical point of view, see
        e.g.~\cite[Sect.~3.2.2]{waldmann:2007a}.
    \item \index{Christoffel symbol}%
        \label{item:christoffels-locally} Using the inverse matrix
        $g^{ij}$ we have the following local Christoffel symbols of
        the Levi-Civita connection
        \begin{equation}
            \label{eq:christoffel-local}
            \Gamma_{ij}^k = \frac{1}{2} g^{k\ell}
            \left(
                \frac{\partial g_{\ell i}}{\partial x^j}
                + \frac{\partial g_{\ell j}}{\partial x^i}
                - \frac{\partial g_{ij}}{\partial x^\ell}
            \right).
        \end{equation}
    \end{remarklist}
\end{remark}
Since by Corollary~\ref{corollary:ric-symmetry} and
Proposition~\ref{proposition:levi-civita-connection-etc},
\refitem{item:Riem-densitiy-is-cov-const} for a semi-Riemannian
manifold $(M, g)$ the \Index{Ricci tensor} $\Ric$ is in fact symmetric
\begin{equation}
    \label{eq:ric-symmetric-tensor}
    \Ric \in \Secinfty(\Sym^2 T^*M),
\end{equation}
we can compute a further ``trace'' by using the metric $g$. Note that
while $\Ric$ can be defined for every covariant derivative this
further contraction requires $g$. One calls the function
\begin{equation}
    \label{eq:scalar-curvature}
    \scal = \SP{g^{-1}, \Ric} \in \Cinfty(M)
\end{equation}
the \emph{scalar curvature}. Locally, $\scal$ is just
\begin{equation}
    \label{eq:scalar-curvature-local}
    \index{Scalar curvature}%
    \scal \at{U} = g^{ij} \Ric_{ij}.
\end{equation}
In the literature, there are many other notations for $\scal$, e.g.
$R$ (without indexes) or $s$ or $S$.

We come now to differential operators defined by means of a
semi-Riemannian metric. We have already seen the divergence operator
$\divergenz$ which acts on vector fields and which can be extended as
in Lemma~\ref{lemma:covariant-divergence-of-multivector-fields} to all
sections $\Secinfty(\Sym^\bullet TM)$. We have two other important
operators.
\begin{definition}[Gradient and d'Alembertian]
    \label{definition:gradient-and-dAlembertian}
    \index{Gradient}%
    \index{dAlembertian@d'Alembertian}%
    \index{Laplacian}%
    On a semi-Riemannian manifold $(M, g)$ the gradient of a function
    is defined by
    \begin{equation}
        \label{eq:gradient}
        \gradient f = (\D f)^\sharp \in \Secinfty(TM)
    \end{equation}
    and the d'Alembertian of a function $f \in \Cinfty(M)$ is
    \begin{equation}
        \label{eq:dAlembertian}
        \dAlembert f = \divergenz (\gradient f) \in \Cinfty(M).
    \end{equation}
    In case of a Riemannian manifold we write $\Laplace f = \divergenz
    (\gradient f)$ instead and call $\Laplace$ the Laplacian.
\end{definition}
\begin{remark}
    \label{remark:laplacian-sign-convention}
    There are different sign conventions in the definition of the
    Laplacian and the d'Alembertian. In particular, sometimes $-
    \Laplace$ is favoured instead of our $\Laplace$ since $\Laplace$
    as we defined it turns out to be a \emph{negative} essentially
    selfadjoint operator on $\Cinfty(M)$ for compact $M$.
\end{remark}
We discuss now a couple of local formulas which allow to handle the
operators $\divergenz$, $\gradient$ and $\dAlembert$ more explicitly.
\begin{proposition}
    \label{proposition:grad-div-box-local-expressions}
    Let $(M, g)$ be a semi-Riemannian manifold and let $(U, x)$ be a
    chart of $M$.
    \begin{propositionlist}
    \item \label{item:grad-local} The gradient of $f \in \Cinfty(M)$
        is locally given by
        \begin{equation}
            \label{eq:grad-local}
            \gradient (f) \at{U}
            = g^{ij} \frac{\partial f}{\partial x^i}
            \frac{\partial}{\partial x^j}.
        \end{equation}
    \item \label{item:div-local} The divergence of $X \in
        \Secinfty(TM)$ is locally given by
        \begin{equation}
            \label{eq:div-local}
            \divergenz (X) \at{U}
            = \frac{\partial X^i}{\partial x^i} + \Gamma_{ki}^k X^i.
        \end{equation}
    \item \label{item:dAlembert-local} The d'Alembertian of $f \in
        \Cinfty(M)$ is locally given by
        \begin{equation}
            \label{eq:dAlembert-local}
            \dAlembert f \at{U} = g^{ij}
            \left(
                \frac{\partial^2 f}{\partial x^i \partial x^j}
                - \Gamma_{ij}^k \frac{\partial f}{\partial x^k}
            \right).
        \end{equation}
    \item \label{item:dAlembert-is-second-order-diffop} The
        d'Alembertian is a second order differential operator with
        leading symbol
        \begin{equation}
            \label{eq:dAlembert-leading-symbol}
            \sigma_2(\dAlembert) = 2 g^{-1} \in \Secinfty(\Sym^2 TM).
        \end{equation}
        Moreover, with respect to the global symbol calculus induced
        by the Levi-Civita connection we have
        \begin{equation}
            \label{eq:dAlembert-as-stdrep}
            \dAlembert
            = \left( \frac{\I}{\hbar} \right)^2 \stdrep(2 g^{-1}),
        \end{equation}
        whence
        \begin{equation}
            \label{eq:dAlembert-as-trace}
            \dAlembert f = \frac{1}{2} \SP{g^{-1}, \SymD^2 f}.
        \end{equation}
    \end{propositionlist}
\end{proposition}
\begin{proof}
    The local formulas \eqref{eq:grad-local} and \eqref{eq:div-local}
    are clear. Then \eqref{eq:dAlembert-local} follows from some
    straightforward computation using the precise form of
    \eqref{eq:christoffel-local} for the Christoffel symbols. Then
    \eqref{eq:dAlembert-leading-symbol} is clear by definition of the
    leading symbol. For \eqref{eq:dAlembert-as-stdrep} and
    \eqref{eq:dAlembert-as-trace} we compute
    \begin{align*}
        \SymD^2 f &= \SymD \D f
        = \D x^i \vee \nabla_{\frac{\partial}{\partial x^i}}
        \left(
            \frac{\partial f}{\partial x^j} \D x^j
        \right) \\
        &= \D x^i \vee \frac{\partial^2 f}{\partial x^i \partial x^j}
        \D x^j
        + \D x^i \vee \frac{\partial f}{\partial x^j}
        \nabla_{\frac{\partial}{\partial x^i}} \D x^j \\
        &=
        \frac{\partial^2 f}{\partial x^i \partial x^j}
        \D x^i \vee \D x^j
        - \Gamma_{ik}^j \frac{\partial f}{\partial x^j}
        \D x^i \vee \D x^k,
    \end{align*}
    which gives
    \[
    \SymD^2 f =
    \left(
        \frac{\partial^2 f}{\partial x^i \partial x^j}
        - \Gamma_{ij}^k \frac{\partial f}{\partial x^k}
    \right)
    \D x^i \vee \D x^j
    \]
    for a general connection $\nabla$. For $g^{-1} = \frac{1}{2}
    g^{ij} \frac{\partial}{\partial x^i} \vee \frac{\partial}{\partial
      x^j}$ we find
    \[
    \SP{g^{-1}, \SymD^2 f}
    = 2 g^{ij}
    \left(
        \frac{\partial^2 f}{\partial x^i \partial x^j}
        - \Gamma_{ij}^k \frac{\partial f}{\partial x^k}
    \right)
    = 2 \dAlembert f.
    \]
\end{proof}
\begin{remark}[Hessian]
    \label{remark:hessian}
    Sometimes $\frac{1}{2} \SymD^2 f \in \Secinfty(\Sym^2 T^*M)$ is
    also called the \emIndex{Hessian}
    \begin{equation}
        \label{eq:hessian}
        \mathrm{Hess}(f)
        = \frac{1}{2} \SymD^2 f \in \Secinfty(\Sym^2 T^*M).
    \end{equation}
    Then the d'Alembertian is the trace of the Hessian with respect to
    $g^{-1}$. Moreover, the gradient $\gradient: \Cinfty(M)
    \longrightarrow \Secinfty(TM)$ is a differential operator of order
    one, the same holds for the divergence $\divergenz: \Secinfty(TM)
    \longrightarrow \Cinfty(M)$.
\end{remark}
\begin{remark}
    \label{remark:leibniz-rules}
    For later use we also mention the following Leibniz rules
    \begin{equation}
        \label{eq:grad-leibniz}
        \gradient (fg) = g \gradient (f) + f \gradient (g),
    \end{equation}
    \begin{equation}
        \label{eq:div-leibniz}
        \divergenz (fX)= f \divergenz (x) + X(f),
    \end{equation}
    \begin{equation}
        \label{eq:dAlembert-leibniz}
        \dAlembert (fg)
        =
        g \dAlembert f + \gradient(g) f + \gradient(f) g
        + f \dAlembert g
        = g \dAlembert f
        + 2 \SP{\gradient(f), \gradient(g)}
        + f \dAlembert g,
    \end{equation}
    for $f,g  \in \Cinfty(M)$ and $X \in \Secinfty(TM)$. They can
    easily be obtained from the definitions.
\end{remark}
\begin{example}[Minkowski spacetime]
    \label{example:grad-div-dAlembert-in-Minkowskispace}
    \index{Minkowski spacetime}%
    \index{Minkowski metric}%
    \index{Special relativity}%
    We consider the $n$-dimensional Minkowski spacetime. As a manifold
    we have $M = \mathbb{R}^n$ with canonical coordinates $x^0, x^1,
    \ldots, x^{n-1}$. Then the Minkowski metric $\eta$ on $M$ is the
    \emph{constant} metric
    \begin{equation}
        \label{eq:minkowski-metric}
        \eta = \frac{1}{2} \eta_{ij} \D x^i \vee \D x^j
    \end{equation}
    with $(\eta_{ij}) = \diag (+1, -1, \ldots, -1)$. One easily
    computes that in this global chart all Christoffel symbols vanish:
    $(M, \eta)$ is \emph{flat}. Moreover, we have for the above
    differential operators
    \begin{equation}
        \label{eq:grad-in-Minkowskispace}
        \gradient f
        = \frac{\partial f}{\partial x^0}\frac{\partial}{\partial x^0}
        - \sum_{i=1}^{n-1}
        \frac{\partial f}{\partial x^i}\frac{\partial}{\partial x^i},
    \end{equation}
    \begin{equation}
        \label{eq:div-in-Minkowskipace}
        \divergenz X
        = \frac{\partial X^0}{\partial x^0}
        + \sum_{i=1}^{n-1} \frac{\partial X^i}{\partial x^i},
    \end{equation}
    \begin{equation}
        \label{eq:dAlembert-in-Minkowskispace}
        \dAlembert f
        = \frac{\partial^2 f}{\partial (x^0)^2}
        - \sum_{i=1}^{n-1} \frac{\partial^2 f}{\partial (x^i)^2}.
    \end{equation}
    This shows that $\dAlembert$ is indeed the usual wave operator or
    d'Alembertian as known from the theory of special relativity, see
    e.g. \cite{roemer.forger:1993a}.  Finally, the Lorentz density
    with respect to $\eta$ is just the usual Lebesgue measure
    \begin{equation}
        \label{eq:lorentz-densitiy-in-Minkowskispace}
        \index{Lorentz density}%
        \index{Lebesgue measure}%
        \mu_\eta = |\D x^0 \wedge \cdots \wedge \D x^{n-1}|.
    \end{equation}
\end{example}

%
%

\subsection{Normally Hyperbolic Differential Operators}
\label{subsec:normally-hyperbolic-diffops}

The aim of this subsection is to generalize the d'Alembertian to more
general fields than scalar fields. As it will turn out later, the most
important feature of $\dAlembert$ is the fact that the leading symbol
is given by the metric. This motivates the following definition:
\begin{definition}[Normally hyperbolic operator]
    \label{definition:normally-hyperbolic-operator}
    \index{Normally hyperbolic operator}%
    \index{Leading symbol}%
    Let $E \longrightarrow M$ be a vector bundle over a Lorentz
    manifold $(M, g)$. A differential operator $D: \Secinfty(E)
    \longrightarrow \Secinfty(E)$ is called normally hyperbolic if it
    is of second order and
    \begin{equation}
        \label{eq:normally-hyperbolic-condition}
        \sigma_2 (D) = 2 g^{-1} \tensor \id_E.
    \end{equation}
\end{definition}
Recall that $\sigma_2 (D) \in \Secinfty(\Sym^2 TM \tensor \End(E))$
which explains the second tensor factor in
\eqref{eq:normally-hyperbolic-condition}. Usually, we simply write
$\sigma_2(D) = 2 g^{-1}$ with some slight abuse of notation. Note also
that, as already for the d'Alembertian itself, the factor $2$ in the
symbol comes from our convention for symbols. Here also other
conventions are used in the literature. However, this will not play
any role later. The important fact is that $D$ has a symbol being just
a \emph{constant nonzero} multiple of $g^{-1}$.

The following construction will always lead to a normally hyperbolic
operator:
\begin{example}[Connection d'Alembertian]
    \label{example:connection-dAlembert}
    \index{Connection dAlembertian@Connection d'Alembertian}%
    Let $\nabla^E$ be a covariant derivative for $E \longrightarrow M$
    and let $\nabla$ be the Levi-Civita connection. This yields a
    global symbol calculus whence by
    \begin{equation}
        \label{eq:connection-dAlembert}
        \dAlembert^\nabla
        = \left( \frac{\I}{\hbar} \right)^2
        \stdrep (2 g^{-1} \tensor \id_E)
        = \frac{1}{2}
        \SP{2 g^{-1} \tensor \id_E, \frac{1}{2}(\SymD^E)^2 \argument}
    \end{equation}
    a second order differential operator is given with leading symbol
    \begin{equation}
        \label{eq:6}
        \sigma_2 (\dAlembert^\nabla)
        = \left( \frac{\I}{\hbar} \right)^2
        \sigma_2 (\stdrep(2 g^{-1} \tensor \id_E))
        = 2 g^{-1} \tensor \id_E
    \end{equation}
    by Theorem~\ref{theorem:Global-Symbol-Calculus}. Thus
    $\dAlembert^\nabla$ is normally hyperbolic for any choice of
    $\nabla^E$. An operator of this type is called the
    \emph{connection d'Alembertian} with respect to $\nabla^E$.
\end{example}
\begin{lemma}[Connection d'Alembertian]
    \label{lemma:connection-dAlembert}
    Let $\nabla^E$ be a covariant derivative for $E \longrightarrow M$
    and $\dAlembert^\nabla$ the corresponding connection
    d'Alembertian.
    \begin{lemmalist}
    \item \label{item:connection-dAlembert-leibniz} For $f \in
        \Cinfty(M)$ and $s \in \Secinfty(E)$ we have
        \begin{equation}
            \label{eq:connection-dAlembert-leibniz}
            \dAlembert^\nabla (fs)
            = (\dAlembert f) s + 2 \nabla^E_{\gradient(f)} s
            + f \dAlembert^\nabla s.
        \end{equation}
    \item \label{item:connection-dAlembert-local} Let $A_{i
          \beta}^\alpha = e^\alpha
        \left(\nabla^E_{\frac{\partial}{\partial x^i}} e_\beta\right)
        \in \Cinfty(U)$ denote the local Christoffel symbols with
        respect to a chart $(U, x)$ and local base sections $e_\alpha
        \in \Secinfty(E \at{U})$. Then locally
        \begin{equation}
            \label{eq:connection-dAlembert-local}
            \dAlembert^\nabla s =
            \left(
                g^{ij}
                \frac{\partial^2 s^\alpha}{\partial x^i \partial x^j}
                + 2 g^{ij} \frac{\partial s^\gamma}{\partial x^i}
                A_{j \gamma}^\alpha
                - g^{ij} \Gamma_{ij}^k
                \frac{\partial s^\alpha}{\partial x^k}
                + g^{ij}
                \left(
                    \frac{\partial A_{i \beta}^\alpha}{\partial x^j}
                    - A_{k \beta}^\alpha \Gamma_{ij}^k
                    + A_{i \beta}^\gamma A_{j \gamma}^\alpha
                \right)
                s^\beta
            \right)
            e_\alpha.
        \end{equation}
    \end{lemmalist}
\end{lemma}
\begin{proof}
    For the first part we use
    Proposition~\ref{proposition:Symmetric-Covariant-Derivative} to
    compute
    \[
    (\SymD^E)^2 (fs)
    = \SymD^E (\D f \tensor s + f \SymD^E s)
    = \SymD \D f \tensor s + 2 \D f \vee \SymD^E s + f (\SymD^E)^2 s.
    \]
    Then for the natural pairing we have
    \begin{align*}
        \dAlembert^\nabla (fs)
        &=
        \frac{1}{2} \SP{g^{-1}, (\SymD^E)^2 (f \cdot s)} \\
        &=
        \frac{1}{2} \SP{g^{-1}, \SymD \D f} \cdot s
        +  \SP{g^{-1}, \D f \vee \SymD^E s}
        + \frac{1}{2} f \SP{g^{-1}, (\SymD^E)^2 s} \\
        &=
        \dAlembert f \cdot s
        + 2 g^{ij} \frac{\partial f}{\partial x^i}
        \nabla_{\frac{\partial}{\partial x^j}}^E s
        + f \dAlembert^\nabla s \\
        &=
        \dAlembert f \cdot s
        + 2 \nabla_{\gradient(f)}^E s + f \dAlembert^\nabla s,
    \end{align*}
    proving the first part. For the second, let $A_{i\beta}^\alpha$ be
    the local Christoffel symbols. Then first we have
    \begin{align*}
        \SymD^E s
        & = \D x^i \tensor \nabla_{\frac{\partial}{\partial x^i}}^E s
        = \D x^i \tensor
        \left(
            \frac{\partial s^\alpha}{\partial x^i} e_\alpha
            + s^\alpha \nabla_{\frac{\partial}{\partial x^i}}^E
            e_\alpha
        \right)
        = \D x^i \tensor \frac{\partial s^\alpha}{\partial x^i}
        e_\alpha
        + \D x^i \tensor s^\alpha A_{i\alpha}^\beta e_\beta \\
        & = \left(
            \frac{\partial s^\beta}{\partial x^i}
            + s^\alpha A_{i \alpha}^\beta
        \right)
        \D x^i \tensor e_\beta.
    \end{align*}
    Consequently, we have
    \begin{align*}
        (\SymD^E)^2 s
        & =
        \D x^j \vee
        \nabla_{\frac{\partial}{\partial x^j}} ^{E \tensor T^*M}
        \left(
            \left(
                \frac{\partial s^\beta}{\partial x^i}
                + s^\alpha A_{i \alpha}^\beta
            \right)
            \D x^i \tensor e_\beta
        \right) \\
        & =
        \D x^j \vee
        \left(
            \frac{\partial^2 s^\beta}{\partial x^i \partial x^j}
            + \frac{\partial s^\alpha}{\partial x^j}
            A_{i \alpha}^\beta
            + s^\alpha
            \frac{\partial A_{i \alpha}^\beta}{\partial x^j}
        \right)
        \D x^i \tensor e_\beta \\
        & \quad + \D x^j \vee
        \left(
            \frac{\partial s^\beta}{\partial x^i}
            + s^\alpha A_{i \alpha}^\beta
        \right)
        \left(
            -\Gamma_{jk}^i \D x^k \tensor e_\beta
            + A_{j \beta} ^\gamma \D x^i \tensor e_\gamma
        \right) \\
        & =
        \frac{\partial^2 s^\beta}{\partial x^i \partial x^j}
        \D x^i \vee \D x^j \tensor e_\beta
        + 2 \frac{\partial s^\alpha}{\partial x^i} A_{j \alpha}^\beta
        \D x^i \vee \D x^j \tensor e_\beta
        - \frac{\partial s^\beta}{\partial x^i} \Gamma_{ji}^i
        \D x^j \vee \D x^k \tensor e_\beta \\
        & \quad + s^\alpha \frac{\partial A_{i \alpha}^\beta}{\partial x^j}
        \D x^i \vee \D x^j \tensor e_\beta
        - s^\alpha A_{i \alpha}^\beta \Gamma_{jk}^i
        \D x^j \vee \D x^k \tensor e_\beta
        + s^\alpha A_{i \alpha}^\gamma A_{j \gamma}^\beta
        \D x^i \D x^j \tensor e_\beta.
    \end{align*}
    The natural pairing with $g^{-1}$ means replacing $\frac{1}{2}\D
    x^i \vee \D x^j$ with $g^{ij}$ everywhere. This gives the result.
\end{proof}

We now prove that every normally hyperbolic operator is actually a
connection d'Alembertian up to a $\Cinfty(M)$-linear operator. We have
the following result, sometimes called a generalized Weitzenböck
formula, see e.g.~\cite[Prop.~3.1]{baum.kath:1996a}:
\begin{proposition}[Weitzenböck formula]
    \label{proposition:weitzenboeck-formula}
    \index{Weitzenboeck formula@Weitzenböck formula}%
    Let $D \in \Diffop^2(E)$ be a normally hyperbolic differential
    operator. Then there exists a unique covariant derivative
    $\nabla^E$ for $E$ and a unique $B \in \Secinfty(\End(E))$ such
    that
    \begin{equation}
        \label{eq:weitzenboeck-formula}
        D = \dAlembert^\nabla + B.
    \end{equation}
\end{proposition}
\begin{proof}
    First we show uniqueness. Assume that $\nabla^E$ and $B$ exist
    such that \eqref{eq:weitzenboeck-formula} holds. Then from
    Lemma~\ref{lemma:connection-dAlembert} we know that
    \[
    D(f \cdot s) - f D(s)
    = \dAlembert^\nabla (f \cdot s) + B(f \cdot s)
    - f \dAlembert^\nabla (s) - f B(s)
    = (\dAlembert f) \cdot s + 2 \nabla_{\gradient(f)}^E s,
    \]
    since $B$ is $\Cinfty(M)$-linear. Thus we have
    \[
    \nabla^E_{\gradient (f)} s
    = \frac{1}{2}
    \left(
        D (f \cdot s) - f D (s)
        - (\dAlembert f) \cdot s
    \right)
    \tag{$*$}
    \]
    for all $f \in \Cinfty(M)$ and $s \in \Secinfty(E)$. Since
    gradients of functions span every $T_p M$ for all $p \in M$, the
    covariant derivative $\nabla^E$ is uniquely determined by $D$ via
    ($*$). But then also $B = D - \dAlembert^\nabla$ is uniquely
    determined. Let us now turn to the existence: to this end we
    compute the right hand side of ($*$) locally in order to show that
    it actually \emph{defines} a connection. Let locally
    \[
    D s \at{U}
    = g^{ij} \frac{\partial^2 s^\alpha}{\partial x^i \partial x^j}
    e_\alpha
    + D^i{}_\alpha^\beta
    \frac{\partial s^\alpha}{\partial x^i} e_\beta
    + D^\beta_\alpha s^\alpha e_\beta
    \]
    with local coefficients $D^i{}_\alpha^\beta, D_\alpha^\beta \in
    \Cinfty(U)$. Then we have
    \begin{align*}
        & \frac{1}{2} (D(f \cdot s) - f D(s)  - (\dAlembert f) \cdot s) \\
        & =  \frac{1}{2}
        \left(
            g^{ij}
            \frac{\partial^2 (f s^\alpha)}{\partial x^i \partial x^j}
            e_\alpha
            + D^i{}_\alpha^\beta
            \frac{\partial (f s^\alpha)}{\partial x^i} e_\beta
            + f D_\alpha^\beta s^\alpha e_\beta
            -f g^{ij}
            \frac{\partial^2 s^\alpha}{\partial x^i \partial x^j}
            e_\alpha
            - f D^i{}_\alpha^\beta
            \frac{\partial s^\alpha}{\partial x^i} e_\beta
            - f D_\alpha^\beta s^\alpha e_\beta
        \right. \\
        & \quad \left.  - g^{ij}
            \frac{\partial^2 f}{\partial x^i \partial x^j} s^\alpha
            e_\alpha
            + g^{ij} \Gamma_{ij}^k \frac{\partial f}{\partial x^k}
            s^\alpha e_\alpha
        \right) \\
        & = \frac{1}{2}
        \left(
            2 g^{ij} \frac{\partial f}{\partial x^i}
            \frac{\partial s^\alpha}{\partial x^j}e_\alpha
            + D^i{}_\alpha^\beta \frac{\partial f}{\partial x^i}
            s^\alpha e_\beta
            + g^{ij} \Gamma_{ij}^k \frac{\partial f}{\partial x^k}
            s^\alpha e_\alpha
        \right) \\
        & = (\gradient f)^j \frac{\partial s^\alpha}{\partial x^j}
        e_\alpha
        + \frac{1}{2}
        \left(
            D^i{}_\alpha^\beta g_{ij} (\gradient f)^j s^\alpha e_\beta
            + g^{ij} \Gamma_{ij}^k g_{k\ell} (\gradient f)^\ell
            s^\alpha e_\alpha
        \right).
    \end{align*}
    On one hand we know that the right hand side of ($*$) is globally
    defined. On the other hand, we see from the local expression that
    replacing $\gradient f$ by an arbitrary vector field $X$ defines
    locally a connection with connection one-forms
    \[
    A_{i \alpha}^\beta
    = \frac{1}{2} D^j{}_\alpha^\beta g_{ij}
    + \frac{1}{2} g^{rs} \Gamma_{rs}^j g_{ij} \delta_\alpha^\beta
    \tag{$**$},
    \]
    i.e. a connection $\nabla^E$ such that on $U$
    \[
    \nabla^E_X s
    = (\Lie_X s^\alpha) e_\alpha
    + A_{i \alpha}^\beta X^i s^\alpha e_\beta.
    \]
    This is clear from the local expression. Together we see that we
    indeed have a global connection $\nabla^E$ with local connection
    one-forms $A_{i \alpha}^\beta$ as in ($**)$. It remains to show
    that this connection $\nabla^E$ yields
    \eqref{eq:weitzenboeck-formula}. So we have to show that $D -
    \dAlembert^\nabla$ is $\Cinfty(M)$-linear. Using the explicit
    expression ($**$) for $A_{i \alpha}^\beta$ together with
    Lemma~\ref{lemma:connection-dAlembert},
    \refitem{item:connection-dAlembert-local} this is a
    straightforward computation. We have
    \begin{align*}
        Ds - \dAlembert^\nabla s
        & = D^i{}_\alpha^\beta \frac{\partial s^\alpha}{\partial x^i}
        e_\beta +
        D_\alpha^\beta s^\alpha e_\beta
        - 2 g^{ij} \frac{\partial s^\alpha}{\partial x^i}
        A_{i \alpha}^\beta e_\beta
        + g^{ij} \Gamma_{ij}^k \frac{\partial s^\alpha}{\partial x^k}
        e_\alpha \\
        & \quad - g^{ij}
        \left(
            \frac{\partial A_{i \alpha}^\beta}{\partial x^j} s^\alpha
            - A_{k \alpha}^\beta \Gamma_{ij}^k s^\alpha
            + A_{i \alpha}^\gamma A_{j \gamma} ^\beta s^\alpha
        \right)
        e_\beta \\
        & = \frac{\partial s^\alpha}{\partial x^i}
        \left(
            D^i{}_\alpha^\beta - 2 g^{ij}
            \left(
                \frac{1}{2} D^\ell{}_\alpha^\beta g_{j \ell}
                +\frac{1}{2} g^{rs} \Gamma_{rs}^\ell g_{j\ell}
                \delta_\alpha^\beta
            \right)
            + g^{jk} \Gamma_{jk}^i \delta_\alpha^\gamma
        \right)
        e_\beta \\
        & \quad + D_\alpha^\beta s^\alpha e_\beta
        - g^{ij}
        \left(
            \frac{\partial A_{i \alpha}^\beta}{\partial x^j} s^\alpha
            - A_{k \alpha}^\beta \Gamma_{ij}^k s^\alpha
            + A_{i \alpha}^\gamma A_{j \gamma}^\beta s^\alpha
        \right)
        e_\beta \\
        & = \left(
            D_\alpha^\beta
            - g^{ij} \frac{\partial A_{i \alpha}^\beta}{\partial x^j}
            + g^{ij} A_{k \alpha}^\beta \Gamma_{ij}^k
            - g^{ij} A_{i \alpha}^\gamma A_{j \gamma}^\beta
        \right)
        s^\alpha e_\beta.
    \end{align*}
    This is clearly $\Cinfty(M)$-linear and hence the local expression
    for an endomorphism field $B \in \Secinfty(\End(E))$. Since $D -
    \dAlembert^\nabla$ is globally defined, $B$ is indeed a globally
    defined section. Of course, taking the explicit but complicated
    transformation laws for coefficients of second order differential
    operators, connection one-forms and Christoffel symbols, this can
    also be checked by hand (though it is not very funny).
\end{proof}
\begin{remark}[Normally hyperbolic operators]
    \label{remark:NormallyHyperbolic}
    ~
    \begin{remarklist}
    \item \index{Normally hyperbolic operator!Leibniz rule}%
        \index{Covariant derivative!D-compatible@$D$-compatible}%
        \label{item:NormallyHyperbolicLeibniz}
        If $D$ is normally hyperbolic and $\nabla^E$ is the
        corresponding covariant derivative then $D$ satisfies the
        Leibniz rule
        \begin{equation}
            \label{eq:NormallyHyperbolicLeibniz}
            D (f \cdot s)
            = f D(s) + 2 \nabla^E_{\gradient(f)} s
            +(\dAlembert f) \cdot s
        \end{equation}
        for all $f \in \Cinfty(M)$ and $s \in \Secinfty(E)$. This
        follows from the above proof. The connection $\nabla^E$ is
        also called the \emph{$D$-compatible connection}. In the
        following, we can safely assume that $D$ is of the form
        $\dAlembert^\nabla + B$ as above.
    \item While in general every $B \in \Secinfty(\End(E))$ gives a
        normally hyperbolic $\dAlembert^\nabla + B$, in specific
        contexts there are sometimes more geometrically motivated
        choices for both, the connection $\nabla^E$ and the additional
        tensor field $B$.
    \item Even though we formulated the above proposition and the
        definition of normally hyperbolic differential operators with
        respect to a Lorentz signature, it is clear that the above
        considerations apply also to the general semi-Riemannian
        case. In the Riemannian case, the corresponding operators are
        called connection Laplacians and normally elliptic operators,
        respectively.
    \end{remarklist}
\end{remark}


%% file: causal.tex
%
%

While most of the material up to now was applicable for general
semi-Riemannian manifolds we shall now discuss the causal structure
referring to the Lorentz signature exclusively.

%
%

\subsection{Some Motivation from General Relativity}
\label{subsec:motivation-from-GR}

In general relativity\index{General relativity} the
spacetime\index{Spacetime} is described by a four-dimensional manifold
$M$ equipped with a Lorentz metric $g$ subject to Einstein's
equation. One defines the \emIndex{Einstein tensor}
\begin{equation}
    \label{eq:einstein-tensor}
    G = \Ric - \frac{1}{2} \scal \cdot g,
\end{equation}
which is a symmetric covariant tensor field
\begin{equation}
    \label{eq:einstein-tensor-field}
    G \in \Secinfty(\Sym^2 T^*M).
\end{equation}
It can be shown that the covariant divergence of $G$ vanishes,
\begin{equation}
    \label{eq:einstein-tensor-divless}
    \divergenz G = 0,
\end{equation}
while $G$ itself needs not to be covariant constant at all.
Physically, \eqref{eq:einstein-tensor-divless} is interpreted as a
\emph{conservation law}.  Einstein's equation is then given by
\begin{equation}
    \label{eq:einstein-equation}
    \index{Einstein equation}%
    G = \kappa T,
\end{equation}
where $T \in \Secinfty(\Sym^2 T^*M)$ is the so-called
\emph{energy-momentum tensor} of all matter and interaction fields on
$M$ \emph{excluding} gravity. The constant $\kappa$ is up to numerical
constants Newton's constant of gravity. The precise form of $T$ is
complicated and depends on the concrete realization of the matter
content of the spacetime under consideration.  More generally,
Einstein's equation with \emph{cosmological constant} are
\begin{equation}
    \label{eq:einstein-equation-with-cosmological-constant}
    \index{Cosmological constant}%
    G + \lambda g = \kappa T,
\end{equation}
where $\lambda \in \mathbb{R}$ is a constant, additional parameter of
the theory, the cosmological constant.

The nature of these equations is that for a given functional
expression for $T$ usually coming from a variational principle, the
metric $g$ has to be found in such a way that
\eqref{eq:einstein-equation} or
\eqref{eq:einstein-equation-with-cosmological-constant} is satisfied.
However, this is rather complicated as \eqref{eq:einstein-equation}
and \eqref{eq:einstein-equation-with-cosmological-constant} turn out
to be quadratic partial differential equations of second order in the
coefficients of the metric which are of a rather complicated type. On
one hand they are ``hyperbolic'' and therefor ask for an ``initial
value problem''. On the other hand, when formulating
\eqref{eq:einstein-equation} or
\eqref{eq:einstein-equation-with-cosmological-constant} as initial
value problem for a metric on a $3$-dimensional submanifold, the
Equations~\eqref{eq:einstein-equation} or
\eqref{eq:einstein-equation-with-cosmological-constant} have a certain
gauge freedom thanks to the diffeomorphism invariance of the condition
\eqref{eq:einstein-equation} and
\eqref{eq:einstein-equation-with-cosmological-constant}, respectively.
This yields ``constraints'' which have to be satisfied. For more
details on this initial value problem in general relativity see e.g.
\cite{choquet-bruhat.geroch:1969a, fischer.marsden:1979a,
  choquet-bruhat.york:1980a}.

All this makes general relativity quite complicated, both from the
conceptual and practical point of view. We refer to textbooks on
general relativity for a more detailed and sophisticated discussion,
see e.g.  \cite{sexl.urbantke:1987a, straumann:1988a,
  hawking.ellis:1973a, beem.ehrlich.easley:1996a}.

For $T = 0$ one speaks of a \emph{vacuum solution} to Einstein's
equation: only those degrees of freedom are relevant which come
directly from geometry and hence from gravity. Already this particular
case is very complicated as it is still highly non-linear.
Nevertheless, there are solutions which look like propagating waves or
black holes.

In the following, we take the point of view that a certain energy and
momentum content of the spacetime results in a certain metric $g$.
Then we assume that there is a slight perturbation by some additional
field $\phi$ on $M$ which on one hand has only a minor contribution to
$T$ and thus does not influence $g$. On the other hand, the field is
subject to field equations determined by $g$. With other words, we
neglect the \emph{back-reaction}\index{Back-reaction} of the field on
$g$ but investigate the field equations in a \emph{fixed} background
metric $g$.

Thus we arrive at field equations for $\phi$ on a given spacetime $(M,
g)$. It turns out that the question whether $g$ is a solution to
Einstein's equation or not, is of minor importance when we want to
understand the field equations for $\phi$. In fact, the geometric
features of $g$ which guarantee a ``good behaviour'' of $\phi$ are
rather independent of Einstein's equation.

In general, physically relevant field equations for $\phi$ can be
quite complicated: if we are interested in ``interacting fields'' then
the field equations are non-linear. Thus all the technology of
distributions etc. does \emph{not} apply, at least not in a naive way.
For this reason we restrict to \emph{linear} field equations: one
motivation is that even if the original field equations for $\phi$ are
non-linear, a linearization\index{Linearization} around a solution
$\phi_0$ might be interesting. Assuming that $\phi_0$ is a solution
one considers $\phi = \phi_0 + \psi$ and rewrites the (non-linear)
equations for $\phi$ as field equations for the perturbation $\psi$
and \emph{neglects} higher order terms in $\psi$. This way one obtains
an approximation in form of a linear field equation for $\psi$.

We shall now discuss some typical examples. The prototype of a field
equation is the \emIndex{Klein-Gordon equation} for a scalar field
$\phi \in \Cinfty(M)$ of mass $m \in \mathbb{R}$
\begin{equation}
    \label{eq:klein-gordon}
    \dAlembert \phi + m^2 \phi = 0.
\end{equation}
On non-trivial geometries there are physical arguments suggesting that
the Klein-Gordon equation should be modified in a way incorporating
the scalar curvature, i.e. one considers
\begin{equation}
    \label{eq:klein-gordon-with-curvature}
    \dAlembert \phi + \xi \scal \phi + m^2 \phi = 0,
\end{equation}
where $\xi \in \mathbb{R}$ is a parameter. While
\eqref{eq:klein-gordon-with-curvature} is still linear, a
self-interacting modification of the Klein-Gordon equation is e.g.
\begin{equation}
    \label{eq:klein-gordon-with-interaction}
    \dAlembert \phi + m^2 \phi + \lambda \phi^2 + \mu \phi^3 = 0,
\end{equation}
where again $\lambda, \mu \in \mathbb{R}$ are parameters of the
theory. If $\phi_0$ is a solution of
\eqref{eq:klein-gordon-with-interaction} then a linearized version of
\eqref{eq:klein-gordon-with-interaction} for $\phi = \phi_0 + \psi$ is
given by
\begin{equation}
    \label{eq:klein-gordon-interact-linearization}
    \dAlembert \psi + m^2 \psi + 2 \lambda \phi_0 \psi
    + 3 \phi_0^2 \psi
    = 0.
\end{equation}
By this procedure we obtain a rather general linear equation with
leading symbol being the metric but fairly general $\Cinfty(M)$-linear
part, in our case either $\xi \scal +m^2$ or $m^2 + 2 \lambda \phi_0 +
3 \mu \phi_0^2$ or a combination of both.

This motivates that one should consider linear second order
differential equations of normal hyperbolic type, i.e.
\begin{equation}
    \label{eq:2nd-order-diffop-normal-hyperbolic}
    (\dAlembert + B) \phi = 0,
\end{equation}
with $B \in \Cinfty(M)$. Finally, the step towards general vector
bundles and sections $\phi \in \Secinfty(E)$ is only a mild
generalization: in many physical field theories the fields have more
than one component. This way we arrive at field equations of the form
\begin{equation}
    \label{eq:2nd-order-diffop-normal-hyperbolic-for-sections}
    (\dAlembert^\nabla + B) \phi = 0
\end{equation}
for $\phi \in \Secinfty(E)$ with a connection d'Alembertian
$\dAlembert^\nabla$ and some $B \in \Secinfty(\End(E))$. Note once
more that in our approximation to general relativity we have a fixed
background metric $g$ used in the definition of $\dAlembert^\nabla$.

%
%

\subsection{Future and Past on a Lorentz Manifold}
\label{subsec:future-past-on-lorentz}

\begin{figure}
    \centering
    \input{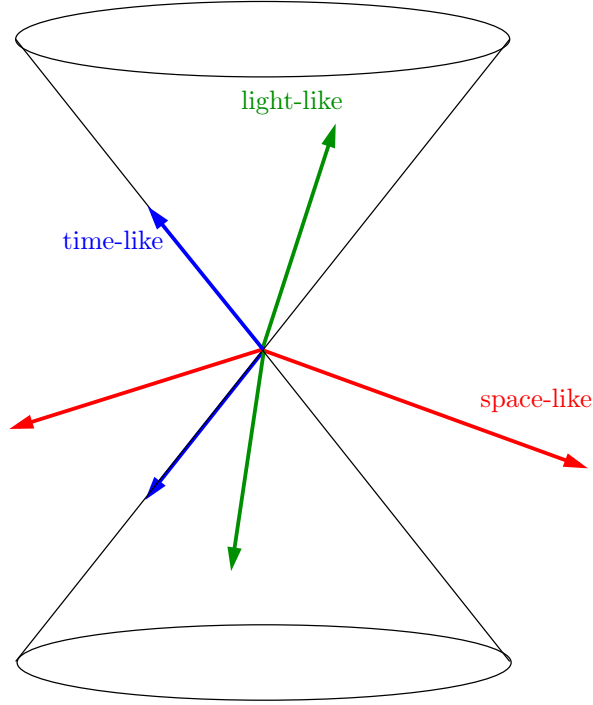}
    \caption{\label{fig:light-cone}%
      Light cone structure in Minkowski spacetime
    }
\end{figure}

Having a fixed Lorentz metric $g$ on a spacetime manifold $M$ we can
now transfer the notions of special relativity, see e.g.
\cite{roemer.forger:1993a}, to $(M, g)$. In fact, each tangent space
$(T_p M, g_p)$ is isometrically isomorphic to Minkowski spacetime
$(\mathbb{R}^n, \eta)$ with $\eta = \diag (+1, -1, \ldots, -1)$, by
choosing a Lorentz frame: there exist tangent vectors $e_i \in T_p M$
with $i = 1, \ldots, n$ such that
\begin{equation}
    \label{eq:lorentz-frame}
    g_p(e_i, e_j) = \eta_{ij} = \pm \delta_{ij}.
\end{equation}
\begin{remark}[Local Lorentz frame]
    \label{remark:local-lorentz-frame}
    \index{Frame!Lorentz}%
    The pointwise isometry from $(T_p M, g_p)$ to $(\mathbb{R}^n,
    \eta)$ can be made to depend smoothly on $p$ at least in a local
    neighborhood: For every $p \in M$ there exists a small open
    neighborhood $U$ of $p$ and local sections $e_1, \ldots, e_n \in
    \Secinfty(E \at{U})$ such that for all $q \in U$
    \begin{equation}
        \label{eq:local-lorentz-frame}
        g_q( e_i(q), e_j(q) ) = \eta_{ij}.
    \end{equation}
    In general, the frame $\{ e_i \}_{i=1, \ldots, n}$ can \emph{not}
    be chosen to come from a chart $x$ on $U$, i.e. $e_i$ is not
    $\frac{\partial}{\partial x^i}$. Here the curvature of $g$ is the
    obstruction. Nevertheless, such \emph{local Lorentz frames} will
    simplify certain computations. We note that for two local Lorenz
    frames $\{ e_i \}_{i=1, \ldots, n}$ and $\{\widetilde{e}_i
    \}_{i=1, \ldots, n}$ on $U$ there exists a unique smooth function
    $\Lambda: U \longrightarrow \group{O}(1, n-1)$ such that
    \begin{equation}
        \label{eq:lorentz-transformation-between-local-lorentz-frames}
        e_i(p) = \Lambda_i^j(p) \widetilde{e}_j(p),
    \end{equation}
    since the Lorentz transformations $\group{O}(1, n-1)$ are
    precisely the linear isometries of $(\mathbb{R}^n, \eta)$.
\end{remark}
As in special relativity one states the following definition:
\begin{definition}
    \label{definition:timelike-lightlike-spacelike}
    \index{Timelike vector}%
    \index{Spacelike vector}%
    \index{Lightlike vector}%
    \index{Causal vector}%
    Let $(M, g)$ be a Lorentz manifold and $v_p \in T_p M$ a non-zero
    vector. Then $v_p$ is called
    \begin{definitionlist}
    \item \label{item:timelike} timelike if $g_p (v_p, v_p) > 0$,
    \item \label{item:lightlike} lightlike or null if $g_p (v_p, v_p)
        = 0$,
    \item \label{item:spacelike} spacelike if $g_p (v_p, v_p) < 0$.
    \end{definitionlist}
    Non-zero vectors with $g_p(v_p, v_p) \geq 0$ are sometimes also
    called causal. To the zero vector, no attribute is assigned.
\end{definition}
In a fixed tangent space we have two open convex cones of timelike
vectors whose boundaries consists of the lightlike vectors together
with the zero vector, see Figure~\ref{fig:light-cone}.
Already in Minkowski spacetime there are Lorentz transformations
which exchange the two connected components of the timelike vectors.
Thus there is no intrinsic definition of ``future-'' and
``past-directed'' vectors in $(\mathbb{R}^n, \eta)$. Clearly, for
physical purposes it is crucial to have such a distinction: we choose
once and for all a \emph{time-orientation} on Minkowski spacetime
$(\mathbb{R}^n, \eta)$, i.e. a choice of one of the interiors of the
light-cones to be future directed. We symbolize this choice by
$(\mathbb{R}^n, \eta, \uparrow)$. Now only the \emph{orthochronous
  Lorentz transformations}
\begin{equation}
    \label{eq:orthocrone-lorentz-trans}
    \index{Minkowski spacetime!time-orientation}%
    \index{Lorentz transformation!orthochronous}%
    \group{L}^\uparrow (1, n-1) =
    \left\{
        \Lambda \in \group{O}(1, n,1)
        \; \big| \;
        \Lambda^0_0 > 0
    \right\}
\end{equation}
preserve the time-orientation $(\mathbb{R}^n, \eta,
\uparrow)$. Clearly, $\group{L}^\uparrow (1, n-1)$ is a closed
subgroup of $\group{O}(1, n-1)$ of the same dimension.

\index{Minkowski spacetime!space-orientation}%
\index{Lorentz transformation!proper}%
Analogously, Lorentz transformations do not preserve the
space-orientation in general. For a spacelike sub vector space $\Sigma
\subseteq \mathbb{R}^n$ (of dimension $n-1$), there are orientation
preserving and reversing Lorentz transformations. Choosing one
orientation of $\Sigma$ we obtain an additional structure on Minkowski
spacetime which we symbolize as $(\mathbb{R}^n, \eta, +)$ or
$(\mathbb{R}^n, \eta, \uparrow, +)$ in the case where we have chosen a
time-orientation as well. One can check that ``$+$'' does not depend
on the particular choice of $\Sigma$. The subgroups preserving $+$ or
$+$ and $\uparrow$ are the \emph{proper} and the \emph{proper and
  orthochronous} Lorentz transformations denoted by $\group{L}_+ (1,
n-1)$ and $\group{L}_+^\uparrow (1, n-1)$, respectively. It is a
standard fact that $\group{L}_+^\uparrow (1, n-1)$ is the connected
component of the identity and hence a normal subgroup. The discrete
resulting quotient group is
\begin{equation}
    \label{eq:lorentz-modulo-propor-orthocrone}
    \group{L}(1, n-1) / \group{L}_+^\uparrow (1, n-1)
    = \{ \id, \mathrm{P}, \mathrm{T}, \mathrm{PT} \}
\end{equation}
with relations $\mathrm{P}^2 = \mathrm{T}^2 = \id$ and $\mathrm{PT} =
\mathrm{TP}$. Then $\mathrm{T}$ is the
time-reversal\index{Time-reversal} while $\mathrm{P}$ is the
parity\index{Parity} operation.

We want to use now the time- and space-oriented Minkowski spacetime
$(\mathbb{R}^n ,\eta, \uparrow, +)$ in order to obtain time and space
orientations for $(M, g)$ as well. Here we meet the usual obstructions
analogously to the obstructions for orientability in general. In the
following the time-orientability will be crucial while the
space-orientability is not that important. Thus we focus on the
time-orientability. Here one has the following result:
\begin{proposition}
    \label{proposition:time-orientability}
    Let $(M, g)$ be a Lorentz manifold. Then the following statements
    are equivalent:
    \begin{propositionlist}
    \item \label{item:time-or-and-timelike-vf} There exists a
        timelike vector field $X \in \Secinfty(TM)$, i.e. $X(p)$ is
        timelike for all $p \in M$.
    \item \label{item:time-or-and-transition-matrices} There exists an
        open cover $\{U_\alpha\}_\alpha$ of $M$ with local Lorentz
        frames $\{e_{\alpha i}\}_{i = 1, \ldots, n} \in
        \Secinfty(TU_\alpha)$ such that on $U_\alpha \cap U_\beta \neq
        \emptyset$ the transition matrix $\Lambda_{\alpha\beta} \in
        \group{O}(1, n-1)$ with
        \begin{equation}
            \label{eq:transition-mat-for-lorentz-frames}
            e_{\alpha i} = \Lambda_{\alpha \beta} {}^j_i e_{\beta j}
        \end{equation}
        takes values in $\group{L}^\uparrow (1, n-1)$.
    \end{propositionlist}
\end{proposition}
\begin{proof}
    Assume that $X \in \Secinfty(TM)$ is timelike. Then we choose an
    open cover $\{U_\alpha\}$ of $M$ with locally defined Lorentz
    frames $\{ e_{\alpha i} \}$ on $U_\alpha$. Without restriction we
    can choose the $U_\alpha$ to be connected. Then on $U_\alpha$
    either the timelike vector $e_{\alpha 1}$ or the timelike vector
    $-e_{\alpha 1}$ is in the same connected component of the
    timelike vectors as $X$. Changing $e_{\alpha 1}$ to $-e_{\alpha
      1}$ if necessary yields a local Lorentz frame on $U_\alpha$ with
    $e_{\alpha 1}$ in the same connected component as $X$. Since this
    holds for all $\alpha$ we obtain transition matrices
    $\Lambda_{\alpha \beta}$ in $\Cinfty( U_\alpha \cap U_\beta,
    \group{L}^\uparrow)$ as wanted.

    Conversely, let such an open cover and local Lorentz frames be
    given. We choose a partition of unity $\chi_\alpha$ subordinate to
    $U_\alpha$ with $\chi_\alpha \geq 0$. Then we define
    \[
    X = \sum_\alpha \chi_\alpha e_{\alpha 1}
    \tag{$*$}
    \]
    which is clearly a globally defined smooth vector field $X \in
    \Secinfty(TM)$. At $p \in M$ only finitely many $\alpha$
    contribute to ($*$). Moreover, since by
    \eqref{eq:transition-mat-for-lorentz-frames} all the $e_{\alpha
      1}(p)$ are in the \emph{same} connected component of the
    timelike vectors and since this connected component is
    \emph{convex}, also $X(p)$ is in this connected component. It
    follows that $X(p)$ is timelike.
\end{proof}

There are still alternative formulations of the property described by
\refitem{item:time-or-and-timelike-vf} and
\refitem{item:time-or-and-transition-matrices} in
Proposition~\ref{proposition:time-orientability}. However, for the
time being we take the result of
Proposition~\ref{proposition:time-orientability} as definition of
time-orientability:
\begin{definition}[Time-orientability]
    \label{definition:time-orientability}
    \index{Time-orientable}%
    \index{Time-orientation}%
    \index{Future directed}%
    \index{Past directed}%
    Let $(M, g)$ be a Lorentz manifold.
    \begin{definitionlist}
    \item \label{item:time-orientable} $(M, g)$ is called
        time-orientable if there exists a timelike vector field $X
        \in \Secinfty(TM)$.
    \item \label{item:time-orientation} The choice of a timelike
        vector field $X \in \Secinfty(TM)$ is called a
        time-orientation.
    \item \label{item:future-and-past-directed} With respect to a
        time-orientation, a timelike vector $v_p \in T_p M$ is called
        future directed if $v_p$ is in the same connected component as
        $X(p)$. It is called past directed if $-v_p$ is future
        directed.
    \end{definitionlist}
\end{definition}
\begin{remark}[Time-orientability]
    \label{remark:time-orientability}
    Note that time-orientability of $(M, g)$ is rather independent of
    (topological) orientability of $M$. One can find easily a Lorentz
    metric on the Möbius strip which is time-orientable and,
    conversely, a Lorentz metric on the cylinder $\mathbb{S}^1 \times
    \mathbb{R}$ which is \emph{not} time-orientable. We leave it as an
    exercise to figure out the details of these examples.
\end{remark}

\begin{figure}
    \centering
    \input{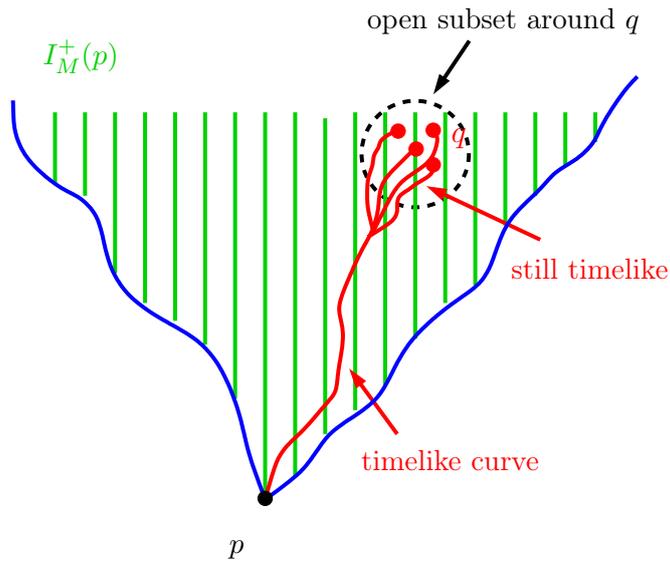}
    \caption{\label{fig:chronological-future-is-open}%
      The chronological future is open.
    }
\end{figure}

\begin{figure}
    \centering
    \input{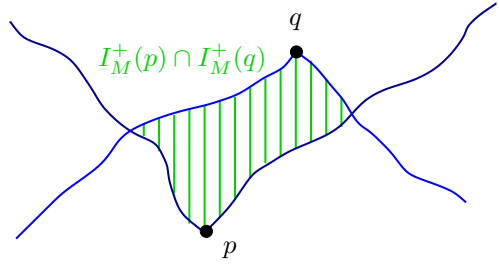}
    \caption{\label{fig:open-diamond}%
      An chronological open diamond in a spacetime.
    }
\end{figure}

\begin{figure}
    \centering
    \input{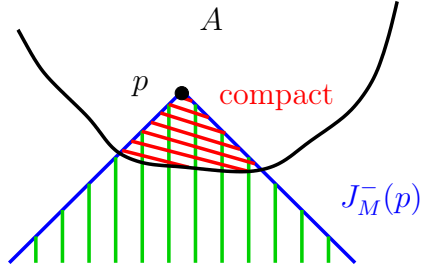}
    \caption{\label{fig:past-compact-in-Minkowski}%
      A past compact subset $A$ in Minkowski spacetime.
    }
\end{figure}

In the following, we shall always assume that $(M, g)$ is time
orientable. Moreover, we assume that a time-orientation has been
chosen once and for all. This will be important for a consistent
interpretation of $(M, g)$ as a spacetime manifold. Using the time
orientation we can define the future and past of a given point in $M$.
More precisely, one calls a curve $\gamma: I \subseteq \mathbb{R}
\longrightarrow M$ \emph{timelike}, \emph{lightlike}, \emph{spacelike}
or \emph{causal} if $\dot{\gamma}(t)$ is timelike, lightlike,
spacelike, or causal for all $t \in I$, respectively. A causal vector
$v_p \in T_p M$ is called future or past directed if it is contained
in the closure of the future or past directed timelike vectors at
$p$. Then a curve $\gamma$ is called future or past directed if
$\dot{\gamma}(t)$ is causal and future or past directed at every
$t$. By continuity we see that a causal curve is either future or past
directed. In a time-oriented spacetime it cannot change the causal
direction. Clearly a $\Fun[1]$-curve is sufficient for this argument.
\begin{definition}
    \label{definition:past-relations}
    \index{Future directed}%
    \index{Past directed}%
    Let $(M, g)$ be a time-oriented Lorentz manifold and $p, q \in
    M$. The we define
    \begin{definitionlist}
    \item \label{item:strict-earlier-relation} $p \ll q$ if there
        exists a future directed, timelike smooth curve from $p$ to
        $q$.
    \item \label{item:earlier-or-equal-relation} $p \leq q$ if either
        $p = q$ or there exists a future directed, causal smooth curve
        from $p$ to $q$.
    \item \label{item:earlier-relation} $p < q$ if $p \leq q$ but $p
        \neq q$.
    \end{definitionlist}
\end{definition}
Clearly the relations $\ll$ and $\leq$ are transitive. We use these
relations to define the chronological and causal future and past of a
point:
\begin{definition}[Chronological and causal future and past]
    \label{definition:chron-causal-future-and-past}
    \index{Chronological future}%
    \index{Chronological past}%
    \index{Causal future}%
    \index{Causal past}%
    Let $(M, g)$ be a time-oriented Lorentz manifold and $p \in M$.
    \begin{definitionlist}
    \item \label{item:chronological-future} The chronological future
        of $p$ is
        \begin{equation}
            \label{eq:chronological-future}
            I^+(p) = \left\{ q \in M \; \big| \; p \ll q \right\}.
        \end{equation}
    \item \label{item:chronological-past} The chronological past of
        $p$ is
        \begin{equation}
            \label{eq:chronological-past}
            I^-(p) = \left\{ q \in M \; \big| \; q \ll p \right\}.
        \end{equation}
    \item \label{item:causal-future} The causal future of $p$ is
        \begin{equation}
            \label{eq:causal-future}
            J^+(p) = \left\{ q \in M \; \big| \; p \leq q \right\}.
        \end{equation}
    \item \label{item:causal-past} The causal past of $p$ is
        \begin{equation}
            \label{eq:causal-past}
            J^-(p) = \left\{ q \in M \; \big| \; q \leq p \right\}.
        \end{equation}
    \end{definitionlist}
\end{definition}
Sometimes we indicate the ambient spacetime $M$ in the definitions
by $I_M^\pm (p)$ and $J_M^\pm (p)$ since they will play a crucial
role. The definitions of $I_M^\pm (p)$ and $J_M^\pm (p)$ reflect
\emph{global} properties of $M$ which are not necessarily preserved
under isometric embeddings. We illustrate the meaning of $I_M^\pm (p)$
and $J_M^\pm (p)$ by some examples:
\begin{example}
    \label{example:past-and-future-for-subsets-of-Minkowski-space}
    The spacetime $(M, g)$ in
    Figure~\ref{fig:future-and-past-for-convex-spacetime} and the
    following pictures are open subsets of the usual Minkowski
    spacetime $(\mathbb{R}^2, \eta)$ with future direction being
    ``upward''.
    \begin{figure}
        \centering
        \input{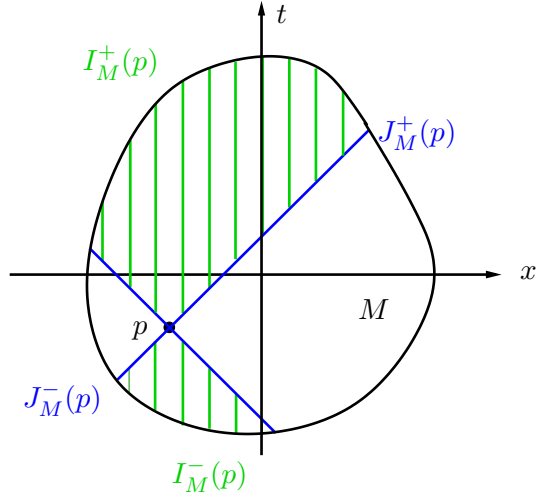}
        \caption{\label{fig:future-and-past-for-convex-spacetime}%
          Future and past for a convex subset of Minkowski spacetime.
        }
    \end{figure}
    \begin{figure}
        \centering
        \input{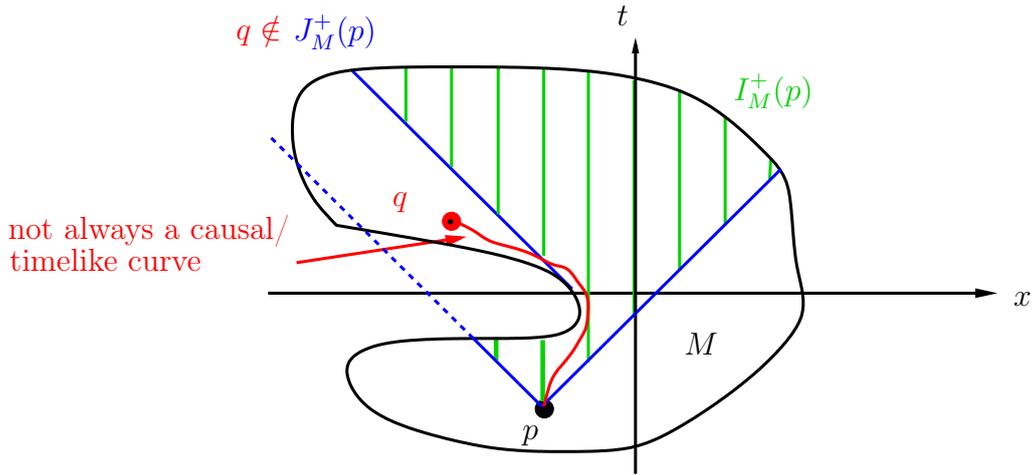}
        \caption{\label{fig:future-and-past-for-nonconvex-spacetime}%
          Future and past for a spacetime $M$ with ``notch''.
        }
    \end{figure}
    \begin{figure}
        \centering
        \input{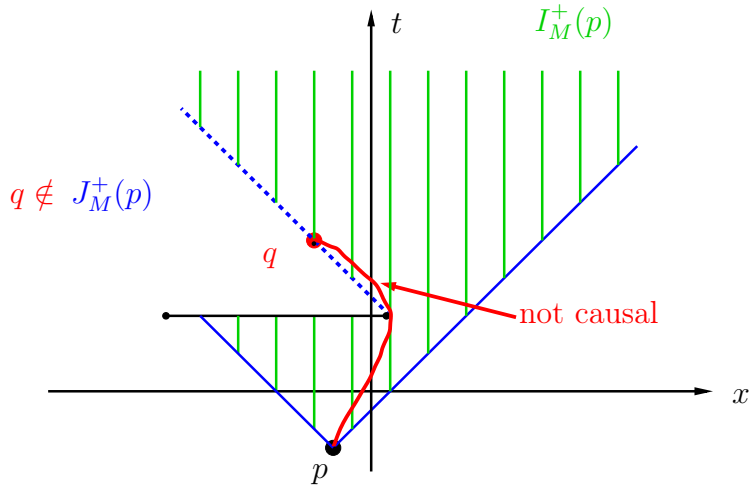}
        \caption{\label{fig:future-for-spacetime-without-line-segment}%
          Future and past for a spacetime with an excluded line
          segment.
        }
    \end{figure}
    Figure~\ref{fig:future-and-past-for-nonconvex-spacetime} shows
    that $I_M^+ (p)$ and $J_M^+ (p)$ are not just the intersections of
    $M$ with $I_{\mathbb{R}^2}^+ (p)$ and $J_{\mathbb{R}^2}^+ (p)$,
    but actually smaller.
    Figure~\ref{fig:future-for-spacetime-without-line-segment}
    illustrates that $J_M^+(p)$ needs not to be the closure of
    $I_M^+(p)$. In fact, $J_M^+(p)$ is not closed at all in this
    example.
\end{example}
Without proof we state the following result, see e.g.
\cite[Chap.~14]{oneill:1983a}:
\begin{proposition}
    \label{proposition:chronological-future-past-is-open}
    Let $(M, g)$ be a time-oriented Lorentz manifold. Then for every
    $p \in M$ the chronological future and past $I_M^\pm (p)$ of $p$
    is an open subset of $M$.
\end{proposition}
The intuition behind this proposition is clear and is visualized
in Figure~\ref{fig:chronological-future-is-open}.
Since the sets $I_M^\pm (p)$ are open, we can use them to define a
collection of open subsets of $M$. In particular, we consider the
intersections $I_M^+(p) \cap I_M^-(q)$ for $p, q \in M$. These subsets
are sometimes called (chronological) open
\emph{diamonds}\index{Diamond!open} as Figure~\ref{fig:open-diamond}
suggests. In flat Minkowski space the sets $I_M^+(p) \cap I_M^-(q)$
are diamond-shaped.
These open diamonds can be used to define a \emph{new} topology on
$M$: they form a basis of a topology sometimes called the
\emIndex{Alexandrov topology} of $(M, g)$. By
Proposition~\ref{proposition:chronological-future-past-is-open} it is
coarser than the original topology. We will come back to the question
whether the Alexandrov topology actually coincides with the usual one;
a case which is of course physically interesting: in this case the
topological structure of $M$ is determined by the causal structure.
Analogously to the chronological open diamonds, we define the diamonds
\begin{equation}
    \label{eq:diamonds}
    \index{Diamond!causal}%
    J_M(p, q) = J_M^+(p) \cap J_M^-(q).
\end{equation}
Finally, we can extend
Definition~\ref{definition:chron-causal-future-and-past} to arbitrary
subsets $A \subseteq M$. One defines the chronological future and past
as well as the causal future and past of $A$ by
\begin{equation}
    \label{eq:chronological-future-past-of-subset}
    I_M^\pm (A) = \bigcup_{p \in A} I_M^\pm (p)
\end{equation}
and
\begin{equation}
    \label{eq:causal-future-past-of-subset}
    J_M^\pm (A) = \bigcup_{p \in A} J_M^\pm (p),
\end{equation}
respectively. Again, $J_M^\pm (A)$ needs not to be closed but is
contained in the closure of $I_M^\pm (A)$ which is always open by
Proposition~\ref{proposition:chronological-future-past-is-open}.
\begin{figure}
    \centering
    \input{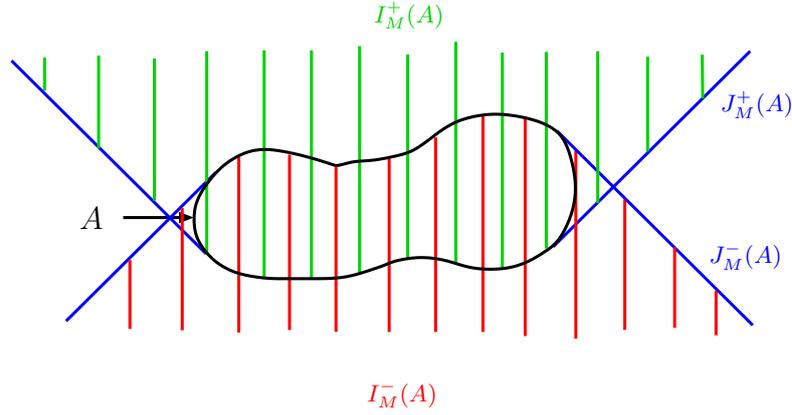}
    \caption{\label{fig:chron-caus-future-past-of-subset-in-Minkowski-space}%
      Chronological and causal future and past of $A$ in Minkowski
      spacetime $(\mathbb{R}^2, \eta)$.
    }
\end{figure}
\begin{definition}[Future and past compactness]
    \label{definition:future-past-compactness}
    \index{Future compact}%
    \index{Past compact}%
    Let $(M, g)$ be a time-oriented Lorentz mani\-fold. Then a subset
    $A \subseteq M$ is called future compact if $J_M^+(p) \cap A$ is
    compact for all $p \in M$ and past compact if $J_M^-(p) \cap A$ is
    compact for all $p \in M$.
\end{definition}
The geometric interpretation is clear and can be visualized again in
Minkowski spacetime as in Figure~\ref{fig:past-compact-in-Minkowski}.
Clearly, $A$ needs not be compact in the topological sense. However,
if all the $J_M^\pm (p)$ are \emph{closed} then every compact subset
$A \subseteq M$ is future and past compact.

The phenomenon in
Figure~\ref{fig:future-and-past-for-nonconvex-spacetime} motivates the
following definition:
\begin{definition}[Causal compatibility]
    \label{definition:causal-compatibility}
    \index{Causal compatibility}%
    Let $(M, g)$ be a time-oriented Lorentz manifold and $U \subseteq
    M$ open. Then $U$ is called causally compatible if for all $p \in
    M$ we have
    \begin{equation}
        \label{eq:causal-compatible}
        J_U^\pm(p) = J_M^\pm (p) \cap U.
    \end{equation}
    More generally, a time-orientation preserving isometric embedding
    $\iota: (N, h) \hookrightarrow (M, g)$ of a time-oriented Lorentz
    manifold $(N, h)$ into $(M, g)$ is called causally compatible if
    $\iota(N) \subseteq M$ is causally compatible.
\end{definition}
\begin{remark}
    \label{remark:causal-compatible}
    Let $(M, g)$ be a time-oriented Lorentz manifold.
    \begin{remarklist}
    \item \label{item:causal-compatible-and-causal-curves} $U
        \subseteq M$ is causally compatible if for every causal curve
        from $p \in U$ to $q \in U$ in $M$ one also finds a causal
        curve from $p$ to $q$ which lies entirely in $U$. In
        Figure~\ref{fig:future-and-past-for-nonconvex-spacetime} this
        is not the case for the subset $M \subseteq \mathbb{R}^2$.
    \item \label{item:causal-compatible-is-transitive} If $V \subseteq
        U \subseteq M$ are open subset such that $V \subseteq U$ is
        causally compatible in the Lorentz manifold $(U, g\at{U})$ and
        $U$ is causally compatible in $M$, then also $V \subseteq M$
        is causally compatible.
    \item \label{item:causal-future-past-of-subset-of-caus-comp-set}
        If $U \subseteq M$ is causally compatible and $A \subseteq U$
        the clearly
        \begin{equation}
            \label{eq:causal-future-past-of-subset-of-caus-comp-set}
            J_U^\pm (A) = J_M^\pm (A) \cap U.
        \end{equation}
    \item \label{item:lorentz-category} Since the relation ``causally
        compatible'' is transitive with respect to inclusion, we
        obtain a category of $n$-dimensional time-oriented Lorentz
        manifolds $\mathsf{Lorentz}_n$ as follows: the objects will be
        $n$-dimensional time-oriented Lorentz manifolds and the
        morphisms $\iota: (N, h) \hookrightarrow (M, g)$ will be
        isometric embeddings preserving the time-orientations which
        are causally compatible. Even though there are usually not
        many morphisms between two objects in this category, it will
        turn out to be a very useful notion. In recent approaches to
        axiomatic quantum field theory on generic spacetimes this
        point of view becomes important, see e.g.
        \cite{brunetti.fredenhagen.verch:2003a, hollands.wald:2010a}
        and references therein.
    \end{remarklist}
\end{remark}

%
%

\subsection{Causality Conditions and Cauchy-Hypersurfaces}
\label{subsec:caus-cond-cauchy-hypersurfaces}

\begin{figure}
    \centering
    \input{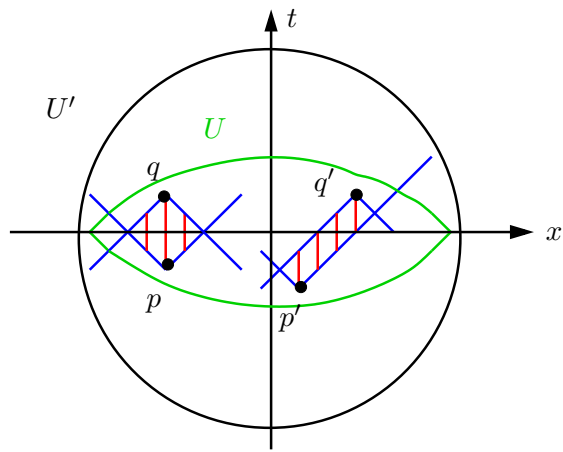}
    \caption{\label{fig:convex-and-causal}%
      A subset $U$ which is both convex and causal.
    }
\end{figure}

\begin{figure}
    \centering
    \input{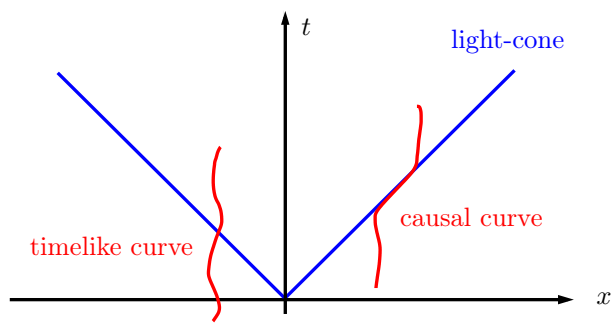}
    \caption{\label{fig:lightcone-achronal-but-not-acausal}%
      The light cones in Minkowski spacetime are achronal but not
      acausal.
    }
\end{figure}

We continue our investigation of the causality structure of a
time-oriented Lorentz manifold $(M, g)$. We start with the following
definition:
\begin{definition}[Causal subsets]
    \label{definition:causal-subsets}
    \index{Causal subset}%
    \index{Geodesically convex}%
    Let $U \subseteq M$ be an open subset. Then $U$ is called causal
    if there is a geodesically convex open subset $U' \subseteq M$
    such that $U^\cl \subseteq U'$ and for any two points $p, q \in
    U^\cl$ the diamond $J_{U'}(p, q)$ is compact and contained in
    $U^\cl$.
\end{definition}
Figure~\ref{fig:convex-but-not-causal} to
Figure~\ref{fig:convex-and-causal} show the relations between the
notions of geodesically convex and causal subsets. Again, the ambient
spacetime is the Minkowski spacetime $(\mathbb{R}^2, \eta)$. Since the
geodesics are still the straight lines, open convex subsets $U'
\subseteq \mathbb{R}^2$ in the usual sense coincide with the
geodesically convex subsets.
\begin{figure}
    \centering
    \input{convex-not-causal.\pictype}
    \caption{\label{fig:convex-but-not-causal}%
      A subset $U$ which is convex but not causal.
    }
\end{figure}
\begin{figure}
    \centering
    \input{causal-not-convex.\pictype}
    \caption{\label{fig:causal-but-not-convex}%
      A subset $U$ which is causal but not convex.
    }
\end{figure}
The ``opposite'' of a causal domain are the acausal subsets of $M$.
\begin{definition}[Acausal and achronal subsets]
    \label{definition:acausal-and-achronal}
    \index{Acausal subset}%
    \index{Achronal subset}%
    Let $A \subseteq M$ be a subset of a time-oriented Lorentz
    manifold. Then $A$ is called
    \begin{definitionlist}
    \item \label{item:achronal} achronal if every timelike curve
        intersects $A$ in at most one point.
    \item \label{item:acausal} acausal if every causal curve
        intersects $A$ in at most one point.
    \end{definitionlist}
\end{definition}
Clearly, acausal subsets are achronal but the reverse is not true.
Already the light cones in Minkowski spacetime are achronal but not
acausal, as Figure~\ref{fig:lightcone-achronal-but-not-acausal}
illustrates.
Using the causal structure of $(M, g)$ we obtain a refined notion of
boundary and closure of a subset $A \subseteq M$. One defines $p \in
A^\cl$ to be an \emph{edge point}\index{Edge point} if for all open
neighborhoods $U$ of $p$ there exists a timelike curve from
$I_U^-(p)$ to $I_U^+(p)$ which does not meet $A$.
\begin{figure}
    \centering
    \input{edge-points-2-dim.\pictype}
    \caption{
      \label{fig:edge-points-2-dim}%
      Examples of edge points of a line segment in $2$-dimensional
      Minkowski spacetime.
    }
\end{figure}
In Figure~\ref{fig:edge-points-2-dim} the point $q$ is an edge point
of the segment while $p$ is not. In
Figure~\ref{fig:edge-points-3-dim}, the line segment $A$ is considered
as subset of $3$-dimensional Minkowski spacetime $(\mathbb{R}^3,
\eta)$. Then all points in $A^\cl$ are edge points. Thus the notion of
edge points is finer than the notion of a (topological) boundary
point.
\begin{figure}
    \centering
    \input{edge-points-3-dim.\pictype}
    \caption{\label{fig:edge-points-3-dim}%
      Examples of edge points of a line segment $A$ in $3$-dimensional
      Minkowski spacetime.
    }
\end{figure}
We want to get as large achronal or acausal subsets as possible: they
will be good candidates for Cauchy hypersurfaces where we can impose
initial conditions. The following theorem states that we can expect at
least $\Fun[0]$-submanifolds.
\begin{theorem}[Achronal hypersurfaces]
    \label{theorem:achronal-hypersurfaces}
    \index{Achronal hypersurface}%
    Let $(M, g)$ be a time-oriented Lorentz manifold and $A \subseteq
    M$ achronal. Then $A$ is a topological hypersurface in $M$ if and
    only if $A$ does not contain any of its edge points.
\end{theorem}
Recall that a topological hypersurface $\Sigma$ of $M$ is a
$\Fun[0]$-manifold $\Sigma$ together with a $\Fun[0]$-embedding $i:
\Sigma \hookrightarrow M$ with codimension one. In general, we can not
expect more than a $\Fun[0]$-hypersurface as the example of the light
cone shows. For a proof we refer to \cite[Prop.~24 in
Chap~14]{oneill:1983a}. The following corollary is a straightforward
consequence of Theorem~\ref{theorem:achronal-hypersurfaces}.
\begin{corollary}
    \label{corollary:edgeless-achronal-subset}
    An achronal subset $A$ is a closed topological hypersurface if and
    only if $A$ is edgeless.
\end{corollary}
The extreme case of an achronal hypersurface will be a Cauchy
hypersurface. First recall that a timelike curve $\gamma: I \subseteq
\mathbb{R} \longrightarrow M$ is called
\emph{inextensible}\index{Inextensible} if there is no
``reparametrization'' $\widetilde{\gamma}: J \subset \mathbb{R}
\longrightarrow M$ of $\gamma$ such that $\widetilde{\gamma}(J)
\supseteq \gamma(I)$ is strictly larger. Then we can formulate the
following definition:
\begin{definition}[Cauchy hypersurface]
    \label{definition:cauchy-hypersurface}
    \index{Cauchy hypersurface}%
    Let $(M, g)$ be a time-oriented Lorentz manifold. A subset $\Sigma
    \subseteq M$ is called a Cauchy hypersurface if every inextensible
    timelike curve meets $\Sigma$ in exactly one point.
\end{definition}
\begin{figure}
    \centering
    \input{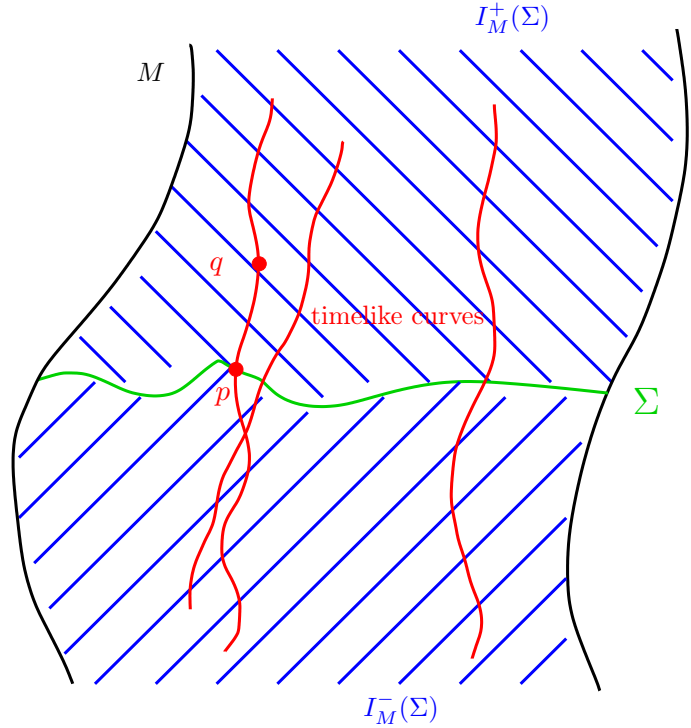}
    \caption{\label{fig:cauchy-hypersurface}%
      A Cauchy hypersurface $\Sigma$ in a spacetime $M$.
    }
\end{figure}
\begin{remark}[Cauchy hypersurface]
    \label{remark:cauchy-hypersurface}
    Clearly, a Cauchy hypersurface $\Sigma$ is achronal. Moreover, by
    the very definition of an edge point, $\Sigma$ has no edge points.
    Thus $\Sigma$ is a closed topological hypersurface by
    Theorem~\ref{theorem:achronal-hypersurfaces}. Finally, if $q \in
    M$ there exists a timelike curve through $q$, say a timelike
    geodesic. Thus such a timelike curve has an extension which meets
    $\Sigma$ in one point $p \in \Sigma$. It follows that either $q
    \ll p$, $p = q$, or $p \ll q$. Thus $M$ is the disjoint union of
    the non-empty open subsets $I_M^\pm (\Sigma)$ and $\Sigma$. Hence
    $\Sigma$ is the topological boundary of $I_M^\pm (\Sigma)$, i.e.
    we have the disjoint union
    \begin{equation}
        \label{eq:cauchy-hypersurface-separates-M}
        M = I_M^+ (\Sigma)
        \; \dot{\cup} \;
        \Sigma
        \; \dot{\cup} \;
        I_M^-(\Sigma).
    \end{equation}
    Furthermore, on can show that a Cauchy hypersurface is met by
    every inextensible causal curve at least once, see
    e.g.~\cite[Lem.~29 in Chap.~14]{oneill:1983a}.
\end{remark}
The physical interpretation of a Cauchy hypersurface is that the whole
future of the spacetime, viewed from $\Sigma$ is
predictable\index{Predictable} in the sense that every particle or
light ray being in the future $I_M^+(\Sigma)$ of $\Sigma$ has passed
through $\Sigma$ at earlier times. Analogously, viewed from $\Sigma$,
the whole past of $M$ is already known.

For an arbitrary subset $A \subseteq M$ we can still ask which part of
$M$ is predictable from $A$. This motivates the following definition
of the Cauchy development of $A$:
\begin{definition}[Cauchy development]
    \label{definition:cauchy-development}
    \index{Cauchy development}%
    Let $A \subseteq M$ be a subset. The future Cauchy development
    $D_M^+(A) \subseteq M$ of $A$ is the set of all those points $p
    \in M$ for which every past-inextensible causal curve through $p$
    also meets $A$. Analogously, one defines the past Cauchy
    development $D_M^-(A)$ and we call
    \begin{equation}
        \label{eq:cauchy-development}
        D_M(A) = D_M^+(A) \cup D_M^-(A)
    \end{equation}
    the Cauchy development of $A$.
\end{definition}
\begin{remark}[Cauchy development]
    \label{remark:cauchy-development}
    Let $A \subseteq M$ be a subset. The physical interpretation of
    $D_M^+(A)$ is that $D_M^+(A)$ is predictable from
    $A$. Analogously, $D_M^-(A)$ consists of those points which
    certainly influence $A$ in their future. We have $A \subseteq
    D_M^\pm (A)$.
\end{remark}
\begin{figure}
    \centering
    \input{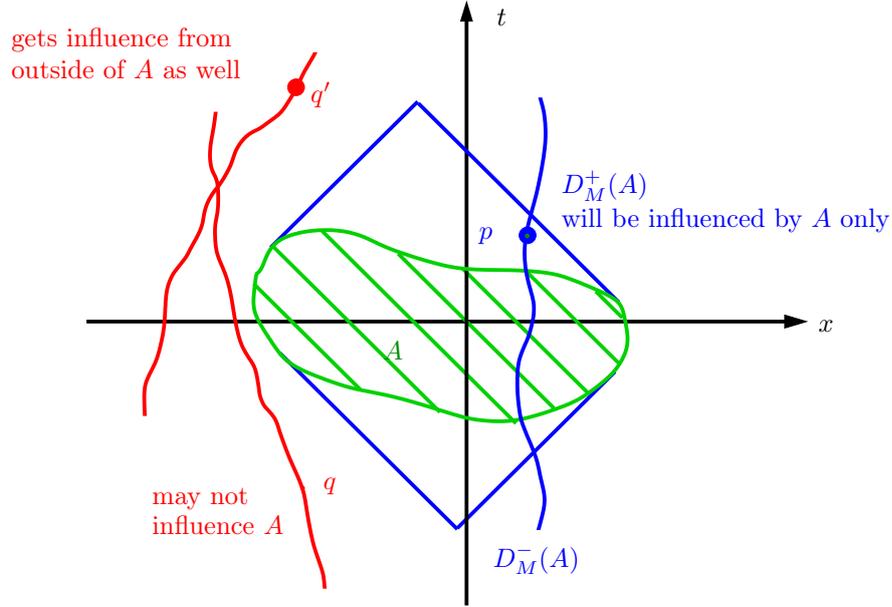}
    \caption{\label{fig:cauchy-development}%
      Cauchy development of a subset $A$ in the 2-dimensional
      Minkowski spacetime.
    }
\end{figure}
\begin{remark}
    \label{remark:properties-of-cauchy-development-operations}
    For $A \subseteq M$ we clearly have
    \begin{equation}
        \label{eq:8}
        D_M^\pm ( D_M^\pm(A) ) = D_M^\pm (A)
    \end{equation}
    and hence
    \begin{equation}
        \label{eq:9}
        D_M ( D_M(A) ) = D_M (A).
    \end{equation}
    Moreover, for $A \subseteq B \subseteq M$ we have
    \begin{equation}
        \label{eq:10}
        D_M^\pm (A) \subseteq D_M^\pm (B)
    \end{equation}
    and
    \begin{equation}
        \label{eq:11}
        D_M(A) \subseteq D_M(B).
    \end{equation}
    Thus the three operations $D_M^\pm(\argument)$ and
    $D_M(\argument)$ behave similar as the topological closure $A
    \mapsto A^\cl$.
\end{remark}
\begin{remark}
    \label{remark:4}
    If $A \subseteq M$ is achronal then $A$ is a Cauchy hypersurface
    if and only if $D_M(A) = M$. Thus for an achronal hypersurface,
    $D_M(A)$ can be viewed as the largest subset of $M$ for which $A$
    is a Cauchy hypersurface. In fact, one can show that $D_M(A)$ is
    open for an acausal topological hypersurface, see
    e.g.~\cite[Lem.~43 in Chap.~14]{oneill:1983a}.
\end{remark}
While the existence of a Cauchy hypersurface is from the physical
point of view very appealing, it is by far not evident. In fact, not
every time-oriented Lorentz manifold has a Cauchy hypersurface. Quite
contrary to the existence of a Cauchy hypersurface is the following
example:
\begin{example}
    \label{example:no-cauchy-hypersurface}
    \index{Geodesic!periodic}%
    We consider the cylinder $M = \mathbb{S}^1 \times \mathbb{R}$ with
    Lorenz metric $\D t^2 - \D x^2$ where the time variable is in
    $\mathbb{S}^1$-direction. The global vector field
    $\frac{\partial}{\partial t}$ is timelike and defines the
    time-orientation. Then through every point $p \in M$ there is a
    timelike geodesic which is \emph{periodic}. Thus there cannot be
    any Cauchy
    hypersurface. Figure~\ref{fig:periodic-geodesics-no-cauchy-hypersurface}
    illustrates this situation.
    \begin{figure}
        \centering
        \input{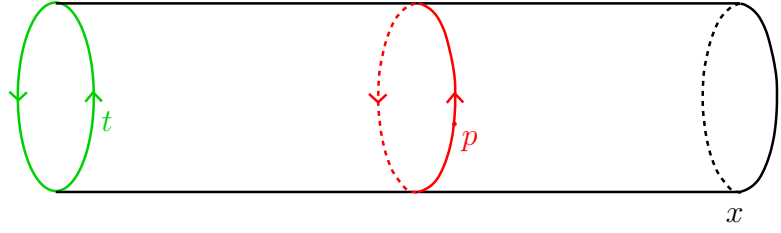}
        \caption{\label{fig:periodic-geodesics-no-cauchy-hypersurface}%
          Periodic timelike geodesic on a cylinder.
        }
    \end{figure}
    A slight variation is obtained by removing two lines in
    Figure~\ref{fig:almost-periodic-geodesic}.
    \begin{figure}
        \centering
        \input{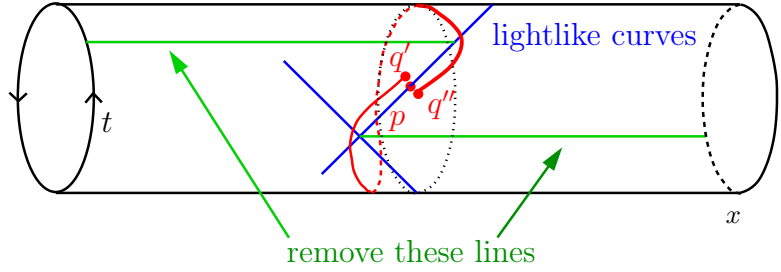}
        \caption{\label{fig:almost-periodic-geodesic}%
          Almost periodic timelike curves on a cylinder with removed
          line segments.
        }
    \end{figure}
    Then there are no longer closed timelike curves. However, starting
    arbitrarily close to the point $p$ at $q'$ there is a timelike
    curve (no longer geodesic of course) which ends again arbitrarily
    close to $p$ in $q''$.
\end{example}
Both situations are of course very bad for physical interpretations:
in the first case one could travel into ones own past with all the
funny paradoxa appearing. In the second case one could do so at least
approximately. This motivates the following definition:
\begin{definition}[Causality condition]
    \label{definition:causality-condition}
    \index{Causality condition}%
    \index{Causally convex}%
    \index{Spacetime!strongly causal}%
    \index{Spacetime!causal}%
    \index{Causal loop}%
    Let $(M, g)$ be a time-oriented Lorentz manifold.
    \begin{definitionlist}
    \item \label{item:causal} $M$ is called causal if there are no closed
        causal curves in $M$.
    \item \label{item:causally-convex} An open subset $U \subseteq M$
        is called causally convex if no causal curve intersects with
        $U$ in a disconnected subset of $U$.
    \item \label{item:strongly-causal-at-p} $M$ is called strongly
        causal at $p \in M$ if every open neighborhood of $p$
        contains an open causally convex neighborhood.
    \item \label{item:strongly-causal} $M$ is called strongly causal
        if $M$ is strongly causal at every point $p \in M$.
    \end{definitionlist}
\end{definition}
Without proof we mention the following interpretation of the strong
causality condition, see
e.g.~\cite[Prop.~3.11]{beem.ehrlich.easley:1996a}:
\begin{theorem}[Kronheimer, Penrose]
    \label{theorem:strongly-causal-and-Alexandrov-topology}
    \index{Alexandrov topology}%
    A time-oriented Lorentz manifold $(M, g)$ is strongly causal if
    and only if the Alexandrov topology coincides with the original
    topology of $M$.
\end{theorem}

The last ingredient we need is the following: In
Example~\ref{example:past-and-future-for-subsets-of-Minkowski-space}
we have seen examples of time-oriented spacetimes where the sets
$J_M^\pm (p)$ are not closed and hence not the closure of the $I_M^\pm
(p)$. To cure this effect one demands that the diamonds $J_M(p,q) =
J_M^+(p) \cap J_M^-(q)$ are \emph{compact} for all $p,q \in M$. Here
one has the following nice consequence, see
e.g. \cite{minguzzi.sanchez:2006a:pre}:
\begin{proposition}
    \label{proposition:compact-diamonds}
    \index{Diamond!compact}%
    Assume that $J_M(p,q) = J_M^+(p) \cap J_M^-(q)$ is compact for all
    $p,q \in M$ on a time-oriented spacetime $(M, g)$. Then the causal
    past and future $J_M^\pm (p)$ of any point $p \in M$ are closed
    subsets of $M$.
\end{proposition}
\begin{remark}
    \label{remark:1}
    In this section we only introduced some of the characteristic
    features of a time-oriented Lorentz manifold. There are many other
    notions of causality with increasing strength. Remarkably, many
    fundamental insights have been obtained only recently. We refer to
    the very nice review article of Minguzzi and S\'anchez
    \cite{minguzzi.sanchez:2006a:pre} for an additional discussion.
\end{remark}

%
%

\subsection{Globally Hyperbolic Spacetimes}
\label{subsec:glob-hyperb-spacetimes}

We are now in the position to define a globally hyperbolic spacetime
according to \cite{bernal.sanchez:2007a}:
\begin{definition}[Globally hyperbolic spacetime]
    \label{definition:globally-hyperbolic-spacetime}
    \index{Globally hyperbolic}%
    A time-oriented Lorentz manifold $(M, g)$ is called globally
    hyperbolic if
    \begin{definitionlist}
    \item \label{item:gl-hyp-is-causal} $(M, g)$ is causal,
    \item \label{item:gl-hyp-diamonds-are-compact} all diamonds
        $J_M(p, q)$ are compact for $p, q \in M$.
    \end{definitionlist}
\end{definition}
Note that in earlier works the notion of globally hyperbolic
spacetimes involved a strongly causal $(M, g)$ instead of just a
causal one. It was observed only recently that these two notions
actually coincide, see \cite{bernal.sanchez:2007a}.

The relevance of this condition comes from the relation to Cauchy
hypersurfaces. To this end, we first introduce the notion of a time
function:
\begin{definition}[Time function]
    \label{definition:time-function}
    \index{Time function}%
    \index{Temporal function}%
    \index{Cauchy hypersurface}%
    Let $(M, g)$ be a time-oriented Lorentz manifold and $t: M
    \longrightarrow \mathbb{R}$ a continuous function. Then $t$ is
    called a
    \begin{definitionlist}
    \item \label{item:time-function} time function if $t$ is strictly
        increasing along all future directed causal curves.
    \item \label{item:temporal-function} temporal function if $t$ is
        smooth and $\gradient t$ is future directed and timelike.
    \item \label{item:cauchy-time-function} Cauchy time function if
        $t$ is a time function whose level sets are Cauchy
        hypersurfaces.
    \item \label{item:cauchy-temporal-function} Cauchy temporal
        function if $t$ is a temporal function such that all level
        sets are Cauchy hypersurfaces.
    \end{definitionlist}
\end{definition}
\begin{remark}[Time functions]
    \label{remark:temporal-function}
    ~
    \begin{remarklist}
    \item \label{item:sign-convention-and-temporal-function} With the
        other sign convention for the metric a temporal function has
        \emph{past} directed gradient.
    \item \label{item:temporal-level-sets-are-submanifolds} If $t$ is
        temporal, its level sets are (if nonempty) embedded smooth
        submanifolds since the gradient is non-zero everywhere and
        hence every value is a regular value. Note that they do not
        need to be Cauchy hypersurfaces at all: In fact, remove a
        single point from Minkowski spacetime then the usual time
        function is temporal but there is no Cauchy hypersurface at
        all.
    \item \label{item:gradient-flow-of-temp.function} The gradient
        flow of $t$ gives a diffeomorphism between the different level
        sets of $t$. Since every timelike curve intersects a Cauchy
        hypersurface precisely once we see that this gives a
        diffeomorphism
        \begin{equation}
            \label{eq:spacetime-splitting-with-tempfunction}
            M \simeq t(M) \times \Sigma_{t_0},
        \end{equation}
        and all Cauchy hypersurfaces are diffeomorphic to a given
        reference Cauchy hypersurface $\Sigma_{t_0}$. This gives a
        very strong implication on the structure of $M$.
        \begin{figure}
            \centering
            \input{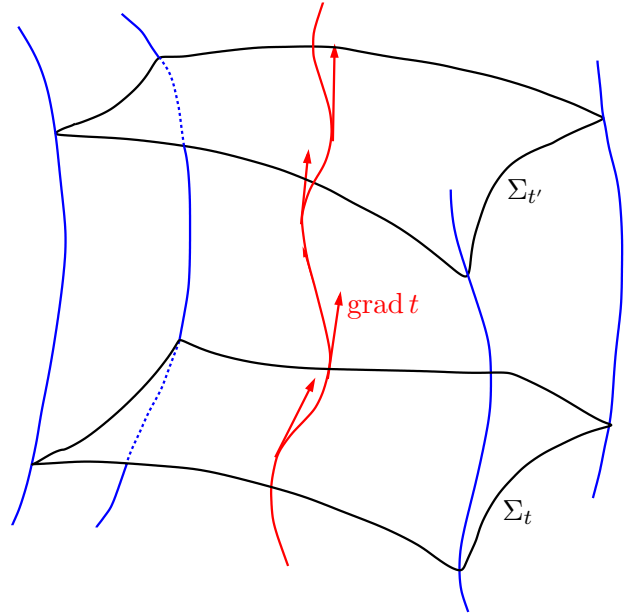}
            \caption{\label{fig:gradient-flow-of-cauchy-temp-function}%
              The gradient flow of a Cauchy temporal function.
            }
        \end{figure}
    \item \label{item:image-of-temp-function} By rescaling $t$ we can
        always assume that the image of $t$ is the whole real line
        $\mathbb{R}$. This follows as the image of $t$ is necessarily
        open and connected (for connected $M$).
    \end{remarklist}
\end{remark}
The following celebrated and non-trivial theorem brings together the
notions of globally hyperbolic spacetimes and the existence of Cauchy
temporal functions.
\begin{theorem}
    \label{theorem:globalhyp-and-cauchy-surface}
    \index{Globally hyperbolic}%
    \index{Cauchy hypersurface}%
    \index{Temporal function}%
    Let $(M, g)$ be a connected time-oriented Lorentz manifold. Then
    the following statements are equivalent:
    \begin{theoremlist}
    \item \label{item:globalhyp} $(M, g)$ is globally hyperbolic.
    \item \label{item:cauchy-existence} There exists a topological
        Cauchy hypersurface.
    \item \label{item:smooth-spacelike-cauchy-existence} There exists
        a smooth spacelike Cauchy hypersurface.
    \end{theoremlist}
    In this case there even exists a Cauchy temporal function $t$ and
    $(M, g)$ is isometrically diffeomorphic to the product manifold
    \begin{equation}
        \label{eq:cauchy-surface-splitting-of-M-and-g}
        \mathbb{R} \times \Sigma
        \quad \textrm{with metric} \quad
        g = \beta \D t^2 - g_t,
    \end{equation}
    where $\beta \in \Cinfty(\mathbb{R} \times \Sigma)$ is positive
    and $g_t \in \Secinfty(\Sym^2 T^*\Sigma)$ is a Riemannian metric
    on $\Sigma$ depending smoothly on $t$. Moreover, each level set
    \begin{equation}
        \label{eq:levelset-of-cauchy-temporal-function}
        \Sigma_t
        = \{ (t, \sigma) \in \mathbb{R} \times \Sigma \}
        \subseteq M
    \end{equation}
    of the temporal function $t$ is a smooth spacelike Cauchy
    hypersurface.
\end{theorem}
\begin{remark}
    \label{remark:globalhyp-and-cauchy-surface}
    The equivalence of \refitem{item:globalhyp} and
    \refitem{item:cauchy-existence} is the celebrated theorem of
    Geroch\index{Geroch's Theorem} \cite{geroch:1970a}. The
    enhancement to the smooth setting is due to Bernal and S\'anchez
    \cite{bernal.sanchez:2003a, bernal.sanchez:2005a,
      bernal.sanchez:2006a, bernal.sanchez:2007a}. Conversely, having
    a metric of the form $\beta \D t^2 -g_t$ on $\mathbb{R} \times
    \Sigma$ it is trivial to see that all level sets $\Sigma_t$ are
    spacelike hypersurfaces diffeomorphic to $\Sigma$. Note however,
    that the form \eqref{eq:cauchy-surface-splitting-of-M-and-g} alone
    does not guarantee that the $\Sigma_t$ are Cauchy hypersurfaces.
\end{remark}
\begin{example}[Minkowski strip]
    \label{example:minkowski-strip}
    \index{Minkowski strip}%
    We consider $\Sigma= (a,b)$ an open interval with $- \infty < a <
    b < + \infty$ and $M = \mathbb{R} \times \Sigma \subseteq
    \mathbb{R}^2$ as open subset of Minkowski space.
    \begin{figure}
        \centering
        \input{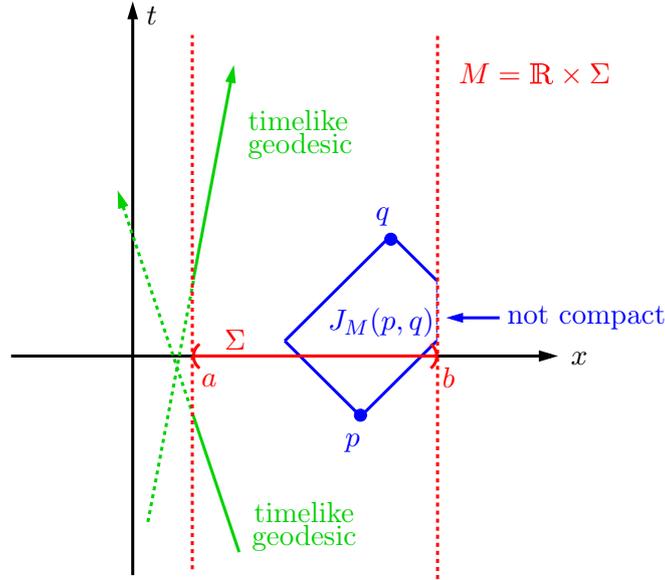}
        \caption{\label{fig:minkowski-space}%
          The Minkowski strip
        }
    \end{figure}
    Then $\Sigma_t$ is not a Cauchy hypersurface for any $t$.  This is
    clear from the observation that there are inextensible timelike
    geodesics not passing through $\Sigma_t$. In fact, $M$ is not
    globally hyperbolic at all: while $M$ is causal (and even strongly
    causal) it fails to satisfy the second condition of global
    hyperbolicity: there are diamonds $J_M(p,q)$ which are not
    compact, see Figure~\ref{fig:minkowski-space}. Thus by
    Theorem~\ref{theorem:globalhyp-and-cauchy-surface} there cannot
    exist any Cauchy hypersurface. Nevertheless, the metric is of the
    very simple form
    \begin{equation}
        \label{eq:metric-on-minkowski-strip}
        g = \D t^2 - \D x^2.
    \end{equation}
\end{example}

The problem with this example comes from the geometric feature of the
open interval  $\Sigma = (a,b) \subseteq \mathbb{R}$ of being ``too
short''. The following proposition gives now a sufficient condition
such that this can not happen:
\begin{proposition}
    \label{proposition:1}
    \index{Geodesically complete}%
    Let $M = \mathbb{R} \times \Sigma$ with Lorentz metric
    \begin{equation}
        \label{eq:metric-splitting-on-product-space}
        g = \frac{1}{2} \D t \vee \D t - f(t) g_\Sigma,
    \end{equation}
    where $g_\Sigma$ is a Riemannian metric on $\Sigma$ and $f \in
    \Cinfty(\mathbb{R})$ is positive. The time-orientation is such
    that $\frac{\partial}{\partial t}$ is future directed. Then $(M,
    g)$ is globally hyperbolic if and only if $g_\Sigma$ is
    geodesically complete.
\end{proposition}
For a proof see e.g. \cite[Lem.~A.5.14]{baer.ginoux.pfaeffle:2007a}.
Many of the physically interesting examples of spacetimes from general
relativity can be brought to the form
\eqref{eq:metric-splitting-on-product-space} whence the above
Proposition can be used to discuss the global hyperbolicity of $(M,
g)$.

For later use we mention the following result which still enhances
Theorem~\ref{theorem:globalhyp-and-cauchy-surface}, see
\cite[Thm.~1.2]{bernal.sanchez:2006a}.
\begin{theorem}
    \label{theorem:1}
    \index{Temporal function}%
    Let $(M, g)$ be globally hyperbolic and let $\Sigma \subseteq M$
    be a smooth spacelike Cauchy hypersurface. Then there exists a
    Cauchy temporal function $t$ such that the $t = 0$ Cauchy
    hypersurface coincides with $\Sigma$.
\end{theorem}


%% file: cauchy.tex
%
%

Having the notion of a Cauchy hypersurface we are now in the position
to formulate the Cauchy problem for a normally hyperbolic differential
operator. Here we still be rather informal only fixing the principal
ideas. The precise formulation of the Cauchy problem will be given and
discussed in detail in Section~\ref{satz:cauchy-problem}.

Thus let $(M, g)$ be globally hyperbolic and $\Sigma \subseteq M$ a
smooth Cauchy hypersurface which we assume to be spacelike throughout
the following. At a given point $p \in \Sigma \subseteq M$ the tangent
plane $T_p\Sigma \subseteq T_pM$ is spacelike whence there exists a
unique vector $\mathfrak{n}_p \in T_pM$ which satisfies
\begin{equation}
    \label{eq:future-directed-normal-vector-normal}
    g_p(\mathfrak{n}_p, T_p\Sigma) = 0,
\end{equation}
\begin{equation}
    \label{eq:future-directed-normal-vector-normalized}
    g_p (\mathfrak{n}_p, \mathfrak{n}_p) = 1,
\end{equation}
\begin{equation}
    \label{eq:future-directed-normal-vector-future}
    \mathfrak{n}_p \; \textrm{is future directed}.
\end{equation}
\index{Normal vector}%
This vector is called the \emph{future directed normal vector} of
$\Sigma$ at $p$. Taking all points $p \in \Sigma$ we obtain the future
directed normal vector field of $\Sigma$, i.e. the vector field
\begin{equation}
    \label{eq:future-directed-normal-field}
    \index{Normal vector field}%
    \mathfrak{n} \in \Secinfty\left(TM \at{\Sigma}\right),
\end{equation}
such that \eqref{eq:future-directed-normal-vector-normal},
\eqref{eq:future-directed-normal-vector-normalized}, and
\eqref{eq:future-directed-normal-vector-future} hold for every $p \in
\Sigma$. Since $\Sigma$ is a smooth submanifold, $\mathfrak{n}$ is
smooth itself.
\begin{figure}
    \centering
    \input{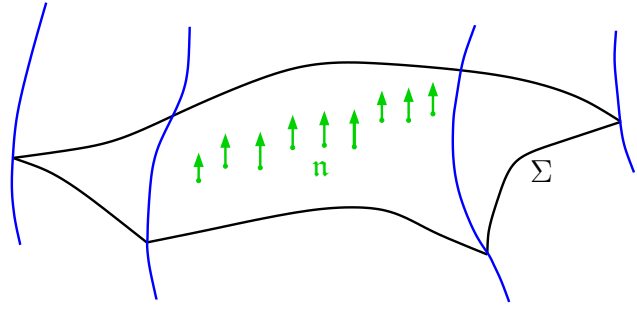}
    \caption{\label{fig:future-directed-normal-field}%
      The future directed normal vector field of a Cauchy hypersurface
      $\Sigma$.
    }
\end{figure}
We consider now a normally hyperbolic differential operator $D \in
\Diffop(E)$ on some vector bundle $E \longrightarrow M$. Then this
operator gives the homogeneous wave equation
\begin{equation}
    \label{eq:homogeneous-field-equation}
    \index{Wave equation!homogeneous}%
    D u = 0,
\end{equation}
or more generally
\begin{equation}
    \label{eq:inhomogeneous-field-equation}
    \index{Wave equation!inhomogeneous}%
    \index{Source term}%
    \index{Inhomogeneity}%
    D u = v,
\end{equation}
where $v \in \Secinfty(E)$ is a given \emph{inhomogeneity} and $u \in
\Secinfty(E)$ is the field we are looking for. Having specified the
inhomogeneity which physically corresponds to a \emph{source term}, we
can try to find a solution $u$ which has specified initial values and
initial velocities on $\Sigma$. More precisely, we want
\begin{equation}
    \label{eq:initial-values}
    \index{Inital values}%
    u \At{\Sigma} = u_0 \in \Secinfty(E \at{\Sigma})
\end{equation}
and
\begin{equation}
    \label{eq:initial-velocities}
    \index{Inital velocity}%
    \nabla^E_{\mathfrak{n}} u \At{\Sigma}
    = \dot{u}_0 \in \Secinfty(E \at{\Sigma})
\end{equation}
with a priori given $u_0$ and $\dot{u}_0$. The hope is that this
\emIndex{Cauchy problem} has a unique solution, probably after
considering compactly supported $v$, $u_0$, and $\dot{u}_0$. Moreover,
one hopes that the solution $u$ depends in a reasonably continuous way
on the initial values $u_0$ and $\dot{u}_0$ and perhaps also on $v$.

More generally, one can try to find solutions $u \in \Sec[-\infty](E)$
for distributional initial values $u_0, \dot{u}_0 \in
\Sec[-\infty](E\at{\Sigma})$ and distributional $v \in
\Sec[-\infty](E)$. In general, however, we meet difficulties with this
Cauchy problem. Namely, we can not just restrict a distribution $u$ to
a submanifold $\Sigma$ in order to make sense out of
\eqref{eq:initial-values} and \eqref{eq:initial-velocities}: this is
only possible if $u$ behaves nicely enough around $\Sigma$. Clearly,
the restriction is not problematic as soon as $u$ is at least
$\Fun[1]$.

As a last comment we note that the Cauchy problem still makes sense if
$\Sigma$ is just a spacelike hypersurface which is not necessarily a
Cauchy hypersurface. In this case we still can hope to get a solution
to the Cauchy problem but we have to expect non-uniqueness for obvious
reasons.

\begin{figure}
    \centering
    \input{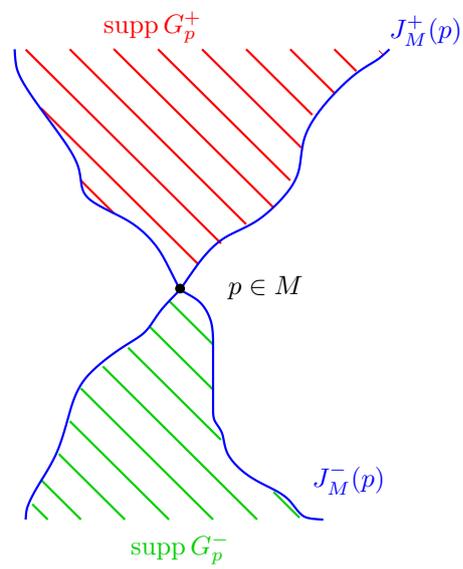}
    \caption{\label{fig:supp-of-greens-function}%
      The support of a Green function at a point $p \in M$.
    }
\end{figure}

The main idea to attack this problem is to construct particular
distributional solutions, the \emph{fundamental solutions} $F_p \in
\Sec[-\infty](E) \tensor E_p^* \tensor \Dichten T_p^*M$ such that
\begin{equation}
    \label{eq:fundamental-solution}
    \index{Fundamental solution}%
    \index{Delta-Functional@$\delta$-Functional}%
    D F_p = \delta_p,
\end{equation}
where $\delta_p$ is the $\delta$-distribution at $p \in M$ viewed as
$E_p^* \tensor \Dichten T_p^*M$-valued generalized section of $E$,
i.e. for a test section $\mu \in \Secinfty_0(E^* \tensor \Dichten
T^*M)$ we have
\begin{equation}
    \label{eq:delta-distribution}
    \delta_p(\mu) = \mu(p) \in E_p^* \tensor \Dichten T_p^*M.
\end{equation}
\begin{definition}[Green function]
    \label{definition:green-function}
    \index{Fundamental solution}%
    \index{Green function!advanced}%
    \index{Green function!retarded}%
    Let $p \in M$. A generalized section $F_p$ of $E$ which satisfies
    \eqref{eq:fundamental-solution} is called fundamental solution of
    $D$ at $p$. If a fundamental solution $F_p^\pm$ in addition
    satisfies
    \begin{equation}
        \label{eq:supp-condition-of-green-function}
        \supp F_p^\pm \subseteq J_M^\pm (p),
    \end{equation}
    then $F_p^\pm$ is called advanced or retarded Green function of
    $D$ at $p$, respectively.
\end{definition}
\begin{remark}[Green function]
    \label{remark:green-function}
    Note that the notion of a fundamental solution makes sense for
    every differential operator on any manifold. The notion of
    advanced and retarded Green functions makes sense for any
    differential operator on a time-oriented Lorentz manifold, see
    also Figure~\ref{fig:supp-of-greens-function}.
\end{remark}
The remaining part of these notes are now devoted to the study of
existence and uniqueness of Green functions $F_p^\pm$. Moreover, we
have to relate the Green functions to the Cauchy problem for $D$.
Here it will be important not only to have a Green function $F_p$ for
every $p \in M$. We also will need a reasonable dependence of $F_p$ on
$p$.


%% file: chap3.tex
%
%

\chapter{The Local Theory of Wave Equations}
\label{cha:LocalTheory}

The purpose of this chapter is to discuss the existence and uniqueness
of fundamental solutions for the wave equation determined by a
normally hyperbolic differential operator at least on small enough
open subsets of $M$. Thus the global structure of $M$ does not yet
play a role in this chapter. Nevertheless, already locally the
geometry enters in form of non-trivial curvature terms and resulting
non-trivial parallel transports. Thus already at this stage we will be
beyond the usual flat situation of the wave equation in
$\mathbb{R}^{2n}$.

We basically follow \cite{baer.ginoux.pfaeffle:2007a} and construct
the fundamental solution first in the flat case of Minkowski
spacetime. Here we use the approach of Riesz \cite{riesz:1949a} by
specifying the fundamental solutions using holomorphic function
techniques. Then one constructs a formal solution on a domain as a
series with certain coefficients, the \emph{Hadamard
  coefficients}. This solution will be a series with no good control
of convergence and in fact, no convergence in general. Thus an
additional step is needed to find the ``true'' fundamental
solutions. To this end certain cut-off parameters are introduced
yielding a convergent series which is however no longer a fundamental
solution but only a parametrix. With some convolution tricks this can
be cured in the last step. The fundamental solution will have nice
causal properties allowing to find solutions to the inhomogeneous wave
equation with good causal properties as well.

%
%

\section{The d'Alembert Operator on Minkowski Spacetime}
\label{sec:dAlembertOnMinkowskiSpaceTime}

\input{minkowski}

%
%

\section{The Riesz Distributions on a Convex Domain}
\label{sec:RieszOnConvexDomain}

\input{riesz}

%
%

\section{The Hadamard Coefficients}
\label{sec:HadamardCoefficients}

\input{hadamard}

%
%

\section{The Fundamental Solution on Small Neighborhoods}
\label{satz:LocalFundamentalSolution}

\input{fundamental}

%
%

\section{Solving the Wave Equation Locally}
\label{satz:SolvingWaveEqLocally}

\input{local}


%% file: minkowski.tex
%
%

As warming up we consider the most simple case of a normally
hyperbolic differential operator, the \emph{d'Alembert operator} on
flat \Index{Minkowski spacetime}.

%
%

\subsection{The Riesz Distributions}
\label{subsec:riesz-distributions}

We shall not only construct the fundamental solutions of the
d'Alembert operator
\begin{equation}
    \label{eq:dAlembertOp}
    \index{dAlembertian@d'Alembertian}%
    \dAlembert = \frac{\partial^2}{\partial t^2} - \Laplace
\end{equation}
in $n$ dimensions but a local family of distributions associated to
$\dAlembert$. Sometimes we will set $t = x^0$ and $\vec{x} = (x^1,
\ldots, x^{n-1})$ for abbreviation. In more physical terms, we set the
speed of light $c$ to $1$ by choosing appropriate units. Here we
follow essentially the approach of Riesz \cite{riesz:1949a}, see
\cite[Sect.~1.2]{baer.ginoux.pfaeffle:2007a} for a modern presentation
of this approach.

\begin{figure}
    \centering
    \includegraphics{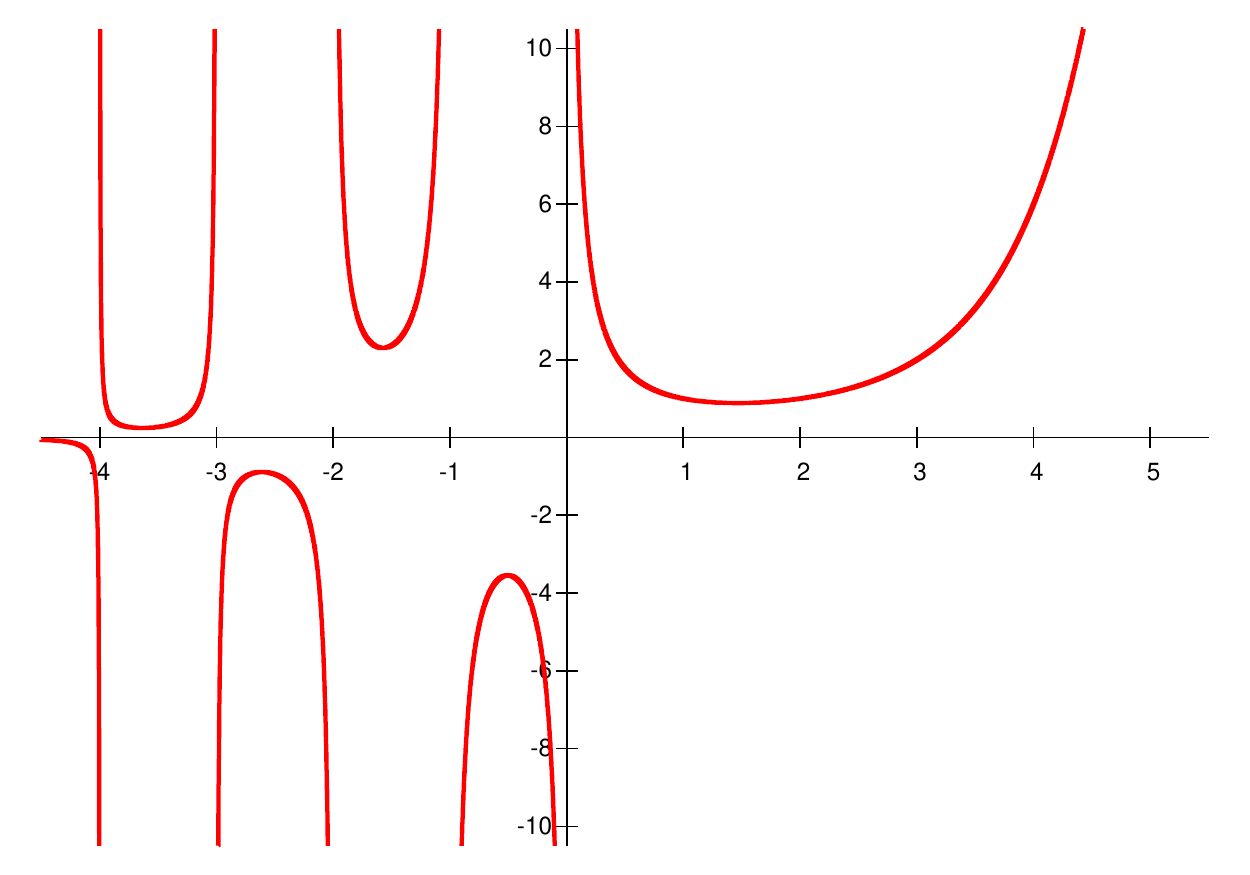}
    \caption{
      \label{fig:gamma-function}%
      The Gamma function along the real axis
    }
\end{figure}

Using the Minkowski metric $\eta$ we have the following function, also
denoted by $\eta$,
\begin{equation}
    \label{eq:etanorm}
    \eta(x) = \eta(x, x)
\end{equation}
on $\mathbb{R}^n$. Clearly $\eta \in \Pol^2(\mathbb{R}^2)$ is a
homogeneous quadratic polynomial. Explicitly, in the standard
coordinates we have
\begin{equation}
    \label{eq:etanorm-standbasis}
    \eta (x^0, \ldots, x^{n-1})
    = (x^0)^2 - \sum_{i=1}^{n-1} (x^i)^2
    = t^2 - (\vec{x})^2.
\end{equation}
We consider the following family of continuous functions on Minkowski
spacetime:
\begin{definition}
    \label{definition:riesz-function}
    Let $\alpha \in \mathbb{C}$ have $\RE (\alpha) > n$. Then one
    defines
    \begin{equation}
        \label{eq:riesz-function}
        R^\pm (\alpha) (x)
        =
        \begin{cases}
            c(\alpha, n) \eta(x)^{\frac{\alpha-n}{2}}
            & \mathrm{for} \, \alpha \in I^\pm(0) \\
            0 & \mathrm{else},
        \end{cases}
    \end{equation}
    where the coefficient is
    \begin{equation}
        \label{eq:riesz-coefficient}
        c(\alpha ,n)
        = \frac{2^{1-\alpha} \pi^{\frac{2-n}{2}}}
        {\Gamma(\frac{\alpha}{2}) \Gamma(\frac{\alpha-n}{2} +1)}.
    \end{equation}
\end{definition}
\begin{remark}[Gamma function]
    \label{remark:gamma-function}
    The \emIndex{Gamma function}
    \begin{equation}
        \label{eq:gamma-function}
        \Gamma: \mathbb{C} \backslash \{0, -1, -2, \ldots \}
        \longrightarrow \mathbb{C}
    \end{equation}
    is known to be a holomorphic function with simple poles at $-n$
    for $n \in \mathbb{N}_0$. One has the following properties:
    \begin{remarklist}
    \item \label{item:gamma-residue} The residue at $-n \in
        \mathbb{N}_0$ is given by
        \begin{equation}
            \label{eq:gamma-residue}
            \mathrm{res}_{-n} \Gamma
            = \frac{(-1)^n}{n!}.
        \end{equation}
    \item \label{item:gamma-recursive-identiy} For $z \in \mathbb{C}
        \setminus \{0, -1, -2, \ldots \}$ one has the functional
        equation
        \begin{equation}
            \label{eq:gamma-recursive-identity}
            \Gamma (z+1) = z \Gamma (z)
            \quad
            \textrm{with}
            \quad
            \Gamma(1) = 1.
        \end{equation}
    \item \label{item:gamma-is-faculty} For $n \in \mathbb{N}_0$ one
        obtains from \eqref{eq:gamma-recursive-identity} immediately
        \begin{equation}
            \label{eq:1}
            \Gamma(n+1) = n!.
        \end{equation}
    \item \label{item:gamma-euler-formula} For $\RE(z) > 0$ one has
        Euler's integral formula
        \begin{equation}
            \label{eq:gamma-euler-fomula}
            \index{Euler's integral formula}%
            \Gamma(z) = \int_0^\infty t^{z-1} \E^{-t} \D t
        \end{equation}
        in the sense of an improper Riemann integral.
    \item \label{item:gamma-legrende-duplication} For all $z \in
        \mathbb{C} \setminus \{0, -1, -2, \ldots \}$ one has
        \emIndex{Legendre's duplication formula}
        \begin{equation}
            \label{eq:gamma-legendre-duplication}
            \Gamma(z) \Gamma(z + \frac{1}{2})
            = 2^{1-2z} \sqrt{\pi} \: \Gamma(2z).
        \end{equation}
    \end{remarklist}
\end{remark}
For more details and proofs of the above properties of $\Gamma$ we
refer to any textbook on complex function theory like
e.g.~\cite[Chap.~2]{remmert:1995a}. The graph of the Gamma function
along the real axis can be seen in Figure~\ref{fig:gamma-function}.

Since the Gamma function $\Gamma$ has no zeros we conclude that the
prefactor $c(\alpha, n)$ is holomorphic for all $\alpha \in
\mathbb{C}$: indeed, for those $\alpha \in \mathbb{C}$ where
$\Gamma(\frac{\alpha}{2})$ or $\Gamma(\frac{\alpha-n}{2} + 1)$ has a
pole the inverse is well-defined and has a zero of the same (first)
order as the pole of the $\Gamma$ function. This happens for
\[
\frac{\alpha}{2} = 0, -1, -2, \ldots
\quad
\textrm{and}
\quad
\frac{\alpha-n}{2} +1 =0, -1, -2, \ldots
\]
Thus we conclude
\begin{equation}
    \label{eq:riesz-coefficients-zeros}
    c(\alpha, n) = 0
    \quad
    \textrm{iff}
    \quad
    \alpha \in
    \left\{
        -2k \; \big| \; k \in \mathbb{N}_0
    \right\}
    \cup
    \left\{
        n- 2k \; \big| \; k \in \mathbb{N}_0
    \right\},
\end{equation}
since the nominator has clearly no zeros. For $\alpha$ not being in
the above special set but still with $\RE(\alpha) > n$, the function
$R^\pm (\alpha)$ is continuous but not smooth on $\mathbb{R}^n$:
\begin{lemma}
    \label{lemma:riesz-function-continuity}
    For $\RE(\alpha) > n$ the function $R^\pm(\alpha)$ is continuous
    on $\mathbb{R}^n$.  It is smooth in $I^\pm(0)$ and in
    $\mathbb{R}^n \setminus J^\pm(0)$.
\end{lemma}
\begin{proof}
    The function $x \mapsto c(\alpha, n)\eta(x)^{\frac{\alpha-n}{2}}$
    is clearly smooth for $x \in I^\pm(0)$ since here $\eta(x) > 0$.
    Conversely, on the open subset $\mathbb{R}^n \setminus J^\pm(0)$
    the function $R^\pm(\alpha)$ is zero and hence smooth, too. The
    continuity follows as $\eta(x) \longrightarrow 0$ for $x \in
    I^\pm(0)$ with $x \longrightarrow \partial I^\pm(0)$ and $\alpha >
    n$ guarantees that the function $0 \leq \xi \mapsto
    \xi^{\frac{\alpha-n}{2}}$ is at least continuous at $0$.
\end{proof}
The next lemma clarifies the behaviour under Lorentz transformations.
\begin{lemma}
    \label{lemma:riesz-functions-lorenz-trafos}
    \index{Lorentz transformation!orthochronous}%
    \index{Time-reversal}%
    Let $\Lambda \in \group{L}^\uparrow (1, n-1)$ be an orthochronous
    Lorentz transformation and $\RE(\alpha) > n$. Then
    \begin{equation}
        \label{eq:ries-invariant-under-lorentz-trafo}
        \Lambda^* R^\pm(\alpha) = R^\pm(\alpha).
    \end{equation}
    If $\mathrm{T} \in \group{L}(1, n-1)$ is the time-reversal $x^0
    \mapsto -x^0$ then
    \begin{equation}
        \label{eq:riesz-time-reversal}
        \mathrm{T}^*R^\pm(\alpha) = R^\mp(\alpha).
    \end{equation}
\end{lemma}
\begin{proof}
    If $\Lambda \in \group{L}(1,n-1)$ is an arbitrary Lorentz
    transformation then by the very definition of $\group{L}(1,n-1)$
    we have
    \[
    (\Lambda^* \eta) (x) = \eta (\Lambda x)
    = \eta (\Lambda x, \Lambda x) = \eta(x, x) = \eta (x),
    \]
    and thus $\Lambda^* \eta = \eta$. But then
    \eqref{eq:ries-invariant-under-lorentz-trafo} and
    \eqref{eq:riesz-time-reversal} are obvious since the light cones
    $J^\pm(0)$ are mapped to $J^\pm(0)$ and to $J^\mp(0)$ under
    $\Lambda \in \group{L}^\uparrow(1, n-1)$ and under $\mathrm{T}$,
    respectively.
\end{proof}
In particular, it would be sufficient to consider $R^+(\alpha)$ alone
since we can recover every information about $R^-(\alpha)$ from
$R^+(\alpha)$ via \eqref{eq:riesz-time-reversal}.

Since $R^\pm(\alpha) \in \Fun[0](\mathbb{R}^n)$ we can consider
$R^\pm(\alpha)$ also as a distribution (of order zero) via the usual
identification, i.e.
\begin{equation}
    \label{eq:riesz-distribution}
    R^\pm(\alpha): \varphi \; \mapsto \;
    \int_{\mathbb{R}^n} \varphi(x) R^\pm(\alpha)(x) \D^n x
\end{equation}
\index{Lorentz density}%
for test functions $\varphi \in \Cinfty_0(\mathbb{R}^n)$. Here and in
the following we use the \Index{Lebesgue measure} $\D^n x$ for
integration.  Note that this coincides with the Lorentz density
induced by $\eta$.
\begin{lemma}
    \label{lemma:ries-distribution-is-homolorphic}
    Let $\RE(\alpha) > n$.
    \begin{lemmalist}
    \item \label{item:riesz-function-holomorphic} For every $x \in
        \mathbb{R}^n$ the function
        \begin{equation}
            \label{eq:ries-function-holomorphic}
            \alpha \; \mapsto \; R^\pm(\alpha)(x)
        \end{equation}
        is holomorphic.
    \item \label{item:ries-distribution-holomorphic} For every test
        function $\varphi \in \Cinfty_0(\mathbb{R}^n)$ the function
        \begin{equation}
            \label{eq:riesz-distribution-holomorphic}
            \alpha \; \mapsto \; R^\pm(\alpha) (\varphi)
        \end{equation}
        is holomorphic.
    \end{lemmalist}
\end{lemma}
\begin{proof}
    The first part is clear as the Gamma function and hence the
    coefficient $c(\alpha, n)$ is holomorphic. Moreover, for $x \in
    J^\pm(0)$ the map $\alpha \mapsto \eta(x)^{\frac{\alpha-n}{2}}$ is
    holomorphic. However, this pointwise holomorphy of $R^\pm(\alpha)$
    is not the relevant feature for the following.  Instead, we need
    the second part. To prove this, we consider $\varphi \in
    \Cinfty_0(\mathbb{R})$. Then
    \[
    R^\pm(\alpha) (\varphi)
    = \int_{\mathbb{R}^n} \varphi(x) R^\pm(\alpha)(x) \D^n x
    = \int_{\supp \varphi} \varphi(x) R^\pm(\alpha)(x) \D^n x.
    \]
    Since $\supp \varphi$ is compact we can exchange the orders of
    integration for every closed triangle path $\Delta$ in $\{ \alpha
    \in \mathbb{C} \; | \; \RE(\alpha) > n \}$ by Fubini's
    theorem. Thus
    \[
    \int_\Delta R^\pm(\alpha) (\varphi) \D \alpha
    = \int_\Delta \int_{\supp \varphi}
    \varphi(x) R^\pm(\alpha)(x) \D^n x \D \alpha
    = \int_{\supp \varphi} \varphi(x)
    \int_\Delta R^\pm(\alpha)(x) \D \alpha \D^n x
    = 0,
    \]
    since $R^\pm(\alpha)(x)$ is holomorphic for every $x \in
    \mathbb{R}^n$. It follows by Morera's theorem that
    \eqref{eq:riesz-distribution-holomorphic} is holomorphic, too.
\end{proof}

In this sense we have a holomorphic map
\begin{equation}
    \label{eq:holomorphic-riesz-distribution}
    \{ \alpha \in \mathbb{C} \; | \; \RE(\alpha) > n \} \ni \alpha
    \; \mapsto \;
    R^\pm(\alpha) \in \mathcal{D}'(\mathbb{R}^n)
\end{equation}
with values in the distributions.  The key idea is now to investigate
\eqref{eq:holomorphic-riesz-distribution} in detail to show that, as a
holomorphic map, it has a unique extension to the whole complex plane
$\mathbb{C}$. To this end we need the following technical lemma:
\begin{lemma}
    \label{lemma:identities-for-riesz-distribution}
    In the sense of continuous functions we have:
    \begin{lemmalist}
    \item \label{item:riesz-multiplication-by-eta-identity} For
        $\RE(\alpha) > n$ we have
        \[
        \eta R^\pm(\alpha) = \alpha (\alpha - n + 2) R^\pm(\alpha+2).
        \]
    \item \label{item:riesz-differentiation-rule} For $\RE(\alpha) > n
        + 2k$ the function $R^\pm(\alpha)$ is $\Fun$ and we have
        \begin{equation}
            \label{eq:riesz-differentiation-rule}
            \frac{\partial}{\partial x^i} R^\pm(\alpha)
            = \frac{1}{\alpha - 2} R^\pm(\alpha-2) \eta_{ij} x^j.
        \end{equation}
    \item \label{item:riesz-and-gradient} For $\RE(\alpha) > n$ we
        have
        \begin{equation}
            \label{eq:riesz-and-gradient-of-eta}
            \gradient \eta \cdot R^\pm(\alpha)
            = 2 \alpha \gradient R^\pm(\alpha + 2).
        \end{equation}
    \item \label{item:dAlembert-of-riesz} For $\RE(\alpha) > n+2$ we
        have
        \begin{equation}
            \label{eq:dAlembert-of-riesz}
            \dAlembert R^\pm(\alpha + 2) = R^\pm(\alpha).
        \end{equation}
    \end{lemmalist}
\end{lemma}
\begin{proof}
    The first part is a simple calculation. We have
    \begin{align*}
        \alpha (\alpha + 2 - n) R^\pm(\alpha + 2)
        & = \alpha (\alpha + 2 - n) c(\alpha+2 , n)
        \eta^{\frac{\alpha+2-n}{2}} \\
        & = \alpha (\alpha + 2 - n) c(\alpha+2 , n)
        \eta^{\frac{\alpha-n}{2}} \eta \\
        &= \frac{\alpha (\alpha+2-n) c(\alpha+2, n)}{c(\alpha,n)}
        \eta R^\pm(\alpha),
    \end{align*}
    and
    \begin{align*}
        \frac{c(\alpha+2, n)}{c(\alpha, n)}
        = \frac{
          2^{1-2-\alpha} \pi^{\frac{2-n}{2}} \Gamma(\frac{\alpha}{2})
          \Gamma(\frac{\alpha-n}{2}+1)
        }{
          \Gamma(\frac{\alpha+2}{2}) \Gamma(\frac{\alpha+2-n}{2}+1)
          2^{1-\alpha} \pi^{\frac{2-n}{2}}
        }
        = \frac{
          2^{-2} \Gamma(\frac{\alpha}{2}) \Gamma(\frac{\alpha-n}{2}+1)
        }{
          \frac{\alpha}{2} \Gamma(\frac{\alpha}{2})
          \frac{\alpha+2-n}{2} \Gamma(\frac{\alpha-n}{2}+1)
        }
        = \frac{1}{\alpha (\alpha+2-n)}.
        \tag{$*$}
    \end{align*}
    For the second part we recall that in $I^\pm(0)$ the function
    $R^\pm(\alpha)$ is smooth as well as in $\mathbb{R}^n \setminus
    J^\pm(0)$. On the latter, the function and hence all its
    derivatives are zero. In $I^\pm(0)$ we compute
    \begin{align*}
        \frac{\partial}{\partial x^i} R^\pm(\alpha) \At{I^\pm(0)}
        & = c(\alpha, n) \frac{\partial}{\partial x^i}
        \eta(x)^{\frac{\alpha-n}{2}}
        = c(\alpha, n) \frac{\alpha-n}{2}
        \eta(x)^{\frac{\alpha-n}{2}-1} \frac{\partial}{\partial x^i}
        \eta(x) \\
        & = c(\alpha, n) \frac{\alpha-n}{2}
        \eta(x)^{\frac{\alpha-2-n}{2}} 2 \eta_{ij} x^j
        = c(\alpha, n) (\alpha-n) \eta^{\frac{\alpha-2-n}{2}}
        \eta_{ij} x^j \\
        & = \frac{c(\alpha, n)}{c(\alpha-2, n)} (\alpha-n)
        R^\pm(\alpha-2) \eta_{ij} x^j
        \stackrel{\mathclap{(*)}}{=}
        \frac{1}{(\alpha-2)(\alpha-n)} (\alpha-n)
        R^\pm(\alpha-2) \eta_{ij} x^j \\
        & = \frac{1}{(\alpha-2)} R^\pm(\alpha-2) \eta_{ij} x^j.
    \end{align*}
    Now if $\RE(\alpha) > n + 2k$ then $\RE(\alpha-2) > n + 2k - 2$ is
    still larger than $n$ for positive $k \in \mathbb{N}$. Thus the
    partial derivative $\frac{\partial}{\partial x^i} R^\pm(\alpha)
    \at{I^\pm(0)}$ is the continuous function $\frac{1}{(\alpha-2)}
    R^\pm(\alpha-2) \eta_{ij} x^j$ in $I^\pm(0)$ which continuously
    extends to $\mathbb{R}^n$ by setting it zero outside of
    $I^\pm(0)$. Indeed, since $R^\pm(\alpha-2)$ has this as continuous
    extension, we obtain a continuous extension of
    $\frac{\partial}{\partial x^i} R^\pm(\alpha) \at{I^\pm(0)}$. But
    this matches the partial derivative of $R^\pm(\alpha)$ outside of
    $J^\pm(0)$. Thus we obtain a continuous partial derivative
    $\frac{\partial}{\partial x^i} R^\pm(\alpha)$ on all of Minkowski
    space $\mathbb{R}^n$ which shows that $R^\pm(\alpha)$ is at least
    $\Fun[1]$. By induction we can proceed as long as $\alpha - 2k >
    n$. The third part is now a simple consequence of the first and
    second part. We have
    \begin{align*}
        \gradient \eta
        = \left(
            \frac{\partial \eta}{\partial x^i} \D x^i
        \right)^\sharp
        = \frac{\partial \eta}{\partial x^i} \eta^{ij}
        \frac{\partial}{\partial x^j}
        = 2 \eta_{ik} x^k \eta^{ij} \frac{\partial}{\partial x^j}
        = 2 x^j \frac{\partial}{\partial x^j}
        = 2 \xi.
    \end{align*}
    Thus $\gradient \eta$ is twice the \emIndex{Euler vector field} on
    $\mathbb{R}^n$, which, remarkably, does not depend on the metric
    $\eta$ but only on the vector space structure. Using
    \eqref{eq:riesz-differentiation-rule} we compute for $\RE(\alpha)
    > n$
    \begin{align*}
        2 \alpha \gradient R^\pm(\alpha+2)
        & = 2 \alpha \eta^{ij}
        \frac{\partial R^\pm(\alpha+2)}{\partial x^i}
        \frac{\partial}{\partial x^j}
        = 2 \alpha \eta^{ij} \frac{1}{\alpha+2-2} R^\pm(\alpha)
        \eta_{ik} x^k \frac{\partial}{\partial x^j} \\
        & = R^\pm(\alpha) 2 x^k \frac{\partial}{\partial x^k}
        = R^\pm(\alpha) \gradient \eta.
    \end{align*}
    For the last part we use \eqref{eq:riesz-differentiation-rule}
    twice and obtain
    \begin{align*}
        \dAlembert R^\pm(\alpha+2)
        & = \eta^{ij} \frac{\partial}{\partial x^i}
        \frac{\partial}{\partial x^j} R^\pm(\alpha+2) \\
        & = \eta^{ij} \frac{\partial}{\partial x^i}
        \left(
            \frac{1}{\alpha+2-2} R^\pm(\alpha+2-2) \eta_{jk} x^k
        \right) \\
        & = \eta^{ij} \frac{1}{\alpha}
        \left(
            \frac{\partial}{\partial x^i} R^\pm(\alpha)
        \right)
        \eta_{jk} x^k
        + \eta^{ij} \frac{1}{\alpha} R^\pm(\alpha) \eta_{jk}
        \frac{\partial}{\partial x^i} x^k \\
        & = \frac{1}{\alpha} \eta^{ij} \frac{1}{\alpha-2} \eta_{il} x^l
        R^\pm(\alpha-2) \eta_{jk} x^k
        + \frac{1}{\alpha} R^\pm(\alpha) \eta^{ij} \eta_{jk}
        \delta^k_i \\
        & = \frac{1}{\alpha} \frac{1}{\alpha-2} \eta_{il} x^l x^i
        R^\pm(\alpha-2) + \frac{n}{\alpha} R^\pm(\alpha) \\
        & = \frac{1}{\alpha (\alpha-2)} \eta \cdot R^\pm(\alpha-2)
        + \frac{n}{\alpha} R^\pm(\alpha) \\
        &\stackrel{\mathclap{\refitem{item:riesz-multiplication-by-eta-identity}}}{=} \;
        \frac{(\alpha-2) (\alpha-2-n+2)}{\alpha (\alpha-2)}
        R^\pm(\alpha)
        + \frac{n}{\alpha} R^\pm(\alpha) \\
        & = \frac{\alpha-n+n}{\alpha} R^\pm(\alpha)
        = R^\pm(\alpha).
    \end{align*}
\end{proof}

The above relations hold in the ``strong sense'', i.e. they are
equalities of continuous or even $\Fun$-functions valid point by
point. Since $\Fun(\mathbb{R}^n) \hookrightarrow
\mathcal{D}'(\mathbb{R}^n)$ is injectively embedded via
\eqref{eq:riesz-distribution} we conclude that the above relations
also hold in the sense of distributions. This gives us now the idea
how one can \emph{define} $R^\pm(\alpha)$ for arbitrary $\alpha \in
\mathbb{C}$ at least in the sense of distributions. On one hand, we
want to obtain a holomorphic family of distributions $R^\pm(\alpha)$
for all $\alpha \in \mathbb{C}$ extending the already given ones as in
Lemma~\ref{lemma:ries-distribution-is-homolorphic},
\refitem{item:ries-distribution-holomorphic}. Since a holomorphic
function is already determined by its values on the non-empty open
half space of $\RE(\alpha) > n$, such an extension is necessarily
unique if it exists at all. On the other hand, we can make use of the
relations in Lemma~\ref{lemma:identities-for-riesz-distribution}, in
particular the one in \refitem{item:dAlembert-of-riesz}, to define
such an extension. Indeed, we can express $R^\pm(\alpha)$ as the
d'Alembert operator acting on $R^\pm(\alpha+2)$ for $\RE(\alpha) >
n+2$. Now if $\RE(\alpha) > n$ we \emph{define} $R^\pm(\alpha)$ as
\emph{distribution} by
\begin{equation}
    \label{eq:riesz-distribution-new-definition}
    R^\pm(\alpha) = \dAlembert R^\pm(\alpha+2).
\end{equation}
Since $\alpha \mapsto R^\pm(\alpha)$ is a holomorphic family of
distributions for $\RE(\alpha) > n$ by
Lemma~\ref{lemma:ries-distribution-is-homolorphic}
\refitem{item:ries-distribution-holomorphic} the definition
\eqref{eq:riesz-distribution-new-definition} and the previous
Definition~\ref{definition:riesz-function} coincide as they coincide
for $\RE(\alpha) > n+2$ by
Lemma~\ref{lemma:identities-for-riesz-distribution},
\refitem{item:dAlembert-of-riesz}. Thus we can define inductively for
$\RE(\alpha + 2k) > n$
\begin{equation}
    \label{eq:riesz-inductive-definition}
    R^\pm(\alpha) = \dAlembert R^\pm(\alpha+2k)
\end{equation}
for $k \in \mathbb{N}$. We need the following Lemma:
\begin{lemma}
    \label{lemma:riesz-holomorphic-extension}
    Let $\alpha \in \mathbb{C}$ and define $R^\pm(\alpha)$ by
    \begin{equation}
        \label{eq:riesz-inductive-k}
        R^\pm(\alpha) = \dAlembert {}^k R^\pm(\alpha+2k),
    \end{equation}
    where $k \in \mathbb{N}_0$ is such that $\RE(\alpha+2k) > n$. Then
    \eqref{eq:riesz-inductive-k} does not depend on the choice of $k$
    and yields an entirely holomorphic family of distributions which
    extends the family $\{ R^\pm(\alpha) \}_{\RE(\alpha) > n}$.
\end{lemma}
\begin{proof}
    First we note that \eqref{eq:riesz-inductive-k} yields a
    well-defined distribution as $R^\pm(\alpha+2k)$ is even a
    continuous function for all $k \in \mathbb{N}_0$ with
    $\RE(\alpha+2k) > n$ and derivatives of distributions yield
    distributions. Thus $R^\pm(\alpha) \in \mathcal{D}'(\mathbb{R}^n)$
    is well-defined. If $k' \in \mathbb{N}_0$ is another number with
    $\RE(\alpha+2k') > n$, say $k' > k$, then $\dAlembert^k
    R^\pm(\alpha+2k) = \dAlembert^{k'} R^\pm(\alpha+2k')$ since by
    Lemma~\ref{lemma:identities-for-riesz-distribution}
    \refitem{item:dAlembert-of-riesz} we have $R^\pm(\alpha+2k) =
    \dAlembert^{k'-k} R^\pm(\alpha+2k')$. This shows that
    \eqref{eq:riesz-inductive-k} does not depend on $k$. In
    particular, if already $\RE(\alpha) > n$ then $k = 0$ would
    suffice and $R^\pm(\alpha)$ coincides with the previous definition
    in this case. Thus \eqref{eq:riesz-inductive-k} extends our
    previous definition. Finally, let $\varphi \in
    \Cinfty_0(\mathbb{R}^n)$ be a test function, then
    \[
    R^\pm(\alpha)(\varphi)
    = \left(
        \dAlembert {}^k R^\pm(\alpha+2k)
    \right)
    (\varphi)
    = R^\pm(\alpha+2k) (\dAlembert {}^k \varphi)
    \]
    depends holomorphically on $\alpha$ since $\dAlembert^k \varphi
    \in \Cinfty_0(\mathbb{R}^n)$ is again a test function and
    $R^\pm(\alpha+2k)$ depends holomorphically on $\alpha$ by
    Lemma~\ref{lemma:ries-distribution-is-homolorphic}
    \refitem{item:ries-distribution-holomorphic} in the distributional
    sense. Thus \eqref{eq:riesz-inductive-k} is a holomorphic
    extension of our previous definition.
\end{proof}
\begin{corollary}
    \label{corollary:riesz-holomorphic-extension-is-unique}
    The family $\{ R^\pm(\alpha) \}_{\alpha \in \mathbb{C}}$ of
    distributions as in \eqref{eq:riesz-inductive-k} is the unique
    holomorphic family of distributions extending the family from
    \eqref{eq:riesz-distribution-holomorphic}.
\end{corollary}
After these preparations we are now in the position to state the main
definition of this section:
\begin{definition}[Riesz distributions]
    \label{definition:riesz-distributions}
    \index{Riesz distribution!advanced}%
    \index{Riesz distribution!retarded}%
    For $\alpha \in \mathbb{C}$ the distributions $R^+(\alpha)$ are
    called the advanced Riesz distributions and the $R^-(\alpha)$ are
    called the retarded Riesz distributions.
\end{definition}

%
%

\subsection{Properties of the Riesz Distributions}
\label{subsec:properties-of-riesz-distributions}

Having a definition of $R^\pm(\alpha)$ for all complex numbers $\alpha
\in \mathbb{C}$ we can start to collect some properties of the Riesz
distributions. In particular, they will turn out to provide Green
functions for $\dAlembert$ on Minkowski spacetime.  We start with the
following observation:
\begin{proposition}
    \label{proposition:riesz-distribution-identities}
    \index{Lorentz transformation!orthochronous}%
    \index{Time-reversal}%
    Let $\alpha \in \mathbb{C}$. Then we have:
    \begin{propositionlist}
    \item \label{item:riesz-distribution-and-lorentz-trafo} For all
        orthochronous Lorentz transformations $\Lambda \in
        \group{L}^\uparrow(1, n-1)$ we have
        \begin{equation}
            \label{eq:riesz-dist-invariant-under-orthochronous}
            \Lambda^* R^\pm(\alpha) = R^\pm(\alpha),
        \end{equation}
        and for the time-reversal $\mathrm{T} \in \group{L}(1, n-1)$
        we have
        \begin{equation}
            \label{eq:riesz-dist-and-time-reversal}
            \mathrm{T}^* R^\pm(\alpha) = R^\mp(\alpha).
        \end{equation}
    \item \label{item:riesz-distribution-times-eta} One has
        \begin{equation}
            \label{eq:riesz-distribution-times-eta}
            \eta R^\pm(\alpha) = \alpha (\alpha-n+2) R^\pm(\alpha+2).
        \end{equation}
    \item \label{item:riesz-distribution-differentiation} For all $i=
        1, \ldots, n$ one has
        \begin{equation}
            \label{eq:riesz-distribution-differentiation}
            (\alpha-2) \frac{\partial}{\partial x^i} R^\pm(\alpha)
            = R^\pm(\alpha-2) \eta_{ij} x^j.
        \end{equation}
    \item \label{item:riesz-distribution-is-homogeneous} Let $\lambda
        > 0$. Then for all $\varphi \in \Cinfty_0(\mathbb{R}^n)$ one
        has
        \begin{equation}
            \label{eq:riesz-distribution-homogeneous}
            \lambda^{\alpha-n} R^\pm(\alpha) (\varphi_\lambda)
              = R^\pm(\alpha) (\varphi),
        \end{equation}
        where $\varphi_\lambda(x) = \lambda^n \varphi(\lambda x)$.
        Infinitesimally, this means for the Lie derivative with
        respect to the \Index{Euler vector field}
        \begin{equation}
            \label{eq:riesz-distribution-homogeneous-infintesimal}
            \Lie_\xi R^\pm(\alpha) = (\alpha-n) R^\pm(\alpha),
        \end{equation}
        i.e. $R^\pm(\alpha)$ is homogeneous of degree $\alpha-n$.
    \item \label{item:riesz-distribution-dAlembert} One has
        \begin{equation}
            \label{eq:riesz-distribution-gradient}
            \gradient \eta \cdot R^\pm(\alpha)
            = 2 \alpha \gradient R^\pm(\alpha+2)
        \end{equation}
        and
        \begin{equation}
            \label{eq:riesz-distribution-dAlembert}
            \dAlembert R^\pm(\alpha+2) = R^\pm(\alpha).
        \end{equation}
    \end{propositionlist}
\end{proposition}
\begin{proof}
    For the first part we first note that the Jacobi determinant of
    the diffeomorphism $x \mapsto \Lambda x$ is $\pm 1$ for $\Lambda
    \in \group{L}^\uparrow(1, n-1)$ whence it preserves the Lorentz
    volume density $|\D x^1 \wedge \cdots \wedge \D x^n| = \D^n
    x$. Thus the general definition of $\Lambda^* R^\pm(\alpha)$
    simplifies in this case and is compatible with
    \eqref{eq:riesz-distribution} for $\RE(\alpha) > n$. In fact, we
    have for $\RE(\alpha) > n$ and $\varphi \in
    \Cinfty_0(\mathbb{R}^n)$
    \begin{align*}
        \int \varphi(x) (\Lambda^* R^\pm(\alpha))(x) \D^n x
        & = \int \varphi(x) R^\pm(\alpha)(\Lambda x) \D^n x \\
        & = \int \varphi(\Lambda^{-1}y) R^\pm(\alpha)(y) \D^n y \\
        & = \int (\Lambda_* \varphi)(y) R^\pm(\alpha)(y) \D^n y.
    \end{align*}
    Since the continuous function $R^\pm(\alpha)$ for $\RE(\alpha) >
    n$ is $\group{L}^\uparrow(1,n-1)$-invariant by
    Lemma~\ref{lemma:riesz-functions-lorenz-trafos}, and since
    \[
    \alpha \; \mapsto \;
    (\Lambda^* R^\pm(\alpha))(\varphi)
    = R^\pm(\alpha) (\Lambda \varphi)
    \]
    as well as $\alpha \mapsto R^\pm(\alpha)(\varphi)$ are both
    holomorphic for all $\alpha \in \mathbb{C}$, these holomorphic
    functions coincide for all $\alpha \in \mathbb{C}$.  The second
    and third part follow by the same arguments as both sides are
    holomorphic functions of $\alpha$ when evaluated on $\varphi \in
    \Cinfty_0(\mathbb{R}^n)$ and they coincide for $\RE(\alpha)$
    sufficiently large by
    Lemma~\ref{lemma:ries-distribution-is-homolorphic}.  Now let
    $\lambda > 0$. Then $\alpha \mapsto \lambda^\alpha$ is holomorphic
    on $\mathbb{C}$ and thus $\alpha \mapsto \lambda^\alpha
    R^\pm(\alpha)(\varphi_\lambda)$ is holomorphic on $\mathbb{C}$ for
    any fixed $\varphi \in \Cinfty_0(\mathbb{R}^n)$. Thus we have to
    show \eqref{eq:riesz-distribution-homogeneous} only for
    sufficiently large $\RE(\alpha)$ in order to apply the uniqueness
    arguments. But for $\RE(\alpha) > n$ we have
    \begin{align*}
        R^\pm(\alpha)(\lambda x)
        & = \begin{cases}
            c(\alpha, n) \eta(\lambda x)^{\frac{\alpha-n}{2}}
            \qquad & x \in I^\pm(0) \\
            0 & \mathrm{else}
        \end{cases} \\
        & = \begin{cases}
            c(\alpha, n) (\lambda^2)^{\frac{\alpha-n}{2}}
            \eta(x)^{\frac{\alpha-n}{2}} \qquad & x \in I^\pm(0) \\
            0 & \mathrm{else}
        \end{cases} \\
        & = \lambda^{\alpha-n} R^\pm(\alpha)(x)
    \end{align*}
    for all $x \in \mathbb{R}^n$. Then, in the sense of
    distributions,
    \begin{align*}
        \lambda^{\alpha-n} R^\pm(\alpha) (\varphi_\lambda)
        & = \lambda^{\alpha-n} \int R^\pm(\alpha)(x) \lambda^n
        \varphi(\lambda x) \D^n x \\
        & = \int R^\pm(\alpha)(\lambda x) \varphi(\lambda x) \lambda^n
        \D^n x \\
        & = \int R^\pm(\alpha)(y) \varphi(y) \D^n y \\
        & = R^\pm(\alpha)(\varphi).
    \end{align*}
    Thus we conclude that \eqref{eq:riesz-distribution-homogeneous}
    holds for all $\alpha \in \mathbb{C}$. To prove the infinitesimal
    version \eqref{eq:riesz-distribution-homogeneous-infintesimal} one
    can either use \eqref{eq:riesz-distribution-differentiation} and
    \eqref{eq:riesz-distribution-times-eta} or differentiate
    \eqref{eq:riesz-distribution-homogeneous}: Indeed, since
    $(\lambda, x) \mapsto \lambda^n \varphi(\lambda x)$ is smooth and
    compactly supported in $x$ ``locally uniform in $\lambda$'', a
    slight variation of
    Lemma~\ref{lemma:parameter-differentiation-under-integral} shows
    that $\lambda \mapsto \lambda^{\alpha-n}
    R^\pm(\alpha)(\varphi_\lambda)$ is smooth in $\lambda$ and the
    derivatives can be computed by differentiating ``under the
    integral sign'' as in
    Lemma~\ref{lemma:parameter-differentiation-under-integral}. We
    find
    \begin{align*}
        \frac{\partial}{\partial \lambda}
        \left(
            \lambda^{\alpha-n} R^\pm(\alpha) (\varphi_\lambda)
        \right)
        & =
        (\alpha-n) \lambda^{\alpha-n-1} R^\pm(\alpha)
        (\varphi_\lambda)
        + \lambda^{\alpha-n} R^\pm(\alpha)
        \left(
            \frac{\partial}{\partial \lambda}
            (x \; \mapsto \; \lambda^n \varphi(\lambda x))
        \right) \\
        & =
        (\alpha-n) \lambda^{\alpha-n-1} R^\pm(\alpha)
        (\varphi_\lambda)
        + \lambda^{\alpha-n} \lambda^n R^\pm(\alpha)
        \left(
            \frac{\partial \varphi}{\partial x^i} (\lambda x) x^i
        \right) \\
        & \quad + \lambda^{\alpha-n} n \lambda^{n-1} R^\pm(\alpha)
        (x \; \mapsto \; \varphi(\lambda x)).
    \end{align*}
    Since the left hand side does \emph{not} depend on $\lambda$, this
    has to vanish for all $\lambda > 0$. Setting $\lambda = 1$ yields
    \begin{align*}
       0 & = (\alpha-n) R^\pm(\alpha)(\varphi) + R^\pm(\alpha)
       \left(
       x^i \frac{\partial \varphi}{\partial x^i}
       \right)
       + n R^\pm(\alpha)(\varphi) \\
       & = \alpha R^\pm(\alpha)(\varphi) + R^\pm(\alpha)
       \left(
           \frac{\partial}{\partial x^i}
           \left(x \; \mapsto \; x^i \varphi(x)\right)
       \right)
       - R^\pm(\alpha)(n \varphi) \\
       & = (\alpha-n) R^\pm(\varphi)
       -
       \frac{\partial R^\pm(\alpha)}{\partial x^i} (x^i \varphi) \\
       & = (\alpha-n) R^\pm(\alpha)(\varphi)
       - (\Lie_\xi R^\pm(\alpha))(\varphi),
   \end{align*}
   and thus \eqref{eq:riesz-distribution-homogeneous-infintesimal}.
   The last part again follows from
   Lemma~\ref{lemma:identities-for-riesz-distribution},
   \refitem{item:riesz-and-gradient} and
   \refitem{item:dAlembert-of-riesz} as well as the uniqueness
   argument: clearly both sides evaluated on a test function give
   holomorphic functions of $\alpha$ which coincide for large
   $\RE(\alpha)$.
\end{proof}
\begin{remark}[Homogeneous distributions]
    \label{remark:homogeneous-distributions}
    \index{Distribution!homogeneous}%
    In general, a distribution $u \in \mathcal{D}'(\mathbb{R}^n
    \setminus \{0\})$ is called \emph{homogeneous of degree $\alpha
      \in \mathbb{C}$} if for all test functions $\varphi \in
    \Cinfty_0(\mathbb{R}^n \setminus \{0\})$ one has
    \begin{equation}
        \label{eq:homogeneous-distribution}
        \lambda^\alpha u(\varphi_\lambda) = u( \varphi)
    \end{equation}
    for all $\lambda > 0$, where $\varphi_\lambda(x) = \lambda^n
    \varphi(\lambda x)$ as before. By the same argument as in the
    proof one can show that \eqref{eq:homogeneous-distribution}
    implies
    \begin{equation}
        \label{eq:homogeneous-distribution-infinitesimal}
        \Lie_\xi u = \alpha u.
    \end{equation}
    In fact, \eqref{eq:homogeneous-distribution-infinitesimal} turns
    out to be equivalent to its integrated form
    \eqref{eq:homogeneous-distribution}. It is then a non-trivial but
    interesting question whether a homogeneous distribution $u \in
    \mathcal{D}'(\mathbb{R}^n \setminus \{0\})$ of some degree
    $\alpha$ can be extended to a distribution $u \in
    \mathcal{D}'(\mathbb{R}^n)$ such that the homogeneity is
    preserved. A detailed discussion of homogeneous distributions can
    be found in \cite[Sect.~3.2]{hoermander:2003a}. As a final remark
    we mention that many problems in renormalization theory of quantum
    field theories can be reformulated mathematically as the question
    whether certain homogeneous distributions on $\mathbb{R}^n
    \setminus \{0\}$ have homogeneous extensions to $\mathbb{R}^n$,
    see e.g.~\cite{steinmann:2000a, scharf:1995a}.
\end{remark}

In a next step we want to understand the support and singular support
of the Riesz distributions $R^\pm(\alpha)$. Here we can build on the
results from Lemma~\ref{lemma:riesz-function-continuity} and
\ref{lemma:riesz-functions-lorenz-trafos}: the support and singular
support have to be Lorentz invariant subsets under the orthochronous
Lorentz group $\group{L}^\uparrow(1, n-1)$. We denote by
\begin{equation}
    \label{eq:boundary-of-lightcone}
    C^\pm(0)
    = \left\{
        x \in \mathbb{R}^n
        \; \big| \;
        x \in J^\pm(0) \; \textrm{and} \; \eta(x,x) = 0
    \right\}
\end{equation}
the boundary of $I^\pm(0)$. The particular values $\alpha \in
\mathbb{C}$ where $c(\alpha, n)$ vanishes play an exceptional role for
the support of $R^\pm(\alpha)$. We call them \emph{exceptional}, i.e.
$\alpha \in \mathbb{C}$ is exceptional if
\begin{equation}
    \label{eq:exceptional-alphas}
    \alpha \in \{ n-2k, -2k \; | \; k \in \mathbb{N}_0 \}.
\end{equation}
Then we have the following result:
\begin{proposition}[Support of $R^\pm(\alpha)$]
    \label{proposition:support-of-riesz-distribution}
    \index{Riesz distribution!support}%
    \index{Riesz distribution!singular support}%
    Let $\alpha \in \mathbb{C}$.
    \begin{propositionlist}
    \item \label{item:riesz-supp-nonexceptional} If $\alpha$ is not
        exceptional then
        \begin{equation}
            \label{eq:riesz-supp-nonexceptional}
            \supp R^\pm(\alpha) = J^\pm(0),
        \end{equation}
        and the singular support
        \begin{equation}
            \label{eq:riesz-singsupp-nonexceptional}
            \singsupp R^\pm(\alpha) \subseteq
            \partial I^\pm(0) = C^\pm(0)
        \end{equation}
        is either $\{0\}$ or $C^\pm(0)$.
    \item \label{item:riesz-support-exceptional} If $\alpha$ is
        exceptional then
        \begin{equation}
            \label{eq:riesz-supp-and-singsupp-exceptional}
            \supp R^\pm(\alpha) = \singsupp R^\pm(\alpha)
            \subseteq C^\pm(0).
        \end{equation}
    \item \label{item:riesz-supp-in-high-dimension-for-special-alpha}
        Let $n \geq 3$. For $\alpha \in \{ n-2k \; | \; k \in
        \mathbb{N}_0, k < \frac{n}{2} \}$ we have
        \begin{equation}
            \label{eq:riesz-supp-singsupp-higgdim-special-alpha}
            \supp R^\pm(\alpha) = \singsupp R^\pm(\alpha) = C^\pm(0).
        \end{equation}
    \end{propositionlist}
\end{proposition}
\begin{proof}
    Let $\alpha \in \mathbb{C}$ be arbitrary. Since by definition of
    $R^\pm(\alpha)$ we have
    \[
    R^\pm(\alpha) = \dAlembert {}^k R^\pm(\alpha + 2k)
    \]
    for $k$ sufficiently large such that $\RE(\alpha+2k) > n$, we have
    by Theorem~\ref{theorem:differentiation-of-gensec},
    \refitem{item:diffop-for-gensec-is-local}
    \begin{align*}
        R^\pm(\alpha) \At{\mathbb{R}^n \setminus C^\pm(0)}
        & = \dAlembert {}^k
        \left(
            R^\pm(\alpha+2k) \At{\mathbb{R}^n \setminus C^\pm(0)}
        \right) \\
        & = \begin{cases}
            \dAlembert {}^k c(\alpha+2k, n)
            \eta^{\frac{\alpha+2k-n}{2}} & \mathrm{on} \quad I^\pm(0) \\
            0 & \mathrm{else},
        \end{cases} \\
        & = \begin{cases}
            c(\alpha, n) \eta^{\frac{\alpha-n}{2}}
            & \mathrm{on} \quad I^\pm(0) \\
            0 & \mathrm{else},
        \end{cases}
    \end{align*}
    using the explicit computation of $\dAlembert \eta$ as in the
    proof of Lemma~\ref{lemma:riesz-function-continuity}. Thus on the
    open subset $\mathbb{R}^n \setminus C^\pm(0) = I^\pm(0) \cup
    (\mathbb{R}^n \setminus J^\pm(0))$ we have a smooth function
    \[
    R^\pm(\alpha) \At{\mathbb{R}^n \setminus C^\pm(0)}
    = \begin{cases}
        c(\alpha, n) \eta^{\frac{\alpha-n}{2}}
        & \mathrm{on} \quad I^\pm(0) \\
        0 & \mathrm{else},
    \end{cases}
    \]
    for all $\alpha \in \mathbb{C}$.  From this we immediately
    conclude that for all $\alpha \in \mathbb{C}$
    \[
    \supp R^\pm(\alpha) \subseteq J^\pm(0)
    \tag{$*$}
    \]
    and
    \[
    \singsupp R^\pm(\alpha) \subseteq C^\pm(0),
    \tag{$**$}
    \]
    since $J^\pm(0)$ and the light cone $C^\pm(0)$ are already closed.
    Then the Lorentz invariance $\Lambda_* R^\pm(\alpha) =
    R^\pm(\alpha)$ for all $\Lambda \in \group{L}^\uparrow(1, n-1)$
    yields that the support and the singular support have to be
    Lorentz invariant subsets. Indeed, in general one has
    \[
    \supp \Lambda_* R^\pm(\alpha) = \Lambda (\supp R^\pm(\alpha))
    \]
    \[
    \singsupp \Lambda_* R^\pm(\alpha)
    = \Lambda (\singsupp R^\pm(\alpha))
    \]
    for every diffeomorphism $\Lambda$. Thus $\supp R^\pm(\alpha)$ and
    $\singsupp R^\pm(\alpha)$ are closed Lorentz invariant subsets of
    Minkowski space. In particular, $\singsupp R^\pm(\alpha)$ is
    either $\{0\}$ or $C^\pm(0)$ as these are the only Lorentz
    invariant subsets of $C^\pm(0)$. Now let $\alpha$ be not
    exceptional. Then $c(\alpha, n)$ is non-zero and hence
    $R^\pm(\alpha)\at{I^\pm(0)}$ is non-zero and even smooth. Thus
    \[
    \supp R^\pm(\alpha) \supseteq I^\pm (0).
    \]
    On the other hand, by ($*$), we note $\supp R^\pm(\alpha)
    \subseteq J^\pm(0)$, whence
    \eqref{eq:riesz-singsupp-nonexceptional} follows. This shows the
    first part.  For the second part, let $\alpha$ be exceptional.
    Then $c(\alpha, n)= 0$ whence $R^\pm(\alpha)\at{\mathbb{R}^n
      \setminus C^\pm(0)}$ vanishes identically. Thus
    \[
    \supp R^\pm(\alpha) \subseteq C^\pm(0)
    \]
    follows. Now $C^\pm(0)$ has empty open interior whence the support
    of $R^\pm(\alpha)$ is either empty or necessarily entirely
    singular. Thus
    \[
    \supp R^\pm(\alpha) = \singsupp R^\pm(\alpha)
    \]
    follows, proving the second part.  For the last part we follow
    \cite[Prop.~1.2.4.]{baer.ginoux.pfaeffle:2007a} and prove first
    the following technical statement. We consider a test function
    $\psi \in \Cinfty_0(\mathbb{R})$ with $\supp \psi \subseteq [a,
    b]$ and a bump function $\chi \in \Cinfty_0(\mathbb{R}^{n-1})$
    such that $\chi \at{B_r(0)} = 1$ for some $r > b$.
    \begin{figure}
        \centering
        \input{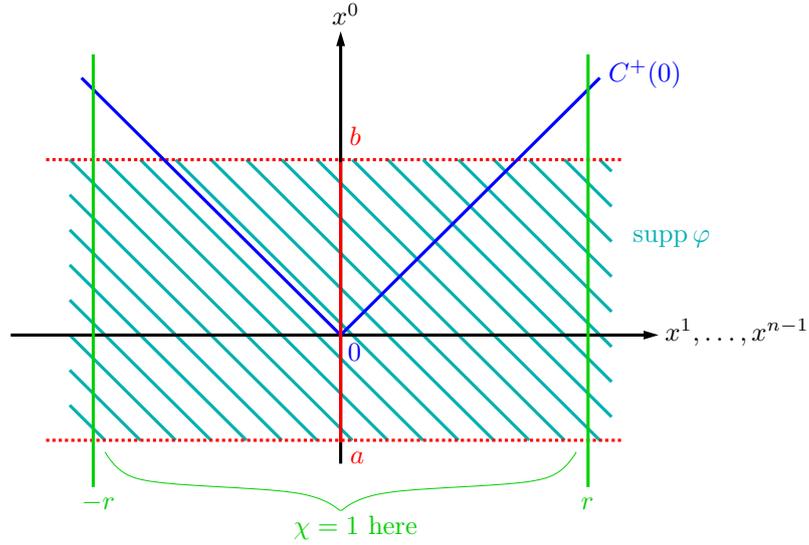}
        \caption{
          \label{fig:special-test-function}%
          The test function constructed in the proof of
          Proposition~\ref{proposition:support-of-riesz-distribution}.
        }
    \end{figure}
    Then the test function $\varphi(x^0, x^1, \ldots, x^{n-1}) =
    \psi(x^0) \chi(x^1, \ldots, x^{n-1})$ has the property that
    \[
    \varphi \At{J^+(0)} (x^0, \ldots, x^{n-1}) = \psi (x^0)
    \tag{$*$$**$}
    \]
    for all $x \in J^+(0)$, see also
    Figure~\ref{fig:special-test-function}. Then the claim is that for
    all $\RE(\alpha) > 0$ one has
    \[
    R^+(\alpha)(\varphi)
    = \frac{1}{\Gamma(\alpha)}
    \int_0^\infty (x^0)^{\alpha-1} \psi(x^0) \D x^0.
    \tag{$**$$**$}
    \]
    Indeed, we first note that both sides are holomorphic in $\alpha$.
    For the left hand side this is true for all $\alpha \in
    \mathbb{C}$ and for the right hand side this follows as
    $\frac{1}{\Gamma(\alpha)}$ is entire and the integral is
    holomorphic by the same Morera type argument as in the proof of
    Lemma~\ref{lemma:ries-distribution-is-homolorphic}. Thus it will
    be sufficient to show ($**$$**$) for $\RE(\alpha) > n$ where we can
    use the explicit form of $R^+(\alpha)$ as continuous function. We
    compute
    \begin{align*}
        R^+(\alpha)(\varphi)
        &= c(\alpha, n) \int_{J^+(0)} \eta(x)^{\frac{\alpha-n}{2}}
        \varphi(x) \D^n x \\
        &= c(\alpha, n) \int_0^\infty \D x^0
        \int_{|\vec{x}| \leq x^0}
        \left(
            (x^0)^2 - (\vec{x})^2
        \right)^{\frac{\alpha-n}{2}} \psi(x^0) \chi(\vec{x})
        \D^{n-1} x \\
        &\stackrel{\mathclap{(*{*}*)}}{=}
        \;\;
        c(\alpha, n)
        \int_0^\infty \D x^0 \psi(x^0)
        \int_{|\vec{x}| \leq x^0}
        \left(
            (x^0)^2 - (\vec{x})^2
        \right)^{\frac{\alpha-n}{2}} \D^{n-1} x.
    \end{align*}
    For the $\vec{x}$-integration we use $(n-1)$-dimensional polar
    coordinates $r$ and $\vec{\Omega}$, i.e. the radius $r =
    |\vec{x}|$ and the remaining point $\vec{\Omega} =
    \frac{\vec{x}}{|\vec{x}|}$ on the unit sphere
    $\mathbb{S}^{n-2}$. We evaluate for fixed $x^0$ the inner
    integral
    \begin{align*}
        \int_{|\vec{x}| \leq x^0}
        \left(
            (x^0)^2 - (\vec{x})^2
        \right)^{\frac{\alpha-n}{2}} \D^{n-1} x
        &=
        \int_0^{x^0} r^{n-2} \D r
        \int_{\mathbb{S}^{n-2}}
        \left(
            (x^0)^2 - r^2
        \right)^{\frac{\alpha-n}{2}} \D \Omega \\
        &=
        \vol(\mathbb{S}^{n-2}) \cdot
        \int_0^{x^0}
        \left(
            (x^0)^2 - r^2
        \right)^{\frac{\alpha-n}{2}}
        r^{n-2} \D r.
    \end{align*}
    The remaining integral can be brought to the following form. First
    we substitute $\rho = \frac{r}{x^0}$ and then $\rho = \cos
    \theta$. This yields
    \begin{align*}
        \int_0^{x^0}
        \left(
            (x^0)^2 - r^n
        \right)^{\frac{\alpha-n}{2}}
        r^{n-2} \D r
        & =
        \int_0^1 (1-\rho^2)^{\frac{\alpha-n}{2}} \rho^{n-2}
        (x^0)^{\alpha-n+n-2+1} \D \rho \\
        & =
        (x^0)^{\alpha-1}
        \int_0^1 (1-\rho^2)^{\frac{\alpha-n}{2}} \rho^{n-2} \D \rho \\
        & = (x^0)^{\alpha-1}
        \int_0^{\pi / 2}
        (\sin^2 \theta)^{\frac{\alpha-n}{2}} (\cos \theta)^{n-2}
        \sin \D \theta \\
        & =(x^0)^{\alpha-1} \int_0^{\pi / 2}
        (\sin \theta)^{\alpha-n+1} (\cos \theta)^{n-2} \D \theta.
    \end{align*}
    The last integral is \emph{Bronstein-integrable}, see
    e.g. \cite[Sect.~1.1.3.4,
    Integral~10]{bronstein.semendjajew:1989a} and gives
    \begin{align*}
        (x^0)^{\alpha-1} \int_0^{\pi / 2}
        (\sin \theta)^{\alpha-n+1} (\cos \theta)^{n-2} \D \theta
        =
        \frac{1}{2}
        \frac{\Gamma(\frac{\alpha-n}{2}+1) \Gamma(\frac{n-3}{2}+1)}
        {\Gamma(\frac{\alpha-n}{2} + \frac{n-3}{2} - 2)}
        =
        \frac{1}{2}
        \frac{\Gamma(\frac{\alpha-n}{2}+1) \Gamma(\frac{n-3}{2}+1)}
        {\Gamma(\frac{\alpha-1}{2} + 1)}.
    \end{align*}
    Since finally the surface of the $(n-2)$-dimensional unit sphere
    is known to be
    \[
    \vol(\mathbb{S}^{n-2})
    = \frac{2 \pi^{\frac{n-1}{2}}}{\Gamma(\frac{n-1}{2})},
    \]
    see e.g. \cite[p.~142]{gallot.hulin.lafontaine:1990a}, we obtain
    \begin{align*}
        R^+(\alpha)(\varphi)
        & =
        c(\alpha, n)
        \cdot
        \frac{2 \pi^{\frac{n-1}{2}}}{\Gamma(\frac{n-1}{2})}
        \cdot
        \frac{1}{2}
        \frac{\Gamma(\frac{\alpha-n}{2}+1) \Gamma(\frac{n-3}{2}+1)}
        {\Gamma(\frac{\alpha-1}{2} + 1)}
        \cdot \int_0^\infty (x^0)^{\alpha-1} \psi(x^0) \D x^0 \\
        & =
        \frac{2^{1-\alpha} \pi^{\frac{2-n}{2}}}
        {\Gamma(\frac{\alpha}{2}) \Gamma(\frac{\alpha-n}{2}+1)}
        \frac{2 \pi^{\frac{n-1}{2}}}{\Gamma(\frac{n-1}{2})}
        \cdot
        \frac{1}{2}
        \frac{\Gamma(\frac{\alpha-n}{2}+1) \Gamma(\frac{n-3}{2}+1)}
        {\Gamma(\frac{\alpha-1}{2} + 1)}
        \cdot \int_0^\infty (x^0)^{\alpha-1} \psi(x^0) \D x^0 \\
        & =
        \frac{2^{1-\alpha} \sqrt{\pi}}
        {\Gamma(\frac{\alpha}{2})
          \Gamma(\frac{\alpha}{2}+\frac{1}{2})}
        \cdot \int_0^\infty (x^0)^{\alpha-1} \psi(x^0) \D x^0 \\
        & =
        \frac{1}{\Gamma(\alpha)}
        \cdot \int_0^\infty (x^0)^{\alpha-1} \psi(x^0) \D x^0 ,
    \end{align*}
    \index{Legendre's duplication formula}%
    where the last equality is valid thanks to Legendre's duplication
    formula~\eqref{eq:gamma-legendre-duplication}. This finally
    establishes the claim ($*$$*$$*$$*$).  In particular, for $\alpha
    = 2$ we obtain
    \begin{align*}
        R^+(2)(\varphi)
        = \frac{1}{\Gamma(2)}
        \int_0^\infty (x^0)^{\alpha-1} \psi(x^0) \D x^0
        = \int_0^\infty (x^0)^{\alpha-1} \psi(x^0) \D x^0,
    \end{align*}
    from which it follows that the support of $R^+(2)$ cannot be $0
    \in \mathbb{R}^n$ alone as we get a non-trivial result for a
    $\varphi$ with $0 \notin \supp \varphi$ by taking a $\psi$ with
    support away from zero. Thus by the previous arguments the support
    is at least $C^+(0)$. So if $n$ is even then $2$ is an exceptional
    value whence $\supp R^+(2) = \singsupp R^+(2) \subseteq C^+(0)$
    and thus
    \[
    \supp R^+(2) = \singsupp R^+(2) = C^+(0)
    \]
    follows. Since in this case also $2, 4, \ldots, n-2, n$ are
    exceptional and
    \[
    R^+(2) = \dAlembert {}^k R^+(2+2k)
    \]
    for all $k \in \mathbb{N}_0$ we conclude from the locality
    \eqref{eq:supp-of-Ds} of differential operators by
    Theorem~\ref{theorem:differentiation-of-gensec} that
    \[
    C^+(0) = \supp R^+(2) \subseteq \supp R^+(2+2k) \subseteq C^+(0)
    \]
    for all those $k$ with $2+2k \leq n$. But then again $\supp
    R^+(2+2k) = C^+(0)$ follows. Now let $n$ be odd. Since
    $R^+(\alpha)(\varphi)$ is holomorphic for all $\alpha$ and since
    the limit $\alpha \longrightarrow 1$ of ($**$$**$) exists, we
    conclude
    \[
    R^+(1)(\varphi) = \int_0^\infty \psi(x^0) \D x^0,
    \]
    whence the support of $R^+(1)$ is again not only $\{0\}$. Thus we
    can repeat the argument with $R^+(1)$ instead of $R^+(2)$ and
    obtain \eqref{eq:riesz-supp-singsupp-higgdim-special-alpha} also
    in this case. Of course the result for $R^-(\alpha)$ is completely
    analogous or can be deduced from the time reversal symmetry
    \eqref{eq:riesz-dist-and-time-reversal}.
\end{proof}

The following counting of the order of the Riesz distributions
$R^\pm(\alpha)$ is straightforward:
\begin{proposition}[Order of $R^\pm(\alpha)$]
    \label{proposition:order-of-riesz-distribution}
    \index{Riesz distribution!order}%
    Let $\alpha \in \mathbb{C}$.
    \begin{propositionlist}
    \item \label{item:riesz-order-for-big-alpha} If $\RE(\alpha) > n$
        then the global order of $R^\pm(\alpha)$ is zero
        \begin{equation}
            \label{eq:riesz-order-for-big-alpha}
            \ord (R^\pm(\alpha)) = 0.
        \end{equation}
    \item \label{item:riesz-global-order} The global order of
        $R^\pm(\alpha)$ is bounded by $2k$ where $k \in \mathbb{N}_0$
        is such that $\RE(\alpha) + 2k > n$.
    \item \label{item:riesz-order-if-alpha-nonnegative} If
        $\RE(\alpha) > 0$ then the global order of $R^\pm(\alpha)$ is
        bounded by $n$ if $n$ is even and by $n+1$ if $n$ is odd.
    \end{propositionlist}
\end{proposition}
\begin{proof}
    The first part is clear since for $\RE(\alpha) > n$ the
    distribution $R^\pm(\alpha)$ is even a continuous function. For
    the second part let $k \in \mathbb{N}_0$ be such that
    $\RE(\alpha)+2k > n$. Then
    \begin{align*}
        \ord( R^\pm(\alpha) )
        = \ord( \dAlembert {}^k R^\pm(\alpha+2k) )
        \leq \ord( R^\pm(\alpha+2k) ) +2k
        = 0 + 2k,
    \end{align*}
    since $\ord(R^\pm(\alpha+2k))=0$ by the first part. Finally, let
    $\RE(\alpha) > 0$ and $n=2k$ be even. Then by the second part
    $\ord(R^\pm(\alpha)) \leq 2k = n$ since $\RE(\alpha)+n > n$. If on
    the other hand $n=2k+1$ is odd then by the second part
    $\ord(R^\pm(\alpha)) \leq 2(k+1) = n+1$ since $\RE(\alpha)+2(k+1)
    > n$.
\end{proof}

The next statement is on the reality of $R^\pm(\alpha)$ for real
$\alpha \in \mathbb{R}$. In fact, one has the following statement:
\begin{proposition}[Reality of $R^\pm(\alpha)$]
    \label{proposition:realitiy-of-riesz-distribution}
    \index{Riesz distribution!reality}%
    Let $\alpha \in \mathbb{C}$. Then one has
    \begin{equation}
        \label{eq:conjugation-of-riesz-distribution}
        \cc{R^\pm(\alpha)} = R^\pm(\cc{\alpha}).
    \end{equation}
    In particular, for $\alpha \in \mathbb{R}$ one has
    \begin{equation}
        \label{eq:riesz-real-for-alpha-real}
        \cc{R^\pm(\alpha)} = R^\pm(\alpha).
    \end{equation}
\end{proposition}
\begin{proof}
    First we consider $\RE(\alpha) > n$. Then we have
    \begin{align*}
        \cc{R^\pm(\alpha)(x)}
        & = \begin{cases}
            \cc{c(\alpha,n)} \cc{\eta(x)^{\frac{\alpha-n}{2}}}
            & \textrm{for} \; x \in I^\pm(0) \\
            0 & \textrm{else}
        \end{cases} \\
        & = \begin{cases}
            c(\cc{\alpha},n) \eta(x)^{\frac{\cc{\alpha}-n}{2}}
            & \textrm{for} \; x \in I^\pm(0) \\
            0 & \textrm{else}
        \end{cases} \\
        & = R^\pm(\cc{\alpha})(x)
    \end{align*}
    for all $x \in \mathbb{R}^n$ since $\cc{\Gamma(\alpha)} =
    \Gamma(\cc{\alpha})$ and hence $\cc{c(\alpha, n)} = c(\cc{\alpha},
    n)$. For arbitrary $\alpha \in \mathbb{C}$ let $k \in
    \mathbb{N}_0$ be such that $\RE(\alpha)+2k > n$. Then
    \begin{align*}
        \cc{R^\pm(\alpha)}
        = \cc{\dAlembert {}^k R^\pm(\alpha+2k)}
        = \dAlembert {}^k \cc{R^\pm(\alpha+2k)}
        = \dAlembert {}^k R^\pm(\cc{\alpha}+2k)
        = R^\pm(\cc{\alpha}),
    \end{align*}
    since $\dAlembert$ is a real differential operator and
    $\cc{R^\pm(\alpha+2k)} = R^\pm(\cc{\alpha}+2k)$ for $\RE(\alpha +
    2k) > n$.
\end{proof}

The next statement is the key observation why the Riesz distributions
are actually what we are looking for.
\begin{proposition}
    \label{proposition:riesz-zero-is-delta}
    One has
    \begin{equation}
        \label{eq:ries-zero-is-delta}
        R^\pm(0) = \delta_0.
    \end{equation}
\end{proposition}
\begin{proof}
    We have to compute $R^\pm(0)(\varphi)$ for $\varphi \in
    \Cinfty_0(\mathbb{R}^n)$. Let $\varphi \in
    \Cinfty_K(\mathbb{R}^n)$ with some compact $K \subseteq
    \mathbb{R}^n$ and choose $\chi \in \Cinfty_0(\mathbb{R}^n)$ with
    $\chi \at{K} =  1$. Then we have $\varphi = \chi
    \varphi$. Moreover, by the usual Hadamard trick we have smooth
    functions $\varphi_i \in \Cinfty(\mathbb{R}^n)$ such that
    \[
    \varphi(x) = \varphi(0) + \sum\nolimits_{i=1}^n x^i \varphi_i(x).
    \]
    In fact,
    \[
    \varphi_i(x)
    = \int_0^1 \frac{\partial \phi}{\partial x^i} (tx) \D t
    \]
    will do the job. Note that $\supp \varphi_i$ is \emph{not}
    compact. In any case, we have
    \[
    \varphi = \chi \varphi
    = \chi \varphi(0) + \sum\nolimits_{i=1}^n x^i \chi \varphi_i
    \]
    with compactly supported $\chi \varphi(0)$ and $x^i \chi
    \varphi_i$. Only now we can apply the distribution $R^\pm(0)$ to
    both terms giving
    \begin{align*}
        R^\pm(0)(\varphi)
        = R^\pm(0)
        \left(\chi \varphi(0) + \sum_i x^i \chi \varphi_i\right)
        = \varphi(0) R^\pm(0)(\chi)
        + \sum_i\left(x^i R^\pm(0)\right) (\chi \varphi_i).
    \end{align*}
    Now $2 x^i$ is the $i$-th component of $\gradient \eta$ whence by
    Proposition~\ref{proposition:riesz-distribution-identities},
    \refitem{item:riesz-distribution-dAlembert} for $\alpha = 0$ we
    obtain
    \[
    2 \left(x^i R^\pm(0)\right)(\chi \varphi_i)
    =
    2 \cdot 0 \cdot \eta^{ij} \frac{\partial R^\pm(2)}{\partial x^j}
    (\chi \varphi_i)
    = 0.
    \]
    This shows
    \[
    R^\pm(0)(\varphi) = \varphi(0) R^\pm(0) (\chi).
    \]
    Since $R^\pm(0)(\varphi)$ does not depend on the choice of the
    cut-off function $\chi$ the constant $R^\pm(0)(\chi)$ does
    neither. However, it might still depend on the chosen compactum
    $K$ which is easy to see to be not the case. This shows that
    \[
    R^\pm(0) = R^\pm(0)(\chi) \delta_0 = c \cdot \delta_0
    \]
    is a multiple of the $\delta$-functional at zero. We are left with
    the computation of $c = R^\pm(0)(\chi)$. To this end it is
    obviously sufficient to compute $R^\pm(0)(\varphi)$ for one
    function with $\varphi(0) \neq 0$. Thus we again use a factorizing
    function
    \[
    \varphi(x) = \psi(x^0) \chi(x^1, \ldots, x^{n-1})
    \]
    with $\psi \in \Cinfty_0(\mathbb{R})$ and $\chi \in
    \Cinfty_0(\mathbb{R}^{n-1})$ such that $\chi$ is equal to $1$ on a
    large enough ball around $0$ in order to have
    \[
    \varphi \At{J^+(0)} (x^0, \ldots, x^{n-1}) = \psi(x^0).
    \]
    Recall that we constructed such a function in the proof of
    Proposition~\ref{proposition:support-of-riesz-distribution},
    \refitem{item:riesz-supp-in-high-dimension-for-special-alpha}. Then
    we have
    \[
    \dAlembert \varphi \At{J^+(0)} (x) = \ddot{\psi}(x^0),
    \]
    since the $x^1$-, \ldots, $x^{n-1}$-derivatives in $\dAlembert$ do
    not contribute. From the above proof we know that
    \begin{align*}
        R^\pm(0)(\varphi)
        = R^\pm(2) (\dAlembert \varphi)
        = \int_0^\infty x^0 \ddot{\psi}(x^0) \D x^0
        = - \int_0^\infty \dot{\psi}(x^0) \D x^0
        = \psi(0)
        = \varphi(0),
    \end{align*}
    by integration by parts and using that $\psi$ has compact support.
    Thus $R^\pm(0)(\varphi) = \varphi(0)$ whence the multiple is $1$
    and the proof is finished for dimensions $n \ge 3$. The two
    remaining cases $n = 1, 2$ are indeed much simpler. Either, one
    can modify the above argument to work also in this simpler
    situation. Or, as we shall do in
    Subsection~\ref{subsec:riesz-distribution-for-1dim-and-2dim}, one
    uses a direct computation.
\end{proof}

The last proposition allows us to formulate the following main result
of this subsection: we have found the advanced and retarded Green
functions of the scalar wave equation on Minkowski spacetime.
\begin{theorem}[Green function of $\dAlembert$]
    \label{theorem:green-function-of-dAlembert}
    \index{Green function!dAlembertian@d'Alembertian}%
    \index{Wave equation!scalar}%
    \index{Minkowski spacetime}%
    The Riesz distributions $R^\pm(2)$ are advanced and retarded Green
    functions for the scalar d'Alembert operator $\dAlembert$ on
    Minkowski spacetime.
\end{theorem}
\begin{proof}
    First we know by
    Proposition~\ref{proposition:riesz-distribution-identities},
    \refitem{item:riesz-distribution-dAlembert} that $\dAlembert
    R^\pm(2) = R^\pm(0)$ which is $\delta_0$ by
    Proposition~\ref{proposition:riesz-zero-is-delta}. Thus the
    $R^\pm(2)$ are fundamental solutions of $\dAlembert$. Moreover, by
    Proposition~\ref{proposition:support-of-riesz-distribution} we
    know that $\supp R^\pm(2) \subseteq J^\pm(0)$ whence we indeed
    have advanced and retarded Green functions.
\end{proof}

\begin{remark}
    \label{remark:riesz-on-Ck-functions}
    For the later use we mention that for $\varphi \in
    \Fun_0(\mathbb{R}^n)$ the distribution $R^\pm(\alpha)$ can still
    be applied to $\varphi$ as long as $\ord( R^\pm(\alpha) ) \leq k$
    by Remark~\ref{remark:order-of-generalized-section}. This is the
    case for
    \begin{equation}
        \label{eq:condition-for-riesz-on-Ck}
        \RE(\alpha) > n - 2 \cdot \left[ \frac{k}{2} \right]
    \end{equation}
    by Proposition~\ref{proposition:order-of-riesz-distribution},
    \refitem{item:riesz-global-order}. In this case
    $R^\pm(\alpha)(\varphi) = R^\pm(\alpha+2\ell) (\dAlembert {}^\ell
    \varphi)$ for $2 \ell \leq k$ and since $R^\pm(\alpha)(\varphi)$
    is still holomorphic for $\RE(\alpha) > n$ and $\varphi \in
    \Fun[0]_0(\mathbb{R}^n)$ we obtain the result that
    $R^\pm(\alpha)(\varphi)$ is holomorphic for $\RE(\alpha) > n - 2
    \cdot \left[ \frac{k}{2} \right]$ and $\varphi \in
    \Fun_0(\mathbb{R}^n)$.
\end{remark}

%
%

\subsection{The Riesz Distributions in Dimension $n = 1, 2$}
\label{subsec:riesz-distribution-for-1dim-and-2dim}

In this small section we compute the Riesz distributions
$R^\pm(\alpha)$ and in particular $R^\pm(2)$ for low dimensions
explicitly.

We start with the most trivial case $n = 1$. In this case
$\mathbb{R}^1$ is equipped with the \emph{Riemannian} metric $\eta =
\D t^2$, where we denote the canonical coordinate simply by $t$.
Though we do not even have an honest Lorentz spacetime in this case
the results from the preceding sections are nevertheless valid.

In this case, the advanced and retarded Green functions $R^\pm(2)$ are
even defined as \emph{continuous functions} since $\RE(2) = 2 > 1 =
n$.
\begin{proposition}
    \label{proposition:green-functions-in-1dim}
    \index{Riesz distribution!in dimension 1}%
    Let $n = 1$. Then the advanced and retarded Green functions of
    $\dAlembert = \frac{\partial^2}{\partial t^2}$ are explicitly
    given as the continuous functions
    \begin{equation}
        \label{eq:advanced-green-function-in-1dim}
        R^+(2)(t)
        = \begin{cases}
            t & \textrm{for}\quad t > 0 \\
            0 & \textrm{else}
        \end{cases}
    \end{equation}
    and
    \begin{equation}
        \label{eq:retarded-green-function-in-1dim}
        R^-(2)(t)
        = \begin{cases}
            |t| & \textrm{for}\quad t < 0 \\
            0 & \textrm{else}.
        \end{cases}
    \end{equation}
    Moreover, for $\RE(\alpha) > 1$ we have
    \begin{equation}
        \label{eq:riesz-in-1dim}
        R^\pm(\alpha)(t)
        = \begin{cases}
            \frac{1}{\Gamma(\alpha)} |t|^{\alpha-1}
            & \textrm{for} \quad t \in \mathbb{R}^\pm \\
            0 & \textrm{else}.
        \end{cases}
    \end{equation}
\end{proposition}
\begin{proof}
    First we compute for all $\alpha \in \mathbb{C}$
    \begin{align*}
        c(\alpha, 1)
        = \frac{2^{1-\alpha} \pi^{\frac{2-1}{2}}}
        {\Gamma(\frac{\alpha}{2}) \Gamma(\frac{\alpha-1}{2}+1)}
        = \frac{2^{1-\alpha} \sqrt{\pi}}
        {\Gamma(\frac{\alpha}{2}) \Gamma(\frac{\alpha}{2} +
          \frac{1}{2})}
        = \frac{1}{\Gamma(\alpha)}
    \end{align*}
    by \Index{Legendre's duplication formula}. Since $\eta(t) = t^2$
    and $I^\pm(0) = \mathbb{R}^\pm$ we have
    \eqref{eq:riesz-in-1dim}. Finally, $\Gamma(2) = 1$ whence
    \eqref{eq:advanced-green-function-in-1dim} and
    \eqref{eq:retarded-green-function-in-1dim} follow.
\end{proof}

\begin{remark}[Riesz distribution in one dimension]
    \label{remark:riesz-distributions-in-1dim}
    ~
    \begin{remarklist}
    \item \label{item:riesz-is-green-explicit} It is an easy exercise
        to compute $\frac{\partial^2}{\partial t^2} R^\pm(2)$ in the
        sense of distributions directly to show that
        \begin{equation}
            \label{eq:riesz-is-green}
            \frac{\partial^2}{\partial t^2} R^\pm(2) = \delta_0.
        \end{equation}
        In fact, we have done this implicitly in the proof of
        Proposition~\ref{proposition:riesz-zero-is-delta}.
    \item \label{item:riesz-in-1dim-in-hormander} The functions
        $R^\pm(\alpha)$ for $\RE(\alpha) > 1$ coincide with the
        functions $\chi^{\alpha-1}_\pm$ of Hörmander in
        \cite[Sect.~3.2., (3.2.17)]{hoermander:2003a}.  In fact, even
        though the function $R^\pm(\alpha)$ defined by
        \eqref{eq:riesz-in-1dim} is no longer continuous for
        $\RE(\alpha) > 0$, it is still \emph{locally integrable}. Thus
        it defines a distribution also in this case, depending
        holomorphically on $\alpha$. Hence we conclude
        \begin{equation}
            \label{eq:riesz-for-alpha-greater-zero}
            R^\pm(\alpha)(t) =
            \begin{cases}
                \frac{1}{\Gamma(\alpha)} |t|^{\alpha-1}
                & \textrm{for} \; t \in \mathbb{R}^\pm\\
                0 & \textrm{else}
            \end{cases}
        \end{equation}
        \index{Distribution!homogeneous}%
        is valid for $\RE(\alpha) > 0$ in the sense of locally
        integrable functions. The functions $\chi_\pm^\alpha$ are at
        the heart of the study of homogeneous distributions and can be
        used to obtain fundamental solutions of much more general
        second order differential operators with constant coefficients
        than just for $\dAlembert$, see
        \cite[Sect.~3.2]{hoermander:2003a}.
    \end{remarklist}
\end{remark}

\begin{figure}
    \centering
    \input{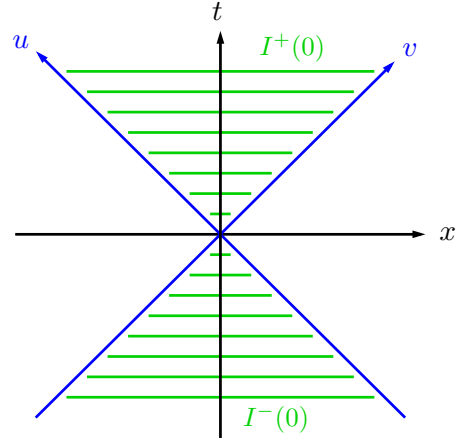}
    \caption{\label{fig:light-cone-coordinates}%
      Light cone coordinates.
    }
\end{figure}

We turn now to the case $n = 1+1$. Here it is convenient to use the
coordinates $(t, x) \in \mathbb{R}^2$ with
\begin{equation}
    \label{eq:eta-in-2dim-in-stand-coord}
    \eta(t, x) = t^2 -x^2.
\end{equation}
First we compute the prefactor $c(\alpha, n)$ for $n = 2$. We have
\begin{equation}
    \label{eq:prefactor-in-2dim}
    c(\alpha, 2) = \frac{2^{1-\alpha}}{\Gamma(\frac{\alpha}{2})^2}
\end{equation}
as one immediately obtains from the definition. In order to evaluate
$\eta^{\frac{\alpha-2}{2}}$ we introduce new coordinates on
$\mathbb{R}^2$. We pass to the \emph{light cone coordinates}
\begin{equation}
    \label{eq:light-cone-coordinates}
    \index{Light cone coordindates}%
    u = \frac{1}{\sqrt{2}} (t-x)
    \qquad
    \textrm{and}
    \qquad
    v = \frac{1}{\sqrt{2}} (t+x),
\end{equation}
i.e.
\begin{equation}
    \label{eq:lightcone-coord-inversion}
    t = \frac{1}{\sqrt{2}} (u+v)
    \qquad
    \textrm{and}
    \qquad
    x = \frac{1}{\sqrt{2}} (v-u).
\end{equation}
Since this is clearly a global diffeomorphism we can evaluate
$R^\pm(\alpha)$ in these new coordinates. The prefactors are chosen in
such a way that the diffeomorphism is orientation preserving and has
Jacobi determinant equal to one: It is just the counterclockwise
rotation by $45^\circ$ in the $(t, x)$-plane, see
Figure~\ref{fig:light-cone-coordinates}. First we note that the
function $\eta$ in these coordinates is
\begin{align}
    \label{eq:eta-in-lightcone-coord}
    \eta(u, v) = \frac{1}{2} (u+v)^2 - \frac{1}{2}(v-u)^2
    =\frac{1}{2} (u^2 + 2uv + v^2 - u^2 + 2uv - v^2)
    = 2uv.
\end{align}
Moreover, the future and past $I^\pm(0)$ of $0$ can be described by
\begin{equation}
    \label{eq:future-in-lightcone-coord}
    I^+(0) =
    \left\{(u,v) \in \mathbb{R}^2 \; \big| \; u,v > 0 \right\}
\end{equation}
and
\begin{equation}
    \label{eq:past-in-lightcone-coord}
    I^-(0) =
    \left\{(u,v) \in \mathbb{R}^2 \; \big| \; u,v < 0 \right\},
\end{equation}
see again Figure~\ref{fig:light-cone-coordinates}.
Thus we have for $\RE(\alpha) > 2$
\begin{align}
    \label{eq:riesz-in-lightcone-coordinates}
    R^\pm(\alpha)(u,v)
    & = \begin{cases}
        \frac{2^{1-\alpha}}{\Gamma(\frac{\alpha}{2})^2}
        (2uv)^{\frac{\alpha-2}{2}}
        & \textrm{for} \; u,v \in \mathbb{R}^\pm \\
        0 & \textrm{else}
    \end{cases} \\
    & = \begin{cases}
        \frac{2^{1-\alpha}}{\Gamma(\frac{\alpha}{2})^2}
        |\sqrt{2} u|^{\frac{\alpha-2}{2}}
        |\sqrt{2} v|^{\frac{\alpha-2}{2}}
        & \textrm{for} \; u,v \in \mathbb{R}^\pm \\
        0 & \textrm{else},
    \end{cases}
\end{align}
whence $R^\pm(\alpha)$ is \emph{factorizing} in these coordinates.
This suggests to consider the following functions
\begin{align}
    \label{eq:special-ansatz-for-lightcone-coordinates}
    r^\pm(\alpha)(u)
    = \begin{cases}
        \frac{2^{-\frac{\alpha}{4}}}{\Gamma(\frac{\alpha}{2})}
        |u|^{\frac{\alpha-2}{2}}
        & \textrm{for} \; u \in \mathbb{R}^\pm \\
        0 & \mathrm{else}
    \end{cases}
\end{align}
for $\RE(\alpha) > 2$. Since the prefactor is still holomorphic for
all $\alpha \in \mathbb{C}$ and since $|u|^z$ is locally integrable
for $\RE(z) > -1$ we can extend this definition to the case
$\RE(\alpha) > 0$.
\begin{proposition}
    \label{proposition:properties-of-special-ansatz-functions}
    Let $\RE(\alpha) > 0$.
    \begin{propositionlist}
    \item \label{item:ansatz-function-is-loc-integrable} The functions
        $r^\pm(\alpha)$ on $\mathbb{R}$ are locally integrable and
        thus define distributions of order zero with
        \begin{equation}
            \label{eq:supp-of-specia-ansatz-functions}
            \supp r^\pm(\alpha) = \mathbb{R}^\pm \cup \{0\}.
        \end{equation}
    \item \label{item:special-ansatz-defines-holomorphic-distribution}
        For every $\varphi \in \Cinfty_0(\mathbb{R})$ the function
        \begin{equation}
            \label{eq:special-ansatz-distribution}
            \alpha \; \mapsto \; \int r^\pm(\alpha)(u) \varphi(u) \D u
        \end{equation}
        is holomorphic for $\RE(\alpha) > 0$.
    \item \label{item:special-ansatz-funciton-2-is-heavyside} For
        $\alpha = 2$ we have
        \begin{equation}
            \label{eq:special-ansatz-funciton-2-is-heavyside}
            r^\pm(2)(u)
            = \begin{cases}
                \frac{1}{\sqrt{2}} & \textrm{for} \; u \in \mathbb{R}^\pm \\
                0 & \mathrm{else},
            \end{cases}
        \end{equation}
        i.e. a multiple of the \Index{Heaviside distribution}.
    \end{propositionlist}
\end{proposition}
\begin{proof}
    Let $\RE(\alpha) > 0$. Clearly, the only interesting thing about
    the local integrability of $r^\pm(\alpha)$ is around zero since on
    $\mathbb{R} \setminus \{0\}$ the function is clearly smooth and
    thus in $L^1_\loc$. Thus we consider $\alpha = \beta + \I \gamma$
    with $\beta > 0$ and $\gamma \in \mathbb{R}$. Then
    \begin{align*}
        |r^\pm(\alpha) (u) |
        = \begin{cases}
            \left|
                \frac{2^{\frac{\beta+\I \gamma}{4}}}
                {\Gamma(\frac{\beta+\I \gamma}{2})}
            \right|
            \left|
                |u|^{\frac{\beta+\I \gamma - 2}{2}}
            \right|
            & u \in \mathbb{R}^\pm \\
            0 & \mathrm{else}
        \end{cases}
        \quad
        \leq c_{\beta, \gamma} |u|^{\frac{\beta-2}{2}}.
    \end{align*}
    But the function $u \mapsto |u|^{\frac{\beta-2}{2}}$ is locally
    integrable for $\beta > 0$. This shows the first part as
    \eqref{eq:supp-of-specia-ansatz-functions} is obvious. Now let
    $\varphi \in \Cinfty_0(\mathbb{R})$ or, which would be sufficient,
    $\varphi \in \Fun[0]_0(\mathbb{R})$. Let $\supp \varphi \subseteq
    [a,b]$ and without restriction $b > 0$, then
    \begin{align*}
        \int_{\mathbb{R}} r^+(\alpha)(u) \varphi(u) \D u
        = \int_0^\infty
        \frac{2^{-\frac{\alpha}{4}}}{\Gamma(\frac{\alpha}{2})}
        u^{\frac{\alpha-2}{2}} \varphi(u) \D u
        = \frac{2^{-\frac{\alpha}{4}}}{\Gamma(\frac{\alpha}{2})}
        \int_0^b u^{\frac{\alpha-2}{2}} \varphi(u) \D u.
    \end{align*}
    Since the function $u \mapsto u^{\frac{\alpha-2}{2}} \varphi(u)$
    is integrable over $[a,b]$ we can again exchange the integration
    over a triangle path $\int_\Delta \D \alpha$ in the complex half
    space with $\RE(\alpha) > 0$ and the integral $\int_0^b \D u$. Thus
    Morera's theorem again yields the statement that
    \eqref{eq:special-ansatz-distribution} is holomorphic. The last
    part is clear.
\end{proof}
\begin{lemma}
    \label{lemma:dAlembert-in-lightcone-coordinates}
    \index{dAlembertian@d'Alembertian!in light cone coordinates}%
    In the lightcone coordinates the d'Alembert operator is
    \begin{equation}
        \label{eq:dAlembert-in-lightcone-coordinates}
        \dAlembert = 2 \frac{\partial^2}{\partial u \partial v}.
    \end{equation}
\end{lemma}
\begin{proof}
    This is a trivial computation.
\end{proof}
\begin{proposition}
    \label{proposition:green-function-in-2dim}
    \index{Riesz distribution!in dimension 2}%
    Let $u, v$ be the light cone coordinates on $\mathbb{R}^2$. Then
    the distributions
    \begin{equation}
        \label{eq:green-function-in-2dim}
        R^\pm(2)(u,v) = r^\pm(2)(u) r^\pm(2)(v)
    \end{equation}
    are advanced and retarded Green functions of $\dAlembert$ of order
    zero.
\end{proposition}
\begin{proof}
    Of course, we know this from the general
    Theorem~\ref{theorem:green-function-of-dAlembert}, but here we can
    give a more elementary proof. We consider the $+$ case where we
    have
    \[
    r^+(2)(u) = \frac{1}{\sqrt{2}} \Theta(u)
    \]
    with the Heaviside function $\Theta$. Since $r^+(2)$ is locally
    integrable, we have a distribution $r^+(2) \in
    \mathcal{D}'(\mathbb{R})$ of order zero. Moreover, one knows
    \[
    \frac{\partial}{\partial u} \Theta = \delta_0.
    \]
    The same holds for the $v$-dependence. Thus we can interpret
    \eqref{eq:green-function-in-2dim} as \emph{external} tensor
    product
    \[
    R^\pm(2) = \frac{1}{2} \Theta_u \extensor \Theta_v,
    \]
    whence
    \begin{align*}
        2 \frac{\partial^2}{\partial u \partial v} R^\pm(2)
        = 2 \frac{\partial}{\partial u} \extensor
        \frac{\partial}{\partial v}
        \left(
            \frac{1}{2} \Theta_u \extensor \Theta_v
        \right)
        = \frac{\partial}{\partial u} \Theta_u
        \extensor
        \frac{\partial}{\partial v} \Theta_v
        = \delta_{(0,0)}.
    \end{align*}
    Since the Jacobi determinant of the coordinate change is one, the
    $\delta$-distribution in $(u,v)$ is the same as the one in
    $(t,x)$. Thus the claim follows. Note that this formulation is of
    course more elementary and can almost be ``guessed''.
\end{proof}
\begin{remark}
    \label{remark:riesz-in-2dim}
    In $n = 1+1$ all the Riesz distributions are factorizing as
    external tensor products of the distributions $r^\pm(\alpha)$ of
    one variable. This simplifies the discussion considerably. Note
    however, that this is a particular feature of $n = 1+1$ and no
    longer true in higher dimensions. Note also that from
    Proposition~\ref{proposition:order-of-riesz-distribution} we only
    get the estimate $\ord (R^\pm(2)) \leq 2$ which is clearly not
    optimal: The Riesz distributions $R^\pm(2)$ in $n = 1+1$ are
    locally integrable and hence of order zero.
\end{remark}

It is a good exercise to work out the cases $n = 1+2$ and $n = 1+3$
explicitly.


%% file: riesz.tex
%
%

We pass now from Minkowski spacetime to a general Lorentz manifold
$(M,g)$ and try to find analogs of the Riesz distributions at least
locally around a point $p \in M$. The main idea is to use the Riesz
distributions on the tangent space $T_pM$, which is isometric to
Minkowski space, and push forward the Riesz distributions via the
exponential map.

%
%

\subsection{The Functions $\varrho_p$ and $\eta_p$}
\label{subsec:FunctionsRhopEtap}

Since on $M$ we have a canonical positive density, namely the Lorentz
volume density\index{Lorentz density} $\mu_g$ from
Proposition~\ref{proposition:levi-civita-connection-etc},
\refitem{item:Riem-density}, we can use this density to identify
functions and densities once and for all. In particular, this results
in an identification of the generalized sections $\Sec[-\infty](E)$ of
a vector bundle $E \longrightarrow M$ with the topological dual of
$\Secinfty_0(E^*)$ and not of $\Secinfty_0(E \tensor \Dichten T^*M)$
as we did before. In more detail, for $s \in \Sec[-\infty](E)$ and a
test section $\varphi \in \Secinfty_0(E^*)$ we first map $\varphi$ to
$\varphi \tensor \mu_g \in \Secinfty_0(E^* \tensor \Dichten T^*M)$ and
then apply $s$, i.e. we set
\begin{equation}
    \label{eq:distributions-for-fixed-lorentz-density}
    s(\varphi) = s(\varphi \tensor \mu_g),
\end{equation}
and drop the explicit reference to $\mu_g$ to simplify our
notation. Since
\begin{equation}
    \label{eq:identication-functions-densities-with-lorentz-dens}
    \Secinfty_0(E^*)
    \ni \varphi \; \mapsto \; \varphi \tensor \mu_g \in
    \Secinfty_0(E^* \tensor \Dichten T^*M)
\end{equation}
is indeed an isomorphism of LF spaces as discussed in
Remark~\ref{remark:generalized-sections-for-fixed-density}, we have an
induced isomorphism of the topological duals which is
\eqref{eq:distributions-for-fixed-lorentz-density}.

If we now want to push forward the $R^\pm(\alpha)$ from $T_pM$ to $M$
we have to take care of the two different notion of volume densities.
On $T_pM$ we have the \emph{constant} density coming from the
Minkowski scalar product $g_p$ while one $M$ we have $\mu_g$. In
general, the push-forward of $\mu_g(p)$ via $\exp_p$ to $M$ \emph{does
  not} coincide with $\mu_g$ whence we need a way to compare the two
densities. This is done by the following construction. Let $V_p
\subseteq T_pM$ be a suitable open star-shaped neighborhood of $0_p$
and let $U_p = \exp_p(V_p) \subseteq M$ be the corresponding open
neighborhood of $p$ such that
\begin{equation}
    \label{eq:starshaped-domain-with-exp-diffeo}
    \exp_p: V_p \longrightarrow U_p
\end{equation}
is a diffeomorphism. Then we define the function
\begin{equation}
    \label{eq:density-compare-function}
    \varrho_p = \frac{\mu_g \at{U_p}}{\exp_{p*}(\mu_g(p))\at{U_p}}.
\end{equation}
\begin{lemma}
    \label{lemma:density-dompare-function}
    The function $\varrho_p$ is well-defined and smooth on $U_p$. We
    have $\varrho_p > 0$ and
    \begin{equation}
        \label{eq:density-compare-function-alternative}
        \varrho_p \exp_{p*} (\mu_g(p)) = \mu_g
    \end{equation}
    on $U_p$.
\end{lemma}
\begin{proof}
    Since on $V_p$ the exponential map is a diffeomorphism, the
    push-forward of the constant density $\mu_g(p) \in \Dichten
    T^*_pM$ gives a smooth density on $U_p$. Clearly, it is still
    positive whence the quotient \eqref{eq:density-compare-function}
    is well-defined and a smooth function. Since also $\mu_g > 0$ it
    follows that $\varrho_p > 0$ everywhere.
\end{proof}

Sometimes it will be convenient to work on $V_p$ instead of
$U_p$. Thus we can pull-back everything to $V_p$ by $\exp_p$ and
obtain
\begin{equation}
    \label{eq:pull-back-of-lorentz-density}
    \exp_p^*(\varrho_p)
    \underbrace{\exp_p^*(\exp_{p*} \mu_g(p))}_{\mu_g(p)}
    = \exp_p^*(\varrho_p) \mu_g(p)
    = \exp_p^*(\mu_g)
\end{equation}
on $V_p$. To simplify our notation we abbreviate
\begin{equation}
    \label{eq:density-compare-function-on-TpM}
    \widetilde{\varrho}_p = \exp_p^*(\varrho_p) \in \Cinfty(V_p)
\end{equation}
and have
\begin{equation}
    \label{eq:densitiy-compare-func-on-TpM-equation}
    \widetilde{\varrho}_p \mu_g(p) = \exp_p^*(\mu_g).
\end{equation}
Thus $\widetilde{\varrho}_p$ is the function which measures how much
$\exp_p^*(\mu_g)$ is \emph{not} constant.

To effectively compute $\varrho_p$ or $\widetilde{\varrho}_p$ one
proceeds as follows. Let $e_1, \ldots, e_n \in T_pM$ be a basis. Then
we can evaluate both densities on $e_1, \ldots, e_n$ to get
$\varrho_p$ and $\widetilde{\varrho}_p$. More precisely, by the
definition of the pull-back we have for $v \in V_p$
\begin{equation}
    \label{eq:dens-comp-func-by-evaluation}
    \widetilde{\varrho}_p (v)
    =
    \frac{\mu_g (\exp_p(v))
      \left(
          T_v \exp_p e_1, \ldots, T_v \exp_p e_n
      \right)}
    {\mu_g(p)(e_1, \ldots, e_n)}.
\end{equation}
\index{Jacobi vector field}%
Thus we have to compute ``determinants'' of the tangent map of
$\exp_p$ in order to obtain $\widetilde{\varrho}_p$. This can indeed
be done rather explicitly by using Jacobi vector fields at least in a
formal power series expansion in $v$. We give here the result without
going into details, but refer to
Appendix~\ref{sec:jacobi-determinants-exp} for more background
information.
\begin{proposition}
    \label{proposition:taylor-expansion-of-dens-compare-function}
    The Taylor expansion of $\widetilde{\varrho}_p$ up to second order
    is explicitly given by
    \begin{equation}
        \label{eq:taylor-expansion-of-dens-compare-function}
        \widetilde{\varrho}_p (v)
        = 1 - \frac{1}{6} \Ric_p(v,v) + \cdots,
    \end{equation}
    where $\Ric_p$ is the \Index{Ricci tensor} at $p$ and $v \in
    T_pM$.
\end{proposition}
\begin{proof}
    First we note that $\mu_g$ is covariantly constant and hence
    parallel along all curves. Thus we do not get any contributions of
    covariant derivatives of $\mu_g$. This simplifies the general
    result from
    Theorem~\ref{theorem:taylor-expansion-of-dens-comp-func} and
    yields the result
    \eqref{eq:taylor-expansion-of-dens-compare-function}.
\end{proof}

\begin{corollary}
    \label{corollary:dAlembert-of-densitiy-comp-func}
    \index{Scalar curvature}%
    At $p \in M$ we have
    \begin{equation}
        \label{eq:dAlembert-of-densitiy-comp-func}
        \dAlembert \varrho_p \at{p}
        = - \frac{1}{3} \scal(p).
    \end{equation}
\end{corollary}
\begin{proof}
    By general results from
    Appendix~\ref{sec:taylor-parallel-transport} we know that for any
    function $f \in \Cinfty(M)$ one has the formal Taylor expansion
    \[
    (\exp_p^* f)(v) \sim_{v \rightarrow 0}
    \sum_{r=0}^\infty \frac{1}{r!} \frac{1}{r!}
    \SymD^r f \at{p}(v, \ldots, v),
    \]
    where $\SymD$ is the symmetrized covariant derivative. By
    Proposition~\ref{proposition:taylor-expansion-of-dens-compare-function}
    we have for $\varrho_p = \exp_{p*}(\widetilde{\varrho}_p)$
    \[
    \frac{1}{4} \SymD^2 \varrho_p \at{p}(v,v)
    = - \frac{1}{6} \Ric_p(v,v),
    \]
    whence the Hessian of $\varrho_p$ at $p$ is given by
    \[
    \mathrm{Hess} \varrho_p
    = \frac{1}{2} \SymD^2 \varrho_p
    = - \frac{1}{3} \Ric_p.
    \]
    Thus we conclude
    \begin{align*}
        \dAlembert \varrho_p \at{}
        = \frac{1}{2} \SP{g^{-1}, \SymD^2 \varrho_p} \at{p}
        = - \frac{1}{3} \SP{g^{-1}, \Ric} \at{p}
        = - \frac{1}{3} \scal(p)
    \end{align*}
    by the definition of the scalar curvature as in
    \eqref{eq:scalar-curvature} as well as by
    Proposition~\ref{proposition:grad-div-box-local-expressions},
    \refitem{item:dAlembert-is-second-order-diffop}.
\end{proof}

With the general techniques from the appendix it is also possible to
obtain the higher orders in the Taylor expansion of $\varrho_p$ in a
rather explicit and systematic way. They turn out to be universal
algebraic combinations of the curvature tensor and its covariant
derivatives. However, we shall not need this here. Instead, we mention
that by the usual expansion $\sqrt{1+x} = 1 + \frac{1}{2}x + \cdots$
we immediately find
\begin{equation}
    \label{eq:dAlembert-of-sqrt-of-dens-comp-funcI}
    \dAlembert \sqrt{\varrho_p} \At{p} = - \frac{1}{6} \scal(p)
\end{equation}
and
\begin{equation}
    \label{eq:dAlembert-of-sqrt-of-dens-comp-funcII}
    \dAlembert \frac{1}{\sqrt{\varrho_p}} \At{p}
    = \frac{1}{6} \scal(p).
\end{equation}
For the Riesz distributions we needed the quadratic function $\eta(x)
= \eta(x,x)$ as basic ingredient. Clearly, we have this on every
tangent space whence we can define
\begin{equation}
    \label{eq:length-function-of-vector-on-TpM}
    \widetilde{\eta}_p(v) = g_p(v,v)
\end{equation}
for every $v \in T_pM$. Analogously to the relation between
$\varrho_p$ and $\widetilde{\varrho}_p$ we set
\begin{equation}
    \label{eq:distance-function-on-M}
    \eta_p(q) = \widetilde{\eta}_p(\exp_p^{-1}(q))
\end{equation}
for $q \in U_p$. With other words, $\eta_p \in \Cinfty(U_p)$ is the
function with
\begin{equation}
    \label{eq:relation-of-distance-and-lengthonTpM-functions}
    \exp_p^*(\eta_p) = \widetilde{\eta}_p.
\end{equation}
We collect now some properties of the functions $\eta_p$ and
$\widetilde{\eta}_p$.
\begin{proposition}
    \label{proposition:distance-function-on-M}
    Let $(M,g)$ be a time-oriented Lorentz manifold and $p \in
    M$. Moreover, let $U \subseteq M$ be geodesically star-shaped with
    respect to $p$.
    \begin{propositionlist}
    \item\label{item:gradient-of-distance-function} The gradient of
        $\eta_p \in \Cinfty(N)$ is given by
        \begin{equation}
            \label{eq:gradient-of-distance-function}
            \gradient \eta_p \at{q}
            = 2 T_{\exp_p^{-1}(q)} \exp_p (\exp_p^{-1}(q))
        \end{equation}
        for $q \in U$.
    \item\label{item:length-of-gradient-of-dist-func} One has
        \begin{equation}
            \label{eq:length-of-gradient-of-dist-func}
            g (\gradient \eta_p, \gradient \eta_p)
            = 4 \eta_p.
        \end{equation}
    \item \label{item:gradient-of-dist-is-future-past-directed} On
        $I^\pm_U(p)$ the gradient of $\eta_p$ is a future resp. past
        directed timelike vector field.
    \item\label{item:dAlembert-of-distant-function} One has
        \begin{equation}
            \label{eq:dAlembert-of-distant-function}
            \dAlembert \eta_p
            = 2n + g (\gradient \log \varrho_p, \gradient \eta_p).
        \end{equation}
    \end{propositionlist}
\end{proposition}
\begin{proof}
    For the first part we need the \Index{Gauss Lemma} which says
    \[
    g_{\exp_p(v)}
    \left(
        T_v \exp_p (v), T_v \exp_p (w)
    \right)
    = g_p(v,w)
    \]
    for $v \in V_p \subseteq T_pM$ and $w \in T_pM$ arbitrary, see
    Proposition~\ref{proposition:gauss-lemma}. Using this we compute
    for $q \in U$ and $w_q \in T_qM$
    \begin{align*}
        \D \eta_p \at{q} (w_q)
        &= w_q (\widetilde{\eta}_p \circ \exp_p^{-1}) \\
        &= \D \widetilde{\eta}_p \At{\exp_p^{-1}(q)}
        \left(
            T_q \exp_p^{-1} (w_q)
        \right)\\
        &= 2 g_p
        \left(
            \exp_p^{-1}(q), T_q \exp_p^{-1} (w_q)
        \right) \\
        &= 2 g_q
        \left(
            T_{\exp_p^{-1}(q)} \exp_p (\exp_p^{-1}(q)),
            T_{\exp_p^{-1}(q)} \exp_p T_q \exp_p^{-1} (w_q)
        \right) \\
        &= 2 g_q
        \left(
            T_{\exp_p^{-1}(q)} \exp_p (\exp_p^{-1}(q)), w_q
        \right)
    \end{align*}
    by the Gauss Lemma for $v = \exp_p^{-1}(q)$ and the chain rule. By
    the very definition of the gradient this gives
    \eqref{eq:gradient-of-distance-function}. For the second part we
    again use the Gauss Lemma and get with $v = \exp_p^{-1}(q)$ for $q
    \in U$
    \begin{align*}
        g_q (\gradient \eta_p \at{q}, \gradient \eta_p \at{q})
        &= 4 g_{\exp_p(v)}
        \left( T_v \exp_p(v), T_v \exp_p(v) \right) \\
        &= 4 g_p(v,v) \\
        &= 4 g_p (\exp_p^{-1}(q), \exp_p^{-1}(q)) \\
        &= 4 \eta_p(q)
    \end{align*}
    as claimed. For the third part we first notice that the points in
    $I^\pm(0_p) \subseteq T_pM$ are mapped under $\exp_p$ to points in
    $I^\pm_U(p)$ since there is a timelike curve joining $p$ and such
    a point $q = \exp_p(v)$, namely the geodesic $t \mapsto
    \exp_p(tv)$. This is indeed a timelike curve for all $t$ thanks to
    the Gauss Lemma. Thus for $q \in I^\pm_U(p)$ we have $q =
    \exp_p(v)$ with $v \in I^\pm(0_p) \subseteq T_pM$ whence
    $\eta_p(q) = g_p (\exp_p^{-1}(q), \exp_p^{-1}(q)) > 0$. This shows
    $\eta_p > 0$ on $I^\pm_U(p)$. By the second part we conclude
    \[
    g (\gradient \eta_p, \gradient \eta_p) > 0
    \]
    on $I^\pm_U(p)$ whence $\gradient \eta_p$ is timelike on
    $I^\pm_U(p)$. Now let $v \in I^+(0_p) \subseteq T_pM$ be future
    directed. Then $t \mapsto \exp_p(tv)$ is a future directed
    geodesic with tangent vector
    \[
    \frac{\D}{\D t} \exp_p(tv)
    = T_{tv} \exp_p (v)
    = \frac{1}{2t} \gradient \eta_p \At{\exp_p(tv)}.
    \]
    Thus for $t > 0$ the gradient of $\eta_p$ is a positive multiple
    of the tangent vector of $\exp_p(tv)$ and hence future directed
    itself at $\exp_p(tv)$. Since every point in $I^+_U(p)$ can be
    reached this way, $\gradient \eta_p$ is future directed on all of
    $I^+_U(p)$. With the same argument we see that $\gradient \eta_p$
    is past directed on $I^-_U(p)$. The last part is again a
    computation. First we note that thanks to $\varrho_p > 0$
    everywhere, we have a smooth real-valued logarithm $\log \varrho_p
    \in \Cinfty(U)$. The Leibniz rule \eqref{eq:divergence-product
      rule} for $\divergenz$ gives
    \begin{align*}
        \divergenz \left(\frac{1}{\varrho_p} \gradient \eta_p \right)
        = \Lie_{\gradient \eta_p} \left(\frac{1}{\varrho_p} \right)
        + \frac{1}{\varrho_p} \divergenz \gradient \eta_p
        = g
        \left( \gradient \frac{1}{\varrho_p}, \gradient \eta_p \right)
        + \frac{1}{\varrho_p} \dAlembert \varrho_p
    \end{align*}
    and thus
    \begin{align*}
        \dAlembert \varrho_p
        & = \varrho_p \divergenz
        \left( \frac{1}{\varrho_p} \gradient \eta_p \right)
        - g
        \left(
            \varrho_p \gradient \frac{1}{\varrho_p}, \gradient \eta_p
        \right) \\
        & = \varrho_p \divergenz
        \left( \frac{1}{\varrho_p} \gradient \eta_p \right)
        - g
        \left(
            \gradient \log \frac{1}{\varrho_p}, \gradient \eta_p
        \right) \\
        & = \varrho_p \divergenz
        \left( \frac{1}{\varrho_p} \gradient \eta_p \right)
        + g
        \left(
            \gradient \log \varrho_p, \gradient \eta_p
        \right).
    \end{align*}
    We still have to compute the first divergence. Since $\divergenz$
    here is always the divergence with respect to $\mu_g$ we consider
    on $U$
    \begin{align*}
        \divergenz_{\mu_g} \left( \frac{1}{\varrho_p} X \right)
        & = \divergenz_{\varrho_p \exp_{p*}(\mu_g(p))}
        \left( \frac{1}{\varrho_p} X \right)\\
        & = \divergenz_{\exp_{p*}(\mu_g(p))}
        \left( \frac{1}{\varrho_p} X \right)
        +
        \Lie_X (\log \varrho_p) \\
        & = \Lie_X \left( \frac{1}{\varrho_p} \right)
        + \frac{1}{\varrho_p} \divergenz_{\exp_{p*}(\mu_g(p))}(X)
        + \frac{1}{\varrho_p} \Lie_X \varrho_p \\
        & = \frac{1}{\varrho_p} \divergenz_{\exp_{p*}(\mu_g(p))} (X)
    \end{align*}
    by the chain rule and the behaviour of the divergence operator
    under the change of the reference density, see e.g.
    \cite[Lemma~2.3.45]{waldmann:2007a}. Thus we have for a general
    vector field $X$
    \[
    \varrho_p \divergenz \left( \frac{1}{\varrho_p} X \right)
    = \divergenz_{\exp_{p*}(\mu_g(p))} (X)
    \]
    on $U$. Since the definition of the divergence operator is natural
    with respect to diffeomorphisms we have
    \[
    \divergenz_{\exp_{p*}(\mu_g(p))} (X)
    = \exp_{p*} \left( \divergenz_{\mu_g(p)} (\exp_p^* X) \right).
    \]
    Now we consider again $X = \gradient \eta_p$ whence
    \begin{align*}
        \exp_p^*(\gradient \eta_p) \At{v}
        &= T_{\exp_p(v)} \exp_p^{-1}
        \left( \gradient \eta_p \At{\exp_p(v)} \right) \\
        & = 2 T_{\exp_p(v)} \exp_p^{-1}
        \left(
            T_{\exp_p^{-1} (\exp_p(v))} \exp_p (
            \exp_p^{-1}(\exp_p(v)))
        \right) \\
        & = 2v.
    \end{align*}
    With other words
    \[
    \exp_p^*(\gradient \eta_p) = 2 \xi_{T_pM}
    \]
    is twice the \Index{Euler vector field} on the tangent space
    $T_pM$. But the divergence of $\xi_{T_pM}$ with respect to the
    constant density is easily seen to be $n = \dim M$. Thus we end up
    with $\dAlembert \eta_p = 2n + g (\gradient \log \varrho_p,
    \gradient \eta_p)$, finishing the proof.
\end{proof}
\begin{remark}
    \label{remark:dAlembert-of-distance}
    In fact, it will be the last statement of the last proposition
    which causes new complications compared to the trivial, flat
    case. Here we have of course
    \begin{equation}
        \label{eq:dAlembert-of-distance-in-flat-case}
        \dAlembert_{\textrm{flat}} \eta_p^{\textrm{flat}} = 2n
    \end{equation}
    \emph{without} the additional term as in
    \eqref{eq:dAlembert-of-distant-function}. Clearly
    $\varrho_p^{\textrm{flat}}=1$ whence this additional contribution
    vanishes. However, \eqref{eq:dAlembert-of-distance-in-flat-case}
    was essential for the correct functional equation of the (flat)
    Riesz distributions in
    Section~\ref{sec:dAlembertOnMinkowskiSpaceTime}.
\end{remark}

%
%

\subsection{Construction of the Riesz Distributions $R^\pm_U(\alpha, p)$}
\label{subsec:ConstructionRieszDistribution}

For $\RE(\alpha) > n$ the Riesz distributions $R^\pm(\alpha)$ are even
continuous functions on Minkowski space. As such we can simply
push-forward them via $\exp_p$, at least on the star-shaped $V
\subseteq T_pM$, to a continuous function on $U \subseteq M$. There, a
continuous function defines a distribution after multiplying with the
density $\mu_g$.
\begin{remark}
    \label{remark:continuous-functions-on-TpM-or-M-as-distributions}
    Let $f \in \Fun[0](T_pM)$ be a continuous function on the tangent
    space of $p$. We view $f$ as a distribution as usual via
    \begin{equation}
        \label{eq:continuous-function-on-TpM-as-distribution}
        f(\varphi) = \int_{T_pM} f(v) \varphi(v) \: \mu_g(p)
    \end{equation}
    for $\varphi \in \Cinfty_0(T_pM)$. Using $\exp_p$ we can write
    this as follows. Let $\varphi \in \Cinfty_0(M)$ with $\supp
    \varphi \subseteq U$ then the continuous function $\exp_{p*}
    (f\at{V}) \in \Fun[0](U)$ can be viewed as a distribution on $U$
    \begin{equation}
        \label{eq:push-forward-of-function-as-distribution}
        \exp_{p*}(f\at{V}) (\varphi)
        = \int_M \exp_{p*}(f\at{V})(q) \varphi(q) \: \mu_g(q)
    \end{equation}
    according to our convention. This equals
    \begin{align}
        \label{eq:cont-func-on-TpM-as-distrib-on-M}
        \exp_{p*}(f\at{V}) (\varphi)
        &= \int_M \exp_{p*} (f\at{V}) \exp_{p*} (\exp_p^* \varphi)(q)
        \varrho_p(q) (\exp_{p*} \mu_g(p))(q) \\
        & = \int_M \exp_{p*}
        \left(
            f\at{V} \exp_p^* \varphi \exp_p^* \varrho_p \mu_g(p)
        \right)(q) \\
        & = \int_{T_pM} f \exp_p^* \varphi \widetilde{\varrho}_p
        \: \mu_g(p) \\
        & = (\widetilde{\varrho}_pf) (\exp_p^* \varphi).
    \end{align}
    Thus, if we want to have a consistent definition of the
    push-forward of a distribution on $T_pM$ to a distribution on $M$
    we should include the prefactor $\widetilde{\varrho}_p$: let $u
    \in \mathcal{D}'(T_pM) = \Cinfty_0(T_pM)'$ be a distribution. Then
    one defines $\exp_{p*}u$ as the distribution $\exp_{p*}(u\at{V})
    \in \mathcal{D}'(U) = \Cinfty_0(U)'$ via
    \begin{equation}
        \label{eq:push-forward-of-distribution-from-TpM-to-M}
        \exp_{p*}(u\at{V}) (\varphi)
        = u (\widetilde{\varrho}_p \exp_p^* \varphi),
    \end{equation}
    which is a well-defined distribution as the restriction of $u$ to
    $V$ is a well-defined distribution on $V$ and $\supp
    (\widetilde{\varrho}_p \exp_p^* \varphi) \subseteq V$ thanks to
    $\supp \varphi \subseteq U$. Note that this definition
    \emph{differs} from the entirely intrinsic definition of the
    push-forward of distributions in
    Proposition~\ref{proposition:push-forward-of-distributions} in so
    far as we have modified our notion of distributions itself.
\end{remark}

We apply this construction of the push-forward now to the Riesz
distributions $R^\pm(\alpha)$. First we note that $R^\pm(\alpha)$ is
intrinsically defined on $T_pM$ \emph{without} specifying a particular
isometric isomorphism $(T_pM, g_p) \simeq (\mathbb{R}^n, \eta)$. The
reason is that $R^\pm(\alpha)$ on Minkowski spacetime is invariant
under orthochronous Lorentz transformations. We still denote the Riesz
distribution on $T_pM$ by $R^\pm(\alpha)$. Then the following
definition makes sense:
\begin{definition}[Riesz distributions on $U$]
    \label{definition:riesz-distribution-on-U}
    \index{Riesz distribution!on domain}%
    Let $p \in M$ and let $U \subseteq M$ be a geodesically
    star-shaped open neighborhood of $p$. Moreover, let $V =
    \exp_p^{-1}(U) \subseteq T_pM$ be the corresponding star-shaped
    open neighborhood of $0 \in T_pM$. Then the advanced and retarded
    Riesz distributions $R^\pm_U(\alpha,p) \in \Cinfty_0(U)'$ are
    defined by
    \begin{equation}
        \label{eq:riesz-distribution-on-U}
        R^\pm_U(\alpha,p)(\varphi)
        = \exp_{p*}\left(R^\pm(\alpha) \at{V}\right) (\varphi)
        = R^\pm(\alpha) \at{V}
        \left(\widetilde{\varrho}_p \exp_p^* \varphi\right)
    \end{equation}
    for $\alpha \in \mathbb{C}$ and $\varphi \in \Cinfty_0(U)$.
\end{definition}
We collect now the properties of $R^\pm(\alpha,p)$ in complete analogy
to those of $R^\pm(\alpha)$. In fact, most properties can be
transferred immediately using
\eqref{eq:riesz-distribution-on-U}. However, when it comes to
differentiation, the additional prefactor $\widetilde{\varrho}_p$ has
to be taken into account properly.
\begin{proposition}
    \label{proposition:properties-of-riesz-on-U}
    Let $U \subseteq M$ be geodesically star-shaped around $p \in M$.
    Then the Riesz distributions $R^\pm_U(\alpha,p)$ have the
    following properties:
    \begin{propositionlist}
    \item\label{item:riesz-on-U-is-cont-for-large-alpha} If
        $\RE(\alpha) > n$ then $R^\pm_U(\alpha,p)$ is continuous on
        $U$ and given by
        \begin{equation}
            \label{eq:riesz-on-U-for-large-alpha}
            R^\pm_U(\alpha,p)(q)
            = \begin{cases}
                c(\alpha,n) \left(\eta_p(q)\right)^{\frac{\alpha-n}{2}} &
                \textrm{for}\quad q \in I^\pm_U(p) \\
                0 & \textrm{else}.
            \end{cases}
        \end{equation}
    \item\label{item:differentiable-of-riesz-on-U-for-large-alpha} For
        $\RE(\alpha) > n+2k$ the function $R^\pm_U(\alpha,p)$ is even
        $\Fun$ on $U$.
    \item\label{item:riesz-inside-and-oustide-lightcone} For all
        $\alpha$ we have $R^\pm_U(\alpha,p) \at{I^\pm_U(p)} =
        c(\alpha,n) \eta_p^{\frac{\alpha-n}{2}} \in
        \Cinfty\left(I^\pm_U(p)\right)$ and $0 = R^\pm(\alpha,p) \at{U
          \setminus J^\pm_U(p)} \in \Cinfty\left(U \setminus
            J^\pm_U(p)\right)$.
    \end{propositionlist}
\end{proposition}
\begin{proof}
    By definition of $\exp_{p*} \left(R^\pm(\alpha)\at{V}\right)$ the
    singularities of $R^\pm(\alpha)$ correspond one-to-one to the
    singularities of $R^\pm_U(\alpha,p)$ under $\exp_p$ since $\exp_p$
    is a \emph{diffeomorphism} and the function
    $\widetilde{\varrho}_p$ is smooth and nonzero on $V$. In
    particular, for $q \notin J^\pm_U(p)$ we have $\exp_p^{-1}(q)
    \notin J^\pm(0) \subseteq T_pM$. Thus on this open subset,
    $R^\pm(\alpha)$ coincides with the smooth function being
    identically zero. This shows $R^\pm_U(\alpha,p) \at{U \setminus
      J^\pm_U(p)} = 0$. Inside the light cone, i.e. for $q \in
    I^\pm_U(p)$ and hence $\exp_p^{-1}(q) \in I^\pm(0)$, we have that
    $R^\pm(\alpha)$ is the smooth function $c(\alpha,n)
    \widetilde{\eta}_p^{\frac{\alpha-n}{2}}$. Thus by
    \eqref{eq:push-forward-of-distribution-from-TpM-to-M} we have for
    $\varphi \in \Cinfty_0\left(I^\pm_U(p)\right)$
    \begin{align*}
        R^\pm_U(\alpha,p)(\varphi)
        &= R^\pm_U(\alpha,p) \At{I^\pm_U(p)} (\varphi) \\
        &= R^\pm(\alpha) \At{I^\pm(0)}
        \left(\widetilde{\varrho}_p \exp_p^* \varphi\right) \\
        & = \int_{I^\pm(0)} R^\pm(\alpha)(v) \widetilde{\varrho}_p(v)
        (\exp_p^* \varphi)(v) \D^n v \\
        & = \int_{I^\pm_U(p)} c(\alpha,n)
        \left(
            \widetilde{\eta}_p \circ \exp_p^{-1}
        \right)^{\frac{\alpha-n}{2}}(q)
        \varphi(q) \: \mu_g(q) \\
        & = \int_{I^\pm_U(p)} c(\alpha,n)
        \left(\eta_p(q)\right)^{\frac{\alpha-n}{2}}
        \varphi(q) \: \mu_g(q) \\
        & = \left(
            c(\alpha,n) \eta_p^{\frac{\alpha-n}{2}}
        \right)
        (\varphi),
    \end{align*}
    since $\widetilde{\eta}_p \circ \exp_p^{-1} = \eta_p$ by
    definition of $\eta_p$. This shows the third part. The first and
    second part follow from the continuity properties of
    $R^\pm(\alpha)$ as in Lemma~\ref{lemma:riesz-function-continuity}
    and Lemma~\ref{lemma:identities-for-riesz-distribution},
    \refitem{item:riesz-differentiation-rule}.
\end{proof}

The analogue of Lemma~\ref{lemma:ries-distribution-is-homolorphic} and
Lemma~\ref{lemma:riesz-holomorphic-extension} is the following
statement:
\begin{proposition}
    \label{proposition:riesz-on-U-is-holomorphic}
    Let $U \subseteq M$ be star-shaped around $p \in M$. Then for
    every fixed test function $\varphi \in \Cinfty_0(U)$ the map
    $\alpha \mapsto R^\pm_U(\alpha,p)(\varphi)$ is entirely
    holomorphic on $\mathbb{C}$.
\end{proposition}
\begin{proof}
    Since for $\varphi \in \Cinfty_0(U)$ the function
    $\widetilde{\varrho}_p \exp_p^* \varphi$ is a test function on $V
    \subseteq T_pM$ and hence on $T_pM$,
    Lemma~\ref{lemma:ries-distribution-is-homolorphic} and
    Lemma~\ref{lemma:riesz-holomorphic-extension} guarantee that
    $\alpha \mapsto R^\pm(\alpha)(\widetilde{\varrho}_p \exp_p^*
    \varphi)$ is holomorphic.
\end{proof}
\begin{proposition}
    \label{proposition:differentiation-of-riesz-on-U}
    Let $U \subseteq M$ be geodesically star-shaped around $p \in M$.
    \begin{propositionlist}
    \item\label{item:identity-for-eta-times-riesz} For all $\alpha \in
        \mathbb{C}$ we have
        \begin{equation}
            \label{eq:identity-for-eta-times-riesz}
            \eta_p R^\pm_U(\alpha,p)
            = \alpha (\alpha-n+2) R^\pm_U(\alpha+2,p).
        \end{equation}
    \item\label{item:gradient-of-eta-times-riesz} For all $\alpha \in
        \mathbb{C}$ we have
        \begin{equation}
            \label{eq:gradient-of-eta-times-riesz}
            \gradient \eta_p \cdot R^\pm_U(\alpha,p)
            = 2 \alpha \gradient R^\pm_U(\alpha+2,p).
        \end{equation}
    \item\label{item:dAlembert-of-riesz-on-U-alpha-not-zero} For all
        $\alpha \in \mathbb{C} \setminus \{0\}$ we have
        \begin{equation}
            \label{eq:dAlembert-of-riesz-on-U-alpha-not-zero}
            \dAlembert R^\pm_U(\alpha+2,p)
            = \left(
                \frac{\dAlembert \eta_p - 2n}{2\alpha}+1
            \right)
            R^\pm_U(\alpha,p).
        \end{equation}
    \item\label{item:riesz-zero-is-delta} For $\alpha = 0$ we have
        \begin{equation}
            \label{eq:riesz-zero-is-delta}
            R^\pm_U(0,p) = \delta_p.
        \end{equation}
    \end{propositionlist}
\end{proposition}
\begin{proof}
    The first part is the literal translation of
    Proposition~\ref{proposition:riesz-distribution-identities}
    \refitem{item:riesz-distribution-times-eta} together with the fact
    that $\eta_p = \widetilde{\eta}_p \circ \exp_p^{-1}$. For the
    second part we have to be slightly more careful: in general, the
    gradient operator $\gradient$ on $M$ with respect to $g$ is
    \emph{not} intertwined into the gradient operator on $T_pM$ with
    respect to the flat metric $g_p$ via $\exp_p$. This is only true
    for arbitrary functions if the metric $g$ is \emph{flat}.
    Nevertheless we have for $\RE(\alpha) > n$ on $I^\pm_U(p)$
    \begin{align*}
        2 \alpha \gradient R^\pm_U(\alpha+2,p)
        &= 2 \alpha c(\alpha+2,n) \gradient
        \left(
            \eta_p^{\frac{\alpha+2-n}{2}}
        \right) \\
        &= 2 \alpha c(\alpha+2,n) \frac{\alpha+2-n}{2}
        \eta_p^{\frac{\alpha-n}{2}} \gradient \eta_p \\
        &= c(\alpha,n) \eta_p^{\frac{\alpha-n}{2}} \gradient \eta_p \\
        &= \gradient \eta_p \cdot R^\pm_U(\alpha,p).
    \end{align*}
    Since for $\RE(\alpha) > n$ the distribution $R^\pm_U(\alpha+2,p)$
    is actually a $\Fun[1]$-function and since on $U \setminus
    I^\pm_U(p)$ the relation \eqref{eq:gradient-of-eta-times-riesz} is
    trivially fulfilled, \eqref{eq:gradient-of-eta-times-riesz} holds
    on $U$ in the sense of $\Fun[0]$-functions and thus also in the
    sense of distributions. The usual holomorphy argument shows that
    \eqref{eq:gradient-of-eta-times-riesz} holds for all $\alpha \in
    \mathbb{C}$. For the third part we repeat our considerations from
    Lemma~\ref{lemma:identities-for-riesz-distribution},
    \refitem{item:riesz-and-gradient}. We first consider $\RE(\alpha)
    > n+2$ whence $R^\pm_U(\alpha+2,p)$ is $\Fun[2]$,
    $R^\pm_U(\alpha,p)$ is $\Fun[1]$, and we can compute $\dAlembert$
    in the sense of functions. On $I^\pm_U(p)$ we have
    \begin{align*}
        \dAlembert R^\pm_U(\alpha+2,p)
        &= \divergenz ( \gradient R^\pm_U(\alpha+2,p)) \\
        &\stackrel{\mathclap{\eqref{eq:gradient-of-eta-times-riesz}}}{=}
        \quad
        \divergenz
        \left(
            \frac{1}{2\alpha} \gradient \eta_p \cdot R^\pm_U(\alpha,p)
        \right) \\
        &= \frac{1}{2\alpha}
        g( \gradient R^\pm_U(\alpha,p), \gradient \eta_p)
        +
        \frac{1}{2\alpha} R^\pm_U(\alpha,p) \dAlembert \eta_p \\
        &\stackrel{\mathclap{\eqref{eq:gradient-of-eta-times-riesz}}}{=}
        \quad
        \frac{1}{2\alpha} g \left(
            \frac{1}{2(\alpha-2)}
            \gradient \eta_p \cdot R^\pm_U(\alpha-2,p),
            \gradient \eta_p
        \right)
        +
        \frac{1}{2\alpha} \dAlembert \eta_p \cdot R^\pm_U(\alpha,p) \\
        &\stackrel{\mathclap{\eqref{eq:length-of-gradient-of-dist-func}}}{=}
        \quad
        \frac{1}{2\alpha} \frac{1}{2(\alpha-2)}
        4 \eta_p R^\pm_U(\alpha-2,p)
        +
        \frac{1}{2\alpha} \dAlembert \eta_p \cdot R^\pm_U(\alpha,p) \\
        &\stackrel{\mathclap{\eqref{eq:identity-for-eta-times-riesz}}}{=}
        \quad
        \frac{1}{\alpha(\alpha-2)} (\alpha-2) (\alpha-2-n+2)
        R^\pm_U(\alpha,p)
        +
        \frac{1}{2\alpha} R^\pm_U(\alpha,p) \dAlembert \eta_p \\
        &= \left(
            \frac{\alpha-n}{\alpha} + \frac{1}{2\alpha}
            \dAlembert \eta_p
        \right)
        R^\pm_U(\alpha,p) \\
        &= \left(
            \frac{\dAlembert \eta_p - 2n}{2\alpha} + 1
        \right)
        R^\pm_U(\alpha,p).
    \end{align*}
    Since for $\RE(\alpha) > n+2$
    Equation~\eqref{eq:dAlembert-of-riesz-on-U-alpha-not-zero} is an
    equality between at least continuous functions, we have shown
    \eqref{eq:dAlembert-of-riesz-on-U-alpha-not-zero} since on $U
    \setminus J^\pm_U(p)$ we trivially have
    \eqref{eq:dAlembert-of-riesz-on-U-alpha-not-zero} as both sides
    are identically zero. Thus
    \eqref{eq:dAlembert-of-riesz-on-U-alpha-not-zero} holds for
    $\RE(\alpha) > n+2$ and by the obvious holomorphy in $\alpha
    \in \mathbb{C} \setminus \{0\}$ of both sides it holds for all
    $\alpha \neq 0$. Finally, we have
    \begin{align*}
        R^\pm_U(0,p) (\varphi)
        & = R^\pm(0) (\widetilde{\varrho}_p \exp_p^* \varphi)
        = \delta_0 (\widetilde{\varrho}_p \exp_p^* \varphi)
        = \widetilde{\varrho}_p(0) \cdot \varphi (\exp_p(0))
        = 1 \cdot \varphi(p)
        = \delta_p(\varphi),
    \end{align*}
    since $\widetilde{\varrho}_p(0) = 1$.
\end{proof}

Note that in the flat case we have $\dAlembert \eta_p = 2n$ whence
\eqref{eq:dAlembert-of-riesz-on-U-alpha-not-zero} simplifies to
$\dAlembert_{\textrm{flat}} R^\pm_{\textrm{flat}} (\alpha+2,p) =
R^\pm_{\textrm{flat}}(\alpha,p)$ from which we deduced that
$R^\pm_{\textrm{flat}}(2,p)$ is the Green function to
$\dAlembert_{\textrm{flat}}$ in
Theorem~\ref{theorem:green-function-of-dAlembert}. However, in the
general situation we have
\begin{equation}
    \label{eq:dAlembert-of-distance-p}
    \dAlembert \eta_p
    = 2n + g (\gradient \log \varrho_p, \gradient \eta_p )
\end{equation}
by our computation in
Proposition~\ref{proposition:properties-of-riesz-on-U},
\refitem{item:dAlembert-of-distant-function}. This additional term is
responsible for the failure of $R^\pm_U(2,p)$ to be a Green function
at $p$.

In order to determine the support and singular support of
$R^\pm_U(\alpha,p)$ we recall that under $\exp_p$ the chronological
future and past $I^\pm(0)$ of $0 \in T_pM$ are mapped to
$I^\pm_U(p)$. The same holds for $J^\pm(0)$ and $J^\pm_U(p)$ since
$\exp_p$ is assumed to be a diffeomorphism on the neighborhood $U$ of
$p$. Then the following statement is again a direct consequence of
Proposition~\ref{proposition:support-of-riesz-distribution}.
\begin{proposition}[Support and singular support of $R^\pm_U(\alpha,p)$]
    \label{proposition:support-of-riezs-on-U}
    \index{Riesz distribution!on domain!support}%
    \index{Riesz distribution!on domain!singular support}%
    Let $U \subseteq M$ be star-shaped around $p \in M$ and let
    $\alpha \in \mathbb{C}$.
    \begin{propositionlist}
    \item\label{item:singsupp-of-riesz-on-U-nonexceptional} If
        $\alpha$ is not exceptional then
        \begin{equation}
            \label{eq:supp-of-riesz-on-U-nonexc}
            \supp R^\pm_U(\alpha,p) = J^\pm_U(p)
        \end{equation}
        and
        \begin{equation}
            \label{eq:singsupp-of-riesz-on-U-nonexc}
            \singsupp R^\pm_U(\alpha,p) \subseteq \partial I^\pm_U(p).
        \end{equation}
    \item\label{item:singsupp-of-riesz-on-U-exceptional} If $\alpha$
        is exceptional then
        \begin{equation}
            \label{eq:singsupp-supp-of-U-excpetional}
            \singsupp R^\pm_U(\alpha,p) = \supp R^\pm_U(\alpha,p)
            \subseteq \partial I^\pm_U(p).
        \end{equation}
    \item\label{item:singsupp-riesz-on-U-big-dimension} If $n \geq 3$
        and $\alpha \in \left\{ n-2k \big| k \in \mathbb{N}_0, k <
        \frac{n}{2} \right\}$ we have
    \begin{equation}
        \label{eq:singsupp-riesz-on-U-big-dimension}
        \singsupp R^\pm_U(\alpha,p) = \supp R^\pm_U(\alpha,p)
        = \partial I^\pm_U(p).
    \end{equation}
    \end{propositionlist}
\end{proposition}
\begin{proof}
    This follows from
    Proposition~\ref{proposition:support-of-riesz-distribution} and
    the general behaviour of $\supp$ and $\singsupp$ under
    push-forwards with diffeomorphisms and multiplication with
    positive smooth functions.
\end{proof}

\begin{proposition}[Order of $R^\pm_U(\alpha,p)$]
    \label{proposition:order-of-riesz-on-U}
    \index{Riesz distribution!on domain!order}%
    Let $U \subseteq M$ be star-shaped around $p \in M$ and let
    $\alpha \in \mathbb{C}$.
    \begin{propositionlist}
    \item \label{item:order-for-large-alpha} If $\RE(\alpha) > n$ then
        $\ord_U(R^\pm_U(\alpha,p)) = 0$.
    \item \label{item:boundary-for-order} The global order of
        $R^\pm_U(\alpha,p)$ is bounded by $2k$ where $k \in
        \mathbb{N}_0$ is such that $\RE(\alpha) + 2k > n$.
    \item \label{item:boundary-for-order-for-alpha-greater-zero} If
        $\RE(\alpha) > 0$ then the global order of $R^\pm_U(\alpha,p)$
        is bounded by $n$ if $n$ is even and by $n+1$ if $n$ is odd.
    \end{propositionlist}
\end{proposition}
\begin{proof}
    The order of a distribution does not change under push-forwards
    with diffeomorphisms and multiplication with positive smooth
    functions. Thus the result follows directly from
    Proposition~\ref{proposition:order-of-riesz-distribution}.
\end{proof}

\begin{proposition}[Reality of $R^\pm_U(\alpha,p)$]
    \label{proposition:reality-of-riesz-on-U}
    \index{Riesz distribution!on domain!reality}%
    Let $U \subseteq M$ be star-shaped around $p \in M$ and let
    $\alpha \in \mathbb{C}$. Then we have
    \begin{equation}
        \label{eq:reality-of-riesz-on-U}
        \cc{R^\pm_U(\alpha, p)} = R^\pm_U(\cc{\alpha}, p).
    \end{equation}
\end{proposition}
\begin{proof}
    Since $\widetilde{\varrho}_p = \cc{\widetilde{\varrho}_p} > 0$
    this follows from
    Proposition~\ref{proposition:realitiy-of-riesz-distribution}.
\end{proof}

In a next step we need to understand how the Riesz distribution
$R^\pm_U(\alpha,p)$ depends on the point $p \in M$. To this end we
have to be slightly more specific with our definition of
$R^\pm_U(\alpha,p)$. In order to compare
\eqref{eq:riesz-distribution-on-U} for different $p$ it is convenient
to choose a common reference Minkowski spacetime. Thus we consider
the following situation: assume that $U$ is not only star-shaped with
respect to $p$ but also with respect to $p' \in O$ where $O \subseteq
U$ is a small open neighborhood of $p$. In particular, if $U$ is even
geodesically convex then we can choose $O=U$. Moreover, let $e_1,
\ldots, e_n$ be a smooth Lorentz frame on $U$ inducing isometric
isomorphisms
\begin{equation}
    \label{eq:isometric-lorentz-isomorphisms}
    I_{p'}: (T_{p'}M, g_{p'}, \uparrow) \longrightarrow
    (\mathbb{R}^n, \eta, \uparrow)
\end{equation}
preserving the time orientation. Clearly, $I_{p'}$ depends smoothly on
$p'$ in this case. Then for $\varphi \in \Cinfty_0(U)$ we have for all
$p' \in O$
\begin{equation}
    \label{eq:riesz-distribution-at-p-prime}
    R^\pm_U(\alpha,p') (\varphi)
    = R^\pm(\alpha)
    \left(
        I_{p'*} (\widetilde{\varrho}_{p'} \exp_{p'}^* \varphi)
    \right)
\end{equation}
with $R^\pm(\alpha)$ being the Riesz distributions on $\mathbb{R}^n$,
independent of $p'$.
\begin{lemma}
    \label{lemma:support-of-isometric-pullback-of-varphi}
    Let $K \subseteq U$ be compact. Then for every compact subset $L
    \subseteq O$ there exists a compactum $\widetilde{K} \subseteq
    \mathbb{R}^n$ such that
    \begin{equation}
        \label{eq:support-of-isometric-pullback-of-varphi}
        \supp \left(
            I_{p'*} (\widetilde{\varrho}_{p'} \exp_{p'}^* \varphi)
        \right)
        \subseteq \widetilde{K}
    \end{equation}
    for all $\varphi \in \Cinfty_K(U)$ and all $p' \in L$.
\end{lemma}
\begin{proof}
    For all $p' \in O$ the function $x \mapsto \left( I_{p'*}
        (\widetilde{\varrho}_{p'} \exp_{p'}^* \varphi) \right) (x)$ is
    a compactly supported smooth function on $\mathbb{R}^n$. Since
    $I_{p'}$ is a linear isomorphism and $\widetilde{\varrho}_{p'}$ is
    strictly positive,
    \[
    K_{p'} = \supp \left(
        I_{p'*} (\widetilde{\varrho}_{p'} \exp_{p'}^* \varphi)
    \right)
    = I_{p'} ( \exp_{p'}^{-1} (\supp \varphi) )
    \]
    by the general behaviour of supports under diffeomorphisms. The
    various compacta $K_{p'}$ depend on $p'$ in a continuous way. More
    precisely, there is a map $\Phi_{p'}: \widetilde{V} \subseteq
    \mathbb{R}^n \longrightarrow \mathbb{R}^n$ such that $K_{p'} =
    \Phi_{p'} (K_p)$ which depends continuously on $p'$. In fact,
    define
    \[
    \Phi(p',x)
    = I_{p'} ( \exp_{p'}^{-1} ( \exp_p (I_p^{-1}(x))))
    \]
    for $x \in I_p(V) \subseteq \mathbb{R}^n$. Then $\Phi: O \times V
    \longrightarrow \mathbb{R}^n$ is even smooth. Now for $\{p'\}
    \subseteq O$ compact we have $K_{p'} \subseteq \Phi (\{p'\} \times
    K_p)$ and thus $\bigcup_{p'} K_{p'} \subseteq \Phi \left(
        \bigcup_{p'} \{p'\} \times K_p \right)$. If $p' \in L$ runs
    through a compact subset $L \subseteq O$ then the union of the
    $K_{p'}$ is contained in a compactum itself since $\Phi$ is
    continuous. This is the $\widetilde{K}$ we are looking for.
\end{proof}

Using this lemma we see that the support of $I_{p'*}
(\widetilde{\varrho}_{p'} \exp_{p'} \varphi)$ is uniformly contained
in some compactum in $\mathbb{R}^n$. This allows to use the continuity
of the distributions $R^\pm(\alpha)$ to obtain the following result:
\begin{proposition}
    \label{proposition:riesz-dependence-on-base-point-p-prime}
    Let $U \subseteq M$ be star-shaped around $p \in M$ and let $O
    \subseteq U$ be an open neighborhood of $U$ such that $U$ is
    star-shaped around every  $p' \in O$.
    \begin{propositionlist}
    \item \label{item:uniform-cont-estimate-for-riesz-on-U} For every
        compacta $K \subseteq U$ and $L \subseteq O$ there exists a
        constant $c_{K,L,\alpha} > 0$ such that
        \begin{equation}
            \label{eq:uniform-cont-estimate-for-riesz-on-U}
            |R^\pm_U(\alpha,p') (\varphi)|
            \leq c_{K,L,\alpha} \seminorm[K,2k] (\varphi)
        \end{equation}
        for all $\varphi \in \Cinfty_K(U)$ and $p' \in L$ where $k \in
        \mathbb{N}_0$ is such that $\RE(\alpha) + 2k > n$.
    \item \label{item:uniform-riesz-estimate-2} In particular, for
        $\RE(\alpha) > 0$ and every compacta $K \subseteq U$ and $L
        \subseteq O$ there exists a constant $c_{K,L,\alpha} > 0$ such
        that
        \begin{equation}
            \label{eq:uniform-riesz-estimate-2}
            |R^\pm_U(\alpha,p')(\varphi)|
            \leq c_{K,L,\alpha} \seminorm[K,n+1] (\varphi)
        \end{equation}
        for all $\varphi \in \Cinfty_K(U)$.
    \item \label{item:riesz-dependence-on-p-prime} Let $k \in
        \mathbb{N}_0$ satisfy $\RE(\alpha) + 2k > n$. Then for every
        $\Phi \in \Fun[2k+\ell]_0(O \times U)$ the map
        \begin{equation}
            \label{eq:riesz-dependence-on-p-prime}
            O \ni p'
            \; \mapsto \;
            R^\pm_U(\alpha,p') (\Phi(p', \argument))
            \in \mathbb{C}
        \end{equation}
        is $\Fun[\ell]$ on $O$.
    \item \label{item:riesz-dependence-on-p-prime-2} Again, for
        $\RE(\alpha) > 0$ and $\Phi \in \Fun[n+1+\ell](O \times U)$
        the corresponding map \eqref{eq:riesz-dependence-on-p-prime}
        is $\Fun[\ell]$ on $O$.
    \item \label{item:riesz-is-holomorphic} Let $\varphi \in
        \Fun_0(U)$ then the map
        \begin{equation}
            \label{eq:ries-is-holomorphic}
            \alpha \; \mapsto \; R^\pm_U(\alpha,p)(\varphi)
        \end{equation}
        is holomorphic for $\RE(\alpha) > n - 2
        \left[\frac{k}{2}\right]$.
    \item \label{item:riesz-smooth-on-p-prime} If $\Phi \in \Cinfty(O
        \times U)$ is even smooth and has support $\supp \Phi
        \subseteq O \times K$ with some compact $K$, then the function
        \begin{equation}
            \label{eq:riesz-smooth-on-p-prime}
            O \ni p'
            \; \mapsto \;
            R^\pm_U(\alpha,p') (\Phi(p', \argument))
        \end{equation}
        is smooth on $O$.
    \end{propositionlist}
\end{proposition}
\begin{proof}
    By Lemma~\ref{lemma:support-of-isometric-pullback-of-varphi} we
    have a compact subset $\widetilde{K} \subseteq \mathbb{R}^n$ such
    that
    \[
    \supp \left(
        I_{p'*} (\widetilde{\varrho}_{p'} \exp_{p'}^* \varphi)
    \right)
    \subseteq \widetilde{K}
    \]
    for all $p' \in L$ and $\varphi \in \Cinfty_K(U)$. Thus by
    continuity of $R^\pm(\alpha)$ and the fact that $R^\pm(\alpha)$
    has order $\leq 2k$ whenever $\RE(\alpha) + 2k > n$, see
    Proposition~\ref{proposition:order-of-riesz-distribution},
    \refitem{item:riesz-global-order}, we have
    \[
    \left|
        R^\pm_U(\alpha,p')(\varphi)
    \right|
    =
    \left|
        R^\pm(\alpha)
        (I_{p'*} (\widetilde{\varrho}_{p'} \exp_{p'}^* \varphi))
    \right|
    \leq
    c \seminorm[\widetilde{K},2k]
    (I_{p'*} (\widetilde{\varrho}_{p'} \exp_{p'}^* \varphi))
    =
    c' \seminorm[K,2k](\varphi),
    \]
    since $I_{p'*} \widetilde{\varrho}_{p'}$ is bounded with all its
    derivatives on the compactum $\widetilde{K}$ as it is smooth
    anyway, and $I_{p'*} \exp_{p'}^* \varphi$ is also smooth on
    $\widetilde{K}$. Since the exponential map $\exp_{p'}$ also
    depends smoothly on $p'$ all its derivatives up to order $2k$ are
    bounded as long as $p' \in L$, the same holds for $I_{p'}$. This
    gives the new constant $c'$ independent of $p'$ but only depending
    on $L$. This proves the first part. The second follows since for
    $\RE(\alpha) > 0$ the order of $R^\pm(\alpha)$ is bounded by $n+1$
    by Proposition~\ref{proposition:order-of-riesz-distribution},
    \refitem{item:riesz-order-if-alpha-nonnegative}. The third part
    follows immediately from the technical
    Lemma~\ref{lemma:support-of-isometric-pullback-of-varphi} and a
    careful counting of the number of derivatives needed in the proof
    of that lemma, see also
    Proposition~\ref{proposition:geometric-parameter-differentiation}. The
    fourth part is a particular case thereof. The holomorphy follows
    immediately from Remark~\ref{remark:riesz-on-Ck-functions}.  For
    the last part note that by definition of $R^\pm_U(\alpha,p')$ we
    have
    \[
    R^\pm_U(\alpha,p') (\Phi(p', \argument))
    = R^\pm(\alpha)
    \left(
        I_{p'*}
        (\widetilde{\varrho}_{p'} \exp_{p'}^* \Phi(p', \argument))
    \right),
    \]
    and the function
    \[
    (p', x) \; \mapsto \;
    I_{p'*} (\widetilde{\varrho}_{p'} \exp_{p'}^* \Phi(p', \argument))
    \At{x}
    \]
    has support in $O \times \widetilde{K}$ with $\widetilde{K}
    \subseteq \mathbb{R}^n$ compact. Moreover, by the smooth choice of
    $I_{p'}$ and the smoothness of $\widetilde{\varrho}$ and $\exp$ we
    conclude that it is smooth in \emph{both} variables. Thus we can
    apply Lemma~\ref{lemma:parameter-differentiation-under-integral}
    to obtain the smoothness of \eqref{eq:riesz-smooth-on-p-prime}.
\end{proof}

In particular, it follows from the fourth part that the map $p'
\mapsto (\id \tensor R^\pm_U)(\alpha,p')(\Phi)$ is smooth on $O$ for
$\Phi \in \Cinfty_0(O \times U)$.

Let us now discuss an additional symmetry property of the Riesz
distributions. In the flat case the exponential map
\begin{equation}
    \label{eq:exp-map}
    \exp_p: T_pM \longrightarrow M
\end{equation}
is just the translation, i.e. for $(M, g) = (\mathbb{R}^n, \eta)$ we
have
\begin{equation}
    \label{eq:exp-flat-case}
    \exp_p(v) = p + v.
\end{equation}
Thus in this case for $\RE(\alpha) > 0$ we have
\begin{equation}
    \label{eq:riesz-flat-depends-only-on-difference}
    R^\pm(\alpha,p)(q)
    = (\exp_{p*} R^\pm(\alpha))(q)
    = R^\pm(\alpha) \left( \exp_p^{-1}(q) \right)
    = R^\pm(\alpha) (q-p).
\end{equation}
In particular,
\begin{equation}
    \label{eq:riesz-flat-symmetry}
    R^\pm(\alpha,p)(q) = R^\mp(\alpha,q)(p)
\end{equation}
follows since $q-p \in I^+(0)$ iff $p-q \in I^-(0)$ and the function
$\eta$ is invariant under total inversion $x \mapsto -x$. While the
phrase ``$R^\pm(\alpha,p)(q)$ depends only on the difference $q-p$''
clearly only makes sense on a vector space, the symmetry feature
\eqref{eq:riesz-flat-symmetry} remains to be true also in the
geometric context. Of course, now we have to take care that the points
$p$ and $q$ enter equally in \eqref{eq:riesz-flat-symmetry} whence the
domain $U$ has to be star-shaped with respect to both. But then we
have the following statement:
\begin{proposition}[Symmetry of $R^\pm_U(\alpha,p)$]
    \label{proposition:symmetry-of-riesz-on-U}
    \index{Riesz distribution!on domain!symmetry}%
    Let $U \subseteq M$ be geodesically convex and $\alpha \in
    \mathbb{C}$.
    \begin{propositionlist}
    \item \label{item:symmetry-of-riesz-on-U-as-function} If
        $\RE(\alpha) > n$ then
        \begin{equation}
            \label{eq:symmetry-of-riesz-on-U-as-function}
            R^\pm_U(\alpha,p)(q) = R^\mp_U(\alpha,q)(p)
        \end{equation}
        for all $p,q \in U$.
    \item \label{item:symmetry-of-riesz-on-U-as-distribution} For all
        $\Phi \in \Cinfty_0(U \times U)$ one has
        \begin{equation}
            \label{eq:symmetry-of-riesz-on-U-as-distribution}
            \int_U R^\pm(\alpha,p) (\Phi(p, \argument)) \: \mu_g(p)
            =
            \int_U R^\mp(\alpha,q) (\Phi(\argument, q)) \: \mu_g(q).
        \end{equation}
    \end{propositionlist}
\end{proposition}
\begin{proof}
    First we note that thanks to the convexity of $U$ the Riesz
    distributions $R^\pm_U(\alpha,p)$ are defined for all $p \in
    U$. For $\RE(\alpha) > n$ the Riesz distributions are continuous
    functions explicitly given
    by~\eqref{eq:riesz-on-U-for-large-alpha} in
    Proposition~\ref{proposition:properties-of-riesz-on-U},
    \refitem{item:riesz-on-U-is-cont-for-large-alpha}. We compute
    \begin{align*}
        \eta_p(q)
        &= g_p \left(
            \exp_p^{-1}(q), \exp_p^{-1}(q)
        \right) \\
        &= g_{\exp_p(\exp_p^{-1}(q))} \left(
            T_{\exp_p^{-1}(q)} \exp_p (\exp_p^{-1}(q)),
            T_{\exp_p^{-1}(q)} \exp_p (\exp_p^{-1}(q))
        \right) \\
        &= g_q \left(
            T_{\exp_p^{-1}(q)} \exp_p (\exp_p^{-1}(q)),
            T_{\exp_p^{-1}(q)} \exp_p (\exp_p^{-1}(q))
        \right)
    \end{align*}
    by the Gauss Lemma. Now $v = \exp_p^{-1}(q)$ is the tangent vector
    of the geodesic $t \mapsto \exp_p(tv)$ which starts at $p$ and
    reaches $q$ at $t=1$. Reversing the time the curve $\tau \mapsto
    \exp_p((1-\tau)v)$ is still a geodesic which now starts at $q$ for
    $\tau=0$ and reaches $p$ at $\tau=1$. Thus the tangent vector of
    this geodesic is uniquely fixed to be $\exp_p^{-1}(q)$ since in
    the convex $U$ the exponential map $\exp_p$ is a diffeomorphism.
    On the other hand, by the chain rule it follows that
    \[
    \frac{\D}{\D \tau} \At{\tau=0} \exp_p((1-\tau)v)
    = T_v \exp_p(-v) = -T_v \exp_p(v),
    \]
    whence we have shown
    \[
    \exp_q^{-1}(p) = -T_{\exp_p^{-1}(q)} \exp_p (\exp_p^{-1}(q)).
    \tag{$*$}
    \]
    \begin{figure}
        \centering
        \input{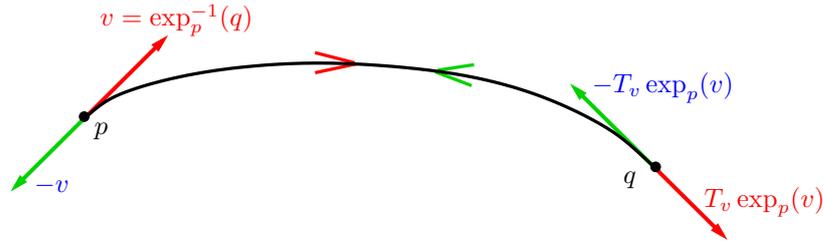}
        \caption{\label{fig:geodesic-backwards}%
          A geodesic running backwards.
        }
    \end{figure}
    It follows that
    \[
    \eta_p(q)
    = g_q ( \exp_q^{-1}(p), \exp_q^{-1}(p) )
    = \eta_q(p).
    \]
    \index{Lorentz distance square}%
    Since $\eta_p(q)$ is something like the ``Lorentz distance
    square'' it is not surprising that this quantity is symmetric in
    $p$ and $q$: everything else would be rather disturbing. Since we
    have a relative sign in ($*$) we see that if $\exp_p^{-1}(q) \in
    I^+(0_p)$ then the geodesic $t \mapsto \exp_p (t \exp_p^{-1}(q))$
    is future directed for all times whence $\exp_p^{-1}(q) \in
    I^-(0_q)$ is past directed. From
    Figure~\ref{fig:geodesic-backwards} this is clear. But then
    \eqref{eq:symmetry-of-riesz-on-U-as-function} follows directly
    from \eqref{eq:riesz-on-U-for-large-alpha} since the prefactors
    $c(\alpha,n)$ are the same for the advanced and retarded Riesz
    distributions. For the second part we first consider $\RE(\alpha)
    > n$. Then $R^\pm_U(\alpha,p)(q)$ is a continuous function on $U
    \times U$ since $\eta_p(q)$ is smooth in both variables. Thus
    $R^\pm_U(\alpha,p)(\argument)$ is locally integrable and hence the
    function
    \[
    (p,q) \; \mapsto \; R^\pm_U(\alpha,p)(q) \Phi(p,q)
    \]
    has compact support and is continuous. Thus we apply Fubini's
    theorem and interchange the $q$- and $p$-integrations
    \begin{align*}
        \int_U R^\pm_U(\alpha,p) (\Phi(p, \argument)) \: \mu_g(p)
        &=
        \int_U \int_U R^\pm_U(\alpha,p)(q) \Phi(p,q)
        \: \mu_g(q) \mu_g(p) \\
        &=
        \int_U \int_U R^\pm_U(\alpha,p)(q) \Phi(p,q)
        \: \mu_g(p) \mu_g(q) \\
        &\stackrel{\mathclap{\eqref{eq:symmetry-of-riesz-on-U-as-function}}}{=}
        \quad
        \int_U \int_U R^\mp_U(\alpha,q)(p) \Phi(p,q)
        \: \mu_g(p) \mu_g(q)
        \\
        &=
        \int_U R^\mp_U(\alpha,q) (\Phi(\argument, q)) \: \mu_g(q),
    \end{align*}
    which proves \eqref{eq:symmetry-of-riesz-on-U-as-distribution} for
    $\RE(\alpha) > n$. For general $\alpha \in \mathbb{C}$ we notice
    that the integrands of both sides are compactly supported smooth
    function on $U$ thanks to
    Proposition~\ref{proposition:geometric-parameter-differentiation}
    and Remark~\ref{remark:smooth-reduced-section}. Thus the usual
    Morera type argument shows that both sides are holomorphic
    functions of $\alpha$ since the integrands are holomorphic in
    $\alpha$ and we exchange the integrations $\int_U$ and
    $\int_\Delta \D \alpha$ as usual: by holomorphy we conclude that
    the equality \eqref{eq:symmetry-of-riesz-on-U-as-distribution}
    holds for all $\alpha$ as it holds for $\RE(\alpha) > n$.
\end{proof}


%% file: hadamard.tex
%
%

Differently from the flat situation, the Riesz distribution
$R^\pm_U(2,p)$ does not yield a fundamental solution for
$\dAlembert$. Indeed, we cannot evaluate $\dAlembert R^\pm_U(2,p)$ as
we did in the flat case since in
Proposition~\ref{proposition:differentiation-of-riesz-on-U} we had to
exclude the value of $\alpha$ needed for $\dAlembert R^\pm_U(2,p)$
explicitly. Instead, from
\begin{equation}
    \label{eq:dAlembert-of-riesz-on-U-again}
    \dAlembert R^\pm_U(\alpha+2,p)
    = \left(
        \frac{\dAlembert \eta_p - 2n}{2\alpha} + 1
    \right)
    R^\pm_U(\alpha,p),
\end{equation}
valid for $\alpha \neq 0$ we only see the following: The limit $\alpha
\longrightarrow 0$ of the right hand side, which would be the
interesting point, is problematic. One has $R^\pm_U(0,p) = \delta_p$
but the prefactor itself is singular, at least on first
sight. However, the simple pole in $\frac{\dAlembert \eta_p -
  2n}{2\alpha}$ is not as dangerous as it seems. In fact, we
\emph{know} that $\alpha \mapsto R^\pm_U(\alpha+2,p)(\dAlembert
\varphi)$ is holomorphic on the whole complex plane. Hence the limit
$\alpha \longrightarrow 0$ of the left hand side certainly
exists. Thus we \emph{do} have an analytic continuation of the right
hand side for $\alpha = 0$, the singularity was not present after
all. However, the precise value at $\alpha = 0$ is hard to obtain and
not just $\delta_p$. Of course, we know it is $\dAlembert
R^\pm_U(2,p)$, but this does not help.

Thus one proceeds differently. The Ansatz is to use all Riesz
distributions $R^\pm_U(2+2k,p)$ and approximate the true Green
function by a series in the $R^\pm_U(2+2k,p)$ for $k \in \mathbb{N}_0$
with appropriate coefficients. These coefficients are the Hadamard
coefficients we are going to determine now. The expansion we obtain
can be thought of as an expansion of the Green functions in increasing
regularity as the $R^\pm_U(2+2k,p)$ become more and more regular for
$k \longrightarrow \infty$.

%
%

\subsection{The Ansatz for the Hadamard Coefficients}
\label{subsec:AnsatzHadamardCoefficients}

The setting will the following. We consider a normally hyperbolic
differential operator $D = \dAlembert^\nabla + B$ on some vector
bundle $E \longrightarrow M$ over $M$ with induced connection
$\nabla^E$ and $B \in \Secinfty(\End(E))$ as in
Section~\ref{subsec:normally-hyperbolic-diffops}. Moreover, for $p \in
M$ we choose a geodesically star-shaped open neighborhood $U \subseteq
M$ on which $R^\pm_U(\alpha,p)$ is defined as before. According to our
convention for distributions, the Green functions are now generalized
sections
\begin{equation}
    \label{eq:green-functions-as-gensecs}
    \mathcal{R}^\pm(p) \in \Sec[-\infty](E) \tensor E_p^*,
\end{equation}
as we take care of the density part using $\mu_g$. The pairing with a
test section $\varphi \in \Secinfty_0(E^*)$ yields then an element in
$E_p^*$. The equation to solve is
\begin{equation}
    \label{eq:green-function-equation}
    D \mathcal{R}^\pm(p) = \delta_p,
\end{equation}
where $\delta_p$ is viewed as $E_p^*$-valued distribution on
$\Secinfty_0(E^*)$ and $D \mathcal{R}^\pm(p)$ is defined as usual.

The Ansatz for $\mathcal{R}^\pm(p)$ is now the following. Since the
$R^\pm_U(\alpha,p)$ have increasing regularity for increasing
$\RE(\alpha)$ we try a series
\begin{equation}
    \label{eq:Ansatz-for-Green-function}
    \index{Ansatz!Green function}%
    \mathcal{R}^\pm(p)
    = \sum_{k=0}^\infty V_p^k R^\pm_U(2+2k, p)
\end{equation}
with \emph{smooth} sections
\begin{equation}
    \label{eq:Hadamard-coefficients}
    V_p^k \in \Secinfty\left(E\at{U}\right) \tensor E^*_p.
\end{equation}
Then \eqref{eq:Ansatz-for-Green-function} should be thought of as an
expansion with respect to regularity. The starting point for $k=0$
will be the most singular term coming from $R^\pm_U(2,p)$. Of course,
such an Ansatz can hardly be expected to work just like that. Even if
we can find reasonable $V_p^k$ such that
\eqref{eq:green-function-equation} holds ``in each order of
regularity'', the series \eqref{eq:Ansatz-for-Green-function} has to
be shown to converge: In fact, this will not be the case (except for
some very particular cases) whence we have to go a step beyond
\eqref{eq:Ansatz-for-Green-function}. However, for the time being we
shall investigate the Ansatz \eqref{eq:Ansatz-for-Green-function}.

First we note that a scalar distribution like $R^\pm_U(\alpha,p)$ can
be multiplied with a smooth section like $V_p^k$ and yields a
distributional section
\begin{equation}
    \label{eq:riesz-times-hadamard-is-distributional-section}
    V_p^k R^\pm_U(2+2k, p) 
    \in \Secinfty_0(E^*)' \tensor E_p^*
    = 
    \Sec[-\infty](E) \tensor E_p^*.
\end{equation}
In Remark~\ref{remark:module-structure-on-generalized-section} it is
only necessary that one factor of the product is actually smooth. We
compute now \eqref{eq:green-function-equation}. First we assume that
the series \eqref{eq:Ansatz-for-Green-function} converges at least in
the weak$^*$ topology so that we can apply $D$ componentwise. This
yields
\begin{equation}
    \label{eq:componentwise-D-on-Ansatz}
    \begin{split}
        D \mathcal{R}^\pm(p)
        & = D \sum_{k=0}^\infty V_p^k R^\pm_U(2+2k,p) \\
        & = \sum_{k=0}^\infty D \left(
            V_p^k R^\pm_U(2+2k, p)
        \right) \\
        & = \sum_{k=0}^\infty \left(
            D(V_p^k) R^\pm_U(2+2k,p)
            + 2 \nabla^E_{\gradient R^\pm_U(2+2k,p)} V_p^k
            + V_p^k \dAlembert R^\pm_U(2+2k,p)
        \right)
    \end{split}
\end{equation}
by the Leibniz rule of a normally hyperbolic differential operator as
in Remark~\ref{remark:NormallyHyperbolic},
\refitem{item:NormallyHyperbolicLeibniz}. Note that in
\eqref{eq:NormallyHyperbolicLeibniz} it is sufficient that one of the
factors is smooth. Inserting the properties of $R^\pm_U(\alpha,p)$
from Proposition~\ref{proposition:properties-of-riesz-on-U} yields
then
\begin{equation}
    \label{eq:D-on-Ansatz-2}
    \begin{split}
        D \mathcal{R}^\pm(p)
        & = D(V_p^0) R^\pm_U(2,p)
        + 2 \nabla^E_{\gradient R^\pm_U(2,p)} V_p^0
        + V_p^0 \dAlembert R^\pm_U(2,p) \\
        & \quad + \sum_{k=1}^\infty \left(
            D(V_p^k) R^\pm_U(2+2k,p)
            + 2 \nabla^E_{\frac{1}{4k} R^\pm_U(2k,p) \gradient \eta_p}
            V_p^k
            + V_p^k \left(
                \frac{\dAlembert \eta_p - 2n}{4k} + 1
            \right) R^\pm_U(2k,p)
        \right) \\
        & = 2 \nabla^E_{\gradient R^\pm_U(2,p)} V_p^0 
        + V_p^0 \dAlembert R^\pm_U(2,p)
        + \sum_{k=0}^\infty D(V_p^k) R^\pm_U(2+2k,p) \\
        & \quad + \sum_{k=1}^\infty \left(
            2 \nabla^E_{\frac{1}{4k} \gradient \eta_p} V_p^k
            + V_p^k \left(
                \frac{\dAlembert \eta_p - 2n}{4k} + 1
            \right)
        \right) R^\pm_U(2k,p) \\
        &  = 2 \nabla^E_{\gradient R^\pm_U(2,p)} V_p^0 
        + V_p^0 \dAlembert R^\pm_U(2,p) \\
        & \quad +\sum_{k=1}^\infty \left(
            D(V_p^{k-1})
            + 2 \nabla^E_{\frac{1}{4k} \gradient \eta_p} V_p^k
            + \left(
                \frac{\dAlembert \eta_p - 2n}{4k} + 1
            \right) V_p^k
        \right) R^\pm_U(2k,p).
    \end{split}
\end{equation}
We view \eqref{eq:D-on-Ansatz-2} as an expansion with respect to
regularity. Thus, we ask for \eqref{eq:componentwise-D-on-Ansatz} in
each ``order'', i.e. \eqref{eq:componentwise-D-on-Ansatz} should be
fulfilled for each component in front of the $R^\pm_U(2k,p)$. This
yields the following equations. In lowest order we have for $V_p^0$
the equation
\begin{equation}
    \label{eq:lowest-hadamard-equation}
    2 \nabla^E_{\gradient R^\pm_U(2,p)} V_p^0 
    + V_p^0 \dAlembert R^\pm_U(2,p)
    = \delta_p,
\end{equation}
while for $k \geq 1$ we have the recursive equations
\begin{equation}
    \label{eq:higher-hadamard-equations}
    \frac{1}{2k} \nabla^E_{\gradient \eta_p} V_p^k
    + \left(
        \frac{\dAlembert \eta_p - 2n}{4k} + 1
    \right) V_p^k
    = - D(V_p^{k-1})    
\end{equation}
for $V_p^k$. Equivalently, we can write this for $k \geq 1$ as
\begin{equation}
    \label{eq:higher-hadamard-equations-2}
    \nabla^E_{\gradient \eta_p} V_p^k
    + \left(
        \frac{1}{2} \dAlembert \eta_p - n +2k
    \right) V_p^k
    = -2k D(V_p^{k-1}).
\end{equation}
Since \eqref{eq:higher-hadamard-equations-2} also makes sense for
$k=0$ it seems tempting to unify \eqref{eq:lowest-hadamard-equation}
and \eqref{eq:higher-hadamard-equations-2}. To this end, we take
\eqref{eq:higher-hadamard-equations-2} for $k=0$ and multiply this by
$R^\pm_U(\alpha,p)$ yielding
\begin{equation}
    \label{eq:lowest-hadamard-alternative}
    \nabla^E_{\gradient \eta_p R^\pm_U(\alpha,p)} V_p^0
    + \left(
        \frac{1}{2} \dAlembert \eta_p - n
    \right) V_p^0 R^\pm_U(\alpha,p)
    = 0,
\end{equation}
which is equivalent to
\begin{equation}
    \label{eq:lowest-hadamard-alternative-2}
    \nabla^E_{2\alpha \gradient R^\pm_U(\alpha+2,p)} V_p^0
    + \alpha \left(
        \dAlembert R^\pm_U(\alpha+2,p) - R^\pm_U(\alpha,p)
    \right) V_p^0
    = 0,
\end{equation}
by Proposition~\ref{proposition:differentiation-of-riesz-on-U}. Now we
can divide by $\alpha$ and obtain the condition
\begin{equation}
    \label{eq:lowest-hadamard-alternative-3}
    2 \nabla^E_{\gradient R^\pm_U(\alpha+2,p)} V_p^0
    + \left(
        \dAlembert R^\pm_U(\alpha+2,p) - R^\pm_U(\alpha,p)
    \right) V_p^0
    = 0,
\end{equation}
whose limit $\alpha \longrightarrow 0$ exists and is given by
\begin{equation}
    \label{eq:lowest-hadamard-alternative-4}
    2 \nabla^E_{\gradient R^\pm_U(2,p)} V_p^0
    + \left(
        \dAlembert R^\pm_U(2,p) - R^\pm_U(0,p)
    \right) V_p^0
    = 0,
\end{equation}
since $R^\pm_U(\alpha,p)$ is holomorphic in $\alpha$ for all $\alpha
\in \mathbb{C}$. Since moreover $R^\pm_U(0,p) = \delta_p$ we can
evaluate the condition \eqref{eq:lowest-hadamard-alternative-4}
further and obtain
\begin{equation}
    \label{eq:lowest-hadamard-alternative-5}
    2 \nabla^E_{\gradient R^\pm_U(2,p)} V_p^0
    + V_p^0 \dAlembert R^\pm_U(2,p)
    = V_p^0 \delta_p.
\end{equation}
Thus we conclude that \eqref{eq:higher-hadamard-equations-2} for $k=0$
implies \eqref{eq:lowest-hadamard-equation} iff $V_p^0(p) =
\id_{E_p}$. This motivates that we want to solve
\eqref{eq:lowest-hadamard-equation} with the additional requirement
\begin{equation}
    \label{eq:hadamard-initial-condition}
    V_p^0(p) = \id_{E_p},
\end{equation}
which we can view as an \emph{initial condition}. Indeed, all the
gradients $\gradient R^\pm_U(\alpha,p)$ are pointing in ``radial''
direction parallel to $\gradient \eta_p$ by
Proposition~\ref{proposition:differentiation-of-riesz-on-U}. Thus a
differential equation like \eqref{eq:lowest-hadamard-equation} should
have a unique solution once the value is fixed in the ``center'',
i.e. at $p$. Then one has just to follow the flow of $\gradient
\eta_p$ in order to determine the value elsewhere. Of course, this
geometric intuition has to be justified more carefully. In any case,
we take these heuristic considerations as motivation for the following
definition:
\begin{definition}[Transport equations]
    \label{definition:transport-equations}
    \index{Transport equation}%
    Let $k \in \mathbb{N}_0$ and let $D \in \Diffop^2(E)$ be normally
    hyperbolic. Then the recursive equations
    \begin{equation}
        \label{eq:transport-equations}
        \nabla^E_{\gradient \eta_p} V_p^k
        + \left(
            \frac{1}{2} \dAlembert \eta_p - n + 2k
        \right) V_p^k
        = - 2k D V_p^{k-1}
    \end{equation}
    together with the initial condition
    \begin{equation}
        \label{eq:transport-equations-initial-condition}
        V_p^0(p) = \id_{E_p}
    \end{equation}
    are called the transport equations for $V_p^k \in
    \Secinfty(E\at{U}) \tensor E^*_p$ corresponding to $D$.
\end{definition}
\begin{remark}[Transport equations]
    \label{remark:transport-equations}
    Let $D \in \Diffop^2(E)$ be normally hyperbolic.
    \begin{remarklist}
    \item \label{item:lowest-hadamard-equation} According to our above
        computation, the transport equation for $k=0$ implies
        \begin{equation}
            \label{eq:lowest-hadamard-equationInTheRemark}
            2 \nabla^E_{\gradient R^\pm_U(2,p)} V_p^0
            + V_p^0 \dAlembert R^\pm_U(2,p)
            = \delta_p.
        \end{equation}
    \item \label{item:transp-eq-same-for-adv-ret-green} The transport
        equations are the same for the advanced and retarded
        $\mathcal{R}^\pm(p)$. Thus we only have to solve them once and
        can us the \emph{same} coefficients $V_p^k$ for both Green
        functions.
    \end{remarklist}
\end{remark}
\begin{definition}[Hadamard coefficients]
    \label{definition:hadamard-coefficients}
    \index{Hadamard coefficients}%
    \index{Normally hyperbolic operator}%
    Let $D \in \Diffop^2(E)$ be normally hyperbolic and $U \subseteq
    M$ geodesically star-shaped around $p \in M$ as before. Solutions
    $V_p^k \in \Secinfty(E\at{U}) \otimes E^*_p$ of the transport
    equations are then called Hadamard coefficients for $D$ at the
    point $p$.
\end{definition}

In the following we shall now explicitly construct the Hadamard
coefficients and show their uniqueness. Note however, that even having
the $V_p^k$ does not yet solve the problem of finding a Green function
since the convergence of \eqref{eq:Ansatz-for-Green-function} is still
delicate.

%
%

\subsection{Uniqueness of the Hadamard Coefficients}
\label{subsec:UniquenessHadamardCoefficients}

We shall now prove that the Hadamard coefficients are necessarily
unique. To this end we need the parallel transport in $E$ with respect
to the covariant derivative $\nabla^E$ induced by $D$. Since on $U$ we
have unique geodesics joining $p$ with any other point $q \in U$,
namely
\begin{equation}
    \label{eq:geodesics-in-U}
    \gamma_{p \rightarrow q}(t) = \exp_p (t \exp_p^{-1}(q)),
\end{equation}
we shall always use these paths for parallel transport. For
abbreviation, we set
\begin{equation}
    \label{eq:parallel-transport-in-U}
    P_{p \rightarrow q} 
    = P_{\gamma_{p \rightarrow q}, 0 \rightarrow 1}:
    E_p \longrightarrow E_q.
\end{equation}
From the explicit definition of the parallel transport we find the
following technical statement:
\begin{lemma}
    \label{lemma:partrans-as-section}
    \index{Parallel transport}%
    The parallel transport along geodesics in $U$ yields a smooth map
    \begin{equation}
        \label{eq:partrans-as-section}
        U \ni q \mapsto P_{p \rightarrow q} \in E_q \tensor E^*_p,
    \end{equation}
    which we can view as a smooth section
    \begin{equation}
        \label{eq:partrans-as-section-2}
        P_{p \rightarrow \argument} \in \Secinfty(E\at{U}) \tensor E^*_p.
    \end{equation}
\end{lemma}
\begin{proof}
    Let $e_\alpha \in \Secinfty(E\at{U})$ be a locally defined smooth
    frame and let $A_\alpha^\beta$ be the corresponding smooth
    connection one-forms. Then the parallel transport is determined by
    the equation
    \[
    \dot{s}^\beta(t) + A_\alpha^\beta
    \left(
        \frac{\D}{\D t} \exp_p(t \exp_p^{-1}(q))
    \right)
    s^\alpha(t)
    = 0.
    \tag{$*$}
    \]
    Since the map $(q,t) \mapsto \exp_p ( t \exp_p^{-1}(q))$ is smooth
    on an open neighborhood of $U \times [0,1] \subseteq U \times
    \mathbb{R}$ the solutions to ($*$) also depend smoothly on $q$ and
    $t$ on this neighborhood. Thus, the solutions depend smoothly on
    $q$ when evaluated at $t=1$, which implies the smoothness of
    \eqref{eq:partrans-as-section-2}.
\end{proof}

Using this smoothness of the parallel transport we can obtain the
following result:
\begin{theorem}[Uniqueness of the Hadamard coefficients]
    \label{theorem:uniqueness-of-hadamard}
    \index{Hadamard coefficients!uniqueness}%
    Let $U \subseteq M$ be geodesically star-shaped around $p$ and let
    $D \in \Diffop^2(E)$ be normally hyperbolic. Then the Hadamard
    coefficients for $D$ at $p$ are necessarily unique. In fact, they
    satisfy
    \begin{equation}
        \label{eq:lowest-hadamard-is-partrans}
        V_p^0 = \frac{1}{\sqrt{\varrho_p}} P_{p \rightarrow \argument}
    \end{equation}
    and for $k \geq 1$ and $q \in U$
    \begin{equation}
        \label{eq:higher-hadamard-recursive-formula}
        V_p^k(q)
        = -k \frac{1}{\sqrt{\varrho_p(q)}} P_{p \rightarrow q}
        \left(
            \int_0^1 \sqrt{\varrho_p} (\gamma_{p \rightarrow q}(\tau))
            \tau^{k-1}
            P_{\gamma_{p \rightarrow q}, 0 \rightarrow \tau}^{-1}
            \left(
                D(V_p^{k-1})(\gamma_{p \rightarrow q})(\tau)
            \right)
        \right) \D \tau.
    \end{equation}
\end{theorem}
\begin{proof}
    We consider the ``\Index{Lorentz radius}'' function $r_p =
    \sqrt{|\eta_p|} \in \Fun[0](U)$ which is continuous but not
    differentiable. However, on $U \setminus C_U(p)$ where $C_U(p) =
    C^+_U(p) \cup C^-_U(p)$ with
    \[
    C^\pm_U(p) = \exp_p (C^\pm(0) \cap V),
    \]
    the function $\eta_p$ is non-zero and hence $r_p \in \Cinfty(U
    \setminus C_U(p))$ is smooth. On $U \setminus C_U(p)$ we have
    \[
    \eta_p = \epsilon r_p^2
    \]
    with $\epsilon(q) = +1$ for $\exp_p^{-1}(q)$ timelike and
    $\epsilon(q) = -1$ for $\exp_p^{-1}(q)$ spacelike, respectively.
    Using our results from
    Proposition~\ref{proposition:distance-function-on-M} we find
    \begin{align*}
        \frac{1}{2} \dAlembert \eta_p - n
        &= \frac{1}{2} g(\gradient \log \varrho_p, \gradient \eta_p)
        = \frac{1}{2} \Lie_{\gradient \eta_p} \log \varrho_p
        = \Lie_{\gradient \eta_p} \log \sqrt{\varrho_p},
    \end{align*}
    valid on $U$ since $\varrho_p > 0$. Moreover, by
    \eqref{eq:length-of-gradient-of-dist-func} we get on $U \setminus
    C_U(p)$
    \begin{align*}
        \Lie_{\gradient \eta_p} (\log r_p^k)
        &= k \Lie_{\gradient \eta_p} (\log r_p)
        = k \frac{1}{r_p} \Lie_{\gradient \eta_p} (r_p)
        = k \frac{1}{\sqrt{\epsilon \eta_p}}
        \Lie_{\gradient \eta_p} \left(\sqrt{\epsilon \eta_p}\right) \\
        &= k \frac{1}{\sqrt{\epsilon \eta_p}}
        \frac{\epsilon}{2 \sqrt{\epsilon \eta_p}}
        \Lie_{\gradient \eta_p} \eta_p
        = \frac{k}{2 \eta_p} \SP{\gradient \eta_p, \gradient \eta_p}
        = 2 k.
    \end{align*}
    Since $\sqrt{\varrho_p} r_p^k > 0$ on $U \setminus C_U(p)$ we can
    rewrite the transport equation \eqref{eq:transport-equations}
    equivalently as
    \begin{align*}
        -2k D(V_p^{k-1})
        = \nabla^E_{\gradient \eta_p} V_p^k
        + \left(
            \frac{1}{2} \dAlembert \eta_p - n +2k
        \right) V_p^k
        = \nabla^E_{\gradient \eta_p} V_p^k
        + \frac{1}{\sqrt{\varrho_p} r_p^k}
        \Lie_{\gradient \eta_p} \left(\sqrt{\varrho_p} r_p^k\right) V_p^k,
    \end{align*}
    and thus as
    \[
    \nabla^E_{\gradient \eta_p} \left(\sqrt{\varrho_p} r_p^k V_p^k\right)
    = - \sqrt{\varrho_p} r_p^k 2k D(V_p^{k-1}).
    \tag{$*$}
    \]
    Since on $U \setminus C_U(p)$ the additional factor
    $\sqrt{\varrho_p} r_p^k$ is both positive and smooth, ($*$) is
    equivalent to the transport equation on $U \setminus C_U(p)$.

    Now we consider first $k=0$. Then ($*$) means that
    \[
    \nabla^E_{\gradient \eta_p} \left(\sqrt{\varrho_p} V_p^0\right)
    \At{U \setminus C_U(p)}
    = 0.
    \tag{$**$}
    \]
    Since the gradient $\gradient \eta_p$ is at every point $q$ just
    twice the tangent vector of the geodesic $\gamma_{p \rightarrow
      q}(t) = \exp_p(t \exp_p^{-1}(q))$ we conclude from ($**$) that
    the local section $\sqrt{\varrho_p} V_p^0 \in \Secinfty(E\at{U})
    \tensor E_p^*$ is covariantly constant in direction of all
    geodesics $\gamma_{p \rightarrow q}$ as long as $q \in U \setminus
    C_U(p)$, i.e. as long as $\exp_p^{-1}(q)$ is either timelike or
    spacelike. But $\sqrt{\varrho_p} V_p^0$ is smooth and thus by
    continuity we conclude that ($**$) holds on all of $U$. But this
    shows that $\sqrt{\varrho_p} V_p^0$ is parallel along all
    geodesics starting at $p$ whence it is given by means of the
    parallel transport, i.e.
    \[
    \sqrt{\varrho_p} V_p^0 \at{q}
    = P_{p \rightarrow q}
    \left(\sqrt{\varrho_p} V_p^0 \at{p} \right)
    = P_{p \rightarrow q} \left( 1 \cdot \id_{E_p} \right)
    = P_{p \rightarrow q},
    \]
    since by assumption $V_p^0(p) = \id_{E_p}$ and $\varrho_p(p) = 1$
    by
    Proposition~\ref{proposition:taylor-expansion-of-dens-compare-function}.
    Indeed, if $e_\alpha \in E_p$ is a basis then $P_{p \rightarrow q}
    (\id_{E_p}) = P_{p \rightarrow q} (e_\alpha \tensor e^\alpha) =
    P_{p \rightarrow q}(e_\alpha) \tensor e^\alpha = P_{p \rightarrow
      q}$ since the parallel transport only acts on the $E_p$-part of
    $\id_{E_p}$ and not on the $E_p^*$-part which is considered as
    values in all of our considerations up to now. But this shows
    \eqref{eq:lowest-hadamard-is-partrans} and hence the uniqueness of
    $V_p^0$.

    Now let $k \geq 1$. Then we again consider ($*$) on $U \setminus
    C_U(p)$. To this end we first note that since $\gradient \eta_p$
    is twice the push-forward of the Euler vector field $\xi_{T_pM}$
    on $T_pM$ its flow can be computed explicitly. In fact, let $c(t)
    = \exp_p(\E^{2t} \exp_p^{-1}(q))$ then for small $t$ around $0$ we
    have by Proposition~\ref{proposition:distance-function-on-M}
    \begin{align*}
        \dot{c}(t)
        = T_{\exp_p(\E^{2t} \exp_p^{-1}(q))} \exp_p 
        \left(
            2 \E^{2t} \exp_p^{-1}(q)
        \right)
        = 2 T_{c(t)} \exp_p (\exp_p^{-1}(c(t)))
        = \gradient \eta_p \at{c(t)},
    \end{align*}
    whence $c(t)$ is the integral curve of $\gradient \eta_p$ through
    $c(0)=q$. Thus ($*$) implies
    \[
    \nabla^\#_{\frac{\partial}{\partial t}}
    \left(c^\# \sqrt{\varrho_p} r_p^k V_p^k\right)
    = -2k c^\# \left(\sqrt{\varrho_p} r_p^k D(V_p^{k-1})\right),
    \tag{$*$$**$}
    \]
    where $\nabla^\#$ is the pull-back connection with respect to the
    curve $c$. Thus $\sqrt{\varrho_p} r_p^k V_p^k$ satisfies the
    perturbed parallel transport equation along $c$ with
    perturbation given by the right hand side of ($*$$**$). The
    solutions of such equations are obtained in terms of the parallel
    transport as follows:
    \begin{lemma}
        \label{lemma:variation-of-constants}
        \index{Parallel transport!perturbation}%
        Let $\gamma: I \subseteq \mathbb{R} \longrightarrow M$ be a
        smooth curve on an open interval and let $f \in
        \Secinfty(\gamma^\# E)$ be a smooth section. Then the
        perturbed parallel transport equation
        \begin{equation}
            \label{eq:perturbated-partrans-equation}
            \nabla^\#_{\frac{\partial}{\partial t}} s = f
        \end{equation}
        has
        \begin{equation}
            \label{eq:perturbated-partrans-solution}
            s(t)
            = P_{\gamma, t_0 \rightarrow t}
            \left(
                s(t_0) + \int_{t_0}^t
                P_{\gamma, t_0 \rightarrow \tau}^{-1} (f(\tau))
                \D \tau
            \right)
        \end{equation}
        as unique and smooth solution $s \in \Secinfty(\gamma^\# E)$
        with initial condition $s(t_0) \in E_{\gamma(t_0)}$ for $t_0
        \in I$.
    \end{lemma}
    \begin{subproof}
        We choose a frame $e_\alpha(t_0) \in E_{\gamma(t_0)}$ at
        $\gamma(t_0)$ and parallel transport it to $e_\alpha(t) =
        P_{\gamma, t_0 \mapsto t} (e_\alpha(t_0))$. This yields a
        covariantly constant frame $e_\alpha \in \Secinfty(\gamma^\#
        E)$, i.e. we have $\nabla^\#_{\frac{\partial}{\partial t}}
        e_\alpha = 0$, see also the proof of
        Lemma~\ref{lemma:time-derivative-parallel-transport}. Then
        $s(t) = s^\alpha(t) e_\alpha(t)$ for any section $s \in
        \Secinfty(\gamma^\# E)$ with $s^\alpha \in \Cinfty(I)$. Thus
        \eqref{eq:perturbated-partrans-equation} becomes
        \[
        \dot{s}^\alpha(t) = f^\alpha(t)
        \]
        for all $\alpha$ with initial conditions $s^\alpha(t_0)$ for
        $t=t_0$. The unique solution to this system of ordinary first
        order differential equations is
        \[
        s^\alpha(t) 
        = s^\alpha(t_0) + \int_{t_0}^t f^\alpha(\tau) \D \tau.
        \tag{$\diamond$}
        \]
        Now we compute
        \begin{align*}
            P_{\gamma, t_0 \rightarrow t} \circ
            P_{\gamma, t_0 \rightarrow \tau}^{-1}
            (f^\alpha(\tau) e_\alpha(\tau))
            = f^\alpha(\tau) P_{\gamma, t_0 \rightarrow t}
            (e_\alpha(t_0))
            = f^\alpha(\tau) e_\alpha(t),
        \end{align*}
        whence
        \begin{align*}
            P_{\gamma, t_0 \rightarrow t}
            \left(
                s(t_0) + \int_{t_0}^t
                P_{\gamma, t_0 \rightarrow \tau}^{-1} (f(\tau))
                \D \tau
            \right)
            &= P_{\gamma, t_0 \rightarrow t}
            (s^\alpha(t_0) e_\alpha(t_0))
            + \int_{t_0}^t f^\alpha(\tau) e_\alpha(t) \D \tau \\
            &= \left(
                s^\alpha(t_0) + \int_{t_0}^t f^\alpha(\tau) \D \tau
            \right) e_\alpha(t) \\
            &= s^\alpha(t) e_\alpha(t) \\
            &= s(t)
        \end{align*}
        as wanted. The uniqueness is clear from specifying the initial
        conditions and smoothness follows from the smoothness of the
        $f^\alpha$ and the explicit form ($\diamond$).
    \end{subproof}

    \medskip
    \noindent
    We apply the lemma to the curve $c(t) = \exp_p \left(\E^{2t}
        \exp_p^{-1}(q)\right)$ where $v = \exp_p^{-1}(q)$ is either
    spacelike or timelike and $t \in (-\infty, \epsilon)$ with some
    small $\epsilon > 0$ such that $\E^{2t} v$ is still in the domain
    $V \subseteq T_pM$ of $\exp_p$. Then the \emph{homogeneous}
    transport equation for $k \geq 1$ implies
    \[
    \nabla^\#_{\frac{\partial}{\partial t}}
    c^\# \left(\sqrt{\varrho_p} r_p^k V_p^k\right)
    = 0,
    \tag{\frownie}
    \]
    and hence
    \[
    \sqrt{\varrho_p} r_p^k V_p^k \at{c(t)}
    = P_{c, t_0 \rightarrow t}
    \left(\sqrt{\varrho_p} r_p^k V_p^k \at{c(t_0)}\right).
    \]
    Taking e.g. $t_0 = 0$ we obtain for all $t \in (-\infty,
    \epsilon)$
    \[
    V_p^k (\exp_p (\E^{2t} v))
    =
    \frac{1}{
      \left(\sqrt{\varrho_p} r_p^k\right)
      \left(\exp_p(\E^{2t} v)\right)
    }
    P_{c, 0 \rightarrow t} 
    \left(\sqrt{\varrho_p} r_p^k V_p^k \at{q}\right).
    \]
    Suppose $V_p^k(q) \neq 0$ then also $\sqrt{\varrho_p} r_p^k V_p^k
    \at{q} \neq 0$. Since the parallel transport is reparametrization
    invariant we can write this equally well as
    \[
    V_p^k (\exp_p(tv))
    = \frac{1}
    {
      \left(\sqrt{\varrho_p} r_p^k\right) 
      \left(\exp_p(\E^{2t} v)\right)
    } 
    P_{\gamma, t_0 \rightarrow t}
    \left(\sqrt{\varrho_p} r_p^k V_p^k \at{q}\right),
    \tag{\smiley}
    \]
    with $\gamma(t) = \exp_p(tv)$ being the geodesic reparametrization
    by the ``arc length''. Now the limit $t \longrightarrow 0$ of
    $P_{\gamma, t_0 \rightarrow t}$ exists and is given by $P_{p
      \rightarrow q}^{-1}$. Thus the limit of $t \longrightarrow 0$ of
    $P_{\gamma, t_0 \rightarrow t} \left(\sqrt{\varrho_p} r_p^k V_p^k
        \at{q}\right)$ exists and is a certain non-zero vector. But
    for $k \geq 1$ the limit of $r_p^k(\exp_p(tv)$) for $t
    \longrightarrow 0$ is $0$ whence the prefactor in (\smiley)
    becomes singular. Thus $V_p^k$ can not be continuous at $p$. Thus
    $V_p^k(q) \neq 0$ for \emph{some} $q \in U \setminus C_U(p)$
    already implies that $V_p^k$ is non-continuous at $p$. We conclude
    that the homogeneous equation (\frownie) has only the trivial
    solution $V_p^k = 0$ as everywhere \emph{smooth} solution.  This
    implies that the inhomogeneous transport equation ($**$$*$) can
    have at most \emph{one} everywhere smooth solution which shows
    uniqueness of the $V_p^k$ for $k \geq 1$. It remains to shows that
    they necessarily satisfy
    Equation~\eqref{eq:higher-hadamard-recursive-formula}.  According
    to Lemma~\ref{lemma:variation-of-constants}, a particular solution
    along the curve $c$ is given by
    \[
    \left(\sqrt{\varrho_p} r_p^k V_p^k)(c(t)\right)
    = P_{c, t_0 \rightarrow t}
    \left(
        \left(\sqrt{\varrho_p} r_p^k V_p^k\right)(c(t_0))
        + \int_{t_0}^t
        \left( P_{\gamma,t_0 \rightarrow \tau} \right)^{-1}
        \left(
            -2k \sqrt{\varrho_p} r_p^k D(V_p^{k-1})
        \right) (c_\tau)
        \D \tau
    \right),
    \]
    where $t, t_0 \in (-\infty, \epsilon)$ with some suitable
    $\epsilon > 0$. Since we are interested in a solution $V_p^k$
    which is still defined at $q = p$, the limit of $V_p^k (c(t_0))$
    for $t_0 \longrightarrow -\infty$ should exist. But then the limit
    $t_0 \longrightarrow -\infty$ yields $(\sqrt{\varrho_0}
    r_p^k)(c(t_0)) \longrightarrow (\sqrt{\varrho_0} r_p^k)(p) = 0$
    whence the first term on the right hand side does not contribute
    in this limit. Thus under the regularity assumption we have
    \begin{align*}
        (\sqrt{\varrho_p} r_p^k V_p^k)(c(t))
        &= P_{c, -\infty \rightarrow t}
        \left(
            \int_{-\infty}^t
            \left( P_{\gamma, -\infty \rightarrow \tau} \right)^{-1}
            \left(
                -2k \sqrt{\varrho_p} r_p^k D(V_p^{k-1})
            \right) (c_\tau)
            \D \tau
        \right) \\
        &= -2k P_{c, -\infty \rightarrow t}
        \int_{-\infty}^t \sqrt{\varrho_p}(c(\tau)) r_p^k(c(\tau))
        P_{c, -\infty \rightarrow \tau}^{-1} \left(D(V_p^{k-1})(c(\tau))\right)
        \D \tau.
    \end{align*}
    Now we have by the isometry properties of the exponential map in
    radial direction
    \begin{align*}
        r_p(c(t)) 
        &=
        \sqrt{|\eta_p|} \exp_p \left(\E^{2t} \exp_p^{-1}(q)\right) \\
        &=
        \sqrt{
          \left|
              g_p\left(
                  \E^{2t} \exp_p^{-1}(q), \E^{2t} \exp_p^{-1}(q)
              \right)
          \right|
        } \\
        &=
        \E^{2t}
        \sqrt{
          \left|
              g_p
              \left(\exp_p^{-1}(q), \exp_p^{-1}(q)\right)
          \right|
        },
    \end{align*}
    and by assumption $g_p\left(\exp_p^{-1}(q), \exp_p^{-1}(q)\right)
    \neq 0$.  After dividing by $\sqrt{\left|g_p\left(\exp_p^{-1}(q),
              \exp_p^{-1}(q)\right)\right|}$ we obtain
    \begin{align*}
        \E^{2kt} \sqrt{\varrho_p}(c(t)) V_p^k(c(t))
        & = -2k P_{c, -\infty \rightarrow t}
        \int_{-\infty}^t \sqrt{\varrho_p}(c(\tau)) \E^{2k\tau}
        P_{c, -\infty \rightarrow \tau}^{-1} 
        \left( D(V_p^{k-1})(c(\tau)) \right)
        \D \tau \\
        & = -2k P_{\alpha, 0 \rightarrow \E^{2t}}
        \int_0^{\E^{2t}} \sqrt{\varrho_p}(c(\sigma)) \sigma^k
        P_{\gamma, 0 \rightarrow \sigma}^{-1}
        \left( D(V_p^{k-1})(\gamma(\sigma)) \right)
        \frac{\D \sigma}{2 \sigma},
    \end{align*}
    with the substitution $\sigma = \E^{2\tau}$ and thus $\D \tau =
    \frac{\D \sigma}{2 \sigma}$. Note that with this substitution
    $c(\tau(\sigma)) = \exp_p \left(\E^{2\tau(\sigma)}
        \exp_p^{-1}(q)\right) = \exp_p (\sigma \exp_p^{-1}(q)) =
    \gamma(\sigma)$ is indeed the geodesic from $p$ to $q$. By the
    invariance under reparametrization of the parallel transport we
    get $P_{c, -\infty \rightarrow \tau} = P_{\gamma, 0 \rightarrow
      \sigma}$ which explains the above formula. Taking $t=0$ we find
    $c(0) = \gamma(1) = q$ and thus
    \[
    \sqrt{\varrho_p(q)} V_p^k(q)
    = -k P_{p \rightarrow q}
    \int_0^1 \sqrt{\varrho_p}(\gamma_{p \rightarrow q}(\tau))
    \tau^{k-1}
    P_{\gamma_{p \rightarrow q}, 0 \rightarrow \tau}^{-1}
    \left(
        D(V_p^{k-1}) (\exp_p(\tau \exp_p^{-1}(q)))
    \right)
    \D \tau
    \]
    after replacing $\sigma$ by $\tau$ again. Since $\sqrt{\varrho_p}
    > 0$ this gives
    \eqref{eq:higher-hadamard-recursive-formula}. Indeed,
    \eqref{eq:higher-hadamard-recursive-formula} only follows for $q
    \in U \setminus C_U(p)$ but the continuity of the right hand side
    makes \eqref{eq:higher-hadamard-recursive-formula} correct
    everywhere.
\end{proof}

\begin{remark}
    \label{remark:uniqueness-of-hadamards}
    Note that the additional $r_p^k$ in the higher transport equations
    yields a completely different behaviour of the solution for $q
    \longrightarrow p$. While for $k=0$ no singularities arise the
    case $k \geq 1$ behaves much more singular. In fact, only one
    solution is everywhere smooth. This is the reason why for $k=0$ we
    have to specify an initial condition $V_p^0(p) = \id_{E_p}$ while
    for $k \geq 1$ the boundary condition of being smooth at $q = p$
    fixes the solution.
\end{remark}

%
%

\subsection{Construction of the Hadamard Coefficients}
\label{subsec:ConstructionHadamardCoefficients}

In Theorem~\ref{theorem:uniqueness-of-hadamard} we have not only shown
the uniqueness of the Hadamard coefficients which was essentially a
consequence of the desired smoothness at $p$ but we also obtained a
rather explicit recursive formula for the $V_p^k$. Using
\eqref{eq:lowest-hadamard-is-partrans} and
\eqref{eq:higher-hadamard-recursive-formula} we recursively
\emph{define} $V_p^k$ for $k \geq 0$ by
\begin{equation}
    \label{eq:lowest-hadamard}
    V_p^0(q) 
    = \frac{1}{\sqrt{\varrho_p(q)}} P_{p \rightarrow q}      
\end{equation}
and
\begin{equation}
    \label{eq:higher-hadamards}
    V_p^k(q)
    = -\frac{k}{\sqrt{\varrho_p(q)}} P_{p \rightarrow q}
    \int_0^1 \sqrt{\varrho_p}(\gamma_{p \rightarrow q}(\tau))
    \tau^{k-1}
    P_{\gamma_{p \rightarrow q}, 0 \rightarrow \tau}^{-1}
    \left(
        D(V_p^{k-1})\exp_p(\tau \exp_p^{-1}(q))
    \right)
    \D \tau
\end{equation}
for $q \in U$. Thus it remains to show that these $V_p^k$ indeed
define smooth sections satisfying the transport equations. The
smoothness is guaranteed from the following proposition which even
handles the smooth dependence on $p$. We again formulate it for a
situation as in
Proposition~\ref{proposition:riesz-dependence-on-base-point-p-prime}.
\begin{proposition}[Smoothness of $V^k$]
    \label{proposition:smoothness-of-hadamard}
    \index{Hadamard coefficients!smoothness}%
    Let $O \subseteq U \subseteq M$ be open subsets such that $U$ is
    geodesically star-shaped around all $p \in O$. Then the recursive
    definitions \eqref{eq:lowest-hadamard} and
    \eqref{eq:higher-hadamards} yield smooth sections
    \begin{equation}
        \label{eq:hadamard-as-smooth-section}
        V^k \in \Secinfty(E^* \extensor E \at{O \times U})
    \end{equation}
    via the definition
    \begin{equation}
        \label{eq:hadamard-dependence-on-p-and-q}
        V^k(p,q) = V_p^k(q)
    \end{equation}
    for $(p,q) \in O \times U$ and $k \geq 0$.
\end{proposition}
\begin{proof}
    First we note that $\varrho(p,g) = \varrho_p(q)$ is actually a
    smooth function $\varrho \in \Cinfty(O \times U)$ with $\varrho >
    0$ everywhere. This follows from
    Lemma~\ref{lemma:density-comparison-function}. From
    Lemma~\ref{lemma:partrans-as-section} we deduce that the
    dependence of $P_{p \rightarrow q}$ on $q$ is smooth and a similar
    argument shows that also the dependence on $p$ is smooth. In fact,
    the parallel transport depends smoothly on $(p,q) \in O \times U$
    yielding thereby a smooth section
    \[
    P \in \Secinfty\left(E^* \extensor E \at{O \times U}\right).
    \]
    It follows that $V^0$ is smooth on $O \times U$. We rewrite the
    recursive definition \eqref{eq:higher-hadamards} in terms of
    $\varrho$,
    \[
    \gamma(p,q,\tau) = \gamma_{p \rightarrow q}(\tau)
    = \gamma_\tau(p,q)
    \]
    and the $V^k$. Then \eqref{eq:higher-hadamards} becomes
    \[
    V^k
    = \frac{k}{\sqrt{\varrho}} P \int_0^1 \sqrt{\varrho} \circ
    \gamma_\tau P \circ \gamma_\tau 
    \left(
        (\id \extensor D(V^{k-1})) \circ \gamma_\tau
    \right)
    \tau^{k-1} \D \tau.
    \]
    By induction we assume that $V^{k-1}$ is smooth. Now $\gamma$ is
    smooth on $O \times U \times [0,1]$ and thus the integrand is
    smooth with a compact domain of integration. This results in a
    smooth $V^k$.
\end{proof}

As already in
Proposition~\ref{proposition:riesz-dependence-on-base-point-p-prime}
we can e.g. take a \emph{convex} $U \subseteq M$ and set $O = U$ in
order to meet the conditions of
Proposition~\ref{proposition:smoothness-of-hadamard}. It remains to
show that the $V_p^k$ actually satisfy the transport equations with
the correct initial condition.
\begin{proposition}
    \label{proposition:hadamard-solve-transport}
    \index{Transport equation}%
    Let $U \subseteq M$ be geodesically star-shaped around $p \in
    M$. Then the sections $V_p^k \in \Secinfty(E \at{U}) \tensor
    E_p^*$ defined by \eqref{eq:lowest-hadamard} and
    \eqref{eq:higher-hadamards} satisfy the transport equations
    \eqref{eq:transport-equations} with initial condition $V_p^0(p) =
    \id_{E_p}$.
\end{proposition}
\begin{proof}
    Clearly $V_p^0(p) = \id_{E_p}$ since $\varrho_p(p) = 1$. In the
    proof of Theorem~\ref{theorem:uniqueness-of-hadamard} we have seen
    that \eqref{eq:transport-equations} is equivalent to
    \[
    \nabla^E_{\gradient \eta_p} 
    \left(\sqrt{\varrho_p} r_p^k V_p^k\right)
    = 
    -2k \sqrt{\varrho_p} r_p^k D(V_p^{k-1})
    \tag{$*$}
    \]
    on the open subset $U \setminus C_U(p)$. Since we already know
    that the section $V^k$ are smooth on $U$ by
    Proposition~\ref{proposition:smoothness-of-hadamard} we know that
    they satisfy \eqref{eq:transport-equations} on $U$ iff they
    satisfy \eqref{eq:transport-equations} on $U \setminus C_U(p)$ by
    a continuity argument. Thus it suffices to show ($*$) on $U
    \setminus C_U(p)$. In the proof of
    Theorem~\ref{theorem:uniqueness-of-hadamard} we have shown that
    ($*$) implies
    \[
    \nabla^\#_{\frac{\partial}{\partial t}}
    \left(
        c^\# \left(\sqrt{\varrho_p} r_p^k V_p^k\right)
    \right)
    = -2k c^\#
    \left(
        \sqrt{\varrho_p} v_p^k D(V_p^{k-1})
    \right)
    \tag{$**$}
    \]
    for the curve $c(t) = \exp_p (\E^{2t} \exp_p^{-1}(q))$ with $q \in
    U$ and $t \in (-\infty, \epsilon)$ and $\epsilon > 0$ sufficiently
    small. But if we have ($**$) for \emph{all} such curves $c$ then
    we get back ($*$) since $\gradient \eta_p \at{q} = \dot{c}(0)$ and
    the left hand side of ($*$) can be evaluated point by point as
    $\nabla^E_{\gradient \eta_p}$ is tensorial in $\gradient
    \eta_p$. Thus ($**$) for \emph{all} such curves is equivalent to
    ($*$). But $V_p^k(q)$ was precisely the solution of ($**$) at
    $t=0$ by Lemma~\ref{lemma:variation-of-constants}. But this means
    at $q$ we have
    \[
    \nabla^E_{\gradient \eta_p} 
    \left(\sqrt{\varrho_p} r_p^k V_p^k\right) \At{q}
    = -2k \sqrt{\varrho_p} r_p^k D(V_p^{k-1}) \at{q}.
    \]
    Since $q \in U \setminus C_U(p)$ was arbitrary, ($*$) follows
    which completes the claim.
\end{proof}

\begin{theorem}[Hadamard Coefficients]
    \label{theorem:hadamard-coefficients}
    \index{Hadamard coefficients!existence}%
    Let $O \subseteq U \subseteq M$ be open subsets such that $U$ is
    geodesically star-shaped around all $p \in O$. Let $D \in
    \Diffop^2(E)$ be normally hyperbolic. Then for each $p \in O$ the
    operator $D$ has unique Hadamard coefficients $V_p^k \in
    \Secinfty\left(E\at{U}\right) \tensor E_p^*$ explicitly given by
    $V_p^k(q) = V^k(p,q)$ where $V^k \in \Secinfty\left(E^* \extensor
        E \at{O \times U}\right)$ is recursively determined by
    \begin{equation}
        \label{eq:lowest-hadamard-no-p}
        V^0 = \frac{1}{\sqrt{\varrho}} P
    \end{equation}
    and
    \begin{equation}
        \label{eq:higher-hadamard-no-p}
        V^k = - \frac{k}{\sqrt{\varrho}} P \int_0^1
        \left(
            \sqrt{\varrho} P \left(\id \extensor D (V^{k-1})\right)
        \right) \circ \gamma_\tau \: \tau^{k-1}
        \D \tau,
    \end{equation}
    where $P \in \Secinfty\left(E^* \extensor E \at{O \times
          U}\right)$ is the parallel transport $P(p,q) = P_{p
      \rightarrow q}$ along $\gamma_{p \rightarrow q}(\tau) =
    \gamma_\tau(p,q) = \exp_p(\tau \exp_p^{-1}(q))$. On the diagonal
    we explicitly have the simplified recursion
    \begin{equation}
        \label{eq:higher-hadamard-on-diagonal}
        V^k(p,p) = - \left((\id \extensor D) (V^{k-1})\right) (p,p).
    \end{equation}
\end{theorem}
\begin{proof}
    It remains to show the simplified recursion
    \eqref{eq:higher-hadamard-on-diagonal} for $p = q$. For $k \geq 1$
    we have
    \begin{align*}
        V^k(p,p)
        &= - \frac{k}{\sqrt{\varrho}(p,p)} P_{p \rightarrow p} \int_0^1
        \left(
            \sqrt{\varrho} P \left((\id \extensor D) (V^{k-1})\right)
        \right)(\gamma_\tau(p,p)) \tau^{k-1}
        \D \tau \\
        &= -k \int_0^1 (\id \extensor D) (V^{k-1}) \at{(p,p)}
        \tau^{k-1} \D \tau \\
        &= - k (\id \extensor D) (V^{k-1}) \at{(p,p)}
        \int_0^1 \tau^{k-1} \D \tau \\
        &= - (\id \extensor D) (V^{k-1}) \at{(p,p)}.
    \end{align*}
\end{proof}

We illustrate the recursion formula by computing the first non-trivial
Hadamard coefficient along the diagonal.
\begin{example}[First Hadamard coefficient]
    \label{example:first-hadamards}
    \index{Hadamard coefficient!first}%
    Let $D = \dAlembert^\nabla + B$ be normally hyperbolic as
    usual. Thus let $s_p \in E_p$ be a vector in $E_p$ and let
    \begin{equation}
        \label{eq:partransed-vector-field}
        s(q) = P_{p \rightarrow q} (s_p),
    \end{equation}
    which defines a vector field $s \in
    \Secinfty\left(E\at{U}\right)$. We compute the covariant
    derivatives of $s$ at $p$. At general points $q \in U$ this might
    be very complicated but at $p$ we have by
    Proposition~\ref{proposition:taylor-coefficients-of-partrans} the
    formal Taylor expansion
    \begin{equation}
        \label{eq:derivatives-of-partransed-vectorfield}
        \ins(e_{i_1}) \cdots \ins(e_{i_k}) \frac{1}{k!}
        \left(\SymD^E\right)^k s
        \At{p}
        =
        \frac{\partial^k}{\partial v^{i_1} \cdots \partial v^{i_k}}
        \left( P_{\gamma_v, 0 \rightarrow 1} \right)^{-1}
        (s(\gamma_v(1))),
    \end{equation}
    with a basis $e_1, \ldots, e_n \in T_pM$ and $\gamma_v(t) =
    \exp_p(tv)$ as usual. But
    \[
    \left( P_{\gamma, 0 \rightarrow 1} \right)^{-1}
    (s(\gamma_v(1)))
    = \left( P_{\gamma_v, 0 \rightarrow 1} \right)^{-1}
    P_{p \rightarrow q = \exp_p(v)} (s_p)
    = s_p
    \]
    is independent of $v$. Thus all partial derivatives vanish and we
    conclude $\left(\SymD^E\right)^k s \at{p} = 0$.  But then
    $\dAlembert^\nabla s \at{p} = \frac{1}{2} \SP{g^{-1},
      \left(\SymD^E\right)^2 s} \at{p} = 0$ follows as well. From this
    we conclude by \eqref{eq:higher-hadamard-on-diagonal}
    \begin{align*}
        &V^1(p,p) \\
        &\quad= - (\id \extensor D) (V^0) (p,p)
        = - (\id \extensor D) \left(
            \frac{1}{\sqrt{\varrho}} P
        \right)(p,p)
        = - D \left(
            \frac{1}{\sqrt{\varrho_p}} P_{p \rightarrow \argument}
            (e_\alpha)
        \right) \tensor e^\alpha \At{p} \\
        &\quad=
        - \frac{1}{\sqrt{\varrho_p}} 
        D \left( P_{p \rightarrow \argument}(e_\alpha) \right)
        \tensor e^\alpha \At{p} 
        -
        2 \left(
            \nabla^E_{\gradient \frac{1}{\sqrt{\varrho_p}}}
            P_{p \rightarrow \argument}(e_\alpha)
        \right) \tensor e^\alpha \At{p} 
        -
        \dAlembert \frac{1}{\sqrt{\varrho_p}} \At{p}
        P_{p \rightarrow \argument}(e_\alpha) \tensor e^\alpha \At{p} \\
        &\quad=
        - (\dAlembert^\nabla + B)
        \left(
            P_{p \rightarrow \argument}(e_\alpha)
        \right) \At{p} \tensor e^\alpha
        +
        0
        -
        \dAlembert \frac{1}{\sqrt{\varrho_p}} \At{p}
        e_\alpha \tensor e^\alpha \\
        &\quad=
        - B\At{p} (e_\alpha) \tensor e^\alpha
        - \frac{1}{6} \scal(p) \id_{E_p},
    \end{align*}
    by \eqref{eq:dAlembert-of-sqrt-of-dens-comp-funcII} and $\id_{E_p}
    = e_\alpha \tensor e^\alpha$ with a basis $e_\alpha$ of $E_p$.
    Thus we have
    \begin{equation}
        \label{eq:first-hadamard-on-diagonal}
        V^1(p,p)
        = -\frac{1}{6} \scal(p) \id_{E_p} - B(p).
    \end{equation}
\end{example}

%
%

\subsection{The Klein-Gordon Equation}
\label{subsec:klein-gordon-equation}

Even though in general the convergence of
\eqref{eq:Ansatz-for-Green-function} is hard to control and may even
fail in general there is one example where we can compute the Hadamard
coefficients explicitly and show weak$^*$ convergence of
\eqref{eq:Ansatz-for-Green-function}.

We consider again the flat Minkowski spacetime $(\mathbb{R}^n, \eta)$
but now the \Index{Klein-Gordon equation}
\begin{equation}
    \label{eq:klein-gordon-equation}
    \left( \dAlembert + m^2 \right) \phi = 0
\end{equation}
instead of $\dAlembert$ alone. As usual $m^2$ denotes a positive
constant. The physical meaning in quantum field theory of $m$ is that
of the mass of the particle described by
(\ref{eq:klein-gordon-equation}).

Since the metric $\eta$ is translation invariant and the operator
$\dAlembert + m^2$ is translation invariant as well, we only have to
compute the Hadamard coefficients at a single point $p \in
\mathbb{R}^n$ and can then translate everything. Thus we can choose
$p=0$. As already mentioned before, $\exp_p$ is just the addition with
$p$ whence
\begin{equation}
    \label{eq:flat-exp-at-0}
    \exp_0: T_0\mathbb{R}^n = \mathbb{R}^n 
    \longrightarrow \mathbb{R}^n
\end{equation}
is simply the identity map. Also the density function $\varrho_p$
becomes very simple as we have
\begin{equation}
    \label{eq:density-function-flat-case-is-constant-one}
    \varrho_p = 1
\end{equation}
for all $p$. Thus the recursion for the Hadamard coefficients
simplifies drastically. Finally, we note that the Klein-Gordon
operator $\dAlembert + m^2$ has already the normal form with $B =
m^2$. Thus the covariant derivative is the flat one and the parallel
transport is the identity. Therefor we have
\[
V_p^0 = 
\frac{1}{\sqrt{\varrho_p}} P_{p \rightarrow \argument} = \id
\]
and
\begin{align*}
    V_p^k(q) 
    &= - \frac{k}{\sqrt{\varrho_p(q)}} P_{p \rightarrow q}
    \int_0^1 \sqrt{\varrho_p} (\gamma_{p \rightarrow q}(\tau))
    P_{\gamma_{p \rightarrow q}, 0 \rightarrow \tau}^{-1}
    \left(
        D(V_p^{k-1})(\gamma_{p \rightarrow q}(\tau))
    \right)
    \tau^{k-1} \D \tau \\
    &= -k \int_0^1 
    D(V_p^{k-1})(p + \tau (q-p)) \tau^{k-1} \D \tau.
\end{align*}
Now $V_p^0$ is constant. We claim that, since $m^2$ is constant as
well, all Hadamard coefficients are constant, too. Indeed, assuming
this for $k-1$ shows that
\begin{align*}
    V_p^k(q) 
    &= -k \int_0^1
    D(V_p^{k-1})(p + \tau (q-p)) \tau^{k-1} \D \tau \\
    &= -k D(V_p^{k-1}) \int_0^1 \tau^{k-1} \D \tau \\
    &= - D(V_p^{k-1}) \\
    &= -m^2 V_p^{k-1},
\end{align*}
which is again constant. Thus by induction we conclude the following:
\begin{lemma}
    \label{lemma:hadamards-for-klein-gordon}
    \index{Klein-Gordon equation!Hadamard coefficients}%
    The Hadamard coefficients for the Klein-Gordon operator
    $\dAlembert + m^2$ on Minkowski spacetime are constant and
    explicitly given by
    \begin{equation}
        \label{eq:hadamards-for-klein-gordon}
        V_p^k = (-m^2)^k
    \end{equation}
    for $k \in \mathbb{N}_0$ and all points $p \in \mathbb{R}^n$.
\end{lemma}
This particularly simple form allows to determine the convergence of
\eqref{eq:Ansatz-for-Green-function} explicitly. We consider large $k
\in \mathbb{N}_0$ such that $R^\pm(2+2k)$ is actually a continuous
function. More precisely, we fix $N \in \mathbb{N}_0$ then for $2k
\geq n-2+2N$ the distribution $R^\pm(2+2k)$ is actually a $\Fun[N]$
function according to Lemma~\ref{lemma:riesz-function-continuity},
explicitly given by
\begin{equation}
    \label{eq:riesz-functions-large-k}
    R^\pm(2+2k)(x)
    = \frac{2^{1-(2+2k)} \pi^{\frac{2-n}{2}}}
    {\Gamma \left( \frac{2+2k}{2} \right)
      \Gamma \left( \frac{2+2k-n}{2}+1 \right)
    }
    \:
    \eta(x)^{\frac{2+2k-n}{2}}
    = \frac{\pi^{\frac{2-n}{2}}}
    {2^{2k-1} k! \Gamma \left( k+2-\frac{n}{2} \right)}
    \:
    \eta(x)^{k+1-\frac{n}{2}}
\end{equation}
for $x \in I^\pm(0)$ and $0$ elsewhere. We want to estimate
$R^\pm(2k)$ and its derivatives over a compactum $K \subseteq
\mathbb{R}^n$. To this end we compute the first partial derivatives of
$R^\pm(\alpha)$ explicitly. We know already
\begin{equation}
    \label{eq:first-derivative-of-riesz}
    \frac{\partial}{\partial x^{i_1}} R^\pm(\alpha)
    = \frac{1}{\alpha-2} R^\pm(\alpha-2) \eta_{i_1 j} x^j
    = \frac{1}{\alpha-2} R^\pm(\alpha-2) x_{i_1},
\end{equation}
where we use the notation
\begin{equation}
    \label{eq:index-pull-down}
    x_i = \eta_{ij} x^j.
\end{equation}
Thus we get
\begin{equation}
    \label{eq:second-derivatives-of-riesz}
    \frac{\partial^2}{\partial x^{i_1} \partial x^{i_2}}
    R^\pm(\alpha)
    = \frac{R^\pm(\alpha-4)}{(\alpha-2)(\alpha-4)}
    x_{i_1} x_{i_2}
    + \frac{R^\pm(\alpha-2)}{\alpha-2} \eta_{i_1 i_2},
\end{equation}
since clearly $\frac{\partial}{\partial x^{i_2}} x_{i_1} = \eta_{i_1
  i_2}$. Moving on from this we get
\begin{equation}
    \label{eq:third-derivatives-of-riesz}
    \begin{split}
        \frac{\partial^3}
        {\partial x^{i_1} \partial x^{i_2} \partial x^{i_4}}
        R^\pm(\alpha)
        &=
        \frac{R^\pm(\alpha-6)}{(\alpha-2)(\alpha-4)(\alpha-6)}
        x_{i_1} x_{i_2} x_{i_3} \\
        &\quad+
        \frac{R^\pm(\alpha-4)}{(\alpha-2)(\alpha-4)}
        \left(
            \eta_{i_1 i_3} x_{i_2} + \eta_{i_2 i_3} x_{i_1}
            + \eta_{i_1 i_2} x_{i_4}
        \right)
    \end{split}
\end{equation}
and
\begin{equation}
    \label{eq:fourth-derivatives-of-riesz}
    \begin{split}
        &\frac{\partial^4}{\partial x^{i_1} \partial x^{i_2}
          \partial x^{i_4} \partial x^{i_4}}
        R^\pm(\alpha) \\
        &\quad=
        \frac{ R^\pm(\alpha-8)}{(\alpha-2)(\alpha-4)(\alpha-6)(\alpha-8)}
        x_{i_1} x_{i_2} x_{i_3} x_{i_4} \\
        &\quad+
        \frac{ R^\pm(\alpha-6)}{(\alpha-2)(\alpha-4)(\alpha-6)}
        \left(
            \eta_{i_1 i_4} x_{i_2} x_{i_3}
            + \eta_{i_2 i_4} x_{i_1} x_{i_3}
            + \eta_{i_3 i_4} x_{i_1} x_{i_2}
            + \eta_{i_1 i_2} x_{i_3} x_{i_4}
            + \eta_{i_1 i_3} x_{i_2} x_{i_4}
        \right) \\
        &\quad+
        \frac{R^\pm(\alpha-4)}{(\alpha-2)(\alpha-4)}
        \left(
            \eta_{i_1 i_2} \eta_{i_3 i_4}
            + \eta_{i_1 i_3} \eta_{i_2 i_4}
            + \eta_{i_1 i_4} \eta_{i_2 i_3}
        \right).
    \end{split}
\end{equation}
Now we see how one can guess the general formula: For $\ell = 2r$
derivatives we have contributions of $\frac{1}{(\alpha-2) \cdots
  (\alpha-2(r+s))} R^\pm(\alpha-2(r+s))$ with coefficients consisting
of symmetrizations of $s$ factors $\eta$ and $2r - 2s$ factors of $x$
where only those symmetrizations are done which are not automatic,
i.e. $x_{i_1} x_{i_2}$ only occurs \emph{once} and not twice. For
$\ell = 2r+1$ we have the analogous statement. Summarizing this in a
more formalized way gives the following result:
\begin{proposition}[Taylor coefficients of $R^\pm(\alpha)$]
    \label{proposition:taylor-coefficients-or-riesz}
    \index{Riesz distribution!Taylor coefficients}%
    Let $\ell \in \mathbb{N}_0$ and set $r = \left[ \frac{\ell}{2}
    \right]$ whence $\ell = 2r$ or $\ell = 2r+1$ depending on $\ell$
    being even or odd. Then the partial derivatives of the Riesz
    distribution $R^\pm(\alpha)$ for $\alpha \notin \{2, 4, \ldots,
    4r\}$ are given by
    \begin{equation}
        \label{eq:partial-derivatives-or-riesz}
        \begin{split}
            &\frac{\partial^\ell}
            {\partial x^{i_1} \cdots \partial x^{i_\ell}}
            R^\pm(\alpha) \\
            &\quad=
            \sum_{s=0}^r
            \frac{R^\pm(\alpha-2\ell+2(r+s))}
            {(\alpha-2) \cdots (\alpha-2\ell+2(r-s))}
            \sum_{\sigma \in \Sym_{r,s}}
            \eta_{i_{\sigma(1)} i_{\sigma(2)}} \cdots 
            \eta_{i_{\sigma(2r-2s-1)} i_{\sigma(2(r-s))}}
            x_{i_{\sigma(2(r-s)+1)}} \cdots x_{i_{\sigma(\ell)}},
        \end{split}
    \end{equation}
    where $\Sym_{r,s}$ denotes those permutations of $\{1, \ldots,
    \ell\}$ such that
    \begin{equation}
        \label{eq:special-permutations}
        \begin{split}
            &\sigma(1) < \sigma(2), \ldots, \sigma(2(r-s)) \\
            &\sigma(3) < \sigma(4), \ldots, \sigma(2(r-s)) \\
            &\vdots \\
            &\sigma(2(r-s)-1) < \sigma(2(r-s)) \\
            \textrm{and}\quad  & \sigma(2(r-s)+1) < \sigma(2(r-s)+2) < 
            \ldots < \sigma(\ell).
        \end{split}
    \end{equation}
\end{proposition}
\begin{proof}
    The proof consists in a rather boring and tedious understanding of
    the above symmetrization procedure. Since we only need some
    qualitative consequences of
    \eqref{eq:partial-derivatives-or-riesz} we leave it as an
    exercise.
\end{proof}

\begin{remark}
    \label{remark:derivatives-of-riesz}
    The above result has again two possible interpretations. On one
    hand, \eqref{eq:partial-derivatives-or-riesz} holds for all
    $\alpha$ except for the poles in the sense of distributions. Even
    for the singular $\alpha$, the right hand side has an analytic
    continuation by the left hand side. On the other hand, for
    $\RE(\alpha)$ large enough, $R^\pm(\alpha)$ is a
    $\Fun[\ell]$-function and \eqref{eq:partial-derivatives-or-riesz}
    holds \emph{pointwise} in the sense of functions. By
    Lemma~\ref{lemma:riesz-function-continuity} this is the case for
    $\RE(\alpha) > n + 2\ell$.
\end{remark}

We consider now the case $\RE(\alpha) > n + 2\ell$ and want to use
\eqref{eq:partial-derivatives-or-riesz} to estimate the $\ell$-th
derivatives of the \emph{function} $R^\pm(\alpha)$ over a compactum $K
\subseteq \mathbb{R}^n$. Thus let $R > 0$ be large enough such that
\begin{equation}
    \label{eq:R-ball-around-compactum}
    K \subseteq B_R(0)
\end{equation}
for some \emph{Euclidean} ball around zero. The following is then
obvious from the definition of $R^\pm(\alpha)$ and gives a (rather
rough) estimate on the sup-norm of $R^\pm(\alpha)$ over $K$.
\begin{lemma}
    \label{lemma:supnorm-of-riesz}
    Let $K \subseteq  \mathbb{R}^n$ be compact and let $R > 0$ with $K
    \subseteq B_R(0)$. Then for $\RE(\alpha) > n$ we have
    \begin{equation}
        \label{eq:supnorm-of-riesz}
        \seminorm[K,0] (R^\pm(\alpha))
        \leq |c(\alpha,n)| R^{\RE(\alpha) - n}.
    \end{equation}
\end{lemma}
\begin{proof}
    For those $x \in I^\pm(0)$ we have
    \[
    |\eta(x,x)|
    = \left|
        (x^0)^2 - \sum_{i=1}^{n-1} (x^i)^2
    \right|
    \leq 
    | (x^0)^2 + \cdots + (x^n)^2|
    = R^2,
    \]
    and outside of $I^\pm(0)$, the function $R^\pm(\alpha)$ vanishes
    anyway.
\end{proof}

Taking derivatives into account we have the following estimate for
large $\RE(\alpha)$:
\begin{proposition}
    \label{proposition:seminorms-of-riesz}
    Let $K \subseteq \mathbb{R}^n$ be compact and let $R \geq 1$ with
    $K \subseteq B_R(0)$. Then for $\RE(\alpha) > n + 2\ell$ we have
    \begin{equation}
        \label{eq:seminorms-of-riesz}
        \seminorm[K,\ell] (R^\pm(\alpha))
        \leq \ell \cdot \ell! \cdot R^{\RE(\alpha)-n} \cdot
        \max \left\{
            |c(\alpha)|, \frac{|c(\alpha-2)|}{|\alpha-2|}, \ldots,
            \frac{|c(\alpha-2\ell)|}{|(\alpha-2) \cdots
              (\alpha-2\ell)|}
        \right\},
    \end{equation}
    with $c(\alpha) = c(\alpha,n)$ for abbreviation.
\end{proposition}
\begin{proof}
    From Proposition~\ref{proposition:taylor-coefficients-or-riesz} we
    know that for precisely $\ell'$ derivatives we have for $x \in K$
    \begin{align*}
        \left|
            \frac{\partial^{\ell'}}
            {\partial x^{i_1} \cdots \partial x^{i_{\ell'}}}
            R^\pm(\alpha)(x)
        \right|
        & \leq \sum_{s=0}^{r=\left[\frac{\ell'}{2}\right]}
        \frac{
          |R^\pm(\alpha - 2\ell' + 2(r-s))(x)|
        }
        {|(\alpha-2) \cdots (\alpha - 2 \ell' + 2(r-s))|}
        \sum | \eta \cdots \eta \cdot x \cdots x| \\
        & \leq \sum_{s=0}^{r=\left[\frac{\ell'}{2}\right]}
        \frac{
          |c(\alpha - 2 \ell' + 2(r-s))|
          R^{\RE(\alpha) - 2\ell' + 2(r-s) - n}
        }
        {|(\alpha-2) \cdots (\alpha - 2 \ell' + 2(r-s))|} \:
        \ell'! R^{\ell'},
    \end{align*}
    since in the sum over all allowed permutations we have at most
    $\ell'!$ factors (In fact, we always have much less, but a rough
    estimate will do the job). Moreover, every factor in $\eta \cdots
    \eta \cdot x \cdots x$ is clearly $\leq R$ in absolute value. Now
    since we assumed $R > 1$ we have for
    $r=\left[\frac{\ell'}{2}\right]$ and $s=0, \ldots, r$
    \[
    R^{\RE(\alpha) - 2\ell' + 2(r-s) - n} R^{\ell'}
    \leq R^{\RE(\alpha) - \ell' + 2\left[\frac{\ell'}{2}\right] - n}
    \leq R^{\RE(\alpha)-n},
    \]
    since $-\ell' + 2 \left[\frac{\ell'}{2}\right]$ is either $-1$ or
    $0$. Thus we can simplify this to 
    \begin{align*}
        \left|
            \frac{\partial^{\ell'}}
            {\partial x^{i_1} \cdots \partial x^{i_{\ell'}}}
            R^\pm(\alpha)(x)
        \right|
        & \leq R^{\RE(\alpha) - n} \ell'!
        \sum_{s=0}^{r=\left[\frac{\ell'}{2}\right]}
        \frac{
          |c(\alpha - 2 \ell' + 2(r-s))|
        }
        {|(\alpha-2) \cdots (\alpha - 2 \ell' + 2(r-s))|} \\
        & \leq \ell'! R^{\RE(\alpha) - n} \ell'
        \max_{s=0}^{r=\left[\frac{\ell'}{2}\right]}
        \left\{
            \frac{
              |c(\alpha - 2 \ell' + 2(r-s))|
            }
            {|(\alpha-2) \cdots (\alpha - 2 \ell' + 2(r-s))|}
        \right\}.
    \end{align*}
    For $\seminorm[K,\ell]$ we finally have to take the maximum of
    this expression over all $\ell' = 0, \ldots, \ell$. In the maximum
    over $s$ we can then simply take the largest of \emph{all},
    resulting in
    \[
    \seminorm[K,\ell] (R^\pm(\alpha))
    \leq \ell \cdot \ell! \cdot R^{\RE(\alpha)-n} \cdot
    \max \left\{
        |c(\alpha)|, \frac{|c(\alpha-2)|}{|\alpha-2|}, \ldots,
        \frac{|c(\alpha-2\ell)|}{|(\alpha-2) \cdots
          (\alpha-2\ell)|}
    \right\},
    \]
    which is what we wanted to show.
\end{proof}

Note that we only gave a rather rough estimate, which will
nevertheless be sufficient for the following. We specialize this now
to the case $\alpha = 2+2k$ with $k$ large enough such that $2+2k >
n+2\ell$. Then an even rougher estimate specializes
\eqref{eq:seminorms-of-riesz} to the following estimate:
\begin{corollary}
    \label{corollary:seminorm-for-special-rieszes}
    Let $\ell \in \mathbb{N}$ be fixed and $k \in \mathbb{N}$ such
    that $2+2k > 2\ell + n$ whence $R^\pm(2+2k)$ is $\Fun[\ell]$. Then
    we have for any compactum $K \subseteq \mathbb{R}^n$ with $K
    \subseteq B_R(0)$ for a sufficiently large $R > 1$
    \begin{equation}
        \label{eq:seminorm-for-special-rieszes}
        \seminorm[K,\ell] (R^\pm(2k))
        \leq \ell \ell! R^{2+2k-n}
        \frac{\pi^{\frac{1-n}{2}}}{2^{k-2} k!}.
    \end{equation}
\end{corollary}
\begin{proof}
    We compute explicitly by \eqref{eq:riesz-functions-large-k}
    \begin{align*}
        \frac{c(2+2k-2\ell)}{(2+2k-2) \cdots (2+2k-2\ell)}
        &=
        \frac{\pi^{\frac{2-n}{2}}}
        {
          2^{2(k-\ell)-1} (k-\ell)!
          \:
          \Gamma \left( k-\ell+2-\frac{n}{2} \right) \cdot
          2k \cdots 2(k - \ell +1)
        } \\
        &= \frac{\pi^{\frac{2-n}{2}}}{
          2^{2k-\ell-1} k!
          \:
          \Gamma \left( k-\ell+2-\frac{n}{2} \right)
        }.
    \end{align*}
    Now by assumption $2+2k > 2\ell+n$ whence on one hand $\ell \leq
    k$ since $n \geq 1$. Thus $2^{2k - \ell -1} \geq 2^{k-1}$.
    Moreover, $k - \ell +2 - \frac{n}{2} > 1$ whence by the monotonous
    growth of the $\Gamma$ function, see
    Figure~\ref{fig:gamma-function}, the smallest contribution of
    $\Gamma \left( k-\ell+2-\frac{n}{2} \right)$ occurs at $\Gamma
    \left( \frac{3}{2} \right) = \frac{1}{2} \sqrt{\pi}$.
    Thus we have
    \[
    \left|
        \frac{c(2+2k-2\ell)}{(2+2k-2) \cdots (2+2k-2\ell)}
    \right|
    \leq
    \frac{\pi^{\frac{2-n}{2}}}{2^{k-1} k! \frac{1}{2} \sqrt{\pi}}
    \]
    for all $\ell$. Inserting this into
    \eqref{eq:seminorms-of-riesz} gives the result.
\end{proof}

Again, estimating $\frac{1}{\Gamma \left( k-\ell+2-\frac{n}{2}
  \right)}$ by $\frac{1}{2} \sqrt{\pi}$ is very rough, in particular
as we are interested for fixed $\ell$ in the asymptotic behaviour for
$k \longrightarrow \infty$. The additional $\Gamma$-factor behaves
essentially like a $\frac{1}{k!}$ therefor improving the estimate
\eqref{eq:seminorm-for-special-rieszes} significantly. However, for
the following theorem, already \eqref{eq:seminorm-for-special-rieszes}
is sufficient.
\begin{theorem}[Green function of the Klein-Gordon operator]
    \label{theorem:green-function-for-klein-gordon}
    \index{Klein-Gordon equation!Green function}%
    Let $p \in \mathbb{R}^n$. Then the series
    \begin{equation}
        \label{eq:green-for-klein-gordon}
        \mathcal{R}^\pm(p)
        = \sum_{k=0}^\infty (-m^2)^k R^\pm(2+2k,p)
    \end{equation}
    converges in the weak$^*$ topology to the advanced and retarded
    Green function of the Klein-Gordon operator $\dAlembert + m^2$,
    respectively. More precisely, for $2+2k > 2\ell +n$ the series
    \begin{equation}
        \label{eq:Cl-tail-series}
        \sum_{2+2k > 2\ell+n} (-m^2)^k R^\pm(2+2k,p)
    \end{equation}
    converges in the $\Fun[\ell]$-topology to a $\Fun[\ell]$-function
    on $\mathbb{R}^n$. Finally, on $I^\pm(0)$ the series
    \eqref{eq:green-for-klein-gordon} converges in the
    $\Cinfty$-topology to a smooth function given by
    \begin{equation}
        \label{eq:smooth-green-in-inner-lightcone}
        \mathcal{R}^\pm(0)\At{I^\pm(0)}
        =
        \sum_{k=0}^\infty
        \frac{ \pi^{\frac{2-n}{2}} (-m^2)^k}{2^{2k-1} k!
          \:
          \Gamma \left( k+2-\frac{n}{2} \right)}
        \:
        \eta^{k+1-\frac{n}{2}}
    \end{equation}
    for $p=0$ from which the other $\mathcal{R}^\pm(p)$ can be
    obtained by translation.
\end{theorem}
\begin{proof}
    Clearly it suffices to show the convergence of
    \eqref{eq:Cl-tail-series} in the $\Fun[\ell]$ topology: since
    $\Fun[\ell](\mathbb{R}^n) \hookrightarrow
    \mathcal{D}'(\mathbb{R}^n)$ is continuously embedded, we can
    deduce the weak$^*$ convergence of
    \eqref{eq:green-for-klein-gordon} from that at once. To show
    \eqref{eq:Cl-tail-series}, we even show absolute convergence: let
    $K \subseteq \mathbb{R}^n$ be compact with $K \subseteq B_R(0)$
    for sufficiently large $R > 1$. Then
    \begin{align*}
        \sum_{2+2k > 2\ell+n} \seminorm[K,\ell]
        \left(
            (-m^2)^k R^\pm(2+2k,0)
        \right)
        & \leq \sum_{2+2k > 2\ell+n} (m^2)^k \ell \ell! R^{2+2k-n}
        \frac{\pi^{\frac{1-n}{2}}}{2^{k-2}k!}
        \leq c \sum_{2+2k > 2\ell+n} \frac{(m^2 R^2)^k}{k!}
    \end{align*}
    with some constant $c > 0$ depending on $\ell, R$. Since the
    series on the right is dominated by $\E^{m^2 R^2}$ we see that we
    indeed have absolute convergence with respect to
    $\seminorm[K,\ell]$ for all $K$. This shows
    $\Fun[\ell]$-convergence everywhere and hence weak$^*$
    convergence. Finally, on $I^\pm(0)$ the functions
    $R^\pm(2+2k)\at{I^\pm(0)}$ are always smooth whence the above
    result shows that they converge in \emph{all}
    $\Fun[\ell]$-topologies. But this means convergence in the
    $\Cinfty$-topology, establishing the last claim
    \eqref{eq:smooth-green-in-inner-lightcone}. By translation
    invariance, the convergence results also hold for any other $p \in
    \mathbb{R}^n$.
\end{proof}

\begin{remark}
    \label{remark:green-function-for-klein-gordon}
    Of course, there are much more straightforward techniques to
    obtain the Green functions for $\dAlembert + m^2$ on Minkowski
    spacetime. The standard approach is to use Fourier transformation
    techniques and to construct $\mathcal{R}^\pm(0)$ as even
    \emph{tempered} distribution on $\mathbb{R}^n$. In fact, for most
    applications in quantum field theory the momentum space
    representation of $\mathcal{R}^\pm(0)$ is needed anyway. However,
    our approach here is intrinsically geometric in the following
    sense: on a general spacetime Fourier transformation is not
    available, at least not in the naive way. Also, the above
    construction shows that $\mathcal{R}^\pm(0)$ depends
    \emph{analytically} on $m^2$: the series
    \eqref{eq:green-for-klein-gordon} being precisely the weak$^*$
    convergent Taylor expansion in the variable $m^2$ which may even
    be taken to be complex. This gives an entirely holomorphic family
    of distributions for $m^2 \in \mathbb{C}$.  Finally, the series
    \eqref{eq:smooth-green-in-inner-lightcone} can actually be
    expressed in terms of known transcendental functions, depending on
    the dimension $n$.
\end{remark}


%% file: fundamental.tex
%
%

In this section we construct out of the local Riesz distributions
$R^\pm(\alpha,p)$ and the corresponding Hadamard coefficients a
fundamental solution on a small neighborhood of $p \in M$. One
proceeds in two steps, first the formal series $\mathcal{R}^\pm(p)$ is
made to converge by brutally modifying the higher order terms. The
price paid is that the result is not yet a fundamental solution but
differs from the fundamental solution by a ``smoothing'' kernel, i.e.
one gets a \emph{parametrix} for $D$. In a second step one shows how
the parametrix can be changed to a fundamental solution by using an
appropriate geometric series of the smooth kernel. Again, we follows
essentially \cite{baer.ginoux.pfaeffle:2007a}.

In the following we fix a geodesically convex open subset $U'
\subseteq M$ and use the corresponding Riesz distributions
$R^\pm_{U'}(\alpha,p)$ which are now available for all $p \in U'$.
Moreover, by Theorem~\ref{theorem:hadamard-coefficients} the Hadamard
coefficients are now smooth sections
\begin{equation}
    \label{eq:smooth-hadamard}
    V^k \in \Secinfty\left(E^* \extensor E \at{U' \times U'}\right),
\end{equation}
out of which we obtain the \emph{formal} fundamental solution
\begin{equation}
    \label{eq:formal-fundamental-solution}
    \index{Fundamental solution!formal}%
    \mathcal{R}^\pm(p)
    = \sum_{k=0}^\infty V_p^k R^\pm_{U'}(2+2k,p)
\end{equation}
on $U'$. Of course, there is no reason to believe that
\eqref{eq:formal-fundamental-solution} converges in general, even not
in the weak$^*$ sense. However, the Riesz distributions
$R^\pm_{U'}(2+2k,p)$ are \emph{continuous} functions if $k$ is large
enough. In fact, by
Proposition~\ref{proposition:properties-of-riesz-on-U} we know that
$R^\pm_{U'}(2+2k,p)$ is at least continuous if $k > \frac{n}{2}$. Thus
we fix $N \in \mathbb{N}_0$ with $N > \frac{n}{2}$ and split the sum
\eqref{eq:formal-fundamental-solution} at $k=N$.

%
%

\subsection{The Approximate Fundamental Solution}
\label{satz:approximate-fundamental-solution}

The idea is now that the finite sum
\begin{equation}
    \label{eq:distributional-part-of-ansatz}
    \sum_{k=0}^{N-1} V_p^k R^\pm_{U'}(2+2k,p) \in
    \Secinfty_0\left(E^*\at{U'}\right)'
\end{equation}
is a well-defined distribution. On the other hand, this contribution
is believed to yield the most singular contribution to the yet to be
found fundamental solution responsible for the $\delta$-distribution
in $D \mathcal{R}^\pm(p) = \delta_p$. Thus the hope is that the
remaining, \emph{infinite} sum can be modified and made to converge
but yielding a less singular contribution than $\delta_p$, in fact
only a smooth one.

For technical reasons we will need a cutoff function $\chi \in
\Cinfty_0(\mathbb{R})$ with
\begin{equation}
    \label{eq:special-cutoff-function}
    \supp \chi \subseteq [-1,1],
    \quad
    0 \leq \chi \leq 1,
    \quad \textrm{and} \quad
    \chi \at{[-\frac{1}{2}, \frac{1}{2}]} = 1.
\end{equation}
For every choice of such a cutoff function, we have the following
technical lemma:
\begin{lemma}
    \label{lemma:seminorms-of-cutoff-function}
    Let $\ell \in \mathbb{N}$ and $\ell' > \ell + 1$. Then there are
    universal constants $c(\ell,\ell')$ such that for all $0 <
    \epsilon \leq 1$ one has
    \begin{equation}
        \label{eq:supnorm-of-rescaled-cutoff-function}
        \seminorm[K,0]
        \left(
            \frac{\D^\ell}{\D t^\ell}
            \left(
                \chi \left( \frac{t}{\epsilon} \right) t^{\ell'}
            \right)
        \right)
        \leq \epsilon c(\ell,\ell') \seminorm[K,\ell](\chi),
    \end{equation}
    where $K$ is any compactum containing $[-1,1]$.
\end{lemma}
\begin{proof}
    First note that $\chi \left( \frac{t}{\epsilon} \right) = 0$ for
    $\left| \frac{t}{\epsilon} \right| > 1$ and hence $|t| >
    \epsilon$. Thus the support of $t \mapsto \chi \left(
        \frac{t}{\epsilon} \right)$ is contained in $[-\epsilon,
    \epsilon] \subseteq [-1,1]$. It follows that in
    \eqref{eq:supnorm-of-rescaled-cutoff-function} we can safely
    replace the supremum over $K$ by a supremum over $\mathbb{R}$
    everywhere. In any case, we have by the Leibniz rule and the chain
    rule
    \begin{align*}
        \frac{\D^\ell}{\D t^\ell}
        \left(
            \chi \left( \frac{t}{\epsilon} \right) t^{\ell'}
        \right)
        &= \sum_{m=0}^\ell \binom{\ell}{m}
        \frac{\D^m}{\D t^m}\left(
            \chi \left(\frac{t}{\epsilon}\right)
        \right)
        \frac{\D^{\ell-m} t^{\ell'}}{\D t^{\ell-m}} \\
        &= \sum_{m=0}^\ell \binom{\ell}{m} \frac{1}{\epsilon^m}
        \frac{\D^m \chi}{\D t^m} \left( \frac{t}{\epsilon} \right)
        \ell' (\ell'-1) \cdots (\ell'-\ell+m+1) t^{\ell'-\ell+m}.
    \end{align*}
    Now for $|t| > \epsilon$ the factor $\frac{\D^m \chi}{\D t^m}
    \left( \frac{t}{\epsilon} \right)$ vanishes whence we find
    \begin{align*}
        \seminorm[K,0]
        \left(
            \frac{\D^\ell}{\D t^\ell}
            \left(
                \chi \left( \frac{t}{\epsilon} \right) t^{\ell'}
            \right)
        \right)
        & \leq \sup_t \sum_{m=0}^\ell \binom{\ell}{m}
        \ell' (\ell'-1) \cdots (\ell'-\ell+m+1)
        \epsilon^{\ell'-\ell}
        \left|
            \frac{\D^m \chi}{\D t^m} \left( \frac{t}{\epsilon} \right)
        \right| \\
        & \leq \epsilon \sum_{m=0}^\ell \binom{\ell}{m}
        \ell' (\ell'-1) \cdots (\ell'-\ell+m+1)
        \seminorm[K,\ell](\chi),
    \end{align*}
    since only $|t| \leq \epsilon$ contribute and
    $\epsilon^{\ell'-\ell} \leq \epsilon$ for $\ell' \geq \ell+1$.
\end{proof}

Since $U'$ is assumed to be convex, the Lorentz distance square is
defined on $U' \times U'$ and gives a smooth function $\eta \in
\Cinfty(U' \times U')$ by setting
\begin{equation}
    \label{eq:lorentz-distance-square}
    \eta(p,q)
    = \eta_p(q)
    = g_p \left( \exp_p^{-1}(q), \exp_p^{-1}(q) \right).
\end{equation}
We know from the proof of
Proposition~\ref{proposition:symmetry-of-riesz-on-U} that $\eta$ is
even a symmetric function
\begin{equation}
    \label{eq:symmetry-of-lorentz-distance-suare}
    \eta(p,q) = \eta(q,p).
\end{equation}
Finally, since $U'$ is assumed to be geodesically convex the geodesics
joining $p, q \in U'$ in $U'$ are unique. Thus we see that $\eta(p,q)
= 0$ iff the geodesic joining $p$ and $q$ is \emph{lightlike}. Since
the points $q$ which are in the image of $C(0) \subseteq T_pM$ under
$\exp_p$ are just $C_{U'}(p)$ we see that
\begin{equation}
    \label{eq:points-with-eta-is-zero}
    \eta^{-1}(\{0\})
    = \bigcup_{p \in U'} \{p\} \times C_{U'}(p).
\end{equation}

The idea is now to keep the series
\eqref{eq:formal-fundamental-solution} unchanged in a small, and in
fact only infinitesimal, neighborhood of the singular support,
i.e. the light cones $\eta^{-1}(\{0\})$, and modify it outside to
ensure convergence. To this end we will choose a sequence $\epsilon_j
\in (0,1]$ of cutoff parameters and consider the series
\begin{equation}
    \label{eq:cutoffed-series}
    \index{Cutoff parameters}%
    \begin{split}
        (p,q) \; \mapsto \; &\sum_{j=N}^\infty
        \chi \left( \frac{\eta(p,q)}{\epsilon_j} \right) V^j(p,q)
        R^\pm_{U'}(2+2j,p)(q) \\
        &= \begin{cases}
            \sum_{j=N}^\infty
            \chi \left( \frac{\eta(p,q)}{\epsilon_j} \right) V^j(p,q)
            c(2+2j,n) \eta(p,q)^{j+1-\frac{n}{2}}
            & \textrm{for} \quad q \in I^\pm_{U'}(p) \\
            0 & \textrm{else}.
        \end{cases}
    \end{split}
\end{equation}
Since $N \geq \frac{n}{2}$ all the terms in the modified (and
truncated) series are at least $\Fun[0]$. In fact, the $j$-th term is
at least $(j-N)$-times continuously differentiable by
Proposition~\ref{proposition:properties-of-riesz-on-U},
\refitem{item:differentiable-of-riesz-on-U-for-large-alpha} and by our
choice of $N$. For estimating the derivatives of $\chi \left(
    \frac{\eta}{\epsilon_j} \right)$ in a suitable way, we first
recall the following version of the chain rule:
\begin{lemma}
    \label{lemma:multi-index-chain-rule}
    Let $g: U \subseteq \mathbb{R}^n \longrightarrow \mathbb{R}$ and
    $f: \mathbb{R} \longrightarrow \mathbb{R}$ be smooth, then for
    every multi-index $I \in \mathbb{N}_0^n$
    \begin{equation}
        \label{eq:mutli-index-chain-rule}
        \frac{\partial^{|I|}}{\partial x^I} (f \circ g)
        =
        \sum_{\substack{r=1, \ldots, |I| \\ J_1, \ldots, J_r \leq I}}
        c_{J_1 \cdots J_r}^r \frac{\D^r f}{\D t^r} \circ g
        \;
        \frac{\partial^{|J_1|} g}{\partial x^{J_1}} \cdots
        \frac{\partial^{|J_r|} g}{\partial x^{J_r}}
    \end{equation}
    with some universal constants $c_{J_1 \cdots J_r}^r \in
    \mathbb{Q}$.
\end{lemma}
\begin{proof}
    This is clear by iterating the chain rule $|I|$ times. In fact,
    most of the $c_{J_1 \cdots J_r}^r$ are zero anyway.
\end{proof}

We shall now use an exhausting series of compacta for $U'$, i.e. we
choose compact subsets
\begin{equation}
    \label{eq:exhausting-series-of-compacta}
    K_0 \subseteq \ldots K_\ell \subseteq
    \mathring{K}_{\ell+1} \subseteq \ldots \subseteq U'
\end{equation}
with $U' = \bigcup_{\ell \geq 0} K_\ell$. This choice will give us
seminorms $\seminorm[K_\ell,k]$ for all involved bundles satisfying a
good estimate for natural pairings and
\begin{equation}
    \label{eq:filration-of-seminorms}
    \seminorm[K_\ell,k] \leq \seminorm[K_{\ell'}, k']
\end{equation}
for $\ell \leq \ell'$ and $k \leq k'$, see
Remark~\ref{remark:SymbolicSeminorms}. The filtration property
\eqref{eq:filration-of-seminorms} will turn out to be crucial. We
shall use the same exhausting sequence of compacta
\eqref{eq:exhausting-series-of-compacta} to obtain an exhausting
sequence $K_\ell \times K_\ell$ of $U' \times U'$ as well.

We consider now the function
\begin{equation}
    \label{eq:cutoff-times-lorentz-square}
    \chi_j^\pm(p,q) = \begin{cases}
        \chi \left( \frac{\eta(p,q)}{\epsilon_j} \right)
        \eta(p,q)^{j+1-\frac{n}{2}}
        & \textrm{for} \quad q \in I_{U'}^\pm(p) \\
        0 & \textrm{else},
    \end{cases}
\end{equation}
which is $\Fun[j-N]$ for $j \geq N > \frac{n}{2}$ according to the
properties of $\eta$ as in
Proposition~\ref{proposition:properties-of-riesz-on-U},
\refitem{item:differentiable-of-riesz-on-U-for-large-alpha}. We apply
now Lemma~\ref{lemma:seminorms-of-cutoff-function} and
Lemma~\ref{lemma:multi-index-chain-rule} to obtain the following
estimate:
\begin{lemma}
    \label{lemma:seminorms-of-cutoff-times-lorentz-square}
    Let $\ell, k \in \mathbb{N}_0$ and $j$ large enough such that $j-N
    \geq k$. Then we have
    \begin{equation}
        \label{eq:seminorms-of-cutoff-times-lorentz-square}
        \seminorm[K_\ell \times K_\ell, k](\chi_j^\pm)
        \leq \epsilon_j c(k,\ell,j),
    \end{equation}
    with constants $c(k,\ell,j) > 0$ independent of $\epsilon_j$
    satisfying
    \begin{equation}
        \label{eq:estimate-constants-filtration}
        c(k,\ell,j) \leq c(k',\ell',j)
    \end{equation}
    for $\ell \leq \ell'$ and $k \leq k'$.
\end{lemma}
\begin{proof}
    We have by the chain rule as in
    Lemma~\ref{lemma:multi-index-chain-rule}
    \begin{align*}
        \seminorm[K_\ell \times K_\ell,k](\chi_j^\pm)
        & \leq \sup_{\substack{ (p,q) \in K_\ell \times K_\ell \\ |I|
            \leq k}}
        \sum_{\substack{r \leq |I| \\ J_1, \ldots, J_r \leq I}}
        c_{J_1 \cdots J_r}^r
        \left|
            \frac{\D^r \left(
                  \chi \left(\frac{t}{\epsilon_j}\right)
                  t^{j+1-\frac{n}{2}}\right)}
            {\D t^r}
        \right|
        \left|
            \frac{\partial^{|J_1|} \eta}{\partial x^{J_1}}
        \right| \cdots
        \left|
            \frac{\partial^{|J_r|} \eta}{\partial x^{J_r}}
        \right| \\
        & \leq \epsilon_j \sup_{|I| \leq k}
        \sum_{\substack{r \leq |I| \\ J_1, \ldots, J_r \leq I}}
        c_{J_1 \cdots J_r}^r
        c\Big( r,j+1-\frac{n}{2} \Big)
        c_r \seminorm[K_\ell \times K_\ell, k](\eta)^r
    \end{align*}
    with $c_r = \max_{t \in \mathbb{R}} \left| \frac{\D^r \chi}{\D
          t^r} (t) \right| < \infty$. The maximum over $r \leq k$ is
    denoted by $\widetilde{c}_k$. The finitely many coefficients
    $c_{J_1 \cdots J_r}^r$ have a maximum depending only on $k$ and
    the sum has a certain maximal number of terms, again depending
    only on $k$. Thus there is a $\widetilde{\widetilde{c}}_k$ with
    \[
    \seminorm[K_\ell \times K_\ell,k](\chi_j^\pm)
    \leq \epsilon_j \widetilde{\widetilde{c}}_k \widetilde{c}_k
    \widetilde{c}
    \Big(k, j+1-\frac{n}{2}\Big)
    \max_{r \leq k} \seminorm[K_\ell \times K_\ell, k](\eta)^r,
    \]
    where $\widetilde{c}\big(k,j+1-\frac{n}{2} \big) = \max_{r \leq k}
    c(r,j+1-\frac{n}{2})$. But this is already the desired form since
    clearly $\widetilde{\widetilde{c}}_k$ increases with $k$,
    $\widetilde{c}_k$ increases with $k$ and so does
    $\widetilde{c}(k,j+1-\frac{n}{2})$. Finally, the last maximum also
    increases with $k$ and $\ell$ whence we can set
    \[
    c(k,\ell,j)
    = \widetilde{\widetilde{c}}_k \widetilde{c}_k
    \widetilde{c}\big(k,j+1-\frac{n}{2} \big)
    \max_{r \leq k} \seminorm[K_\ell \times K_\ell, k](\eta)^r,
    \]
    which will do the job.
\end{proof}

Together with the usual product rule for the seminorms
$\seminorm[K_\ell \times K_\ell, k]$ we obtain the following result:
\begin{lemma}
    \label{lemma:seminorms-of-cutoff-terms}
    Let $k, \ell \in \mathbb{N}_0$ and $j \geq N + k$. Then the $j$-th
    term of the series \eqref{eq:cutoffed-series} satisfies the
    estimate
    \begin{equation}
        \label{eq:seminorms-of-cutoff-terms}
        \seminorm[K_\ell \times K_\ell, k]
        \left(
            \chi \left(\frac{\eta}{\epsilon_j}\right) V^j
            R^\pm_{U'}(2+2j,\argument)
        \right)
        \leq \epsilon_j c(k,\ell,j) c(2+2j,n)
        \seminorm[K_\ell \times K_\ell,k](V^j).
    \end{equation}
\end{lemma}
\begin{proof}
    This is now easy from the product rule of the seminorms which
    gives a $k$-depending universal constant absorbed into the
    definition of $c(k, \ell, j)$ and the formula
    \eqref{eq:cutoffed-series} for the $j$-th term.
\end{proof}

Choosing the $\epsilon_j$ appropriately, this can be made arbitrarily
small in the following way:
\begin{proposition}
    \label{proposition:convergence-of-tail-series}
    ~
    \begin{propositionlist}
    \item\label{item:choice-of-epsilon} For any $j \ge N$ and every
        $\epsilon_j \in (0,1]$ such that
        \begin{equation}
            \label{eq:choice-of-epsilon}
            \epsilon_j \max_k
            \left\{
                c(k,j,j) c(2+2j,n) \seminorm[K_j \times K_j, k](V^j)
            \right\}
            \leq \frac{1}{2^j}
        \end{equation}
        the series \eqref{eq:cutoffed-series} converges absolutely in
        the $\Fun[0]$-topology to a continuous section of $E^*
        \extensor E \at{U' \times U'}$.
    \item\label{item:Ck-convergence-of-cutoffed-series} The series
        \eqref{eq:cutoffed-series} starting at $j \geq N+k$ converges
        absolutely in the $\Fun$-topology to a $\Fun$-section.
    \item\label{item:Cinfty-convergence-of-cutoffed-on-restriction}
        The series \eqref{eq:cutoffed-series} restricted to the open
        subset $U' \times U' \setminus \eta^{-1}(\{0\})$ converges in
        the $\Cinfty$-topology to a smooth section of $E^* \extensor E
        \at{U' \times U' \setminus \eta^{-1}(\{0\})}$.
    \end{propositionlist}
\end{proposition}
\begin{proof}
    For a fixed $j \geq N$ there are only finitely many $k \in
    \mathbb{N}_0$ with $j-N \geq k$ whence the maximum over the $k$'s
    in \eqref{eq:choice-of-epsilon} is well-defined. Thus we clearly
    can choose $\epsilon_j \in (0,1]$ to satisfy
    \eqref{eq:choice-of-epsilon}. Since we can take $k=0$, the second
    part implies the first as well. Thus let $k \in \mathbb{N}_0$ be
    arbitrary and consider the truncated series for $j \geq
    N+k$. First we note that every term is $\Fun$ whence we have to
    estimate their $\seminorm[K_\ell \times K_\ell,k]$-seminorms. We
    have for every $\ell \ge N+k$
    \begin{align*}
        \seminorm[K_\ell \times K_\ell, k]
        \left(
            \sum_{j \geq N+k} \chi_j V^j c(2+2j,n)
        \right)
        &\leq \sum_{j \geq N+k} \epsilon_j c(k,\ell,j) c(2+2j,n)
        \seminorm[K_\ell \times K_\ell, k](V^j) \\
        &\le \sum_{N+k \leq j \leq \ell} \epsilon_j c(k,\ell,j)
        c(2+2j,n) \seminorm[K_\ell \times K_\ell, k](V^j) \\
        &\quad+ \sum_{j > \ell} \epsilon_j c(k,\ell,j)
        c(2+2j,n) \seminorm[K_\ell \times K_\ell, k](V^j) \\
        &\le \textrm{const.} + \sum_{j > \ell} \epsilon_j
        c(k,j,j) c(2+2j,n) \seminorm[K_j \times K_j, k](V^j)
        \\
        &\le \textrm{const.} + \sum_{j > \ell} \frac{1}{2^j}
        < \infty,
    \end{align*}
    by the choice \eqref{eq:choice-of-epsilon} and the fact that for
    $j \ge \ell$ we can replace $\seminorm[K_\ell \times
    K_\ell,k](V^j)$ by $\seminorm[K_j \times K_j, k](V^j)$ as well as
    $c(k,\ell,j) \leq c(k,j,j)$ according to
    \eqref{eq:estimate-constants-filtration}. This shows absolute
    convergence with respect to $\seminorm[K_\ell \times K_\ell, k]$
    for \emph{all} $\ell \ge N+k$. But the compacta are increasing
    whence this shows absolute convergence in the $\Fun$-topology by
    the completeness of $\Sec\left(E^* \extensor E \at{U' \times
          U'}\right)$. Finally, we note that every term in
    \eqref{eq:cutoffed-series} is smooth on $U' \times U' \setminus
    \eta^{-1}(\{0\})$. Then we have $\Fun$-convergence by the second
    part for these restrictions, since omitting the first $k$ terms
    does not change the convergence behaviour of the series. But this
    means that we have convergence in the $\Cinfty$-topology.
\end{proof}

We can thus define an \emph{approximate fundamental solution}
$\widetilde{\mathcal{R}}^\pm_{U'}(p)$ by taking
\begin{equation}
    \label{eq:approximate-fundamental-solution}
    \index{Fundamental solution!approximate}%
    \widetilde{\mathcal{R}}^\pm_{U'}(p)
    = \sum_{j=0}^{N-1} V_p^j R^\pm_{U'}(2+2j,p)
    + \sum_{j=N}^\infty \chi \left(\frac{\eta_p}{\epsilon_j}\right)
    V^j_p R^\pm_{U'}(2+2j,p),
\end{equation}
after choosing the $\epsilon_j$ as in
Proposition~\ref{proposition:convergence-of-tail-series}. From the
support properties of the $R^\pm_{U'}(2+2j,p)$ and the above
convergence statement, we obtain the following result:
\begin{corollary}
    \label{corollary:order-and-supps-of-approx-fund-solution}
    Let the $\epsilon_j \in (0,1]$ be chosen to satisfy
    \eqref{eq:choice-of-epsilon}. Then
    \eqref{eq:approximate-fundamental-solution} is weak$^*$ convergent
    to a distributional section
    \begin{equation}
        \label{eq:approx-fund-solution-is-distributionl-section}
        \widetilde{\mathcal{R}}^\pm_{U'}(p)
        \in \Sec[-(n+1)]\left(E\at{U'}\right) \tensor E_p^*
    \end{equation}
    of global order $\leq n+1$ with
    \begin{equation}
        \label{eq:support-of-approx-fund-solution}
        \supp \widetilde{\mathcal{R}}^\pm_{U'}(p) \subseteq
        J^\pm_{U'}(p),
    \end{equation}
    \begin{equation}
        \label{eq:singsupp-of-approx-fund-solution}
        \singsupp \widetilde{\mathcal{R}}^\pm_{U'}(p) \subseteq
        C^\pm_{U'}(p).
    \end{equation}
\end{corollary}
\begin{proof}
    By Proposition~\ref{proposition:convergence-of-tail-series} the
    series converges in the $\Fun$-topology and hence also in the
    weak$^*$ topology.  Since the series is a continuous section it is
    of order $0$, the finitely many extra terms for $j \leq N-1$ are
    all of order $\leq n+1$ by
    Proposition~\ref{proposition:order-of-riesz-on-U},
    \refitem{item:boundary-for-order-for-alpha-greater-zero}. This
    shows \eqref{eq:approx-fund-solution-is-distributionl-section}.
    Since each term in \eqref{eq:approximate-fundamental-solution} has
    support in $J^\pm_{U'}(p)$ also the limit has support in
    $J^\pm_{U'}(p)$ as this is already a \emph{closed} subset of $U'$
    as we assume $U'$ to be geodesically convex. Moreover, the
    singular support of the first terms with $j \leq N-1$ is in
    $C^\pm_{U'}(p)$. By
    Proposition~\ref{proposition:convergence-of-tail-series}
    \refitem{item:Cinfty-convergence-of-cutoffed-on-restriction}, the
    series is smooth inside $I^\pm_{U'}(p)$ whence
    \eqref{eq:singsupp-of-approx-fund-solution} follows as well.
\end{proof}

Let us now determine in which sense
$\widetilde{\mathcal{R}}^\pm_{U'}(p)$ is an approximate solution.
Since the series converges in the weak$^*$ sense we can apply $D =
\dAlembert^\nabla + B$ term by term thanks to the continuity of
differential operators, see
Theorem~\ref{theorem:differentiation-of-gensec},
\refitem{item:diffops-are-continuous}. In our situation we can even
argue in the sense of functions if we start the series at $N+2$
because then we have $\Fun[2]$-convergence for which $D$ is continuous
as well. In any case we get
\begin{align}
    D \widetilde{\mathcal{R}}^\pm_{U'}(p)
    & = \sum_{j=0}^{N-1} D
    \left(
        V^j_p R^\pm_{U'}(2+2j,p)
    \right)
    + \sum_{j=N}^\infty D
    \left(
        \chi \left(\frac{\eta_p}{\epsilon_j}\right) V_p^j
        R^\pm_{U'}(2+2j,p)
    \right) \nonumber \\
    & = \delta_p + D(V_p^{N-1}) R^\pm_{U'}(2N,p)
    + \sum_{j=N}^\infty D
    \left(
        \chi \left(\frac{\eta_p}{\epsilon_j}\right) V_p^j
        R^\pm_{U'}(2+2j,p)
    \right),
    \label{eq:D-on-approx-fund-solution}
\end{align}
thanks to the transport equations for $V_p^j$. Indeed, the transport
equations, by their very construction, yield Hadamard coefficients
$V_p^j$ such that
\begin{align}
    \sum_{j=0}^{N-1}
    & D \left( V_p^j R^\pm_{U'}(2+2j,p) \right) \nonumber \\
    &= D(V^0_p) R^\pm_{U'}(2,p)
    + 2 \nabla^E_{\gradient R^\pm_{U'}(2,p)} V^0_p
    + V^0_p \dAlembert R^\pm_{U'}(2,p) \nonumber \\
    &\quad+ D(V^1_p) R^\pm_{U'}(4,p)
    + 2 \nabla^E_{\gradient R^\pm_{U'}(4,p)} V^1_p
    + V^1_p \dAlembert R^\pm_{U'}(4,p) \nonumber \\
    &\quad+ \cdots \nonumber \\
    &\quad+ D(V^{N-1}_p) R^\pm_{U'}(2N,p)
    + 2 \nabla^E_{\gradient R^\pm_{U'}(2N,p)} V^{N-1}_p
    + V^{N-1}_p \dAlembert R^\pm_{U'}(2N,p) \nonumber \\
    &= \delta_p \nonumber \\
    & \quad+ 2 \nabla^E_{\gradient R^\pm_{U'}(4,p)} V^1_p
    + V^1_p \dAlembert R^\pm_{U'}(4,p)
    + D(V^0_p) R_{U'}(2,p) \nonumber \\
    &\quad+ \cdots \nonumber \\
    &\quad + 2 \nabla^E_{\gradient R^\pm_{U'}(2N,p)} V^{N-1}_p
    + V^{N-1}_p \dAlembert R^\pm_{U'}(2N,p)
    + D(V^{N-2}_p) R^\pm_{U'}(2N-2,p) \nonumber \\
    &\quad+ D(V^{N-1}_p) R^\pm_{U'}(2N,p) \nonumber \\
    &= \delta_p + 0 + \cdots + 0 + D(V^{N-1}_p) R^\pm_{U'}(2N,p)
    \label{eq:D-on-first-part-of-approx-solution}
\end{align}
for arbitrary $N$ by \eqref{eq:lowest-hadamard-equation} and
\eqref{eq:higher-hadamard-equations}. We consider now the remaining
sum over $j$ in \eqref{eq:D-on-approx-fund-solution} and get by the
Leibniz rule for $D$
\begin{align}
    D \left(
        \chi \left(\frac{\eta_p}{\epsilon_j}\right) V_p^j
        R^\pm_{U'}(2+2j,p)
    \right)
    &= \dAlembert
    \left(
        \chi \left(\frac{\eta_p}{\epsilon_j}\right)
    \right)
    V_p^j R^\pm_{U'}(2+2j,p)
    +
    2 \nabla^E_{
      \gradient \chi \big(\frac{\eta_p}{\epsilon_j}\big)
    }
    \left(
        V_p^j R^\pm_{U'}(2+2j,p)
    \right) \nonumber \\
    &\quad + \chi \left(\frac{\eta_p}{\epsilon_j}\right)
    D \left(
        V_p^j R^\pm_{U'}(2+2j,p)
    \right).
    \label{eq:D-on-terms-in-second-part-of -approx-solution}
\end{align}
By the transport equations we have
\begin{align}
    D \left(
        V_p^j R^\pm_{U'}(2+2j,p)
    \right)
    &= D(V_p^j) R^\pm_{U'}(2+2j,p)
    + 2 \nabla^E_{\gradient R^\pm_{U'}(2+2j,p)} V_p^j
    + V_p^j \dAlembert R^\pm_{U'}(2+2j,p) \nonumber \\
    &= D(V^j_p) R^\pm_{U'}(2+2j,p)
    - D(V^{j-1}_p) R^\pm_{U'}(2j,p).
    \label{eq:D-on-hadamard-times-riesz}
\end{align}
By shifting the summation index appropriately, we get
\begin{align}
    D \widetilde{\mathcal{R}}^\pm_{U'}(p) - \delta_p
    & = D(V_p^{N-1}) R^\pm_{U'}(2N,p)
    + \sum_{j=N}^\infty \dAlembert
    \left(
        \chi \left(\frac{\eta_p}{\epsilon_j}\right)
    \right)
    V_p^j R^\pm_{U'}(2+2j,p) \nonumber \\
    &\quad+
    \sum_{j=N}^\infty
    2 \nabla^E_{
      \gradient \chi \big(\frac{\eta_p}{\epsilon_j}\big)
    }
    \left(
        V_p^j R^\pm_{U'}(2+2j,p)
    \right)
    + \sum_{j=N}^\infty
    \chi \left(\frac{\eta_p}{\epsilon_j}\right)
    D(V^j_p) R^\pm_{U'}(2+2j,p) \nonumber \\
    &\quad-
    \sum_{j=N}^\infty
    \chi \left(\frac{\eta_p}{\epsilon_j}\right)
    D(V^{j-1}_p) R^\pm_{U'}(2j,p) \nonumber \\
    & = \left(
        1 - \chi \left(\frac{\eta_p}{\epsilon_N}\right)
    \right)
    D(V_p^{N-1}) R^\pm_{U'}(2N,p)
    + \sum_{j=N}^\infty \dAlembert
    \left(
        \chi \left(\frac{\eta_p}{\epsilon_j}\right)
    \right)
    V_p^j R^\pm_{U'}(2+2j,p) \nonumber \\
    &\quad+
    \sum_{j=N}^\infty
    2 \nabla^E_{
      \gradient \chi \big(\frac{\eta_p}{\epsilon_j}\big)
    }
    \left(
        V_p^j R^\pm_{U'}(2+2j,p)
    \right) \nonumber \\
    &\quad+
    \sum_{j=N}^\infty
    \left(
        \chi \left(\frac{\eta_p}{\epsilon_j}\right)
        - \chi \left(\frac{\eta_p}{\epsilon_{j+1}}\right)
    \right)
    D(V_p^j) R^\pm_{U'}(2+2j,p) \nonumber \\
    & = \left(
        1 - \chi \left(\frac{\eta_p}{\epsilon_N}\right)
    \right)
    D (V_p^{N-1}) R^\pm_{U'}(2N,p)
    + \Sigma_1 + \Sigma_2 + \Sigma_3,
    \label{eq:D-on-approx-is-delta-plus-this-right-here}
\end{align}
where we abbreviated the last three series with $\Sigma_1, \Sigma_2,$
and $\Sigma_3$, respectively. In order to investigate these three
series we need the following technical lemma:
\begin{lemma}
    \label{lemma:vanishing-of-terms-in-D-approx-solution}
    Let $\epsilon_j \in (0,1]$ be chosen as in
    \eqref{eq:choice-of-epsilon}.
    \begin{lemmalist}
    \item \label{item:one-minus-chi-vanishes} The function $(p,q)
        \mapsto 1 - \chi \left(\frac{\eta(p,q)}{\epsilon_N}\right)$
        vanishes on an open neighborhood of $\eta^{-1}(\{0\})$.
    \item \label{item:gradient-chi-vanishes} The vector field $U'
        \times U' \ni (p,q) \mapsto (\id \extensor \gradient) \left(
            \chi \left(\frac{\eta(p,q)}{\epsilon_j}\right)\right) \in
        T_qU'$ vanishes on an open neighborhood of $\eta^{-1}(\{0\})$.
    \item \label{item:dAlemebrt-of-chi-vanishes} The function $(\id
        \extensor \dAlembert) \left(\chi
            \left(\frac{\eta}{\epsilon_j}\right)\right)$ vanishes on
        an open neighborhood of $\eta^{-1}(\{0\})$.
    \item \label{item:difference-in-last-sum-vanishes} The function
        $\chi \left(\frac{\eta}{\epsilon_j}\right) - \chi
        \left(\frac{\eta}{\epsilon_{j+1}}\right)$ vanishes on an open
        neighborhood of $\eta^{-1}(\{0\})$.
    \item \label{item:remainder-is-smooth} The section $\left(1 - \chi
            \left(\frac{\eta}{\epsilon_N}\right)\right) D(V^{N-1})
            R^\pm_{U'}(2N,\argument)$ as well as all the sections
        in the three series $\Sigma_1, \Sigma_2,$ and $\Sigma_3$ are
        smooth on $U' \times U'$.
    \end{lemmalist}
\end{lemma}
\begin{proof}
    We consider the open neighborhood
    \[
    S_j = \left\{
        (p,q) \in U' \times U'
        \; \Big| \;
        -\frac{\epsilon_j}{2} < \eta(p,q) < \frac{\epsilon_j}{2}
    \right\} \subseteq U' \times U'
    \]
    of $\eta^{-1}(\{0\})$. Clearly, by continuity of $\eta$ this is an
    open neighborhood, see
    Figure~\ref{fig:open-neighborhood-of-lightcone} for the flat
    analogue.
    \begin{figure}
        \centering
        \input{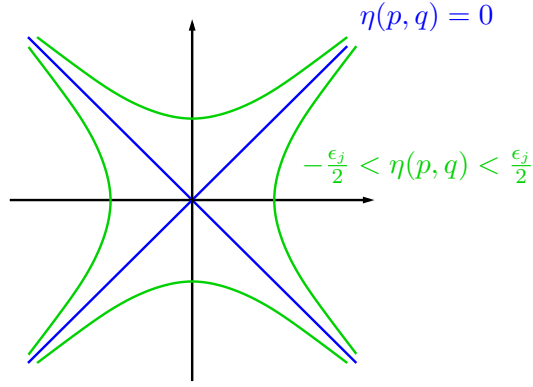}
        \caption{The open neighborhood $S_j$ of $\eta^{-1}(\{0\})$ in
          the flat case.}
        \label{fig:open-neighborhood-of-lightcone}
    \end{figure}
    Since the cutoff function $\chi$ is constant and equal to one on
    $\left[ -\frac{1}{2}, \frac{1}{2} \right]$, we see that the
    function $\chi \big( \frac{\eta}{\epsilon_j} \big)$ is equal to
    one on the \emph{open} $S_j$. From this
    \refitem{item:one-minus-chi-vanishes} follows at once. Thus also
    the gradient vanishes on $S_j$ whence
    \refitem{item:gradient-chi-vanishes} and
    \refitem{item:dAlemebrt-of-chi-vanishes} follow. Since $S_j \cap
    S_{j+1}$ is still an open neighborhood of $\eta^{-1}(\{0\})$, we
    get \refitem{item:difference-in-last-sum-vanishes}. But this means
    that the prefactors in all the above terms vanish on an open
    neighborhood of $\eta^{-1}(\{0\})$ which was the only place where
    the Riesz distributions $R^\pm_{U'}(2+2j,\argument)$ were
    non-smooth. Thus \refitem{item:remainder-is-smooth} follows, too.
\end{proof}

This lemma suggests that the weak$^*$ convergence of all the three
sums $\Sigma_1, \Sigma_2$, and $\Sigma_3$, which we already know, can
be sharpened to a $\Cinfty$-convergence: in this case the defect of
$\widetilde{\mathcal{R}}^\pm_{U'}(p)$ of being a fundamental solution
would be just a \emph{smooth section} and not a general,
distributional section. After possibly redefining the $\epsilon_j$
this can indeed be achieved as we shall see now.

\begin{figure}
    \centering
    \input{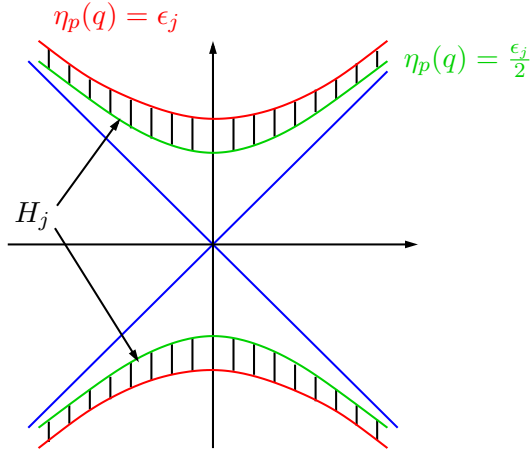}
    \caption{The set $H_j \cap \{p\} \times U'$.}
    \label{fig:the-mass-shell}
\end{figure}

First we note that the functions $\chi
\left(\frac{\eta}{\epsilon_j}\right)$ are only interesting in the
following subset
\begin{equation}
    \label{eq:the-subset-Hj}
    H_j = \left\{
        (p,q) \in U' \times U'
        \; \Big| \;
        \frac{\epsilon_j}{2} \leq \eta(p,q) \leq \epsilon_j
    \right\}.
\end{equation}
Indeed, for $\eta(p,q) > \epsilon_j$ the cutoff function produces a
zero, for $\eta(p,q) < \frac{\epsilon_j}{2}$ the function is
identically one until $\eta(p,q) < - \frac{\epsilon_j}{2}$. But for
negative $\eta(p,q)$ the definition of $R^\pm_{U'}(2+2j,p)(q)$ gives
already zero. Thus we only get contributions to each of the series
$\Sigma_1$ and $\Sigma_2$ from $H_j$ for the $j$-th
term. Geometrically, $H_j \cap \{p\} \times U'$ looks like a thick
mass shell, see Figure~\ref{fig:the-mass-shell}.
It follows that for the $j$-th term in $\Sigma_1$ or $\Sigma_2$ we get
only contributions from the compactum $K_\ell \times K_\ell \cap H_j$
for the seminorm $\seminorm[K_\ell \times K_\ell, k]$.

We start now estimating the $\seminorm[K_\ell \times K_\ell, k]$ of
the $j$-th term in the sum $\Sigma_2$. To this end we first estimate
the function $\eta$ on $K_\ell \times K_\ell \cap H_j$ as follows.
\begin{lemma}
    \label{lemma:seminorm-estimate-of-eta-on-Hj}
    Let $j \geq N$ and $k, \ell \in \mathbb{N}_0$ arbitrary. Then
    \begin{equation}
        \label{eq:seminorm-estimate-of-eta-on-Hj}
        \seminorm[K_\ell \times K_\ell \cap H_j, k+1]
        \left(
            \eta^{j+1-\frac{n}{2}}
        \right)
        \leq
        d(k,\ell,j) \epsilon_j^{j-\frac{n}{2}-k},
    \end{equation}
    with some constants $d(k,\ell,j) > 0$ such that
    \begin{equation}
        \label{eq:filtration-of-d-constants}
        d(k,\ell,j) \leq d(k', \ell', j)
    \end{equation}
    for $k \leq k'$ and $\ell \leq \ell'$.
\end{lemma}
\begin{proof}
    By the chain rule as in Lemma~\ref{lemma:multi-index-chain-rule}
    we have
    \begin{align*}
        &\seminorm[K_\ell \times K_\ell \cap H_j, k+1]
        \left(
            \eta^{j+1-\frac{n}{2}}
        \right) \\
        & \leq
        \sup_{\substack{(p,q) \in K_\ell \times K_\ell \cap H_j
            \\ |I| \leq k+1}}
        \sum_{\substack{r \leq |I| \\ J_1, \ldots, J_r \leq I}}
        c_{J_1 \cdots J_r}^r
        \left|
            \frac{\D t^{j+1-\frac{n}{2}}}{\D t^r}
        \right|_{t=\eta(p,q)}
        \left|
            \frac{\partial^{|J_1|} \eta}{\partial x^{J_1}}
        \right|
        \cdots
        \left|
            \frac{\partial^{|J_1|} \eta}{\partial x^{J_1}}
        \right| \\
        & \leq
        \sup_{\substack{|I| \leq k+1 \\
            \frac{\epsilon_j}{2} \leq t \leq \epsilon_j}}
        \sum_{\substack{r \leq |I| \\ J_1, \ldots, J_r \leq I}}
        c_{J_1 \cdots J_r}^r
        \left|
            \big( j+1-\frac{n}{2} \big) \cdots
            \big( j+1-\frac{n}{2}-r+1 \big)
            t^{j+1-\frac{n}{2}-r}
        \right|
        \left(
            \seminorm[K_\ell \times K_\ell \cap H_j, k+1](\eta)
        \right)^r \\
        & \leq
        \sup_{|I| \leq k+1}
        \sum_{\substack{r \leq |I| \\ J_1, \ldots, J_r \leq I}}
        c_{J_1 \cdots J_r}^r
        \left(
            \frac{\epsilon_j}{2}
        \right)^{j+1-\frac{n}{2}-(k+1)}
        \left(
            \seminorm[K_\ell \times K_\ell \cap H_j, k+1](\eta)
        \right)^r \\
        & \leq
        \epsilon_j^{j-\frac{n}{2}+k}
        \underbrace{
          \max_{|I| \leq k+1}
          \sum_{\substack{r \leq |I| \\ J_1, \ldots, J_r \leq I}}
          c_{J_1 \cdots J_r}^r
          \frac{1}{2^{j+1-\frac{n}{2}-(k+1)}}
          \left(
              \seminorm[K_\ell \times K_\ell \cap H_j, k+1](\eta)
          \right)^r
        }_{d(k,\ell,j)}.
    \end{align*}
    Note that the supremum over $t$ and $r \leq k+1$ of
    $t^{j+1-\frac{n}{2}-r}$ is obtained for the smallest $t =
    \frac{\epsilon_j}{2}$ and the largest $r=k+1$. The constants
    $d(k,\ell,j)$ clearly grow if the compactum $K_\ell$ is replaced
    by the bigger one $K_{\ell'}$. They also grow if we allow larger
    $k$.
\end{proof}

This can now be used to estimate the $j$-th term of the series
$\Sigma_2$. We have the following result:
\begin{lemma}
    \label{lemma:seminorms-of-terms-in-second-sum}
    Let $k, \ell \in \mathbb{N}_0$ and $j \geq N$. Then we have
    \begin{equation}
        \label{eq:seminorms-of-terms-in-second-sum}
        \begin{split}
            & \seminorm[K_\ell \times K_\ell, k]
            \left(
                \nabla^E_{
                  \gradient \chi \big(\frac{\eta}{\epsilon_j}\big)}
                \left(
                    V^j R^\pm_{U'}(2+2j,\argument)
                \right)
            \right) \\
            & \leq
            c_k c(2+2j,n) d(k+1,\ell,j)
            \seminorm[K_\ell \times K_\ell,k+1](V^j)
            \max_{r \leq k+1}
            \seminorm[K_\ell \times K_\ell, k+1](\eta)^r
            \cdot \epsilon_j^{j-\frac{n}{2}-2k-1}.
        \end{split}
    \end{equation}
\end{lemma}
\begin{proof}
    We simply compute
    \begin{align*}
        \seminorm[K_\ell \times K_\ell, k]&
        \left(
            \nabla^E_{
              \gradient \chi \big(\frac{\eta}{\epsilon_j}\big)}
            \left(
                V^j R^\pm_{U'}(2+2j,\argument)
            \right)
        \right) \\
        &=
        \seminorm[K_\ell \times K_\ell \cap H_j, k]
        \left(
            \nabla^E_{
              \gradient \chi \big(\frac{\eta}{\epsilon_j}\big)}
            \left(
                V^j R^\pm_{U'}(2+2j,\argument)
            \right)
        \right) \\
        &\leq
        c_{k,\ell}
        \seminorm[K_\ell \times K_\ell \cap H_j, k+1]
        \left(
            \chi \left(\frac{\eta}{\epsilon_j}\right)
        \right)
        \seminorm[K_\ell \times K_\ell \cap H_j, k+1](V^j)
        \seminorm[K_\ell \times K_\ell \cap H_j, k+1]
        (R^\pm_{U'}(2+2j,\argument)),
    \end{align*}
    since we need one order of differentiation for the gradient and
    one for the covariant derivative. In the constant $c_{k,\ell}$ the
    estimates of the derivatives of the metric, the connection, the
    Leibniz rule, etc. enter. Note that since these quantities are
    smooth everywhere, we can take the supremum over $K_\ell \times
    K_\ell$ whence $c_{k,\ell}$ does not depend on $j$. Now by the
    chain rule as in Lemma~\ref{lemma:multi-index-chain-rule} we have
    \begin{align*}
        \seminorm[K_\ell \times K_\ell \cap H_j, k+1]
        \left(
            \chi \left(\frac{\eta}{\epsilon_j}\right)
        \right)
        &\leq
        \sup_{\substack{(p,q) \in K_\ell \times K_\ell \cap H_j \\
            |I| \leq k+1}}
        \sum_{\substack{r \leq |I| \\ J_1, \ldots, J_r \leq I}}
        c_{J_1 \cdots J_r}^r
        \left|
            \frac{\D^r \chi}{\D t^r}
        \right|_{t=\frac{\eta}{\epsilon_j}}
        \frac{1}{\epsilon_j^r}
        \left|
            \frac{\partial^{|J_1|} \eta}{\partial x^{J_1}}
        \right|
        \cdots
        \left|
            \frac{\partial^{|J_r|} \eta}{\partial x^{J_r}}
        \right| \\
        &\leq
        \frac{1}{\epsilon_j^{k+1}} c_k
        \max_{r \leq k+1}
        \seminorm[K_\ell \times K_\ell, k+1](\eta)^r,
    \end{align*}
    where the sum over the $c_{J_1 \cdots J_r}^r$ as well as the
    supremum over the $r$-th derivatives of $\chi$ are combined into
    the constant $c_k$. For the seminorm of $R^\pm_{U'}(2+2j,\argument)$
    we get
    \begin{align*}
        \seminorm[K_\ell \times K_\ell \cap H_j, k+1]
        (R^\pm_{U'}(2+2j,\argument))
        &\leq
        \seminorm[K_\ell \times K_\ell \cap H_j, k+1]
        \left(
            c(2+2j,n) \eta^{j+1-\frac{2}{2}}
        \right) \\
        &\leq
        c(2+2j,n) d(k+1,\ell,j)
        \epsilon_j^{j-\frac{n}{2}-k}
    \end{align*}
    by Lemma~\ref{lemma:seminorm-estimate-of-eta-on-Hj}. Putting
    things together we obtain
    \begin{align*}
        \seminorm[K_\ell \times K_\ell, k]&
        \left(
            \nabla^E_{
              \gradient \chi \big(\frac{\eta}{\epsilon_j}\big)
            }
            V^j R^\pm_{U'}(2+2j,\argument)
        \right) \\
        & \leq
        c_k c(2+2j,n) d(k+1,\ell,j)
        \epsilon_j^{j-\frac{n}{2}-2k-1}
        \seminorm[K_\ell \times K_\ell,k+1](V^j)
        \max_{r \leq k+1}
        \seminorm[K_\ell \times K_\ell, k+1](\eta)^r.
    \end{align*}
\end{proof}
\begin{lemma}
    \label{lemma:absolute-convergence-of-second-sum}
    Let $j \geq N$. Choose $\epsilon_j \in (0,1]$ such that in
    addition to \eqref{eq:choice-of-epsilon}
    \begin{equation}
        \label{eq:new-choice-of-epsilon}
        \epsilon_j
        \max_{\substack{\ell \leq j \\
            k \leq \frac{1}{2}\left(j-\frac{n}{2}-1\right)}}
        c_k c(2+2j,n) d(k+1,\ell,j)
        \seminorm[K_\ell \times K_\ell,k+1](V^j)
        \max_{r \leq k+1}
        \seminorm[K_\ell \times K_\ell, k+1](\eta)^r
        < \frac{1}{2^j}.
    \end{equation}
    Then the sum $\Sigma_2$ converges absolutely in the
    $\Cinfty$-topology to some $\Sigma_2 \in \Secinfty\left(E^*
        \extensor E \at{U' \times U'}\right)$.
\end{lemma}
\begin{proof}
    First we note that we can indeed find $\epsilon_j \in (0,1]$
    meeting the requirement \eqref{eq:new-choice-of-epsilon}. Then we
    have for fixed $k, \ell$ the estimate
    \begin{align*}
        &\seminorm[K_\ell \times K_\ell, k]
        \left(
            \sum_{j \geq N}
            2 \nabla^E_{
              \gradient \chi \big(\frac{\eta}{\epsilon_j}\big)
            }
            V^j R^\pm_{U'}(2+2j,\argument)
        \right) \\
        &\leq
        \seminorm[K_\ell \times K_\ell, k]
        \left(
            \sum_{j=N}^{j_0-1}
            2 \nabla^E_{
              \gradient \chi \big(\frac{\eta}{\epsilon_j}\big)
            }
            V^j R^\pm_{U'}(2+2j,\argument)
        \right) \\
        &\quad+ 2 \sum_{j \geq j_0}
        c_k c(2+2j,n) d(k+1,\ell,j)
        \seminorm[K_\ell \times K_\ell,k+1](V^j)
        \max_{r \leq k+1}
        \seminorm[K_\ell \times K_\ell, k+1](\eta)^r
        \cdot \epsilon_j^{j-\frac{n}{2}-2k-1} \\
        &\leq
        \textrm{const.} + 2 \sum_{j \geq j_0} \frac{1}{2^j}
        < \infty,
    \end{align*}
    provided we set $j_0$ larger than $\ell$ and such that $j_0 -
    \frac{n}{2} - 2k - 1 \geq 1$, which is clearly possible. In this
    case $\epsilon_j^{j-\frac{n}{2}-2k-1} \leq \epsilon_j$ for $j \geq
    j_0$, and we can use \eqref{eq:new-choice-of-epsilon} to get the
    estimate. But this shows absolute convergence in the seminorm
    $\seminorm[K_\ell \times K_\ell,k]$ as the finitely many terms
    with $N \leq j \leq j_0-1$ do not matter. Since $\ell$ and $k$
    were arbitrary we get $\Cinfty$-convergence. Note that it is
    crucial that each term of $\Sigma_2$ is already smooth, quite
    differently from the ideas in
    Proposition~\ref{proposition:convergence-of-tail-series}.
\end{proof}

By a completely analogous argument one can estimate the terms in the
sum $\Sigma_1$ and show that again finitely many conditions on each
$\epsilon_j \in (0,1]$ yield $\Cinfty$-convergence also of
$\Sigma_1$. We do not write down the explicit condition but leave this
as an exercise. The result is the following:
\begin{lemma}
    \label{lemma:absolute-convergence-of-first-sum}
    There are choices of $\epsilon_j \in (0,1]$ analogous to
    \eqref{eq:new-choice-of-epsilon} such that the sum $\Sigma_1$
    converges absolutely in the $\Cinfty$-topology to some section
    $\Sigma_1 \in \Secinfty\left(E^* \extensor E \at{U' \times
          U'}\right)$.
\end{lemma}

Finally, we consider the third sum $\Sigma_3$. Here the argument is
slightly different leading nevertheless to the same consequences.
\begin{lemma}
    \label{lemma:seminorm-estimate-for-terms-in-third-sum}
    Let $\ell, k \in \mathbb{N}_0$ and let $j \geq N$ satisfy $j \geq
    2k + \frac{n}{2}$. Then we have
    \begin{equation}
        \label{eq:seminorm-estimate-for-terms-in-third-sum}
        \seminorm[K_\ell \times K_\ell, k]
        \left(
            \left(
                \chi \left(\frac{\eta}{\epsilon_j}\right)
                -
                \chi \left(\frac{\eta}{\epsilon_{j+1}}\right)
            \right)
            D(V^j) R^\pm_{U'}(2+2j,\argument)
        \right)
        \leq
        (\epsilon_j + \epsilon_{j+1}) f(k,\ell,j),
    \end{equation}
    with some constants $f(k,\ell,j)$ not depending on the choices of
    the $\epsilon_j$.
\end{lemma}
\begin{proof}
    We estimate
    \begin{align*}
        &\seminorm[K_\ell \times K_\ell, k]
        \left(
            \left(
                \chi \left(\frac{\eta}{\epsilon_j}\right)
                -
                \chi \left(\frac{\eta}{\epsilon_{j+1}}\right)
            \right)
            D(V^j) R^\pm_{U'}(2+2j,\argument)
        \right) \\
        &=
        \seminorm[K_\ell \times K_\ell, k]
        \left(
            \left(
                \chi \left(\frac{\eta}{\epsilon_j}\right)
                -
                \chi \left(\frac{\eta}{\epsilon_{j+1}}\right)
            \right)
            \eta^{k+1} D(V^j) c(2+2j,n) \eta^{j-\frac{n}{2}-k}
        \right) \\
        &\leq
        c_k c(2+2j,n) \\
        &\qquad
        \left(
            \seminorm[K_\ell \times K_\ell,k]
            \left(
                \chi \left(\frac{\eta}{\epsilon_j}\right)
                \eta^{k+1}
            \right)
            +
            \seminorm[K_\ell \times K_\ell,k]
            \left(
                \chi \left(\frac{\eta}{\epsilon_{j+1}}\right)
                \eta^{k+1}
            \right)
        \right)
        \seminorm[K_\ell \times K_\ell,k]
        \left(
            D(V^j) \eta^{j-\frac{n}{2}-k}
        \right) \\
        &\leq
        c_k c(2+2j,n)
        \left(
            \epsilon_j e(k,\ell,j)
            + \epsilon_{j+1} e(k,\ell,j)
        \right)
        \seminorm[K_\ell \times K_\ell,k]
        \left(
            D(V^j) \eta^{j-\frac{n}{2}-k}
        \right),
    \end{align*}
    with some constants $e(k,\ell,j)$ obtained from a Leibniz rule and
    arguments as in the proof of
    Lemma~\ref{lemma:seminorms-of-cutoff-times-lorentz-square} and
    Lemma~\ref{lemma:seminorms-of-cutoff-function}. Note that for $j
    \geq 2k - \frac{n}{2}$ the function $\eta^{j-k-\frac{n}{2}}$ is
    still $\Fun$ whence the last seminorm is still finite. Putting all
    the constants together, we get the desired estimate.
\end{proof}

Again, we can turn \eqref{eq:seminorm-estimate-for-terms-in-third-sum}
into a condition on the $\epsilon_j$ in order to make the seminorm
smaller than $\frac{1}{2^j}$.
\begin{lemma}
    \label{lemma:absolute-convergence-of-third-sum}
    Let the $\epsilon_j \in (0,1]$ be chosen such that in addition to
    \eqref{eq:choice-of-epsilon} we have
    \begin{equation}
        \label{eq:third-choice-of-epsilon}
        \epsilon_j \cdot
        \max \left\{
            \max_{\substack{\ell \leq j \\ 2k+\frac{n}{2} \leq j}}
            f(k,\ell,j),
            \max_{\substack{\ell \leq j-1 \\ 2k+\frac{n}{2} \leq j-1}}
            f(k,\ell,j-1)
        \right\}
        \leq \frac{1}{2^j}.
    \end{equation}
    Then the sum $\Sigma_3$ converges absolutely with respect to the
    $\Cinfty$-topology and yields a smooth section $\Sigma_3 \in
    \Secinfty\left(E^* \extensor E \at{U' \times U'}\right)$.
\end{lemma}
\begin{proof}
    Note that \eqref{eq:third-choice-of-epsilon} are again finitely
    many condition on each $\epsilon_j$ whence we indeed can find an
    $\epsilon_j \in (0,1]$ satisfying
    \eqref{eq:third-choice-of-epsilon}. Now
    Lemma~\ref{lemma:seminorm-estimate-for-terms-in-third-sum} yields
    the estimate
    \begin{align*}
        &\seminorm[K_\ell \times K_\ell, k]
        \sum_{j \geq N}
        \left(
            \left(
                \chi \left(\frac{\eta}{\epsilon_j}\right)
                -
                \chi \left(\frac{\eta}{\epsilon_{j+1}}\right)
            \right)
            D(V^j) R^\pm_{U'}(2+2j,\argument)
        \right) \\
        &\leq
        \seminorm[K_\ell \times K_\ell, k]
        \sum_{j=N}^{j_0-1}
        \left(
            \left(
                \chi \left(\frac{\eta}{\epsilon_j}\right)
                -
                \chi \left(\frac{\eta}{\epsilon_{j+1}}\right)
            \right)
            D(V^j) R^\pm_{U'}(2+2j,\argument)
        \right)
        + \sum_{j=j_0}^\infty (\epsilon_j + \epsilon_{j+1})
        f(k,\ell,j) \\
        &\leq
        \textrm{const.} + 2 \sum_{j=j_0}^\infty \frac{1}{2^j}
        < \infty,
    \end{align*}
    if we take $j_0 \geq N$ such that $j_0 \geq 2k + \frac{n}{2}$ and
    $j_0 \geq \ell$. Indeed, in this case we have
    \[
    \epsilon_j \cdot
    \max_{\substack{\ell \leq j \\
        2k+\frac{n}{2} \leq j}}
    f(k,\ell,j)
    \leq \frac{1}{2^j}
    \quad
    \textrm{and}
    \quad
    \epsilon_{j+1} \cdot
    \max_{\substack{\ell \leq j \\
        2k+\frac{n}{2} \leq j}}
    f(k,\ell,j)
    \leq \frac{1}{2^j},
    \]
    both by \eqref{eq:third-choice-of-epsilon}. But then the absolute
    convergence of $\Sigma_3$ is clear as the finitely many terms $N
    \leq j \leq j_0-1$ do not change the convergence.
\end{proof}

Collecting the results of the previous lemmas we arrive at the
following result:
\begin{proposition}
    \label{proposition:D-on-approx-solution}
    There is a choice of $\epsilon_j \in (0,1]$ such that the
    approximate solution $\widetilde{\mathcal{R}}^\pm_{U'}(p)$
    satisfies in addition to the properties described in
    Proposition~\ref{proposition:convergence-of-tail-series} and
    Corollary~\ref{corollary:order-and-supps-of-approx-fund-solution}
    \begin{equation}
        \label{eq:D-on-approx-solution-is}
        D \widetilde{\mathcal{R}}^\pm_{U'}(p)
        = \delta_p + K^\pm_{U'}(p,\argument)
    \end{equation}
    with some smooth section $K^\pm_{U'} \in \Secinfty\left(E^*
        \extensor E \at{U' \times U'}\right)$ for $p \in U'$.
\end{proposition}
\begin{proof}
    Indeed, the section $K^\pm_{U'}$ is obtained from the computation
    in \eqref{eq:D-on-approx-is-delta-plus-this-right-here} as
    \[
    K^\pm_{U'}
    = \left(1 - \chi \left(\frac{\eta}{\epsilon_N}\right)\right)
    D(V^{N-1}) R^\pm_{U'}(2N,\argument)
    + \Sigma_1 + \Sigma_2 + \Sigma_3.
    \]
    The convergence results on the series $\Sigma_1, \Sigma_2$, and
    $\Sigma_3$ yield $K^\pm_{U'} \in \Secinfty\left(E^* \extensor E
        \at{U' \times U'}\right)$ as we wanted. Note that in total, we
    only have to impose finitely many conditions on each $\epsilon_j$
    according to
    Proposition~\ref{proposition:convergence-of-tail-series},
    \refitem{item:choice-of-epsilon},
    Lemma~\ref{lemma:absolute-convergence-of-second-sum}, the analogue
    condition from $\Sigma_1$, and
    Lemma~\ref{lemma:absolute-convergence-of-third-sum}.
\end{proof}
\begin{remark}[Parametrix]
    \label{remark:parametrix}
    \index{Parametrix}%
    The proposition just says that we have constructed a parametrix
    $\widetilde{\mathcal{R}}^\pm_{U'}(p)$ of $D$ for every $p \in U'$,
    see also \cite[Sect.~7.1]{hoermander:2003a} for more information
    on parametrices.
\end{remark}

In
Proposition~\ref{proposition:riesz-dependence-on-base-point-p-prime}
we had some estimates for $|R^\pm_{U'}(p)(\varphi)|$ \emph{locally
  uniform} in $p$. Since $\widetilde{\mathcal{R}}^\pm_{U'}(p)$ is
build out of the $R^\pm_{U'}(\alpha,p)$ we can expect a similar
feature also for $\widetilde{\mathcal{R}}^\pm_{U'}(p)$. Indeed, this
is the case:

For a fixed $\varphi \in \Secinfty_0(E^* \at{U'})$ we can view $U' \ni
p \mapsto \widetilde{\mathcal{R}}^\pm_{U'}(p)(\varphi) \in E_p^*$ as a
section of $E^*$ defined on $U'$. This section has nice features, it
will be smooth again. More precisely, we have the following
statements:
\begin{proposition}
    \label{proposition:estimates-and-cont-of-pairing-with-approx-solution}
    Let $\widetilde{\mathcal{R}}^\pm_{U'}(p)$ be the approximate
    fundamental solution. Moreover, let $k \in \mathbb{N}_0$ and $K, L
    \subset U'$ be compact. Then we have:
    \begin{propositionlist}
    \item \label{item:estimates-and-order-of-approx-solution} There is
        a constant $c_{K,L} > 0$ such that
        \begin{equation}
            \label{eq:cont-estimate-for-approx-solution}
            \left|
                \widetilde{\mathcal{R}}^\pm_{U'}(p)(\varphi)
            \right|
            \leq
            c_{K,L} \seminorm[K,n+1](\varphi)
        \end{equation}
        for all $p \in L$ and $\varphi \in
        \Secinfty_K\left(E^*\at{U'}\right)$. In particular, the
        distribution $\widetilde{\mathcal{R}}^\pm_{U'}(p)$ is of
        global order $\leq n+1$.
    \item \label{item:aprox-fundamental-paired-is-smooth} The section
        $\widetilde{\mathcal{R}}^\pm_{U'}(\argument)(\varphi)$ of
        $E^*\at{U'}$ is smooth for all $\varphi \in
        \Secinfty_0(E^*\at{U'})$.
    \item \label{item:cont-estimate-for-pairing-with-approx-solution}
        There are constants $c_{K,L,k} > 0$ such that
        \begin{equation}
            \label{eq:cont-estimate-for-pairing-with-approx-solution}
            \seminorm[L,k](\widetilde{\mathcal{R}}^\pm_{U'}(\argument)(\varphi))
            \leq c_{K,L,k} \seminorm[K,k+n+1](\varphi)
        \end{equation}
        for all $\varphi \in \Secinfty_K(E^*\at{U})$.
    \item \label{item:pairing-operator-is-cont} The operator
        \begin{equation}
            \label{eq:pairing-operator}
            \mathcal{R}^\pm_{U'}:
            \Secinfty_0(E^*\at{U'}) \ni \varphi
            \; \mapsto \;
            \left(
                p \; \mapsto \; \mathcal{R}^\pm_{U'}(p)(\varphi)
            \right)
            \in \Secinfty(E^*\at{U'})
        \end{equation}
        is continuous in the $\Cinfty_0$- and $\Cinfty$-topology.
    \end{propositionlist}
\end{proposition}
\begin{proof}
    Clearly, the estimate \eqref{eq:cont-estimate-for-approx-solution}
    is a particular case of the more general situation in
    \eqref{eq:cont-estimate-for-pairing-with-approx-solution} for
    $k=0$. Thus fix $k \in \mathbb{N}_0$. Then we have
    \begin{align*}
        \widetilde{\mathcal{R}}^\pm_{U'}(p)
        = \sum_{j=0}^{N+1} V_p^j R^\pm_{U'}(2+2j,p)
        + \sum_{j=N}^{N+k-1} \chi \left(\frac{\eta_p}{\epsilon_j}\right)
        V_p^j R^\pm_{U'}(2+2j,p)
        + \sum_{j=N+k}^\infty \chi \left(\frac{\eta_p}{\epsilon_j}\right)
        V_p^j R^\pm_{U'}(2+2j,p),
        \tag{$*$}
    \end{align*}
    and we know that the third contribution converges in the
    $\Fun$-topology to
    \[
    f_k(p,q)
    = \sum_{j=N+k}^\infty
    \chi \left(\frac{\eta_p}{\epsilon_j}\right)(q)
    V_p^j(q) R^\pm_{U'}(2+2j,p)(q),
    \]
    which is a section $f_k \in \Fun(E^* \extensor E \at{U' \times
      U'})$. Now let $\varphi \in \Secinfty_K(E^*\at{U'})$ then the
    pairing of $f_k$ with $\varphi$ is
    \begin{align*}
        p \; \mapsto \; f_k(p, \argument) \varphi
        = \int_{U'}f(p,q) \cdot \varphi(q) \: \mu_g(q)
        = \int_K f_k(p,q) \cdot \varphi(q) \: \mu_g(q),
        \tag{$**$}
    \end{align*}
    which still yields a $\Fun$-section. In fact, we immediately
    obtain an estimate of the form
    \[
    \seminorm[L,k](f_k(\argument)\varphi)
    \leq \vol(K) \seminorm[L \times K,k](f_k) \seminorm[K,0](\varphi)
    \]
    by differentiating into the integral ($**$), which is legal as the
    compactly supported integrand is $\Fun$ in $p$ and all first
    derivatives in $p$-direction yield still a continuous integrand in
    $p$ and $q$. The first and second contribution in ($*$) are
    slightly more complicated. First we note that the sums are all
    finite and each term is of the form $\Phi^k(p,\argument)
    R^\pm_{U'}(2+2j,p)$ with a smooth section $\Phi^j \in
    \Secinfty(E^* \extensor E \at{U' \times U'})$. Thus applying this
    to a fixed test section $\varphi \in \Secinfty_K(E^*\at{U'})$
    gives by the very definition of the Riesz distributions the map
    \begin{align*}
        p \; \mapsto \; \Phi^j(p, \argument) R^\pm_{U'}(2+2j,p)(\varphi)
        &=
        R^\pm_{U'}(2+2j,p)
        \left(\Phi^j(p,\argument) \varphi(\argument)\right)
        \\
        &= R^\pm(2+2j) \left(
            \widetilde{\varrho}_p(\argument)
            \exp_p^*( \Phi^j(p,\argument) \varphi(\argument))
        \right).
        \tag{$*{*}*$}
    \end{align*}
    If we want now to estimate the $p$-dependence we can rely on
    Lemma~\ref{lemma:parameter-differentiation-under-integral}: The
    function $(p,q) \mapsto \widetilde{\varrho}_p(q) \exp_p^*(
    \Phi^j(p,q) \varphi(q))$ is smooth in both variables and has
    support in $U' \times K$ thanks to the support condition on
    $\varphi$. Thus the lemma applies and yields a smooth function of
    $p$. Moreover, we can differentiate into the application of
    $R^\pm(2+2j)$ and have for the $p$-derivatives of ($*{*}*$)
    \begin{align*}
        &\frac{\partial^{|I|}}{\partial x^I}
        \left(
            p \; \mapsto \;
            \Phi^j(p, \argument) R^\pm_{U'}(2+2j,p)(\varphi)
        \right) \\
        &\quad= R^\pm(2+2j)
        \left(
            \frac{\partial^{|I|}}{\partial x^I}
            \left(
                p \; \mapsto \;
                \Phi^j(p, \argument) R^\pm_{U'}(2+2j,p)(\varphi)
            \right)
        \right),
        \tag{\smiley}
    \end{align*}
    where $x$ are some generic coordinates for the $p$-variable. Now
    we know that for $j \geq 0$ the Riesz distribution $R^\pm(2+2j)$
    is of order $\leq n+1$. In fact, the order is much less for some
    $j$, see also the low dimensional discussion in
    Section~\ref{subsec:riesz-distribution-for-1dim-and-2dim}, but the
    above estimate on the order will do the job. Thus for each term we
    get an estimate of the form
    \[
    \seminorm[L,k]\left(\Phi^k R^\pm_U(\argument) \varphi\right)
    \leq c^j_{K,L} \seminorm[K,k+n+1](\varphi),
    \]
    as we need the $n+1$ derivatives of $\varphi$ for $R^\pm(2+2j)$
    and up to $k$ derivatives from the differentiation and the chain
    rule coming from (\smiley). In the constant $c^j_{K,L}$ we get
    contributions of the first $k$ derivatives of $\Phi^j$, $\exp_p$
    and $\widetilde{\varrho}_p$ as well as from the continuity of
    $R^\pm(2+2j)$. Thus we arrive at finitely many estimates for the
    finitely many terms in ($*$) which can be combined into
    \eqref{eq:cont-estimate-for-pairing-with-approx-solution}. This
    shows the third part. But then the fourth part is clear as well.
\end{proof}

\begin{remark}
    \label{remark:extension-of-approx-pairing-to-nonsmooth}
    The estimate in
    \eqref{eq:cont-estimate-for-pairing-with-approx-solution} also
    shows that we can apply the operator
    $\widetilde{\mathcal{R}}^\pm_{U'}$ to less regular sections than
    smooth ones. In fact, $\widetilde{\mathcal{R}}^\pm_{U'}$ extends
    to a well-defined continuous linear operator
    \begin{equation}
        \label{eq:extension-of-approx-pairing-to-nonsmooth}
        \widetilde{\mathcal{R}}^\pm_{U'}:
        \Sec[k+n+1]_0(E^*\at{U'})
        \longrightarrow \Sec(E^*\at{U'})
    \end{equation}
    for all $k \geq 0$ with respect to the $\Fun[k+n+1]_0$- and
    $\Fun$-topology, respectively. This will sometimes be a useful
    extension.
\end{remark}

The last features we will need are some support properties of the
``defect'' $K^\pm_{U'}$ of $\widetilde{\mathcal{R}}^\pm_{U'}$ being a
fundamental solution.
\begin{lemma}
    \label{lemma:support-feature-of-defect}
    The smooth section $K^\pm_{U'} \in \Secinfty\left(E^* \extensor E
        \at{U' \times U'}\right)$ satisfies
    \begin{equation}
        \label{eq:support-of-defect-is-ligthcone-like}
        (p,q) \in \supp K^\pm_{U'} \subseteq U' \times U'
        \Longrightarrow q \in J^\pm_{U'}(p).
    \end{equation}
\end{lemma}
\begin{proof}
    Assume that $K^\pm_{U'}(p,q)$ is non-zero. From
    \[
    K^\pm_{U'}(p,q)
    = \left(
        1 - \chi \left(\frac{\eta_p(q)}{\epsilon_j}\right)
    \right)
    D \left(
        V_p^{N-1} R^\pm_{U'}(2N,p)
    \right)(q)
    + \Sigma_1(p,q)
    + \Sigma_2(p,q)
    + \Sigma_3(p,q)
    \]
    and the fact that each series $\Sigma_1, \Sigma_2, \Sigma_3$ has
    only terms involving $R^\pm_{U'}(2+2j,p)(q)$, to have a non-zero
    contribution we necessarily need $q \in J^\pm_{U'}(p)$. Thus
    $K^\pm_{U'}(p,q) \neq 0$ implies $q \in J^\pm_{U'}(p)$. Since the
    support of $K^\pm_{U'}$ is the closure of all those point with
    $K^\pm_{U'}(p,q) \neq 0$ it is contained in the closure of those
    points $(p,q) \in U' \times U'$ with $q \in J^\pm_{U'}(p)$, all
    closures taken with respect to $U' \times U'$. Since $U'$ is
    assumed to be geodesically convex, one can show that the causal
    relation\index{Causal relation}
    \[
    J^\pm_{U'}
    = \left\{
        (p,q) \in U' \times U'
        \; \big| \;
        q \in J^\pm_{U'}(p)
    \right\}
    \subseteq U' \times U'
    \]
    is actually closed. Note that this is a stronger statement than
    all $J^\pm_{U'}(p)$ being closed in $U'$, see e.g.
    \cite[Prop.~2.10]{minguzzi.sanchez:2006a:pre} or \cite[Lemma~2 in
    Chap.~14]{oneill:1983a}. But then
    \eqref{eq:support-of-defect-is-ligthcone-like} follows at once.
\end{proof}

\begin{remark}[Future and past stretched subsets]
    \label{remark:future-and-past-stretched}
    \index{Future stretched}%
    \index{Past stretched}%
    A subset $S \subseteq U' \times U'$ with the feature that $(p,q)
    \in S$ implies $q \in J^\pm_{U'}(p)$ is also called future or past
    stretched, respectively. Thus the support of $K^\pm_{U'}$ is
    future and past stretched with respect to $U'$, respectively.
\end{remark}

We are now in the position to collect all the features of the
approximate fundamental solution $\widetilde{\mathcal{R}}^\pm_{U'}$ we
shall need in the following:
\begin{theorem}[Approximate fundamental solution]
    \label{theorem:approx-fundamental-solution}
    \index{Fundamental solution!approximate}%
    \index{Hadamard coefficients}%
    Let $U' \subseteq M$ be geodesically convex and let $V^j \in
    \Secinfty\left(E^* \extensor E \at{U' \times U'}\right)$ be the
    Hadamard coefficients with respect to the normally hyperbolic
    operator $D \in \Diffop^2(E)$. Then there exists a sequence
    $\epsilon_j \in (0,1]$ for $j \geq N > \frac{n}{2}$ such that
    \begin{equation}
        \label{eq:the-approximate-solution}
        \widetilde{\mathcal{R}}^\pm_{U'}(p)
        = \sum_{j=0}^{N-1} V_p^j R^\pm_{U'}(2+2j,p)
        + \sum_{j=N}^\infty
        \chi \left( \frac{\eta_p}{\epsilon_j} \right)
        V_p^j R^\pm_{U'}(2+2j,p)
    \end{equation}
    converges in the weak$^*$ topology to a distribution
    $\widetilde{\mathcal{R}}^\pm_{U'}(p) \in \Sec[-\infty]\left(E
        \at{U'}\right) \tensor E^*_p$ with the following properties:
    \begin{theoremlist}
    \item \label{item:supp-and-singsupp-of-approx-solution} For the
        support and singular support we have
        \begin{equation}
            \label{eq:supp-of-approx-solution}
            \supp \widetilde{\mathcal{R}}^\pm_{U'}(p) \subseteq
            J^\pm_{U'}(p),
        \end{equation}
        \begin{equation}
            \label{eq:singsupp-of-approx-solution}
            \singsupp \widetilde{\mathcal{R}}^\pm_{U'}(p) \subseteq
            C^\pm_{U'}(p).
        \end{equation}
    \item \label{item:defect-of-approx-solution} We have
        \begin{equation}
            \label{eq:defect-of-approx-solution}
            D \widetilde{\mathcal{R}}^\pm_{U'}(p)
            = \delta_p + K^\pm_{U'}(p, \argument)
        \end{equation}
        with a smooth section $K^\pm_{U'} \in \Secinfty\left(E^*
            \extensor E \at{U' \times U'}\right)$.
    \item \label{item:supp-of-defect} The support of $K^+_{U'}$ is
        future stretched and the support of $K^-_{U'}$ is past
        stretched.
    \item \label{item:approx-solution-smooth-in-p} For a test section
        $\varphi \in \Secinfty_0\left(E^* \at{U'}\right)$ the section
        $p \mapsto \widetilde{\mathcal{R}}^\pm_{U'}(p)(\varphi)$ is
        smooth.
    \item \label{item:uniform-cont-est-and-order-of-approx-solution}
        For compact subsets $K,L \subseteq U'$ there exist constants
        $c_{K,L} > 0$ such that
        \begin{equation}
            \label{eq:uniform-cont-estimate-of-approx-solution}
            \left|
                \widetilde{\mathcal{R}}^\pm_{U'}(p) (\varphi)
            \right|
            \leq c_{K,L} \seminorm[K,n+1](\varphi)
        \end{equation}
        for all $p \in L$ and $\varphi \in \Secinfty_K\left(E^*
            \at{U'}\right)$. In particular, for the global order of
        $\widetilde{\mathcal{R}}^\pm_{U'}(p)$ we have
        \begin{equation}
            \label{eq:order-of-the-approx-solution}
            \ord\left(\widetilde{\mathcal{R}}^\pm_{U'}(p)\right)
            \leq
            n+1.
        \end{equation}
    \end{theoremlist}
\end{theorem}

%
%

\subsection{Construction of the Local Fundamental Solution}
\label{satz:construction-local-solution}

Having a (well-behaved) parametrix to a differential operator there is
a more or less standard procedure of how one can obtain a fundamental
solution from it. Roughly speaking, the defect in having a fundamental
solution is so small that one can use a geometric series to resolve
this problem.

We will choose now an open subset $U \subseteq U'$ such that
\begin{equation}
    \label{eq:choose-of-precompact-subset}
    U^\cl \subseteq U'
\end{equation}
is compact. Later on, we will need additional properties of $U$ but
for the time being the compactness of $U^\cl$ will suffice. Then we
consider the following integral operator build out of $K^\pm_{U'} \in
\Secinfty\left(E^* \extensor E \at{U' \times U'}\right)$. Let
$\varphi$ be a section of $E^*$ defined at least on $U^\cl$ then we
can naturally pair $K^\pm_{U'}(p,q) \cdot \varphi(q)$ and integrate.
This gives
\begin{equation}
    \label{eq:the-integral-operator}
    \left(
        \mathcal{K}^\pm_{U} \varphi
    \right)(p)
    = \int_{U^\cl} K^\pm_{U'}(p,q) \cdot \varphi(q) \: \mu_g(q).
\end{equation}
Depending on the properties of $\varphi$ the integral will be
well-defined and yields a rather nice section of $E^*$ defined on
$U'$. One rather general scenario is the following:
\begin{definition}
    \label{definition:measurable-bounded-sections}
    \index{Bounded section}%
    With respect to some auxiliary positive fiber metric on $E^*$ we
    define
    \begin{equation}
        \label{eq:measurable-bounded-sections}
        \Gamma_b\left(E^*\at{U}\right)
        = \left\{
            \varphi: U \longrightarrow E^*
            \; \big| \;
            \varphi(q) \in E^*_q
            \; \textrm{and} \;
            \varphi
            \; \textrm{is bounded and measurable}
        \right\}.
    \end{equation}
\end{definition}
Here the fiber metric is used to define a norm on each fiber. With
respect to these norms we want $\varphi$ to be bounded over $U$.
The following technical lemma is well-known and obtained in a
completely standard way:
\begin{lemma}[The Banach space $\Gamma_b\left(E^*\at{U}\right)$]
    \label{lemma:measurable-bounded-sections}
    Let $U \subseteq M$ be open with compact closure.
    \begin{definitionlist}
    \item \label{item:measurable-bounded-welldefined} The definition
        of $\Gamma_b\left(E^*\at{U}\right)$ does not depend on the
        auxiliary smooth fiber metric.
    \item \label{item:measurable-bounded-banach-space} The vector
        space $\Gamma_b\left(E^*\at{U}\right)$ becomes a Banach space
        via the norm
        \begin{equation}
            \label{eq:banach-norm-on-measurable-bounded}
            \seminorm[U,0](\varphi)
            = \sup_{q \in U} \norm{\varphi(q)}_{E^*_q}.
        \end{equation}
    \item \label{item:different-norms} Different choices of positive
        fiber metrics on $E^*$ yield equivalent Banach norms
        \eqref{eq:banach-norm-on-measurable-bounded}.
    \item \label{item:Gammak-in-meas-bounded} The restriction map
        \begin{equation}
            \label{eq:Gammak-in-measu-bounded}
            \Sec(E^*) \ni \varphi \; \mapsto \;
            \varphi \at{U} \in \Gamma_b(E^*\at{U})
        \end{equation}
        is continuous for all $k \in \mathbb{N}_0 \cup \{+\infty\}$.
    \end{definitionlist}
\end{lemma}
\begin{proof}
    The measurability of a section is intrinsically defined and refers
    only to the Borel $\sigma$-algebra of the topological space $M$.
    Clearly, the boundedness does not depend on the choice of the
    fiber metric. Only the numerical value of the bound depends on
    this choice. Obviously,
    \eqref{eq:banach-norm-on-measurable-bounded} is a norm and
    different choices of the fiber metric yield equivalent norms in
    \eqref{eq:banach-norm-on-measurable-bounded}. This can entirely be
    copied from our considerations in
    Theorem~\ref{theorem:Ck-Topologie}. We have to show completeness
    of $\Gamma_b(E^*\at{U})$. Thus let $\varphi_n \in
    \Gamma_b(E^*\at{U})$ be a Cauchy sequence with respect to
    $\seminorm[U,0]$. Then we have uniform convergence of
    $\varphi_n(q) \longrightarrow \varphi(q)$ on $U^\cl$. Since every
    $\varphi_n$ is bounded the limit is bounded as well. Finally,
    already the pointwise limit of measurable functions (and hence by
    local triviality: of sections) is known to be measurable again,
    see e.g. \cite[Satz X.1.11]{amann.escher:2001a}. Thus $\varphi \in
    \Gamma_b\left(E^*\at{U}\right)$ is the desired limit of
    $\varphi_n$.  Finally, if $\varphi \in \Sec(E^*)$ then
    $\varphi\at{U} \in \Gamma_b\left(E^* \at{U}\right)$ since over a
    compactum $U^\cl$ any continuous section is bounded and
    measurable. Moreover, by elementary features of the supremum we
    have
    \[
    \seminorm[U,0] \left(\varphi\at{U}\right)
    = \seminorm[U^\cl,0] (\varphi),
    \]
    with our previous definition of the seminorm $\seminorm[K,\ell]$.
    This gives the continuity of \eqref{eq:Gammak-in-measu-bounded}.
\end{proof}

We claim that the operator $\mathcal{K}^\pm_{U}$ is well-defined on
$\Gamma_b\left(E^*\at{U}\right)$ and maps into the smooth sections in
a continuous manner.
\begin{lemma}
    \label{lemma:integral-operator}
    Let $k \in \mathbb{N}_0$ and $U \subseteq U'$ open with compact
    closure $U^\cl \subseteq U'$.
    \begin{lemmalist}
    \item \label{item:integral-operator-makes-smooth} For $\varphi \in
        \Gamma_b\left(E^*\at{U}\right)$ we have $\mathcal{K}^\pm_{U}
        \varphi \in \Secinfty\left(E^*\at{U'}\right)$.
    \item \label{item:integral-operator-is-cont} We have an estimate
        of the form
        \begin{equation}
            \label{eq:cont-estimate-of-integral-operator}
            \seminorm[K,k] \left(\mathcal{K}^\pm_{U} \varphi\right)
            \leq
            \vol(U^\cl) \seminorm[K \times U^\cl,k](K^\pm_{U'})
            \seminorm[U,0](\varphi)
        \end{equation}
        for all $\varphi \in \Gamma_b\left(E^*\at{U}\right)$ and
        compact $K \subseteq U'$.
    \end{lemmalist}
\end{lemma}
\begin{proof}
    We first proof continuity. Thus let $p \in U'$ be fixed and
    consider $p_n \longrightarrow p$. Since the integrand
    $K^\pm_{U}(p_n,q) \cdot \varphi(q)$ is bounded by some integrable
    function, namely by the constant function $\seminorm[K \times
    U^\cl, 0](K^\pm_{U'}) \seminorm[U,0](\varphi)$ where $K$ is any
    compactum containing the convergent sequence $p_n$, we can apply
    Lebesgue's dominated convergence and find
    \begin{align*}
        \lim\limits_{n \rightarrow \infty}
        (\mathcal{K}^\pm_{U} \varphi)(p_n)
        &=
        \lim\limits_{n \rightarrow \infty} \int_{U^\cl}
        K^\pm_{U'}(p_n,q) \cdot \varphi(q)
        \: \mu_g(q) \\
        &\stackrel{\mathclap{\textrm{Lebesgue}}}{=}
        \quad
        \int_{U^\cl} \lim\limits_{n \rightarrow \infty}
        K^\pm_{U'}(p_n,q) \cdot \varphi(q)
        \: \mu_g(q) \\
        &=
        \int_{U^\cl} K^\pm_{U'}(p,q) \cdot \varphi(q)
        \: \mu_g(q) \\
        &= (\mathcal{K}^\pm_{U} \varphi)(p),
    \end{align*}
    which is the continuity of $\mathcal{K}^\pm_{U} \varphi$. By an
    analogous argument we can also exchange the partial
    differentiation with the integration yielding a continuous partial
    derivative
    \[
    \frac{\partial}{\partial x^i} \mathcal{K}^\pm_{U} \varphi
    =
    \int_{U^\cl}
    \frac{\partial K^\pm_{U'}(p,q)}{\partial x^i} \cdot \varphi(q)
    \: \mu_g(q),
    \tag{$*$}
    \]
    all with respect to some local trivialization of $E^*$. Thus
    $\mathcal{K}^\pm_{U} \varphi$ turns out to be $\Fun[1]$ and by
    induction we get $\mathcal{K}^\pm_{U} \varphi \in
    \Secinfty\left(E^*\at{U'}\right)$. This shows the first part. For
    the second, we use a local trivialization and ($*$) to obtain
    \[
    \frac{\partial^{|I|}}{\partial x^I}
    ( \mathcal{K}^\pm_{U} \varphi) \At{p}
    =
    \int_{U^\cl} \frac{\partial^{|I|} K^\pm_{U'}}{\partial x^I}(p,q)
    \cdot \varphi(q)
    \: \mu_g(q),
    \]
    from which we get
    \begin{align*}
        \seminorm[U^\cl, k](\mathcal{K}^\pm_{U} \varphi)
        &\leq
        \sup_{\substack{p \in U^\cl \\ |I| \leq k}}
        \int_{U^\cl}
        \norm{\frac{\partial^{|I|} K^\pm_{U'}}{\partial x^I}(p,q)}
        \norm{\varphi(q)}
        \: \mu_g(q) \\
        &\leq
        \vol(U^\cl)
        \seminorm[K \times U^\cl, k](K^\pm_{U'})
        \seminorm[U,0](\varphi).
    \end{align*}
\end{proof}

With other words, the integral operator behaves like a convolution
integral: the result inherits the better properties concerning
smoothness of both factors under the integral.

The problem is now that the operator $\mathcal{K}^\pm_{U}$ is far from
being ``local'': it changes and typically enlarges the support
strictly. Thus it is slightly tricky to define powers of
$\mathcal{K}^\pm_{U}$. However, as we did not insist on $\varphi$
being continuous at all we can proceed as follows: For $\varphi \in
\Gamma_b\left(E^* \at{U}\right)$ the section $\mathcal{K}^\pm_{U}
\varphi$ is smooth and defined on the \emph{larger} open subset $U'$.
Thus restricting $\mathcal{K}^\pm_{U} \varphi$ back to $U^\cl$ yields
a section which is clearly measurable and bounded \emph{and} still
smooth on the interior $U$ of $U^\cl$. Thus we have
\begin{equation}
    \label{eq:restriction-of-integral-op-result}
    \Gamma_b\left(E^*\at{U}\right) \ni \varphi
    \; \mapsto \;
    \mathcal{K}^\pm_{U} \varphi \at{U} \in
    \Gamma_b\left(E^*\at{U}\right).
\end{equation}
By some slight abuse of notation we denote the composition $\varphi
\mapsto \mathcal{K}^\pm_{U} \varphi \mapsto \mathcal{K}^\pm_{U}
\varphi \at{U}$ again simply by $\mathcal{K}^\pm_{U}$.
\begin{lemma}
    \label{lemma:restricted-integral-operator}
    \index{Operator norm}%
    The linear operator
    \begin{equation}
        \label{eq:restricted-integral-operator}
        \mathcal{K}^\pm_{U}: \Gamma_b\left(E^*\at{U}\right)
        \ni
        \varphi \; \mapsto \; \mathcal{K}^\pm_{U} \varphi \at{U}
        \in
        \Gamma_b\left(E^*\at{U}\right)
    \end{equation}
    is continuous with operator norm
    \begin{equation}
        \label{eq:operator-norm-of-restricted-integral-operator}
        \norm{\mathcal{K}^\pm_{U}}
        \leq \vol(U^\cl) \seminorm[U^\cl \times U^\cl,0](K^\pm_{U'}).
    \end{equation}
\end{lemma}
\begin{proof}
    From Lemma~\ref{lemma:integral-operator} we know that for all
    $\varphi \in \Gamma_b\left(E^*\at{U}\right)$ we have
    \[
    \seminorm[U,0](\mathcal{K}^\pm_{U} \varphi)
    = \seminorm[U^\cl,0](\mathcal{K}^\pm_{U} \varphi)
    \leq \vol(U^\cl)
    \seminorm[U^\cl \times U^\cl,0](\mathcal{K}^\pm_{U})
    \seminorm[U,0](\varphi),
    \]
    which gives the continuity as well as the estimate on the operator
    norm \eqref{eq:operator-norm-of-restricted-integral-operator}.
\end{proof}
\begin{corollary}
    \label{corollary:geometric-series-as-cont-inverse-of-id-plus-intop}
    If the open subset $U \subseteq U'$ is sufficiently small in the
    sense that
    \begin{equation}
        \label{eq:small-enough-subset}
        \vol(U^\cl)
        \seminorm[U^\cl \times U^\cl,0](\mathcal{K}^\pm_{U})
        < 1,
    \end{equation}
    then the operator
    \begin{equation}
        \label{eq:id-plus-int-operator}
        \id + \mathcal{K}^\pm_{U}:
        \Gamma_b\left(E^*\at{U}\right) \longrightarrow
        \Gamma_b\left(E^*\at{U}\right)
    \end{equation}
    is invertible with continuous inverse given by the absolutely
    norm-convergent geometric series
    \begin{equation}
        \label{eq:geometric-series-as-cont-inverse-of-id-plus-intop}
        \index{Geometric series}%
        \left(
            \id + \mathcal{K}^\pm_{U}
        \right)^{-1}
        = \sum_{j=0}^\infty (-\mathcal{K}^\pm_{U})^j.
    \end{equation}
\end{corollary}
\begin{proof}
    Since the operator norm of $\mathcal{K}^\pm_{U}$ is smaller or
    equal to $\vol(U^\cl) \seminorm[U^\cl \times
    U^\cl,0](\mathcal{K}^\pm_{U})$ the statement follows from general
    arguments on the geometric series and the fact that bounded
    operators on a Banach space form a Banach space themselves with
    respect to the operator norm.
\end{proof}

Note that since $\seminorm[U^\cl \times U^\cl,0](\mathcal{K}^\pm_{U})$
is only getting smaller for smaller $U^\cl$, there always exists a
small enough $U \subseteq U'$ around a given point in $U'$.

The idea is now to use the inverse $(\id + \mathcal{K}^\pm_{U})^{-1}$
to correct the approximate solution $\widetilde{\mathcal{R}}^\pm_{U'}$
at least on some small enough $U \subseteq U'$. There are now two
problems: the inverse a priori maps into
$\Gamma_b\left(E^*\at{U}\right)$ but we want some smooth section
instead of a bounded and measurable one. Moreover, we want to control
the support of the result at least in so far that we get ``causal
behaviour''.

The first problem is solved by a more careful investigation of the
geometric series: indeed the operator $\mathcal{K}^\pm_{U}$ already
maps into much nicer sections than just bounded and measurable ones.
By Lemma~\ref{lemma:integral-operator} they are restrictions of smooth
sections on $U'$.

The second problem will persist unless we make some additional
assumptions on the subset $U$. It has to be causal, see
Section~\ref{subsec:caus-cond-cauchy-hypersurfaces}. We will postpone
this investigation to
Section~\ref{satz:causal-props-of-fund-solution}.

We start to discuss the smoothness properties. For continuous sections
things are still very simple as there is a good and easy notion of a
continuous section over a compact subset. In fact, the continuous
sections over $U^\cl$ form a \emph{closed subspace}
\begin{equation}
    \label{eq:continous-closed-in-meas-bounded}
    \Sec[0]\left(E^*\at{U^\cl}\right)
    \subseteq \Gamma_b\left(E^*\at{U}\right)
\end{equation}
with respect to the norm $\seminorm[U,0] = \seminorm[U^\cl,0]$.
Clearly, restricting a continuous section $\varphi \in
\Sec[0]\left(E^*\at{U'}\right)$ to $U^\cl$ yields $\varphi \at{U^\cl}
\in \Sec[0]\left(E^*\at{U^\cl}\right)$. From
Lemma~\ref{lemma:integral-operator},
\refitem{item:integral-operator-makes-smooth} we obtain
\begin{equation}
    \label{eq:intop-on-cont-on-U-closure}
    \mathcal{K}^\pm_{U}:
    \Sec[0]\left(E^*\at{U^\cl}\right) \longrightarrow
    \Sec[0]\left(E^*\at{U^\cl}\right)
\end{equation}
in a continuous way. Moreover, the operator norm estimate
\eqref{eq:operator-norm-of-restricted-integral-operator} for the
restriction \eqref{eq:intop-on-cont-on-U-closure} of
$\mathcal{K}^\pm_{U}$ to continuous sections is still valid. Since
$\Sec[0]\left(E^*\at{U^\cl}\right)$ is a Banach space by its own, we
get a continuous invertible operator
\begin{equation}
    \label{eq:id-plus-intop-on-cont-sect-on-U-cl-is-cont-and-inv}
    \left( \id + \mathcal{K}^\pm_{U} \right)^{-1}
    = \sum_{j=0}^\infty (-\mathcal{K}^\pm_{U})^j:
    \Sec[0]\left(E^*\at{U^\cl}\right)
    \longrightarrow \Sec[0]\left(E^*\at{U^\cl}\right)
\end{equation}
with absolutely norm-convergent geometric series analogously to
Corollary~\ref{corollary:geometric-series-as-cont-inverse-of-id-plus-intop}.

In order to control the smoothness properties of the inverse of $\id +
\mathcal{K}^\pm_{U}$ we introduce the following subspaces of
$\Sec[0]\left(E^*\at{U^\cl}\right)$. The tricky point is to define
smoothness on a \emph{closed} subset $U^\cl$ instead of an open one in
such a way that we still get a good functional space.
\begin{definition}[The space $\Sec\left(E^*\at{U^\cl}\right)$]
    \label{definition:differentiable-on-closed-set}
    Let $k \in \mathbb{N}_0$, then a section $\varphi \in
    \Sec[0]\left(E^*\at{U^\cl}\right)$ is called $\Fun$ on $U^\cl$ if
    it can be approximated by sections $\varphi_n\at{U^\cl}$, with
    $\varphi_n \in \Sec\left(E^*\at{U_n}\right)$ with respect to the
    norm $\seminorm[U^\cl,k]$, where $U_n \supseteq U^\cl$ is open.
    The set of all such section is denoted by
    \begin{equation}
        \label{eq:differentiable-on-closed-set}
        \Sec\left(E^*\at{U^\cl}\right)
        = \left\{
            \varphi \in \Sec[0]\left(E^*\at{U^\cl}\right)
            \; \Big| \;
            \varphi
            \; \textrm{is} \;
            \Fun
        \right\}.
    \end{equation}
\end{definition}
\begin{remark}
    \label{remark:differentiable-on-closed-set}
    For sections in $\Sec[0]\left(E^*\at{U^\cl}\right)$ which are
    $\Fun$ in $U$ and have bounded derivatives the seminorm
    $\seminorm[U^\cl,k]$ is actually a norm with $\seminorm[U^\cl,0]
    \leq \seminorm[U^\cl,k]$. We obtain a norm topology on the subset
    of sections $\varphi \in \Sec[0]\left(E^*\at{U^\cl}\right)$ which
    are restrictions of $\Fun$-sections defined on an (arbitrarily
    small) open neighborhood of $U^\cl$. By definition,
    $\Sec\left(E^*\at{U^\cl}\right)$ is the Banach space completion of
    these sections. Note however that for $\varphi \in
    \Sec\left(E^*\at{U^\cl}\right)$ it is \emph{not clear} whether
    there is a section $\widetilde{\varphi} \in
    \Sec\left(E^*\at{\widetilde{U}}\right)$ with
    \begin{equation}
        \label{eq:restriction-of-Ck-sections}
        \varphi = \widetilde{\varphi} \at{U^\cl}
    \end{equation}
    for some open $\tilde{U} \supseteq U^\cl$.  In fact, the existence
    of such a $\Fun$-section $\widetilde{\varphi}$ depends very much
    of the form of the boundary $\partial U^\cl$ of $U^\cl$ which can
    be very ``wild''.
\end{remark}
Though this is a difficult question in general, we shall not be
bothered by it too much as in the end we are only interested in
$\varphi \at{U}$ for $\varphi \in \Sec\left(E^*\at{U^\cl}\right)$
which is $\Fun$ on $U$. In fact, we have that
\begin{equation}
    \label{eq:restriction-of-Ck-on-closure-to-interior}
    \Sec\left(E^*\at{U^\cl}\right)
    \ni \varphi \; \mapsto \; \varphi\at{U} \in
    \Sec\left(E^*\at{U}\right)
\end{equation}
is a continuous injective linear map with
\begin{equation}
    \label{eq:cont-estimate-of-restriction-to-interior}
    \seminorm[K,k](\varphi) \leq \seminorm[U^\cl,k](\varphi)
\end{equation}
for all compact $K \subseteq U^\cl$. This is obvious. Note however,
that in general \eqref{eq:restriction-of-Ck-on-closure-to-interior} is
far from being \emph{surjective}.
\begin{remark}
    \label{remark:diffop-on-Ck-on-closure}
    Let $D \in \Diffop^k(E^*)$ be a differential operator of order $k$
    and $\ell \geq k$. Then there is a canonical extension of
    $D\at{U}$ to $\Sec[\ell](E^*\at{U^\cl})$ such that for $\varphi
    \in \Sec[\ell]\left(E^*\at{U^\cl}\right)$ we have $D \varphi \in
    \Sec[\ell-k]\left(E^*\at{U^\cl}\right)$ and
    \begin{equation}
        \label{eq:diffop-on-Ck-on-closure}
        D: \Sec[\ell]\left(E^*\at{U^\cl}\right) \longrightarrow
        \Sec[\ell-k]\left(E^*\at{U^\cl}\right)
    \end{equation}
    is continuous. Indeed, let $\widetilde{\varphi} \in
    \Sec[\ell]\left(E^*\at{\widetilde{U}}\right)$ then
    $\seminorm[U^\cl, \ell-k](D \widetilde{\varphi}) \leq c
    \seminorm[U^\cl,k](\widetilde{\varphi})$ for some $c > 0$
    depending on $D$ by Theorem~\ref{theorem:Continuity-of-Diffops}.
    Since the restrictions of such $\widetilde{\varphi}$ to $U^\cl$
    form a dense set in the Banach space
    $\Sec[\ell]\left(E^*\at{U^\cl}\right)$ we obtain the result.
\end{remark}
\begin{lemma}
    \label{lemma:int-operator-on-Ck-on-closure}
    The operator $\mathcal{K}^\pm_{U}:
    \Sec[0]\left(E^*\at{U^\cl}\right) \longrightarrow
    \Sec[0]\left(E^*\at{U^\cl}\right)$ restricts to a continuous
    linear operator
    \begin{equation}
        \label{eq:int-op-on-Ck-on-closure}
        \mathcal{K}^\pm_{U}: \Sec\left(E^*\at{U^\cl}\right)
        \longrightarrow
        \Sec\left(E^*\at{U^\cl}\right)
    \end{equation}
    for all $k \in \mathbb{N}_0$ whose image are restrictions of
    smooth sections of $E^*$ defined on $U'$. The operator norm of
    \eqref{eq:int-op-on-Ck-on-closure} is bounded by
    \begin{equation}
        \label{eq:boundary-on-op-norm-of-intop-on-Ck-on-closure}
        \norm{\mathcal{K}^\pm_{U}}
        \leq \vol(U^\cl)
        \seminorm[U^\cl \times U^\cl,k](K^\pm_{U'}).
    \end{equation}
\end{lemma}
\begin{proof}
    Since $\Sec\left(E^*\at{U^\cl}\right) \subseteq
    \Sec[0]\left(E^*\at{U^\cl}\right) \subseteq
    \Gamma_b\left(E^*\at{U^\cl}\right)$ we can use
    Lemma~\ref{lemma:integral-operator} to get the estimate
    \[
    \seminorm[U^\cl,k](\mathcal{K}^\pm_{U} \varphi)
    \leq \vol(U^\cl)
    \seminorm[U^\cl \times U^\cl,k](K^\pm_{U'})
    \seminorm[U^\cl,0](\varphi)
    \]
    and $\mathcal{K}^\pm_{U} \varphi \in
    \Secinfty\left(E^*\at{U'}\right)$. Since in general
    $\seminorm[U^\cl,0](\varphi) \leq \seminorm[U^\cl,k](\varphi)$ for
    $\varphi \in \Sec\left(E^*\at{U^\cl}\right)$ we have the
    continuity and also the estimate for the operator norm of
    $\mathcal{K}^\pm_{U}$ as in
    \eqref{eq:boundary-on-op-norm-of-intop-on-Ck-on-closure}.
\end{proof}

If we want to repeat the argument of invertibility of
$\mathcal{K}^\pm_{U}: \Sec\left(E^*\at{U^\cl}\right) \longrightarrow
\Sec\left(E^*\at{U^\cl}\right)$ we face the following problem: for a
fixed $k$ we can certainly shrink $U$ in such a way that the operator
norm \eqref{eq:boundary-on-op-norm-of-intop-on-Ck-on-closure} becomes
less than one, but as we are interested in \emph{all} $k \in
\mathbb{N}$ the countable intersection of all shrinkings of $U$ might
be empty. Thus we have to proceed differently. The idea is that we
influence the \emph{numerical value} of the operator norm of
$\mathcal{K}^\pm_{U}$ by passing to a different but equivalent Banach
norm for $\Sec\left(E^*\at{U^\cl}\right)$.
\begin{lemma}
    \label{lemma:new-norm-on-Ck-on-closure-sections}
    Let $U \subseteq U'$ be small enough such that
    \begin{equation}
        \label{eq:new-small-enough-U}
        \delta
        = \vol(U^\cl)
        \seminorm[U^\cl \times U^\cl,0](\mathcal{K}^\pm_{U})
        < 1,
    \end{equation}
    and let $k \in \mathbb{N}_0$. Then
    \begin{equation}
        \label{eq:new-norm-on-Ck-on-closure-sections}
        \widetilde{\seminorm}_{U^\cl,k} (\varphi)
        = \seminorm[U^\cl,0](\varphi)
        + \frac{1-\delta}
        {2 \vol(U^\cl) \seminorm[U^\cl \times U^\cl,k](K^\pm_{U'}) +1}
        \seminorm[U^\cl,k](\varphi)
    \end{equation}
    defines a norm on $\Sec\left(E^*\at{U^\cl}\right) \subseteq
    \Sec[0]\left(E^*\at{U^\cl}\right)$ which is equivalent to
    $\seminorm[U^\cl,k]$. With respect to this Banach norm the
    operator $\mathcal{K}^\pm_{U}$ has operator norm
    \begin{equation}
        \label{eq:operator-norm-of-intop-for-new-norm-on-Ck-on-cl}
        \norm{\mathcal{K}^\pm_{U}}^{\widetilde{}}
        \leq \frac{1+\delta}{2} < 1.
    \end{equation}
\end{lemma}
\begin{proof}
    We know that $1-\delta > 0$ by assumption. Thus it is an easy task
    to see that the two norms $\widetilde{\seminorm}_{U^\cl,k}$ and
    $\seminorm[U^\cl,k]$ are equivalent, since
    $\seminorm[U^\cl,0](\varphi) \leq \seminorm[U^\cl,k](\varphi)$ as
    well as $\seminorm[U^\cl,0](\varphi) <
    \widetilde{\seminorm}_{U^\cl,k}(\varphi)$. Moreover, by
    \eqref{eq:boundary-on-op-norm-of-intop-on-Ck-on-closure} we have
    for the operator norm of $\mathcal{K}^\pm_{U}$ the following
    estimate
    \begin{align*}
        \widetilde{\seminorm}_{U^\cl,k}(\mathcal{K}^\pm_{U} \varphi)
        &= \seminorm[U^\cl,0](\mathcal{K}^\pm_{U} \varphi)
        + \frac{1-\delta}
        {2 \vol(U^\cl) \seminorm[U^\cl \times U^\cl,k](K^\pm_{U'}) +1}
        \seminorm[U^\cl,k](\mathcal{K}^\pm_{U} \varphi) \\
        & \leq \delta \seminorm[U^\cl,0](\varphi)
        + \frac{1-\delta}
        {2 \vol(U^\cl) \seminorm[U^\cl \times U^\cl,k](K^\pm_{U'}) +1}
        \vol(U^\cl) \seminorm[U^\cl \times U^\cl](K^\pm_{U'})
        \seminorm[U^\cl,0](\varphi) \\
        & \leq
        \left(
            \delta + \frac{1-\delta}{2}
        \right)
        \seminorm[U^\cl,0](\varphi) \\
        & \leq \frac{1+\delta}{2}
        \widetilde{\seminorm}_{U^\cl,k}(\varphi)
    \end{align*}
    for $\varphi \in \Sec\left(E^*\at{U^\cl}\right)$. Since $0 \leq
    \delta < 1$ by assumption \eqref{eq:new-small-enough-U} we
    conclude $\frac{1+\delta}{2} < 1$ as desired.
\end{proof}
\begin{corollary}
    \label{corollary:inverse-of-id-plus-int-op}
    Let $k \in \mathbb{N}_0$. Then the operator
    \begin{equation}
        \label{eq:int-op-is-cont-on-Ck}
        \id + \mathcal{K}^\pm_{U}: \Sec\left(E^*\at{U^\cl}\right)
        \longrightarrow
        \Sec\left(E^*\at{U^\cl}\right)
    \end{equation}
    is linear, continuous, and bijective with continuous inverse given
    by the absolutely norm-convergent series
    \begin{equation}
        \label{eq:inverse-of-id-plus-int-op}
        \left( \id + \mathcal{K}^\pm_{U} \right)^{-1}
        = \sum_{j=0}^\infty \left(-\mathcal{K}^\pm_{U}\right)^j.
    \end{equation}
\end{corollary}
\begin{proof}
    This is now obvious by the lemma.
\end{proof}

We shall now compute the inverse of $\id + \mathcal{K}^\pm_{U}$
slightly more explicit: in fact, it is again an integral operator with
a nice kernel. The $j$-th power of $\mathcal{K}^\pm_{U}$ is explicitly
given by
\begin{equation}
    \label{eq:jth-power-of-int-op}
    \begin{split}
        \left(
            (\mathcal{K}^\pm_{U})^j \varphi
        \right)(p)
        &= \int_{U^\cl} K^\pm_{U'}(p,z_1)
        \left(
            (\mathcal{K}^\pm_{U})^{j-1} \varphi
        \right)(z_1)
        \: \mu_g(z_1) \\
        &= \int_{U^\cl} \cdots \int_{U^\cl}
        K^\pm_{U'}(p,z_1) \cdots K^\pm_{U'}(z_{j-1},z_j)
        \varphi(z_j)
        \: \mu_g(z_1) \cdots \mu_g(z_j) \\
        &= \int_{U^\cl}
        \left(
            \int_{U^\cl} \cdots \int_{U^\cl}
            K^\pm_{U'}(p,z_1) \cdots K^\pm_{U'}(z_{j-1},q)
            \: \mu_g(z_1) \cdots \mu_g(z_{j-1})
        \right)
        \varphi(q)
        \: \mu_g(q)
    \end{split}
\end{equation}
by Fubini's theorem. Thus $(\mathcal{K}^\pm_{U})^j$ has again a nice
kernel given by
\begin{equation}
    \label{eq:kernel-of-jth-power-of-intop}
    K^{\pm (j)}_{U}(p,q)
    =
    \int_{U^\cl} \cdots \int_{U^\cl}
    K^\pm_{U'}(p,z_1) \cdots K^\pm_{U'}(z_{j-1},q)
    \: \mu_g(z_1) \cdots \mu_g(z_{j-1}).
\end{equation}
For this kernel we have the following properties:
\begin{lemma}
    \label{lemma:kernel-of-jth-power-of-intop}
    Let $j \in \mathbb{N}$. Then the $j$-th power of
    $\mathcal{K}^\pm_{U}$ has again a smooth integral kernel $K^{\pm
      (j)}_{U'} \in \Secinfty\left(E^* \extensor E \at{U' \times
          U'}\right)$ explicitly given by
    \begin{equation}
        \label{eq:the-kernel-of-jth-power-of-intop}
        K^{\pm (j)}_{U}(p,q)
        =
        \int_{U^\cl} \cdots \int_{U^\cl}
        K^\pm_{U'}(p,z_1) \cdots K^\pm_{U'}(z_{j-1},q)
        \: \mu_g(z_1) \cdots \mu_g(z_{j-1}),
    \end{equation}
    satisfying the estimate
    \begin{equation}
        \label{eq:pKl-estimate-of-power-kernel}
        \seminorm[K \times K,k] \left(K^{\pm (j)}_{U'}\right)
        \leq
        \vol(U^\cl)
        \seminorm[(K \cup U^\cl) \times (K \cup U^\cl), k]
        (K^\pm_{U'}) \delta^{j-2},
    \end{equation}
    with $\delta$ as in \eqref{eq:new-small-enough-U} where $K
    \subseteq U'$ is compact.
\end{lemma}
\begin{proof}
    The above computation \eqref{eq:jth-power-of-int-op} shows that
    \eqref{eq:the-kernel-of-jth-power-of-intop} is indeed the kernel
    of $(\mathcal{K}^\pm_{U})^j: \Sec[0]\left(E^*\at{U^\cl}\right)
    \longrightarrow \Sec[0]\left(E^*\at{U^\cl}\right)$. From the
    explicit formula \eqref{eq:the-kernel-of-jth-power-of-intop} and
    an argument analogous to the one in the proof of
    Lemma~\ref{lemma:integral-operator} we see that $K^{\pm (j)}_U$
    has a continuation for all $(p,q) \in U' \times U'$ to a
    \emph{smooth} section by the very same expression
    \eqref{eq:the-kernel-of-jth-power-of-intop}. Moreover, we can
    differentiate $K^{\pm (j)}_U$ by differentiating under the
    integral. This yields
    \begin{align*}
        \seminorm[K \times K,k]\left(K^{\pm (j)}_U\right)
        &\leq
        \int_{U^\cl} \cdots \int_{U^\cl}
        \seminorm[K \times U^\cl,k](K^\pm_{U'})
        \underbrace{
          \seminorm[U^\cl \times U^\cl,0](K^\pm_{U'}) \cdots
          \seminorm[U^\cl \times U^\cl,0](K^\pm_{U'})
        }_{j-1 \textrm{ times}} \\
        &\qquad \qquad \qquad \quad
        \seminorm[U^\cl \times K,k](K^\pm_{U'})
        \: \mu_g(z_1) \cdots \mu_g(z_{j-1}) \\
        &\leq
        \vol(U^\cl)^{j-1}
        \seminorm[U^\cl \times U^\cl,0](K^\pm_{U'})^{j-2}
        \seminorm[(K \cup U^\cl) \times (K \cup U^\cl), l]
        (K^\pm_{U'})^2
        \\
        &= \vol(U^\cl) \delta^{j-2}
        \seminorm[(K \cup U^\cl) \times (K \cup U^\cl), l]
        (K^\pm_{U'})^2,
    \end{align*}
    since only the first and last $K^\pm_{U'}$ in
    \eqref{eq:the-kernel-of-jth-power-of-intop} depend on the points
    $p,q \in K \subseteq U'$ which are used for differentiation in
    $\seminorm[K \times K, k]$. Thanks to the factorization of the
    variables, we do not get extra ($k$-dependent) constants from the
    Leibniz rule.  Thus \eqref{eq:pKl-estimate-of-power-kernel}
    follows.
\end{proof}

\begin{corollary}
    \label{corollary:kernel-of-geometric-series}
    The operator $(\id + \mathcal{K}^\pm_{U})^{-1} \circ
    \mathcal{K}^\pm_{U}$ has a smooth kernel explicitly given by the
    series $\sum_{j=1}^\infty (-1)^{j-1} K^{\pm (j)}_{U'}$, which
    converges in the $\Cinfty$-topology of $\Secinfty\left(E^*
        \extensor E \at{U' \times U'}\right)$.
\end{corollary}
\begin{proof}
    By the lemma, each $K^{\pm (j)}_{U'}$ is smooth on $U' \times U'$.
    Moreover, with respect to a given seminorm $\seminorm[K \times
    K,k]$, the above series converges since $\delta < 1$ by assumption
    on $U^\cl$. This shows that the series $\sum_{j=1}^\infty
    (-1)^{j-1} K^{\pm (j)}_{U'}$ converges (even absolutely) with
    respect to $\seminorm[K \times K,k]$. Since $K \subseteq U'$ and
    $k \in \mathbb{N}_0$ are arbitrary, we have $\Cinfty$-convergence.
    Clearly, when restricting to $U^\cl \times U^\cl$, the series is
    the kernel of $\left(\id + \mathcal{K}^\pm_{U}\right)^{-1} \circ
    \mathcal{K}^\pm_{U}$.
\end{proof}

\begin{lemma}
    \label{lemma:geometric-series-on-section}
    Let $\varphi \in \Secinfty\left(E^*\at{U'}\right)$ be smooth. Then
    $\left(\id + \mathcal{K}^\pm_{U}\right)^{-1}
    \left(\varphi\at{U^\cl}\right)$ is in
    $\Sec\left(E^*\at{U^\cl}\right)$ for all $k \in \mathbb{N}_0$.
    Moreover,
    \begin{equation}
        \label{eq:geometric-series-maps-smooth-to-smooth}
        \left(
            \id + \mathcal{K}^\pm_{U}
        \right)^{-1} \left(\varphi \at{U^\cl}\right) \At{U}
        \in \Secinfty\left(E^*\at{U}\right)
    \end{equation}
    and the map
    \begin{equation}
        \label{eq:geometric-series-is-continuous-map}
        \Secinfty\left(E^*\at{U'}\right)
        \ni \varphi \; \mapsto \;
        \left(
            \id + \mathcal{K}^\pm_{U}
        \right)^{-1} \left(\varphi \at{U^\cl}\right) \At{U}
        \in \Secinfty\left(E^*\at{U}\right)
    \end{equation}
    is continuous. The image is even in the subset of those smooth
    sections on $U$ which are restrictions of smooth sections of $E^*$
    on $U'$.
\end{lemma}
\begin{proof}
    First we note that $\varphi\at{U^\cl} \in
    \Sec\left(E^*\at{U^\cl}\right)$ by the very definition as in
    Definition~\ref{definition:differentiable-on-closed-set}.
    Moreover, since
    \[
    \seminorm[U^\cl,k]\left(\varphi\at{U^\cl}\right)
    = \seminorm[U^\cl,k](\varphi),
    \]
    the restriction map is a continuous map
    \[
    \Secinfty(E^*\at{U'}) \longrightarrow \Sec(E^*\at{U^\cl})
    \]
    for any $k \in \mathbb{N}_0$. Now $\left(\id +
        \mathcal{K}^\pm_{U}\right)^{-1}(\varphi\at{U^\cl}) \in
    \Sec(E^*\at{U^\cl})$ by
    Corollary~\ref{corollary:inverse-of-id-plus-int-op} and applying
    $\left(\id + \mathcal{K}^\pm_{U}\right)^{-1}$ is again
    continuous. Finally, restricting a section in
    $\Sec\left(E^*\at{U^\cl}\right)$ to $U$ gives a $\Fun$-section in
    the usual sense by
    \eqref{eq:restriction-of-Ck-on-closure-to-interior} in
    Remark~\ref{remark:differentiable-on-closed-set}. Moreover, this
    restriction is again continuous whence finally
    \[
    \Secinfty\left(E^*\at{U'}\right)
    \ni \varphi
    \; \mapsto \;
    \left(\id + \mathcal{K}^\pm_{U}\right)^{-1}
    \left(\varphi \at{U^\cl}\right) \At{U}
    \in \Sec\left(E^*\at{U}\right)
    \]
    is continuous for all $k \in \mathbb{N}_0$. In particular, it
    follows that $\left(\id + \mathcal{K}^\pm_{U}\right)^{-1}
    \left(\varphi\at{U^\cl}\right) \At{U} \in
    \Secinfty\left(E^*\at{U}\right)$.  Since the inverse is given by
    the geometric series we see that
    \[
    \left(\id + \mathcal{K}^\pm_{U}\right)^{-1}
    \left(\varphi\at{U^\cl}\right)
    = \varphi\at{U^\cl}
    - \left(
        \left( \id + \mathcal{K}^\pm_{U} \right)^{-1} \circ
        \mathcal{K}^\pm_{U}
    \right)
    \left(\varphi\at{U^\cl}\right).
    \]
    Now $\varphi\at{U^\cl}$ is the restriction of the smooth section
    $\varphi$ on $U'$ to $U^\cl$. Also the operator $( \id +
    \mathcal{K}^\pm_{U} )^{-1} \circ \mathcal{K}^\pm_{U}$ has a smooth
    integral kernel defined even on $U' \times U'$ by
    Corollary~\ref{corollary:inverse-of-id-plus-int-op}. Hence the
    result $\left(\left(\id + \mathcal{K}^\pm_{U}\right)^{-1} \circ
        \mathcal{K}^\pm_{U} \right) \left(\varphi\at{U^\cl}\right)$
    can also be viewed as the smooth section
    \[
    \left(
        \left(\id + \mathcal{K}^\pm_{U}\right)^{-1}
        \circ
        \mathcal{K}^\pm_{U}
    \right)
    \left(\varphi\at{U^\cl}\right) (p)
    = \int_{U^\cl}
    \left(
        \sum_{j=1}^\infty (-1)^{j-1} K^{\pm (j)}_U(p,q)
    \right) \varphi(q) \: \mu_g(q),
    \tag{$*$}
    \]
    defined even for $p \in U'$. Since the kernel of ($*$) is smooth
    it follows easily that
    \[
    \left(
        \id + \mathcal{K}^\pm_{U}
    \right)^{-1} \circ
    \mathcal{K}^\pm_{U} :
    \Sec\left(E^*\at{U^\cl}\right) \longrightarrow
    \Secinfty\left(E^*\at{U'}\right)
    \]
    is a continuous linear map: this can be done analogously to the
    argument in Lemma~\ref{lemma:integral-operator} where we only have
    to replace $K^\pm_U$ by the smooth kernel of ($*$) in
    \eqref{eq:cont-estimate-of-integral-operator}. This shows that
    $\left(\left(\id + \mathcal{K}^\pm_{U}\right)^{-1} \circ
        \mathcal{K}^\pm_{U}\right) \left(\varphi\at{U^\cl}\right) \in
    \Secinfty\left(E^*\at{U'}\right)$ and hence
    \eqref{eq:geometric-series-is-continuous-map}. Moreover, the
    composition of all the involved maps including the last
    restriction to $U$ are continuous. Thus
    \eqref{eq:geometric-series-is-continuous-map} is continuous as
    well.
\end{proof}
Note that $\left(\id + \mathcal{K}^\pm_U\right)^{-1}$ is defined even
on $U'$ via the integral formula. But here it is no longer the inverse
of the operator $\id + \mathcal{K}^\pm_U$.

We can now use the inverse of $\id + \mathcal{K}^\pm_{U}$ to build a
true fundamental solution as follows:
\begin{definition}[Local fundamental solution]
    \label{definition:local-fundamental-solution}
    \index{Fundamental solution!local}%
    Let $U' \subseteq M$ be geodesically convex and $U \subseteq U'$
    be open with compact closure $U^\cl \subseteq U'$ such that the
    volume of $U^\cl$ is small enough. For $p \in U$ we define
    \begin{equation}
        \label{eq:local-fundamental-solution}
        F^\pm_U(p): \Secinfty_0\left(E^*\at{U}\right)
        \ni \varphi \; \mapsto \;
        \left( \id + \mathcal{K}^\pm_{U} \right)^{-1}
        \left(\widetilde{\mathcal{R}}^\pm_U(\argument) (\varphi)\right)
        \At{p}
        \in E^*_p.
    \end{equation}
\end{definition}
\begin{theorem}[Local fundamental solution]
    \label{theorem:local-fundamental-solution}
    Let $U' \subseteq M$ be geodesically convex and let $U \subseteq
    U'$ be open with compact closure $U^\cl \subseteq U'$ such that
    the volume of $U^\cl$ is small enough. Then for $p \in U$ the map
    \begin{equation}
        \label{eq:loc-fund-solution}
        F^\pm_U(p): \Secinfty_0\left(E^*\at{U}\right)
        \longrightarrow E^*_p
    \end{equation}
    is a local fundamental solution of $D$ at $p$ such that for every
    $\varphi \in \Secinfty_0\left(E^*\at{U}\right)$
    \begin{equation}
        \label{eq:local-fund-solution--paired-with-varphiis-smooth}
        F^\pm_U(\argument)\varphi: p \; \mapsto \; F^\pm_U(p)(\varphi)
    \end{equation}
    is a smooth section of $E^*$ over $U$. In fact,
    \begin{equation}
        \label{eq:local-fund-solution-is-cont}
        F^\pm_U: \Secinfty_0\left(E^*\at{U}\right)
        \ni \varphi \; \mapsto \;
        F^\pm_U(\argument)(\varphi) \in \Secinfty\left(E^*\at{U}\right)
    \end{equation}
    is a continuous linear map.
\end{theorem}
\begin{proof}
    From Theorem~\ref{theorem:approx-fundamental-solution},
    \refitem{item:approx-solution-smooth-in-p} we know that
    $\widetilde{\mathcal{R}}^\pm_{U'}(\argument)(\varphi)$ defines a
    smooth section of $E^*$ over $U'$.  By
    Proposition~\ref{proposition:estimates-and-cont-of-pairing-with-approx-solution}
    we know that $\widetilde{\mathcal{R}}^\pm_{U'}:
    \Secinfty_0(E^*\at{U'}) \ni \varphi \mapsto
    \widetilde{\mathcal{R}}^\pm_{U'}(\argument) \varphi \in
    \Secinfty(E^*\at{U'})$ is continuous in the $\Cinfty_0$- and
    $\Cinfty$-topology, respectively.  By
    Lemma~\ref{lemma:geometric-series-on-section}, also the map
    \[
    \Secinfty\left(E^*\at{U'}\right)
    \ni \varphi \; \mapsto \;
    \left(
        \id + \mathcal{K}^\pm_{U}
    \right)^{-1} \left(\varphi\at{U^\cl}\right) \At{U}
    \in \Secinfty\left(E^*\at{U}\right)
    \]
    is continuous, whence it follows that
    \eqref{eq:local-fund-solution-is-cont} is continuous and
    linear. This also implies
    \eqref{eq:local-fund-solution--paired-with-varphiis-smooth}. Thus
    it remains to shows that $F^\pm_U(p)$ is indeed a fundamental
    solution of $D$ at $p$. We compute
    \begin{align*}
        \left(
            D F^\pm_U(p)
        \right) (\varphi)
        &= F^\pm_U(p) (D^\Trans \varphi) \\
        &=
        \left(
            \left(
                \id + \mathcal{K}^\pm_{U}
            \right)^{-1}
            \left(
                \widetilde{\mathcal{R}}^\pm_{U}(\argument)
                (D^\Trans \varphi)
            \right)
        \right)
        \At{p}\\
        &=
        \left(
            \left(
                \id + \mathcal{K}^\pm_{U}
            \right)^{-1}
            \left(
                \left(
                    D \widetilde{\mathcal{R}}^\pm_U(\argument)
                \right) (\varphi)
            \right)
        \right)
        \At{p} \\
        &=
        \left(
            \left(
                \id + \mathcal{K}^\pm_{U}
            \right)^{-1}
            (\varphi + \mathcal{K}^\pm_U \varphi)
        \right)
        \At{p}\\
        &= \varphi(p)
    \end{align*}
    by \eqref{eq:D-on-approx-solution-is}. But this is precisely the
    defining property of a fundamental solution.
\end{proof}
\begin{corollary}
    \label{corollary:local-fundamental-solution}
    Let $D \in \Diffop^2(E)$ be a normally hyperbolic differential
    operator. Then every point in $M$ has a small neighborhood $U
    \subseteq M$ such that on $U$ we have a fundamental solution
    $F^\pm_U(p)$ for all $p \in U$, i.e.
    \begin{equation}
        \label{eq:local-fund-solution}
        D F^\pm_U(p) = \delta_p,
    \end{equation}
    and such that the linear map
    \begin{equation}
        \label{eq:fund-sol-on-varphi-yields-smooth-section}
        F^\pm_U:
        \Secinfty_0\left(E\at{U}\right)
        \ni \varphi
        \; \mapsto \;
        \left(
            \; p \mapsto \; F^\pm_U(p)(\varphi)
        \right)
        \in \Secinfty\left(E\at{U}\right)
    \end{equation}
    is continuous.
\end{corollary}

%
%

\subsection{Causal Properties of $F^\pm_U$}
\label{satz:causal-props-of-fund-solution}

The construction of the integral operator $\mathcal{K}^\pm_U$ and the
invertibility of $\id + \mathcal{K}^\pm_U$ works for arbitrary small
enough $U \subseteq U'$. However, since $\mathcal{K}^\pm_U$ is
\emph{non-local} the nice support properties of
$\widetilde{\mathcal{R}}^\pm_{U'}$ are typically destroyed. To
guarantee good causal behaviour we need to put some extra conditions
on $U$.

\begin{remark}
    \label{remark:causal-choice-of-U}
    \index{Causal subset}%
    \index{Causal compatibility}%
    Let $U \subseteq U'$ be causal, i.e. for $p,q \in U^\cl \subseteq
    U'$ we have $J^\pm_{U'}(p,q) \subseteq U^\cl$ and the diamond is
    compact. Then $U^\cl$ is causally compatible with $U'$. Indeed, if
    say $q \in J^+_{U'}(p)$ then we can join $p$ and $q$ by a unique
    future directed geodesic which is entirely in
    $J^+_{U'}(p,q)$. Thus this curve is also entirely in $U^\cl$
    whence $q \in J^+_{U^\cl}(p)$ proving that $U^\cl$ is causally
    compatible with $U'$.
\end{remark}
In the following, we assume that $U \subseteq U'$ is in addition a
causal subset. As a first consequence we have
\begin{equation}
    \label{eq:causally-compatible-choice}
    J^\pm_{U^\cl}(p)
    = J^\pm_{U'}(p) \cap U^\cl
\end{equation}
for $p \in U^\cl$.
\begin{lemma}
    \label{lemma:causal-choice-of-U}
    Let $U \subseteq U'$ be in addition causal. Then for $\varphi \in
    \Sec[0]\left(E^*\at{U^\cl}\right)$ we have
    \begin{equation}
        \label{eq:supp-of-intop-varphi-for-causal-choice-of-U}
        \supp (\mathcal{K}^\pm_U \varphi)
        \subseteq J^\mp_{U^\cl} (\supp \varphi).
    \end{equation}
\end{lemma}
\begin{proof}
    We know that $(p,q) \in \supp K^\pm_{U'}$ implies $q \in
    J^\pm_{U'}(p)$ by Lemma~\ref{lemma:support-feature-of-defect}.
    Thus for $p \in U^\cl$ and
    \[
    \left(
        \mathcal{K}^\pm_U \varphi
    \right)(p)
    = \int_{U^\cl} K^\pm_{U'}(p,q) \cdot \varphi(q) \: \mu_g(q)
    \]
    we get $(\mathcal{K}^\pm_U \varphi)(p) = 0$ if the
    integrand vanishes identically. But if
    $K^\pm_{U'}(p,q) \cdot \varphi(q) \neq 0$ for some $(p,q)$ then on
    one hand $q \in J^\pm_{U'}(p)$ by the support features of
    $K^\pm_{U'}(p,q)$ and $q \in \supp \varphi$ on the other
    hand. Thus $q \in J^\pm_{U'}(p) \cap \supp \varphi$ follows. We
    conclude that necessarily $(\mathcal{K}^\pm_U \varphi)(p) = 0$ if
    $J^\pm_{U'}(p) \cap \supp \varphi = \emptyset$.
    \begin{figure}
        \centering
        \input{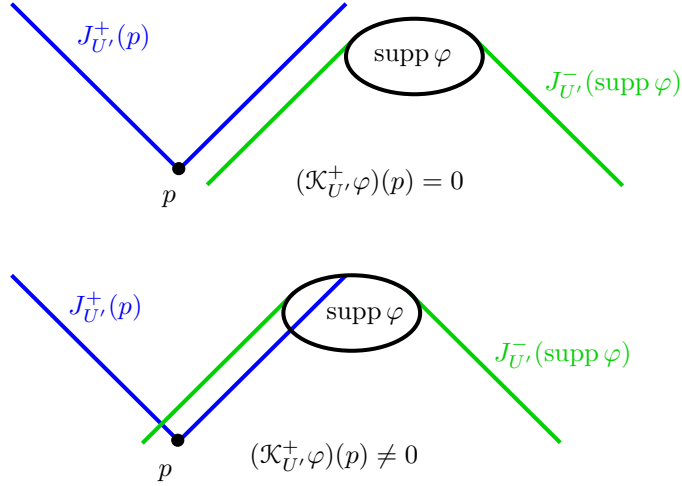}
        \caption{\label{fig:illustration-for-supps}%
          The relation between the supports in the proof of
          Lemma~\ref{lemma:causal-choice-of-U}.
        }
    \end{figure}
    From this we conclude
    \[
    \supp (\mathcal{K}^\pm_U \varphi) \subseteq
    J^\mp_{U'}(\supp \varphi) \cap U^\cl
    = J^\mp_{U^\cl}(\supp \varphi)
    \]
    see also Figure~\ref{fig:illustration-for-supps}.
\end{proof}

To compute the support of $\left(\id + \mathcal{K}^\pm_U\right)^{-1}
\varphi$ one may have the idea that with
\eqref{eq:supp-of-intop-varphi-for-causal-choice-of-U} also the finite
powers of $\mathcal{K}^\pm_U$ have the property
\eqref{eq:supp-of-intop-varphi-for-causal-choice-of-U}. This is indeed
correct as by induction and
\eqref{eq:supp-of-intop-varphi-for-causal-choice-of-U}
\begin{equation}
    \label{eq:supp-of-finite-power-intop-on-varphi}
    \begin{split}
        \supp \left(
            (\mathcal{K}^\pm_U)^j \varphi
        \right)
        \subseteq J^\mp_{U} \left(
            \supp \left(
                (\mathcal{K}^\pm_U)^{j-1} \varphi
            \right)
        \right)
        \subseteq J^\mp_U \left(
            J^\mp_U \cdots J^\mp_U
            (\supp \varphi)
        \right)
        = J^\mp_U(\supp \varphi),
    \end{split}
\end{equation}
since clearly $J^\mp_U(A) = J^\mp_U(J^\mp_U(A))$ for arbitrary $A
\subseteq U$. However, taking the geometric series for $\left(\id +
    \mathcal{K}^\pm_U\right)^{-1}$ would require to take the closure
of the union of countably many closed subsets of $J^\mp_U(\supp
\varphi)$. Now $J^\mp_U(\supp \varphi)$ need not be closed at all,
even though $\supp \varphi$ is closed. Thus we can \emph{not} conclude
by this argument that the support of $\left(\id +
    \mathcal{K}^\pm_U\right)^{-1} \varphi$ lies in $J^\mp_U(\supp
\varphi)$. However, we can proceed as follows:
\begin{lemma}
    \label{lemma:supps-of-int-kernels-are-future-past-stretched}
    \index{Future stretched}%
    \index{Past stretched}%
    For all $j \in \mathbb{N}$ the supports of the integral kernels
    $K^{\pm (j)}_{U'}$ of $(\mathcal{K}^\pm_U)^j$ are future
    respectively past stretched, i.e.
    \begin{equation}
        \label{eq:supp-of-int-kernels-are-future-past-stretched}
        (p,q) \in \supp K^{\pm (j)}_{J'} \subseteq U' \times U'
        \Longrightarrow q \in J^\pm_{U'}(p).
    \end{equation}
    Moreover, the support of the integral kernel of $(\id +
    \mathcal{K}^\pm_U)^{-1} \circ \mathcal{K}^\pm_U$ is also future
    respectively past stretched.
\end{lemma}
\begin{proof}
    Assume that $K^{\pm (j)}_{U'}(p,q) \neq 0$. Then the integrand in
    \eqref{eq:the-kernel-of-jth-power-of-intop} can not be identically
    zero whence there have to be $z_1, \ldots, z_{j-1} \in U^\cl$ with
    $z_1 \in J^\pm_{U'}(p), \ldots, z_{j-1} \in J^\pm_{U'}(z_{j-2}), q
    \in J^\pm_{U'}(z_{j-1})$. But this means $q \in J^\pm_{U'}(p)$
    proving \eqref{eq:supp-of-int-kernels-are-future-past-stretched}
    with the same closure argument as in the proof of
    Lemma~\ref{lemma:supps-of-int-kernels-are-future-past-stretched}.
    Now we consider the $\Cinfty$-convergent sum of the $K^{\pm
      (j)}_{U'}$. If $\sum_{j=1}^\infty (-1)^{j-1} K^{\pm
      (j)}_{U'}(p,q) \neq 0$ for some $(p,q)$ then at least for one
    $j$ we have $K^{\pm (j)}_{U'}(p,q) \neq 0$. Thus $q \in
    J^\pm_{U'}(p)$ and we can proceed as before.
\end{proof}

\begin{corollary}
    \label{corollary:supp-of-geometric-on-varphi}
    For $\varphi \in \Sec[0]\left(E^*\at{U^\cl}\right)$ we have
    \begin{equation}
        \label{eq:supp-of-geometric-on-varphi}
        \supp \left(
            \left(\id + \mathcal{K}^\pm_U\right)^{-1} \varphi
        \right)
        \subseteq J^\mp_{U^\cl}(\supp \varphi).
    \end{equation}
\end{corollary}
\begin{proof}
    Clearly, we have
    \[
    \left(\id + \mathcal{K}^\pm_U\right)^{-1} \varphi
    = \varphi
    - \left(\id + \mathcal{K}^\pm_U\right)^{-1}
    \circ
    \mathcal{K}^\pm_U \varphi.
    \]
    With $\supp \varphi \subseteq J^\mp_{U^\cl} (\supp \varphi)$ and
    the above lemma the statement follows at once as in
    Lemma~\ref{lemma:causal-choice-of-U}.
\end{proof}

Using this property of $\left(\id + \mathcal{K}^\pm_U\right)^{-1}$ for
\emph{causal} $U$ we arrive at the following statement:
\begin{theorem}[Local Green functions]
    \label{theorem:local-green-functions}
    \index{Green function!local}%
    Let $U \subseteq U'$ be small enough and causal. Then the
    fundamental solutions $F^\pm_U(p)$ from
    Theorem~\ref{theorem:local-fundamental-solution} are advanced and
    retarded Green functions, i.e. we have
    \begin{equation}
        \label{eq:local-green-functions-supp}
        \supp F^\pm_U(p)
        \subseteq J^\pm_U(p).
    \end{equation}
\end{theorem}
\begin{proof}
    Let $\varphi \in \Secinfty_0\left(E^* \at{U}\right)$ be a test
    section. Then
    \begin{align*}
        \supp \left(
            F^\pm_U(\argument) (\varphi)
        \right)
        &= \supp \left(
            \left(\id + \mathcal{K}^\pm_U\right)^{-1}
            \widetilde{\mathcal{R}}^\pm_{U'}(\argument) (\varphi)
        \right) \\
        &\subseteq J^\mp_{U^\cl} \left(
            \widetilde{\mathcal{R}}^\pm_{U'}(\argument) (\varphi)
        \right) \\
        &\subseteq J^\mp_{U^\cl} \left(
            J^\mp_{U^\cl} (\supp \varphi)
        \right)
        = J^\mp_{U^\cl} (\supp \varphi),
        \tag{$*$}
    \end{align*}
    since $\supp \widetilde{\mathcal{R}}^\pm_{U'}(p) \at{U} \subseteq
    J^\pm_{U^\cl}(p)$ whence for $\supp \varphi \cap J^\pm_{U^\cl}(p)
    = \emptyset$ we conclude $\widetilde{\mathcal{R}}^\pm_{U'}(p)
    (\varphi) = 0$. Thus $p \notin J^\mp_{U^\cl}(\supp \varphi)$
    implies $\widetilde{\mathcal{R}}^\pm_{U'}(p) (\varphi) = 0$. Since
    for compactly supported $\varphi$ we have a \emph{closed}
    $J^\mp_{U^\cl}(\supp \varphi)$ by $U$ being causal we conclude
    that $\supp \widetilde{\mathcal{R}}^\pm_{U'}(\argument) (\varphi)
    \subseteq J^\mp_{U^\cl}(\supp \varphi)$. This shows ($*$). Thus if
    $\supp \varphi \cap J^\pm_{U^\cl}(p) = \emptyset$ for $p \in
    U^\cl$ then $p \notin J^\mp_U(\supp \varphi)$ and thus $p \notin
    \supp (F^\pm_U(\argument)(\varphi))$ whence $F^\pm_U(p)(\varphi) =
    0$ follows. But this implies \eqref{eq:local-green-functions-supp}
    as $J^\pm_{U^\cl}(p)$ is closed thanks to $U$ being causal.
\end{proof}

Since every point in a time-oriented Lorentz manifold has an
arbitrarily small causal neighborhood we finally arrive at the
following result:
\begin{corollary}
    \label{corollary:local-green-functions}
    \index{Green function!local existence}%
    Let $D \in \Diffop^2(E)$ be normally hyperbolic. Then every point
    in $M$ has a small enough causal neighborhood $U \subseteq M$ such
    that on $U$ we have advanced and retarded Green functions
    $F^\pm_U(p)$ at $p \in U$, i.e.
    \begin{equation}
        \label{eq:local-green}
        D F^\pm(p) = \delta_p
    \end{equation}
    and
    \begin{equation}
        \label{eq:local-green-supp}
        \supp F^\pm_U(p) \subseteq J^\pm_U(p),
    \end{equation}
    such that in addition
    \begin{equation}
        \label{eq:local-green-is-cont}
        F^\pm_U:
        \Secinfty_0\left(E^*\at{U}\right)
        \ni \varphi \; \mapsto \;
        \left(
            p \; \mapsto \; F^\pm_U(p)(\varphi)
        \right) \in \Secinfty\left(E^*\at{U}\right)
    \end{equation}
    is a continuous linear map.
\end{corollary}


%% file: local.tex
%
%

In this section we show how the Green functions $F_U^\pm(p)$ can be
used to obtain solutions to the wave equation
\begin{equation}
    \label{eq:wave-equation}
    D u = v
\end{equation}
with a prescribed source term $v$. The main idea is that a suitable
$v$ can be written as a superposition\index{Superposition} of
$\delta$-functionals. Since $F^\pm_U(p)$ solves
\eqref{eq:wave-equation} for $v = \delta_p$ we get a solution to
\eqref{eq:wave-equation} for arbitrary $v$ by taking the corresponding
superposition of the fundamental solutions $F^\pm_U(p)$. Of course, at
the moment we are restricted to $v$ having compact support in $U$.

Then we are interested in two extreme cases: for a distributional $v$
we can only expect to obtain distributions $u$ as solutions. However,
if $v$ has good regularity then we can expect $u$ to be regular as
well.

%
%

\subsection{Local Solutions for Distributional Inhomogeneity}
\label{satz:local-solution-distributional-inhomog}

Let $v \in \Sec[-\infty]_0(E\at{U})$ be a generalized section of $E$
with compact support in $U$. We want to solve
\begin{equation}
    \label{eq:distributional-wave-equation}
    \index{Wave equation!inhomogeneous}%
    Du^\pm = v
\end{equation}
with some $u^\pm \in \Sec[-\infty](E\at{U})$.
\begin{remark}
    \label{remark:support-of-solution-not-compact}
    \index{Source term}%
    \index{Propagating wave}%
    Since a normally hyperbolic differential operator $D$ describes a
    wave equation we expect from physical considerations that a source
    term $v$ causes \emph{propagating} waves whence the support of
    $u^\pm$ is expected to be non-compact: In fact, the best we can
    hope for is that in spatial directions the support stays compact
    while in time directions we will have non-compact support at least
    in either the future or the past. Up to now we are dealing with
    the local situation $U \subseteq M$ where thanks to the simple
    geometry those questions are rather harmless. Later on this issue
    will become more subtle.
\end{remark}
\begin{lemma}
    \label{lemma:dualazinig-of-fund-solution}
    \index{Inhomogeneity!distributional}%
    Let $U \subseteq M$ be a small enough open subset such that the
    construction of $F^\pm_U$ as in
    Section~\ref{satz:LocalFundamentalSolution} applies.
    \begin{lemmalist}
    \item \label{item:transposing-of-fund-solution} The map $F^\pm_U:
        \Secinfty_0(E^*\at{U}) \longrightarrow \Secinfty(E^*\at{U})$
        induces a linear map
        \begin{equation}
            \label{eq:transposed-fund-solution}
            (F^\pm_U)': \Sec[-\infty]_0(E\at{U}) \longrightarrow
            \Sec[-\infty](E\at{U})
        \end{equation}
        by dualizing, i.e. for $v \in \Sec[-\infty]_0(E\at{u})$ and
        $\varphi \in \Secinfty_0(E^*\at{U})$ one defines
        \begin{equation}
            \label{eq:dualizing-of-fund-solution}
            \left(
                (F^\pm_U)' (v)
            \right) (\varphi)
            = v \left(
                F^\pm_U(\varphi)
            \right).
        \end{equation}
    \item \label{item:transposed-fund-solution-is-weakstar-cont} The
        map $(F^\pm_U)'$ is weak$^*$ continuous.
    \item We have
        \begin{equation}
            \label{eq:D-and-transposed-fund-solution}
            D (F^\pm_U)'(v) = v
        \end{equation}
        for all $v \in \Sec[-\infty]_0(E\at{U})$.
    \end{lemmalist}
\end{lemma}
\begin{proof}
    For the first part we recall that we have the identification
    \[
    \Secinfty_0(E^*\at{U}) \ni \varphi
    \; \mapsto \;
    \varphi \tensor \mu_g \in
    \Secinfty_0(E^*\at{U} \tensor \Dichten T^*M)
    \]
    from which we obtain the identification
    \[
    \Sec[-\infty](E\at{U}) \ni u \; \mapsto \;
    ( \varphi \mapsto u(\varphi \tensor \mu_g))
    \in \Secinfty_0(E^*\at{U})'.
    \tag{$*$}
    \]
    Since tensoring with $\mu_g > 0$ does not change the supports we
    can dualize the continuous map
    \[
    F^\pm_U: \Secinfty_0(E^*\at{U}) \longrightarrow
    \Secinfty(E^*\at{U})
    \]
    to a map
    \[
    (F^\pm_U)': \Secinfty(E^*\at{U})' \longrightarrow
    \Secinfty_0(E^*\at{U})'.
    \tag{$**$}
    \]
    Using ($*$) and the fact that the dual space of all test sections
    are the compactly supported generalized sections, see
    Theorem~\ref{theorem:gensec-with-compact-support}, we get
    \[
    \Sec[-\infty]_0(E\at{U})
    \stackrel{(*)}{\longrightarrow}
    \Secinfty(E^*\at{U})'
    \stackrel{(F^\pm_U)'}{\longrightarrow}
    \Secinfty_0(E^*\at{U})'
    \stackrel{(*)}{\longrightarrow}
    \Sec[-\infty](E\at{U}),
    \]
    whose composition we denote by $(F^\pm_U)'$ as well. This is the
    map \eqref{eq:transposed-fund-solution}. Dualizing yields a
    weak$^*$ continuous map in ($**$). Finally, the identifications
    ($*$) are weak$^*$ continuous as well, hence it results in a
    weak$^*$ continuous map \eqref{eq:transposed-fund-solution}. Note
    that in \eqref{eq:dualizing-of-fund-solution} we have hidden the
    aspect of the reference density $\mu_g$ in the pairing of $v$ and
    $F^\pm_U(\varphi)$. This shows the first and second part. For the
    third part we unwind the definition of $D F^\pm_U{}'$. Let
    $\varphi \in \Secinfty_0(E^*\at{U})$ be a test section and compute
    \begin{align*}
        \left(D \left((F^\pm_U)'(v)\right)\right)(\varphi)
        &= \left((F^\pm_U)'(v)\right) (D^\Trans \varphi) \\
        &= v\left(
            p \; \mapsto \; F^\pm_U (D^\Trans \varphi)\at{p}
        \right) \\
        &= v\left(
            p \; \mapsto \; (F^\pm_U(p))(D^\Trans \varphi)
        \right) \\
        &= v\left(
            p \; \mapsto \; \varphi(p)
        \right) \\
        &= v(\varphi),
    \end{align*}
    using the definition of the dualized map and the feature $D
    F^\pm_U(p) = \delta_p$. But this means
    \eqref{eq:D-and-transposed-fund-solution}.
\end{proof}

\begin{figure}
    \centering
    \input{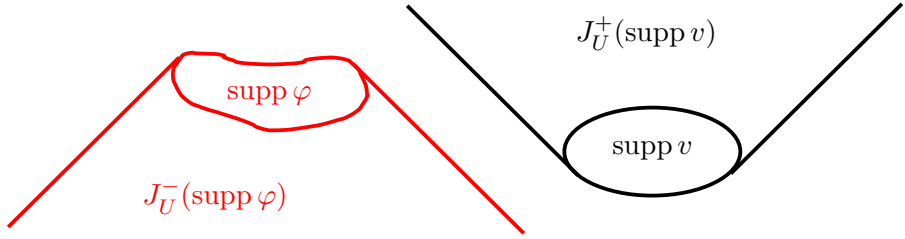}
    \caption{\label{fig:supp-of-solution}%
      The supports of $\varphi$ and $v$ in
      Lemma~\ref{lemma:supp-of-solution}.
    }
\end{figure}

\begin{remark}[Fundamental solutions]
    \label{remark:fund-solutions}
    \index{Fundamental solution}%
    We note that in the above proof we have not used any details of
    the properties of $D$ or $F^\pm_U$. The only thing we needed was
    the property that
    \begin{equation}
        \label{eq:pairing-with-fund-is-cont}
        F^\pm_U: \Secinfty_0(E^*\at{U}) \ni \varphi
        \; \mapsto \;
        \left(p \; \mapsto \; F^\pm_U(p)(\varphi)\right)
        \in \Secinfty(E\at{U})
    \end{equation}
    is continuous in the $\Cinfty_0$- and $\Cinfty$-topology in order
    to dualize \eqref{eq:pairing-with-fund-is-cont} to a map
    \eqref{eq:transposed-fund-solution} \emph{and} the fundamental
    solution property
    \begin{equation}
        \label{eq:fund-solution-prop}
        D F^\pm_U(p) = \delta_p
    \end{equation}
    in order to compute $D (F^\pm_U)'(v)$ as in
    \eqref{eq:D-and-transposed-fund-solution}. Thus the above argument
    shows one principle usage of fundamental solutions: they allow to
    solve the inhomogeneous equations in a distributional sense. Of
    course, up to now we have just found on particular solution for
    each inhomogeneity $v$ but no uniqueness. In fact, for our wave
    equations we expect to have many solutions as we expect traveling
    waves for trivial inhomogeneity $v=0$. Thus we have to specify
    boundary conditions in order to get more specific solutions. In
    order to control the ``boundary conditions'' in our case, we use
    the fundamental solutions $F^\pm_U(p)$ as in
    Theorem~\ref{theorem:local-green-functions}, i.e. on a causal $U
    \subseteq M$.
\end{remark}
\begin{lemma}
    \label{lemma:supp-of-solution}
    Let $U \subseteq M$ be small enough and causal and let
    $F^\pm_U(p)$ be the corresponding fundamental solutions as in
    Theorem~\ref{theorem:local-green-functions}. For $v \in
    \Sec[-\infty]_0(E \at{U})$ we have
    \begin{equation}
        \label{eq:supp-of-solution}
        \supp (F^\pm_U)'(v) \subset J^\pm_U (\supp v).
    \end{equation}
\end{lemma}
\begin{proof}
    We use the causality property $\supp F^\pm_U(p) \subseteq
    J^\pm_U(p)$ for all $p \in U$ of the fundamental solution. Thus
    let $\varphi \in \Secinfty_0(E^*\at{U})$ be a test section with
    $\supp \varphi \cap J^\pm_U(\supp v) = \emptyset$. We have to show
    $((F^\pm_U)'(v))(\varphi) = 0$ for all such $\varphi$. We compute
    \[
    (F^\pm_U{}'(v))(\varphi)
    = v\left(F^\pm_U (\varphi)\right)
    = v\left(p \; \mapsto \; F^\pm_U(p)(\varphi)\right).
    \]
    From the proof of Theorem~\ref{theorem:local-green-functions} we
    know that $\supp (p \mapsto F^\pm_U(p)(\varphi)) \subset
    J^\mp_{U^\cl}(\supp \varphi)$. But $\supp u \cap J^\pm_U(\supp
    \varphi) = \emptyset$ by assumption whence $v( p \mapsto
    F^\pm_U(p)(\varphi)) = 0$ follows, see also
    Figure~\ref{fig:supp-of-solution}.
\end{proof}

\begin{remark}
    \label{remark:2}
    Even though we do not yet have the uniqueness properties, already
    at this stage we see some very nice features familiar from our
    physically motivated expectations:
    \begin{remarklist}
    \item \index{Light speed}%
        Using the solution $u^+ = (F^+_U)'(v)$ of the inhomogeneous
        wave equation we see that the influence of the source term $v$
        is only in the future of $v$. This is a physically reasonable
        behaviour. The interpretation is that at some time one
        switches on a source term, e.g. an oscillating dipole, and
        observes emitted waves $u^+$ in the future of $v$. In
        particular, the signals emitted by $v$ can not propagate
        faster than with light speed. The solution $u^-$ is the other
        extreme which for physical reasons is not acceptable.
    \item In the flat situation of the Minkowski spacetime
        $(\mathbb{R}^n, \eta)$ we can take $U = \mathbb{R}^n$ and
        obtain $F^\pm_{\mathbb{R}^n}(0) = R^\pm(2)$ and
        $F^\pm_{\mathbb{R}^n}(p)$ is the translated Riesz distribution
        for arbitrary $p \in \mathbb{R}^n$. Then the construction of
        the solutions $(F^\pm_{\mathbb{R}^n})'(v)$ for a given $v$ is
        the well-known solution procedure as known e.g. from
        electrodynamics \cite{scheck:2005a, jackson:1975a}.
    \item Of particular interest is the following situation: a charged
        pointlike particle with charge $e$ moves along a trajectory $t
        \mapsto \vec{x}(t)$ in Minkowski spacetime with velocity
        $|\vec{v}(t)| = |\dot{\vec{x}}(t)| < 1$. As usual, we set the
        speed of light $c=1$ by choosing an appropriate unit
        system. Then the charge density is $\varrho(t,\vec{x})= e
        \delta_{\vec{x}(t)}$ while the current density is
        $\vec{j}(t,\vec{x}) = e \vec{v}(t) \delta_{{x}(t)}$, viewed
        both as distributions on the spatial $\mathbb{R}^{n-1}$ inside
        Minkowski spacetime. They combine into an
        $\mathbb{R}^n$-valued distribution on $\mathbb{R}^n$ denoted
        by $j$. The corresponding solution $A = (F^\pm_U)'(j)$ of
        $\dAlembert A = j$ is then known as the
        \emIndex{Lienhard-Wiechert potential}. It describes the
        electromagnetic potential of the radiation emitted by the
        moving charge, see e.g. \cite[Sect.~3.6]{scheck:2005a} or
        \cite[Sect.~14.1]{jackson:1975a}.
    \item From our construction, $(F^\pm_U)'$ is only defined on the
        distributional sections with \emph{compact} support. However,
        the example of the moving charge gives an inhomogeneity with
        non-compact support, at least in timelike directions: Here
        only the support in spatial directions is compact for all
        times. Thus for physical applications it will be necessary to
        extend the domain of $(F^\pm_U)'$ to more general
        distributions.
    \end{remarklist}
\end{remark}

%
%

\subsection{Local Solution for Smooth Inhomogeneity}
\label{satz:local-solution-smooth-inhomogeneity}

In a next step we want to discuss the additional properties of the
solutions $(F^\pm_U)'(v)$ of the inhomogeneous wave equation $Du = v$
for distributional $v$ having some kind of regularity. Of particular
interest is the case where $v$ is actually \emph{smooth} and hence a
test section $v \in \Secinfty_0(E\at{U})$.

To this end we first collect some more specific properties of the
operator $\left(\id + \mathcal{K}^\pm_U\right)^{-1}$. It will be
advantageous to consider integral operators with smooth kernel in
general. Thus we consider the following situation: Let $U \subset M$
be open with $U^\cl$ compact and let $U^\cl \subseteq U'$ with $U'$
open. Moreover, let $K \in \Secinfty(E^* \extensor E \at{U' \times
  U'})$ be a smooth kernel on the larger open subset $U' \times
U'$. For sections $\varphi \in \Sec[]_b(E^*\at{U})$ we consider the
integral operator
\begin{equation}
    \label{eq:integral-operator}
    \index{Integral operator}%
    \index{Smooth kernel}%
    (\mathcal{K} \varphi)(p)
    = \int_{U^\cl} K(p,q) \cdot \varphi(q) \: \mu_g(q)
\end{equation}
analogously to \eqref{eq:the-integral-operator}, where $p \in
U'$. Repeating the arguments from Lemma~\ref{lemma:integral-operator}
and Lemma~\ref{lemma:int-operator-on-Ck-on-closure} we obtain the
following general result:
\begin{lemma}
    \label{lemma:properties-of-intop-with-smooth-kernel}
    Let $U \subseteq U^\cl \subseteq U'$ with $U, U'$ open and $U^\cl$
    compact. For the integral operator $\mathcal{K}$ corresponding to
    a smooth kernel $K \in \Secinfty(E^* \extensor E \at{U' \times
      U'})$ as in \eqref{eq:integral-operator} the following
    statements are true:
    \begin{lemmalist}
    \item \label{item:int-op-smoothes} For $\varphi \in
        \Sec[]_b(E^*\at{U})$ one has $\mathcal{K} \varphi \in \Sec(E^*
        \at{U^\cl})$ for all $k \in \mathbb{N}_0$ and $\mathcal{K}
        \varphi \at{U} \in \Secinfty(E^*\at{U})$.
    \item \label{item:intop-is-continuous} The maps (all denoted by
        $\mathcal{K}$)
        \begin{equation}
            \label{eq:intop-to-Ck}
            \mathcal{K}: \Sec[]_b(E^*\at{U}) \ni \varphi
            \; \mapsto \;
            \mathcal{K} \varphi \in \Sec(E^* \at{U^\cl})
        \end{equation}
        and
        \begin{equation}
            \label{eq:intop-to-smooth}
            \mathcal{K}: \Sec[]_b(E^*\at{U}) \ni \varphi
            \; \mapsto \;
            \mathcal{K} \varphi \at{U} \in \Secinfty(E^*\at{U})
        \end{equation}
        are continuous. In fact, for $k \in \mathbb{N}_0$ one even has
        \begin{equation}
            \label{eq:intop-cont-estimate}
            \seminorm[U^\cl, k] (\mathcal{K} \varphi)
            \leq c \seminorm[U^\cl, 0](\varphi)
        \end{equation}
        for some $c > 0$ depending on $k$.
    \end{lemmalist}
\end{lemma}
\begin{proof}
    For the first part we can copy the proof of
    Lemma~\ref{lemma:integral-operator},
    \refitem{item:integral-operator-makes-smooth} and show that
    \eqref{eq:integral-operator} yields a smooth section $\mathcal{K}
    \varphi \in \Secinfty(E^*\at{U'})$. Its restriction to $U^\cl$ is
    then in $\Sec(E^*\at{U^\cl})$ by the very definition, see
    Definition~\ref{definition:differentiable-on-closed-set}. Moreover,
    the restriction to the open $U$ is of course still smooth. For the
    second part it suffices to show \eqref{eq:intop-to-smooth}. But
    clearly
    \[
    \seminorm[K,k](\mathcal{K} \varphi)
    \leq \vol(U^\cl) \seminorm[K \times U^\cl, k](K)
    \seminorm[U^\cl,0](\varphi)
    \]
    as in Lemma~\ref{lemma:integral-operator},
    \refitem{item:integral-operator-is-cont}. But then the continuity
    is clear by the definition of the locally convex and Banach
    topologies of $\Sec[]_b(E^*\at{U})$, $\Sec(E^*\at{U^\cl})$ and
    $\Secinfty(E^*\at{U})$, respectively.
\end{proof}

We apply this lemma now to the Green functions $F^\pm_U(p)$ in the
following way.
\begin{lemma}
    \label{lemma:seminorm-estimates-for-inverted-intop}
    Let $U \subseteq U^\cl \subseteq U' \subseteq M$ be as in
    Section~\ref{satz:LocalFundamentalSolution} with $U$ small enough
    and let $\mathcal{K}^\pm_U$ be the integral operator from
    \eqref{eq:the-integral-operator}.
    \begin{lemmalist}
    \item \label{item:cont-estimate-for-inverted-intop-1} For every $k \in
        \mathbb{N}_0$ there is a $c > 0$ such that for $\varphi \in
        \Sec[]_b(E^*\at{U'})$ we have
        \begin{equation}
            \label{eq:cont-estimate-for-inverted-intop-1}
            \seminorm[U^\cl,k]
            \left(
                \left(
                    \left(
                        \id + \mathcal{K}^\pm_U
                    \right)^{-1}
                    \circ
                    \mathcal{K}^\pm_U
                \right)
                (\varphi)
            \right)
            \leq c \seminorm[U^\cl,0](\varphi).
        \end{equation}
    \item \label{item:cont-estimate-for-inverted-intop-2} For $\varphi
        \in \Sec(E^*\at{U'})$ there is a $\widetilde{c} > 0$ such that
        \begin{equation}
            \label{eq:cont-estimate-for-inverted-intop-2}
            \seminorm[U^\cl,k]
            \left(
                \left(\id + \mathcal{K}^\pm_U\right)^{-1}
                (\varphi\at{U^\cl})
            \right)
            \leq
            \widetilde{c} \seminorm[U^\cl,k](\varphi).
        \end{equation}
    \end{lemmalist}
\end{lemma}
\begin{proof}
    From Corollary~\ref{corollary:kernel-of-geometric-series} we know
    that the operator $\left(\id + \mathcal{K}^\pm_U\right)^{-1} \circ
    \mathcal{K}^\pm_U$ has a smooth kernel in $\Secinfty(E^* \extensor
    E \at{U' \times U'})$. Thus the previous
    Lemma~\ref{lemma:properties-of-intop-with-smooth-kernel},
    \refitem{item:intop-is-continuous} applies and
    \eqref{eq:intop-cont-estimate} gives
    \eqref{eq:cont-estimate-for-inverted-intop-1}. For the second part
    we note that
    \[
    \left(\id + \mathcal{K}^\pm_U\right)^{-1}
    \left(\varphi \at{U^\cl}\right) \At{U^\cl}
    = \varphi \at{U^\cl}
    - \left(\id + \mathcal{K}^\pm_U\right)^{-1}
    \circ \mathcal{K}^\pm_U (\varphi)
    \at{U^\cl},
    \]
    as we already argued in the proof of
    Lemma~\ref{lemma:integral-operator}. But then
    \begin{align*}
        \seminorm[U^\cl,k]
        \left(
            \left(\id + \mathcal{K}^\pm_U\right)^{-1}
            \left(\varphi \at{U^\cl}\right)
        \right)
        =  \seminorm[U^\cl,k]
        \left(
            \varphi
            - \left(\id + \mathcal{K}^\pm_U\right)^{-1}
            \circ \mathcal{K}^\pm_U
            (\varphi)
        \right)
        \leq \seminorm[U^\cl,k] (\varphi)
        + c \seminorm[U^\cl,0](\varphi)
    \end{align*}
    with $c > 0$ from
    \eqref{eq:cont-estimate-for-inverted-intop-1}. Since
    $\seminorm[U^\cl,k](\varphi) \geq \seminorm[U^\cl,0](\varphi)$ we
    take $\widetilde{c} = 1+c$ to obtain
    \eqref{eq:cont-estimate-for-inverted-intop-2}.
\end{proof}

The importance in the above estimates is that we can control the
``loss of derivatives'': the operator $\left(\id +
    \mathcal{K}^\pm_U\right)^{-1}$ is not loosing orders of
differentiation while $\left(\id + \mathcal{K}^\pm_U\right)^{-1} \circ
\mathcal{K}^\pm_U$ is even gaining smoothness in
\eqref{eq:cont-estimate-for-inverted-intop-1}. We combine this now
with the properties of $\widetilde{\mathcal{R}}^\pm_{U'}$ from
Proposition~\ref{proposition:estimates-and-cont-of-pairing-with-approx-solution}
to obtain the following property of the operator $F^\pm_U$:
\begin{proposition}
    \label{proposition:cont-estimate-of-green-function}
    Let $U \subseteq U^\cl \subseteq U^\cl$ be as before and let
    $F^\pm_U = \left(\id + \mathcal{K}^\pm_U\right)^{-1} \circ
    \widetilde{\mathcal{R}}^\pm_{U'}(\argument)$ be the operator as in
    Definition~\ref{definition:local-fundamental-solution}. Then for
    all compacta $K \subseteq U$ and all $k \in \mathbb{N}_0$ we have
    a $c_{K,k} > 0$ such that
    \begin{equation}
        \label{eq:cont-estimate-of-green-function}
        \seminorm[U^\cl,k](F^\pm_U(\varphi))
        \leq c_{K,k} \seminorm[K,k+n+1](\varphi)
    \end{equation}
    for all $\varphi \in \Secinfty_K(E^*\at{U})$.
\end{proposition}
\begin{proof}
    We know already from the proof of
    Theorem~\ref{theorem:local-fundamental-solution} that the operator
    $F^\pm_U$ is continuous but
    \eqref{eq:cont-estimate-of-green-function} gives a more precise
    statement of this. We have by
    \eqref{eq:cont-estimate-for-inverted-intop-2} and
    \eqref{eq:cont-estimate-for-pairing-with-approx-solution}
    \begin{align*}
        \seminorm[U^\cl,k](F^\pm_U(\varphi))
        &= \seminorm[U^\cl,k]
        \left(
            \left(\id + \mathcal{K}^\pm_U\right)^{-1}
            \left(\widetilde{\mathcal{R}}^\pm_{U'}\right) (\varphi)
        \right) \\
        &\leq \widetilde{c} \seminorm[U^\cl, k]
        \left(\widetilde{\mathcal{R}}^\pm_{U'} (\varphi)\right) \\
        &\leq \widetilde{c}\, c_{K, U^\cl,k+n+1}
        \seminorm[K,k+n+1](\varphi),
    \end{align*}
    which is \eqref{eq:cont-estimate-of-green-function}.
\end{proof}

\begin{corollary}
    \label{corollary:green-function-on-Ck}
    The operator $F^\pm_U$ has a continuous extension to an operator
    \begin{equation}
        \label{eq:LocalGreen-function-on-Ck}
        F^\pm_U: \Sec[k+n+1]_0(E^* \at{U}) \longrightarrow
        \Sec(E^*\at{U})
    \end{equation}
    for all $k \geq 0$, and the estimate
    \eqref{eq:cont-estimate-of-green-function} also holds for $\varphi
    \in \Sec[k+n+1]_K(E^* \at{U})$.
\end{corollary}
\begin{proof}
    The estimate \eqref{eq:cont-estimate-of-green-function} for all
    compact subsets $K \subseteq U$ is just the continuity of
    $F^\pm_U$ in the $\Fun[k+n+1]_0$- and $\Fun$-topology. Thus by the
    usual density argument we have a unique continuous extension
    \eqref{eq:LocalGreen-function-on-Ck} still obeying the estimate
    \eqref{eq:cont-estimate-of-green-function}.
\end{proof}

As usual we can also dualize \eqref{eq:LocalGreen-function-on-Ck} and
get a weak$^*$ continuous map
\begin{equation}
    \label{eq:dualizing-of-green-Ck-extension}
    (F^\pm_U)': \Sec[-k]_0(E\at{U}) \longrightarrow
    \Sec[-k-n-1](E\at{U}),
\end{equation}
again for all $k \geq 0$. Recall that by
Remark~\ref{remark:generalized-sections-for-fixed-density} the
topological dual spaces of $\Sec(E^*\at{U})$ and $\Sec_0(E^*\at{U})$
can be identified with $\Sec[-k]_0(E\at{U})$ and $\Sec[-k](E\at{U})$,
respectively. Note again, that $\Sec[-0](E\at{U})$ are \emph{not} just
the continuous sections $\Sec[0](E\at{U})$. The importance of
Proposition~\ref{proposition:cont-estimate-of-green-function} and
Corollary~\ref{corollary:green-function-on-Ck} is that we only loose a
fixed amount of derivatives under $F^\pm_U$. In this sense the order
of the map $F^\pm_U$ is globally bounded by $n+1$.

In general, a continuous operator $A: \Secinfty_0(E^*) \longrightarrow
\Secinfty(E^*)$ gives a dual operator $A': \Sec[-\infty]_0(E)
\longrightarrow \Sec[-\infty](E)$ as we did this above for $A =
F^\pm_U$. Now this operator $A'$ does not necessarily map
$\Secinfty_0(E) \subseteq \Sec[-\infty]_0(E)$ into $\Secinfty(E)
\subseteq \Sec[-\infty](E)$. For this additional property, $A$ needs
to be a ``symmetric'' operator for the natural pairing. We will now
show this feature for $F^\pm_U$. We consider the following
situation. Let $v$ be a distributional section of $E$ with compact
support in $U$ as before but we assume that $v$ is actually a
$\Fun[\ell]$-section with $\ell \in \mathbb{N}_0$. Then for a test
section $\varphi \in \Secinfty_0(E^*\at{U})$ we have
\begin{align}
    \label{eq:dualized-green-function-on-Cl-distribution}
    (F^\pm_U)'(v) (\varphi)
    = v( F^\pm_U (\varphi))
    = \int_U v(p) \cdot F^\pm_U(\varphi)\at{p} \: \mu_g(p)
    = \int_U v(p) \cdot F^\pm_U(p)(\varphi) \: \mu_g(p),
\end{align}
according to our convention for the pairing of
$\Sec[-\infty]_0(E\at{U})$ and $\Secinfty(E^*\at{U})$. For the Riesz
distributions we already had some symmetry properties as explained in
Proposition~\ref{proposition:symmetry-of-riesz-on-U}. Thus the
question is whether we can extend this to $F^\pm_U$ as well and move
$F^\pm_U$ to the other side in the natural pairing
\eqref{eq:dualized-green-function-on-Cl-distribution}. We start with
the corresponding symmetry property of
$\widetilde{\mathcal{R}}^\pm_{U'}$.
\begin{lemma}
    \label{lemma:1}
    Let $\widetilde{\mathcal{R}}^\pm_{U'}$ be as before and let $k \in
    \mathbb{N}_0$. The for all $u \in \Sec[k+n+1]_0(E\at{U'})$ we have
    \begin{lemmalist}
    \item \label{item:dualizing-of-approx-solution}
        $\widetilde{\mathcal{R}}^\pm_{U'}$ dualizes to a weak$^*$
        continuous linear map
        \begin{equation}
            \label{eq:dualizing-of-approx-solution}
            \left(\widetilde{\mathcal{R}}^\pm_{U'}\right)':
            \Sec[-k]_0(E\at{U'})
            \longrightarrow
            \Sec[-k-n-1](E\at{U'}).
        \end{equation}
    \item \label{item:explicit-formula-for-approx-solution-dual} We
        have $(\widetilde{\mathcal{R}}^\pm_{U'})'(u) \in
        \Sec(E\at{U'})$ explicitly given by
        \begin{equation}
            \label{eq:explicit-formula-for-approx-solution-dual}
            \left(
                \left(\widetilde{\mathcal{R}}^\pm_{U'}\right)'(u)
            \right)
            (q)
            = \sum_{j=0}^\infty
            \left(
                \widetilde{V}^j_q
            \right)^\Trans
            R^\mp_{U'}(2+2j,q)(u),
        \end{equation}
        where $\widetilde{V}^j = V^j$ for $j \leq N-1$ and
        $\widetilde{V}^j = V^j \chi(\frac{\eta}{\epsilon_j})$ for $j
        \geq N$ for abbreviation and
        \begin{equation}
            \label{eq:canonical-transpotisition}
            ^\Trans:
            \Secinfty\left(E^* \extensor E \at{U' \times U'}\right)
            \longrightarrow
            \Secinfty\left(E \extensor E^* \at{U' \times U'}\right)
        \end{equation}
        is the canonical transposition also flipping the arguments.
    \end{lemmalist}
\end{lemma}
\begin{proof}
    The first part is clear since $\widetilde{\mathcal{R}}^\pm_{U'}$
    is a continuous linear map
    \[
    \widetilde{\mathcal{R}}^\pm_{U'}:
    \Sec[k+n+1]_0(E^*\at{U}) \longrightarrow \Sec(E^*\at{U})
    \]
    by Remark~\ref{remark:extension-of-approx-pairing-to-nonsmooth}
    and the duals are just given by $\Sec[-k-n-1](E\at{U'})$ and
    $\Sec[-k]_0(E\at{U'})$ respectively. Thus it remains to evaluate
    $(\widetilde{\mathcal{R}}^\pm_{U'})'(u)$. Since we can interpret
    $u$ as distributional section of any order we want, it is
    sufficient to evaluate the result on smooth test sections $\varphi
    \in \Secinfty_0(E^*\at{U'})$ since they will by dense in every
    other test section space $\Sec[\ell]_0(E^*\at{U'})$. Thus we
    compute
    \begin{align*}
        (\widetilde{\mathcal{R}}^\pm_{U'})'(u)(\varphi)
        &= u \left(
            \widetilde{\mathcal{R}}^\pm_{U'} (\varphi)
        \right) \\
        &= \int_{U'} u(p) \cdot
        \widetilde{\mathcal{R}}^\pm_{U'}(\varphi)\at{p} \: \mu_g(p) \\
        &= \int_{U'} u(p) \cdot
        \widetilde{\mathcal{R}}^\pm_{U'}(p)(\varphi) \: \mu_g(p) \\
        &= \int_{U'} u(p) \cdot
        \sum_{j=0}^{N-1} V_p^j R^\pm_{U'}(2+2j,p)(\varphi) \: \mu_g(p) \\
        &\quad + \int_{U'} u(p) \cdot
        \left(
            \sum_{j=N}^\infty V_p^j
            \chi \left(\frac{\eta_p}{\epsilon_j}\right)
            R^\pm_{U'}(2+2j,p)(\varphi)
        \right) \mu_g(p).
    \end{align*}
    We set $\widetilde{V}^j_p = V^j_p$ for $j \leq N-1$ and
    $\widetilde{V}^j_p = V^j_p \chi(\frac{\eta_p}{\epsilon_j})$ for $j
    \geq N$ to abbreviate the single terms. Then we have
    \begin{align*}
        (\widetilde{\mathcal{R}}^\pm_{U'})'(u)(\varphi)
        &= \sum_{j=0}^{N+k-1} \int_{U'} u(p) \cdot \widetilde{V}^j_p
        R^\pm_{U'}(2+2j,p)(\varphi) \: \mu_g(p) \\
        &\quad + \sum_{j=N+k}^\infty \int_{U'} u(p) \cdot
        \int_{U'} \widetilde{V}^j_p(q) R^\pm_{U'}(2+2j,p)(q) \cdot
        \varphi(q) \: \mu_g(q) \mu_g(p),
        \tag{\smiley}
    \end{align*}
    since in the second series we have $\Fun$-convergence by
    Proposition~\ref{proposition:convergence-of-tail-series},
    \refitem{item:Ck-convergence-of-cutoffed-series} and compact
    support. Thus the series can indeed be taken outside the
    integrals. For the first $N+k$ terms we use
    Proposition~\ref{proposition:symmetry-of-riesz-on-U} in a slightly
    more general setting: the function
    \[
    (p,q)
    \; \mapsto \;
    u(p) \cdot \widetilde{V}^j_p(p,q) \cdot \varphi(q)
    \]
    is compactly supported in $U' \times U'$ but only $\Fun[k+n+1]$
    instead of $\Cinfty$. However, the involved Riesz distributions
    are all of order $\leq n+1$ whence we still can apply
    Proposition~\ref{proposition:symmetry-of-riesz-on-U},
    \refitem{item:symmetry-of-riesz-on-U-as-distribution}, e.g. by
    arguing with the usual density trick. This gives
    \begin{align*}
        \sum_{j=0}^{N+k-1}& \int_{U'} R^\pm_{U'}(2+2j,p)
        \left(
            q \; \mapsto \;
            u(p) \cdot \widetilde{V}^j(p,q) \cdot \varphi(q)
        \right)
        \mu_g(p) \\
        &= \sum_{j=0}^{N+k-1} \int_{U'} R^\mp_{U'}(2+2j,q)
        \left(
            p \; \mapsto \;
            u(p) \cdot \widetilde{V}^j(p,q)
        \right) \cdot \varphi(q))
        \: \mu_g(q) \\
        &= \sum_{j=0}^{N+k-1} \int_{U'} R^\mp_{U'}(2+2j,q)
        \left(
            p \; \mapsto \;
            u(p) \cdot \widetilde{V}^j(p,q)
        \right)
        \cdot \varphi(q) \: \mu_g(q).
    \end{align*}
    Now it is useful to consider the transposition map
    \[
    ^\Trans:
    \Secinfty\left(E* \extensor E \at{U' \times U'}\right)
    \longrightarrow
    \Secinfty\left(E \extensor E^* \at{U' \times U'}\right),
    \]
    defined in the usual way by exchanging the order of arguments
    $(p,q) \leftrightarrow (q,p)$ and the $E$- and $E*$-parts,
    respectively. Thus we have
    \begin{align*}
        \sum_{j=0}^{N+k-1}& \int_{U'} R^\mp_{U'}(2+2j,q)
        \left(
            p \; \mapsto \;
            u(p) \cdot \widetilde{V}^j(p,q)
        \right)
        \cdot \varphi(q) \: \mu_g(q) \\
        &= \sum_{j=0}^{N+k-1} \int_{U'} R^\mp_{U'}(2+2j,q)
        \left(
            p \; \mapsto \;
            \widetilde{V}^j{}^\Trans(p,q) \cdot u(p)
        \right)
        \cdot \varphi(q) \: \mu_g(q) \\
        &= \sum_{j=0}^{N+k-1} \int_{U'}
        \widetilde{V}^j{}^\Trans R^\mp_{U'}(2+2j,\argument)(u) \at{q}
        \cdot \varphi(q) \: \mu_g(q).
    \end{align*}
    By the smoothness of $\widetilde{V}^j{}^\Trans$ and
    Proposition~\ref{proposition:riesz-dependence-on-base-point-p-prime}
    we conclude that the section
    \[
    q \; \mapsto \; \sum_{j=0}^\infty
    \left(
        (\widetilde{V}^j R^\mp_{U'}(2+2j,\argument))^\Trans (u)
    \right) (q)
    \]
    is actually a $\Fun$-section of $E$ on $U'$ since $u$ is
    $\Fun[k+n+1]$. It remains to consider the second part of
    (\smiley). First we again use
    Proposition~\ref{proposition:symmetry-of-riesz-on-U},
    \refitem{item:symmetry-of-riesz-on-U-as-function} to move
    $R^\pm_{U'}(2+2j,p)$ to the other side. Afterwards we exchange the
    order of integration and summation back by the same
    $\Fun$-convergence yielding eventually
    \begin{align*}
        \sum_{j=N+k}^\infty& \int_{U'} \int_{U'} u(p) \cdot
        \widetilde{V}^j_p(q) R^\pm_{U'}(2+2j,p)(q) \cdot \varphi(q)
        \: \mu_g(p) \mu_g(q) \\
        &= \sum_{j=N+k}^\infty \int_{U'} \int_ {U'}
        R^\mp_{U'}(2+2j,q)(p) u(p) \cdot \widetilde{V}^j(p,q) \cdot
        \varphi(q)
        \: \mu_g(p) \mu_g(q) \\
        &= \int_{U'} \varphi(q) \cdot
        \int_{U'}
        \left(
            \sum_{j=N+k}^\infty (\widetilde{V}^k_q)^\Trans
            R^\mp_{U'}(2+2j,q)
        \right)(p) \cdot u(p)
        \: \mu_g(p) \mu_g(q).
    \end{align*}
    The series still converges \emph{in the $\Fun$-topology} as we
    only switched the labels. Thus the inner integrand is a
    $\Fun$-section on $U' \times U'$ being paired with a compactly
    supported $\Fun[k+n+1]$-section $u$. This gives still a
    $\Fun$-section on $U'$ which is then paired with the remaining
    $\varphi$. We conclude that
    \begin{align*}
        \left(
            (\widetilde{\mathcal{R}}^\pm_{U'})'(u)
        \right)(\varphi)
        &= \int_{U'}
        \left(
            \sum_{j=0}^\infty (\widetilde{V}^j)^\Trans
            R^\mp_{U'}(2+2j,\argument)(u)
        \right)(q) \cdot \varphi(q)
        \: \mu_g(q)
    \end{align*}
    with a $\Fun$-section
    \[
    \left(
        (\widetilde{\mathcal{R}}^\pm_{U'})'(u)
    \right)(q)
    = \sum_{j=0}^\infty (\widetilde{V}^j_q)^\Trans
    R^\mp_{U'}(2+2j,q)(u)
    \]
    as claimed.
\end{proof}

\begin{remark}
    \label{remark:dual-of-riesz}
    The Riesz distributions $R^\pm_{U'}(\alpha,p)$ enjoy the symmetry
    property $R^\pm_{U'}(\alpha,p)(q) = R^\mp_{U'}(\alpha,q)(p)$ as
    soon as $\RE(\alpha) > n$. For all $\alpha \in \mathbb{C}$, the
    correct analog of this symmetry was obtained in
    Proposition~\ref{proposition:symmetry-of-riesz-on-U},
    \refitem{item:symmetry-of-riesz-on-U-as-distribution}. Thus
    extending the transposition $^\Trans$ from smooth to continuous or
    even distributional sections we have
    \begin{equation}
        \label{eq:symmetry-of-ries-distri-on-U}
        \left(R^\pm_{U'}\right)^\Trans
        = R^\mp_{U'}
    \end{equation}
    in the sense of
    Proposition~\ref{proposition:symmetry-of-riesz-on-U},
    \refitem{item:symmetry-of-riesz-on-U-as-distribution}. Moreover,
    since in the series
    \eqref{eq:explicit-formula-for-approx-solution-dual} we have the
    ``same'' coefficients as for the original series defining
    $\widetilde{\mathcal{R}}^\pm_{U'}$ only at flipped points, we get
    the same sort of estimates and convergence results. In particular
    we have
    \begin{equation}
        \label{eq:dual-of-ries-is-transpose}
        \left(\widetilde{\mathcal{R}}^\pm_{U'}\right)'
        =
        \left(\widetilde{\mathcal{R}}^\mp_{U'}\right)^\Trans
    \end{equation}
    on distributional sections which are at least $\Fun[n+1]$. This
    allows to efficiently compute
    $\left(\widetilde{\mathcal{R}}^\pm_{U'}\right)'(u)$ for $u \in
    \Sec[n+1]_0(E\at{U'})$ by means of the nicely convergent series
    \eqref{eq:explicit-formula-for-approx-solution-dual} or
    \eqref{eq:dual-of-ries-is-transpose}.
\end{remark}

\begin{corollary}
    \label{corollary:dual-of-riesz}
    Let $u \in \Secinfty_0(E\at{U'})$ then
    $\left(\widetilde{\mathcal{R}}^\pm_{U'}\right)'(u) \in
    \Secinfty(E\at{U'})$.
\end{corollary}
\begin{corollary}
    \label{corollary:dual-of-riesz-on-Fun-sections}
    Let $k \in \mathbb{N}_0 \cup \{ +\infty \}$ and $u \in
    \Sec[k+n+1]_0(E\at{U'})$. Then the series
    \eqref{eq:explicit-formula-for-approx-solution-dual} converges in
    the $\Fun$-topology.
\end{corollary}
\begin{proof}
    This follows analogously to the statements for
    $\widetilde{\mathcal{R}}^\pm_{U'}$ as in
    Proposition~\ref{proposition:convergence-of-tail-series}: the
    finitely many terms with $j \leq N+k-1$ are already $\Fun$ by
    themselves and the remaining sum converges in $\Fun$ \emph{before}
    applying to $u$ on $U' \times U'$. Then the integration over $p$
    together with the compactly supported $u$ can be exchanged with
    the summation by the usual arguments. It gives then the
    $\Fun$-convergence on $U'$.
\end{proof}

\begin{figure}
    \centering
    \input{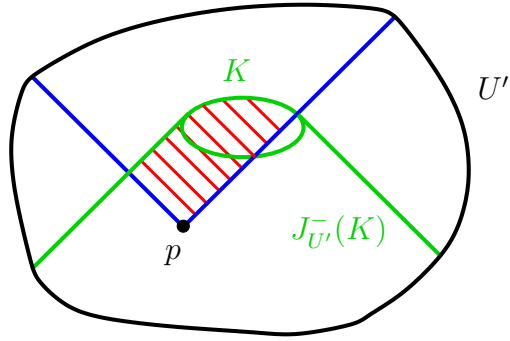}
    \caption{\label{fig:future-of-point-intersect-with-past-of-compact}%
      The intersection of the future of a point $p$ with the past of a
      compactum $K$, all in a geodesically convex $U'$.
    }
\end{figure}

We can use the lemma also to extend $\widetilde{\mathcal{R}}^\pm_{U'}$
as well as its dual $\left(\widetilde{\mathcal{R}}^\pm_{U'}\right)'$
and $\left(\widetilde{\mathcal{R}}^\mp_{U'}\right)^\Trans$ to some
more general test sections and distributions with not necessarily
compact support. We consider the following situation: Let $K \subseteq
U'$ be compact, then the intersection $J^+_{U'}(p) \cap J^-_{U'}(K)$
is still compact since $U'$ is geodesically convex, see
Figure~\ref{fig:future-of-point-intersect-with-past-of-compact}.
In fact, also the intersection $J^+_{U'}(L) \cap J^-_{U'}(K)$ is
compact for another compactum $L \subseteq U'$. Suppose $\supp \varphi
\subseteq J^-_{U'}(K)$ for a test section $\varphi \in
\Sec(E^*\at{U'})$ with not necessarily compact support. Then for every
$j$ and every $p \in U'$ the overlap
\begin{align*}
    \supp \left(
        \widetilde{V}^j_p R^+(2+2j,p)
    \right) \cap \supp \varphi
    &\subseteq
    \supp R^+_{U'}(2+2j,p) \cap \supp \varphi \\
    &\subseteq
    J^+_{U'}(p) \cap \supp \varphi \\
    &\subseteq J^+_{U'}(p) \cap J^-_{U'}(K)
\end{align*}
is compact. Thus $\widetilde{V}^j_p R^+_{U'}(2+2j,p)(\varphi)$ is
defined by Proposition~\ref{proposition:supp-overlap-is-compact} in a
non-ambiguous way. By the same argument, also
$\widetilde{\mathcal{R}}^+_{U'}(\varphi)$ is well-defined. Moreover,
since for $p \in L$ the support of $\widetilde{V}^j_p R^+(2+2j,p)$ has
still compact overlap with $\supp \varphi$ we can replace $\varphi$ by
some $\chi \varphi$ as in the proof of
Proposition~\ref{proposition:supp-overlap-is-compact} and get the same
convergence results of the series
\begin{equation}
    \label{eq:approx-solution-for-supp-overlap}
    \widetilde{\mathcal{R}}^+_{U'}(p)(\varphi)
    = \sum_{j=0}^\infty \widetilde{V}^j_p R^+(2+2j,p)(\varphi)
\end{equation}
as for compactly supported $\varphi$. In conclusion, this gives a
$\Fun$-convergence if $\varphi$ is of class $\Fun[k+n+1]$ for all $k
\in \mathbb{N}_0 \cup \{+\infty\}$. With the same argument, also the
series $\left(\widetilde{\mathcal{R}}^-_{U'}\right)^\Trans$
converges. Here of course we need $u \in \Sec[k+n+1](E\at{U'})$ with
$\supp u \subseteq J^+_{U'}(K)$ to make the series
\begin{equation}
    \label{eq:dual-series-on-supp-overlap}
    (\widetilde{\mathcal{R}}^-_{U'})^\Trans (u)
    = \sum_{j=0}^\infty (\widetilde{V}^j_{\argument})^\Trans
    R^-_{U'}(2+2j,\argument)(u)
\end{equation}
converge in the $\Fun$-topology. We collect these results in the
following lemma:
\begin{lemma}
    \label{lemma:convergence-on-supp-overlap}
    Let $K \subseteq U'$ be compact and $k \in \mathbb{N}_0 \cup
    \{+\infty\}$.
    \begin{lemmalist}
    \item Assume $u \in \Sec[k+n+1](E\at{U'})$ has support in
        $J^\mp_{U'}(K)$. Then
        \begin{equation}
            \label{eq:dual-series-for-supp-oerlap-converges-Ck}
            (\widetilde{\mathcal{R}}^\mp_{U'})^\Trans (u)
            = \sum_{j=0}^\infty (\widetilde{V}^j_{\argument})^\Trans
            R^\mp_{U'}(2+2j,\argument)(u)
        \end{equation}
        converges in the $\Fun$-topology.
    \item Assume $\varphi \in \Sec[k+n+1](E^*\at{U'})$ has support in
        $J^\mp_{U'}(K)$. Then
        \begin{equation}
            \label{eq:approx-sol-on-supp-overlap-converges-Ck}
            \widetilde{\mathcal{R}}^\pm_{U'}(p)(\varphi)
            = \sum_{j=0}^\infty \widetilde{V}^j_p
            R^\pm(2+2j,p)(\varphi)
        \end{equation}
        converges in the $\Fun$-topology.
    \end{lemmalist}
\end{lemma}

We can now study the dual of $F^\pm_U$ under the assumption that $U
\subset U^\cl \subset U'$ is \emph{causal} in order to have good
support properties of the integral operator $\mathcal{K}^\pm_U$.
\begin{lemma}
    \label{lemma:dual-of-fundamental-solution}
    Let $u \in \Secinfty_0(E\at{U})$. Then
    \begin{equation}
        \label{eq:dual-of-fundamental-solution}
        (F^\pm_U)'(u)
        =
        \left(\widetilde{\mathcal{R}}^\mp_U\right)^\Trans
        \left(
            q \; \mapsto \;
            u(q) - \int_{U} u(p) \cdot L^\pm_U(p,q) \: \mu_g(p)
        \right)
    \end{equation}
    with $L^\pm_U$ being the smooth integral kernel of $\left(\id +
        \mathcal{K}^\pm_U\right)^{-1} \circ \mathcal{K}^\pm_{U}$. Thus
    $(F^\pm_U)'(u) \in \Secinfty(E\at{U})$.
\end{lemma}
\begin{proof}
    For $\varphi \in \Secinfty_0(E^*\at{U})$ we have to evaluate the
    pairing
    \begin{align*}
        (F^\pm_U)'(u)(\varphi)
        &= u\left(F^\pm_U(\varphi)\right) \\
        &= \int_U u(p) \cdot F^\pm_U(\varphi)\at{p} \: \mu_g(p) \\
        &= \int_U u(p) \cdot
        \left(\id + \mathcal{K}^\pm_{U}\right)^{-1}
        \left(\widetilde{\mathcal{R}}^\pm_U(\argument)(\varphi)\right)\At{p}
        \: \mu_g(p) \\
        &= \int_U u(p)
        \left(
            \widetilde{\mathcal{R}}^\pm_U(p)(\varphi)
            - \left(\id + \mathcal{K}^\pm_U\right)^{-1}
            \circ \mathcal{K}^\pm_U
            \left(\widetilde{\mathcal{R}}^\pm_U(\argument)(\varphi)\right)\At{p}
        \right)
        \mu_g(p).
    \end{align*}
    Now $\left(\id + \mathcal{K}^\pm_U\right)^{-1} \circ
    \mathcal{K}^\pm_U$ is again an integral operator whose kernel is
    smooth and given by the truncated geometric series as in
    Corollary~\ref{corollary:kernel-of-geometric-series}. Thus denote
    its kernel by $L^\pm_U \in \Secinfty(E^* \extensor E \at{U' \times
      U'})$, noting that even though we only integrate over $U^\cl$
    the kernel has a smooth continuation to $U' \times U'$. Since we
    integrate at least continuous functions and sections over compact
    sets $U^\cl$ and $U^\cl \times U^\cl$, respectively, we can
    exchange the orders of integration and obtain
    \begin{align*}
        (F^\pm_U)'(u)(\varphi)
        &= \int_U u(p) \cdot \widetilde{\mathcal{R}}^\pm_U(\varphi)(p)
        \: \mu_g(p)
        - \int_{U'} \int_{U'} u(p) \cdot L^\pm_U(p,q) \cdot
        \widetilde{\mathcal{R}}^\pm_U(\varphi)(q)
        \: \mu_g(q) \mu_g(p) \\
        &= \int_U
        \left(
            u(q) - \int_U u(p) \cdot L^\pm_U(p,q)
            \: \mu_g(p)
        \right)
        \cdot \widetilde{\mathcal{R}}^\pm_U(\varphi)(q)
        \: \mu_g(q) \\
        &= \int_U v(q) \cdot \widetilde{\mathcal{R}}^\pm_U(\varphi)(q)
        \: \mu_g(q),
        \tag{$*$}
    \end{align*}
    with
    \[
    v(q) = u(q) - \int_U u(p) \cdot L^\pm_U(p,q) \mu_g(p).
    \]
    Now the second term in $v$ is smooth and has a smooth extension to
    $U'$. The first contribution $u$ is compactly supported in $U$ and
    smooth whence it also has a smooth extension to $U'$: we conclude
    $v \in \Secinfty(E\at{U})$. We claim that in ($*$) we are allowed
    to move $\widetilde{\mathcal{R}}^\pm_U$ from $\varphi$ to $v$ on
    the other side of the pairing. Indeed, by the causal properties of
    $L^\pm_U$ according to
    Lemma~\ref{lemma:supps-of-int-kernels-are-future-past-stretched}
    we know
    \[
    \supp L^\pm_U
    \subseteq \left\{ (p,q) \; \big| \; q \in J^+_{U'}(p) \right\}
    \subseteq U' \times U'.
    \]
    Thus when restricting to $U^\cl$ and using that $U$ is causal we
    see that the integrand $u(p) \cdot L^\pm_U(p,q)$ is possibly
    non-trivial only for $p \in \supp v$ and $q \in J^\pm_{U'}(p) \cap
    U^\cl = J^\pm_{U^\cl}(p)$. But this is equivalent to $p \in
    J^\mp_{U^\cl}(q)$ and hence the integrand is possibly non-trivial
    only for $\supp u \cap J^\mp_{U^\cl}(q) \neq \emptyset$. In other
    words, $\supp \left(\left(\id + \mathcal{K}^\pm_U\right)^{-1}
        \circ \mathcal{K}^\pm_U \right)(u) \subseteq J^\pm_{U^\cl}(
    \supp u)$. Hence $\supp v \subseteq J^\pm_{U^\cl}(\supp u)$. Note
    that due to the transposed integration this differs from the
    considerations for $L^\pm_U$ acting on $\varphi \in
    \Secinfty_0(E^*\at{U})$. But then expanding the series over $j$ in
    $\widetilde{\mathcal{R}}^\pm_U(\varphi)$ we get
    \begin{align*}
        \int_U v(q) \cdot \widetilde{\mathcal{R}}^\pm_U(\varphi)(q)
        \: \mu_g(q)
        &=
        \int_U v(q) \cdot
        \sum_{j=0}^\infty \widetilde{V}^j_q R^\pm_U(2+2j,q)(\varphi)
        \: \mu_g(q) \\
        &= \sum_{j=0}^\infty \int_U  v(q)
        \cdot \widetilde{V}^j_q R^\pm_U(2+2j,p)(\varphi)
        \: \mu_g(q) \\
        &= \sum_{j=0}^\infty \int_U
        \left(
            \widetilde{V}^j_p{}^\Trans R^\mp_U(2+2j,p)(v)
        \right)
        \cdot \varphi(p)
        \: \mu_g(p) \\
        &= \int_U \sum_{j=0}^\infty \widetilde{V}^j_p{}^\Trans
        R^\mp_U(2+2j,p)(v) \cdot \varphi(p)
        \: \mu_g(p) \\
        &= \int_U (\widetilde{\mathcal{R}}^\mp_U)^\Trans(v)(p) \cdot
        \varphi(p)
        \: \mu_g(p).
    \end{align*}
    Here we used that $\Fun[0]$-convergent series can be exchanges
    with integration over compacta and $R^\pm_U$ can be transposed as
    in Proposition~\ref{proposition:symmetry-of-riesz-on-U},
    \refitem{item:symmetry-of-riesz-on-U-as-distribution} even though
    $v$ has non-compact support: The main point is that the overlap of
    the supports is compact even though $\supp v \subseteq
    J^\pm_{U^\cl}(\supp u)$ typically is non-compact. But then we know
    that the series still converges in the $\Fun[0]$-topology and can
    be moved inside the integral by
    Lemma~\ref{lemma:convergence-on-supp-overlap}.
\end{proof}

\begin{remark}
    \label{remark:dual-fund-on-Ck}
    A careful counting of derivatives shows that the operator
    $\left(\id + \mathcal{K}^\pm_U\right)^{-1}$ does not eat orders of
    differentiation and $(\widetilde{\mathcal{R}}^\mp_U)^\Trans$ needs
    at most $n+1$. Thus we also obtain the statement that
    \begin{equation}
        \label{eq:dual-fund-on-Ck}
        (F^\pm_U)': \Sec[k+n+1]_0(E\at{U}) \longrightarrow \Sec(E\at{U})
    \end{equation}
    holds for all $k \in \mathbb{N}_0 \cup \{+\infty\}$.
\end{remark}
We summarize the result of this section in the following theorem:
\begin{theorem}
    \label{theorem:dual-of-fund-solution-gives-inhom-solution}
    \index{Inhomogeneity!smooth}%
    \index{Wave equation!local solution}%
    Let $k \in \mathbb{N}_0 \cup \{+\infty\}$ and $u \in
    \Sec[k+n+1]_0(E\at{U}$. Then $(F^\pm_U)'(u)$, explicitly given by
    \eqref{eq:approx-sol-on-supp-overlap-converges-Ck}, is a
    $\Fun$-section of $E\at{U}$ with
    \begin{equation}
        \label{eq:dual-of-fund-solution-gives-inhom-solution}
        \supp (F^\pm_U)'(u)
        \subseteq J^\pm_U (\supp u)
        \quad \textrm{and} \quad
        D(F^\pm_U)'(u) = u.
    \end{equation}
    In particular, we have a smooth local solution of the wave
    equation for a smooth and compactly supported inhomogeneity.
\end{theorem}


%% file: chap4.tex
%
%

\chapter{The Global Theory of Geometric Wave Equations}
\label{satz:global-theory}

Since in a time-oriented Lorentz manifold every point has a causal
neighborhood we see from the results in the last chapter that locally
we have advanced and retarded fundamental solutions, i.e. Green
functions, for a given normally hyperbolic differential
operator. Moreover, we have seen how these fundamental solutions can
be used to construct solutions to the inhomogeneous wave equations for
different kinds of inhomogeneities.

The topic in this chapter is now to globalize these results from the
(small) neighborhoods to the whole Lorentz manifold. Here the global
causal structure yields obstructions of various kinds: in general we
will not be able to find global Green functions. Instead, we will need
some assumptions on the global geometry. Here the best situation will
be obtained for globally hyperbolic Lorentz manifolds. On such
spacetimes we can then also formulate and solve the Cauchy problem for
the wave equation. This nice solutions theory allows to treat the wave
equation essentially as an (infinite-dimensional) Hamiltonian
dynamical system. We will illustrate this point of view by determining
the relevant Poisson algebra of observables.

%
%

\section{Uniqueness Properties of Fundamental Solutions}
\label{satz:uniqueness-properties}

\input{unique}

%
%

\section{The Cauchy Problem}
\label{satz:cauchy-problem}

\input{cauchyprob}

%
%

\section{Global Fundamental Solutions and Green Operators}
\label{satz:global-fund-and-green-operators}

\input{green}

%
%

\section{A Poisson Algebra}
\label{satz:poisson-algebra}

\input{poisson}


%% file: unique.tex
%
%

It will be easier to show uniqueness of fundamental solutions than
their actual existence. In the following we will provide criteria
under which there is at most one advanced and one retarded fundamental
solution. In order to treat a rather general situation we first recall
some more refined techniques for the description of the causal
structure.

%
%

\subsection{Time Separation}
\label{satz:time-separation}

The time separation function $\tau$ on $M$ will be the Lorentz
analogue of the \Index{Riemannian distance} $d$. However, in various
aspects it behaves quite differently. It will help us to formulate
appropriate conditions on $M$ to ensure uniqueness properties for the
fundamental solutions. We recall here its definition and some of the
basic properties.
\begin{definition}[Arc length]
    \label{definition:arc-length}
    \index{Arc lenght}%
    Let $\gamma: [a,b] \longrightarrow M$ be a (piecewise) $\Fun[1]$
    curve in a semi-Riemannian manifold $(M,g)$. Then its arc length
    is defined by
    \begin{equation}
        \label{eq:arc-length}
        L(\gamma)
        = \int_a^b
        \sqrt{\left|
              g_{\gamma(t)}(\dot{\gamma}(t),\dot{\gamma}(t))
          \right|
        } \D t.
    \end{equation}
\end{definition}
Clearly, the definition makes sense for piecewise $\Fun[1]$-curves as
well. The following is obvious:
\begin{lemma}
    \label{lemma:arc-lenght-reparametrization-invariant}
    The arc length of a piecewise $\Fun[1]$ curve $\gamma$ is
    invariant under monotonous piecewise $\Fun[1]$
    reparametrization.
\end{lemma}
Unlike in Riemannian geometry, for different points $p$ and $q$ there
may still be curves $\gamma$ joining $p$ and $q$ which have arc length
$0$, namely if $\dot{\gamma}$ is timelike. This makes the concept of a
``distance'' more complicated. One has the following definition:
\begin{definition}[Time separation]
    \label{definition:time-separation}
    \index{Time separation}%
    The time separation function $\tau: M \times M \longrightarrow
    \mathbb{R} \cup \{+\infty\}$ in a time-oriented Lorentz manifold
    $(M,g)$ is defined by
    \begin{equation}
        \label{eq:time-separation}
        \tau(p,q)
        = \sup \left\{
            L(\gamma)
            \; \big| \;
            \gamma
            \;
            \textrm{is a future directed causal curve from}
            \;
            p
            \;
            \textrm{to}
            \;
            q
        \right\}
    \end{equation}
    if $q \in J^+_M(p)$ and $\tau(p,q) = 0$ if $q \notin J^+_M(p)$.
\end{definition}
In contrast to the Riemannian situation where one uses the infimum
over all arc lengths of curves joining $p$ and $q$ to define the
Riemannian distance, the time separation $\tau$ has some new features:
first it is clear that $\tau(p,q) = 0$ may happen even for $p \neq q$;
this is possible already in Minkowski spacetime. Moreover, in general
$\tau(p,q)$ is \emph{not} a symmetric function as it involves the
choice of the time-orientation. Again, this can easily be seen for
Minkowski spacetime and points $p \neq q$ with $q \in I^+_M(p)$.
\begin{figure}
    \centering
    \input{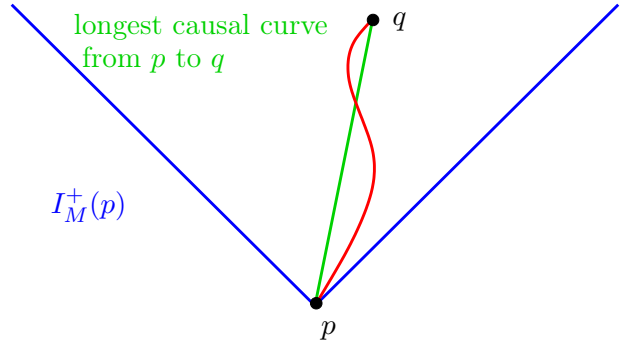}
    \caption{\label{fig:twin-paradoxon}%
      The twin paradoxon}
\end{figure}
\index{Twin paradoxon}%
In this case $\tau(p,q)$ is the Minkowski length of the vector
$\vec{pq} = q - p$. The fact that all other future directed causal
curves from $p$ to $q$ are shorter is the mathematical fact underlying
the so-called \emph{twin paradoxon}. In the more weird examples of
Lorentz manifolds it may happen that $\tau(p,q) = + \infty$ for some
or even all pairs of points: the Lorentz cylinder from
Figure~\ref{fig:periodic-geodesics-no-cauchy-hypersurface} is an
example. By spiralling around the cylinder we find a future directed
timelike geodesic $\gamma$ from $p$ to $q$ of arbitrarily big length
$L(\gamma)$. This already indicates that the points $p$ and $q$ with
$\tau(p,q) = + \infty$ will be responsible for bad behaviour of the
causal structure.

\begin{figure}
    \centering
    \input{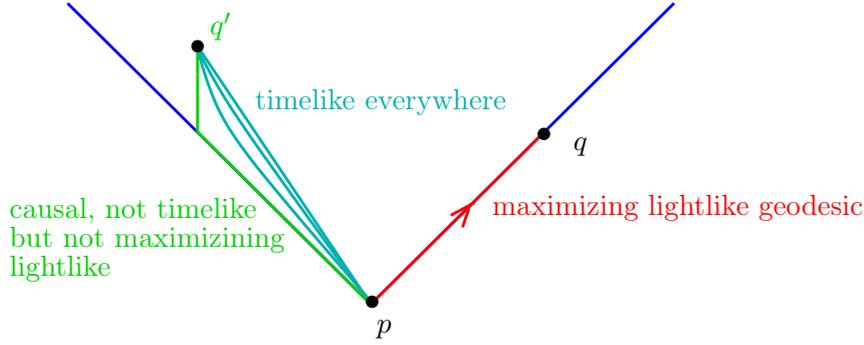}
    \caption{\label{fig:deforming-lightlike-into-timelike}
      Illustration for
      Lemma~\ref{lemma:no-maximizing-lightlike-curve-gives-a-timelike-one}.}
\end{figure}

Recall that a lightlike curve $\gamma$ from $p$ to $q$ is called
\emph{maximizing} if there is no timelike curve from $p$ to $q$. Then
we have the following useful Lemma:
\begin{lemma}
    \label{lemma:no-maximizing-lightlike-curve-gives-a-timelike-one}
    \index{Curve!maximizing lightlike}%
    \index{Curve!timelike}%
    If there is a causal curve $\gamma$ from $p$ to $q$ which is not a
    maximizing lightlike curve then there also exists a timelike curve
    from $p$ to $q$.
\end{lemma}
The proof can be found e.g. in \cite[Thm.~10.51]{oneill:1983a}, see
also the discussion in
\cite[Thm.~2.30]{minguzzi.sanchez:2006a:pre}. The geometric meaning of
this is illustrated in
Figure~\ref{fig:deforming-lightlike-into-timelike}.
In fact, it can be shown that a maximizing lightlike curve is, up to
reparametrization, a lightlike geodesic without conjugate points
between the endpoints. Moreover, one can show that the timelike curve
in the lemma can be chosen arbitrarily close to the original causal
curve $\gamma$. Using this lemma one arrives at the following
properties of the time separation:
\begin{theorem}[Time separation]
    \label{theorem:time-separation}
    \index{Time separation}%
    \index{Timelike loop}%
    Let $(M,g)$ be a time-oriented Lorentz manifold and $p,q,r \in
    M$.
    \begin{theoremlist}
    \item \label{item:timesep-positivity} One has $\tau(p,q) > 0$ iff
        $p \ll q$.
    \item \label{item:timesep-p-p} If there exists a timelike closed
        curve through $p$ then we have $\tau(p,p) = +
        \infty$. Otherwise one has $\tau(p,p) = 0$.
    \item \label{item:timesep-not-symmetric} If $0 < \tau(p,q) <
        +\infty$ then $\tau(q,p) = 0$.
    \item \label{item:reverse-triangle-inequ} For $p \leq q \leq r$
        one has a reverse triangle inequality, i.e.
        \begin{equation}
            \label{eq:revers-triangle-inequ}
            \tau(p, q) + \tau(q, r) \leq \tau(p, r).
        \end{equation}
    \item \label{item:maximizing-geodesic-in-geod-convex-set} Suppose
        $p,q \in U \subseteq M$ with an open geodesically convex
        $U$. If $q \in I^+_U(p)$ then the geodesic $\gamma(t) = \exp_p
        ( t \exp_p^{-1}(q) )$ maximizes the arc length of all causal
        curves from $p$ to $q$ which are entirely in $U$ and
        $\tau_U(p,q) = \sqrt{ g_p( \exp_p^{-1}(q), \exp_p^{-1}(q))}$.
    \item \label{item:timesep-lower-semi-cont} The time separation
        function $\tau$ is lower semi continuous, i.e. for convergent
        sequence $p_n \longrightarrow p$ and $q_n \longrightarrow q$
        one has
        \begin{equation}
            \label{eq:timesep-lower-semi-cont}
            \liminf_{n \rightarrow \infty} \tau(p_n,q_n) = \tau(p,q).
        \end{equation}
    \end{theoremlist}
\end{theorem}
\begin{proof}
    We only sketch the arguments and refer to
    \cite[Chapter~14]{oneill:1983a} or
    \cite[Sect.~2.5]{minguzzi.sanchez:2006a:pre} for details.  If $p
    \ll q$ then there is a timelike future directed curve $\gamma$
    from $p$ to $q$. Thus $L(\gamma) > 0$ and $\tau(p,q) \geq
    L(\gamma)$. Conversely, suppose $\tau(p,q) > 0$ then there is a
    causal future directed curve $\gamma$ from $p$ to $q$ which cannot
    be a lightlike curve as for lightlike curves we have arc length
    $0$. By
    Lemma~\ref{lemma:no-maximizing-lightlike-curve-gives-a-timelike-one}
    we can deform $\gamma$ into a timelike curve whence $p \ll q$
    follows. This gives the first part.  If we have a timelike closed
    loop $\gamma$ through $p$ then clearly $L(\gamma) > 0$. Thus
    winding around more and more often produces $L(\gamma^n) = n
    L(\gamma) \longrightarrow + \infty$, showing $\tau(p,p) =
    +\infty$. Otherwise, there can be at most a maximizing lightlike
    loop through $p$ or $p \notin J^+_M(p)$ at all, by
    Lemma~\ref{lemma:no-maximizing-lightlike-curve-gives-a-timelike-one}. In
    both cases $L(\gamma) = 0$ for all (possibly none at all) curves
    whence $\tau(p,p) = 0$.  The third part is clear since $0 <
    \tau(p,q)$ shows that there is a timelike curve from $p$ to $q$
    and hence $p \ll q$. If also $\tau(q, p) > 0$ then also $q \ll p$
    whence we would obtain a closed timelike loop from $p$ to $p$ with
    non-trivial length $L(\gamma)$. Running around this loop $n$ times
    and then to $q$ gives a timelike curve from $p$ to $q$ with arc
    length at least $n L(\gamma) \longrightarrow + \infty$. This
    contradicts $\tau(p,q) < \infty$, see also
    Figure~\ref{fig:timelike-loop-for-two-points}.
    \begin{figure}
        \centering
        \input{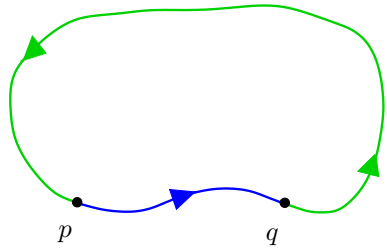}
        \caption{\label{fig:timelike-loop-for-two-points}%
          A timelike loop from $p$ to $q$.
        }
    \end{figure}
    For the fourth part, let $p \leq q \leq r$ be given and let
    $\epsilon > 0$. We find future directed causal curves $\gamma_1$
    from $p$ to $q$ and $\gamma_2$ from $q$ to $p$ with
    \[
    \tau(p,q) < L(\gamma_1) + \epsilon
    \quad
    \textrm{and}
    \quad
    \tau(q,r) < L(\gamma_2) + \epsilon
    \]
    by definition of $\tau$ as supremum. Since $\tau(p,r)$ is clearly
    not less than $L(\gamma_1) + L(\gamma_2)$ as $\gamma_2$ after
    $\gamma_1$ is joining $p$ to $r$, we find
    \[
    \tau(p,r) \geq L(\gamma_1) + L(\gamma_2)
    > \tau(p,q) - \epsilon + \tau(q,r) - \epsilon,
    \]
    whence $\tau(p,r) \geq \tau(p,r) + \tau(q,r) - 2 \epsilon$. Since
    $\epsilon > 0$ was arbitrary, we get the reverse triangle
    inequality.  For the fifth part we refer to e.g. \cite[Lem.~5.33
    and Prop.~5.34]{oneill:1983a}.  Using this we can prove the last
    part as follows: for $\tau(p,q) = 0$ nothing is to be shown. Thus
    consider $0 < \tau(p,q) < +\infty$. Now we fix $\epsilon >
    0$. Then we have to find a neighborhood $U$ of $p$ and a
    neighborhood $V$ of $q$ such that for $p' \in U$ and $q' \in V$ we
    have $\tau(p',q') > \tau(p,q) - \epsilon$. Since $0 < \tau(p,q) <
    +\infty$ we find a timelike curve $\gamma$ from $p$ to $q$ with
    $\tau(p,q) < L(\gamma) + \frac{\epsilon}{3}$ by the first
    part. Now we choose a geodesically convex neighborhood $V'$ of $q$
    and fix a point $q_1 \in V'$ on the curve $\gamma$ such that the
    curve $\gamma$ from $q_1$ to $q$ stays inside $V'$, see
    Figure~\ref{fig:choosing-the-neighborhoods}.
    \begin{figure}
        \centering
        \input{choosing-neighborhoods.\pictype}
        \caption{\label{fig:choosing-the-neighborhoods}%
          Illustration for proof of
          Theorem~\ref{theorem:time-separation},
          \refitem{item:timesep-lower-semi-cont}.
        }
    \end{figure}
    Since the curve $\gamma$ from $q_1$ to $q$ is inside $V'$ and
    timelike, we know from the fifth part that the geodesic segment
    from $q_1$ to $q$ in $V'$ maximizes the arc length and hence it is
    longer (or equal) as the curve $\gamma$ from $q_1$ to $q$. Now we
    fix a smaller neighborhood $V$ of $q$ by the condition that $q'
    \in V$ is in the causal future of $q_1$ and the geodesic
    $c_{q_1,q'}(t) = \exp_{q_1}( t \exp_{q_1}^{-1}(q') )$ from $q_1$
    to $q'$ has arc length
    \[
    L(c_{q_1,q'}) > L(c_{q_1,q}) - \frac{\epsilon}{3}.
    \]
    This is clearly possible as the arc length depends continuously on
    the endpoint. From the two conditions we see that the curve from
    $p$ to $q'$ first along $\gamma$ and then along $c_{q_1, q'}$ has
    arc length $L(\gamma) - \frac{\epsilon}{3}$. An analogous
    construction around $p$ specifies a $p_1$ and the neighborhood
    $U$. Then for $p' \in U$ and $q' \in V$ we have a timelike curve
    by first taking the geodesic from $p'$ to $p_1$ then via $\gamma$
    from $p_1$ to $q_1$ and finally along the geodesic from $q_1$ to
    $q'$. Its arc length is at least $L(\gamma) - 2
    \frac{\epsilon}{3}$. Since $\gamma$ was chosen such that
    $L(\gamma) + \frac{\epsilon}{3} > \tau(p,q)$ we see that the arc
    length of the curve from $p'$ to $q'$ is at least $\tau(p,q) -
    \frac{\epsilon}{3} - 2 \frac{\epsilon}{3} = \tau(p,q) -
    \epsilon$. It follows that for all $p', q'$ in these neighborhoods
    we have $\tau(p',q') \geq \tau(p,q) - \epsilon$. This shows the
    lower semi continuity of $\tau$ for the case $\tau(p,q) <
    \infty$. The construction for $\tau(p,q) = \infty$ proceeds
    analogously by choosing large $L(\gamma)$ and neighborhoods as
    before.
\end{proof}

The following example shows that $\tau$ is \emph{not} continuous in
general:
\begin{example}[Discontinuous time separation]
    \label{example:discont-timesep}
    \index{Time separation!discontinuous}%
    Consider the Minkowski plane with a half axis removed, i.e. $M =
    \mathbb{R}^2 \setminus (-\infty, 0]$, see
    Figure~\ref{fig:discont-time-sep}.
    \begin{figure}
        \centering
        \input{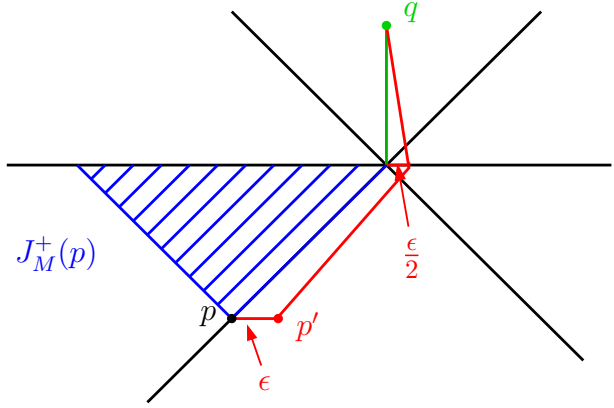}
        \caption{\label{fig:discont-time-sep}%
          A discontinuous time separation.
        }
    \end{figure}
    Let $p = (-1,-1)$ then the causal future $J^+_M(p)$ is the
    triangle under the removed axis. In particular, $q = (1,0)$ is
    \emph{not} in the future of $p$ whence $\tau(p,q) = 0$. However,
    for $p' = (-1, -1+\epsilon)$ with $0 < \epsilon < 1$ the point $q$
    is in $J^+_M(p')$. The broken geodesic from $p'$ to $(0,
    \frac{\epsilon}{2})$ and then from $(0, \frac{\epsilon}{2})$ to
    $q$ are both timelike and the length of the first is
    \[
    L(\gamma_1)
    = \sqrt{1 - (1- \frac{\epsilon}{2})^2}
    = \sqrt{1-1+\epsilon-\frac{\epsilon^2}{4}}
    = \sqrt{\epsilon - \frac{\epsilon^2}{4}}
    \]
    while the length of the second curve is
    \[
    L(\gamma_2) = \sqrt{1 - \frac{\epsilon^2}{4}}.
    \]
    It follows that $\tau(p',q)$ is at least $\sqrt{\epsilon -
      \frac{\epsilon^2}{4}} + \sqrt{1 - \frac{\epsilon^2}{4}}$, whence
    \begin{equation}
        \label{eq:discont-time-sep}
        \limsup_{\epsilon \rightarrow 0} \tau(p',q) \geq 1
    \end{equation}
    follows at once (in fact equality holds). But since $p'
    \longrightarrow p$ for $\epsilon \longrightarrow 0$ we see that
    $\tau$ is \emph{not} upper semi continuous and hence not
    continuous. In fact, moving $q$ further upwards we can make the
    jump arbitrarily high.
\end{example}
The question is now whether we have spacetimes where $\tau$ is
continuous (and finite). Clearly, Minkowski spacetime is an example
where $\tau$ is continuous and finite. More generally, convex
spacetimes have this feature:
\begin{example}[Time separation for convex spacetimes]
    \label{example:convex-spacetimes-timesep-good}
    \index{Spacetime!convex}%
    Suppose that $M$ is geodesically convex, or $U \subseteq M$ is a
    geodesically convex neighborhood. Then the time separation
    $\tau_U$ on $U$ is finite and continuous. Indeed, this follows
    from Theorem~\ref{theorem:time-separation},
    \refitem{item:maximizing-geodesic-in-geod-convex-set} at once.
\end{example}
Slightly less obvious is the following situation of a globally
hyperbolic spacetime: In fact, this statement can be seen as an
additional motivation for the definition of globally hyperbolic
spacetimes as in
Definition~\ref{definition:globally-hyperbolic-spacetime}. However, it
was noted that
Definition~\ref{definition:globally-hyperbolic-spacetime} implies
strong causality as well. Using this observation, we can quote the
following result \cite[Prop.~21 in Chap.~14]{oneill:1983a}:
\begin{example}[Time separation for globally hyperbolic spacetimes]
    \label{example:timesep-globally-hyp}
    \index{Spacetime!globally hyperbolic}%
    Suppose that $(M,g)$ is globally hyperbolic. Then the time
    separation $\tau$ is finite and continuous, see also
    \cite[Thm.~3.83]{minguzzi.sanchez:2006a:pre}.
\end{example}
With these two fundamental examples in mind we conclude this short
subsection on time separation and refer to
\cite[Chap.~14]{oneill:1983a} for additional information.

%
%

\subsection{Uniqueness of Solutions to the Wave Equation}
\label{satz:uniqueness-solutions}

\begin{figure}
    \centering
    \input{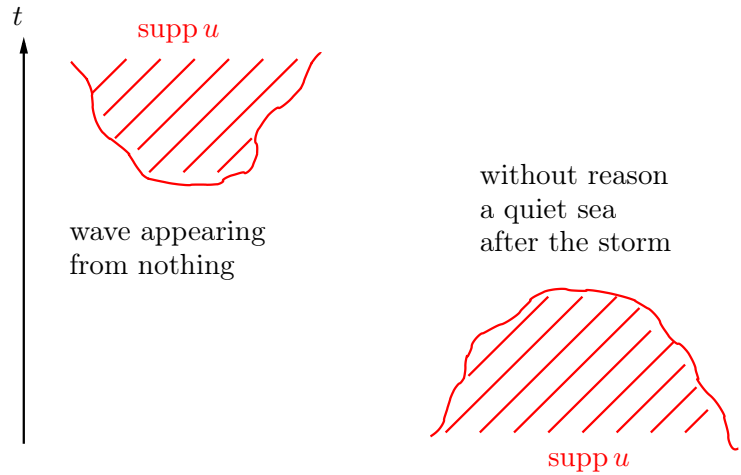}
    \caption{\label{fig:wave-with-future-past-compact-support}%
      Waves with either past or future compact support should not
      exist.
    }
    \index{Freak wave}%
    \index{Storm}%
\end{figure}

In general, the wave equation
\begin{equation}
    \label{eq:hom-wave-equation}
    Du = 0
\end{equation}
has many solutions $u \in \Sec[-\infty](E)$: physically such solutions
correspond to propagating waves without sources. However, also from
our physical intuition we expect that a propagating wave without any
possibility to interact with source terms has to ``travel
forever''. Thus a non-trivial solution of \eqref{eq:hom-wave-equation}
with either future or past compact support should not exist, see
Figure~\ref{fig:wave-with-future-past-compact-support}.
Assuming some (technical) conditions about the causality structure of
the spacetime this is indeed true.

To formulate these conditions first recall that the causal relation
$\leq$ is called \emph{closed} if for any sequence $p_n
\longrightarrow p$ and $q_n \longrightarrow q$ with $p_n \leq q_n$ we
have $p \leq q$ as well. Equivalently, this means that
\begin{equation}
    \label{eq:lying-in-the-future-relation}
    \index{Causal relation!closed}%
    J^+_M =
    \left\{
        (p,q) \in M \times M
        \; \big| \;
        p \leq q
    \right\}
    \subseteq M \times M
\end{equation}
is a closed subset of $M \times M$.

We consider now the following three properties which will turn out to
be sufficient to guarantee the uniqueness of the solutions to
\eqref{eq:hom-wave-equation} with future or past compact support.
\begin{compactenum}
\item \index{Spacetime!causal}%
    \label{item:causally-simple-1} $(M,g)$ is causal, i.e. there are
    no causal loops.
\item \label{item:causally-simple-2} $J^+_M$ is closed.
\item \index{Time separation!finite and continuous}%
    \label{item:timesep-con-and-finite} The time separation $\tau$
    is finite and continuous.
\end{compactenum}
Concerning the relation among these three properties some remarks are
in due:
\begin{remark}[Causally simple spacetimes]
    \label{remark:causally-simple}
    \index{Spacetime!causally simple}%
    \index{Lorentz manifold|see{Spacetime}}%
    A time-oriented Lorentz manifold $(M,g)$ which satisfies the
    causality condition \refitem{item:causally-simple-1} is called
    \emph{causally simple} if in addition $J^\pm_M(p)$ are closed for
    all $p \in M$, see
    e.g. \cite[Sect.~3.10]{minguzzi.sanchez:2006a:pre}. One can show
    that this is equivalent to being causal and $J^+_M$ being closed
    which is equivalent to being causal and $J^\pm_M(K)$ being closed
    for all compact subsets $K \subseteq M$. Thus
    \refitem{item:causally-simple-1} and
    \refitem{item:causally-simple-2} just say that $(M,g)$ is causally
    simple.
\end{remark}
\begin{remark}
    \label{remark:3}
    ~
    \begin{remarklist}
    \item \label{item:timesep-finite} The finiteness of $\tau$ clearly
        implies that there are no timelike loops.
    \item \label{item:timesep-cont-is-more} There are examples of
        causally simple spacetimes which do not satisfy
        \refitem{item:timesep-con-and-finite}. So this is indeed an
        additional requirement.
    \item \label{item:convex-does-it-all} Convex spacetimes satisfy
        all three requirements, see
        Example~\ref{example:convex-spacetimes-timesep-good}.
    \item \label{item:globally-hyp-does-it-all} Also globally
        hyperbolic spacetimes satisfy all three conditions, see
        e.g. the discussion in
        \cite[Thm.~3.83]{minguzzi.sanchez:2006a:pre}.
    \end{remarklist}
\end{remark}
With these conditions we can now prove the following theorem:
\begin{theorem}
    \label{theorem:future-past-comp--hom-solution-is-zero}
    \index{Wave equation!homogeneous}%
    Assume that a time-oriented Lorentz manifold $(M,g)$ satisfies the
    three conditions \refitem{item:causally-simple-1},
    \refitem{item:causally-simple-2},
    \refitem{item:timesep-con-and-finite}. Let $D \in \Diffop^2(E)$ be
    a normally hyperbolic differential operator on some vector bundle
    $E \longrightarrow M$ and let $u \in \Sec[-\infty](E)$ be a
    distributional section. If $u$ has either past or future compact
    support and satisfies the homogeneous wave equation
    \begin{equation}
        \label{eq:hom-wave}
        Du = 0,
    \end{equation}
    then $u = 0$.
\end{theorem}
\begin{proof}
    We follow \cite[Thm.~3.1.1]{baer.ginoux.pfaeffle:2007a} and
    consider the case of a future compact support $\supp u$. We have
    to show $\supp u = \emptyset$. We assume the converse and choose a
    point $\widetilde{q} \in \supp u$. The future compactness of
    $\supp u$ means that for all $p \in M$ the subset $\supp u \cap
    J^+_M(p) \subseteq M$ is compact. Choosing $p \in
    I^-_M(\widetilde{q})$ we obtain a non-empty intersection $\supp u
    \cap J^+_M(p)$, see Figure~\ref{fig:top-of-supp-u}.
    \begin{figure}
        \centering
        \input{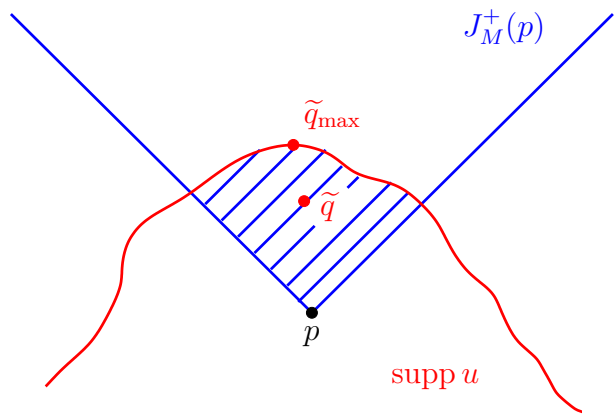}
        \caption{\label{fig:top-of-supp-u}%
          Finding the ``top'' of the support of $u$.
        }
    \end{figure}
    We now want to find the ``top'' of the intersection $\supp u \cap
    J^+_M(p)$: since the time separation $\tau$ is continuous the map
    $q \mapsto \tau(p,q)$ for $q \in \supp u \cap J^+_M(p)$ takes its
    maximal value $\tau(p,\widetilde{q}_{\textrm{max}})$ at some (not
    necessarily unique) $\widetilde{q}_{\textrm{max}} \in \supp u \cap
    J^+_M(p)$ by compactness. We consider now the intersection $\supp
    u \cap J^+_M(\widetilde{q}_{\textrm{max}})$ which is still compact
    and non-empty since $\widetilde{q}_{\textrm{max}} \in \supp u \cap
    J^+_M(\widetilde{q}_{\textrm{max}})$. Figure~\ref{fig:top-of-supp-u}
    suggests that this subset is actually rather small. In fact, for
    $q \in \supp u \cap J^+_M(\widetilde{q}_{\textrm{max}})$ we have
    on one hand $\tau(p,q) \geq \tau(p,\widetilde{q}_{\textrm{max}})$
    since $q \geq \widetilde{q}_{\textrm{max}}$ and $\tau(p,q) \leq
    \tau(p,\widetilde{q}_{\textrm{max}})$ by the maximality of
    $\widetilde{q}_{\textrm{max}}$. Thus
    \[
    \tau(p,q) = \tau(p,\widetilde{q}_{\textrm{max}})
    \]
    for all $q \in \supp u \cap
    J^+_M(\widetilde{q}_{\textrm{max}})$. Among all the $q \in \supp u
    \cap J^+_M(\widetilde{q}_{\textrm{max}})$ we want to find a
    particular $q_{\textrm{max}}$ such that the intersection $\supp u
    \cap J^+_M(q_{\textrm{max}})$ contains \emph{only}
    $q_{\textrm{max}}$ and no other points. In order to find such an
    optimal point we proceed as follows. The compact subset $\supp u
    \cap J^+_M(\widetilde{q}_{\textrm{max}})$ is \emph{partially
      ordered} via $\leq$. Indeed, $p \leq p$ as well as transitivity,
    $p \leq p'$ and $p' \leq p''$ implies $p \leq p''$, are always
    true. Since we do not have causal loops also $p \leq p'$ and $p'
    \leq p$ implies $p = p'$. Now assume that we have an increasing
    chain of elements $\{q_i\}_{i \in I}$, i.e. a subset of points of
    which any two are in relation ``$\leq$''. Our manifold being
    second countable we can find a countable dense subset $\{q_n\}_{n
      \in \mathbb{N}} \subset \{q_i \}_{i \in I}$ which is ordered
    again since it is the subset of an ordered set. We define $Q_n$ to
    be the maximum of $\{q_1, \ldots, q_n \}$ for all $n \in
    \mathbb{N}$. This gives a sequence $(Q_n)$ of elements in
    $\{q_n\}_{n \in \mathbb{N}}$ such that for every $q_k$ there is an
    $n_0$ with $q_k \leq Q_n$ for all $n \geq n_0$. Now the $Q_n$ have
    accumulation points in the compact subset $\supp u \cap
    J^+_M(\widetilde{q}_{\textrm{max}})$. Thus fixing a suitable
    subsequence $Q_{n_m}$ this converges to some $Q_\infty$ which is
    still in $\supp u \cap J^+_M(\widetilde{q}_{\textrm{max}})$. Since
    the relation $\leq$ is closed we see that $Q_\infty$ is an upper
    bound for all the $q_n$, i.e. we have $q_n \leq Q_\infty$ for all
    $n \in \mathbb{N}_0$. Since the $\{q_n \}_{n \in \mathbb{N}_0}
    \subseteq \{q_i \}_{i \in I}$ are dense and ``$\leq$'' is a closed
    relation, we also have
    \[
    q_i \leq Q_\infty
    \]
    \index{Zorn's Lemma}%
    for all indexes $i \in I$. This shows that inside $\supp u \cap
    J^+_M(\widetilde{q}_{\textrm{max}})$ every increasing chain has an
    upper bound. Thus we are in the position to use Zorn's Lemma and
    conclude that there are maximal elements for all of $\supp u \cap
    J^+_M(\widetilde{q}_{\textrm{max}})$. Thus we pick one of these
    not necessarily unique ones and obtain a $q_{\textrm{max}} \in
    \supp u \cap J^+_M(p)$ such that on one hand $q \mapsto \tau(p,q)$
    attains its maximum at $q_{\textrm{max}}$ \emph{and} we have $q
    \leq q_{\textrm{max}}$ for all $q \in \supp u \cap
    J^+_M(\widetilde{q}_{\textrm{max}})$. Thus it follows that
    \[
    \supp u \cap J^+_M(q_{\textrm{max}}) = \{ q_{\textrm{max}} \}
    \]
    by the maximality property with respect to ``$\leq$''. Thus we
    arrive at the following picture, see
    Figure~\ref{fig:point-on-top}, where $q_{\textrm{max}}$ is now on
    the top of $\supp u \cap J^+_M(p)$ and $J^+_M(q_{\textrm{max}})$
    does not intersect $\supp u \cap J^+_M(p)$ except in
    $q_{\textrm{max}}$.
    \begin{figure}
        \centering
        \input{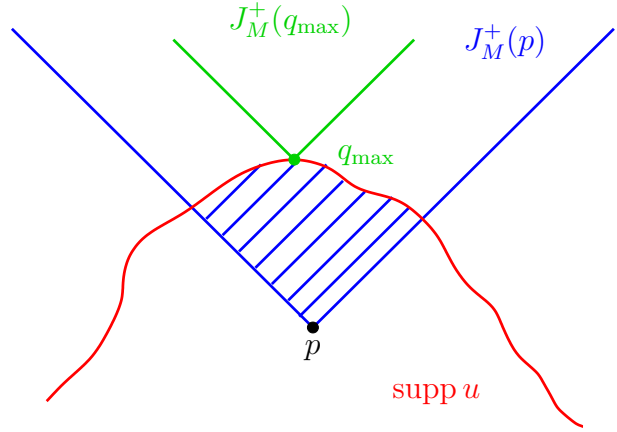}
        \caption{\label{fig:point-on-top}%
          The point $q_{\textrm{max}}$ on top of $\supp u \cap
          J^+_M(p)$.
        }
    \end{figure}
    Now we consider a causal neighborhood $U \subseteq U'$ of
    $q_{\textrm{max}}$ in some convex $U' \subseteq M$ with $U^\cl$
    compact in $U'$, such that the volume of $U^\cl$ is sufficiently
    small.
    \begin{figure}
        \centering
        \input{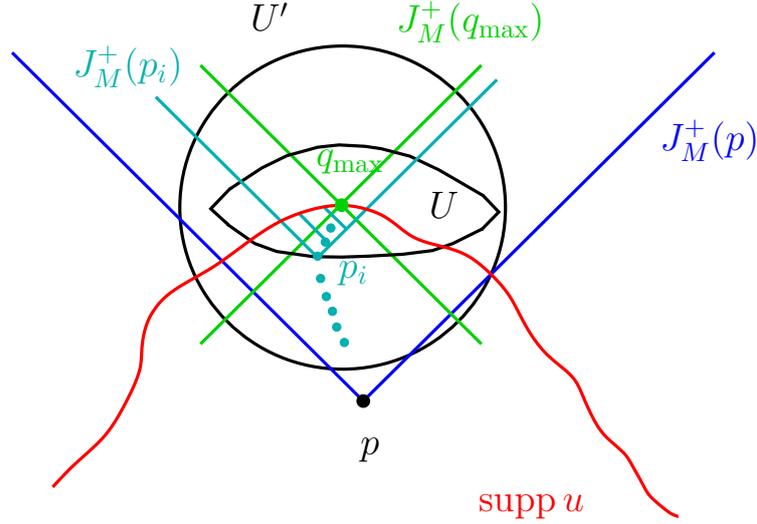}
        \caption{\label{fig:sequence-figure}%
          The sequence $p_i$ approaching $q_{\max}$.
        }
    \end{figure}
    Consider a sequence of points $p_i \in U$ which converge to
    $q_{\textrm{max}}$ and are contained in $I^-_M(q_{\textrm{max}})
    \cap I^+(p)$. Then for large enough $i$ the intersection
    $J^+_M(p_i) \cap \supp u$ is entirely contained in $U$. Indeed,
    assume this is not true. Then for each $i \in \mathbb{N}$ we can
    find a $q_i \in J^+_M(p_i) \cap \supp u$ which is \emph{not} in
    $U$. By the compactness of $J^+_M(p) \cap \supp u$ we can assume
    that $q_i \longrightarrow q$ converges inside $J^+_M(p) \cap \supp
    u$, probably we have to pass to a suitable subsequence. Since $q_i
    \in J^+_M(p_i)$ and $q_i \longrightarrow q$ as well as $p_i
    \longrightarrow q_{\textrm{max}}$ we conclude by the closedness of
    the relation ``$\leq$'' that $q \geq q_{\textrm{max}}$. Thus $q
    \in J^+_M(q_{\textrm{max}}) \cap \supp u = \{ q_{\textrm{max}} \}$
    and hence $q = q_{\textrm{max}}$. On the other hand, $q_i \notin
    U$ implies $q \notin U$ as $U$ is open which gives a contradiction
    to $q_{\textrm{max}} \in U$. Thus we arrive indeed at the
    situation as in Figure~\ref{fig:sequence-figure}. We choose such a
    point $p_i$ and consider the compact subset $K = J^+_M(p_i) \cap
    \supp u \subseteq U$.
    \begin{figure}
        \centering
        \input{sequence2.\pictype}
        \caption{\label{fig:Umgebungqmax}%
          The neighborhood of $q_{\textrm{max}}$.
        }
    \end{figure}
    The open subset $\widetilde{U} = I^+_M(p_i) \cap U$ contains
    $q_{\textrm{max}}$ and is therefor an open neighborhood of
    $q_{\textrm{max}}$, see Figure~\ref{fig:Umgebungqmax}. Now we want
    to show that $u(\varphi) = 0$ for all test sections $\varphi \in
    \Secinfty_0(E^*\at{\widetilde{U}})$. Since with $D \in
    \Diffop^2(E)$ also the transposed operator $D^\Trans \in
    \Diffop^2(E^*)$ is normally hyperbolic we can solve the
    inhomogeneous wave equation
    \[
    D^\Trans \psi = \varphi
    \]
    with some $\psi \in \Secinfty(E^*\at{U})$ by
    Theorem~\ref{theorem:dual-of-fund-solution-gives-inhom-solution}. In
    particular, we know that with $\varphi$ being smooth also $\psi$
    is smooth. Moreover, this theorem also provides us information on
    the support: we can take the advanced solution for which we have
    $\supp \psi \subseteq J^+_U(\supp \varphi) \subseteq J^+_M(p_i)
    \cap U$, see Figure~\ref{fig:SupportOfVarphi}.
    \begin{figure}
       \centering
        \input{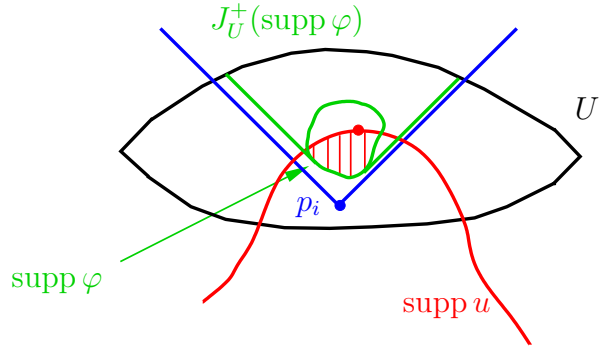}
        \caption{\label{fig:SupportOfVarphi}%
          The support of $\varphi$ and its future.
        }
    \end{figure}
    Thus we get
    \begin{align*}
        \supp  u \cap \supp \psi
        \subseteq \supp u \cap J^+_M(p_i) \cap U
        \subseteq \supp u \cap J^*_M(p_i)
        = K.
    \end{align*}
    This is now the compactness criterion we need for applying $u$ to
    the section $\psi$ according to
    Proposition~\ref{proposition:supp-overlap-is-compact}. Note that
    both have non-compact support in general. But then we have
    \[
    u(\varphi)
    = u( D^\Trans \psi)
    = Du (\psi)
    = 0
    \]
    by $Du = 0$. This shows that $u$ vanishes on all test sections
    $\varphi \in \Secinfty_0(E^*\at{\widetilde{U}})$. Thus the support
    of $u$ is disjoint from $\widetilde{U}$. Now we arrived at the
    desired contradiction as $q_{\textrm{max}} \in \supp u$ but
    $\widetilde{U}$ is an open neighborhood of
    $q_{\textrm{max}}$. Hence $\supp u = \emptyset$ follows and thus
    $u = 0$. The case of past compact support is analogous.
\end{proof}

From this theorem we immediately obtain several statements about the
solutions of the wave equations. Under the same assumptions on the
global structure of $M$, i.e. we require a causally simple spacetime
with finite and continuous time separation, one obtains the following
statement:
\begin{corollary}
    \label{corollary:uniqueness-fundamental-solution-with-pf-compact-supp}
    \index{Fundamental solution!uniqueness}%
    Let $(M,g)$ be a causally simple Lorentz manifold with finite and
    continuous time separation. Then for every normally hyperbolic
    differential operator $D \in \Diffop^2(E)$ there exists at most
    one fundamental solution at $p \in M$ with past compact support
    and at most one with future compact support.
\end{corollary}
\begin{proof}
    Indeed if $DF = \delta_p = D \widetilde{F}$ then $F-\widetilde{F}$
    solves the homogeneous wave equation \emph{and} has still past (or
    future) compact support. Thus $F - \widetilde{F} = 0$ by the
    preceding theorem.
\end{proof}

Now we pass to a globally hyperbolic spacetime $(M,g)$. On one hand we
know from Remark~\ref{remark:3} that $(M,g)$ satisfies the hypothesis
of Theorem~\ref{theorem:future-past-comp--hom-solution-is-zero}. On
the other hand on a globally hyperbolic spacetime the subset
$J^\pm_M(p)$ are always past/future compact: indeed, by the very
definition of global hyperbolicity, $J^+_M(p) \cap J^-_M(q) =
J_M(p,q)$ is a \emph{compact} diamond for all $p,q \in M$. This is
just the statement that $J^+_M(p)$ is past compact and $J^-_M(q)$ is
future compact. This gives immediately the following result:
\begin{corollary}
    \label{corollary:uniqueness-glob-hyp}
    \index{Spacetime!globally hyperbolic}%
    \index{Green function!uniqueness}%
    Let $(M,g)$ be a globally hyperbolic Lorentz manifold. Then for
    every normally hyperbolic differential operator $D \in
    \Diffop^2(E)$ there exists at most one advanced and at most one
    retarded Green function at $p \in M$.
\end{corollary}

\begin{example}[Uniqueness of Green functions]
    \label{example:unique-green-functions}
    Let $(\mathbb{R}^n, \eta)$ be the flat Minkowski spacetime as
    before. Since this is a globally hyperbolic spacetime we have the
    following global and unique Green functions:
    \begin{examplelist}
    \item \index{Riesz distribution}%
        The Riesz distributions $R^\pm(2)$ are the unique advanced and
        retarded Green functions for $\dAlembert$ at $0$. Their
        translates to arbitrary $p \in \mathbb{R}^n$ are the unique
        advanced and retarded Green functions for $\dAlembert$ at $p$.
    \item \index{Klein-Gordon equation}%
        The distributions $\widetilde{\mathcal{R}}^\pm(p) =
        \sum_{k=0}^\infty (-m^2)^k R^\pm(2+2k,p)$ are the unique
        advanced and retarded Green functions at $p \in \mathbb{R}^m$
        of the Klein-Gordon operator $\dAlembert + m^2$ on Minkowski
        spacetime.
    \end{examplelist}
\end{example}

\begin{figure}
    \centering
    \input{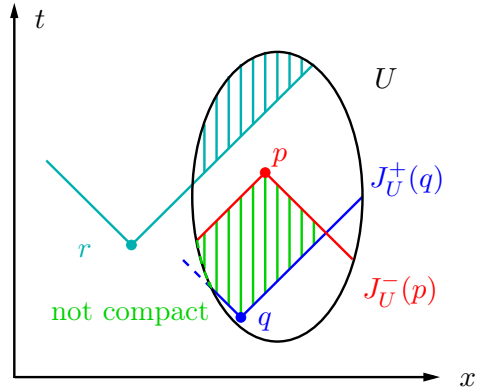}
    \caption{\label{fig:JU-not-compact-on-convex}%
      Convex domain in Minkowski spacetime with non-unique Green
      functions.
    }
\end{figure}

Finally, we mention that on convex domains we can not conclude the
uniqueness of advanced and retarded Green functions without further
assumptions. Even though geodesically convex domains satisfy the
hypothesis of
Theorem~\ref{theorem:future-past-comp--hom-solution-is-zero} it may
\emph{not} be true that $J^+_U(p)$ is past or future compact,
respectively. This is clear from the example in
Figure~\ref{fig:JU-not-compact-on-convex}.
Indeed, if in this situation we take the Green function $R^\pm(2)(p)$
of $\dAlembert$ on $(\mathbb{R}^n,\eta)$ and restrict them to $U$ we
obtain advanced and retarded Green functions $\{ R^\pm(2)(p) \at{U}
\}_{p \in U}$ for all points $p \in U$. Taking now a point $r \in
\mathbb{R}^n$ as in Figure~\ref{fig:JU-not-compact-on-convex} and
adding $R^+(2)(r)\at{U}$ to $R^+(2)(q)\at{U}$ we still have an
advanced Green function since $\dAlembert R^+(2)(r)=0$ on
$U$. However, as $\singsupp R^+(2)(r) = C^+(r)$ by
Proposition~\ref{proposition:support-of-riesz-distribution} for $n$
even, we see that this new advanced Green function differs from
$R^+(2)(q)\at{U}$ on the intersection $C^+(r) \cap U$, even in an
essential way. Thus we cannot hope for uniqueness of advanced and
retarded Green functions in general.


%% file: cauchyprob.tex
%
%

In order to pose the Cauchy problem we have to assume that we have a
Cauchy hypersurface on which we can specify the initial values. Thus
in this section we assume that $(M,g)$ is a globally hyperbolic
spacetime and $\iota: \Sigma \hookrightarrow M$ is a smooth spacelike
Cauchy hypersurface in $M$ whose existence is guaranteed by
Theorem~\ref{theorem:globalhyp-and-cauchy-surface}. Furthermore, the
future directed timelike normal vector field of $\Sigma$ will be
denoted by $\mathfrak{n} \in \Secinfty(TM \at{\Sigma})$ as in
Section~\ref{sec:CauchyProblemGreenFunctions}.
\begin{remark}
    \label{remark:initial-data-distributions}
    \index{Initial conditions!regularity}%
    \index{Initial conditions!distributional}%
    When solving the wave equation $Du = v$ in a distributional sense
    for $u,v \in \Sec[-\infty](E)$ one might be tempted to ask for the
    \emph{initial conditions} of $u$ on $\Sigma$. However, since
    $\iota: \Sigma \hookrightarrow M$ is far from being a submersion
    the restriction $\iota^*u$ is \emph{not at all} well-defined. To
    see the problem one should try to define $\iota^* \delta$ for the
    $\delta$ distribution on $\mathbb{R}$ and $\iota: \{ 0 \}
    \hookrightarrow \mathbb{R}$. Thus for the Cauchy problem to make
    sense we either have to specify conditions on $u$ and $v$ which
    ultimately allow to define $\iota^*u$ etc., or we restrict
    ourselves directly to regular initial conditions and solutions of
    some $\Fun$-regularity. As usual, the most convenient situation
    will be the $\Cinfty$-case.
\end{remark}
In view of the above remark we will therefore focus on regular and
smooth solutions and initial conditions. Thus the Cauchy problem
consists in the following task: Given an inhomogeneity $v \in
\Secinfty(E)$ we want to find a solution $u \in \Secinfty(E)$ of
\begin{equation}
    \label{eq:wave-equation-smooth-solution}
    Du = v
\end{equation}
for given initial conditions $u_0, \dot{u}_0 \in
\Secinfty_0(\iota^\#E)$, i.e.
\begin{align}
    \label{eq:initial-conditions}
    \iota^\#u &= u_0, \\
    \iota^\# \nabla_{\mathfrak{n}}^E u &= \dot{u}_0.
\end{align}
Here $\nabla^E$ will always be the covariant derivative on $E$
determined by $D$ as usual. Note that the left hand side of
\eqref{eq:initial-conditions} is indeed well-defined as for $p \in
\Sigma$ the value $\nabla^E_{\mathfrak{n}(p)}u \in E_p$ is defined as
$\nabla^E$ is function linear in the tangent vector field
argument. Thus we can interpret $p \mapsto
\nabla_{\mathfrak{n}(p)}^Eu$ indeed as a section of $\iota^\#E$.

%
%

\subsection{Uniqueness of the Solution to the Cauchy Problem}
\label{satz:uniq-solut-cauchy}

As for the solutions of the homogeneous wave equation also for the
Cauchy problem the uniqueness will be easier to show than the
existence. We start with some preparatory material on the adjoint
$D^\Trans$ of $D$. Recall from
Theorem~\ref{theorem:agjoint-for-fixed-density} that $D^\Trans \in
\Diffop^2(E^*)$ is determined by
\begin{equation}
    \label{eq:adjoint-of-D}
    \int_M \varphi (Du) \: \mu_g
    = \int_M (D^\Trans \varphi) u \: \mu_g
\end{equation}
for $\varphi \in \Secinfty(E^*)$ and $u \in \Secinfty(E)$ with at
least one of them having compact support. We want to compute now
$D^\Trans$ explicitly.
\begin{lemma}
    \label{lemma:transpose-of-D}
    Let $D \in \Diffop^2(E)$ be a normally hyperbolic differential
    operator written as $D = \dAlembert^\nabla + B$ with $B \in
    \Secinfty(\End(E))$ and the connection d'Alembertian
    $\dAlembert^\nabla$ build out of the connection $\nabla^E$ defined
    by $D$.
    \begin{lemmalist}
    \item \label{item:transpose-of-D} The transposed operator
        $D^\Trans \in \Diffop^2(E^*)$ is given by
        \begin{equation}
            \label{eq:transpose-of-D}
            D^\Trans = \dAlembert^\nabla + B^\Trans
        \end{equation}
        where $\dAlembert^\nabla$ is the connection d'Alembertian with
        respect to the induced connection $\nabla^{E^*}$ for $E^*$
        coming from $\nabla^E$.
    \item \label{item:pairing-rule-dAlembertian} For $s \in
        \Secinfty(E)$ and $\psi \in \Secinfty(E^*)$ we have
        \begin{equation}
            \label{eq:pairing-rule-dAlembertian}
            \dAlembert( \psi(s) )
            = (\dAlembert^\nabla \psi)(s)
            + \psi( \dAlembert s)
            + \SP{g^{-1}, (\SymD^{E^*} \psi) \vee (\SymD^Es)}.
        \end{equation}
    \item \label{item:difference-D-and-Dtrans} For $s \in
        \Secinfty(E)$ and $\psi \in \Secinfty(E^*)$ we have
        \begin{equation}
            \label{eq:difference-D-and-Dtrans}
            (D^\Trans \psi)(s) - \psi (Ds)
            = \divergenz \left(
                \left(
                    (\SymD^{E*} \psi)(s) - \psi (\SymD^E s)
                \right)^\#
            \right).
        \end{equation}
    \end{lemmalist}
\end{lemma}
\begin{proof}
    For the first part we use Theorem~\ref{theorem:Neumaier-Theorem}
    as well as the result from
    Example~\ref{example:connection-dAlembert}. In this example we
    found that $\dAlembert^\nabla = ( \frac{\I}{\hbar} )^2 \stdrep( 2
    g^{-1} \tensor \id_E)$. Since the remaining part $B$ is
    $\Cinfty(M)$-linear it is clear that $B = \stdrep(B)$ in the sense
    that the tensor field $B$ acts pointwise as endomorphism on
    sections of $E$. By Theorem~\ref{theorem:Neumaier-Theorem} we have
    $\stdrep(B)^\Trans = \stdrep(B^\Trans)$ as there are no degrees to
    be lowered by the divergence operator
    $\divergenz_{\mu_g}^{\End(E)}$. In fact, we have $\varphi (Bs) =
    (B^\Trans \varphi)(s)$ by definition of the pointwise
    transposition from which $\stdrep(B)^\Trans = \stdrep(B^\Trans)$
    is immediate. The transpose of $\dAlembert^\nabla$ is more
    involved: here we need to compute the divergence of $2 g^{-1}
    \tensor \id_E$. First we note that the one-form $\alpha$ measuring
    the non-parallelness of the integration density $\mu_g$ is
    vanishing thanks to
    Proposition~\ref{proposition:levi-civita-connection-etc},
    \refitem{item:Riem-densitiy-is-cov-const}. Thus
    $\divergenz_\mu^{\End(E)}$ coincides with the connection
    divergence $\divergenz_\nabla^{\End(E)}$ where we have to use the
    induced connection on $\End(E)$ coming from $\nabla^E$. Thus we
    have to compute
    \begin{align*}
        \divergenz_\nabla^{\End(E)} (g^{-1} \tensor \id_E)
        &= \inss(\D\!x^i)
        \nabla_{\frac{\partial}{\partial x^i}}^{\End(E)}
        (g^{-1} \tensor \id_E) \\
        &= \inss(\D\!x^i)
        \left(
            \nabla_{\frac{\partial}{\partial x^i}} g^{-1} \tensor
            \id_E
            + g^{-1} \tensor \nabla_{\frac{\partial}{\partial
                x^i}}^{\End(E)} \id_E
        \right) \\
        &= 0 + 0,
    \end{align*}
    since on one hand $g^{-1}$ is parallel for the Levi-Civita
    connection and on the other hand $\id_E$ is a parallel section
    with respect to $\nabla^{\End(E)}$. In fact, the latter result is
    just the definition of $\nabla^{\End(E)}$: for $A \in
    \Secinfty(\End(E))$ and $s \in \Secinfty(E)$ the induced
    connection $\nabla^{\End(E)}$ is determined by
    \[
    ( \nabla^{\End(E)} A)(s) = \nabla^E(As) - A( \nabla^E s).
    \]
    Thus $\id_E$ is covariantly constant since the right hand side
    will be zero for $A = \id_E$. We conclude that
    \begin{align*}
        D^\Trans
        &=
        \left(\frac{\I}{\hbar}\right)^2
        \stdrep(2 g^{-1} \tensor \id_E)^\Trans
        +
        \stdrep(B)^\Trans \\
        &=
        \left(\frac{\I}{\hbar}\right)^2
        \stdrep(2 g^{-1} \tensor \id_{E^*})
        + B^\Trans \\
        &= \dAlembert^\nabla + B^\Trans,
    \end{align*}
    where now $\dAlembert^\nabla$ is the connection d'Alembertian on
    $E^*$ with respect to the induced connection $\nabla^{E^*}$. For
    the second part we first show the following Leibniz rule of
    $\dAlembert$ with respect to natural pairings, see also
    Lemma~\ref{lemma:connection-dAlembert}. We compute
    \begin{align*}
        \dAlembert (\psi(s))
        &=
        \frac{1}{2} \SP{g^{-1}, \SymD^2( \psi(s) )} \\
        &=
        \frac{1}{2}
        \SP{g^{-1}, \SymD((\SymD^{E^*} \psi)(s) + \psi(\SymD^E s))} \\
        &=
        \frac{1}{2}
        \SP{g^{-1},
          ((\SymD^{E^*})^2 \psi)(s)
          + 2 (\SymD^{E^*} \psi) \vee (\SymD^E s)
          + \psi ((\SymD^E)^2s)} \\
        &=
        (\dAlembert^\nabla \psi)(s)
        + \SP{g^{-1}, (\SymD^{E^*}\psi) \vee (\SymD^E s)}
        + \psi( \dAlembert^\nabla s),
    \end{align*}
    where we have used the compatibility of the symmetrized covariant
    derivative operators $\SymD, \SymD^E$ and $\SymD^{E^*}$ with
    natural pairings. This compatibility is immediate from the
    definition of these operators, see
    Proposition~\ref{proposition:Symmetric-Covariant-Derivative},
    \refitem{item:SymDLocalForm}. This shows the second part. For the
    last part we know from Theorem~\ref{theorem:Neumaier-Theorem} that
    $(D^\Trans \psi)(s) - \psi(Ds)$ vanishes after integrating over
    $M$ with respect to $\mu_g$. Thus it has to be a divergence of
    some vector field with respect to $\mu_g$. However, this vector
    field is only unique up to a divergence free vector field. Thus
    \eqref{eq:pairing-rule-dAlembertian} gives an explicit
    representative. First we notice that the contribution of $B$
    cancels as $(B^\Trans \psi)(s) - \psi (Bs) = 0$ holds
    pointwise. Thus we only have to consider $(\dAlembert^\nabla
    \psi)(s) - \psi( \dAlembert^\nabla s)$. We compute using the
    compatibility with natural pairing again
    \begin{align*}
        &\SymD \left(\left(\SymD^{E^*} \psi\right)(s)\right)
        - \SymD \left(\psi\left(\SymD^E s\right)\right) \\
        &= \left(\left(\SymD^{E^*}\right)^2 \psi\right)(s)
        + \left(\SymD^{E^*} \psi\right)\left(\SymD^E s\right)
        - \left(\SymD^{E^*} \psi\right)\left(\SymD^E s\right)
        - \psi\left(\left(\SymD^E\right)^2 s\right) \\
        &= \left(\left(\SymD^{E^*}\right)^2 \psi\right)(s)
        - \psi\left(\left(\SymD^E\right)^2 s\right).
    \end{align*}
    Hence we obtain for the left hand side of
    \eqref{eq:pairing-rule-dAlembertian}
    \begin{align*}
        (D^\Trans \psi)(s) - \psi(Ds)
        &=
        \frac{1}{2} \SP{g^{-1}, \left(\SymD^{E^*}\right)^2 \psi}(s)
        - \psi\left(
            \frac{1}{2} \SP{g^{-1}, \left(\SymD^E\right)^2s}
        \right)\\
        &= \frac{1}{2}
        \SP{g^{-1},
          \SymD \left(
              \left(\SymD^{E^*} \psi\right)(s)
              - \psi\left(\SymD^E s\right)
          \right)
        },
    \end{align*}
    since natural pairings commute. Now the one-form in this pairing
    is determined by
    \[
    \left(\left(\SymD^{E^*} \psi\right)(s)
        - \psi\left(\SymD^E s\right)
    \right)
    (\chi)
    = \left(\nabla_\chi^{E^*} \psi\right)(s)
    - \psi\left(\nabla_\chi^E s\right)
    \]
    for $\chi \in \Secinfty(TM)$. Since $g^{-1}$ is covariantly
    constant for the Levi-Civita connection, we have in general
    \[
    \frac{1}{2} \SP{g^{-1}, \SymD \alpha}
    = \divergenz(\alpha^\#)
    \]
    for arbitrary one-forms $\alpha \in \Secinfty(T^*M)$. This
    completes the proof.
\end{proof}

Now we consider again a small convex open subset $U' \subseteq M$ and
a causal open subset $U \subseteq U^\cl \subseteq U'$ of sufficiently
small volume so that we can use our local fundamental solutions from
Chapter~\ref{cha:LocalTheory}. The subset $U$ being causal includes
the diamonds $J_U(p,q)$ being compact and since it is inside the
convex $U'$ there are no causal loops in $U$. Thus $U$ is globally
hyperbolic and by Theorem~\ref{theorem:globalhyp-and-cauchy-surface}
we have a smooth spacelike Cauchy hypersurface $\iota: \Sigma
\hookrightarrow U$ in $U$. In fact, we recall from
\cite[Thm.~2.14]{minguzzi.sanchez:2006a:pre} that every point in $M$
has a neighborhood basis of globally hyperbolic open subsets. Thus we
can safely assume the existence of a smooth Cauchy hypersurface in
$U$. Since $\Sigma$ is spacelike the pull-back of $g$ to $\Sigma$ gives
a \emph{negative} definite metric (beware of our signature convention)
which includes a corresponding volume density. We denote this by
$\mu_\Sigma \in \Secinfty(\Dichten T^* \Sigma)$ and use it for
integration on $\Sigma$. Denote the fundamental solutions of $D^\Trans
\in \Diffop^2(E^*)$ on $U$ as constructed analogously to the ones of
$D$ by $G^\pm_U(p) \in \Sec[-\infty](E^*\at{U}) \tensor E_p$ where $p
\in U$. Then we have operators
\begin{equation}
    \label{eq:fundi-operators}
    G^\pm_U: \Secinfty_0(E\at{U}) \longrightarrow \Secinfty(E\at{U}),
\end{equation}
enjoying properties analogously to the $F^\pm_U$. In particular, we
have a dual map
\begin{equation}
    \label{eq:dual-of-fundi}
    (G^\pm_U)': \Sec[-\infty]_0(E^*\at{U}) \longrightarrow
    \Sec[-\infty](E^*\at{U}),
\end{equation}
which restricts to a map
\begin{equation}
    \label{eq:restriction-of-dual-of-fundi-to-smooth}
    (G^\pm_U)': \Secinfty_0(E^*\at{U}) \longrightarrow
    \Secinfty(E^*\at{U})
\end{equation}
by
Theorem~\ref{theorem:dual-of-fund-solution-gives-inhom-solution}. We
will need the difference between the advanced and retarded fundamental
solutions. We define the map
\begin{equation}
    \label{eq:LocalPropagator}
    \index{Propagator!local}%
    G_U = G^+_U - G^-_U: \Secinfty_0(E\at{U}) \longrightarrow
    \Secinfty(E\at{U}),
\end{equation}
which gives a dual map
\begin{equation}
    \label{eq:dual-of-propagator}
    G_U' = (G^+_U)' - (G^-_U)': \Sec[-\infty](E^*\at{U})
    \longrightarrow \Sec[-\infty](E^*\at{U}).
\end{equation}
On smooth sections $\varphi \in \Secinfty_0(E^*\at{U})$, viewed as
distributional sections, the map $G'_U$ is determined by
\begin{equation}
    \label{eq:dual-of-propagator-on-smooth}
    (G'_U \varphi)(u)
    = \varphi(G_U(u))
    =
    \int_U
    \varphi(p) \cdot
    \left(G^+_U(p)u - G^-_U(p)u\right)
    \: \mu_g(p),
\end{equation}
where $u \in \Secinfty_0(E\at{U})$ is a test section of
$E\at{U}$. Since we know by
Theorem~\ref{theorem:dual-of-fund-solution-gives-inhom-solution} that
$G'_U(\varphi)$ is actually a smooth section of $E^*\at{U}$, it makes
sense to restrict this section to $\Sigma$. Then we obtain the
following lemma:
\begin{lemma}
    \label{lemma:solution-integral-formula}
    Assume $u \in \Secinfty(E\at{U})$ is a solution to the homogeneous
    wave equation $Du=0$ and let $\varphi \in
    \Secinfty_0(E^*\at{U})$. Then we have
    \begin{equation}
        \label{eq:solution-integral-formula}
        \int_U \varphi(p) \cdot u(p) \: \mu_g(p)
        = \int_\Sigma
        \left(
            (\nabla_{\mathfrak{n}}^{E^*} G_U'(\varphi)) \cdot
            u_0(\sigma)
            -
            G_U'(\varphi)(\sigma) \cdot \dot{u}_0(\sigma)
        \right)
        \mu_\Sigma(\sigma),
    \end{equation}
    where $u_0 = \iota^\#u, \dot{u}_0 = \iota^\#
    \nabla_{\mathfrak{n}}^E u \in \Secinfty(i^\#E)$ are the initial
    values of $u$ on $\Sigma$.
\end{lemma}
\begin{proof}
    Let $\varphi \in \Secinfty_0(E^*\at{U})$ be a test section and let
    $\psi^\pm = (G^\pm_U)(\varphi) \in \Sec[-\infty](E^*\at{u})$ which
    is in $\Secinfty(E^*\at{U})$ by
    Theorem~\ref{theorem:dual-of-fund-solution-gives-inhom-solution}. We
    know from this theorem that $D^\Trans \psi^\pm = \varphi$ and
    $\supp \psi^\pm \subseteq J^\pm_U( \supp \varphi)$.
    \begin{figure}
        \centering
        \input{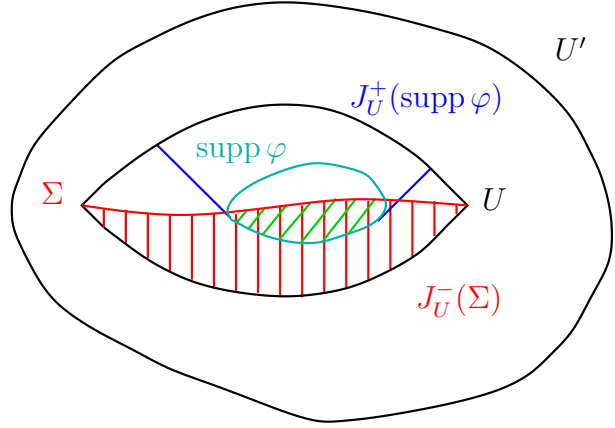}
        \caption{\label{fig:proof-situation} Sketch of the situation
          of the proof for
          Lemma~\ref{lemma:solution-integral-formula}.}
    \end{figure}
    For a Cauchy surface $\Sigma$ and an arbitrary compact subset $K
    \subseteq U$ one knows that $J^\pm_U(K) \cap J^\mp_U(\Sigma)$ is
    again compact, see Figure~\ref{fig:proof-situation}. For a proof
    of this fact we refer to
    \cite[Cor.~A.5.4]{baer.ginoux.pfaeffle:2007a} or
    \cite[p.~44]{minguzzi.sanchez:2006a:pre}. We know that the
    (globally hyperbolic) spacetime $U$ decomposes into the disjoint
    unions
    \[
    U = I^-_U(\Sigma) \dot{\cup} \Sigma \dot{\cup} I^+_U(\Sigma),
    \]
    where $I^\pm_U(\Sigma)$ are open and $\Sigma$ is the common
    boundary of these open subsets, see
    Remark~\ref{remark:cauchy-hypersurface}. Since we have chosen even
    a smooth Cauchy hypersurface, we can apply Gauss' Theorem in the
    form of Theorem~\ref{theorem:gauss-theorem} to the vector field
    \[
    \index{Gauss' Theorem}%
    X^\pm =
    \left(
        \left(\SymD^{E^*} \psi^\pm\right)(u)
        - \psi^\pm\left(\SymD^E u\right)
    \right)^\#
    \in \Secinfty(TU).
    \]
    Indeed, this vector field has support in $J^\pm_U(\supp
    \varphi)$. Thus the integrations over $I^\mp_U(\Sigma)$ and
    $J^\mp_U(\Sigma)$ as well as over $\Sigma$ itself are well defined
    because the integrands all have compact support. We consider first
    the case of $I^-_U(\Sigma)$. Then the future directed normal
    vector $\mathfrak{n}$ on $\Sigma$ points \emph{outwards} whence
    \[
    \int_{I^-_U(\Sigma)} \divergenz(X^+) \: \mu_g
    =
    \int_{\partial I^-_U(\Sigma) = \Sigma}
    g(X^+, \mathfrak{n}) \: \mu_\Sigma
    \tag{$*$}
    \]
    by Theorem~\ref{theorem:gauss-theorem}. We evaluate both sides
    explicitly. First we have
    \[
    \int_{I^-_U(\Sigma)} \divergenz(X^+) \: \mu_g
    = \int_{I^-_U(\Sigma)}
    \left(
        (D^\Trans \psi^+)(u) - \psi^+(Du)
    \right) \mu_g
    = \int_{I^-_U(\Sigma)} \varphi(u) \: \mu_g,
    \]
    by Lemma~\ref{lemma:transpose-of-D} and $Du = 0$ as well as
    $D^\Trans \psi^+ = \varphi$. For the right hand side of ($*$) we
    get
    \begin{align*}
        \int_\Sigma g(X^+, \mathfrak{n}) \: \mu_\Sigma
        &= \int_\Sigma
        \left(
            g\left(
                \left(\SymD^{E*} \psi^+ (u)\right)^\#,
                \mathfrak{n}
            \right)
            -
            g\left(
                \psi^+ \left(\SymD^E u\right)^\#,
                \mathfrak{n}
            \right)
        \right) \mu_\Sigma \\
        &= \int_\Sigma
        \left(
            \left(\SymD^{E^*} \psi^+(u)\right)(\mathfrak{n})
            -
            \left(\psi^+(\SymD u)\right)(\mathfrak{n})
        \right) \mu_\Sigma \\
        &= \int_\Sigma
        \left(
            \left(\nabla^{E^*}_\mathfrak{n} \psi^+\right)(u)
            -
            \psi^+\left(\nabla^E_\mathfrak{n} u\right)
        \right) \mu_\Sigma \\
        &= \int_\Sigma
        \left(
            \left(\nabla^{E^*}_\mathfrak{n} \psi^+\right)(u_0)
             - \psi^+(\dot{u}_0)
         \right) \mu_\Sigma,
    \end{align*}
    where we have omitted the restriction $\iota^\#$ in our notation
    for the sake of simplicity. Analogously, we obtain for
    $I^+_U(\Sigma)$ the result
    \[
    \int_{I^+_U(\Sigma)} \divergenz(X^-) \: \mu_g
    = -
    \int_\Sigma g(X^-, \mathfrak{n}) \: \mu_\Sigma,
    \tag{$**$}
    \]
    since now $\mathfrak{n}$ is pointing \emph{inwards}. Evaluating
    both sides gives
    \[
    \int_{I^+_U(\Sigma)} \divergenz(X^-) \: \mu_g
    = \int_{I^+_U(\Sigma)} \varphi(u) \: \mu_g
    \]
    and
    \[
    - \int_{\Sigma} g(X^-,\mathfrak{n}) \mu_\Sigma
    = - \int_\Sigma
    \left(
        \left(\nabla^{E^*}_\mathfrak{n} \psi^-\right)(u_0)
        - \psi^- (\dot{u}_0)
    \right) \mu_\Sigma.
    \]
    Thus taking the sum of ($*$) and ($**$) gives the equality
    \[
    \int_U \varphi(u) \: \mu_g
    = \int_\Sigma
    \left(
        \nabla^{E^*}_\mathfrak{n}(\psi^+ - \psi^-)(u_0)
        - (\psi^+ - \psi^-)(\dot{u}_0)
    \right) \mu_\Sigma,
    \]
    which is \eqref{eq:solution-integral-formula} by the definition of
    $\psi^+$ and $\psi^-$.
\end{proof}

\begin{lemma}
    \label{lemma:supp-of-solution-and-supp-of-initial-values}
    Assume $u \in \Secinfty(E\at{U})$ is a solution to the homogeneous
    wave equation $Du = 0$ and let $u_0, \dot{u}_0 \in
    \Secinfty(\iota^\# E)$ denote the initial values of $u$ on
    $\Sigma$. Then
    \begin{equation}
        \label{eq:supp-of-solution-and-supp-of-initial-values}
        \supp u \subseteq J_U( \supp u_0 \cup \supp \dot{u}_0).
    \end{equation}
\end{lemma}
\begin{proof}
    We determine the support of $u$ viewed as distributional
    section. This will coincide with the true support thanks to
    Remark~\ref{remark:restriction-and-support-of-gensec},
    \refitem{item:regular-gensec}. Thus let $\varphi \in
    \Secinfty_0(E^*\at{U})$ be a test section. Then we know that $
    \supp (G^\pm_U)'(\varphi) \subseteq J^\pm_U(\supp \varphi)$ by
    Lemma~\ref{lemma:supp-of-solution}. It follows that
    $G'_U(\varphi)$ has its support in $J_U(\supp \varphi)$. Suppose
    that $\supp u_0 \cup \supp \dot{u}_0$ will not intersect
    $J_U(\supp \varphi)$, see
    Figure~\ref{fig:support-of-initial-data}.
    \begin{figure}
        \centering
        \input{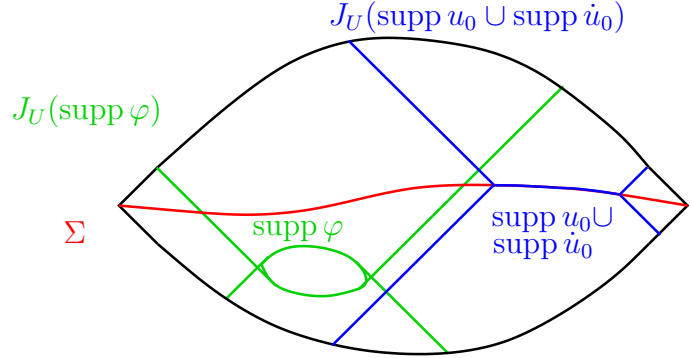}
        \caption{\label{fig:support-of-initial-data}%
          The support of the initial data.
        }
    \end{figure}
    Then this is equivalent to say that $\supp \varphi$ does not
    intersect $J_U(\supp u_0 \cup \supp \dot{u}_0)$. But by
    \eqref{eq:solution-integral-formula} the integral over $\Sigma$ is
    clearly $0$ whence $\int_U \varphi(u) \mu_g = 0$ follows. Thus
    $u$, viewed as distribution, vanishes on all these $\varphi$ where
    $\supp \varphi \cap J_U( \supp u_0 \cup \supp \dot{u}_0) =
    \emptyset$. But this means $\supp u \subseteq J_U( \supp u_0 \cup
    \supp \dot{u}_0)^\cl$. It remains to show that $J_U(\supp u_0 \cup
    \supp \dot{u}_0)$ is closed. In fact, this is true in general as
    we shall sketch now: Let $A \subseteq \Sigma$ be closed and
    consider $J^+_M(A)$ for simplicity. Let $p_n \in J^+_M(A)$ be a
    sequence of points with $p_n \longrightarrow p \in M$.
    \begin{figure}
        \centering
        \input{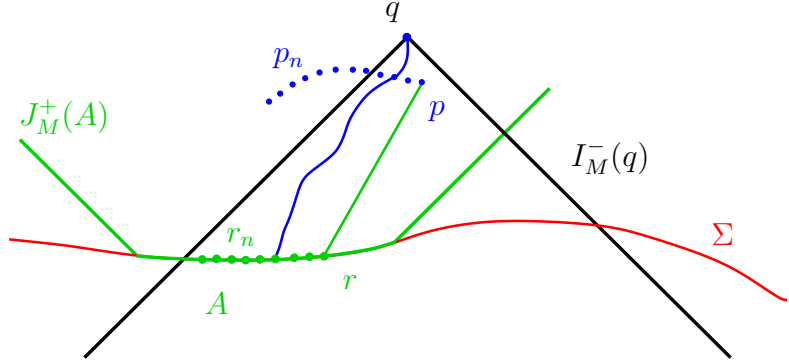}
        \caption{\label{fig:causal-influence-of-closed-set-is-closed}%
          The causal influence of a closed set $A \subseteq \Sigma$ in
          a Cauchy hypersurface is closed again.
        }
    \end{figure}
    Choose a point $q$ in the chronological future of $p$, i.e. we
    have $p \in I^-_M(q)$. Since $I^-_M(q)$ is open, all but finitely
    many $p_n$ are in $I^-_M(q)$ whence $q$ is in the chronological
    future of these $p_n$. Thus in particular $q \in J^+_M(A)$ as we
    can join the curves from $A$ to $p_n$ and then from $p_n$ to
    $q$. Now we find causal curves $\gamma_n$ from $p_n$ through
    $\Sigma$ entirely inside $J^+_M(A)$ giving us a point $r_n \in
    \Sigma \cap A$. Since these curves are in the cone $J^-_M(q)$ we
    have $r_n \in \Sigma \cap J^-_M(q)$. For a Cauchy hypersurface one
    knows that $\Sigma \cap J^-_M(q)$ is always compact. Thus also $A
    \cap \Sigma \cap J^-_M(q)$ is compact and hence the $r_n$ converge
    to some $r \in A$ after passing to a suitable subsequence. But
    then the curves $\gamma_n$ converge to some limiting curve
    $\gamma$ joining $r$ with $p$, see
    \cite[Lemma~14.14]{oneill:1983a} for details on the notion of
    limiting curves. By continuity $\gamma$ is still causal and thus
    $p \in J^+_M(A)$, see
    Figure~\ref{fig:causal-influence-of-closed-set-is-closed}. The
    argument for $J^-_M(A)$ is analogous.
\end{proof}

Later on we will be interested in those $u \in \Secinfty(E\at{U})$
where the initial values $u_0, \dot{u}_0 \in \Secinfty(\iota^\#E)$
have \emph{compact support} in $\Sigma$.

Let us now prove the uniqueness property of the Cauchy
problem. Lemma~\ref{lemma:solution-integral-formula} states that
\emph{locally} on $U$ the solution $u$ of the wave equation is
determined by its initial values $u_0$ and $\dot{u}_0$ on $\Sigma$,
since the left hand side of \eqref{eq:solution-integral-formula}
determines $u$ as a distribution and hence by the injective embedding
according to Remark~\ref{remark:generalized-densities} also as a
section. Thus we need to globalize this uniqueness statement.
\begin{theorem}
    \label{theorem:uniqueness-of-hom-equation-with-initial-values}
    \index{Spacetime!globally hyperbolic}%
    \index{Wave equation!unique solution}%
    \index{Initial conditions}%
    \index{Cauchy problem!uniqueness}%
    Let $(M,g)$ be globally hyperbolic and let $\iota: \Sigma
    \hookrightarrow M$ be a smooth spacelike Cauchy hypersurface with
    future directed normal vector field $\mathfrak{n} \in
    \Secinfty(\iota^\# TM)$. Assume that $u \in \Secinfty(E)$ is a
    solution to the wave equation $Du = 0$ with initial conditions
    \begin{equation}
        \label{eq:trivial-initial-values}
        u_0 = 0 = \dot{u}_0.
    \end{equation}
    Then
    \begin{equation}
        \label{eq:trivial-initial-values-give-trivial-solution}
        u = 0.
    \end{equation}
\end{theorem}
\begin{proof}
    First we note that by
    Theorem~\ref{theorem:globalhyp-and-cauchy-surface} there is a
    Cauchy temporal function $\mathfrak{t}$ on $M$ such that the level
    surface for $\mathfrak{t}=0$ coincides with $\Sigma$. We set
    \[
    \iota_t: \Sigma_t
    = \left\{
        p \in M \; \big| \; \mathfrak{t}(p) = t
    \right\}
    \hookrightarrow M
    \]
    for all times $t \in \mathbb{R}$. The gradient of $\mathfrak{t}$
    is by definition future directed and timelike and for a tangent
    vector $v_p \in T_p \Sigma_t$ we have $\D \mathfrak{t} \at{p}
    (v_p) = 0$ whence the gradient of $\mathfrak{t}$ is orthogonal to
    $T_p \Sigma_t$ at $p \in \Sigma_t$. Normalizing the gradient will
    give a globally defined vector field $\mathfrak{n} \in
    \Secinfty(TM)$ such that for every $t \in \mathbb{R}$ the
    restriction $\mathfrak{n}_t = \iota_t^\# \mathfrak{n} \in
    \Secinfty(\iota_t^\# TM)$ is the future directed normal vector
    field of $\Sigma_t$. Now let $p \in M$ be given and let $t_0 =
    \mathfrak{t}(p)$ be its time value, i.e. $p \in
    \Sigma_{t_0}$. Assume $t_0 > 0$ (the case $t_0 < 0$ is treated
    analogously). Then we define
    \[
    t_{\mathrm{max}}
    = \sup
    \left\{
        t \in [0, t_0]
        \; \big| \;
        u
        \; \textrm{vanishes on} \;
        J^-_M(p) \cap \cup_{0 \leq \tau \leq t} \Sigma_\tau
    \right\}.
    \]
    Since $u$ vanishes on $\Sigma_0$ this is well-defined and we have
    $0 \leq t_{\textrm{max}} \leq t_0$, see also
    Figure~\ref{fig:vanishing-time}.
    \begin{figure}
        \centering
        \input{vanishing-time.\pictype}
        \caption{\label{fig:vanishing-time}%
          The definition of $t_{\textrm{max}}$
        }
    \end{figure}
    The idea is now to show $t_{\textrm{max}} = t_0$ whence by
    continuity $u$ vanishes also at $p$. As $p$ was arbitrary this
    will imply $u=0$ everywhere for positive times. Then the analogous
    argument would give $u=0$ also for negative times. Thus let us
    assume the controversy, i.e. $t_{\textrm{max}} < t_0$. Let $q \in
    J^-_M(p) \cap \Sigma_{t_{\textrm{max}}}$, then we can find a small
    open causal neighborhood $U \subseteq U^\cl \subseteq U'$ of $q$
    such that on one hand we have our local fundamental solutions and
    on the other hand $U \cap \Sigma_{t_{\textrm{max}}}$ is still a
    Cauchy hypersurface. Note that this additional requirement can
    still be achieved, see
    e.g. \cite[Lem.~A.5.6]{baer.ginoux.pfaeffle:2007a}. In fact, the
    \Index{Cauchy development} $D(V)$ of a small enough open
    neighborhood $q \in V \subseteq \Sigma_{t_{\textrm{max}}}$ of $q$
    in $\Sigma_{t_{\textrm{max}}}$ will do the job, see also
    Remark~\ref{remark:4}. We consider the initial values of $u$ on
    this Cauchy hypersurface and denote them by $u_{t_\textrm{max}} =
    \iota_{t_\textrm{max}}^\# u$ and $\dot{u}_{t_\textrm{max}} =
    \iota_{t_\textrm{max}}^\# \nabla^E_{\mathfrak{n}(p)} u$ as usual.
    \begin{figure}
        \centering
        \input{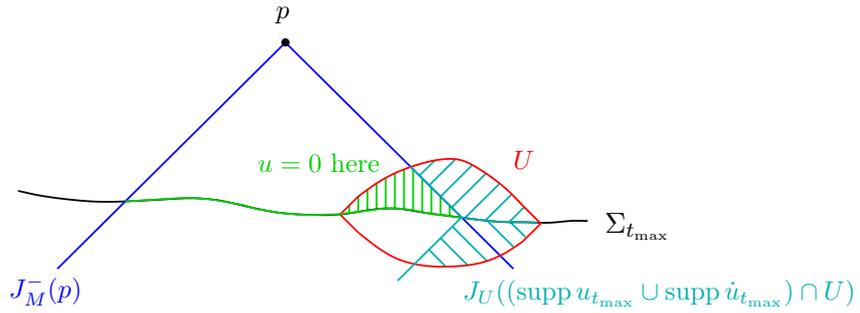}
        \caption{\label{fig:showing-that-u-is-zero-above-locally}%
          Showing that $u$ is zero locally above
          $\Sigma_{\textrm{max}}$.
        }
    \end{figure}
    From Lemma~\ref{lemma:supp-of-solution-and-supp-of-initial-values}
    we know that $u$ restricted to the small open subset $U$ has the
    following property
    \[
    \supp u \subseteq
    J_U \left(
        \supp u_{t_{\textrm{max}}} \cup \supp
        \dot{u}_{t_{\textrm{max}}} \cap U
    \right).
    \]
    Now by continuity and the choice of $t_{\textrm{max}}$ we know
    that $u_{t_{\textrm{max}}} = 0 = \dot{u}_{t_\textrm{max}}$ on
    $\Sigma_{t_{\textrm{max}}} \cap J^-_M(p)$. In particular,
    $u_{t_{\textrm{max}}} = 0 = \dot{u}_{t_{\textrm{max}}}$ in the
    open subset $U \cap \Sigma_{t_{\textrm{max}}} \cap J^-_M(p)$ of
    $\Sigma_{t_{\textrm{max}}}$, see
    Figure~\ref{fig:showing-that-u-is-zero-above-locally}. But then
    Lemma~\ref{lemma:supp-of-solution-and-supp-of-initial-values}
    shows that $u$ still vanishes on $J^-_M(p) \cap J^+_M(
    \Sigma_{t_{\textrm{max}}} \cap U)$, i.e. in this part of $U$ which
    is above $\Sigma_{t_{\textrm{max}}}$ and in the past of $p$. Since
    $J^-_M(p) \cap \Sigma_{t_{\textrm{max}}}$ is compact we can cover
    this part of the Cauchy hypersurface $\Sigma_{t_{\textrm{max}}}$
    with finitely many $U_1, \ldots, U_N$ for which the above argument
    applies.
    \begin{figure}
        \centering
        \input{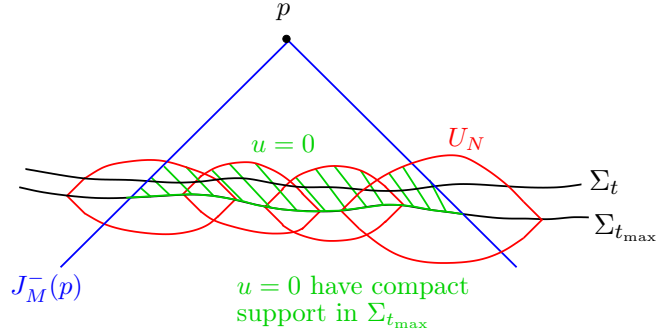}
        \caption{\label{fig:showing-that-u-is-zero-above}%
          Showing that $u$ is zero in a small neighborhood of
          $\Sigma_{t_{\textrm{max}}}$.
        }
    \end{figure}
    Now the union $U_1 \cup \ldots \cup U_N$ is an open neighborhood
    of $J^-_M(p) \cap \Sigma_{t_{\textrm{max}}}$ and hence $u$
    vanishes on this open subset $(U_1 \cup \ldots \cup U_N) \cap
    J^-_M(p) \cap J^+_M(\Sigma_{t_{\textrm{max}}})$ in the future of
    $J^+_M(\Sigma_{t_{\textrm{max}}})$. But this means that there is
    an $\epsilon > 0$ such that on $\Sigma_t \cap J^-_M(p)$ the
    section $u$ still vanishes for all $t \in [t_{\textrm{max}},
    t_{\textrm{max}} + \epsilon)$. This is in contradiction to the
    maximality of $t_{\textrm{max}}$ and hence $t_{\textrm{max}} =
    t_0$ whence $u(p) = 0$ by continuity. This shows that $u=0$ on
    $J^+_M(\Sigma)$ and an analogous argument gives $u=0$ on
    $J^-_M(\Sigma)$.
\end{proof}
As this is one of the central theorems we give an alternative proof of
the uniqueness statement. In particular, it will give some new insight
and an additional technique which turns out to be useful also at other
places.
\begin{altproof}[ of Theorem~\ref{theorem:uniqueness-of-hom-equation-with-initial-values}]
    Again we use a foliation of $M$ by smooth spacelike Cauchy
    hypersurfaces $\Sigma_t$ where for each $t \in \mathbb{R}$ the set
    $\Sigma_t$ is the level hypersurface of a Cauchy temporal function
    as before. We define now
    \[
    u^+(p)
    = \begin{cases}
        u(p) & \textrm{for } t(p) \leq 0 \\
        0    & \textrm{for } t(p) > 0
    \end{cases},
    \]
    and claim that this is a $\Fun[2]$-section still satisfying the
    wave equation $Du^+ = 0$.
    \begin{figure}
        \centering
        \input{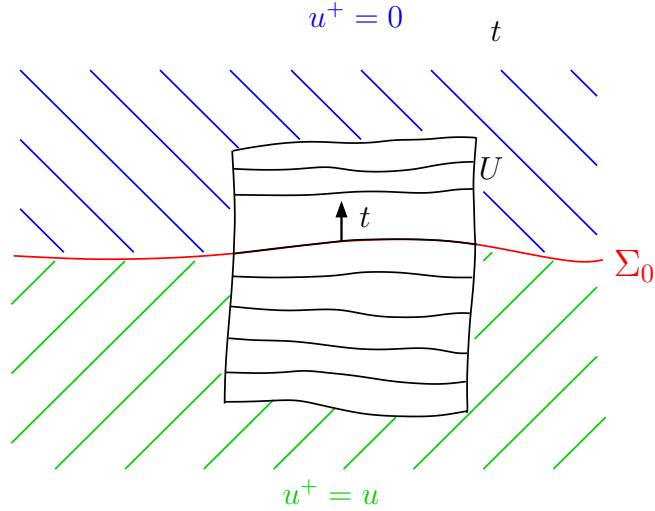}
        \caption{\label{fig:u-plus}%
          The neighborhood $U$.
        }
    \end{figure}
    Since $M = I^+_M(\Sigma_0) \cup \Sigma_0 \cup I^-_M(\Sigma_0)$
    with open $I^\pm_M(\Sigma_0)$ and $\Sigma_0$ the common boundary
    of $I^\pm_M(\Sigma_0)$ we can check the regularity of $u^+$ on
    each piece. Clearly on $I^\pm_M(\Sigma_0)$ we have $u^+
    \at{I^\pm_M(\Sigma_0)} \in \Secinfty(E \at{I^\pm_M(\Sigma_0)})$
    and $Du^+ \at{I^\pm_M(\Sigma_0)} = 0$. Thus we only have to check
    that $u^+$ is $\Fun[2]$ at $\Sigma_0$, then by continuity $D u^+ =
    0$ will follow everywhere. Thus let $p \in \Sigma_0$ and choose a
    small open neighborhood $V \subseteq \Sigma_0$ of $p$ allowing for
    local coordinates $x^1, \ldots, x^{n-1}$ and a trivialization of
    the bundled $E \at{\Sigma_0}$. By the splitting theorem we have an
    open neighborhood $U \subseteq M$ of $p$ such that the time
    function $t$ gives a diffeomorphism $U \simeq (-\epsilon,
    \epsilon) \times V$ and the metric $g\at{U}$ is given by
    \[
    g\at{U}
    = \beta \D\!t^2 - g_t
    \]
    with $\beta \in \Cinfty(U)$ positive and $g_t$ a smooth
    time-dependent metric on $\Sigma_0$, see
    Theorem~\ref{theorem:globalhyp-and-cauchy-surface}. In fact, we
    have this block diagonal structure even globally, see also
    Figure~\ref{fig:u-plus}. Now $u_0 = 0$ implies that $u^+$ is
    continuous at $\Sigma_0$. Moreover, all partial derivatives of $u$
    in $x^1, \ldots, x^{n-1}$ direction vanish on $\Sigma_0$ and hence
    the partial derivative of $u^+$ in $x^1, \ldots, x^{n-1}$
    directions are continuous as well. The block diagonal form of the
    metric shows that $\frac{\partial}{\partial t}$ is parallel to
    $\mathfrak{n}$ at $\Sigma_0$ whence the condition $\dot{u}_0 = 0 $
    means that the partial $\frac{\partial}{\partial t}$-derivative of
    $u$ vanishes at $\Sigma_0$. Indeed this differs (in our
    trivialization) from the covariant derivative by
    $\Cinfty(M)$-linear combinations of the components of $u_0$, which
    vanish by $u_0 = 0$. We conclude that $u^+$ is $\Fun[1]$. For the
    second derivative we first observe that the contributions
    $\frac{\partial^2}{\partial x^i \partial x^j} u $ all vanish on
    $\Sigma_0$ since $u_0 = 0$ is constant. Moreover, since $u$ is
    $\Fun[2]$, the contributions $\frac{\partial}{\partial t}
    \frac{\partial}{\partial x^i} u = \frac{\partial}{\partial x^i}
    \frac{\partial}{\partial t} u$ vanish on $\Sigma_0$ since
    $\frac{\partial}{\partial t} u = 0$ identically on $\Sigma_0$. For
    the last combination $\frac{\partial^2}{\partial t^2} u$ we have
    to use the wave equation. Locally the wave equation reads
    \[
    \left(
        \frac{1}{\beta} \frac{\partial^2}{\partial t^2}
        - g_t^{ij} \frac{\partial^2}{\partial x^i \partial x^j}
    \right) u
    + a \frac{\partial u }{\partial t}
    + b^i \frac{\partial u}{\partial x^i}
    + B u
    = 0,
    \]
    where $g_t^{ij}$ is the inverse metric to the metric $g_t$ on
    $\Sigma_t$, and $a, b^i, B$ are coefficient functions. Evaluating
    this on $\Sigma_0$ using the previous results gives
    $\frac{\partial^2 u }{\partial t^2} = 0$ on $\Sigma_0$. Thus the
    second partial derivatives are also continuous in this local
    chart. It follows that $u^+$ is $\Fun[2]$. By continuity it
    follows that $Du^+=0$ everywhere. But then
    Theorem~\ref{theorem:future-past-comp--hom-solution-is-zero} gives
    immediately $u^+ = 0$ since clearly $u^+$ has future compact
    support, see Figure~\ref{fig:u-plus-future-compact-supp}, and $M$
    being globally hyperbolic fulfills the conditions of
    Theorem~\ref{theorem:future-past-comp--hom-solution-is-zero}. But
    this implies $u \at{I^-_M(\Sigma_0)} = 0$. An analogous argument
    for
    \[
    u^-(p) = \begin{cases}
        0    & p \in I^-_M(\Sigma_0) \\
        u(p) & p \in I^+_M(\Sigma_0)
    \end{cases}
    \]
    shows that $u \at{I^+_M(\Sigma_0)} = 0$ as well.
    \begin{figure}
        \centering
        \input{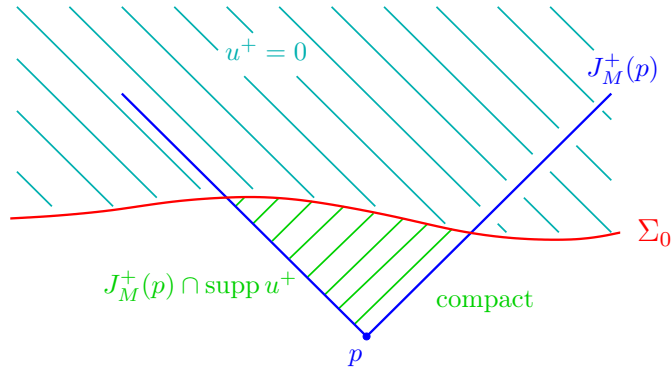}
        \caption{\label{fig:u-plus-future-compact-supp}%
          The section $u^+$ has future compact support.
        }
    \end{figure}
\end{altproof}

\begin{remark}
    \label{remark:alternative-proof}
    The alternative proof gives yet another interpretation of Cauchy
    hypersurfaces. They are the hypersurfaces $\Sigma$ along which
    solutions of the wave equation can be sewed together if they match
    on $\Sigma$. The argument in this approach will be used at several
    instances again.
\end{remark}
In view of the alternative proof we see that the uniqueness of the
solution to the Cauchy Problem is a direct consequence of
Theorem~\ref{theorem:future-past-comp--hom-solution-is-zero}
alone. The considerations in Section~\ref{satz:uniq-solut-cauchy}
before are not needed. Moreover, since
Theorem~\ref{theorem:future-past-comp--hom-solution-is-zero} works
even for distributional sections $u \in \Sec[-\infty](E)$ the
regularity needed for the uniqueness is actually much smaller than
$\Cinfty$:
\begin{theorem}
    \label{theorem:uniqueness-initial-condition-C2-case}
    \index{Spacetime!globally hyperbolic}%
    \index{Wave equation!unique solution}%
    \index{Initial conditions}%
    \index{Cauchy problem!uniqueness}%
    Let $(M,g)$ be globally hyperbolic and let $\iota: \Sigma
    \hookrightarrow M$ be a smooth spacelike Cauchy hypersurface with
    future directed normal vector field $\mathfrak{n} \in
    \Secinfty(\iota^\# TM)$. Let $v \in \Sec[0](E)$ be a continuous
    section and $u \in \Sec[2](E)$ a $\Fun[2]$-section satisfying the
    inhomogeneous wave equation
    \begin{equation}
        \label{eq:inhom-wave-eq-small-differentiability}
        Du = v.
    \end{equation}
    Then $u$ is uniquely determined by its initial conditions $u_0 =
    \iota^\# u$ and $\dot{u}_0 = \iota^\# \nabla^E_{\mathfrak{n}} u$
    on $\Sigma$.
\end{theorem}
\begin{proof}
    Requiring $u \in \Sec[2](E)$ is the minimal requirement to view
    \eqref{eq:inhom-wave-eq-small-differentiability} as a
    \emph{pointwise} equation. In fact, since continuous sections
    still embed into $\Sec[-\infty](E)$ we also have $Du = v$ in the
    sense of distributional sections. Suppose $\widetilde{u} \in
    \Sec[2](E)$ is an alternative solution with the same initial
    conditions. Then $u - \widetilde{u}$ is a $\Fun[2]$-solution of
    the homogeneous wave equation. For this we can repeat the argument
    from the alternative proof of
    Theorem~\ref{theorem:uniqueness-of-hom-equation-with-initial-values}
    since we only needed $\Fun[2]$ there. Thus $u - \widetilde{u} = 0$
    as distributions by
    Theorem~\ref{theorem:uniqueness-of-hom-equation-with-initial-values}
    and hence $u - \widetilde{u} = 0$ as $\Fun[2]$-sections as well.
\end{proof}

%
%

\subsection{Existence of Local Solutions to the Cauchy Problem}
\label{satz:existence-local-solutions}

After the uniqueness we pass to the existence of solutions to the
Cauchy problem. We will assume that the Cauchy data as well as the
inhomogeneity of the wave equation have compact support.

The first statement is still a local result to the Cauchy problem:
\begin{proposition}
    \label{proposition:existence-of-local-solutions-with-initial-values}
    Let $(M,g)$ be a time-oriented Lorentz manifold with a smooth
    spacelike hypersurface $\iota: \Sigma \hookrightarrow M$ with
    future directed normal vector field $\mathfrak{n}$. Moreover, let
    $U \subseteq U^\cl \subseteq U'$ be a sufficiently small causal
    open subset of $M$ such that $\Sigma \cap U \hookrightarrow U$ is
    a Cauchy hypersurface for $U$. Then there exists a unique solution
    $u \in \Secinfty(E\at{U})$ for given initial values $u_0,
    \dot{u}_0 \in \Secinfty_0(\iota^\# E \at{U})$ and given
    inhomogeneity $v \in \Secinfty_0(E\at{U})$ of the inhomogeneous
    wave equation
    \begin{equation}
        \label{eq:inhom-wave-equ-local}
        Du = v
    \end{equation}
    with $\iota^\# u = u_0$ and $\iota^\# \nabla_{\mathfrak{n}}^E u =
    \dot{u}_0$. In addition we have
    \begin{equation}
        \label{eq:supp-of-local-solution}
        \supp u \subseteq
        J_M( \supp u_0 \cup \supp \dot{u}_0 \cup \supp v).
    \end{equation}
\end{proposition}
\begin{proof}
    As usual, sufficiently small means that we have our local
    fundamental solutions and therefor the result of
    Chapter~\ref{cha:LocalTheory}. The uniqueness of $u$ follows
    directly from
    Theorem~\ref{theorem:uniqueness-of-hom-equation-with-initial-values}.
    We can apply the splitting theorem for globally hyperbolic
    manifolds in the form of
    Theorem~\ref{theorem:globalhyp-and-cauchy-surface} to $U$, see
    also \cite[Thm.~2.78]{minguzzi.sanchez:2006a:pre}. Thus we find a
    Cauchy temporal function $t$ on $U$ inducing an isometry of $U$ to
    $\mathbb{R} \times (\Sigma \cap U)$ such that the metric becomes
    $\beta \D\!t^2 - g_t$ with $\beta \in \Cinfty(U)$ positive and
    $g_t$ a time dependent Riemannian metric on $\Sigma \cap U$. Every
    $t$-level surface is Cauchy and we have the normal vector field
    \[
    \mathfrak{n}
    = \frac{1}{\sqrt{\beta}} \frac{\partial}{\partial t}
    \in \Secinfty(TU),
    \]
    which is normal to every level surface. Moreover, since by
    definition $U \subseteq U'$ is contained in a \emph{convex} domain
    $U'$ the vector bundle $E$ is trivializable over $U'$ and hence
    over $U$. Therefore we can choose a frame $\{e_\alpha\}$ over $U$
    of $E \at{U}$ and write $u = u^\alpha e_\alpha$ with smooth
    functions $u^\alpha \in \Cinfty(U)$ for every $u \in
    \Secinfty(E\at{U})$. This allows to identify a section $u$ with a
    collection of scalar function $u^\alpha$. The normally hyperbolic
    operator $D$ is now of the form
    \[
    D = \frac{1}{\beta} \frac{\partial^2}{\partial t^2}
    + \widetilde{D},
    \tag{$*$}
    \]
    where $\widetilde{D}$ contains at most first $t$-derivatives,
    still up to second derivatives in $\Sigma$-directions, and it has
    matrix-valued coefficient functions with respect to our
    trivialization induced by the $e_\alpha$. We claim now that the
    initial conditions together with the wave equation determine all
    $t$-derivatives of a solution along $\Sigma$. The argument is
    similar to the proof of
    Theorem~\ref{theorem:uniqueness-of-hom-equation-with-initial-values}.
    Suppose $u$ is a smooth solution of $Du = v$ with initial
    conditions $u_0$ and $\dot{u}_0$. We already know that $\dot{u}_0$
    is determined by $u_0$ and $\frac{\partial u}{\partial
      t}\at{\Sigma}$ and conversely $\frac{\partial u}{\partial
      t}\at{\Sigma}$ is determined by $\dot{u}_0$ and $u_0$. Using
    ($*$) we see that
    \[
    \frac{\partial^2 u}{\partial t^2}
    = \beta (Du - \widetilde{D}u)
    = \beta (v - \widetilde{D}u).
    \tag{$**$}
    \]
    This shows that $\frac{\partial^2 u}{\partial t^2}\at{\Sigma}$ is
    determined by $u_0$ and $\frac{\partial u}{\partial
      t}\at{\Sigma}$, namely we have
    \[
    \frac{\partial^2 u}{\partial t^2}\At{\Sigma}
    = (\beta v)\At{\Sigma} - ( \beta\widetilde{D}u)\At{\Sigma},
    \]
    where the right hand side uses only $u_0$ and $\frac{\partial
      u}{\partial t}\at{\Sigma}$ since $\widetilde{D}$ is at most of
    first order in the $t$-variable. Moreover, differentiating ($**$)
    $j$-times we get
    \[
    \frac{\partial^{j+2} u}{\partial t^{j+2}}
    = \frac{\partial^j (\beta v)}{\partial t^j}
    - \frac{\partial^j}{\partial t^j} (\beta \widetilde{D}u).
    \]
    Hence on $\Sigma$ we have
    \[
    \frac{\partial^{j+2} u}{\partial t^{j+2}} \At{\Sigma}
    = \frac{\partial^j (\beta v)}{\partial t^j} \At{\Sigma}
    - \frac{\partial^j}{\partial t^j} (\beta \widetilde{D}u)
    \At{\Sigma}.
    \tag{$*{*}*$}
    \]
    We see that the right hand side is a $\Cinfty(\Sigma)$-linear
    combination of the $u_0, \frac{\partial u}{\partial t}\at{\Sigma},
    \ldots, \frac{\partial^{j+1} u}{\partial t^{j+1}} \at{\Sigma}$
    plus an affine term $\frac{\partial^j (\beta v)}{\partial t^j}
    \at{\Sigma}$. Thus by induction we conclude that all
    $t$-derivatives of $u$ on $\Sigma$ are determined by $u_0$ and
    $\frac{\partial u}{\partial t}\at{\Sigma}$, and of course by the
    choice of the inhomogeneity $v$. Moreover, since we have a
    $\Cinfty(\Sigma)$-affine linear combination we conclude that
    \[
    \supp \left(
        \frac{\partial^j u}{\partial t^j}\at{\Sigma}
    \right)
    \subseteq
    \left(
        \underbrace{\supp u_0 \cup \supp \dot{u}_0 \cup \supp v}_{K}
    \right)
    \cap \Sigma
    = K \cap \Sigma
    \]
    is contained in a compact subset $K' = K \cap \Sigma$ of $\Sigma$
    for all $j$. Now we use these recursion formulas to \emph{define}
    sections $u_j \in \Secinfty(E\at{\Sigma})$ by ($*{*}*$) for all $j
    \geq 2$. First we note that we indeed can find a global section
    $\widetilde{u} \in \Secinfty(E)$ whose $t$-derivatives on $\Sigma$
    are given by the $u_j$: this is essentially a consequence of the
    \Index{Borel Lemma} for Fr\'echet spaces, see
    e.g. \cite[Satz~5.3.33]{waldmann:2007a}. For convenience we repeat
    the argument here: We choose a cut-off function $\chi \in
    \Cinfty_0(\mathbb{R})$ with $\supp \chi \subseteq [-1,1]$ and
    $\chi \at{[-\frac{1}{2},\frac{1}{2}]} = 1$. As we did frequently
    in Section~\ref{satz:LocalFundamentalSolution} we consider as
    Ansatz a series
    \[
    \widetilde{u}(t,p)
    = \sum_{j=0}^\infty
    \chi\left(\frac{t}{\epsilon_j}\right)
    \frac{t^j}{j!} \: u_j(p)
    \tag{$\star$}
    \]
    with numbers $0 < \epsilon_j \leq 1$ yet to be chosen. We want to
    choose them in such a way that the series converges in the
    $\Cinfty$-topology of $\Secinfty(E\at{U})$. Clearly, each term has
    support in $[-1,1] \times K'$ whence we only have to
    consider the seminorms of $\Secinfty(E\at{U})$ estimating
    derivatives on this compactum. It is clear from the Ansatz and the
    properties of $\chi$ that if we have $\Cinfty$-convergence then
    $\frac{\partial^j \widetilde{u}}{\partial t^j} \at{t=0} = u_j$ for
    all $j$. Thus let us estimate the $k$-th seminorm
    $\seminorm[{[-1,1] \times K', k}]$ of each term of ($\star$). With
    the usual Leibniz rule and the fact that the seminorms factorize
    on factorizing functions we get from
    Lemma~\ref{lemma:seminorms-of-cutoff-function}
    \[
    \seminorm[ {[-1,1] \times K', k} ]
    \left(
        \chi\left(\frac{t}{\epsilon_j}\right) \frac{t^j}{j!} \: u_j
    \right)
    \leq
    \frac{\epsilon_j}{j!}
    \seminorm[ {[-1,1],k} ](\chi) \seminorm[K', k](u_j).
    \]
    This allows to choose the $\epsilon_j$ such that
    \[
    \epsilon_j \max_{k \leq j}
    \seminorm[{[-1,1],k}](\chi) \seminorm[K' ,k](u_j)
    < 1.
    \]
    Then the series ($\star$) converges in the $\Fun$-norm
    $\seminorm[{[-1,1]\times K', k}]$ absolutely as the first terms do
    not spoil the convergence. Thus we have absolute
    $\Cinfty$-convergence in total. This shows the existence of a
    $\widetilde{u} \in \Secinfty(E\at{U})$ with
    \[
    \frac{\partial^j \widetilde{u}}{\partial t^j}\At{\Sigma}
    = u_j
    \]
    and
    \[
    \supp \widetilde{u}
    \subseteq J_M(K).
    \]
    Indeed, the last claim follows from the fact that $\supp
    \widetilde{u} \subseteq [-1,1] \times K'$ and $\mathbb{R} \times
    K' \subseteq J_M(K)$ since for every $(t,p) \in \mathbb{R} \times
    K'$ the curve $\tau \mapsto (\tau,p)$ connects $(0,p)$ to $(t,p)$
    and the curve is clearly timelike. This follows from the splitting
    of the metric, see also
    Figure~\ref{fig:easy-timelike-curve-since-splitting}.
    \begin{figure}
        \centering
        \input{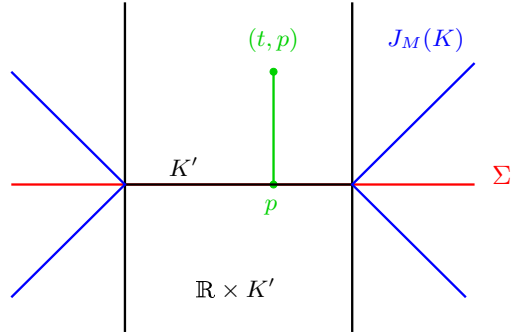}
        \caption{\label{fig:easy-timelike-curve-since-splitting}%
          The splitting yields simple timelike curves.
        }
    \end{figure}
    From the construction of $\widetilde{u}$ we see that $D
    \widetilde{u}$ coincides with $v$ including \emph{all} time
    derivatives on $\Sigma$. In other words, $D \widetilde{u} - v$
    vanishes on $\Sigma$ up to infinite order. Thus we can consider
    the definition
    \[
    w_\pm
    = \begin{cases}
        D \widetilde{u} - v & \textrm{on} \; I^\pm_M(\Sigma) \\
        0                   & \textrm{on} \; J^\mp_M(\Sigma),
    \end{cases}
    \]
    which gives a \emph{smooth} section $w_\pm \in
    \Secinfty(E\at{U})$. Since both $\widetilde{u}$ and $v$ have
    compact support, also $w_\pm$ is compactly supported. Thus we can
    solve the inhomogeneous wave equation
    \[
    D \widetilde{\widetilde{u}}_\pm = w_\pm
    \]
    on the open subset $U$ according to
    Theorem~\ref{theorem:dual-of-fund-solution-gives-inhom-solution}
    with a smooth solution $\widetilde{\widetilde{u}}_\pm \in
    \Secinfty(E\at{U})$ such that $\supp \widetilde{\widetilde{u}}_\pm
    \subseteq J^\pm_U( \supp w_\pm)$. Since $\supp w \subseteq ( \supp
    D\widetilde{u} \cup \supp v ) \cap J^+(\Sigma) \subseteq J^+_M(K)$
    we conclude $J^+_M(\supp w) \subseteq J^+_M(K)$.
    \begin{figure}
        \centering
        \input{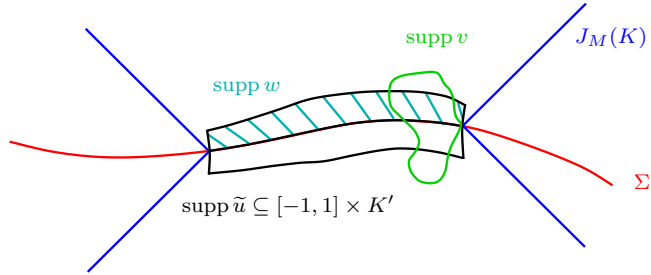}
        \caption{\label{fig:the-supports-of-the-several-sections}%
          The supports of the several sections in the proof of
          Proposition~\ref{proposition:existence-of-local-solutions-with-initial-values}.
        }
    \end{figure}
    This shows that $\supp \widetilde{\widetilde{u}}_\pm \subseteq
    J^\pm_M(K) \cap U = J^\pm_U(K)$. In particular,
    $\widetilde{\widetilde{u}} \at{J^\mp_M(\Sigma)} = 0$. Now we
    consider the smooth section $u_\pm \in \Secinfty(E\at{U})$ defined
    by
    \[
    u_\pm = \widetilde{u} - \widetilde{\widetilde{u}}_\pm.
    \]
    Since $\widetilde{\widetilde{u}}_\pm$ vanishes on
    $J^\mp_M(\Sigma)$ we have $u_\pm = \widetilde{u}$ on
    $J^\mp_M(\Sigma)$. In particular, $u_\mp$ coincides with
    $\widetilde{u}$ up to all orders on $\Sigma$ by continuity of the
    $t$-derivatives. Thus $u_\pm$ satisfies the correct initial
    conditions. Moreover, on $I^\pm_U(\Sigma)$ we have
    \[
    Du_\pm \at{I^\pm_U(\Sigma)}
    = D \widetilde{u} \at{I^\pm_U(\Sigma)}
    - D \widetilde{\widetilde{u}}_\pm \at{U^\pm_U(\Sigma)}
    = (w+v) \at{I^\pm_U(\Sigma)} - w\at{I^\pm_U(\Sigma)}
    = v \at{I^\pm_U}(\Sigma),
    \]
    whence on this open part of $U$ the section $u_\pm$ solves the
    inhomogeneous wave equation. Since both $u_+$ and $u_-$ agree on
    $\Sigma$ up to infinite orders, as they agree with
    $\widetilde{u}$, we can glue them together and set
    \[
    u = \begin{cases}
        u_+ & \textrm{on} \; I^+_U(\Sigma) \\
        u_- & \textrm{on} \; I^-_U(\Sigma).
    \end{cases}
    \]
    On one hand, this yields a smooth section $u \in
    \Secinfty(E\at{U})$ on all of $U$. Moreover, $u$ solves the
    inhomogeneous wave equation on both open parts $I^\pm_U(\Sigma)$
    and hence on all of $U$ by continuity. Finally, we know that
    \[
    \supp (u_\pm)
    \subseteq \supp \widetilde{u} \cup
    \supp \widetilde{\widetilde{u}}_\pm
    \subseteq J_U(K) \cup J^\pm_U(K) = J_U(K),
    \]
    whence also $\supp u \subseteq J_U(K)$. This completes the proof.
\end{proof}

We can refine the above argument for finite order of
differentiability. Here on one hand the Borel-Lemma is not needed as
we can simply take a polynomial in $t$ multiplied by the cut-off
function in order to have compact support. On the other hand, we have
to count orders of differentiation carefully:
\begin{proposition}
    \label{proposition:local-Ck-solutions}
    Let $k \geq 2$. Under the same general assumptions as in
    Proposition~\ref{proposition:existence-of-local-solutions-with-initial-values}
    we assume to have initial values $u_0 \in
    \Sec[2(k+n+1)+2]_0(\iota^\# E\at{U})$, $\dot{u}_0 \in
    \Sec[2(k+n+1)+1]_0(\iota^\# E\at{U})$ and an inhomogeneity $v \in
    \Sec[2(k+n+1)]_0(E\at{U})$. Then there exists a unique solution $u
    \in \Sec(E\at{U})$ of the inhomogeneous wave equation
    \begin{equation}
        \label{eq:Ck-wave-equation}
        Du = v
    \end{equation}
    with initial conditions $\iota^\#\upsilon = u_0$ and $\iota^\#
    \nabla_{\mathfrak{n}}u = \dot{u}_0$. For the support we still have
    \begin{equation}
        \label{eq:supp-Ck-solution}
        \supp u \subseteq
        J_M( \supp u_0 \cup \supp \dot{u}_0 \cup \supp v).
    \end{equation}
\end{proposition}
\begin{proof}
    As in the proof of
    Proposition~\ref{proposition:existence-of-local-solutions-with-initial-values}
    we define the sections $u_j$ recursively by
    \[
    \frac{\partial^{j+1} u}{\partial t^{j+1}} \at{\Sigma}
    = \frac{\partial^j (\beta v)}{\partial t^j} \at{\Sigma}
    - \frac{\partial^j}{\partial t^j} (\beta \widetilde{D}u)
    \at{\Sigma}
    \tag{$*$}
    \]
    with $u_j = \frac{\partial^j u}{\partial t^j}\at{\Sigma}$. Since
    for the right hand side we only have up to $j+1$ time derivatives
    we need $u_0, u_1, \ldots, u_{j+1}$ in order to determine
    $u_{j+2}$. In the local coordinates on $U$ we split the operator
    $\widetilde{D}$ into $\widetilde{D} = D_2 + D_1
    \frac{\partial}{\partial t}$ where $D_2, D_1$ are operators
    differentiating only in spacial directions. The coefficients of
    $D_2, D_1$ depend on all variables and $D_2$ is of order two while
    $D_1$ is of order one. Then the recursion ($*$) for $u = u_0 + t
    u_1 + \frac{t^2}{2} u_2 + \ldots$ can be written as
    \begin{align*}
        u_{j+1}
        &= \frac{\partial^j}{\partial t^j}(\beta v) \At{t=0}
        - \frac{\partial^j}{\partial t^j} \beta D_2
        \sum_{k=0}^j \frac{t^k}{k!} u_k \At{t=0}
        - \frac{\partial^j}{\partial t^j} \beta D_1
        \frac{\partial}{\partial t}
        \sum_{k=0}^j \frac{t^{k+1}}{(k+1)!} u_{k+1} \At{t=0} \\
        &= \frac{\partial^j}{\partial t^j}(\beta v)\At{t=0}
        - \sum_{k=0}^j \binom{j}{k}
        \frac{\partial^{j-k}}{\partial t^{j-k}} (\beta D_2) \At{t=0}
        u_k
        - \sum_{k=0}^j \binom{j}{k}
        \frac{\partial^{j-k}}{\partial t^{j-k}} (\beta D_1) \At{t=0}
        u_{k+1}
        \tag{$**$}.
    \end{align*}
    Note that $\frac{\partial^{j-k}}{\partial t^{j-k}} (\beta D_2)
    \at{t=0}$ is again a differential operator of order two while
    $\frac{\partial^{j-k}}{\partial t^{j-k}} (\beta D_1) \at{t=0}$ is
    of order one. This determines $u_{j+2}$ recursively in terms of
    spacial derivatives of $u_0, \ldots, u_{j+1}$. We claim that
    $u_{j+2}$ contains at most $j+2$ derivatives of $u_0$, at most
    $j+1$ derivatives of $u_1$ and at most $j$ derivatives of
    $v$. Indeed, for $j=0$ we have
    \[
    u_2
    = \beta v \at{\Sigma}
    - \beta D_2 \at{\Sigma} u_0
    - \beta D_1 \at{\Sigma} u_1,
    \]
    which shows the claim for this $j$. By inductions we see from
    ($**$) that $\frac{\partial^{j-k}}{\partial t^{j-k}} (\beta D_2)
    \at{t=0} u_k$ contains at most $k+2$ derivatives of $u_0$ and
    hence at most $j+2$ derivatives since $k=0, \ldots, j$. Moreover,
    it contains at most $k+1$ derivatives of $u_1$ and hence at most
    $k+1 \leq j+1$. Finally it contains at most $k-1$ derivatives of
    $v$ and thus also here things match. For the second sum one
    proceeds analogously. Finally, the first term gives $j$
    derivatives of $v$, which also matches our claim. Now assume we
    are give $u_0$ and $u_1$ of class $\Fun[2(k+n+1)+2]$ and
    $\Fun[2(k+n+1)+1]$, respectively. Moreover, suppose $v \in
    \Sec[2(k+n+1)]_0(E\at{U})$. Then the $u_j$ defined by the
    recursion ($*$) are of class $\Fun[2(k+n+1)+2-j]$. Thus the finite
    sum
    \[
    \widetilde{u}(t,p)
    = \chi(t) \sum_{j=0}^{k+n+1} \frac{t^j}{j!} u_j(p)
    \]
    gives a section of class at least $\Fun[k+n+1+2]$. Moreover, the
    recursion shows that $D \widetilde{u} - v$ vanishes up to order
    $t^{k+n+1}$. Thus gluing this with zero gives a section
    \[
    w_\pm
    = \begin{cases}
        D \widetilde{u} - v & \textrm{on } I^\pm_U(\Sigma) \\
        0                   & \textrm{else},
    \end{cases}
    \]
    which is still of class $\Fun[k+n+1]$ everywhere. Then
    $\widetilde{\widetilde{u}}_\pm$ is of class $\Fun$ by
    Theorem~\ref{theorem:dual-of-fund-solution-gives-inhom-solution}
    and thus $u_\pm$ are both of class $\Fun$. Since
    $\widetilde{\widetilde{u}}_\pm$ is $\Fun$ and vanishes on the open
    subset $I^\mp_U(\Sigma)$, the $u_\pm = \widetilde{u} -
    \widetilde{\widetilde{u}}_\pm$ agree with $\widetilde{u}$ on
    $\Sigma$ up to order $t^k$. Thus also the glued solution $u$ is of
    class $\Fun$ as claimed. The statement about the support is
    analogous to the smooth case.
\end{proof}

\begin{remark}
    \label{remark:integral-formula}
    Having Lemma~\ref{lemma:solution-integral-formula} in mind, it is
    tempting to define the solution of the Cauchy problem (at least in
    the homogeneous case $v=0$) by the formula
    \eqref{eq:solution-integral-formula}: Using instead of a test
    section $\varphi$ a $\delta$-functional at $p$ would directly give
    \begin{equation}
        \label{eq:integral-fomula-wish}
        u(p)
        = \int_{\Sigma}
        \left(
            \nabla^E_{\mathfrak{n}} G_U'(\delta_p) \at{\sigma}
            \cdot u_0(\sigma)
            - G_U'(\delta_p)\at{\sigma} \cdot \dot{u}_0(\sigma)
        \right) \mu_\Sigma(\sigma).
    \end{equation}
    However, here we face two problems. First one has to shows that
    $u$ is indeed a solution of $Du=0$ with the correct initial
    conditions. Second, and more severe, one has to justify the
    restriction of the \emph{distributions} $\nabla^E_{\mathfrak{n}}
    G_U'(\delta_p)$ and $G_U'(\delta_p)$ to the hypersurface, which is
    indeed a nontrivial task. Thus we leave
    \eqref{eq:integral-fomula-wish} as a heuristic formula and stay
    with
    Proposition~\ref{proposition:existence-of-local-solutions-with-initial-values}
    and Proposition~\ref{proposition:local-Ck-solutions}.
\end{remark}

%
%

\subsection{Existence of Global Solutions to the Cauchy Problem}
\label{satz:existence-global-solutions}

\begin{figure}
    \centering
    \input{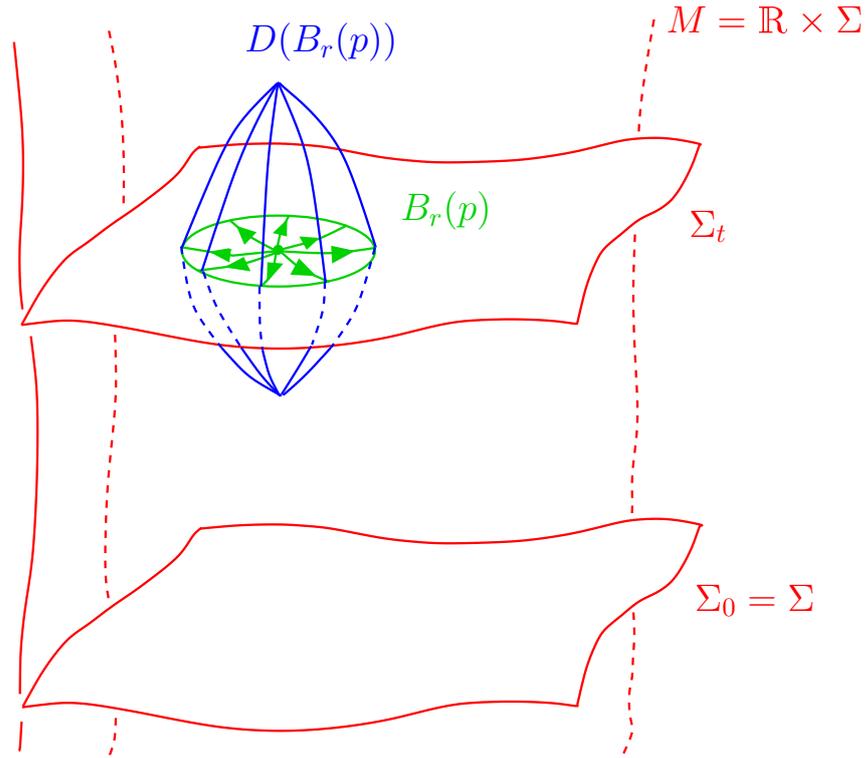}
    \caption{\label{fig:open-balls-in-Sigma}%
      An open ball $B_r(p)$ in a Cauchy hypersurface $\Sigma_t$ and
      its Cauchy development $D(B_r(p))$.
    }
\end{figure}

To approach the global existence of solutions we assume as before that
$M$ is globally hyperbolic with a smooth spacelike Cauchy hypersurface
$\Sigma$. Now we again use the splitting theorem $M \cong \mathbb{R}
\times \Sigma$ with the first coordinate being the Cauchy temporal
function and $\Sigma_t$ the Cauchy hypersurface of constant time $t$
where we shift the origin to $\Sigma_0 = \Sigma$. For every $p \in M$
we have a unique time $t$ with $p \in \Sigma_t$. On each $\Sigma_t$ we
have a Riemannian metric $g_t$ such that $g = \beta \D\!t^2 -
g_t$. This allows to speak of the open balls around $p \in \Sigma_t$
of radius $r > 0$ with respect to this metric $g_t$. We denote these
by $B_r(p)$ without explicit reference to $t$. Note that $B_r(p)
\subseteq \Sigma_t$ is open in $\Sigma_t$ but not in $M$, see also
Figure~\ref{fig:open-balls-in-Sigma}.
Here we use the Riemannian distance $d_{g_t}$ in $\Sigma_t$ with
respect to $g_t$ for defining the ball, i.e.
\begin{equation}
    \label{eq:open-balls-in-cauchyhyp}
    d_{g_t}(p,q)
    = \inf \left\{
        \left.
            \int_a^b g_t(\dot{\gamma}(\tau), \dot{\gamma}(\tau))
            \D\tau
        \right|
        \; \gamma(a) = p, \gamma(b) = q, \gamma(\tau) \in \Sigma_t
    \right\},
\end{equation}
where $\gamma$ is an at least piecewise $\Fun[1]$ curve joining $p,q
\in \Sigma_t$ \emph{inside} $\Sigma_t$. Having such a ball we consider
its \Index{Cauchy development} $D_M(B_r(p)) = D_M^+(B_r(P)) \cup
D_M^-(B_r(p))$ in $M$ according to
Definition~\ref{definition:cauchy-development}, see again
Figure~\ref{fig:open-balls-in-Sigma}. We now want to find $r$ small
enough that $D_M(B_r(p))$ is a nice open neighborhood of $p$ allowing
a local fundamental solution: in this case we call an open
neighborhood a \emph{relatively compact causal open neighborhood of
  small volume} or short \emIndex{RCCSV} for abbreviation. We start
with a couple of technical lemmas, following
\cite{baer.ginoux.pfaeffle:2007a}:
\begin{lemma}
    \label{lemma:sup-over-RCCSV-radi}
    The function $\rho:M \longrightarrow (0,+\infty]$ defined by
    \begin{equation}
        \label{eq:sup-over-RCCSV-radi}
        \rho(p)
        = \sup \left\{
            r > 0
            \; \big| \;
            D(B_r(p)) \textrm{ is RCCSV}
        \right\}
    \end{equation}
    is well-defined and lower semi-continuous.
\end{lemma}
\begin{proof}
    We have to show first that the set of $r>0$ with $D(B_r(p))$ RCCSV
    is non-empty. To this end we choose an RCCSV neighborhood $U
    \subseteq U^\cl \subseteq U'$ as before.
    \begin{figure}
        \centering
        \input{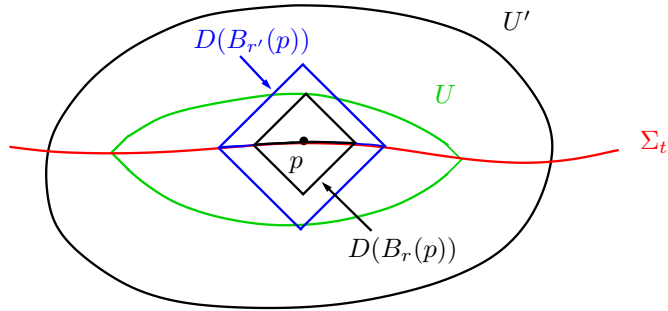}
        \caption{\label{fig:construction-of-RCCSV-radius-1} Illustration
          for the proof of Lemma~\ref{lemma:sup-over-RCCSV-radi}.}
    \end{figure}
    Then $U \cap \Sigma_t$ will be an open neighborhood of $p$ in
    $\Sigma_t$ hence it contains a $B_r(p) \subseteq \Sigma_t$. The
    problem might be that the Cauchy development of $B_r(p)$ may reach
    too far outside of $U$ or even $U'$ such that it is not RCCSV for
    free, see Figure~\ref{fig:construction-of-RCCSV-radius-1}. In
    fact, we have to choose a small enough $r$ such that $D(B_r(p))
    \subseteq U$. In this case it is causal in $U'$ and has small
    enough volume. We choose points $q^\pm \in U$ with $p \in
    J^\pm_U(q^\mp)$. Then we consider the open subset $J_U^+(q-) \cap
    J_U^-(q^+)$ which is an neighborhood of $p$, see
    Figure~\ref{fig:construction-of-RCCSV-radius-2}.
    \begin{figure}
        \centering
        \input{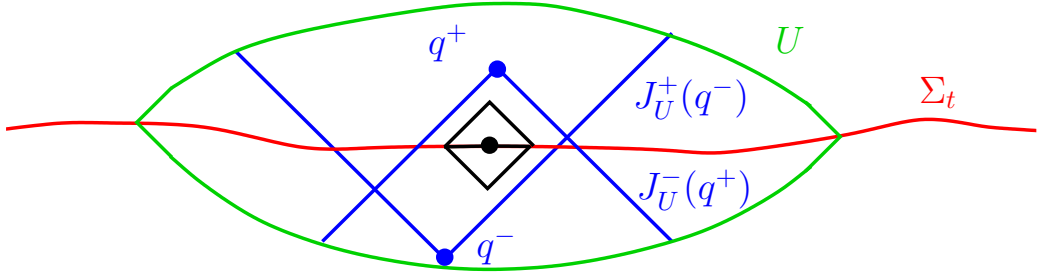}
        \caption{\label{fig:construction-of-RCCSV-radius-2}%
          Constructing a small enough open ball around $p$.
        }
    \end{figure}
    The intersection of this neighborhood of $p$ (in $M$) with
    $\Sigma_t$ gives an open neighborhood of $p$ in $\Sigma_t$. Now we
    choose a $B_r(p)$ contained in this neighborhood. We claim that
    $D_M(B_r(p))$ is in $U$. First we note that $J_U^\pm(q^\mp) = U
    \cap J_M^\pm(q^\mp)$ since $U$ is causally compatible with
    $M$. Now if $q \in D_M^+(B_r(p))$ then every past-inextensible
    causal curve meets $B_r(p)$. We claim that $q \in
    I^-(q^+)$. Assume that this is not the case. Then we have a
    past-inextensible curve from $p$ to $q$ which has to pass through
    the backward light cone of $q^+$. Denote this intersection point
    by $q_0$. Since we are inside a geodesically convex neighborhood
    $U'$, we can take the unique lightlike geodesic from $q^+$ to this
    $q_0$ which is past directed. Since this geodesic is on the light
    cone, it hits the Cauchy hypersurface $\Sigma_t$ \emph{not} in the
    open subset $I^-_U(q^+)$ but on its boundary, say in the point
    $q_1$. Thus it will not intersect the even smaller open ball
    $B_r(p)$. Thus the combined curve from $q$ back to $q_0$ and then
    back to $q_1$ will never hit $B_r(p)$, no matter how we extend it
    further in past directions. This contradicts $q \in D_M^+(p)$
    whence we conclude that $q \in I^-_U(q^+)$. A simpler argument
    shows that $q$ is also in the chronological future of $q^-$ and
    hence in the intersection of the two open subsets $I^+_U(q^-)$ and
    $I^-_U(q^+)$. An analogous argument shows that a point in
    $D_M^-(p)$ is also in this intersection. We finally arrived at the
    desired statement that $D_M(B_r(p))$ is in $U$.

    Now let $p \in M$ and $r > 0$ with $\rho(p) > r$ be given. In
    particular $D_M(B_r(p))$ will be RCCSV. Then we have to show that
    for a given $\epsilon > 0$ we have
    \[
    \rho(p') > r - \epsilon
    \]
    for all $p'$ in an appropriate open neighborhood of $p$. We
    consider the following function defined for $p' \in D_M(B_r(p))$ by
    \[
    \lambda(p')
    = \sup \left\{
        r' > 0 \; \big| \; B_{r'}(p') \subset D_M(B_r(p))
    \right\},
    \]
    i.e. we ask for the balls around $p'$ to be contained in the
    Cauchy development of $B_r(p)$. Note that $p'$ may correspond to a
    different time $t' \neq t$ which has to be taken into account in
    the definition of the radius $r'$, i.e. we use $g_{t'}$. We claim
    that there is an open neighborhood $V$ of $p$ such that for all
    $p' \in V$ we have
    \[
    \lambda(p') > r - \epsilon.
    \]
    Assume that this is not true. Then we can find a sequence $p_n
    \longrightarrow p$ of points in $D_M(B_r(p))$ with $\lambda(p_n)
    \leq r - \epsilon$ for all $n$. Then it follows that for $r' = r -
    \frac{\epsilon}{2}$ the ball $B_{r'}(p_n)$ is \emph{not} entirely
    contained in $D_M(B_r(p))$ for all $n$. This allows to find a
    point $q_n \in B_{r'}(p_n) \setminus D_M(B_r(p))$. Since
    $D_M(B_r(p))$ is RCCSV the closure $D_M(B_r(p))^\cl$ is compact
    and thus also $B_r(p)^\cl \subseteq D_M(B_r(p))^\cl$. Since the
    metric $g_t$ and hence the distance function $d_{g_t}$ depend (at
    least) continuous on $t$ we conclude that with the convergence of
    $p_n \longrightarrow p$ and $r' < r$ we have $B_{r'}(p_n)
    \subseteq [-1,1] \times B_r(p)^\cl$ for all $n \geq n_0$.
    \begin{figure}
        \centering
        \input{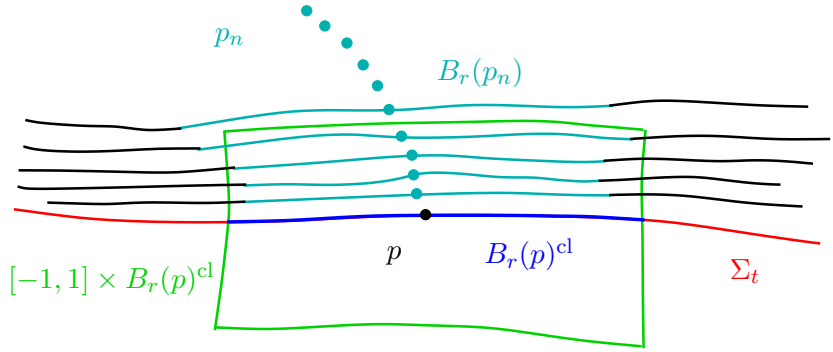}
        \caption{\label{fig:sequence-converging-to-p}%
          Balls around the $p_n$ with radius $r'$ are finally inside
          the box $[-1,1] \times B_r(p)^\cl$.
        }
    \end{figure}
    But then also the points $q_n \in B_{r'}(p_n) \subseteq [-1,1]
    \times B_r(p)^\cl$ are in this compact ``box'', see
    Figure~\ref{fig:sequence-converging-to-p}. Therefore we find a
    convergent subsequence which we denote by $q_n \longrightarrow q$
    as well. Now $p_n \longrightarrow p$ and $q_n \in B_{r'}(p_n)^\cl$
    whence $q \in B_{r'}(p)^\cl$ follows. Since $B_{r'}(p)^\cl
    \subseteq B_r(p)^\cl$ we conclude $q \in B_r(p)$. But
    $D_M(B_r(p))$ is open and hence eventually all sequence elements
    $q_n$ are contained in $D_M(B_r(p))$ which is a
    contradiction. Thus our original claim was in fact true. Thus let
    $p' \in V$ be in this neighborhood and let $r - \epsilon < r' <
    \lambda(p')$. Then by definition we have $B_{r'}(p') \subseteq
    D_M(B_r(p))$ and hence by
    Remark~\ref{remark:properties-of-cauchy-development-operations} we
    have
    \[
    D(B_{r'}(p')) \subseteq D(B_r(p)).
    \]
    Since the larger Cauchy development $D_M(B_r(p))$ is RCCSV this is
    also true for the smaller $D_M(B_{r'}(p'))$. Indeed,
    $D_M(B_{r'}(p'))$ is causal in the surrounding convex $U'$ and has
    smaller volume than $D_M(B_r(p))$ Since $D_M(B_{r'}(p'))^\cl
    \subseteq D_M(B_r(p))^\cl$ it is also pre-compact as wanted. But
    this shows $\rho(p') \geq r' > r - \epsilon$, which is the lower
    semi-continuity.
\end{proof}

\begin{figure}
    \centering
    \input{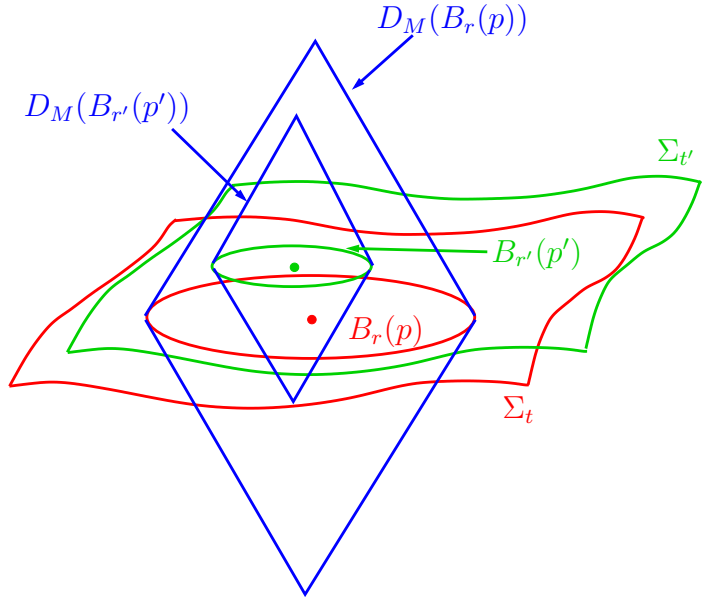}
    \caption{\label{fig:smaller-cauchy-in-bigger-one}%
      For points $p' \in D_M(B_r(p))$ the Cauchy development of a
      smaller ball is included in that of $B_r(p)$.
    }
\end{figure}

\begin{figure}
    \centering
    \input{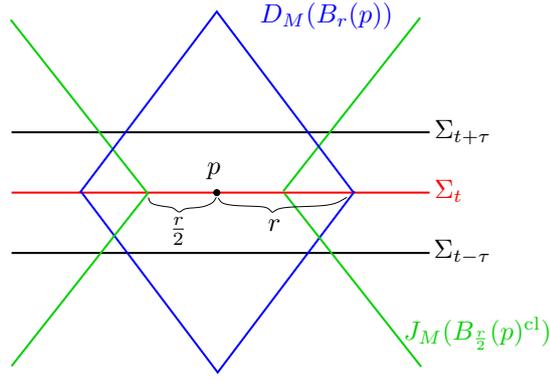}
    \caption{\label{fig:biggest-time-for-box-function}%
      Illustration for the function $\theta_r(p)$.
    }
\end{figure}

Geometrically, this semi-continuity means that for a given $B_r(p)$
around $p$ we can find a ball $B_{r'}(p')$ around $p'$ with only
slightly smaller $r' < r$ such that the Cauchy development of
$B_{r'}(p')$ is still entirely in the one of $B_r(p)$, see also
Figure~\ref{fig:smaller-cauchy-in-bigger-one}.
The next auxiliary function we shall need is the following. We define
for $r>0$ and $p \in M$ (always with respect to the chosen Cauchy
temporal function)
\begin{equation}
    \label{eq:biggest-time-for-box-function}
    \theta_r(p)
    = \sup \left\{
        \tau > 0
        \; \Big| \;
        J_M\left(B_{\frac{r}{2}}(p)^\cl\right) \cap
        \left([t - \tau, t + \tau] \times \Sigma\right)
        \subseteq D_M(B_r(p))
    \right\},
\end{equation}
where $t$ is the time corresponding to the point $p$, i.e. $p \in
\{t\} \times \Sigma \subseteq M$. The picture to have in mind is
sketched in Figure~\ref{fig:biggest-time-for-box-function}.
Again, we first show that this is well-defined, i.e. the subset of
$\tau > 0$ with $J_M(B_{\frac{r}{2}}\cl) \cap ( [t-\tau,t+\tau] \times
\Sigma) \subseteq D_M(B_r(p))$ is non-empty:
\begin{lemma}
    \label{lemma:there-is-a-tau-that-fulfills-condition}
    For every $p \in M$ and $r > 0$ there exists a $\tau > 0$ such
    that
    \begin{equation}
        \label{eq:there-is-a-tau-that-fulfills-condition}
        J_M \left(B_{\frac{r}{2}}(p)^\cl\right)
        \cap
        \left([t-\tau,t+\tau] \times \Sigma\right)
        \subseteq D_M(B_r(p)),
    \end{equation}
    where $t \in \mathbb{R}$ is the unique time with $p \in \Sigma_t$.
\end{lemma}
\begin{proof}
    First we note the following statement: for a compact subset $K
    \subseteq M = \mathbb{R} \times \Sigma$ let $t_{\min}$ and
    $t_{\max}$ be the minimum and maximum of the time function on $K$,
    respectively. Then consider an arbitrary time $t \geq t_{\max}$
    and let $\widetilde{K} = J_M^+(K) \cap \Sigma_t$ which we can
    identify with a subset of $\Sigma$ again since $\Sigma_t \simeq
    \Sigma$.
    \begin{figure}
        \centering
        \input{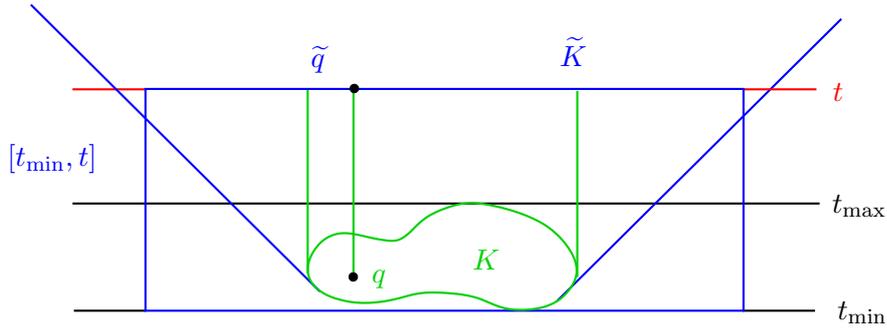}
        \caption{\label{fig:illustration-for-proof-of-existence-of-tau}%
          The compact subset $K$ is in $[t_{\min}, t] \times
          \widetilde{K}$.
        }
    \end{figure}
    Guided by
    Figure~\ref{fig:illustration-for-proof-of-existence-of-tau} we
    claim that $K$ is contained in $[t_{\min}, t] \times
    \widetilde{K}$: indeed, let $p \in K$ be given, then there is a
    timelike curve from $p$ to $q \in \widetilde{K}$ which is just
    $\tau \mapsto (\tau, p)$ where $\tau$ ranges from the time $t(p)
    \geq t_{\min}$ of $p$ to $t$. Thus in the trivialization $p$
    corresponds to $(t(p),q) \in [t_{\min}, t] \times
    \widetilde{K}$. Since $K$ is compact, one knows that $J_M^+(K)
    \cap \Sigma = \widetilde{K}$ is compact as well, see
    e.g. \cite[p.~44]{minguzzi.sanchez:2006a:pre}.  This shows that
    $J_M^+(K) \cap ([t_{\min},t] \times \Sigma) \subseteq [t_{\min},t]
    \times \widetilde{K}$ is compact, too. As we can argue analogously
    for $J_M^-(K)$ we see that for any compact subset $K \subseteq M$
    the subset $J_M(K) \cap ([t_1,t_2] \times \Sigma)$ is compact for
    $t_1 \leq t_{\min}$ and $t_2 \geq t_{\max}$. In particular,
    $J_M(B_{\frac{r}{2}}(p)^\cl) \cap ([t - \frac{1}{2},t +
    \frac{1}{2}] \times \Sigma)$ is compact for all $n \geq 1$ since
    here $t_{\min} = t = t_{\max}$.

    Now assume such a $\tau$ with
    \eqref{eq:there-is-a-tau-that-fulfills-condition} does not
    exist. Then we find $q_n \in ([t - \frac{1}{n}, t + \frac{1}{n}]
    \times \Sigma) \cap J_M(B_{\frac{r}{2}}(p)^\cl)$ which are not in
    $D_M(B_r(p))$. Since the subset $J_M(B_{\frac{r}{2}}(p)^\cl) \cap
    ([t - \frac{1}{n}, t + \frac{1}{n}] \times \Sigma)$ is
    \emph{compact} we can pass to a convergent subsequence, which we
    also denote by $q_n$ converging to some $q$. Clearly, the point
    $q$ has time value $t$. But this means $q \in
    J_M(B_{\frac{r}{2}}(p)^\cl) \cap \{t_0\} \times \Sigma =
    B_{\frac{r}{2}}(p)^\cl$. Now $D_M(B_r(q))$ is an open neighborhood
    of $B_{\frac{r}{2}}(p)^\cl$, thus we have necessarily $q_n \in
    D_M(B_r(q))$ for almost all $n$. This a contradiction and hence we
    have a $\tau > 0$ as wanted.
\end{proof}

\begin{lemma}
    \label{lemma:theta-function-is-lower-semi-continuous}
    The function $\theta_r: M \longrightarrow (0, \infty]$ is
    well-defined and lower semi-continuous.
\end{lemma}
\begin{proof}
    By the last lemma, the function is well-defined.  We consider $p
    \in M$ and $\epsilon > 0$. Then we have to show that for all $p'$
    in a suitable neighborhood of $p$ we still have $\theta_r(p') \geq
    \theta(p) - \epsilon$. Let $t \in \mathbb{R}$ be the time of
    $p$. We assume that there is \emph{no} such open neighborhood of
    $p$. Thus we find a sequence $p_n \longrightarrow p$ of points
    with $\theta_r(p_n) < \theta_r(p) - \epsilon$ for all $n$. Since
    $I_M(B_r(p))$ as well as $(t-\tau,t+\tau) \times \Sigma$ are open
    neighborhoods of $p$ we have
    \[
    p_n \in
    J_M\left(B_{r}(p)^\cl\right)
    \cap
    \left([-T, T] \times \Sigma\right)
    \tag{$*$}
    \]
    for $T$ large enough and $n \geq n_0$. As already argued in the
    proof of Lemma~\ref{lemma:there-is-a-tau-that-fulfills-condition},
    this subset is compact. For the times $t_n$ of $p_n \in \{t_n\}
    \times \Sigma$ we know $t_n \longrightarrow t$ as $p_n
    \longrightarrow p$. Since $\theta_r(p_n) < \theta_r(p) - \epsilon$
    we have
    \[
    J_M\left(B_{\frac{r}{2}}(p)^\cl\right)
    \cap
    \left(
        \left[
            t_n - \theta_r(p) + \epsilon,
            t_n + \theta_r(p) - \epsilon
        \right]
        \times \Sigma
    \right)
    \nsubseteq D_M(B_r(p_n)).
    \]
    Hence we can choose points $q_n \in
    J_M\left(B_{\frac{r}{2}}(p)^\cl\right) \cap \left([t_n -
        \theta_r(p) + \epsilon, t_n + \theta_r(p) - \epsilon] \times
        \Sigma\right)$ which are \emph{not} in $D_M(B_r(p_n))$, see
    Figure~\ref{fig:construction-of-qn-sequence}.
    \begin{figure}
        \centering
        \input{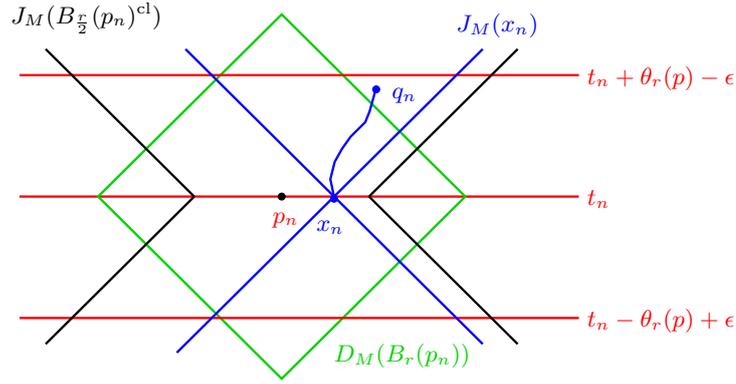}
        \caption{\label{fig:construction-of-qn-sequence}%
          Construction of the points $q_n$ in the proof of
          Lemma~\ref{lemma:theta-function-is-lower-semi-continuous}.
        }
    \end{figure}
    By definition we find $x_n \in B_{\frac{r}{2}}(q_n)^\cl$ with $q_n
    \in J_M(x_n)$. From $p_n \longrightarrow p$ we also conclude that
    for sufficiently large $n$ we have
    \[
    J_M\left(B_{\frac{r}{2}}(p_n)^\cl\right)
    \subseteq
    J_M\left(B_r(p)^\cl\right),
    \]
    see also Figure~\ref{fig:cauch-develop-of-pn}.
    \begin{figure}
        \centering
        \input{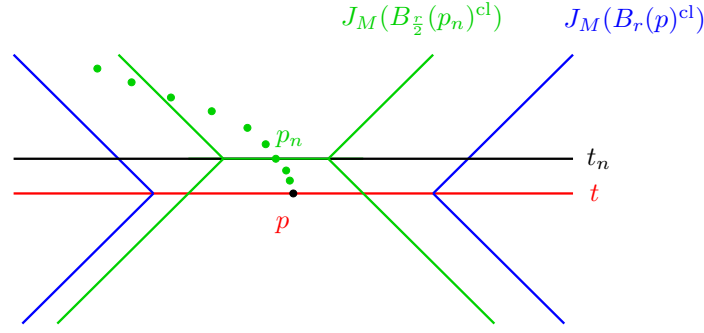}
        \caption{\label{fig:cauch-develop-of-pn}%
          The double cones of the half radius balls around $p_n$ are
          included in the double cone of the full radius ball around
          $p$ for large $n$.
        }
    \end{figure}
    This shows that $q_n \in J_M\left(B_{\frac{r}{2}}(p_n)^\cl\right)
    \subseteq J_M\left(B_r(p)^\cl\right)$ whence together with $q_n
    \in [t_n - \theta_r(p) + \epsilon, t_n + \theta_r(p) - \epsilon]
    \times \Sigma$ we see that all the $q_n$ are in the compact subset
    ($*$). For the $x_n$ this is also true as we have $x_n \in
    B_{\frac{r}{2}}(p_n)^\cl \subseteq
    J_M\left(B_{\frac{r}{2}}(p_n)^\cl\right)$. We may pass to
    convergent subsequences $q_n \longrightarrow q$ and $x_n
    \longrightarrow x$. Since $x_n \in B_{\frac{r}{2}}(p_n)^\cl$ with
    $p_n \longrightarrow p$ we conclude by continuity of the
    Riemannian distance function that $x \in
    B_{\frac{r}{2}}(p)^\cl$. Moreover, since the causal relation
    ``$\leq$'' is closed on a globally hyperbolic spacetime, see
    Remark~\ref{remark:3}, we conclude from $q_n \in J_M(x_n)$ and the
    convergence of the sequences that $q \in J_M(x)$ and hence $q \in
    J_M(B_{\frac{r}{2}}(p)^\cl)$. In addition, since $q_n \in [t_n -
    \theta_r(p) + \epsilon, t_n + \theta_r(p) - \epsilon] \times
    \Sigma$ and $t_n \longrightarrow t$ we conclude that $q \in [t -
    \theta_r(p)+ \epsilon, t + \theta_r(p) - \epsilon] \times
    \Sigma$. Thus we can use the definition of the function
    $\theta_r(p)$ at $p$ and conclude from
    \[
    q \in
    J_M\left(B_{\frac{r}{2}}(p)^\cl\right)
    \cap
    \left(
        \left[
            t - \theta_r(p) + \epsilon,
            t + \theta_r(p) - \epsilon
        \right]
        \times \Sigma
    \right)
    \]
    that $q \in D_M(B_r(p))$. Indeed, $\theta_r(p) > \theta_r(p) -
    \epsilon$ whence we can apply
    \eqref{eq:biggest-time-for-box-function} for $\tau = \theta_r(p) -
    \epsilon$. Since the $q_n$ are not in $D_M(B_r(p_n))$ we have an
    inextensible causal curve $\gamma_n$ through $q_n$ which does not
    meet $B_r(p_n)$, see Figure~\ref{fig:the-curve-gamma-n}. However,
    since $\Sigma_{t_n}$ is a Cauchy hypersurface, it meets $\gamma_n$
    in exactly one point, say $y_n$, see also
    Remark~\ref{remark:cauchy-hypersurface}.
    \begin{figure}
        \centering
        \input{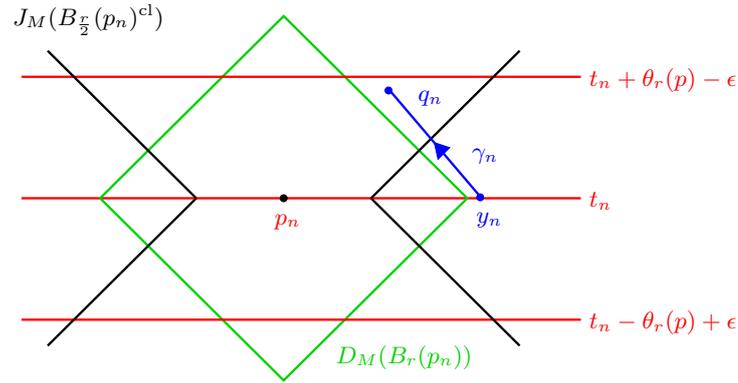}
        \caption{\label{fig:the-curve-gamma-n}%
          The causal curve $\gamma_n$ from the proof of
          Lemma~\ref{lemma:theta-function-is-lower-semi-continuous}
          which does not meet $B_r(p_n)$.
        }
    \end{figure}
    Now we claim that the $y_n$ are also in a compact subset. To this
    end we consider again a large enough $T$ such that all times
    occurring are in $[-T,T]$. First note that it may well happen that
    none of the $y_n$ are in the compact subset ($*$) if $p_n
    \longrightarrow p$ but all the $p_n$ have the same time and come,
    say from the ``right''. In this case, already Minkowski spacetime
    gives us $y_n$ not in ($*$).
    \begin{figure}
        \centering
        \input{last-proof.\pictype}
        \caption{\label{fig:last-proof-figure}%
          The compactum $L$.
        }
    \end{figure}
    However, the intersection $J_M(B_{\frac{r}{2}}(p)^\cl) \cap
    \Sigma_T = L$ is compact and hence the past of $L$ intersected
    with the time interval $[-T,T] \times \Sigma$ is again compact, as
    we argued in the proof of
    Lemma~\ref{lemma:there-is-a-tau-that-fulfills-condition}, see
    Figure~\ref{fig:last-proof-figure}. But now $q_n \in
    J_M\left(B_{\frac{r}{2}}(p)^\cl\right) \subseteq
    J_M\left(B_r(p)^\cl\right)$ shows that $q_n \in J_M^-(L)$. But
    then also the past $J_M^-(q_n)$ is in the past of $J_M^-(L)$ and
    thus $y_n \in J_M^-(q_n) \subseteq J_M^-(L)$. Since the time of
    $y_n$ is $t_n \in [-T,T]$ we conclude that $y_n \in J_M^-(L) \cap
    ([-T,T] \times \Sigma)$ for all $n$. Clearly, if $y_n$ is in the
    future of $q_n$, i.e. the Figure~\ref{fig:last-proof-figure} is
    reversed, the same holds for $J_M^+(L) \cap ([-T,T] \times
    \Sigma)$. Taking the union $J_M(L) \cap ([-T,T] \times \Sigma)$
    will therefore give a compactum for which all $y_n$ are
    inside. Thus we can also here pass to a convergent subsequence
    $y_n \longrightarrow y$. Necessarily $y \in \Sigma_t$ as $t_n
    \longrightarrow t$. Since $y_n \notin B_r(p_n)$ we conclude $y \in
    B_r(p)$ by continuity of the distance function $d_{g_t}$ with
    respect to $t$. Since all the curves $\gamma_n$ are causal, we
    have $q_n \in J_M(y_n)$ and by the closedness of the causal
    relation ``$\leq$'' on a globally hyperbolic spacetime we conclude
    $q \in J_M(y)$. Hence there are inextensible causal curves through
    $y$ and $q$. But since every such curve meets $\Sigma_t$ in only
    one point, namely in $y$, it can not meet $B_r(p)$. However, $q
    \in D_M(B_r(p))$, which is a contradiction.
\end{proof}

The importance of the two lower semi-continuous functions $\rho$ and
$\theta_r$ is that they are bounded from below on every compact
subset: this is an adaption of the statement that a continuous
function is bounded (it takes maximum and minimum) on a compact
subset. Indeed, let $f: K \longrightarrow \mathbb{R}$ be lower
semi-continuous and $K$ compact. Then for all $p \in K$ and $\epsilon
> 0$ we find an open neighborhood $U(p)$ of $p$ such that $f(p') \geq
f(p) - \epsilon$ for all $p' \in U(p)$. Covering $K$ with finitely
many such neighborhoods $U(p_1), \ldots, U(p_n)$ we see that $f(p')
\geq \min_i f(p_i) - \epsilon$ whence $f$ is \emph{bounded} from
below. Let $c = \inf_{p \in K} f(p)$ the infimum of $f$. Then we have
a sequence $p_n \in K$ with $f(p_n) \longrightarrow c$. Now $K$ is
compact whence $p_n$ has a convergent subsequence which we denote also
by $p_n \longrightarrow p$. Thus let $\epsilon > 0$ and choose $U
\subseteq K$ such that $f(p') \geq f(p) - \epsilon$ for all $p' \in
U$. Now all but finitely $p_n$ are in $U$ whence $f(p_n) \geq f(p) -
\epsilon$ for all but finitely many $n$. It follows that also the
limit $\lim_n f(p_n)$ satisfies $c = \lim_n f(p_n) \geq f(p) -
\epsilon$. Thus $c \geq f(p) - \epsilon$ for all $\epsilon > 0$. But
by construction of $c$ we know $c \leq f(p)$ whence $f(p) = c$
follows.

It follows that on a compact subset $K \subseteq M$ the functions
$\rho$ and $\theta_r$ are bounded from zero. We use this in the
following lemma:
\begin{lemma}
    \label{lemma:solution-to-hom-solution-on-timeslice}
    Let $K \subseteq M$ by compact. Then there is a $\delta > 0$ such
    that for all times $t \in \mathbb{R}$ and all $u_t, \dot{u}_t \in
    \Secinfty(\iota_t^\#E)$ on $\Sigma_t$ with support $\supp u_t,
    \supp \dot{u}_t \subseteq K$ we have a smooth solution $u$ of the
    homogeneous wave equation $Du = 0$ on the time slice $(t - \delta,
    t + \delta) \times \Sigma$ with the initial conditions
    $u\at{\Sigma_t} = u_t$ and $\nabla^E_{\mathfrak{n}}u \at{\Sigma_t}
    = \dot{u}_t$. Moreover, for the support one has
    \begin{equation}
        \label{eq:supp-of-timeslice-solution}
        \supp u
        \subseteq J_M \left(\supp u_t \cup \supp \dot{u}_t\right).
    \end{equation}
\end{lemma}
\begin{proof}
    Since $\rho$ is lower semi-continuous according to
    Lemma~\ref{lemma:sup-over-RCCSV-radi} and positive, it admits a
    minimum on the compact subset $K$. Thus we find an $r_0 > 0$ with
    $\rho(p) > 2r_0$ for all $p \in K$. For this radius, the function
    $\theta_{2r_0}$ is lower semi-continuous according to
    Lemma~\ref{lemma:theta-function-is-lower-semi-continuous} and
    positive. Hence we find a $\delta > 0$ with $\theta_{2r_0} >
    \delta$ on $K$. We claim that this $\delta$ will do the job. Thus
    let $t \in \mathbb{R}$ be given. Since $\Sigma_t \cap K$ is again
    compact, we can cover $\Sigma_t \cap K$ with finitely many open
    balls $B_{r_0}(p_1), \ldots, B_{r_0}(p_N)$ of radius $r_0$, where
    as usual the notion of ``ball'' refers to the Riemannian manifold
    $(\Sigma_t, g_t)$. We can find a smooth partition of unity
    $\chi_1, \ldots, \chi_N$ subordinate to the cover $B_{r_0}(p_1)
    \cup \ldots \cup B_{r_0}(p_N)$, i.e. on this open cover of
    $\Sigma_t \cap K$ we have $\chi_1 + \cdots + \chi_N = 1$ and
    $\supp \chi_\alpha \subseteq B_{r_0}(p_\alpha)$ for all $\alpha =
    1, \ldots, N$. It follows that we can decompose the initial
    conditions $u_t$ and $\dot{u}_t$ into smooth pieces having compact
    support in $B_{r_0}(p_\alpha)$ by considering $\chi_\alpha u_t$
    and $\chi_\alpha \dot{u}_t$, respectively. Clearly, we still have
    $\chi_\alpha u_t, \chi_\alpha \dot{u}_t \in \Secinfty_0(\iota_t^\#
    E)$ and $\chi_1 u_t + \cdots + \chi_N u_t = u_t$ as well as
    $\chi_1 \dot{u}_t + \cdots + \chi_N \dot{u}_t = \dot{u}_t$. By
    definition of $\rho$ the Cauchy development
    $D_M(B_{2r_0}(p_\alpha))$ of the balls with twice the radius is
    still RCCSV, see Figure~\ref{fig:covering-of-compact-set}.
    \begin{figure}
        \centering
        \input{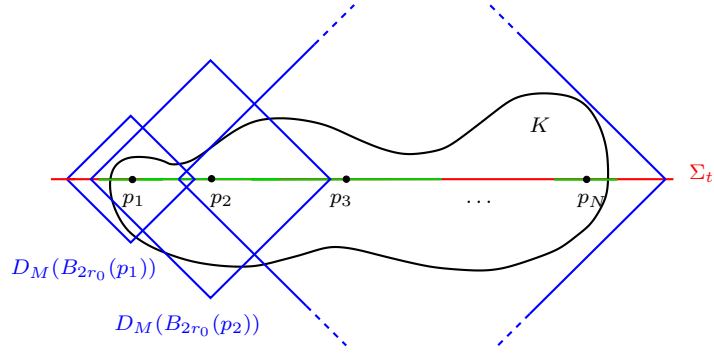}
        \caption{\label{fig:covering-of-compact-set}%
          The covering of the compact subset $K \cap \Sigma_t$ and the
          Cauchy development of the balls.
        }
    \end{figure}
    Thus we can apply
    Proposition~\ref{proposition:existence-of-local-solutions-with-initial-values}
    to these open subsets and obtain a smooth solution $u_\alpha \in
    \Secinfty\left(E \at{D_M(B_{2r_0}(p_\alpha))}\right)$ of the
    homogeneous wave equation $D u_\alpha = 0$ on
    $D_M(B_{2r_0}(p_\alpha))$ for the initial conditions
    \[
    u_\alpha \at{\Sigma_t} = \chi_\alpha u_t
    \quad \textrm{and} \quad
    \nabla^E_{\mathfrak{n}} u_\alpha \at{\Sigma_t}
    = \chi_\alpha \dot{u}_t.
    \]
    Moreover, since we consider the homogeneous wave equation, the
    supports satisfy
    \[
    \supp u_\alpha \subseteq
    J_M \left(
        \supp \chi_\alpha u_t \cup \supp \chi_\alpha \dot{u}_t
    \right).
    \tag{$*$}
    \]
    By definition of the function $\theta_{2r_0}$ and the choice of
    $\delta$ we see that
    \[
    J_M\left(B_{r_0}(p_\alpha)^\cl\right)
    \cap
    \left([t-\delta,t+\delta] \times \Sigma\right)
    \subseteq D_M(B_{2r_0}(p_\alpha)).
    \]
    Hence the solution $u_\alpha$ is defined on the subset
    $J_M\left(B_{r_0}(p_\alpha)^\cl\right) \cap
    \left([t-\delta,t+\delta] \times \Sigma\right)$. Moreover, since
    $\supp \chi_\alpha u_t, \supp \chi_\alpha \dot{u}_t \subseteq
    B_{r_0}(p_\alpha)$ we conclude from ($*$) that
    \[
    \supp u_\alpha \subseteq J_M\left(B_{r_0}(p_\alpha)^\cl\right).
    \]
    \begin{figure}
        \centering
        \input{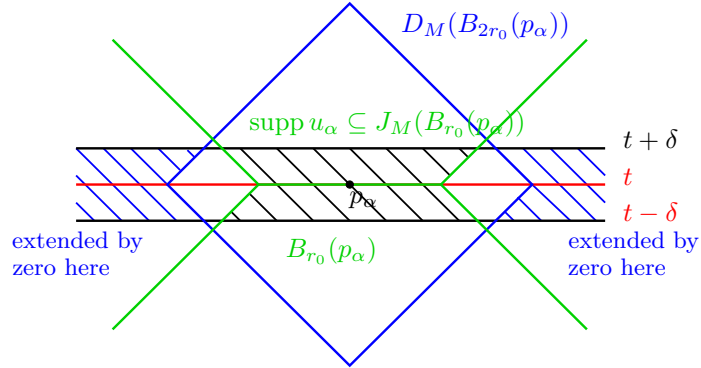}
        \caption{\label{fig:local-solutions-u-alpha}%
          The local solutions $u_\alpha$ constructed in the proof of
          Lemma~\ref{lemma:solution-to-hom-solution-on-timeslice} and
          their support.
        }
    \end{figure}
    Since $u_\alpha$ is smooth on $D_M(B_{2r_0}(p_\alpha))$ we can
    safely extend $u_\alpha$ by zero to $(t-\delta, t+\delta) \times
    \Sigma$, see Figure~\ref{fig:local-solutions-u-alpha}, and have a
    section $u_\alpha \in \Secinfty\left(E\at{(t-\delta,t+\delta)
          \times \Sigma}\right)$ satisfying $\supp u_\alpha \subseteq
    J_M\left(B_{r_0}(p_\alpha)^\cl\right) \cap
    \left([t-\delta,t+\delta] \times \Sigma\right)$ and $D u_\alpha =
    0$ as well as
    \[
    u_\alpha \at{\Sigma_t} = \chi_\alpha u_t
    \quad \textrm{and} \quad
    \nabla^E_{\mathfrak{n}} u_\alpha \at{\Sigma_t}
    = \chi_\alpha \dot{u}_t.
    \]
    Then their sum $u = u_1 + \ldots + u_N$ will still satisfy $Du =
    0$ on $(t-\delta,t+\delta) \times \Sigma$ and
    \[
    u \at{\Sigma_t} = u_t
    \quad \textrm{as well as} \quad
    \nabla^E_{\mathfrak{n}} u \at{\Sigma_t} = \dot{u}_t,
    \]
    since the $\chi_\alpha$ are a partition of unity. Finally,
    \begin{align*}
        \supp u
        &\subseteq
        \supp u_1 \cup \ldots \cup \supp u_N \\
        &\subseteq
        J_M\left(\supp \chi_1 u_t \cup \supp \chi_1 \dot{u}_t\right)
        \cup \ldots \cup
        J_M\left(\supp \chi_N u_t \cup \supp \chi_N \dot{u}_t\right)
        \\
        &\subseteq
        J_M\left(
            \supp \chi_1 u_t \cup \supp \chi_1 \dot{u}_t
            \cup
            \supp \chi_N u_t \cup \supp \chi_N \dot{u}_t
        \right) \\
        &\subseteq
        J_M\left(\supp u_t \cup \supp \dot{u}_t\right),
    \end{align*}
    since on one hand $J_M(A) \cup J_M(B) \subseteq J_M(A \cup B)$ and
    on the other hand $\supp \chi_\alpha u_t \subseteq \supp u_t$ and
    $\supp \chi_\alpha \dot{u}_t \subseteq \supp \dot{u}_t$ for all
    $\alpha$. This completes the proof.
\end{proof}
\begin{remark}
    \label{remark:Ck-timeslice-solutions}
    We see from the proof that we do not loose any differentiability
    by the globalization process. Only for the local solvability of
    the Cauchy problem we need to count orders of differentiation
    carefully. The reason is that the partition of unity can be chosen
    smooth and hence we do not spoil regularity by decomposing
    everything into small pieces. Thus we get from
    Proposition~\ref{proposition:local-Ck-solutions} the analogous
    statement: for initial conditions $u_t \in
    \Sec[2(k+n+1)+2](\iota_t^\# E)$ and $\dot{u}_t \in
    \Sec[2(k+n+2)+1](\iota_t^\# E)$ with the same support conditions
    we get a solution $u \in \Sec\left(E\at{(t-\delta,t+\delta)\times
          \Sigma}\right)$, where of course $k \geq 2$. The statement
    on the support is also still valid.
\end{remark}

Now we come to the existence of global solutions to the Cauchy
problem. As before, $M = \mathbb{R} \times \Sigma$ is globally
hyperbolic with a smooth spacelike Cauchy hypersurface.
\begin{theorem}
    \label{theorem:global-solutions}%
    \index{Cauchy problem!global existence}%
    \index{Wave equation!global solution}%
    \index{Wave equation!inhomogeneous}%
    Let $(M,g)$ be a globally hyperbolic spacetime with smooth
    spacelike Cauchy hypersurface $\iota: \Sigma \hookrightarrow M$.
    \begin{theoremlist}
    \item \label{item:smooth-global-solutions} For $u_0, \dot{u}_0 \in
        \Secinfty_0(\iota^\# E)$ and $v \in \Secinfty_0(E)$ there
        exists a unique global solution $u \in \Secinfty(E)$ of the
        inhomogeneous wave equation $Du = v$ with initial conditions
        $\iota^\# u = u_0$ and $\iota^\# \nabla^E_{\mathfrak{n}} u =
        \dot{u}_0$. We have
        \begin{equation}
            \label{eq:supp-of-smooth-global-solution}
            \supp u
            \subseteq
            J_M \left(
                \supp u_0 \cup \supp \dot{u}_0 \cup \supp v
            \right).
        \end{equation}
    \item \label{item:Ck-global-solution} For $k \geq 2$ and $u_0 \in
        \Sec[2(k+n+1)+2]_0(\iota^\# E)$, $\dot{u}_0 \in
        \Sec[2(k+n+1)+1]_0(\iota^\# E)$ and $v \in
        \Sec[2(k+n+1)]_0(E)$ there exists a unique global solution $u
        \in \Sec(E)$ of the inhomogeneous wave equation $Du = v$ with
        initial conditions $\iota^\# u = u_0$ and $\iota^\#
        \nabla^E_{\mathfrak{n}} u = \dot{u}_0$. It also satisfies
        \eqref{eq:supp-of-smooth-global-solution}.
    \end{theoremlist}
\end{theorem}
\begin{proof}
    Uniqueness follows in both cases from
    Theorem~\ref{theorem:uniqueness-of-hom-equation-with-initial-values}.
    We consider the first case with smooth initial conditions. Since
    all the supports are compact so is their union. Therefore, we can
    cover this compact subset with finitely many RCCSV subsets for
    which we can apply the local existence according to
    Proposition~\ref{proposition:existence-of-local-solutions-with-initial-values}.
    Again, choosing an appropriate partition of unity subordinate to
    this cover, we can decompose the initial conditions and the
    inhomogeneity into pieces having their compact supports inside of
    the RCCSV subsets. If we succeed to show the existence of a global
    solution for such initial conditions and inhomogeneity with
    support in the RCCSV subset, we can afterwards sum up this finite
    number of solutions to get a solution for the arbitrary $u_0,
    \dot{u}_0$ and $v$. This shows that without restriction, we can
    assume that $\supp u_0, \supp \dot{u}_0$ and $\supp v$ are
    contained in a single RCCSV subset $U \subseteq U^\cl \subseteq
    U'$ as required by
    Proposition~\ref{proposition:existence-of-local-solutions-with-initial-values}.
    We set $K = \supp u_0 \cup \supp \dot{u}_0 \cup \supp v \subseteq
    U$, which is still compact. By using a second partition of unity
    argument, we can cut $K$ into even smaller pieces such that we
    have, for the fixed $U$, the properties
    \[
    K \subseteq (-\epsilon, \epsilon) \times \Sigma
    \]
    and
    \[
    J_M(K) \cap ( (-\epsilon,\epsilon) \times \Sigma ) \subseteq U,
    \]
    for an appropriate small $\epsilon > 0$, see
    Figure~\ref{fig:supp-of-K}.
    \begin{figure}
        \centering
        \input{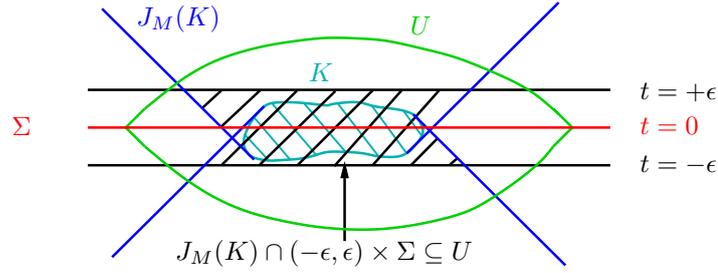}
        \caption{\label{fig:supp-of-K}%
          The compact set $K = \supp u_0 \cup \supp \dot{u}_0 \cup
          \supp v$ of the proof of
          Theorem~\ref{theorem:global-solutions}.
        }
    \end{figure}
    Now let $u \in \Secinfty(E\at{U})$ be the solution according to
    Proposition~\ref{proposition:existence-of-local-solutions-with-initial-values}.
    Since $\supp u \subseteq J_M(K)$ we see that we can smoothly
    extend $u$ to the whole time slice $(-\epsilon,\epsilon) \times
    \Sigma$ by $0$. We have to argue that we can extend this solution
    even further on arbitrarily large time slices $(-T,T) \times
    \Sigma$. Thus we set $T_{\max}$ to be the supremum of all those
    times $T$ for which there exists a smooth extension of $u$ to the
    slice $(-\epsilon, T) \times \Sigma$, still obeying the causality
    condition $\supp u \subseteq J_M(K)$. Since we have at least $T
    \geq \epsilon$ the supremum $T_{\max}$ is positive. Since $K$ is
    in the slice $(-\epsilon,\epsilon) \times \Sigma$ we have $Du=0$
    on $[+\epsilon, T_{\max}) \times \Sigma$ since the inhomogeneity
    has $\supp v \subseteq K$. If we have two extensions, $u$ until
    $T$ and $\widetilde{u}$ until $\widetilde{T}$ with $T <
    \widetilde{T}$, then $\widetilde{u} \at{(-\epsilon, T) \times
      \Sigma} = u$ since the open piece $(-\epsilon, T) \times \Sigma$
    is globally hyperbolic itself. Hence the uniqueness statement from
    Theorem~\ref{theorem:uniqueness-of-hom-equation-with-initial-values}
    applies to $\widetilde{u}\at{(-\epsilon,T) \times \Sigma}$ and
    $u$. Thanks to this uniqueness we only have to show the existence
    of a solution for arbitrary, but fixed finite $T$, i.e. $T_{\max}
    = + \infty$. This will automatically give a solution defined for
    all times $t \in \mathbb{R}^+$ and hence a solution on
    $(-\epsilon, \infty) \times \Sigma$.

    We assume the converse, i.e. $T_{\max} < \infty$. We consider
    $\widetilde{K} = ( [-\epsilon, T_{\max}] \times \Sigma ) \cup
    J_M(K)$ which is compact as we have already argued at the
    beginning of the proof of
    Lemma~\ref{lemma:there-is-a-tau-that-fulfills-condition} in
    greater generality. We can therefore apply
    Lemma~\ref{lemma:solution-to-hom-solution-on-timeslice} to this
    compact subset $\widetilde{K}$ yielding a $\delta > 0$ as
    described there.
    \begin{figure}
        \centering
        \input{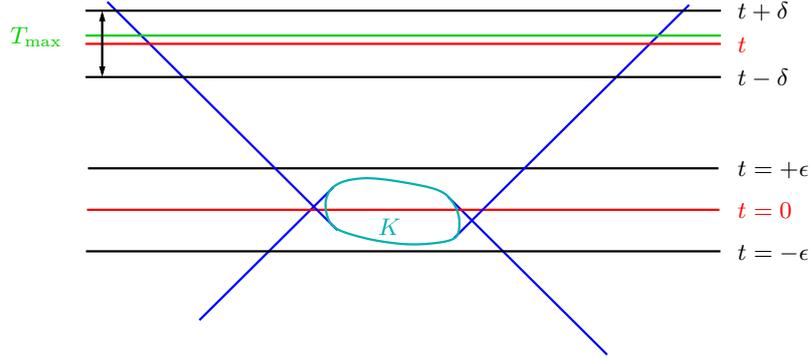}
        \caption{\label{fig:choosing-delta}%
          The choice of $t$ for given $\delta > 0$ in the proof of
          Theorem~\ref{theorem:global-solutions}.}
    \end{figure}
    Now we take a $t < T_{\max}$ with $T_{\max} - t < \delta$ but $K
    \subseteq (-\epsilon,t) \times \Sigma$. Note that since $K
    \subseteq (-\epsilon,\epsilon) \times \Sigma$ and $T_{\max} \geq
    \epsilon$, this is clearly possible no matter what $\delta > 0$
    is, see Figure~\ref{fig:choosing-delta}. On the whole slice
    $(t-\delta,t+\delta) \times \Sigma$ we solve the homogeneous wave
    equation $Dw=0$ for the initial conditions
    \[
    w \at{\Sigma_t} = u \at{\Sigma_t}
    \quad \textrm{and} \quad
    \nabla^E_{\mathfrak{n}} w \at{\Sigma_t}
    = \nabla^E_{\mathfrak{n}} u \at{\Sigma_t},
    \]
    which is possible thanks to
    Lemma~\ref{lemma:solution-to-hom-solution-on-timeslice}. On a
    smaller slice $(t - \eta,t + \eta) \times \Sigma$ the
    inhomogeneity $v$ already vanishes by $\supp v \subseteq K$ since
    $K$ is contained in the \emph{open} slice $(-\epsilon,t) \times
    \Sigma$. Thus on this slice $w$ and $u$ both solve the homogeneous
    wave equation with the same initial conditions on
    $\Sigma_t$. Therefore $w=u$ on $(-\epsilon,t) \times \Sigma$,
    again by the uniqueness theorem. But this shows that $w$ extends
    $u$ to the slice $(-\epsilon, t+\delta) \times \Sigma$ in a smooth
    way. For the support we see that the initial conditions for $w$
    are contained in $J_M(K) \cap \Sigma_t$.
    \begin{figure}
        \centering
        \input{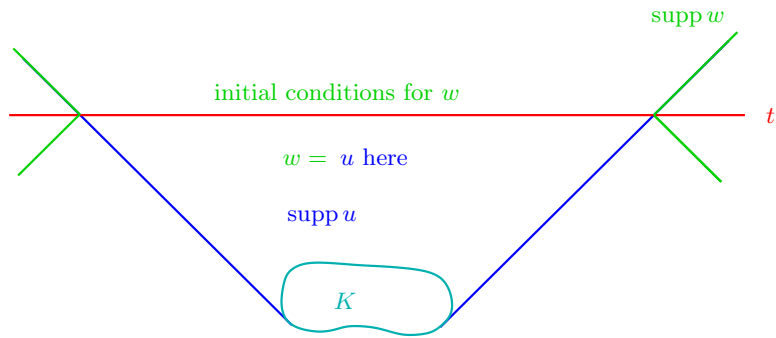}
        \caption{\label{fig:the-extension-w}%
          The extension $w$ of $u$ in the proof of
          Theorem~\ref{theorem:global-solutions}.
        }
    \end{figure}
    For the future of $t$ this means that $\supp w$ is still contained
    in $J_M(K)$, for the past of $t$ we already know that $w = u$
    whence in total $\supp w \subseteq J_M(K)$, see
    Figure~\ref{fig:the-extension-w}. But $T_{\max} < t + \delta$
    whence we get a contradiction since $w$ is a valid extension of
    $u$ with all desired properties. Thus $T_{\max} = + \infty$. An
    analogous argument shows that also in the past directions we can
    extend the solution to $t = - \infty$. This gives the first
    part. The second part proceeds completely analogous, using only
    Proposition~\ref{proposition:local-Ck-solutions} and
    Remark~\ref{remark:Ck-timeslice-solutions} instead.
\end{proof}

%
%

\subsection{Well-Posedness of the Cauchy Problem}
\label{satz:well-posedness}

We have seen that the Cauchy problem for the inhomogeneous wave
equation with smooth initial data and smooth compactly supported
inhomogeneity admits a unique smooth solution. Also in the context of
sufficiently large but finite differentiability we have a unique
solution to the Cauchy problem. A Cauchy problem is called
\emph{well-posed}\index{Cauchy problem!well-posed} if for given
initial data one has a unique solution which depends
\emph{continuously} on the initial data. Of course, this requires to
specify the relevant topologies in detail. In typical situations, the
relevant topologies should be clear from the context. Note also that
for physical applications a continuous dependence on the initial data
is certainly necessary in order to have a physically reasonable
theory: initial data are always subject to (arbitrarily small but
non-zero) uncertainties when measured. Thus a discontinuous dependence
would lead to a physical theory without \index{Predictive power}
predictive power. But even if one has continuous dependence it may
well happen for Cauchy problems that the discrepancy at finite times
between solutions corresponding to very close initial conditions grows
very fast in time, typically in an exponential way when quantified
correctly. Thus it might be of interest to have the continuity even
sharpened by some more quantitative description.

Back to our situation we want to show the well-posedness of the Cauchy
problem with respect to the usual locally convex topologies of smooth
or $\Fun$-sections. The main tool will be the following general
statement from locally convex analysis:
\begin{theorem}[Open mapping theorem]
    \label{theorem:open-mapping}%
    \index{Open mapping theorem}%
    Let $\mathcal{E}, \widetilde{\mathcal{E}}$ be Fr\'echet spaces and
    let $\phi: \mathcal{E} \longrightarrow \widetilde{\mathcal{E}}$ be
    a continuous linear map. If $\phi$ is surjective then $\phi$ is an
    open map.
\end{theorem}
As usual, a map $\phi$ is called \emph{open} if the images of open
subsets are again open. The proof of the open mapping theorem can
e.g. be found in \cite[Thm.~2.11]{rudin:1991a}. We will need the
following corollary of it:
\begin{corollary}
    \label{corollary:continous-linear-bijection-is-homeomorphism}
    Let $\phi: \mathcal{E} \longrightarrow \widetilde{\mathcal{E}}$ be
    a continuous linear bijection between Fr\'echet spaces. Then
    $\phi^{-1}$ is continuous as well.
\end{corollary}
Indeed, let $U \subseteq \mathcal{E}$ be open. Then the set-theoretic
$(\phi^{-1})^{-1}(U)$, i.e. the pre-image of $U$ under $\phi^{-1}$,
coincides simply with $\phi(U)$ which is open by the theorem. Thus
$\phi^{-1}$ is continuous. Note that for general maps between
topological spaces a continuous bijective map needs not have a
continuous inverse at all.

We are now interested in the following situation: the result of
Theorem~\ref{theorem:global-solutions} can be viewed as a map
\begin{equation}
    \label{eq:the-solution-machine}
    \Secinfty_0(\iota^\# E)
    \oplus \Secinfty_0(\iota^\# E)
    \oplus \Secinfty_0(E)
    \longrightarrow \Secinfty(E),
\end{equation}
sending $(u_0, \dot{u}_0, v)$ to the unique solution $u$ of the wave
equation $Du = v$ with initial conditions $u_0$ and
$\dot{u}_0$. Clearly, the map \eqref{eq:the-solution-machine} is
linear which easily follows from the uniqueness statement of
Theorem~\ref{theorem:global-solutions}. Thus continuous dependence on
the initial conditions will refer to the continuity of the map
\eqref{eq:the-solution-machine}. Note that this even includes the
continuous dependence on the inhomogeneity $v$. The relevant
topologies are then the $\Cinfty$-topology on the target side and the
canonical topology of the direct sum of the
$\Cinfty_0$-topologies. Since the direct sum is finite, this is not
problematic and essentially boils down to show $\Cinfty_0$-continuity
for each summand. This way, we arrive at the following theorem:
\begin{theorem}[Well-posed Cauchy problem I]
    \label{theorem:well-posedness-1}%
    \index{Cauchy problem!well-posed}%
    Let $(M,g)$ be a globally hyperbolic spacetime with smooth
    spacelike Cauchy hypersurface $\iota: \Sigma \hookrightarrow
    M$. Then the linear map \eqref{eq:the-solution-machine} sending
    the initial conditions and the inhomogeneity to the corresponding
    solution of the Cauchy problem is continuous.
\end{theorem}
\begin{proof}
    First we note that the ``inverse'' map which evaluates an
    arbitrary section $u \in \Secinfty(E)$ on the Cauchy hypersurface
    and applies $D$ to it is continuous, i.e.
    \[
    \mathcal{P}: \Secinfty(E) \ni u \mapsto
    (\iota^\# u, \iota^\# \nabla^E_{\mathfrak{n}}u, Du)
    \in
    \Secinfty(\iota^\# E)
    \oplus \Secinfty(\iota^\# E)
    \oplus \Secinfty(E)
    \tag{$*$}
    \]
    is continuous in the $\Cinfty$-topologies. This is clear as all
    three components of $\mathcal{P}$ are continuous. Indeed, the
    restriction is continuous by a slight variation of the results
    from
    Proposition~\ref{proposition:pullback-of-sections-is-continous}. The
    application of either $\nabla^E_{\mathfrak{n}}$ or $D$ is
    continuous as well whence the continuity of each of the three
    components of $\mathcal{P}$ follows. However, for a general $u \in
    \Secinfty(E)$ neither the restrictions $\iota^\# u$ and $\iota^\#
    \nabla^E_{\mathfrak{n}}u$ nor $Du$ will have compact support. Thus
    we enforce this by considering a fixed compact subset $K \subseteq
    M$ and the subspaces $\Secinfty_{K \cap \Sigma}(\iota^\# E)$ as
    well as $\Secinfty_K(E)$ of $\Secinfty(\iota^\# E)$ and
    $\Secinfty(E)$ of those sections with compact support in the
    compact subsets $K \cap \Sigma$ and $K$, respectively. By
    Lemma~\ref{lemma:sections-with-compact-support-in-fixed-compacta}
    we know that both spaces are Fr\'echet spaces as they are
    $\Cinfty$-closed subspaces of the Fr\'echet spaces
    $\Secinfty(\iota^\#E)$ and $\Secinfty(E)$, respectively. Hence
    their direct sum is a closed subspace of the target in ($*$)
    whence the pre-image
    \[
    \mathcal{V}_K
    = \mathcal{P}^{-1}
    ( \Secinfty_{K \cap \Sigma}(\iota^\# E) \oplus
    \Secinfty_{K \cap \Sigma}(\iota^\# E) \oplus \Secinfty_K(E) )
    \subseteq \Secinfty(E)
    \]
    is again closed. This way, it becomes a Fr\'echet subspace
    itself. Restricted to $\mathcal{V}_K$, the map $\mathcal{P}_K =
    \mathcal{P} \at{\mathcal{V}_K}$ becomes bijective, this is
    precisely the statement of
    Theorem~\ref{theorem:global-solutions}. Indeed, $\mathcal{P}_K$ is
    surjective since every point in $\Secinfty_{K \cap
      \Sigma}(\iota^\#E) \oplus \Secinfty_{K \cap \Sigma}(\iota^\# E)
    \oplus \Secinfty_K(E)$ has a pre-image. This is just the existence
    of the solutions to the Cauchy problem. However, as the solution
    is unique, we have precisely one pre-image under
    $\mathcal{P}_K$. Since now all involved spaces are Fr\'echet
    themselves and $\mathcal{P}_K$ is obviously continuous, we can
    apply
    Corollary~\ref{corollary:continous-linear-bijection-is-homeomorphism}
    to conclude that $\mathcal{P}_K$ has continuous inverse
    \[
    \mathcal{P}_K^{-1}:
    \Secinfty_{K \cap \Sigma}(\iota^\# E) \oplus
    \Secinfty_{K \cap \Sigma}(\iota^\# E) \oplus \Secinfty(E)
    \longrightarrow \mathcal{V}_K \subseteq \Secinfty(E)
    \]
    for all $K \subseteq M$ compact. By the definition of the
    inductive limit topology this gives us immediately the continuity
    of the map \eqref{eq:the-solution-machine} as claimed. In fact,
    this is again a general feature of LF topologies and this trick
    can be transferred to the general situation, see
    e.g. \cite{jarchow:1981a}.
\end{proof}

With an analogous argument we also obtain the well-posedness of the
Cauchy problem in the following situation of finite differentiability:
\begin{theorem}[Well-posed Cauchy problem II]
    \label{theorem:well-posed-Cauchy-Ck}%
    \index{Cauchy problem!well-posed}%
    Let $(M,g)$ be a globally hyperbolic spacetime with smooth
    spacelike Cauchy hypersurface $\iota: \Sigma \hookrightarrow M$
    and let $k \geq 2$. Then the linear map
    \begin{equation}
        \label{eq:Ck-initial-value-machine-map}
        \Sec[2(k+n+1)+2]_0(\iota^\# E)
        \oplus \Sec[2(k+n+1)+1]_0(\iota^\# E)
        \oplus \Sec[2(k+n+1)]_0(E)
        \longrightarrow \Sec(E)
    \end{equation}
    sending $(u_0, \dot{u}_0, v)$ to the unique solution $u$ of the
    inhomogeneous wave equation $Du = v$ with initial conditions
    $\iota^\# u = u_0$ and $\iota^\# \nabla^E_{\mathfrak{n}} u =
    \dot{u}_0$ is continuous.
\end{theorem}

Thus we have in both cases a well-posed Cauchy problem. There are,
however, some small drawbacks of the above theorems: First, as already
mentioned, we are limited to inhomogeneities $v$ with compact support
in $M$. Physically more appealing would be an inhomogeneity with
compact support only in spacelike direction, i.e. the ``eternally
moving electron''. Note that this is clearly an intrinsic concept on a
globally hyperbolic spacetime. Moreover, the control of derivatives in
Theorem~\ref{theorem:global-solutions} and hence in
Theorem~\ref{theorem:well-posed-Cauchy-Ck} seems not to be optimal. In
particular, it would be nice to show that the map
\eqref{eq:Ck-initial-value-machine-map} has some fixed order
\emph{independent} of $k$.


%% file: green.tex
%
%

While in Chapter~\ref{cha:LocalTheory} we have discussed the local
existence of fundamental solutions as well as their properties we
shall now pass to the global picture. From the uniqueness statements
in Corollary~\ref{corollary:uniqueness-glob-hyp} we see that the local
advanced and retarded fundamental solutions necessarily agree with the
restrictions of the corresponding global ones if the latter exist at
all. Here we have to restrict to such an RCCSV neighborhood which is
globally hyperbolic itself, i.e. a Cauchy development of a small
enough ball in $\Sigma$. Then the question of existence of global
fundamental solutions can be viewed as the question whether the given
local fundamental solutions can be extended to the whole spacetime.

Actually, we shall proceed differently and construct the global
fundamental solutions directly using the global statements on the
Cauchy problem. As before, we assume throughout this section that
$(M,g)$ is globally hyperbolic.

%
%

\subsection{Global Green Functions}
\label{satz:glob-green-functions}

We first consider the smooth version. Here we start with the following
theorem:
\begin{theorem}
    \label{theorem:global-fund-solution}%
    \index{Green function!global}%
    \index{Fundamental solution!global}%
    Let $(M,g)$ be a globally hyperbolic spacetime and $D \in
    \Diffop^2(E)$ a normally hyperbolic differential operator. For
    every point $p \in M$ there is a unique advanced and retarded
    fundamental solution $F^\pm_M(p)$ of $D$ at $p$. Moreover, for
    every test section $\varphi \in \Secinfty_0(E^*)$ the section
    \begin{equation}
        \label{eq:fund-solution-gives-smooth-inhom-solution}
        M \ni p \; \mapsto \; F^\pm_M(p)\varphi \in E^*_p
    \end{equation}
    is a smooth section of $E^*$ which satisfies the equation
    \begin{equation}
        \label{eq:solution-of-tran-inhom-equation}
        D^\Trans F^\pm_M(\argument)\varphi = \varphi.
    \end{equation}
    Finally, the linear map
    \begin{equation}
        \label{eq:solution-depends-continous-on-inhom}
        F^\pm_M: \Secinfty_0(E^*) \ni \varphi
        \; \mapsto \;
        F^\pm_M(\argument)\varphi \in \Secinfty(E^*)
    \end{equation}
    is continuous.
\end{theorem}
\begin{proof}
    The uniqueness was already shown in
    Corollary~\ref{corollary:uniqueness-glob-hyp}. For the existence
    we consider the following construction: we first choose a
    splitting $M \simeq \mathbb{R} \times \Sigma$ with a Cauchy
    temporal function being the first coordinate of the product and
    $\Sigma$ being a smooth spacelike Cauchy hypersurface. We denote
    as usual by $\Sigma_t$ the level set of fixed time $t$,
    i.e. $\Sigma_t = \{t\} \times \Sigma
    \stackrel{\iota_t}{\hookrightarrow} M$, which is again a Cauchy
    hypersurface. Normalizing the gradient of $t$ appropriately we
    obtain the smooth future-directed unit normal vector field
    $\mathfrak{n} \in \Secinfty(TM)$ which, at $\Sigma_t$, is normal
    to $\Sigma_t$ for all times $t$. Now let $\varphi \in
    \Secinfty_0(E^*)$ be a test section of $E^*$. Since $\varphi$ has
    compact support we find a time $t$ such that $\supp \varphi$ is in
    the past of $t$. More precisely, we have $\supp \varphi \subseteq
    I^-_M(\Sigma_{t})$, see
    Figure~\ref{fig:supp-of-inhom-is-in-past-of-hypersurface}.
    \begin{figure}
        \centering
        \input{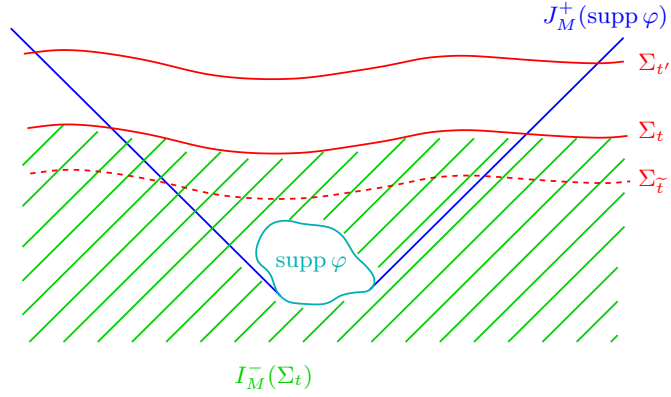}
        \caption{\label{fig:supp-of-inhom-is-in-past-of-hypersurface}%
          The various hypersurfaces chosen in the future of $\supp
          \varphi$ in the proof of
          Theorem~\ref{theorem:global-fund-solution}.
        }
    \end{figure}
    We now apply Theorem~\ref{theorem:global-solutions},
    \refitem{item:smooth-global-solutions} to the transposed operator
    $D^\Trans \in \Diffop^2(E^*)$ which we know to be normally
    hyperbolic as well. Thus we obtain a unique global and smooth
    solution $\psi^+ \in \Secinfty(E^*)$ of the inhomogeneous wave
    equation $D^\Trans \psi^+ = \varphi$ for the initial conditions
    $\iota_t^\# \psi^+ = 0 = \iota_t^\# \nabla^{E^*}_{\mathfrak{n}}
    \psi^+$. First we note that $\psi^+$ does not depend on the
    precise choice of $t$. Indeed, let $t'$ be another time with
    $\supp \varphi \subseteq I^-_M(\Sigma_{t'})$ and assume e.g. $t <
    t'$. Denote by $\psi^+{}' \in \Secinfty(E^*)$ the corresponding
    solution of the Cauchy problem $D^\Trans \psi^+{}' = \varphi$ and
    $\iota_{t'}^\# \psi^+{}' = 0 = \iota_{t'}^\#
    \nabla^{E'}_{\mathfrak{n}} \psi^+{}'$. Then we find a
    $\widetilde{t} < t$ such that $\supp \varphi \subseteq
    I^-_M(\Sigma_{\widetilde{t}})$ since $I_M^-(\Sigma_t)$ is
    \emph{open} while $\supp \varphi$ is closed, see again
    Figure~\ref{fig:supp-of-inhom-is-in-past-of-hypersurface}. The
    open piece $(\widetilde{t}, \infty) \times \Sigma = \bigcup_{t >
      \widetilde{t}} \Sigma_t \subseteq M$ is still a globally
    hyperbolic spacetime on its own. Here $\psi^+{}'$ satisfies
    $D^\Trans \psi^+{}' = 0$ since $\supp \varphi$ is not in this part
    of $M$. Since $\psi^+{}'$ has vanishing initial conditions on
    $\Sigma_{t'}$ we conclude by the uniqueness properties of
    solutions that $\psi^+{}' = 0$ on $(\widetilde{t}, \infty) \times
    \Sigma$. This implies in particular the feature that $\iota_t^\#
    \psi^+{}' = 0 = \iota_t^\# \nabla^{E^*}_{\mathfrak{n}} \psi^+{}'$
    whence both $\psi^+$ and $\psi^+{}'$ have vanishing initial
    conditions on $\Sigma_t$ and satisfy the wave equation $D^\Trans
    \psi^+{}' = \varphi = D^\Trans \psi^+$ on all of $M$. Thus by the
    uniqueness according to
    Theorem~\ref{theorem:uniqueness-initial-condition-C2-case} we
    conclude $\psi^+ = \psi^+{}'$. Hence the section $\psi^+$ does not
    depend on the choice of $t$ as long as $t$ is large
    enough. According to Theorem~\ref{theorem:well-posedness-1} the
    map which assigns $\varphi$ to $\psi^+$ is a continuous linear map
    with respect to the $\Cinfty_0$- and $\Cinfty$-topology. Moreover,
    evaluating $\psi^+$ at a given point $p \in M$ is a $E^*_p$-valued
    continuous linear functional, namely the
    $\delta_p$-functional. Thus the map $\varphi \mapsto \psi^+(p)$ is
    a continuous linear functional for every point $p \in M$. This
    defines the generalized section $F^+_M(p) \in \Sec[-\infty](E)
    \tensor E^*_p$, i.e.
    \[
    F^+_M(p): \varphi \; \mapsto \; \psi^+(p)
    \]
    with $\psi^+$ as above. By definition of $F^+_M(p)$ the map
    \eqref{eq:fund-solution-gives-smooth-inhom-solution} is just the
    map $\varphi \mapsto \psi^+$ which is continuous according to
    Theorem~\ref{theorem:well-posedness-1} and yields a smooth section
    $\psi^+$. This shows
    \eqref{eq:fund-solution-gives-smooth-inhom-solution} and
    \eqref{eq:solution-depends-continous-on-inhom} for the case of
    ``$+$''. We prove that $F^+_M(p)$ is a fundamental solution at
    $p$. For the two test sections $\varphi, D^\Trans \varphi \in
    \Secinfty_0(E^*)$ we have resulting solutions $\psi^+,
    \widetilde{\psi}^+$ as above, i.e. $D^\Trans \psi^+ = \varphi$ and
    $D^\Trans \widetilde{\psi}^+ = D^\Trans \varphi$. Thus $D^\Trans
    (\widetilde{\psi}^+ - \varphi) = 0$ and both $\widetilde{\psi}^+$
    and $\varphi$ have vanishing initial conditions on $\Sigma_t$: the
    section $\varphi$ even vanishes in an open neighborhood of
    $\Sigma_t$ while $\widetilde{\psi}^+$ has vanishing initial
    conditions on $\Sigma_t$ by construction. Thus by uniqueness we
    have $\widetilde{\psi}^+ - \varphi = 0$. Unwinding this gives
    \[
    \left(D F^+_M(p)\right)(\varphi)
    = F^+_M(p) (D^\Trans \varphi)
    = \widetilde{\psi}^+(p)
    = \varphi(p),
    \]
    hence $D F^+_M(p) = \delta_p$ follows as $\varphi$ is an arbitrary
    test section. This gives us a fundamental solution $F^+_M(p)$ at
    every point $p \in M$. It remains to show that $F^+_M(p)$ is
    actually an advanced fundamental solution, i.e. $\supp F^+_M(p)
    \subseteq J^+_M(p)$. Since $J^+_M(p)$ is closed by global
    hyperbolicity of $M$ we have to find an open neighborhood of $q
    \in M \setminus J^+_M(p)$ on which $F^+_M(p)$ vanishes. Thus let
    $q \notin J^+_M(p)$ be such a point. By closedness of $J^+_M(p)$
    there is an open neighborhood of $q$ such that $q' \notin
    J^+_M(p)$ for all $q'$ in this neighborhood. We distinguish two
    cases. If $p \notin J^+_M(q)$ then we also have $p \notin
    J^+_M(q')$ for all $q'$ in a small neighborhood of $q$. Thus we
    can choose $q', q''$ close to $q$ with $q \in I^+_M(q') \cap
    I^+_M(q'')$ but $p \notin J^+_M(q'')$ and $p \notin J^-_M(q')$,
    see Figure~\ref{fig:the-three-q-points}.
    \begin{figure}
        \centering
        \input{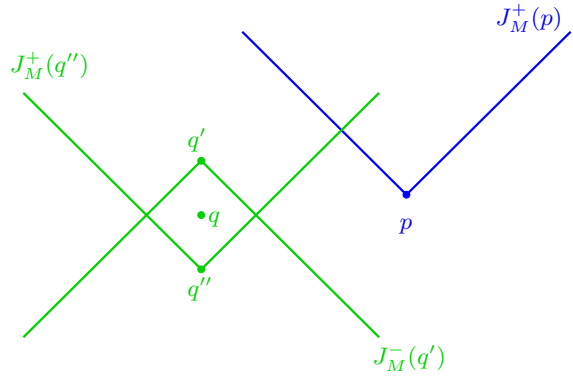}
        \caption{\label{fig:the-three-q-points}%
          Choosing the points $q'$ and $q''$ with $q \in J_M(q'',q')$
          for $q$ and $p$ spacelike.
        }
    \end{figure}
    In this case $p \notin J_M(J_M(q'',q'))$. Since $I_M(q'',q')$ is
    an open neighborhood of $q$ we have for all $\varphi \in
    \Secinfty_0(E^*)$ with $\supp \varphi \subseteq I_M(q'',q')$ by
    Theorem~\ref{theorem:global-solutions},
    \refitem{item:smooth-global-solutions} the property $\supp \psi^+
    \subseteq J_M(\supp \varphi) \subseteq J^+_M(q'') \cup J^-_M(q')$,
    where $\psi^+$ is the section with vanishing initial conditions
    for large times and $D^\Trans \psi^+ = \varphi$. Since $p \notin
    J_M^+(q'') \cup J^-_M(q')$ we have $0 = \psi^+(p) = F^+_M(p)
    \varphi$. However, this simple argument only works for $p$ and $q$
    spacelike. Thus the other case is where $p$ and $q$ are not
    spacelike, but $p$ is in $J^+_M(q)$.
    \begin{figure}
        \centering
        \input{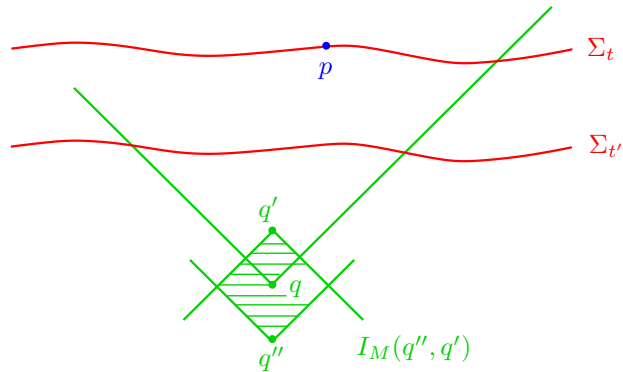}
        \caption{\label{fig:p-q-not-spacelike}%
          The points $p$ and $q$ do not lie spacelike to each other.
        }
    \end{figure}
    But then necessarily the time $t$ of $p$ is strictly larger than
    the one of $q$ as $q \neq p$. We fix a time $t'$ between $t$ and
    the time of $q$ and choose a point $q'$ on $\Sigma_t$ in the
    future $I^+_M(q)$ of $q$, see
    Figure~\ref{fig:p-q-not-spacelike}. Moreover, let $q'' \in
    I^-_M(q)$ be arbitrary. This gives us an open diamond
    $I_M(q'',q')$ which is an open neighborhood of $q$. Let $\varphi
    \in \Secinfty_0(E^*)$ have support in $I_M(q'',q')$ and let
    $\psi^+$ be the solution of $D^\Trans \psi^+ = \varphi$ with
    vanishing initial values for large times as before. Since $t > t'$
    is clearly later than $\supp \varphi$ we have $\psi^+
    \at{\Sigma_t} = 0$. But this gives $\psi^+(p) = 0$ also in this
    case and hence $F^+_M(p) \varphi = 0$ for all such $\varphi$. This
    finally shows that $\supp F^+_M(p) \subseteq J^+_M(p)$ as
    wanted. The retarded case is analogous as usual.
\end{proof}

We can strengthen the above result in the following way. As we have at
least some rough counting of needed derivatives in
Theorem~\ref{theorem:global-solutions},
\refitem{item:Ck-global-solution} for the Cauchy problem we can use
this to estimate the order of the Green functions $F^\pm_M(p)$:
\begin{theorem}
    \label{theorem:order-of-global-green-functions}%
    \index{Green function!global order}%
    Let $(M,g)$ be a globally hyperbolic spacetime and $D \in
    \Diffop^2(E)$ a normally hyperbolic differential operator. Then
    the unique advanced and retarded Green functions $F^\pm_M(p)$ of
    $D$ at $p$ are of global order
    \begin{equation}
        \label{eq:order-of-global-green-functions}
        \ord F^\pm_M(p) \leq 2n + 6.
    \end{equation}
    More precisely, the linear map
    \eqref{eq:solution-depends-continous-on-inhom} extends to a
    continuous linear map
    \begin{equation}
        \label{eq:GlobalGreen-function-on-Ck}
        F^\pm_M: \Sec[2(k+n+1)]_0(E^*) \ni \varphi
        \; \mapsto \;
        F^\pm_M(\argument) \varphi \in \Sec(E^*)
    \end{equation}
    for all $k \geq 2$ such that we still have
    \begin{equation}
        \label{eq:green-function-on-Ck-solves-wave-equation}
        D^\Trans F^\pm_M(\argument)\varphi = \varphi.
    \end{equation}
\end{theorem}
\begin{proof}
    By Theorem~\ref{theorem:global-solutions}
    \refitem{item:Ck-global-solution} we can repeat the whole
    construction in the proof of
    Theorem~\ref{theorem:global-fund-solution} for a test section
    $\varphi \in \Sec[2(k+n+1)]_0(E^*)$. Indeed, the initial
    conditions for $\psi^+$ being zero for large times clearly satisfy
    the differentiability conditions of
    Theorem~\ref{theorem:global-solutions},
    \refitem{item:Ck-global-solution}. Thus we obtain a solution
    $\psi^+ \in \Sec(E)$ of $D^\Trans \psi^+ = \varphi$. With the
    definition $F^+_M(p)\varphi = \psi^+(p)$ and hence $\psi^+ =
    F^+(\argument)\varphi$ we get by
    Theorem~\ref{theorem:well-posed-Cauchy-Ck} the continuity of
    \eqref{eq:solution-depends-continous-on-inhom}. By construction,
    \eqref{eq:green-function-on-Ck-solves-wave-equation} still
    holds. Now let $p \in M$ be given and choose $k=2$ which is the
    minimal one allowed by Theorem~\ref{theorem:global-solutions} and
    Theorem~\ref{theorem:well-posed-Cauchy-Ck}. Then the continuity of
    \eqref{eq:solution-depends-continous-on-inhom} implies that for
    all compact $K \subseteq M$ we find a constant $c > 0$ with
    \[
    \left|
        F^+_M(p) \varphi
    \right|
    = \seminorm[\{p\},0] \left(F^+_M(\argument) \varphi\right)
    \leq c \seminorm[K, 2(k+n+1)](\varphi).
    \]
    But this shows that the local order of $F^+_M(\argument)$ on the
    compactum $K$ is less or equal than $2n+6$, independently on
    $K$. It is clear by the usual density argument that the map
    $F^+_M(p)$ defined here is indeed the unique extension of the
    advanced Green function defined in the previous Theorem. The
    retarded case is analogous.
\end{proof}

\begin{remark}
    \label{remark:order-of-glob-green}
    Again, the estimate on the order is usually very rough and even
    worse than the estimate we found in the local case. Nevertheless,
    the important point is that the order is \emph{globally finite}
    and independent of $p$. Since in the construction of the solution
    $\psi^+$ we only needed the very special initial conditions
    $\iota^\# \psi^+ = 0 = \iota^\# \nabla^{E}_{\mathfrak{n}} \psi^+$
    the proof of the local solution to the Cauchy problem as in
    Proposition~\ref{proposition:local-Ck-solutions} with finite
    differentiability simplifies drastically yielding a simplified
    recursion only involving the inhomogeneity. We leave it as an open
    task to improve the estimate
    \eqref{eq:order-of-global-green-functions} on the global order.
\end{remark}

%
%

\subsection{Green Operators}
\label{satz:green-operators}

The fundamental solutions $F^\pm_M(p)$ were constructed as the map
$\varphi \mapsto (p \mapsto F^\pm_M(p)\varphi)$ being a map
$\Secinfty_0(E^*) \longrightarrow \Secinfty(E)$, i.e. the solution map
from the Cauchy problem. We shall now investigate this map more
closely as it provides almost an inverse to $D$. In general, one
defines the following operators.
\begin{definition}[Green Operators]
    \label{definition:green-operators}%
    \index{Green operator!advanced}%
    \index{Green operator!retarded}%
    Let $(M,g)$ be a time-oriented Lorentz manifold and $D \in
    \Diffop^2(E)$ a normally hyperbolic differential operator. Then a
    continuous linear map
    \begin{equation}
        \label{eq:green-operator}
        G^\pm_U: \Secinfty_0(E) \longrightarrow \Secinfty(E)
    \end{equation}
    with
    \begin{definitionlist}
    \item \label{item:DG-is-id} $D G^\pm_M = \id_{\Secinfty_0(E)}$,
    \item \label{item:GD-is-id} $G^\pm_M D \at{\Secinfty_0(E)} =
        \id_{\Secinfty_0(E)}$,
  \item \label{item:supp-of-green-op} $\supp (G^\pm_M u) \subseteq
      J^\pm_M( \supp u)^\cl$ for all $u \in \Secinfty_0(E)$
    \end{definitionlist}
    is called an advanced and retarded Green operator for $D$,
    respectively.
\end{definition}
Note that if the causal relation is not closed we have to put a
closure in part \refitem{item:supp-of-green-op} by hand. In view of
the local result in \eqref{eq:fundi-operators} one can imagine that a
Green operator for $D$ is linked to the fundamental solutions
$G^\pm_M(p)$ of the dual differential operator $D^\Trans \in
\Diffop^2(E^*)$. In fact, we have the following proposition for
general spacetimes, where we require $\supp G^\pm_M(p) \subseteq
(J^\pm_M(p))^\cl$ in the case when the causal relation is not closed.
\begin{proposition}[Green operators and fundamental solutions]
    \label{proposition:green-op-and-fundi-solution}
    Let $(M,g)$ be a time-oriented Lorentz manifold and $D \in
    \Diffop^2(E)$ a normally hyperbolic differential operator.
    \begin{propositionlist}
    \item \label{item:green-op-from-fundi} Assume $\{G^\pm_M(p)\}$ is
        a family of global advanced or retarded fundamental solutions
        of $D^\Trans$ at every point $p \in M$ with the following
        property: for every test section $u \in \Secinfty_0(E)$ the
        section $p \mapsto G^\pm_M(p) u$ is a smooth section of $E$
        depending continuously on $u$ and satisfying $D
        G^\pm_M(\argument)u = u$. Then
        \begin{equation}
            \label{eq:green-op-from-fundi}
            (G^\pm_M u)(p) = G^\mp_M(p) u
        \end{equation}
        yield advanced and retarded Green operators for $D$,
        respectively.
    \item \label{item:fundi-from-green-op} Assume $G^\pm_M$ are
        advanced or retarded Green operator for $D$,
        respectively. Then $G^\pm_M(p): \Secinfty_0(E) \longrightarrow
        \mathbb{C}$ defined by
        \begin{equation}
            \label{eq:fundi-from-green-op}
            G^\pm_M(p) u = (G^\mp_M u)(p)
        \end{equation}
        defines a family of advanced and retarded fundamental
        solutions of $D^\Trans$ at every point $p \in M$ with the
        properties described in \refitem{item:green-op-from-fundi},
        respectively.
    \end{propositionlist}
\end{proposition}
\begin{proof}
    For the first part we assume to have a family $\{G^\pm_M(p)\}_{p
      \in M}$ of advanced or retarded fundamental solutions of
    $D^\Trans$ with the above properties. By assumption, the resulting
    linear map \eqref{eq:green-op-from-fundi} is continuous. It
    satisfies $D G^\pm_M = \id_{\Secinfty_0(E)}$ also by
    assumption. Since the $G^\pm_M(p)$ are fundamental solutions of
    $D^\Trans$ we have
    \[
    (G^\pm_M D u)(p)
    = G^\mp_M(p) (Du)
    = (D^\Trans G^\mp_M(p))(u)
    = \delta_p (u)
    = u(p)
    \]
    for all $p \in M$ and $u \in \Secinfty_0(E)$. Thus $G^\pm_M D =
    \id_{\Secinfty_0(E)}$ as well. Finally, we have to check the
    support properties thereby explaining the flip from $\pm$ to $\mp$
    in \eqref{eq:green-op-from-fundi}. Thus let $p \in M$ be given
    such that $0 \neq (G^\pm_Mu)(p) = G^\mp(p)u$. Since the support of
    the distributions $G^\mp_M(p)$ is in $J^\mp_M(p)^\cl$ this implies
    that $\supp u$ has to intersect $J^\mp_M(p)^\cl$. Since
    $J^\mp_M(p)^\cl = I^\mp_M(p)^\cl$, see
    \cite[Prop.~2.17]{minguzzi.sanchez:2006a:pre}, and since $\supp u$
    has an open interior which is non-empty, we see that $\supp u$
    also has to intersect $I^\mp_M(p)$. But then $p \in I^\mp_M(\supp
    u)$ whence $\supp (G^\pm_M u) \subseteq I^\pm_M (\supp u)^\cl =
    J^\pm_M(\supp u)^\cl$ follows, proving the first part. For the
    second part assume $G^\pm_M$ is given and define $G^\pm_M(p) =
    \delta_p \circ G^\mp_M$, according to
    \eqref{eq:fundi-from-green-op}. This is clearly a distribution
    since $\delta_p$ is continuous and $G^\mp_M$ is continuous by
    assumption. By construction, the section $p \mapsto G^\pm_M(p) u =
    (G^\mp_M u)(p)$ is smooth and depends continuously on $u$. We have
    \[
    D G^\mp_M(\argument)u
    = D \left(p \mapsto G^\mp_M(p) u\right)
    = D G^\pm_M u
    = u
    \]
    as well as
    \[
    \left(D^\Trans G^\mp_M(p)\right)(u)
    = G^\mp_M(p) (Du)
    = \left(G^\pm_M (Du)\right)(p)
    = u(p),
    \]
    whence $G^\mp_M(p)$ is a fundamental solution satisfying also $D
    G^\pm_M(\argument)u = u$. Finally, for the support we can argue as
    before in part \refitem{item:green-op-from-fundi}.
\end{proof}

\begin{remark}[Green operators]
    \label{remark:green-opps}
    ~
    \begin{remarklist}
    \item \label{item:green-op-causal-closed} If the causal relation
        ``$\le$'' is closed then the definition of a Green operator
        simplifies and also the above proof simplifies. This will be
        the case for globally hyperbolic spacetimes.
    \item \label{item:green-op-inverse-of-D} At first glance, a Green
        operator of $D$ looks like an inverse on the space of
        compactly supported sections. However, this is not quite
        correct as $G^\pm_M$ maps into $\Secinfty(E)$ and not into
        $\Secinfty_0(E)$. Nevertheless, the Green operator behaves
        very much like an inverse of $D\at{\Secinfty_0(E)}$.
    \item \label{item:green-ops-might-not-exist} In general, Green
        operators do not exist: if e.g. $M$ is a compact Lorentz
        manifold and $D = \dAlembert$ is the scalar d'Alembertian then
        the constant function $1$ has compact support but satisfied
        $\dAlembert 1 = 0$. Thus $G \dAlembert 1 = 1$ is impossible
        for a linear map $G$.
    \end{remarklist}
\end{remark}
In the case of a globally hyperbolic spacetime our construction of
advanced and retarded fundamental solutions in
Theorem~\ref{theorem:global-fund-solution} gives immediately advanced
and retarded Green operators:
\begin{corollary}
    \label{corollary:green-ops-on-globally-hyperbolic}%
    \index{Green operator!existence}%
    \index{Green operator!uniqueness}%
    On a globally hyperbolic spacetime any normally hyperbolic
    differential operator has unique advanced and retarded Green
    operators.
\end{corollary}
\begin{proof}
    Indeed, the fundamental solutions were precisely constructed as in
    the proposition with the operator coming from the solvability of
    the Cauchy problem in Theorem~\ref{theorem:global-fund-solution}.
\end{proof}

Having related the Green operators of $D$ to the fundamental solutions
of $D^\Trans$ we can also relate the Green operators of $D$ and
$D^\Trans$ directly. First we notice that, as we already did locally
in Section~\ref{satz:SolvingWaveEqLocally}, the Green operators allow
for dualizing:
\begin{proposition}
    \label{proposition:dualizing-of-green-operators}
    Let $(M,g)$ be globally hyperbolic and let $D \in \Diffop^2(E)$ be
    a normally hyperbolic differential operator with advanced and
    retarded Green operators $G^\pm_M: \Secinfty_0(E) \longrightarrow
    \Secinfty(E)$.
    \begin{propositionlist}
    \item \label{item:dualizing-of-green-operators-1} The dual map
        $(G^\pm_M)': \Sec[-\infty]_0(E^*) \longrightarrow
        \Sec[-\infty](E^*)$ is weak$^*$ continuous and satisfies
      \begin{equation}
          \label{eq:dualized-green-op-and-D-trans}
          D^\Trans (G^\pm_M)' (\varphi) = \varphi
          = (G^\pm_M)' D^\Trans \varphi
      \end{equation}
      for all generalized sections $\varphi \in \Sec[-\infty]_0(E^*)$
      with compact support.
  \item \label{item:supp-property-of-dualized-green-op} For a
      generalized section $\varphi \in \Sec[-\infty]_0(E^*)$ with
      compact support we have
      \begin{equation}
          \label{eq:supp-property-of-dualized-green-op}
          \supp (G^\pm_M)' (\varphi)
          \subseteq J^\mp_M( \supp \varphi).
      \end{equation}
    \end{propositionlist}
\end{proposition}
\begin{proof}
    Since $G^\pm_M: \Secinfty_0(E) \longrightarrow \Secinfty(E)$ is
    linear and continuous we have an induced dual map $(G^\pm_M)':
    \Secinfty(E)' = \Sec[-\infty]_0(E^*) \longrightarrow
    \Secinfty_0(E)' = \Sec[-\infty](E^*)$ where we identify the dual
    spaces as usual by means of the canonical volume density
    $\mu_g$. Then $(G^\pm_M)'$ is automatically weak$^*$
    continuous. To prove \eqref{eq:dualized-green-op-and-D-trans} we
    take a test section $u \in \Secinfty_0(E)$ and compute
    \[
    \left(D^\Trans (G^\pm_M)'(\varphi)\right)(u)
    = (G^\pm_M)'(\varphi)(Du)
    = \varphi\left(G^\pm_M Du\right)
    = \varphi(u)
    \]
    by the very definitions. Since $\Secinfty_0(E) \subseteq
    \Secinfty(E)$ is dense this is sufficient to show the first part
    of \eqref{eq:dualized-green-op-and-D-trans}, which is understood
    as an identity between generalized sections with compact
    support. For the other part we compute
    \[
    \left((G^\pm_M)' D^\Trans \varphi\right)(u)
    = (D^\Trans \varphi) (G^\pm_M u)
    = \varphi\left(D G^\pm_M u\right)
    = \varphi(u).
    \]
    Note that $D^\Trans \varphi$ has again compact support whence the
    above computation is indeed justified. This proves
    \eqref{eq:dualized-green-op-and-D-trans}. For the second statement
    let $u \in \Secinfty_0(E)$ be a test section. Then
    $(G^\pm_M)'(\varphi) u = \varphi(G^\pm_M u)$. Since $\supp
    (G^\pm_M u) \subseteq J^\pm_M( \supp u)$ we see that $\varphi
    (G^\pm_M u)$ vanishes if $\supp \varphi \cap J^\pm_M (\supp u) =
    \emptyset$. But this means $J^\mp_M(\supp \varphi) \cap \supp u =
    \emptyset$. Thus for $\supp u \subseteq M \setminus J^\mp_M( \supp
    \varphi)$ we have $(G^\pm_M)'(\varphi) u = 0$ which implies
    \eqref{eq:supp-property-of-dualized-green-op}, since $J^\pm_M(
    \supp \varphi)$ is already closed.
\end{proof}

As in the local situation we can now apply $(G^\pm_M)'$ to generalized
sections $\varphi$ which are actually smooth, i.e. $\varphi \in
\Secinfty_0(E^*)$. We expect that we obtain the Green operators of
$D^\Trans$. Here we need the following simple Lemma:
\begin{lemma}
    \label{lemma:green-ops-of-transpose-and-green-ops-of-D}
    Let $(M,g)$ be a globally hyperbolic spacetime and let $D \in
    \Diffop^2(E)$ be a normally hyperbolic differential operator with
    advanced and retarded Green operators $G^\pm_M$. Moreover, denote
    the corresponding Green operators of $D^\Trans \in \Diffop^2(E^*)$
    by $F^\pm_M$. Then we have for $\varphi \in \Secinfty_0(E^*)$ and
    $u \in \Secinfty_0(E)$
    \begin{equation}
        \label{eq:green-ops-of-transpose-and-green-ops-of-D}
        \int_M (F^\pm_M \varphi) \cdot u \: \mu_g
        = \int_M \varphi \cdot (G^\mp_M u) \: \mu_g.
    \end{equation}
\end{lemma}
\begin{proof}
    The lemma is a simple integrations by parts argument. First we
    note that $F^\pm_M \varphi$ has (non-compact) support in $J^\pm_M(
    \supp \varphi)$ while $G^\mp_M u$ has (non-compact) support in
    $J^\mp_M( \supp u)$ by the very definition of Green operators. It
    follows from the global hyperbolicity that the overlap $J^\pm_M(
    \supp \varphi) \cap J^\mp_M(\supp \varphi)$ is compact, see
    Figure~\ref{fig:supp-overlap-compact}.
    \begin{figure}
        \centering
        \input{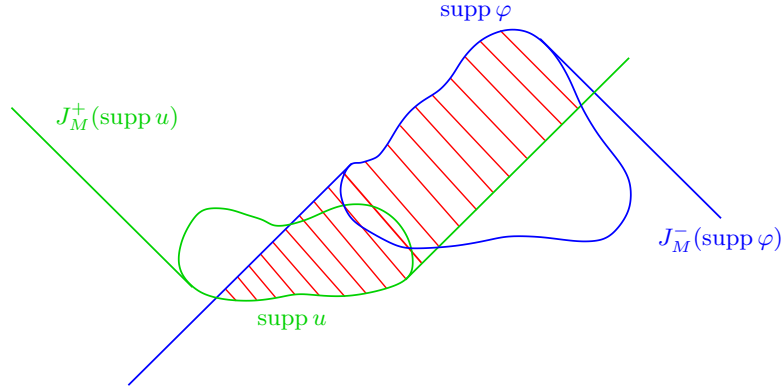}
        \caption{\label{fig:supp-overlap-compact}%
          The compact overlap of $J^+_M( \supp u)$ and $J^-_M(\supp
          \varphi)$.
        }
    \end{figure}
    Thus writing $u = D G^\mp_M u$ we get
    \begin{align*}
        \int_M (F^\pm \varphi) \cdot u \: \mu_g
        &=
        \int_M (F^\pm_M \varphi) \cdot (D G^\mp_M u) \: \mu_g \\
        & \stackrel{(*)}{=}
        \int_M
        (D^\Trans F^\pm_M \varphi) \cdot (G^\mp_M u) \: \mu_g \\
        &=
        \int_M \varphi \cdot (G^\mp_M u) \: \mu_g,
    \end{align*}
    where we have used $D^\Trans F^\pm_M \varphi = \varphi$ and the
    compactness of the overlap to justify the integration by parts in
    ($*$).
\end{proof}

From this lemma we immediately see that the dual operator $(G^\mp_M)':
\Sec[-\infty]_0(E^*) \longrightarrow \Sec[-\infty](E^*)$ applied to a
distributional section which is actually smooth, i.e. to $\varphi \in
\Secinfty_0(E^*) \subseteq \Sec[-\infty]_0(E^*)$ is given by
\begin{equation}
    \label{eq:dual-green-op-on-smooth}
    (G^\mp_M)' \varphi = F^\pm_M \varphi.
\end{equation}
Indeed, this is just the content of
\eqref{eq:green-ops-of-transpose-and-green-ops-of-D} where we
interpret the right hand side as the distributional section $\varphi
\in \Sec[-\infty]_0(E^*)$ evaluated on $G^\mp_M(u)$ as usual. In
particular, the dual map $(G^\mp_M)'$ yields a smooth section and not
just a distributional one when applied to $\varphi \in
\Secinfty_0(E^*)$. Moreover, since $F^\pm_M$ is continuous with
respect to the $\Cinfty_0$- and $\Cinfty$-topology according to
Theorem~\ref{theorem:global-fund-solution} we have also continuity of
the dual operators $(G^\mp_M)'$ on $\Secinfty_0(E^*)$ with respect to
the $\Cinfty_0$- and $\Cinfty$-topology. This way, we obtain the
global analogues of the local results obtained in
Section~\ref{satz:construction-local-solution}. We summarize the
discussion the in the following theorem:
\begin{theorem}
    \label{theorem:dual-of-green-operators}%
    \index{Green operator!dual}%
    Let $(M,g)$ be globally hyperbolic and let $D \in \Diffop^2(E)$ be
    a normally hyperbolic differential operator. Denote the global
    advanced and retarded Green operators of $D$ by $G^\pm_M$ and
    those of $D^\Trans$ by $F^\pm_M$, respectively.
    \begin{theoremlist}
    \item \label{item:dual-of-greenop-and-greenop-of-dual} For the
        dual operators we have
        \begin{equation}
            \label{eq:dual-of-greenop-on-smooth}
            (G^\pm_M)' \at{\Secinfty_0(E^*)}
            = F^\mp_M
        \end{equation}
        \begin{equation}
            \label{eq:dual-of-transposed-greenop-onsmooth}
            (F^\pm_M)' \at{\Secinfty_0(E)}
            = G^\mp_M.
        \end{equation}
    \item \label{item:dual-of-greenops-are-cont} The duals of the
        Green operators restrict to maps
        \begin{equation}
            \label{eq:dual-of-greenop-is-cont}
            (G^\pm_M)': \Secinfty_0(E^*) \longrightarrow \Secinfty(E^*),
        \end{equation}
        \begin{equation}
            \label{eq:dual-of-trans-greenop-is-cont}
            (F^\pm_M)': \Secinfty_0(E) \longrightarrow \Secinfty(E),
        \end{equation}
        which are continuous with respect to the $\Cinfty_0$- and
        $\Cinfty$-topologies, respectively.
    \item \label{item:greenop-on-distributions} The Green operators
        have unique weak$^*$ continuous extensions to operators
        \begin{equation}
            \label{eq:greenop-extension}
            G^\pm_M: \Sec[-\infty]_0(E) \longrightarrow \Sec[-\infty](E)
        \end{equation}
        \begin{equation}
            \label{eq:trans-greenop-extension}
            F^\pm_M: \Sec[-\infty]_0(E^*) \longrightarrow
            \Sec[-\infty](E^*)
        \end{equation}
        satisfying
        \begin{equation}
            \label{eq:supp-property-of-greenop}
            \supp (G^\pm_M u) \subseteq J^\pm_M( \supp u)
        \end{equation}
        \begin{equation}
            \label{eq:supp-property-of-transpose-greenop}
             \supp (F^\pm_M \varphi) \subseteq J^\pm_M (\supp \varphi),
        \end{equation}
        respectively. For these extensions one has
        \begin{equation}
            \label{eq:greenop-is-dual-of-trans-greenop}
            G^\pm_M
            = \left(F^\mp_M \at{\Secinfty_0(E^*)}\right)'
        \end{equation}
        \begin{equation}
            \label{eq:trans-green-op-is-dual-of-greenop}
            F^\pm_M
            = \left(G^\mp_M \at{\Secinfty_0(E)}\right)'.
        \end{equation}
    \end{theoremlist}
\end{theorem}
\begin{proof}
    Indeed, part \refitem{item:dual-of-greenop-and-greenop-of-dual}
    was already discussed and part
    \refitem{item:dual-of-greenops-are-cont} is clear by part
    \refitem{item:dual-of-greenop-and-greenop-of-dual} and the
    continuity of Green operators. The last part is also clear since
    the corresponding dual operators provide us with an extension of
    the Green operators according to
    \refitem{item:dual-of-greenop-and-greenop-of-dual}. The uniqueness
    of the extension is clear as the smooth sections with compact
    support are (sequentially) dense in the distributional sections
    with compact support: this follows analogously to the density
    statement in Theorem~\ref{theorem:gensec-with-compact-support},
    \refitem{item:sec-dense-in-gensec} for the case of arbitrary
    distributional sections. Then \eqref{eq:supp-property-of-greenop}
    and \eqref{eq:supp-property-of-transpose-greenop} are obtained
    from Proposition~\ref{proposition:dualizing-of-green-operators},
    \refitem{item:supp-property-of-dualized-green-op} applied to
    $D^\Trans$ and $D$, respectively. Finally
    \eqref{eq:greenop-is-dual-of-trans-greenop} and
    \eqref{eq:trans-green-op-is-dual-of-greenop} are clear.
\end{proof}

\begin{remark}
    \label{remark:dual-of-greenop-is-trans-greenop}
    With some slight abuse of notation we do not distinguish between
    the Green operators and their canonical extension to generalized
    sections. This gives the short hand version
    \begin{equation}
        \label{eq:dual-of-greenop-and-trans-greenop}
        G^\pm_M = \left(F^\mp_M\right)'
    \end{equation}
    of \eqref{eq:greenop-is-dual-of-trans-greenop} and
    \eqref{eq:trans-green-op-is-dual-of-greenop}. In particular, the
    Green operators of $D^\Trans$ are completely determined by those
    of $D$ and vice versa.
\end{remark}

As a first application of the extended Green operators we obtain a
solution of the wave equation for arbitrary compactly supported
inhomogeneity with good causal behaviour:
\begin{theorem}
    \label{theorem:solution-of-wave-eq-for-distributional-inhom}%
    \index{Cauchy problem!distributional inhomogeneity}%
    Let $(M,g)$ be a globally hyperbolic spacetime and $D \in
    \Diffop^2(E)$ normally hyperbolic with advanced and retarded Green
    operators $G^\pm_M$.
    \begin{theoremlist}
    \item \label{item:extension-of-greenop-are-fundi-for-weak-operator}
        The Green operators $G^\pm_M: \Sec[-\infty]_0(E)
        \longrightarrow \Sec[-\infty](E)$ satisfy
        \begin{equation}
            \label{eq:extension-of-greenop-are-fundi-for-weak-operator}
            D G^\pm_M
            = \id_{\Sec[-\infty]_0(E)}
            = G^\pm_M D \at{\Sec[-\infty]_0(E)}.
        \end{equation}
    \item \label{item:solution-of-wave-eq-for-distributional-inhom}
        For every $v \in \Sec[-\infty]_0(E)$, every smooth spacelike
        Cauchy hypersurface $\iota: \Sigma \hookrightarrow M$ with
        \begin{equation}
            \label{eq:inhomogeneity-supp-condition-of}
            \supp v \subseteq I^+_M(\Sigma),
        \end{equation}
        and all $u_0, \dot{u}_0 \in \Secinfty_0(\iota^\#E)$ there
        exists a unique generalized section $u_+ \in \Sec[-\infty](E)$
        with
        \begin{equation}
            \label{eq:solution-for-distributional-inhom}
            D u_+ = v,
        \end{equation}
        \begin{equation}
            \label{eq:supp-of-distributional-inhom-solution}
            \supp u_+ \subseteq J_M( \supp u_0 \cup \supp \dot{u}_0)
            \cup J_M^+(\supp v),
        \end{equation}
        \begin{equation}
            \label{eq:singsupp-of-distributional-inhom-solution}
            \singsupp u_+ \subseteq J_M^+ (\supp v),
        \end{equation}
        \begin{equation}
            \label{eq:initial-conditions-of-distri-inhom-solution}
            \iota^\# u_+ = u_0
            \quad \textrm{and} \quad
            \iota^\# \nabla^E_{\mathfrak{n}} u = \dot{u}_0.
        \end{equation}
        The section $u_+$ depends weak$^*$ continuously on $v$ and
        continuously on $u_0, \dot{u}_0$.
    \item \label{item:solution-of-wave-eq-for-distributional-inhom-2}
        An analogous statement holds for the case $\supp v \subseteq
        I^-_M(\Sigma)$.
    \end{theoremlist}
\end{theorem}
\begin{proof}
    For the first part we can use the fact that all involved maps are
    weak$^*$ continuous and $\Secinfty_0(E)$ is weak$^*$ dense in
    $\Sec[-\infty]_0(E)$. Then
    \eqref{eq:extension-of-greenop-are-fundi-for-weak-operator} is
    just a consequence of the defining properties of a Green operator
    on $\Secinfty_0(E)$. For the second part we first notice that
    $G^+_M v \in \Sec[-\infty](E)$ is a generalized section with
    support in $J^+_M( \supp v)$ according to
    \eqref{eq:supp-property-of-greenop} and $D G^+_M v = v$ according
    to the first part. Let $w \in \Sec[+\infty](E)$ be the unique
    solution to the Cauchy problem $Dw = 0$ and $\iota^\#w = u_0$ and
    $\iota^\# \nabla^E_{\mathfrak{n}} w = \dot{u}_0$ whose existence
    and uniqueness is guaranteed by
    Theorem~\ref{theorem:global-solutions},
    \refitem{item:smooth-global-solutions}. We set $u = w + G^+_M
    v$. This is a generalized section with $Du = v$ as $w$ solves the
    homogeneous wave equation. Moreover, we have
    \begin{align*}
        \supp u
        &= \supp(w+G^+_Mv) \\
        &\subseteq
        \supp w \cup \supp G^+_M v \\
        &\subseteq
        J_M( \supp u_0 \cup \supp \dot{u}_0) \cup J^+_M(\supp v),
    \end{align*}
    according to \eqref{eq:supp-of-smooth-global-solution} and
    \eqref{eq:supp-property-of-greenop}. Since $w$ is smooth we also
    have $\singsupp u = \singsupp G^+_M v \subseteq J^+_M( \supp
    v)$. Now $\supp v \subseteq J^+_M(\Sigma)$ implies that $M
    \setminus J^+_M(\supp v)$ is an open neighborhood of $\Sigma$, see
    Figure~\ref{fig:supps-of-generalized-solution}.
    \begin{figure}
        \centering
        \input{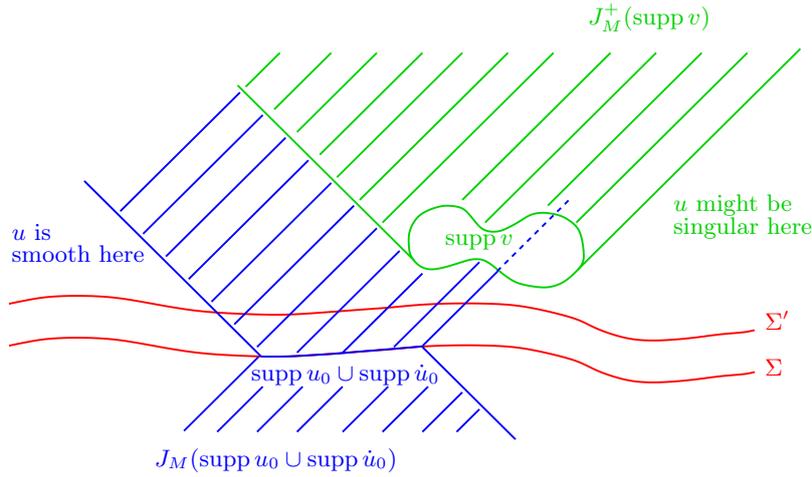}
        \caption{\label{fig:supps-of-generalized-solution}%
          The figure shows the supports of the inhomogeneity $v$, the
          initial conditions $u_0, \dot{u}_0$ and where the solution
          to the inhomogeneous wave equation might be singular.}
    \end{figure}
    Thus $u$ is smooth on an open neighborhood of $\Sigma$ whence the
    restriction of $u$ is well-defined. Note that for a general
    element of $\Sec[-\infty](E)$ this would not be possible. Thus
    \eqref{eq:initial-conditions-of-distri-inhom-solution} is
    meaningful and we have $\iota^\# u = \iota^\#w = u_0$ as well as
    $\iota^\# \nabla^E_{\mathfrak{n}} u = \iota^\#
    \nabla^E_{\mathfrak{n}} w = \dot{u}_0$. Hence $u$ has all required
    properties. Note that $u$ depends weak$^*$ continuously on $v$ as
    $G^+_M$ is weak$^*$ continuous. Moreover, $w$ depends continuously
    on $u_0$, $\dot{u}_0$ with respect to the $\Cinfty$- and
    $\Cinfty_0$-topologies. Finally, suppose that $\widetilde{u}$ is
    another generalized section satisfying the four properties
    \eqref{eq:solution-for-distributional-inhom} -
    \eqref{eq:initial-conditions-of-distri-inhom-solution}. Then $u -
    \widetilde{u}$ solves the homogeneous wave equation and has
    singular support away from $\Sigma$, too. Thus we can speak of
    initial conditions of $u - \widetilde{u}$ on $\Sigma$ which are
    now identically zero. Let $\Sigma'$ be another Cauchy hypersurface
    separating $\Sigma$ and $J^+_M(\supp v)$ as in
    Figure~\ref{fig:supps-of-generalized-solution}, which we clearly
    can find. Then in the globally hyperbolic spacetime
    $I^-_M(\Sigma')$ we have a smooth solution $u - \widetilde{u}$ of
    the homogeneous wave equation with vanishing initial
    conditions. Hence $(u - \widetilde{u}) \at{I^-_M(\Sigma')} = 0$ by
    the uniqueness
    Theorem~\ref{theorem:uniqueness-of-hom-equation-with-initial-values}. But
    this implies that the generalized section $u - \widetilde{u}$
    meets the conditions of
    Theorem~\ref{theorem:future-past-comp--hom-solution-is-zero},
    which gives $u - \widetilde{u} = 0$ everywhere.
\end{proof}

\begin{remark}
    \label{remark:Ck-for-distrib-inhom}
    With other words, we have again a well-posed Cauchy problem in
    this more general context of generalized sections as
    inhomogeneities. Note that due to $u \in \Sec[-\infty](E)$ the
    weak$^*$ continuity is the best we can hope for. Analogously to
    Theorem~\ref{theorem:global-solutions},
    \refitem{item:Ck-global-solution} we can also solve the analogous
    Cauchy problem with finite differentiability of the initial
    conditions. In this case we can have singular support outside of
    $J^+_M( \supp v)$ but only a rather mild one: on $M \setminus
    J^+_M( \supp v)$ the solution $u$ is $\Fun$ whence the
    restrictions to $\Sigma$ still make sense.
\end{remark}

%
%

\subsection{The Image of the Green Operators}
\label{satz:image-of-green-operators}

In this section we want to characterize the image of the Green
operators $G^\pm_M$ in $\Secinfty(E)$ in some more detail. Since
$\supp (G^\pm_M u) \subseteq J_M( \supp u)$ for $u \in \Secinfty_0(E)$
we see already here that in general, the maps $G^\pm_M$ can not be
surjective. In general, $M$ can not be written as $J_M(K)$ for a
compact subset. This would require a \emph{compact} Cauchy
hypersurface $\Sigma$. These considerations motivate the following
definition:
\begin{definition}[The space $\Sec_{\mathrm{sc}}(E)$]
    \label{definition:spacelike-compact-sections}%
    \index{Spacelike compact support}%
    Let $k \in \mathbb{N} \cup \{+\infty\}$. For a time-oriented
    Lorentz manifold we denote by $\Sec_{\mathrm{sc}}(E) \subseteq
    \Sec(E)$ those section $u$ for which there exists a compact subset
    $K \subseteq M$ with $\supp u \subseteq J_M(K)$.
\end{definition}
Of course, we are mainly interested in the globally hyperbolic
case. The notion ``sc'' refers to \emph{spacelike compact support}. We
want to endow the subspace $\Sec_{\mathrm{sc}}(E) \subseteq \Sec(E)$
with a suitable topology analogous to the one of $\Sec_0(E)$. Indeed,
$\Sec_{\mathrm{sc}}(E)$ is dense in $\Sec(E)$ for the
$\Cinfty$-topology as $\Sec_0(E) \subseteq \Sec_{\mathrm{sc}}(E)
\subseteq \Sec(E)$ is already dense. Thus we need a finer topology for
$\Sec_{\mathrm{sc}}(E)$ to have good completeness properties. Since
$J_M(K)$ is closed in $M$ on a globally hyperbolic spacetime we can
use Lemma~\ref{lemma:sections-with-compact-support-in-fixed-compacta}
to construct a LF topology for $\Sec_{\mathrm{sc}}(E)$ as follows: For
$K \subseteq K'$ we have $J_M(K) \subseteq J_M(K')$ whence
\begin{equation}
    \label{eq:inclusion-for-bigger-compact}
    \Sec_{J_M(K)}(E) \hookrightarrow \Sec_{J_M(K')}(E)
\end{equation}
is continuous in the $\Fun_{J_M(K)}$- and $\Fun_{J_M(K')}$-topology
and we have a closed image. Since the induced topology from the
$\Fun_{J(K')}$-topology on the image of
\eqref{eq:inclusion-for-bigger-compact} is again the
$\Fun_{J(K)}$-topology we indeed have a nice embedding. Finally, for
an exhausting sequence $K_n \subseteq M$ of compacta we have
eventually $J_M(K) \subseteq J_M(K_n)$. Thus a countable sequence of
subsets exhausts all $J_M(K)$'s. These are the prerequisites for the
strict inductive limit topology analogously to the case of
$\Secinfty_0(E)$ as formulated in
Theorem~\ref{theorem:inductive-limit-topology}. We call the resulting
topology the \emph{$\Fun_{\mathrm{sc}}$-topology}. Without going into
further details we state the consequences literally translating from
Theorem~\ref{theorem:inductive-limit-topology}.
\begin{theorem}[LF topology for $\Sec_{\mathrm{sc}}(E)$]
    \label{theorem:LF-topology-for-spacelike-compact-sections}%
    \index{Strict inductive limit}%
    \index{Locally convex topology!Cinfinitysc-topology@$\Cinfty_{\mathrm{sc}}$-topology}%
    \index{Locally convex topology!Cksc-topology@$\Fun_{\mathrm{sc}}$-topology}%
    Let $(M,g)$ be a time-oriented Lorentz manifold with closed causal
    relation and let $k \in \mathbb{N}_0 \cup \{+\infty\}$. Endow
    $\Sec_{\mathrm{sc}}(E)$ with the inductive limit topology coming
    from \eqref{eq:inclusion-for-bigger-compact}.
    \begin{theoremlist}
    \item \label{item:sc-secs-are-hausdorff-complete-top-vs}
        $\Sec_{\mathrm{sc}}(E)$ is a Hausdorff locally convex complete
        and sequentially complete topological vector space.
    \item \label{item:closed-subspaces-of-sc-secs} All inclusions
        $\Sec_{J_M(K)}(E) \hookrightarrow \Sec_{\mathrm{sc}}(E)$ are
        continuous and the $\Fun_{\mathrm{sc}}$-topology is the finest
        locally convex topology on $\Sec_{\mathrm{sc}}(E)$ with this
        property. Every $\Sec_{J_M(K)}(E)$ is closed in
        $\Sec_{\mathrm{sc}}(E)$ and the induced topology from the
        $\Fun_{\mathrm{sc}}$-topology is again the
        $\Fun_{J_M(K)}$-topology.
    \item \label{item:sc-secs-cauchy-seqence} A sequence $u_n \in
        \Sec_{\mathrm{sc}}(E)$ is a $\Fun_{\mathrm{sc}}$-Cauchy
        sequence iff there is a compact subset $K \subseteq M$ with
        $u_n \in \Sec_{J_M(K)}(E)$ and $u_n$ is a
        $\Fun_{J_M(K)}$-Cauchy sequence. An analogous statement holds
        for convergent sequences.
    \item \label{item:sc-secs-continuity} If $V$ is a locally convex
        vector space then a linear map $\Phi: \Sec_{\mathrm{sc}}(E)
        \longrightarrow V$ is $\Fun_{\mathrm{sc}}$-continuous iff all
        restrictions $\Phi \at{\Sec_{J_M(K)}}: \Sec_{J_M(K)}(E)
        \longrightarrow V$ are $\Fun_{J_M(K)}$-continuous. It suffices
        to check this for an exhausting sequence of compacta.
    \item \label{item:sc-top-and-Ck-top} If in addition $M$ is
        globally hyperbolic with a smooth spacelike Cauchy
        hypersurface $\Sigma$ then $\Sec_{\mathrm{sc}}(E) = \Sec(E)$
        iff $\Sigma$ is compact in which case the
        $\Fun_{\mathrm{sc}}$- and the $\Fun$-topologies
        coincide. Otherwise the $\Sec_{\mathrm{sc}}$-topology is
        strictly finer. In fact,
        \begin{equation}
            \label{eq:sc-top-restriction}
            \iota^\#: \Sec_{\mathrm{sc}}(E)
            \longrightarrow
            \Sec_0(\iota^\# E)
        \end{equation}
        is a surjective linear map which is continuous in the
        $\Fun_{\mathrm{sc}}$- and $\Fun_0$-topology. It furthermore
        has continuous right inverses.
    \end{theoremlist}
\end{theorem}
\begin{proof}
    First we note that for an exhausting sequence $K_n$ of compacta we
    have $K \subseteq K_n$ for all compacta and $n$ suitably
    large. Thus countably many $K_n$ will suffice to specify the
    inductive limit topology of $\Sec_{\mathrm{sc}}(E)$. Since we have
    the continuous embedding with closed image
    \eqref{eq:inclusion-for-bigger-compact} and the correct induced
    topology on the image, we are indeed in the situation of a
    countable strict inductive limit of Fr\'echet spaces, see again
    e.g. \cite[Sect.~4.6]{jarchow:1981a} for details. In particular,
    the parts \refitem{item:sc-secs-are-hausdorff-complete-top-vs} -
    \refitem{item:sc-secs-continuity} are consequences of the general
    properties of LF topologies. For the last part it is clear that if
    $\Sigma$ is a compact Cauchy hypersurface then $J_M(\Sigma) = M$
    whence the $\Fun_{\mathrm{sc}}$-topology simply coincides with the
    $\Fun$-topology as $\Sec_{J_M(\Sigma)}(E)$ \emph{is} already the
    inductive limit. Thus assume that $\Sigma$ is not
    compact. Moreover, let $K \subseteq \Sigma$ be a compact subset in
    $\Sigma$. Then the restriction of a section $u \in
    \Sec_{J_M(K)}(E)$ to $\Sigma$ yields a section $\iota^\# u \in
    \Sec_{K}(\iota^\#E)$. Moreover, we clearly have that the linear
    map
    \[
    \iota^\#: \Sec_{J_M(K)} \ni u
    \; \mapsto \;
    \iota^\#u \in \Sec_K(i^\#E)
    \tag{$*$}
    \]
    is continuous. This is clear from the concrete form of the
    seminorms defining the $\Fun$-topology on $M$ and $\Sigma$,
    respectively. Here we see that in general
    \[
    \iota^\#: \Secinfty_{\mathrm{sc}}
    \longrightarrow
    \Secinfty_0(\iota^\#E),
    \]
    hence $\Secinfty_{\mathrm{sc}}(E) \subsetneq \Secinfty(E)$ follows
    from $\Secinfty_0(i^\# E) \subsetneq \Secinfty(i^\#E)$ at
    once. Moreover, since ($*$) is continuous for all such $K
    \subseteq M$ we see that also
    \[
    \Secinfty_{J_M(K)}(E)
    \longrightarrow \Secinfty_K(\iota^\#E)
    \hookrightarrow \Secinfty_0(\iota^\#E)
    \tag{$**$}
    \]
    is continuous. Now we use that an exhausting sequence $K_n
    \subseteq \Sigma$ \emph{inside of} $\Sigma$ still provides an
    exhausting sequence $J_M(K_n) \subseteq M$ of $M$. Thus we can use
    ($**$) to conclude the continuity of \eqref{eq:sc-top-restriction}
    by part \refitem{item:sc-secs-continuity}. Conversely, using the
    fact that $M$ is diffeomorphic to $\mathbb{R} \times \Sigma$ we
    can extend a section $u_0 \in \Secinfty_0(\iota^\#E)$ to $M$ by
    using the \emph{prolongation map}\index{Prolongation map}
    \[
    \prol(u_0) \at{(t,\sigma)}
    = u_0(\sigma),
    \tag{$\star$}
    \]
    i.e. $\prol(u_0) = \pr_2^\# u_0$. Note that the vector bundle $E$
    on $M$ can be identified with the pull-back bundle $\pr_2^\#
    \iota^\#E \longrightarrow M$ since the time axis is topologically
    trivial. Here $\pr_2: M = \mathbb{R} \times \Sigma \longrightarrow
    \Sigma$ is the projection onto $\Sigma$ as usual. Note that
    ($\star$) makes use of the diffeomorphism $M \simeq \mathbb{R}
    \times \Sigma$ and is \emph{not} canonical. If $u_0 \in
    \Sec_K(\iota^\#E)$ then $\prol(u_0) \in \Sec_{\pr_2^{-1}(K)}(E)
    \subseteq \Sec_{J_M(K)}(E)$ since clearly $\pr_2^{-1}(K)$ is
    inside $J_M(K)$, see
    Figure~\ref{fig:pre-image-of-K-under-projection}, as the curve $t
    \mapsto (t, \sigma)$ is clearly timelike, see also the proof of
    Proposition~\ref{proposition:existence-of-local-solutions-with-initial-values}.
    \begin{figure}
        \centering
        \input{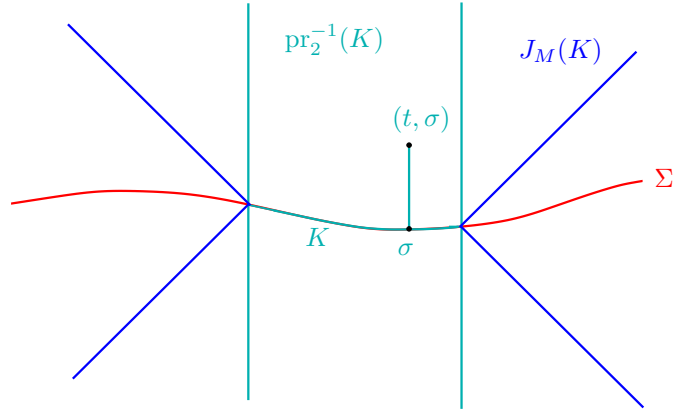}
        \caption{\label{fig:pre-image-of-K-under-projection}%
          The pre-image of a compactum $K \subseteq \Sigma$ under the
          projection $\pr_2$ is inside the causal future of $K$ in a
          globally hyperbolic manifold.
        }
    \end{figure}
    Since $\prol(u_0)$ is ``constant'' in time it is easy to see that
    $\prol: \Sec_K(\iota^\#E) \longrightarrow \Sec_{\pr_2^{-1}(K)}(E)
    \subseteq \Sec_{J_M(K)}(E)$ is continuous. Then also
    \[
    \prol: \Sec_K(\iota^\#E) \longrightarrow \Sec_{\mathrm{sc}}(E)
    \]
    is continuous by part
    \refitem{item:closed-subspaces-of-sc-secs}. Now the
    characterization of the $\Cinfty_0$-topology asserts that $\prol:
    \Sec_0(\iota^\#E) \longrightarrow \Sec_{\mathrm{sc}}(E)$ is
    continuous as well since $K \subseteq \Sigma$ was an arbitrary
    compact subset, see again
    Theorem~\ref{theorem:inductive-limit-topology},
    \refitem{item:LFcontinuity}. Since by construction $\iota^\# \prol
    = \id$ we finally showed the last part. Note that the
    $\Fun_{\mathrm{sc}}$-topology is clearly strictly finer because
    $\Sec_{\mathrm{sc}}(E) \subseteq \Sec(E)$ is dense in the
    $\Sec$-topology but $\Sec_{\mathrm{sc}}(E)$ is complete in the
    $\Sec_{\mathrm{sc}}$-topology.
\end{proof}

\begin{remark}[The $\Fun_{\mathrm{sc}}$-topology]
    \label{remark:spacelike-compact-topology}
    We can repeat the discussion of continuous maps also for the
    $\Fun_{\mathrm{sc}}$-topology in complete analogy to the case of
    the $\Fun_0$-topology as in
    Subsection~\ref{subsec:cont-maps-between-test-section-spaces} and
    Subsection~\ref{subsec:continuity-of-diffops}. In particular, any
    differential operator $D \in \Diffop^k(E;F)$ of order $k$ gives a
    continuous linear map
    \begin{equation}
        \label{eq:diffop-on-sc-is-cont}
        D: \Sec[k+\ell]_{\mathrm{sc}}(E) \longrightarrow
        \Sec[\ell]_{\mathrm{sc}}(F)
    \end{equation}
    with respect to the $\Fun[k+\ell]_{\mathrm{sc}}$- and the
    $\Fun[\ell]_{\mathrm{sc}}$-topology for all $\ell \in \mathbb{N}_0
    \cup \{+\infty\}$. We also have approximation theorems resulting
    from the ones in Subsection~\ref{Approximations}.
\end{remark}

The space $\Secinfty_{\mathrm{sc}}(E) \subseteq \Secinfty(E)$ is the
natural target space for the Green operators $G^\pm_M$ since the
causality requirement
\begin{equation}
    \label{eq:supp-property-of-greenops}
    \supp (G^\pm_M(u)) \subseteq J_M( \supp u)
\end{equation}
immediately implies $G^\pm_M(u) \in \Sec_{\mathrm{sc}}(E)$. The
continuity of $G^\pm_M$ with respect to the $\Cinfty$-topology on
$\Secinfty(E)$ implies also the continuity with respect to the in
general strictly finer $\Cinfty_{\mathrm{sc}}$-topology:
\begin{proposition}
    \label{proposition:greenop-is-cont-for-sc-top-also}
    Let $(M,g)$ be a time-oriented Lorentz manifold with closed causal
    relation. Assume that $G^\pm_M$ are advanced or retarded Green
    operators for a normally hyperbolic differential operator $D \in
    \Diffop^2(E)$. Then
    \begin{equation}
        \label{eq:greenop-is-cont-for-sc-top-also}
        G^\pm_M: \Secinfty_0(E) \longrightarrow
        \Secinfty_{\mathrm{sc}}(E)
    \end{equation}
    is continuous with respect to the $\Cinfty_{\mathrm{sc}}$- and
    $\Cinfty_0$-topology.
\end{proposition}
\begin{proof}
    We know that $G^\pm_M: \Secinfty_0(E) \longrightarrow
    \Secinfty(E)$ is continuous by definition. Thus let $K \subseteq
    M$ be compact then $G^\pm_M: \Secinfty_K(E) \longrightarrow
    \Secinfty(E)$ is continuous in the $\Cinfty_K$- and
    $\Cinfty$-topology be
    Theorem~\ref{theorem:inductive-limit-topology},
    \refitem{item:LFcontinuity}. Since the image is in
    $\Secinfty_{J_M(K)}(E)$ and the $\Cinfty_{J_M(K)}$-topology of
    $\Secinfty_{J_M(K)}$ is the subspace topology inherited from
    $\Secinfty(E)$ we have continuity of
    \[
    G_M^\pm: \Secinfty_K(E) \longrightarrow \Secinfty_{J_M(K)}(E)
    \]
    for all compact subsets $K \subseteq M$. By
    Theorem~\ref{theorem:LF-topology-for-spacelike-compact-sections},
    \refitem{item:closed-subspaces-of-sc-secs} we conclude that also
    \[
    G_M^\pm: \Secinfty_K(E) \longrightarrow \Secinfty_{\mathrm{sc}}(E)
    \]
    is continuous. Since $K$ was arbitrary, by
    Theorem~\ref{theorem:inductive-limit-topology},
    \refitem{item:LFcontinuity} we have the continuity of
    \eqref{eq:greenop-is-cont-for-sc-top-also}.
\end{proof}

Now we come to the main result of this section which describes the
image of the \emph{difference} of the advanced and the retarded Green
operator: as already in the local case we consider the
\emph{propagator}
\begin{equation}
    \label{eq:GlobalPropagator}%
    \index{Propagator!global}%
    G_M = G^+_M - G^-_M: \Secinfty_0(E) \longrightarrow
    \Secinfty_{\mathrm{sc}}(E),
\end{equation}
if $G^\pm_M$ are advanced and retarded Green operators for a normally
hyperbolic differential operator $D$. Here we have the following
statement:
\begin{theorem}
    \label{theorem:propagator-complex}%
    \index{Propagator!image}%
    \index{Propagator!kernel}%
    Let $(M,g)$ be a time-oriented Lorentz manifold with closed causal
    relation. Assume that a normally hyperbolic differential operator
    $D \in \Diffop^2(E)$ has advanced and retarded Green operators
    $G^\pm_M$.
    \begin{theoremlist}
    \item \label{item:propagator-complex} The sequence of linear maps
        \begin{equation}
            \label{eq:propagator-compex}
            0
            \longrightarrow \Secinfty_0(E)
            \stackrel{D}{\longrightarrow} \Secinfty_0(E)
            \stackrel{G_M}{\longrightarrow} \Secinfty_{\mathrm{sc}}(E)
            \stackrel{D}{\longrightarrow} \Secinfty_{\mathrm{sc}}(E)
        \end{equation}
        is a complex of continuous linear maps.
    \item \label{item:propator-complex-is-exact-at-beginning} The
        complex \eqref{eq:propagator-compex} is exact at the first
        $\Secinfty_0(E)$.
    \item \label{item:propagator-complex-exact-on-globally-hyperbolic}
        If $(M,g)$ is globally hyperbolic then
        \eqref{eq:propagator-compex} is exact everywhere.
    \end{theoremlist}
\end{theorem}
\begin{proof}
    The continuity refers to the natural topologies of
    $\Secinfty_0(E)$ and $\Secinfty_{\mathrm{sc}}(E)$, respectively,
    and follows from Remark~\ref{remark:spacelike-compact-topology}
    and
    Proposition~\ref{proposition:greenop-is-cont-for-sc-top-also}. From
    the very definition of Green operators it follows that $G_M \circ
    D = 0 = D \circ G_M$ on $\Secinfty_0(E)$. This shows that
    \eqref{eq:propagator-compex} is a complex. To show exactness at
    the first $\Secinfty_0(E)$ we have to show that $D$ is injective
    on $\Secinfty_0(E)$. Thus let $u \in \Secinfty_0(E)$ with $Du=0$
    be given. Then $0 = G_M^+ Du = u$ shows the injectivity of
    $D$. For the last part assume that $(M,g)$ is globally
    hyperbolic. To show exactness at the second $\Secinfty_0(E)$ we
    have to show $\image D\at{\Secinfty_0(E)} = \ker G_M$. We already
    know ``$\subseteq$'' hence we consider $u \in \Secinfty_0(E)$ with
    $G_M u = 0$. We know that $v = G^+_M u = G^-_M u$ has support in
    $J^+_M(\supp u)$ as well as in $J^-_M( \supp u)$ as $G^\pm_M$ are
    advanced and retarded Green operators. This shows $\supp v
    \subseteq J^+_M(\supp u) \cap J^-_M( \supp u)$ which is
    compact. Indeed, the intersection of $J^+_M$ and $J^-_M$ of
    compact subsets like $\supp u$ is again compact on a globally
    hyperbolic spacetime. This implies $v \in \Secinfty_0(E)$. Since
    in general $D G^+_M u = u$ we see $u = Dv$ with $u$ compactly
    supported. This shows exactness at the second place. To show
    exactness at the third place we have to show that $u \in
    \Secinfty_{\mathrm{sc}}(E)$ with $Du=0$ is actually of the form $u
    = G_M v$ with $v \in \Secinfty_0(E)$. Thus let $u \in
    \Secinfty_{\mathrm{sc}}(E)$ be such a section. For $u \in
    \Secinfty_{\mathrm{sc}}(E)$, the support of $u$ is contained in
    some $J_M(K')$ with $K' \subseteq M$ compact. Choosing an open
    neighborhood of $K'$ with compact closure $K$, i.e. $K' \subseteq
    \mathring{K} \subseteq K$, we see that $\supp u \subseteq I^+_M(K)
    \cup I^-_M(K)$. The two subsets $I^\pm_M(K)$ provide an open cover
    of the open subset $I_M(K) \subseteq M$. Thus we can find a
    subordinate partition of unity $\chi^+, \chi^- \in
    \Cinfty(I_M(K))$ with $\supp \chi^\pm \subseteq I^\pm_M(K)$ and
    $\chi^+ + \chi^- = 1$ on $I_M(K)$. Setting $u^\pm = \chi^\pm u$ we
    have $u = u^+ + u^-$ with $\supp u^\pm \subseteq I^\pm_M(K)
    \subseteq J^\pm_M(K)$. From $Du = 0$ we see $D u^+ = - D u^-$
    which we denote by $v$. Since $\supp D u^\pm \subseteq \supp
    u^\pm$ we conclude $\supp v \subseteq \supp u^+ \cap \supp u^-
    \subseteq J^+_M(K) \cap J^-_M(K)$ which is compact, i.e. $v \in
    \Secinfty_0(E)$. In particular we can apply $G^\pm_M$ to $v$. We
    want to show $G^+_M D u^+ = u^+$: Even though $D u^+ = v$ has
    compact support we can not directly apply the defining property of
    $G^+_M$ since $u^+$ does not have compact support. However, we can
    interpret $u^+$ in a distributional sense and compute for a test
    section $\varphi \in \Secinfty_0(E^*)$
    \begin{align*}
        \int_M \varphi(p) \cdot (G^+_M D u^+)(p) \: \mu_g
        &\stackrel{(*)}{=}
        \int_M (F^-_M \varphi) \cdot (D u^+)(p) \: \mu_g \\
        &\stackrel{(**)}{=}
        \int_M (D^\Trans F^-_M \varphi)(p) \cdot u^+(p) \: \mu_g \\
        &=
        \int_M \varphi(p) \cdot u^+(p) \: \mu_g,
    \end{align*}
    where we have used
    Lemma~\ref{lemma:green-ops-of-transpose-and-green-ops-of-D} in
    ($*$) and integration by parts in ($**$) which is possible since
    $F^-_M \varphi$ has support in $J^-_M( \supp \varphi)$ while $u^+$
    has support in $J^+_M(K)$. Hence the overlap of their supports is
    compact even though their supports are not. Then the above
    computation shows $G^+_M D u^+ = u^+$. Analogously we find $G^-_M
    D u^- = u^-$. Putting these results together gives $G_M v = G^+_M
    v - G^-_M v = G^+_M D u^+ + G^-_M Du^- = u^+ + u^- =
    u$. Therefore, $u$ is in the image of $G_M$ with a pre-image in
    $\Secinfty_0(E)$ as wanted.
\end{proof}

\begin{remark}[Propagator]
    \label{remark:propagator}%
    \index{Propagator}%
    The simple description of the image and kernel of the operator
    $G_M = G^+_M - G^-_M$ has many important consequences. In physics
    in (quantum) field theory this operator is called the
    \emph{propagator} which is one of the most crucial ingredients in
    any perturbative (quantum) field theory. It also appears as the
    kernel of the Poisson bracket in classical field theory which we
    will discuss in Section~\ref{satz:poisson-algebra}.
\end{remark}

As an application of the operator $G_M$ we obtain a global version of
Lemma~\ref{lemma:solution-integral-formula} expressing the solution of
the homogeneous Cauchy problem in terms of the initial data:
\begin{theorem}
    \label{theorem:global-integral-formula-for-solution}%
    \index{Cauchy problem!homogeneous}%
    Let $(M,g)$ be a globally hyperbolic spacetime and let $\iota:
    \Sigma \hookrightarrow M$ be a smooth spacelike Cauchy
    hypersurface. Let $D \in \Diffop^2(E)$ be normally hyperbolic and
    let $F^\pm_M$ be the advanced and retarded Green operators of
    $D^\Trans$. Then the solution $u \in \Secinfty_{\mathrm{sc}}(E)$
    of the homogeneous wave equation $D u = 0$ with initial values
    $\iota^\#u = u_0$ and $\iota^\# \nabla^E_{\mathfrak{n}} u =
    \dot{u}_0$ on $\Sigma$ is determined by
    \begin{equation}
        \label{eq:global-integral-formula-for-solution}
        \int_M \varphi(p) \cdot u(p) \mu_g(p)
        = \int_\Sigma
        \left(
            (\nabla^E_{\mathfrak{n}} F_M(\varphi))(\sigma) \cdot
            u_0(\sigma)
            - F_M(\varphi)(\sigma) \cdot \dot{u}_0(\sigma)
        \right) \mu_\Sigma
    \end{equation}
    for $\varphi \in \Secinfty_0(E^*)$.
\end{theorem}
\begin{proof}
    The proof is literally the same as for
    Lemma~\ref{lemma:solution-integral-formula}. Therefore it will be
    enough to sketch the arguments. We consider the sections
    $F^\pm_M(\varphi) \in \Secinfty_{\mathrm{sc}}(E^*)$ which have
    supports $\supp F^\pm_M(\varphi) \subseteq J^\pm_M( \supp
    \varphi)$. Taking covariant derivatives and pairing with $u$ gives
    the vector field
    \[
    X^\pm
    = \left(
        (\SymD^{E^*} F^\pm_M(\varphi)) \cdot u
        - F^\pm_M(\varphi) \cdot( \SymD^E u)
    \right)^\#
    \in \Secinfty(TM)
    \]
    which has again support in $J^\pm_M( \supp \varphi)$.
    \begin{figure}
        \centering
        \input{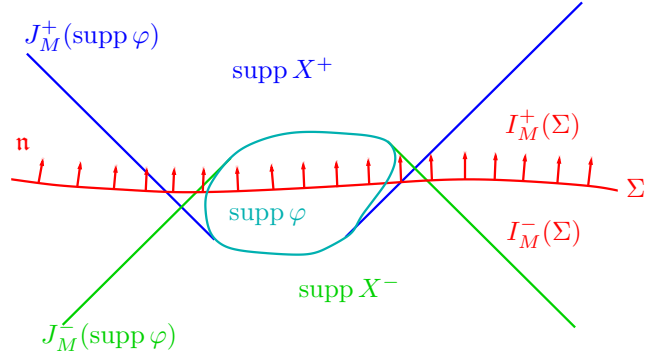}
        \caption{\label{fig:supp-of-X-plus-minus}%
          The supports of the functions $\varphi$ and $X^\pm$.
        }
    \end{figure}
    In particular, $\supp X^\pm \cap I^\mp_M(\Sigma)$ as well as
    $\supp X^\pm \cap \Sigma$ are (pre-) compact and hence the
    following integrations will be well-defined. Integrating over
    $I^-_M(\Sigma)$ the unit normal field $\mathfrak{n}$ is pointing
    outwards as it is future-directed. Conversely, integrating over
    $I^+_M(\Sigma)$ the vector field $- \mathfrak{n}$ is pointing
    outwards. Thus by Theorem~\ref{theorem:gauss-theorem} we get
    \[
    \int_{I^\pm_M(\Sigma)} \divergenz (X^\mp) \: \mu_g
    = \mp \int_\Sigma g(X^\mp\mathfrak{n}) \: \mu_\Sigma,
    \]
    where we of course have restricted $X^\pm$ to $\Sigma$ on the
    right hand side. For the left hand side we obtain
    \begin{align*}
        \int_{I^\pm_M(\Sigma)} \divergenz (X^\mp) \: \mu_g
        &=
        \int_{I^\pm_M(\Sigma)} \divergenz
        \left(
            (\SymD^{E^*} F^\pm_M(\varphi)) \cdot u
            - F^\pm_M(\varphi) \cdot( \SymD^E u)
        \right)^\#
        \mu_g \\
        &\stackrel{\mathclap{\eqref{eq:difference-D-and-Dtrans}}}{=}
        \quad
        \int_{I^\pm_M(\Sigma)}
        \left(
            (D^\Trans D^\mp_M (\varphi)) \cdot u
            - F^\pm(\varphi)(Du)
        \right) \mu_g \\
        &= \int_{I^\pm_M(\Sigma)} \varphi \cdot u \: \mu_g,
    \end{align*}
    since $F^\pm_M$ are the Green operators of $D^\Trans$ and $D u =
    0$. For the right hand side we have
    \begin{align*}
        \mp \int_\Sigma g(X^\mp,\mathfrak{n}) \: \mu_\Sigma
        &=
        \mp \int_\Sigma
        \left(
            (\SymD^{E^*} F^\mp_M(\varphi)) \cdot u
            - F^\mp_M(\varphi) \cdot \SymD^E u
        \right) \cdot \mathfrak{n} \: \mu_\Sigma \\
        &= \mp \int_\Sigma
        \left(
            (\nabla^{E^*}_{\mathfrak{n}} F^\mp_M(\varphi)) \cdot u
            - F^\mp_M(\varphi) \cdot \nabla^E_{\mathfrak{n}} u
        \right) \mu_g \\
        &= \mp \int_\Sigma
        \left(
            (\nabla^{E^*}_{\mathfrak{n}} F^\mp_M(\varphi)) \cdot u_0
            - F^\mp_M(\varphi) \cdot \dot{u}_0
        \right) \mu_\Sigma,
    \end{align*}
    with an analogous computation as in
    Lemma~\ref{lemma:solution-integral-formula}. Putting things
    together gives \eqref{eq:global-integral-formula-for-solution}.
\end{proof}

\begin{remark}
    \label{remark:global-integral-formula-for-solution}
    From this formula we see that the homogeneous Cauchy problem can
    again be encoded completely in terms of the Green operators. Since
    also the inhomogeneous Cauchy problem with vanishing initial
    conditions can be solved by means of the Green operators thanks to
    Theorem~\ref{theorem:solution-of-wave-eq-for-distributional-inhom}
    we see that the Cauchy problem and the construction of the Green
    operators are ultimately the same problem.
\end{remark}


%% file: poisson.tex
%
%

In this section we describe a first attempt to establish a Hamiltonian
picture for the wave equation based on a certain Poisson algebra of
observables coming from the canonical symplectic structure on the
space of initial conditions. Throughout this section, $(M,g)$ will be
globally hyperbolic. For the vector bundle $E \longrightarrow M$ we
have to be slightly more specific: We choose $E$ to be a \emph{real}
vector bundle. The reason will be to get the correct linearity
properties of the Poisson bracket later. From a physical point of
view, many of the complex vector bundles actually arise as
complexifications of real ones. Then the wave operators in question
have the additional property to commute with the complex conjugation
of the sections of the complexified bundles. This will be important in
applications in physics later on, in particular for CPT-like theorems
in quantum field theories, see e.g \cite{streater.wightman:1964a,
  haag:1993a}. For an overview on the geometrical aspects of
(finite-dimensional) classical mechanics we refer to
\cite{abraham.marsden:1985a, marsden.ratiu:1999a, waldmann:2007a}

%
%

\subsection{Symmetric Differential Operators}
\label{satz:symmetric-diffops}

Now we equip the vector bundle $E$ with an additional structure,
namely a fiber metric $h$. In most applications this fibre metric will
be positive definite, a fact which we shall not use though. In any
case, the fibre metric induces a musical isomorphism $\flat: E
\longrightarrow E^*$ with inverse $\#: E^* \longrightarrow E$ as
usual. On sections we have
\begin{equation}
    \label{eq:flatMap}%
    \index{Musical isomorphism}%
    \index{Fiber metric}%
    \flat: \Secinfty(E) \ni u \; \mapsto \; u^\flat = h(u,\argument) \in
    \Secinfty(E^*).
\end{equation}
There should be no confusion with the sharp and flat map coming from
the Lorentz metric $g$. Using this additional structure one can define
symmetric differential operators as usual:
\begin{definition}[Symmetric differential operators]
    \label{definition:symmetric-diffop}%
    \index{Differential operator!symmetric}%
    \index{Differential operator!adjoint}%
    Let $(E,h)$ be a real vector bundle with fibre metric and $D \in
    \Diffop^\bullet(E)$. Then the adjoint of $D$ with respect
    to $h$ is defined to be the unique $D^* \in \Diffop^\bullet(E)$
    with
    \begin{equation}
        \label{eq:adjoint-diffop-wrp-fibre-metric}
        \int_M h(D^*u,v) \mu_g
        = \int_M h(u,Dv) \mu_g
    \end{equation}
    for all $u,v \in \Secinfty_0(E)$. The operator $D$ is called
    symmetric if
    \begin{equation}
        \label{eq:symmetric-diffop}
        D = D^*.
    \end{equation}
\end{definition}
\begin{remark}[Symmetric differential operators]
    \label{remark:symmetric-diffop}
    ~
    \begin{remarklist}
    \item \label{item:adjoint-and-transpose} The definition of the
        adjoint $D^*$ with respect to $h$ is well-defined
        indeed. Namely, if $D \in \Diffop^k(E)$ then one has
        \begin{equation}
            \label{eq:adjoint-and-transpose}
            D^*u = (D^\Trans u^\flat)^\#
        \end{equation}
        with the adjoint operator $D^\Trans \in \Diffop^k(E^*)$ as we
        discussed it before in
        Theorem~\ref{theorem:agjoint-for-fixed-density}. This follows
        from the simple computation
        \begin{equation}
            \label{eq:adjoint-and-tranpose-proof}
            \int_M h((D^\Trans u^\flat)^\#,v) \: \mu_g
            = \int_M (D^\Trans u^\flat) \cdot v \: \mu_g
            = \int_M u^\flat \cdot Dv \: \mu_g
            = \int_M h(u,Dv) \mu_g,
        \end{equation}
        which shows that \eqref{eq:adjoint-and-transpose} solves the
        condition \eqref{eq:adjoint-diffop-wrp-fibre-metric}. It is
        clear that $D^*$ is again a differential operator of the same
        order as $D$ and it is necessarily unique since the inner
        product is non-degenerate.
    \item \label{item:adjoining-is-involution} The adjoint $D^*$
        depends on $h$ but also on the density $\mu_g$ in the
        integration \eqref{eq:adjoint-diffop-wrp-fibre-metric}. The
        map $D \mapsto D^*$ is a linear involutive anti-automorphism,
        i.e. we have
        \begin{equation}
            \label{eq:adjoining-is-involution}
            (D^*)^* = D
            \quad \textrm{and} \quad
            (D \widetilde{D})^* = \widetilde{D}^* D^*
        \end{equation}
        for $D, \widetilde{D} \in \Diffop^\bullet(E)$.
    \item \label{item:hermitian-diffops}%
        \index{Fiber metric!Hermitian}%
        In the case of a complex vector bundle one proceeds similarly:
        for a given (pseudo-) Hermitian fibre metric one defines the
        adjoint $D^*$ by the same condition
        \eqref{eq:adjoint-diffop-wrp-fibre-metric}. Now $D \mapsto
        D^*$ is antilinear in addition to
        \eqref{eq:adjoining-is-involution} and $\Diffop^\bullet(E)$
        becomes a $^*$-algebra over $\mathbb{C}$ by this
        choice. Differential operators with $D = D^*$ are now called
        \emph{Hermitian}. A particular case is obtained for a
        complexified vector bundle $E_{\mathbb{C}} = E \tensor
        \mathbb{C}$. If $h$ is a fibre metric on $E$ then it induces a
        Hermitian fibre metric on $E_{\mathbb{C}}$ by setting
        \begin{equation}
            \label{eq:induced-hermitian-metric-on-complexification}%
            \index{Complexified bundle}%
            h_{\mathbb{C}}(u \tensor z, v \tensor w) = h(u,v) \cc{z}w
        \end{equation}
        for $u,v \in E_p$ and $z,w \in \mathbb{C}$. Then a symmetric
        operator $D \in \Diffop^\bullet(E)$ yields a Hermitian
        operator $D_{\mathbb{C}} \in \Diffop^\bullet(E_{\mathbb{C}})$
        which commutes in addition with the \emph{complex conjugation}
        of sections.
    \end{remarklist}
\end{remark}
In most physically interesting situations the wave operator $D$ will
be symmetric. As a motivation we consider the following example:
\begin{example}[Symmetric connection d'Alembertian]
    \label{example:conn-dAlembertian-is-symmetric}%
    \index{Covariant derivative!metric}%
    \index{Connection dAlembertian@Connection d'Alembertian!symmetric}%
    Let $(E,h)$ be a real vector bundle with fibre metric
    $h$. Moreover, let $\nabla^E$ be a covariant derivative which is
    \emph{metric} with respect to $h$, i.e.
    \begin{equation}
        \label{eq:metric-cov-derivative}
        \Lie_X h(u,v)
        = h( \nabla^E_X u, v) + h(u, \nabla^E_X v)
    \end{equation}
    for all $X \in \Secinfty(TM)$ and $u,v \in \Secinfty(E)$. We claim
    that in this case the connection d'Alembertian is
    symmetric. Indeed, \eqref{eq:metric-cov-derivative} immediately
    implies that the symmetric covariant derivative operators
    $\SymD^E$ and $\SymD^{E^*}$ with respect to $\nabla^E$ and
    $\nabla^{E^*}$ are intertwined by $\#$ and $\flat$ as follows
    \begin{equation}
        \label{eq:flat-and-sharp-intertwine-SymDE-and-SymDEstar}
        \left(\SymD^E u\right)^\flat = \SymD^{E^*} u^\flat
    \end{equation}
    for $u \in \Secinfty(\Sym^k T^*M \tensor E)$ and $\flat $ applied
    to the $E$-part only. This is a simple verification. But then we
    have for $\dAlembert^\nabla$ by Lemma~\ref{lemma:transpose-of-D},
    \refitem{item:transpose-of-D}
    \begin{equation}
        \label{eq:conn-dAlembert-is-symmetric-proof}
        \left(\dAlembert^\nabla\right)^\Trans u^\flat
        = \SP{g^{-1}, \left(\SymD^{E^*}\right)^2 u^\flat}
        = \SP{g^{-1}, \left(\left(\SymD^E\right)^2 u\right)^\flat}
        = \SP{g, \left(\SymD^E\right)^2u }^\flat
        = \left(\dAlembert^\nabla u\right)^\flat,
    \end{equation}
    since the natural pairing of the $\Sym^2 T^*M$ component with
    $g^{-1}$ obviously commutes with the musical isomorphism $\flat$
    acting only on the $E$-component. But this implies
    \begin{equation}
        \label{eq:conn-dAlember5t-is-symmetric}
        \dAlembert^\nabla = \left(\dAlembert^\nabla\right)^*
    \end{equation}
    as claimed. More generally, if $B \in \Secinfty(\End(E))$ is also
    symmetric with respect to $h$, which is now a pointwise criterion,
    then $D = \dAlembert^\nabla + B$ is symmetric as well.

    This construction is also compatible with complexification: if $h$
    is extended to $E_{\mathbb{C}} = E \tensor \mathbb{C}$ as in
    Remark~\ref{remark:symmetric-diffop},
    \refitem{item:hermitian-diffops} then the connection $\nabla^E$
    also extends to $E_{\mathbb{C}}$ yielding a metric connection
    $\nabla^{E_{\mathbb{C}}}$ with respect to $h_{\mathbb{C}}$. The
    condition \eqref{eq:metric-cov-derivative} is then satisfied for
    real tangent vector fields $X \in \Secinfty(TM)$ while we have to
    replace $X$ by $\cc{X}$ in the first term of the right hand side
    of \eqref{eq:metric-cov-derivative} in general. With this (pseudo)
    Hermitian fibre metric $h_{\mathbb{C}}$ and the covariant
    derivative $\nabla^{E_{\mathbb{C}}}$ the property
    \eqref{eq:flat-and-sharp-intertwine-SymDE-and-SymDEstar} still
    holds, resulting in \eqref{eq:conn-dAlember5t-is-symmetric} for
    the connection d'Alembertian $\dAlembert^\nabla$ on
    $E_{\mathbb{C}}$. Again $\dAlembert^\nabla$ is not only Hermitian
    but also commutes with the complex conjugation of sections of
    $E_{\mathbb{C}}$. Note that for general complex vector bundles
    there is no notion of complex conjugation of sections.
    \index{Connection dAlembertian@Connection d'Alembertian!Hermitian}
\end{example}

From now on we shall focus on a symmetric and normally hyperbolic
differential operator $D$. In fact, we shall also assume that
$\nabla^E$ is metric. Then $D = D^*$ means $B = B^*$ for $D =
\dAlembert^\nabla + B$. Since $D = D^*$ essentially means that we can
identify $D$ with $D^\Trans$ via $\flat$ and $\#$ we expect a similar
relation between the Green operators, extending the already found
relations between $F^\pm_M$ and $G^\pm_M$ as in
Theorem~\ref{theorem:dual-of-green-operators}. In fact, one has the
following characterization:
\begin{proposition}[Symmetry of Green operators]
    \label{proposition:symmetry-of-green-operators}%
    \index{Green operator!symmetric}%
    \index{Green operator!reality}%
    Let $(M,g)$ be globally hyperbolic and let $D \in \Diffop^2(E)$ be
    a normally hyperbolic differential operator on the real vector
    bundle $E$. Assume that $D$ is symmetric with respect to a fibre
    metric $h$ on $E$.
    \begin{propositionlist}
    \item \label{item:greenops-and-flat} For the Green operators of
        $D$ and $D^\Trans$ and $u \in \Secinfty_0(E)$ we have
        \begin{equation}
            \label{eq:greenops-and-flat}
            \left(G^\pm_M u\right)^\flat
            = F^\pm_M u^\flat.
        \end{equation}
    \item \label{item:adjoint-of-greenop} For $u,v \in
        \Secinfty_0(E)$ we have
        \begin{equation}
            \label{eq:adjoint-of-greenop}
            \int_M h\left(u, G^\pm_M v\right) \: \mu_g
            = \int_M h\left(G^\mp_M u, v\right) \: \mu_g.
        \end{equation}
    \item \label{item:greenops-and-complexification} The Green
        operators of the canonical $\mathbb{C}$-linear extension of
        $D$ to $E_{\mathbb{C}} = E \tensor \mathbb{C}$ are the
        canonical $\mathbb{C}$-linear extension of the Green operators
        $G^\pm_M$ of $D$. They still satisfy
        \eqref{eq:greenops-and-flat},
        \begin{equation}
            \label{eq:adjoint-of-extended-greenop}
            \int_M h_{\mathbb{C}}\left(u, G^\pm_M v\right) \: \mu_g
            = \int_M h_{\mathbb{C}}\left(G^\mp_M u, v\right) \: \mu_g
        \end{equation}
        for $u,v \in \Secinfty_0(E_{\mathbb{C}})$ and additionally the
        reality condition
        \begin{equation}
            \label{eq:extended-greenops-complex-conjugation}
            \cc{ G^\pm_M u }
            = G^\pm_M \cc{u}.
        \end{equation}
    \end{propositionlist}
\end{proposition}
\begin{proof}
    Clearly, $u \in \Secinfty_0(E)$ has compact support iff $u^\flat$
    has compact support, making \eqref{eq:greenops-and-flat}
    meaningful. We compute for $\varphi \in \Secinfty_0(E^*)$
    \[
    D^\Trans \left(G^\pm_M \varphi^\#\right)^\flat
    \quad
    \stackrel{\mathclap{\eqref{eq:adjoint-and-transpose}}}{=}
    \quad
    \left(DG^\pm_M \varphi^\#\right)^\flat
    = \left(\varphi^\#\right)^\flat
    = \varphi,
    \]
    since $G^\pm_M$ is a Green operator for $D$. Analogously,
    \[
    \left(G^\pm_M \left(D^\Trans \varphi\right)^\#\right)^\flat
    \stackrel{\mathclap{\eqref{eq:adjoint-and-transpose}}}{=}
    \left(G^\pm_M D \varphi^\#\right)^\flat
    = \left(\varphi^\#\right)^\flat
    = \varphi.
    \]
    Now $\varphi \mapsto \left(G^\pm_M \varphi^\#\right)^\flat$ is
    clear linear and continuous since $\#, \flat$ as well as $G^\pm_M$
    are continuous. Finally, since $\#$ and $\flat$ preserve supports
    we have $\supp \left(G^\pm_M \varphi^\#\right)^\flat \subseteq
    J^\pm_M(\supp \varphi)$. This shows that the map $\varphi \mapsto
    \left(G^\pm_M \varphi^\#\right)^\flat$ is indeed an advanced and
    retarded Green operator for $D^\Trans$, respectively. By
    uniqueness according to
    Corollary~\ref{corollary:green-ops-on-globally-hyperbolic} we get
    \eqref{eq:adjoint-of-greenop}. Using this, we compute
    \begin{align*}
        \int_M h\left(u, G^\pm_M v\right) \: \mu_g
        &=
        \int_M u^\flat \cdot \left(G^\pm_M v\right) \: \mu_g \\
        &\stackrel{\mathclap{\eqref{eq:green-ops-of-transpose-and-green-ops-of-D}}}{=}
        \quad
        \int_M \left(F^\mp_M u^\flat\right) \cdot v \: \mu_g \\
        &\stackrel{\mathclap{\eqref{eq:adjoint-of-greenop}}}{=}
        \quad
        \int_M \left(G^\mp_M u\right)^\flat \cdot v \: \mu_g \\
        &= \int_M h\left(G^\mp_Mu, v\right) \: \mu_g
    \end{align*}
    for $u,v \in \Secinfty_0(E)$. Now consider $u,v \in
    \Secinfty_0(E_{\mathbb{C}})$. Then $\cc{Du} = D\cc{u}$ yields the
    hermiticity $D = D^*$ with respect to $h_{\mathbb{C}}$. With the
    same kind of uniqueness argument we see that the Green operators
    $G^\pm_M$ of $D$, canonically extended to $G^\pm_M:
    \Secinfty_0(E_{\mathbb{C}}) \longrightarrow
    \Secinfty(E_{\mathbb{C}})$, yield the Green operators of the
    extension $D \in \Diffop^2(E_{\mathbb{C}})$. Moreover, we clearly
    have \eqref{eq:extended-greenops-complex-conjugation} by
    construction. But then \eqref{eq:adjoint-of-extended-greenop}
    follows from \eqref{eq:extended-greenops-complex-conjugation} and
    \eqref{eq:adjoint-of-greenop} at once.
\end{proof}

\begin{remark}
    \label{remark:propagator-and-adjoining}
    Extending our notation of the adjoint to more general operators we
    can rephrase the result of \eqref{eq:adjoint-of-greenop} or
    \eqref{eq:adjoint-of-extended-greenop} by saying
    \begin{equation}
        \label{eq:adjoing-of-greenops}
        \left(G^\pm_M\right)^* = G^\mp_M.
    \end{equation}
    Note that
    Proposition~\ref{proposition:symmetry-of-green-operators},
    \refitem{item:greenops-and-complexification} still holds for
    arbitrary Hermitian $D = D^*$ on arbitrary complex vector bundles
    except for \eqref{eq:extended-greenops-complex-conjugation}. In
    both cases, it follows that the propagator $G_M = G^+_M - G^-_M$
    is \emph{antisymmetric}
    \begin{equation}
        \label{eq:propagator-is-antisymmetric}%
        \index{Propagator!antisymmetric}%
        \index{Propagator!anti-Hermitian}%
        G_M^* = -G_M
    \end{equation}
    or \emph{anti-Hermitian} in the complex case, respectively. In the
    complex case we can rescale $G_M$ by $\I$ to obtain a
    \emph{Hermitian} operator
    \begin{equation}
        \label{eq:i-times-propagator-is-hermitian}
        (\I G_M)^* = \I G_M.
    \end{equation}
\end{remark}

%
%

\subsection{Interlude: The Lagrangian and the Hamiltonian Picture}
\label{satz:lagrange-and-hamilton}

To put the following construction in the right perspective we briefly
remind on the Lagrangian and Hamiltonian approach to field equations
as it can be found in various textbooks on classical and quantum field
theory. Most of our present considerations should be taken as
heuristic as it would require a lot more effort to justify them on a
mathematically rigorous basis. They serve as a motivation for our
definition of certain Poisson brackets.

Many field equations in physics arise from an action principle where
an action functional\index{Action functional}\index{Action principle}
is defined on the space of all field configurations on the whole
spacetime by means of a Lagrangian density. Such a Lagrangian
density\index{Lagrangian density} $\mathcal{L}$ can be viewed as a
function on the (first) jet bundle $J^1E$ of $E$ which takes values in
the densities $\Dichten T^*M$ on $M$. Roughly speaking, the $k$-th jet
bundle $J^kE$ of $E$ is a fibre bundle over $M$ whose fibre at $p \in
M$ consists of equivalence classes of Taylor expansions of sections of
$E$ around $p$ up to order $k$. Two sections are called equivalent if
they have the same Taylor expansion at $p$ up to order $k$. This is a
coordinate independent statement whence the jet bundles serve the
following purpose: we can make geometrically sense of the statement
that a map $\mathcal{L}: \Secinfty(E) \longrightarrow
\Secinfty(\Dichten T^*M)$ depends at $p \in M$ only on the first $k$
derivatives of $u \in \Secinfty(E)$ at $p$. In our case one typically
has $k = 1$ and symbolically writes $\mathcal{L}(u, \partial u)$ to
emphasize that $\mathcal{L}(u, \partial u) \at{p}$ depends only on
$u(p)$ and $\frac{\partial u}{\partial x^i}(p)$. Having specified such
a Lagrangian density $\mathcal{L}$ the action $S(u)$ is defined by the
(hopefully existing) integral of $\mathcal{L}(u, \partial u)$ over
$M$. Then the stationary points of the action functional are supposed
to be those sections which satisfy the wave equation. With other words
one wants the \Index{Euler-Lagrange equations} for $\mathcal{L}$ to be
the wave equation under consideration. Note that the precise
formulation of an action principle is far from being trivial: on one
hand, one has to require certain integrability conditions on the
sections in order to have a well-defined action. On the other hand, in
deriving the Euler-Lagrange equations one usually neglects certain
boundary terms or considers only variations with compact support. Thus
it is not evident that the Euler-Lagrange equations really describe
the stationary points of $S$. Even worse, in typical situations the
solutions of the Euler-Lagrange equations yield sections $u$ with
\emph{no} good integrability properties at all. Our wave equation is a
good example as here the non-trivial solutions have to have
non-compact support in timelike directions. This way, it may well
happen that none of the solutions of the Euler-Lagrange equation is in
the domain of definition of the action $S$ at all, except for some
trivial solutions like $u=0$. To handle these difficulties a more
sophisticated variational calculus is required which is not within the
reach for us at this stage. Therefore, we take a more pragmatic point
of view and take the Lagrangian density $\mathcal{L}$ and the
corresponding Euler-Lagrange equations as the starting point instead
of the action $S$ itself. These equations and hence the wave equation
are the ultimate goal anyway.

The idea is now to treat the Euler-Lagrange equations for the
Lagrangian density as Euler-Lagrange equations of a suitably defined
\Index{Lagrangian function} defined on the space of initial
conditions: this way we can interpret the field theoretic wave
equations as a classical mechanical system, though of course with
infinitely many degrees of freedom. The idea is roughly as follows:
the initial conditions of the wave equation are specified on a fixed
smooth spacelike Cauchy hypersurface $\iota: \Sigma \hookrightarrow
M$. There we have to specify the value of the section $u_0 \in
\Secinfty_0(\iota^\#E)$ and the normal derivative $\dot{u}_0 \in
\Secinfty_0(\iota^\#E)$. Mechanically speaking, this corresponds to
the initial position\index{Initial position} and the initial
velocity\index{Initial velocity}. Thus the (velocity-) phase space of
the Lagrangian approach is the tangent bundle of the space of initial
positions in complete analogy to Lagrangian mechanics for
finite-dimensional systems. Since the initial positions are described
by the vector space $\Secinfty_0(\iota^\#E)$ the notion of tangent
bundle is simple: we just have to take $\Secinfty_0(\iota^\#E) \times
\Secinfty_0(\iota^\#E)$, i.e. two copies of the configuration
space. The Lagrange function now consists in evaluating the Lagrange
density on $u_0$ and $\dot{u}_0$ on $\Sigma$ and integrating over
$\Sigma$: this indeed makes sense as the Lagrange density
$\mathcal{L}$ can be written relative to the density $\mu_g$ as
$\mathcal{L}(u,\partial u) = \widetilde{\mathcal{L}}(u, \partial u)
\mu_g$ with a function $\widetilde{\mathcal{L}}(u, \partial u)$ on the
first jet bundle. Then we can take this function and evaluate it on
$u_0$ and $\dot{u}_0$ instead of $u$ and $\partial u$ and consider the
density $\widetilde{\mathcal{L}}(u_0, \dot{u}_0) \mu_\Sigma$ on
$\Sigma$. Again, we ignore the technical details which are less severe
as for the action since we are interested in $u_0$ and $\dot{u}_0$
with compact support anyway. The integration over $\Sigma$ is thus
easily defined.

Having a Lagrangian mechanical point of view for our wave equation we
can try to pass to a Hamiltonian description by the usual
\Index{Legendre transform}. This amounts to the passage from the
tangent bundle to the cotangent bundle of the configuration space
$\Secinfty_0(\iota^\# E)$. While the tangent bundle of a vector space
is conceptually easy, the cotangent bundle is more subtle: here the
fact that $\Secinfty_0(\iota^\#E)$ is infinite-dimensional becomes
crucial. Thus we have to decide which dual we want to take. Of course,
the algebraic dual seems inappropriate whence we take the topological
dual which we identify with $\Sec[-\infty](\iota^\#E^*)$ as usual by
means of $\mu_\Sigma$. Then the cotangent bundle of
$\Secinfty_0(\iota^\#E)$ is $\Secinfty_0(\iota^\#E) \times
\Sec[-\infty](\iota^\#E^*)$. The following proposition shows that this
is indeed a symplectic vector space in a very good sense. We formulate
it for a general vector bundle $F \longrightarrow M$ over an arbitrary
manifold.
\begin{proposition}[Symplectic vector space]
    \label{proposition:symplectic-vector-space}%
    \index{Symplectic vector space}%
    \index{Canonical symplectic form}%
    ~
    \begin{propositionlist}
    \item \label{item:symplectic-structure-on-cotangent-space-of-W}
        Let $W$ be a Hausdorff locally convex topological vector space
        with topological dual and consider $V = W \oplus W'$. Then on
        $V$ the two-form
        \begin{equation}
            \label{eq:canonical-symplectic-form-on-W}
            \omega_{\mathrm{can}}\left(
                (w,\varphi), (w', \varphi')
            \right)
            = \varphi'(w) - \varphi(w')
        \end{equation}
        is antisymmetric and non-degenerate.
    \item \label{item:symplectic-form-on-Configspce-cotangent} Let $F
        \longrightarrow M$ be a real vector bundle. Then on
        $\Secinfty_0(F^* \tensor \Dichten T^*M) \oplus
        \Sec[-\infty](F)$ the two-form
        \begin{equation}
            \label{eq:symplectic-form-on-Configspce-cotangent}
            \omega_{\mathrm{can}}\left(
                (\varphi,u), (\varphi',u')
            \right)
            = u'(\varphi) - u(\varphi')
        \end{equation}
        is antisymmetric and non-degenerate.
    \item \label{item:symplectic-form-on-Configspce-cotangent-with-density}
        Let $F \longrightarrow M$ be a real vector bundle and let
        $\mu$ be a positive density on $M$. Then on $\Secinfty_0(F^*)
        \oplus \Sec[-\infty](F)$ the two form
        \begin{equation}
            \label{eq:symplectic-form-on-Configspce-cotangent-with-density}
            \omega_{\mathrm{can}}\left(
                (\varphi,u), (\varphi',u')
            \right)
            = u'(\varphi \tensor \mu) - u(\varphi' \tensor \mu)
        \end{equation}
        is antisymmetric and non-degenerate.
    \end{propositionlist}
\end{proposition}
\begin{proof}
    Clearly, $\omega_{\textrm{can}}$ is bilinear in all three cases
    and antisymmetric on the nose. Assume that $(w, \varphi) \in W
    \oplus W'$ is such that $\omega_{\textrm{can}}( (w,\varphi),
    \argument) = 0$. Then it follows that $\varphi'(\omega) = 0$ for all
    $\varphi' \in W'$ and $\varphi(w') = 0$ for all $w' \in W$. This
    clearly implies $\varphi = 0$. Since $W$ is Hausdorff, by some
    Hahn-Banach-like statements it follows that $W'$ is large enough
    to separate points, see e.g. \cite[Sect.~7.2]{jarchow:1981a}. Thus
    also $w=0$ follows which proves that
    \eqref{eq:canonical-symplectic-form-on-W} is non-degenerate. The
    second and third part are only special cases.
\end{proof}

Since in our situation we have a canonical positive density on
$\Sigma$, namely $\mu_\Sigma$, we can apply the third part and
conclude that $\Secinfty_0(\iota^\#E) \oplus
\Sec[-\infty](\iota^\#E^*)$ is indeed a symplectic vector space.

Without going into the details we can now use the Lagrange function to
define a Legendre transform by which we can pull back the canonical
symplectic form of the cotangent bundle to the tangent bundle. This
constructions boils down to the following simple map, at least in all
cases relevant for us. By means of the fibre metric $h_\Sigma$ on
$\iota^\#E$ coming from $h$ on $E$ we can map a tangent vector
$\dot{u}_0 \in \Secinfty_0(\iota^\#E)$ to a cotangent vector in
$\Sec[-\infty](\iota^\#E^*)$ by taking $\dot{u}^\flat \in
\Secinfty_0(\iota^\# E^*)$ and interpret this smooth section of
$\iota^\#E^*$ as a distributional section $\dot{u}_0^\flat \in
\Sec[-\infty](\iota^\#E^*)$. Clearly, this yields an injective linear
map
\begin{equation}
    \label{eq:flat-into-distributional-sections-is-injective}
    \Secinfty_0(\iota^\#E) \ni \dot{u}_0
    \; \mapsto \;
    \dot{u}_0^\flat \in
    \Secinfty_0(\iota^\#E^*) \subseteq \Sec[-\infty](\iota^\#E^*),
\end{equation}
which allows to pull back $\omega_{\mathrm{can}}$ to the tangent
bundle. This results in the following, still non-degenerate two-form:
\begin{lemma}
    \label{lemma:pullback-of-symplectic-form}
    The pull-back of the symplectic form $\omega_{\mathrm{can}}$ from
    the cotangent bundle of $\Secinfty_0(\iota^\#E)$ to its tangent
    bundle $\Secinfty_0(\iota^\#E) \oplus \Secinfty_0(\iota^\#E)$ via
    \eqref{eq:flat-into-distributional-sections-is-injective} is
    explicitly given by
    \begin{equation}
        \label{eq:pullback-of-symplectic-form}
        \omega_h \left(
            (u_0, \dot{u}_0), (v_0, \dot{v}_0)
        \right)
        = \int_\Sigma
        \left(
            h_\Sigma(u_0, \dot{v}_0) - h_\Sigma(\dot{u}_0,v_0)
        \right) \mu_\Sigma
    \end{equation}
    for $u_0, \dot{u}_0, v_0, \dot{v}_0 \in
    \Secinfty_0(\iota^\#E)$. The two-form $\omega_h$ turns
    $\Secinfty_0(\iota^\#E) \oplus \Secinfty_0(\iota^\#E)$ also into a
    symplectic vector space.
\end{lemma}
\begin{proof}
    Evaluating
    \eqref{eq:symplectic-form-on-Configspce-cotangent-with-density}
    for the distributional sections $\dot{u}_0^\flat, \dot{v}_0^\flat$
    gives immediately \eqref{eq:pullback-of-symplectic-form}. We have
    to check the non-degeneracy: but since $h_\Sigma$ is
    non-degenerate we can always find \emph{smooth} $(v_0, \dot{v}_0)$
    for a given $(u_0, \dot{u}_0) \neq 0$, resulting in a non-trivial
    pairing via $\omega_h$.
\end{proof}

\begin{remark}[Weak vs. strong symplectic]
    \label{remark:weak-and-strong-symplectic}%
    \label{Symplectic vector space!weak vs. strong}%
    The symplectic structure $\omega_{\mathrm{can}}$ on the cotangent
    bundle is even a \emph{strong} symplectic form if one defines the
    topological dual of $\Sec[-\infty](\iota^\#E^*)$ in an appropriate
    way: $\omega_{\mathrm{can}}$ induces an isomorphism from
    $\Secinfty_0(\iota^\#E) \oplus \Sec[-\infty](\iota^\#E^*)$ to its
    topological dual. For $\omega_h$ this is clearly not the case as
    the topological dual of $\Secinfty_0(\iota^\#E) \oplus
    \Secinfty_0(\iota^\#E)$ are two copies of
    $\Sec[-\infty](\iota^\#E^*)$ but with $\omega_h( (u_0, \dot{u}_0),
    \argument)$ we only obtain the (very small) part of smooth
    sections of $\iota^\#E^*$ and not the generalized ones. This is an
    effect of infinite dimension as in finite dimensions an injective
    linear map from a vector space to its dual is necessarily
    bijective. This indicates that for a Hamiltonian description one
    has to expect some (bad!) surprises.
\end{remark}
In any case, we only want to use the symplectic form to define the
Poisson algebra of observables of our ``mechanical'' system. In the
most general approach this algebra consists of smooth functions on the
(co-) tangent bundles. However, we do not want to enter the quite
nontrivial discussion on the appropriate definition of smooth
functions on the LF space $\Secinfty_0(\iota^\#E) \oplus
\Secinfty_0(\iota^\#E)$. There are several competing options which we
do not discuss here. To get a flavour of the complications one should
consult e.g. \cite{kriegl.michor:1997a}. Instead, we focus only on a
very small class of functions, the polynomials on the tangent bundle.

%
%

\subsection{The Poisson Algebra of Polynomials}
\label{satz:poisson-algebra-polynomials}

If $V$ is a Hausdorff locally convex topological vector space over
$\mathbb{R}$, what should the polynomials on $V$ be? Clearly, a
homogeneous polynomial of degree $1$ is just a \emph{linear}
functional $V \longrightarrow \mathbb{R}$ and hence an element of the
dual space of $V$. Having a topological vector space $V$ we require
\emph{continuity} for the homogeneous polynomials of degree $1$ whence
we end up with an element of $V'$.

Passing to homogeneous quadratic polynomials we certainly like to have
expression as
\[
p(v) = \sum_{i=1}^N \varphi_i(v) \psi_i(v)
\]
with $\varphi_i, \psi_i \in V'$ to be part of our observables. Indeed,
if we insist on an \emph{algebra} this is even forced by the algebraic
features: such a $p: V \longrightarrow \mathbb{R}$ is the sum of
products of elements in $V'$. Since we can multiply further we also
have to include functions of the form
\begin{equation}
    \label{eq:homogeneous-polynomial-degree-k}
    p(v) = \sum_{i=1}^N \varphi_1^{(i)}(v) \cdots \varphi_k^{(i)}(v)
\end{equation}
with $\varphi_1^{(i)}, \ldots, \varphi_k^{(i)} \in V'$ and $v \in
V$. Such a function certainly deserves the name ``homogeneous
polynomial of degree $k$''. Taking also linear combinations of such
polynomials of different homogeneity, which is again required if we
want an \emph{algebra} of observables, we end up with functions $p: V
\longrightarrow \mathbb{C}$ of the form
\begin{equation}
    \label{eq:polynomial}
    p(v) = c + \sum_{k=1}^\ell \sum_{i=1}^{N_k}
    \varphi_{1,k}^{(i)}(v) \cdots \varphi_{k,k}^{(i)}(v)
\end{equation}
with $\varphi_{\chi,k}^{(i)} \in V'$ and $v \in V$ and a constant $c
\in \mathbb{C}$.
\begin{definition}[Polynomial functions]
    \label{definition:polynomials}%
    \index{Polynomial function}%
    Let $V$ be a Hausdorff locally convex topological vector
    space. Then the polynomial functions generated by the constants
    and the linear functions $\varphi \in V'$ are denoted by
    $\Pol^\bullet (V)$.
\end{definition}
These functions can be identified with the symmetric algebra over
$V'$.
\begin{proposition}
    \label{proposition:polynomials-and-symmetric-algebra-over-V'}
    Let $V$ be a Hausdorff locally convex topological vector
    space. Then the polynomial functions $p: V \longrightarrow
    \mathbb{R}$ of the form \eqref{eq:polynomial} are in canonical
    bijection with the symmetric algebra $\Sym^\bullet V'$ over
    $V'$. The isomorphism is explicitly given by
    \begin{equation}
        \label{eq:isomorphism-pol-symm-algebra-over-V'}
        \mathcal{J}: \Sym^\bullet V' \ni
        \varphi_1 \vee \cdots \vee \varphi_k
        \; \mapsto \;
        \mathcal{J}(\varphi_1 \vee \cdots \vee \varphi_k)
        =
        \mathcal{J}(\varphi_1) \cdots \mathcal{J}(\varphi_k)
        \in \Pol^\bullet(V),
    \end{equation}
    where for degree $0$ and $1$ we have explicitly
    \begin{equation}
        \label{eq:explicit-form-isomorphism-pol-sym-alg}
        \mathcal{J}(\varphi)(v) = \varphi(v)
        \quad \textrm{and} \quad
        \mathcal{J}(\mathbb{1})(v) = 1.
    \end{equation}
    On arbitrary homogeneous elements $\Phi \in \Sym^k V'$ we have
    \begin{equation}
        \label{eq:isomorphism-on-homogeneous-elements}
        \mathcal{J}(\Phi)(v) = \frac{1}{k!} \Phi(v, \ldots, v).
    \end{equation}
\end{proposition}
\begin{proof}
    This is abstract nonsense on the symmetric algebra. First we
    recall that $\Sym^k V'$ consists of linear combinations of totally
    symmetrized tensor products of $k$ elements $\varphi_1, \ldots,
    \varphi_k \in V'$. We adopt the convection
    \[
    \varphi_1 \vee \cdots \vee \varphi_k
    = \sum_{\sigma \in S_k}
    \varphi_{\sigma(1)} \tensor \cdots \tensor \varphi_{\sigma(k)}
    \]
    \emph{without} prefactors. Then it is well-known that
    $\Sym^\bullet V'$ with $\vee$ is \emph{the} (up to canonical
    isomorphisms) free commutative algebra generated by $\mathbb{1}$
    and $V'$. Since the polynomials \eqref{eq:polynomial} are, by
    construction, also generated by $V'$ and the constants, we get a
    unique algebra homomorphism $\mathcal{J}$ by specifying it on the
    generators by
    \eqref{eq:explicit-form-isomorphism-pol-sym-alg}. Evaluating this
    on higher tensor products gives immediately
    \eqref{eq:isomorphism-on-homogeneous-elements} \emph{with}
    prefactor. It remains to show that $\mathcal{J}$ is injective,
    since the surjectivity is clearly the definition of the
    polynomials. Thus assume that $\Phi = \sum_{k=1}^\ell \Phi_k \in
    \Sym^\bullet V'$ with homogeneous components $\Phi_k \in \Sym^k
    V'$ satisfies $\mathcal{J}(\Phi) = 0$. Then for all $v \in V$ we
    have $\sum_{k=1}^\ell \frac{1}{k!} \Phi_k(v,\ldots,v) =
    0$. Rescaling $v$ to $tv$ with $t \in \mathbb{R}$ we see that the
    polynomial
    \[
    p(t) = \sum_{k=1}^\ell \frac{1}{k!} \Phi_k(v, \ldots, v) t^k
    = 0
    \]
    vanishes identically. Hence $\Phi_k(v,\ldots,v) = 0$ for all $k$
    separately. Now the polarization identities allow to express
    $\Phi_k(v_1, \ldots, v_k)$ in terms of linear combinations of
    terms $\Phi_k(w,\ldots,w)$ with $w$ being certain linear
    combinations of the $v_1, \ldots, v_k$. E.g. for quadratic ones we
    have
    \[
    \Phi_2(v_1,v_2)
    = \frac{1}{2} \left(
        \Phi_2(v_1+v_2,v_1+v_2) - \Phi_2(v_1,v_1) + \Phi_2(v_2,v_2)
    \right)
    \]
    and so on. But then $\Phi_k(v,\ldots,v) = 0$ for all $v \in V$
    implies $\Phi_k = 0$ in $\Sym^k V'$. Thus $\mathcal{J}$ is
    injective.
\end{proof}

\begin{remark}[Polynomial functions]
    \label{remark:polynomials}
    Let again $V$ be a Hausdorff locally convex vector space.
    \begin{remarklist}
    \item \label{item:J-isomorphism} From
        Proposition~\ref{proposition:polynomials-and-symmetric-algebra-over-V'}
        we have that
        \begin{equation}
            \label{eq:J-isomorphism}
            \mathcal{J}: \Sym^\bullet V' \longrightarrow
            \Pol^\bullet(V)
        \end{equation}
        is an isomorphism of commutative, unital, and graded
        algebras.
    \item \label{item:homogeneous-by-pulling-out-scalars} More
        generally, one could define a polynomial function $p: V
        \longrightarrow \mathbb{R}$ on $V$ of degree $k$ to be a
        function with the property
        \begin{equation}
            \label{eq:homogeneous-by-pulling-out-scalars}
            p(tv) = t^k p(v)
        \end{equation}
        for all $v \in V$ and $t \in \mathbb{R}$ \emph{plus} some
        suitable continuity at the origin. This continuity is already
        needed in finite dimensions to exclude functions like
        \begin{equation}
            \label{eq:hom-pol-pull-out-counterex-finite-dim}
            p(v) =
            \begin{cases}
                0 & \textrm{ for } v=0 \\
                \frac{v^iv^jv^k}{\sum_{\ell=1}^{\dim V} (v^\ell)^2} &
                \textrm{ for } v \neq 0
            \end{cases}
        \end{equation}
        to be a ``linear polynomial'' in $v$. Here $v = v^i e_i$
        with a basis $e_i \in V$.
    \item \label{item:completions-of-polynomial-algebras} Since $V'$
        carries a natural Hausdorff locally convex topology, the
        weak$^*$ topology, one can endow $\Sym^k V'$ with a locally
        convex topology as well: in fact, there are several and
        typically inequivalent possibilities. The usage of such
        topologies can be two-fold: on one hand we can complete each
        $\Sym^k V'$ which amounts to obtaining polynomial functions of
        homogeneous degree $k$ of the form
        \begin{equation}
            \label{eq:infinite-linear-combinations-of-hom-polynomials}
            p(v) = \sum_{i=1}^\infty \varphi_1^{(i)}(v) \cdots
            \varphi_k^{(i)}(v),
        \end{equation}
        where the topology on $\Sym^k V'$ is now used to make sense
        out of the limit. But we can also complete into another
        direction: the direct sum $\Sym^\bullet V' =
        \bigoplus_{k=0}^\infty \Sym^k V'$ can be completed to include
        also ``transcendental'' functions and not just polynomials. In
        particular, one would be interested in functions as $f(v) =
        \E^{\varphi(v)}$ with $\varphi \in V'$. This leads to notions
        of holomorphic or real analytic functions on $V$. While the
        first completion does not give anything new in finite
        dimension the second is already interesting in finite
        dimensions. If $V$ is infinite-dimensional, both types of
        completions are typically non-trivial and depend on the
        precise choices of the topologies on the (symmetric) tensor
        products.
    \end{remarklist}
\end{remark}

After these general considerations we come back to our original task:
on the symplectic vector space $\Secinfty_0(\iota^\# E) \oplus
\Secinfty_0(\iota^\#E)$ we want to establish a polynomial algebra with a
\emph{Poisson bracket}.

So the first guess is to use the symmetric algebra over
$\Sec[-\infty](\iota^\# E^*) \oplus \Sec[-\infty](\iota^\#E^*)$, which
is the topological dual of $\Secinfty_0(\iota^\#E) \oplus
\Secinfty_0(\iota^\#E)$ via the usual identification, and endow this
symmetric algebra with a Poisson bracket. The problem is here the
following: Since $\omega_h$ is only a weak symplectic form, not every
linear functional has a Hamiltonian vector field. Thus the Poisson
bracket can not be defined that easily on all linear functionals and
hence on all polynomials. This forces us to proceed differently: we
take as a beginning the subspace
\begin{equation}
    \label{eq:simple-subspace-for-us}
    \Secinfty_0(\iota^\# E^*) \oplus \Secinfty_0(\iota^\#E^*)
    \subseteq
    \Sec[-\infty](\iota^\#E^*) \oplus \Sec[-\infty](\iota^\#E^*)
\end{equation}
as dual space of $\Secinfty_0(\iota^\#E) \oplus
\Secinfty_0(\iota^\#E)$ and consider the symmetric algebra over this
much smaller space. Here the following result is easy to obtain:
\begin{proposition}
    \label{proposition:poisson-structure-on-easy-polynomials}%
    \index{Poisson bracket}%
    \index{Poisson algebra}%
    \index{Hamiltonian vector field}%
    On the symmetric algebra over $\Secinfty_0(i^\#E^*) \oplus
    \Secinfty_0(i^\#E^*)$ exists a unique Poisson bracket
    $\PB{\argument, \argument}_h$ induced by $\omega_h$ with the property
    \begin{equation}
        \label{eq:poisson-bracket-of-linear-polynomials}
        \PB{ \mathcal{J}(\varphi_0, \dot{\varphi}_0),
          \mathcal{J}(\psi_0, \dot{\psi}_0)}_h
        = \int_\Sigma
        \left(
            h_\Sigma^{-1}(\varphi_0, \dot{\psi}_0)
            - h_\Sigma^{-1}(\dot{\varphi}_0, \psi_0)
        \right) \mu_\Sigma
    \end{equation}
    for $(\varphi_0, \dot{\varphi}_0), (\psi_0, \dot{\psi}_0) \in
    \Secinfty_0(\iota^\#E^*) \oplus \Secinfty_0(\iota^\#E^*)$. The
    Hamiltonian vector field of the linear functional
    $\mathcal{J}(\varphi_0, \dot{\varphi}_0)$ with respect to
    $\omega_h$ is the constant vector field
    \begin{equation}
        \label{eq:hamiltonian-vector-field-of-linear-polynomial}
        X_{\mathcal{J}(\varphi_0, \dot{\varphi}_0)}
        = ( - \dot{\varphi}_0^\#, \varphi_0^\# ).
    \end{equation}
\end{proposition}
\begin{proof}
    First we note that any Poisson bracket on a symmetric algebra
    $\Sym^\bullet W$ of any vector space $W$ is uniquely determined by
    its values on $W$ alone: since a Poisson bracket satisfies by
    definition a Leibniz rule in both arguments it is determined by
    its values on a set of generators of the algebra. Since
    necessarily $\PB{\mathbb{1}, \argument} = 0 = \PB{\argument,
      \mathbb{1}}$ for any Poisson bracket it is therefore sufficient
    to specify it on the generators $W \subseteq \Sym^\bullet W$. Thus
    $\PB{\argument, \argument}_h$ will be uniquely determined by
    \eqref{eq:poisson-bracket-of-linear-polynomials}. To motivate the
    formula \eqref{eq:poisson-bracket-of-linear-polynomials} we first
    prove
    \eqref{eq:hamiltonian-vector-field-of-linear-polynomial}. Thus let
    $(\varphi_0, \dot{\varphi}_0)$ be given. Since this is viewed as a
    linear function the \emph{differential} is constant and given by
    $(\varphi_0, \dot{\varphi}_0)$ at every point, i.e.
    \[
    \D \mathcal{J}(\varphi_0, \dot{\varphi}_0) \at{(u_0, \dot{u}_0)}
    = (\varphi_0, \dot{\varphi}_0).
    \tag{$*$}
    \]
    Thus the Hamiltonian vector field, defined by $\omega_h(X_f,
    \argument) = \D f(\argument)$ in general, is determined by
    \begin{align*}
        \int_\Sigma
        \left(
            \varphi_0(v_0) + \dot{\varphi}_0(\dot{v}_0)
        \right) \mu_\Sigma
        &= \D\! \mathcal{J}(\varphi_0, \dot{\varphi}_0)
        \at{(u_0, \dot{u}_0)} (v_0, \dot{v}_0) \\
        &= \omega_h \left(
            X_{\mathcal{J}(\varphi_0, \dot{\varphi}_0)}
            \at{(u_0, \dot{u}_0)}, (v_0, \dot{v}_0)
        \right) \\
        &= \int_\Sigma
        \left(
            h_\Sigma\left(
                X_{\mathcal{J}(\varphi_0, \dot{\varphi}_0)}
                \at{(u_0,\dot{u}_0)},
                \dot{v}_0
            \right)
            -
            h_\Sigma\left(
                \dot{X}_{\mathcal{J}(\varphi_0, \dot{\varphi}_0)}
                \at{(u_0, \dot{u}_0)},
                v_0
            \right)
        \right)
        \mu_\Sigma.
    \end{align*}
    This shows that $X_{\mathcal{J}(\varphi_0,\dot{\varphi}_0)}$ is
    the constant vector field with the two components
    \[
    X_{\mathcal{J}(\varphi_0, \dot{\varphi}_0)} \at{(u_0,\dot{u}_0)}
    = ( - \dot{\varphi}_0^\#, \varphi_0^\# )
    \]
    at every point $(u_0, \dot{u}_0)$,
    i.e. \eqref{eq:hamiltonian-vector-field-of-linear-polynomial}. Now
    the Poisson bracket is, by definition $\PB{f,g} = X_g(f) = \D\!f
    (X_g)$. Hence we get the \emph{constant} function
    \begin{align*}
        \PB{\mathcal{J}(\varphi_0, \dot{\varphi}_0),
          \mathcal{J}(\psi_0, \dot{\psi}_0)} \At{(u_0, \dot{u}_0)}
        &= (\varphi_0, \dot{\varphi}_0) (\dot{\psi}_0^\#, - \psi_0^\#)
        \\
        &= \int_\Sigma
        \left(
            \varphi_0( \dot{\psi}_0^\# ) - \dot{\varphi}_0(\psi_0^\#)
        \right) \mu_\Sigma \\
        &= \int_\Sigma
        \left(
            h_\Sigma^{-1}(\varphi_0, \dot{\psi}_0)
            - h_\Sigma^{-1}(\dot{\varphi}_0, \psi_0)
        \right) \mu_\Sigma,
    \end{align*}
    using the dual fibre metric $h_\Sigma^{-1}$ on $\iota^\#E^*$. This
    explains the statement
    \eqref{eq:poisson-bracket-of-linear-polynomials}. For finite
    dimensional vector spaces (or manifolds) we could now argue with
    the usual calculus of smooth functions that, thanks to the
    closedness of $\omega_h$, the Poisson bracket is indeed a Poisson
    bracket. In infinite dimensions we can not just rely on the
    analogy, in particular since $\omega_h$ is only a weak symplectic
    structure. Instead of establishing an appropriate calculus also in
    this situation, which in principle can be done, we prove the
    \emph{existence} of a Poisson bracket on the polynomials by
    hand. In fact, this follows from the next proposition at once.
\end{proof}

\begin{proposition}
    \label{proposition:poisson-bracket-on-symmetric-algebra}
    Let $W$ be a real vector space and let
    \begin{equation}
        \label{eq:antisymm-bilinear-form-on-W}
        \pi: W \times W \longrightarrow \mathbb{R}
    \end{equation}
    be an antisymmetric bilinear form. Then on $\Sym^\bullet W$ there
    is a unique Poisson bracket $\PB{\argument, \argument}_\pi$ with
    \begin{equation}
        \label{eq:poisson-bracket-on-symmetric-algebra-grading}
        \PB{\argument, \argument}_\pi: \Sym^k W \times \Sym^\ell W
        \longrightarrow \Sym^{k+\ell-2} W,
    \end{equation}
    such that for $v,w \in W = \Sym^1W \subseteq \Sym^\bullet W$ one
    has
    \begin{equation}
        \label{eq:poisson-bracket-on-linear-elements}
        \PB{v,w}_\pi
        = \pi(v,w) \mathbb{1}.
    \end{equation}
\end{proposition}
\begin{proof}
    Again, the uniqueness is clear since by the Leibniz rule, a
    Poisson bracket is determined by its values on the
    generators. Enforcing the Leibniz rule gives us the explicit
    expression
    \[
    \PB{ v_1 \vee \cdots \vee v_k, w_1 \vee \cdots \vee w_\ell}_\pi
    = \sum_{i,j} \pi(v_i,w_j) v_1 \vee \cdots \stackrel{i}{\wedge}
    \cdots \vee v_k \vee w_1 \vee \cdots \stackrel{j}{\wedge} \cdots
    \vee w_\ell
    \]
    as the unique extension of $\pi$ to $\Sym^\bullet W$ which
    satisfies the Leibniz rule in both arguments. Since $\pi$ is
    antisymmetric, $\PB{\argument, \argument}_\pi$ is antisymmetric as
    well. It remains to check the \Index{Jacobi identity}. Thus let
    \[
    \mathrm{Jac}_\pi(f,g,h)
    = \PB{f,\PB{g,h}_\pi}_\pi
    + \PB{g,\PB{h,f}_\pi}_\pi
    + \PB{h,\PB{f,g}_\pi}_\pi
    \]
    be the \Index{Jacobiator} of $\PB{\argument, \argument}_\pi$ for
    arbitrary $f,g,h \in \Sym^\bullet W$. We have to show that
    $\mathrm{Jac}_\pi(f,g,h) = 0$. Now it is a simple algebraic fact
    that $\mathrm{Jac}_\pi$ is a derivation in each argument. Thus
    $\mathrm{Jac}_\pi(f,g,h) = 0$ iff the Jacobiator vanishes on
    \emph{generators} already. In our case $\mathrm{Jac}_\pi(v,w,u) =
    0$ is clear, since $\PB{v,w}_\pi$ is already constant. The grading
    statement \eqref{eq:poisson-bracket-on-symmetric-algebra-grading}
    is clear.
\end{proof}

This way we obtain a Poisson algebra of polynomials modeled by the
symmetric algebra over $\Secinfty_0(\iota^\#E^*) \oplus
\Secinfty_0(\iota^\#E^*)$. Without going into the details we note that
this Poisson bracket has reasonable continuity properties with respect
to the usual LF topology of $\Secinfty_0(\iota^\#E^*) \oplus
\Secinfty_0(\iota^\#E^*)$. To explain these properties we first
rewrite
\begin{equation}
    \label{eq:oplus-of-sec-is-sec-of-oplus}
    \Secinfty_0(\iota^\#E^*) \oplus \Secinfty_0(\iota^\#E^*)
    = \Secinfty_0(\iota^\#(E^* \oplus E^*))
\end{equation}
as usual. Then we have the following lemma:
\begin{lemma}
    \label{lemma:injection-of-symk-into-sec-of-k-times-extensor}
    There is a canonical injection
    \begin{equation}
        \label{eq:injection-of-symk-into-sec-of-k-times-extensor}
        \Sym^k \Secinfty_0(\iota^\#(E^* \oplus E^*)) \hookrightarrow
        \Secinfty_0(\iota^\#(E^* \oplus E^*)
        \underbrace{\extensor \cdots \extensor}_{k\textrm{-times}}
        \iota^\#(E^* \oplus E^*))^{\Sym_k}
    \end{equation}
    of the symmetric power of $\Secinfty_0(\iota^\#(E^* \oplus E^*))$
    of degree $k$ into the sections of the $k$-th external tensor
    product of $\iota^\#(E^* \oplus E^*)$ with itself which are
    totally symmetric under the internal action of the permutations of
    the fibres. Explicitly, we have
    \begin{equation}
        \label{eq:explicit-injection-of-symk-into-sec-of-k-times-extensor}
        (\varphi_1 \vee \cdots \vee \varphi_k)(p_1, \ldots, p_k)
        = \sum_{\sigma \in \Sym_k}
        \varphi_{\sigma(1)}(p_1) \extensor \cdots \extensor
        \varphi_{\varphi(k)}(p_k)
    \end{equation}
    for $p_1, \ldots, p_k \in \Sigma$ and $\Secinfty_0(\iota^\#(E^*
    \oplus E^*))$.
\end{lemma}
\begin{proof}
    Clearly,
    \eqref{eq:explicit-injection-of-symk-into-sec-of-k-times-extensor}
    is injective and well-defined, yielding a totally symmetric
    section with compact support. This follows analogously to
    Theorem~\ref{theorem:sections-on-external-tensor-product}.
\end{proof}

\begin{remark}
    \label{remark:extension-of-poisson-bracket}
    As in Theorem~\ref{theorem:sections-on-external-tensor-product}
    this map is continuous in a very precise way: we have estimates
    analogously to the ones in
    \eqref{eq:cont-estimate-for-tensor-product-of-Ck-functions}. Without
    introducing this notion, we note that
    \eqref{eq:injection-of-symk-into-sec-of-k-times-extensor} is
    continuous with respect to the \emph{projective} tensor product
    topology of $\Sym^k \Secinfty_0(\iota^\#(E^* \oplus E^*))$, see
    e.g. \cite[Chap.~15]{jarchow:1981a} for more details on this
    $\pi$-topology. Moreover, we note that the image of
    \eqref{eq:injection-of-symk-into-sec-of-k-times-extensor} is
    sequentially dense in the totally symmetric sections. This can
    also be shown analogously to
    Theorem~\ref{theorem:sections-on-external-tensor-product}. In
    fact, this even allows to extend the Poisson bracket
    $\PB{\argument, \argument}_h$ to the direct sum over the right
    hand side of
    \eqref{eq:injection-of-symk-into-sec-of-k-times-extensor} for $k
    \in \mathbb{N}_0$ by a continuity argument. However, we shall not
    enter this discussion here.
\end{remark}

From now on, we shall omit the explicit usage of the symbol
$\mathcal{J}$ in \eqref{eq:J-isomorphism} to simplify our notation and
identify elements in $\Sym^\bullet V$ with the polynomials in
$\Pol^\bullet (V)$ directly.

%
%

\subsection{The Covariant Poisson Algebra}
\label{satz:covariant-poisson-algebra}

Up to now the Poisson algebra of observables has certain deficits from
a physical point of view: its definition depends on the \emph{choice}
of a Cauchy hypersurface. In particular, it is not quite clear whether
we get different Poisson algebras for different choices and, if not,
how they are related in detail. In fact, since on a globally
hyperbolic spacetime $M$ all smooth spacelike Cauchy hypersurfaces are
diffeomorphic and since any two positive definite fibre metrics are
isometric, one can cook up an isomorphism of the Poisson algebras
corresponding to $(\Sigma_1, h_{\Sigma_1})$ and $(\Sigma_2,
h_{\Sigma_2})$, respectively. However, this does not seem to be a very
conceptual statement as the isomorphism is just there by ``pure
luck''.

More severe than these aesthetic arguments is the conceptual
disadvantage that all nice symmetries between time- and spacelike
directions will be ``broken'' by the choice of $\Sigma$. As example,
one considers again Minkowski spacetime $(\mathbb{R}^n,\eta)$ with its
Poincare symmetry $O(1,n-1) \ltimes \mathbb{R}^n$. Choosing an
arbitrary smooth spacelike Cauchy hypersurface $\Sigma$ results in
destroying the symmetry: the Poincare group action will \emph{not}
respect the splitting $\mathbb{R}^n \simeq \mathbb{R} \times
\Sigma$, even if $\Sigma$ is a spacelike linear subspace. Thus the
true symmetry of the situation might be hidden after choosing a
splitting $\mathbb{R} \times \Sigma$.

Thus we look for a Poisson algebra isomorphic to the one constructed
in Proposition~\ref{proposition:poisson-structure-on-easy-polynomials}
which is intrinsically defined without reference to $\Sigma$. This
will be accomplished by the following construction, essentially going
back to Peierls \cite{peierls:1952a}, see also
\cite{duetsch.fredenhagen:2001b, duetsch.fredenhagen:2001a,
  duetsch.fredenhagen:2003a, marolf:1994a, marolf:1994b} for a more
modern treatment and applications to the (deformation) quantization of
classical field theories as well as the thesis
\cite{jakobs:2009a}. Note however, that we are only dealing with
rather simple polynomial functions here instead of more general smooth
functions.

We consider $\Secinfty_0(E^*)$ which we can use to evaluate arbitrary
sections $u \in \Secinfty_{\mathrm{sc}}(E)$ on the whole spacetime
$M$. Again, the symmetric algebra $\Sym^\bullet \Secinfty_0(E^*)$
serves as polynomial algebra on \emph{all fields}
$\Secinfty_{\mathrm{sc}}(E)$, whether they are solutions to $Du=0$ or
not. The evaluation is the normal one, i.e. for $\varphi \in
\Secinfty_0(E^*)$ we set
\begin{equation}
    \label{eq:evaluation-with-smooth-section-of-dual-bundle}
    \varphi(u)
    = \int_M \varphi(p) \cdot u(p) \: \mu_g(p)
\end{equation}
and extend this to $\Sym^\bullet \Secinfty_0(E)$ as before. Then these
symmetric tensors become again an observable algebra. However, it
should be emphasized clearly that we are dealing with polynomials on a
much too large space $\Secinfty_{\mathrm{sc}}(E)$ at the
moment. Surprisingly, we will even have a Poisson bracket on this too
large algebra:
\begin{proposition}
    \label{proposition:poisson-algebra-on-polynomials-on-sc-sections}
    Let $(M,g)$ be a globally hyperbolic spacetime and $D \in
    \Diffop^2(E)$ a normally hyperbolic differential operator that is
    symmetric with respect to a fibre metric $h$ on $E$. Then on the
    symmetric algebra $\Sym^\bullet \Secinfty_0(E^*)$ there is a
    unique Poisson bracket $\PB{\argument, \argument}$ determined by
    \begin{equation}
        \label{eq:poisson-bracket-on-generators}
        \PB{ \varphi, \psi }
        = \int_M h^{-1}\left(F_M \varphi, \psi\right) \mu_g
    \end{equation}
    for $\varphi, \psi \in \Secinfty_0(E^*)$, where $F_M = F_M^+ -
    F_M^-$ as before. It satisfies
    \begin{equation}
        \label{eq:grading-of-poisson-bracket}
        \PB{ \Sym^k \Secinfty_0(E^*), \Sym^\ell \Secinfty_0(E^*) }
        \subseteq \Sym^{k+\ell-2} \Secinfty_0(E^*).
    \end{equation}
\end{proposition}
\begin{proof}
    Since $F_M$ is an antisymmetric operator with respect to the
    integration and $h$ according to
    Remark~\ref{remark:propagator-and-adjoining}, see also
    \eqref{eq:propagator-is-antisymmetric}, the right hand side of
    \eqref{eq:poisson-bracket-on-generators} defines an antisymmetric
    bilinear form on $\Secinfty_0(E^*)$. Thus,
    Proposition~\ref{proposition:poisson-bracket-on-symmetric-algebra}
    can be applied.
\end{proof}

\begin{definition}[Covariant Poisson bracket]
    \label{definition:covariant-poisson-bracket}%
    \index{Covariant Poisson bracket}%
    \index{Poisson bracket!covariant}%
    The Poisson bracket on $\Sym^\bullet \Secinfty_0(E^*)$ resulting
    from \eqref{eq:poisson-bracket-on-generators} is called the
    covariant Poisson bracket corresponding to $D$.
\end{definition}
Even though $\Sym^\bullet \Secinfty_0(E^*)$ is enough to separate
points on the too large space of all fields
$\Secinfty_{\mathrm{sc}}(E)$, the covariant Poisson bracket becomes
trivial for elements not sensitive to solutions of the wave
equation. More precisely, we have the following result:
\begin{lemma}
    \label{lemma:degeneration-space-of-poisson-bracket}%
    \index{Casimir element}%
    Let $\varphi \in \Secinfty_0(E^*)$. Then the following statements
    are equivalent:
    \begin{lemmalist}
    \item \label{item:phi-is-casimir} $\varphi$ is a Casimir element
        of the covariant Poisson algebra $( \Sym^\bullet
        \Secinfty_0(E^*), \PB{\argument, \argument})$, i.e. we have
        \begin{equation}
            \label{eq:casimir-element}
            \PB{ \varphi, \argument} = 0.
        \end{equation}
    \item \label{item:phi-vanishes-on-solutions} $\varphi$ vanishes on
        solutions $u \in \Secinfty_{\mathrm{sc}}(E)$ of the wave
        equation $Du = 0$, i.e.
        \begin{equation}
            \label{eq:phi-vanishes-on-solutions}
            \int_M \varphi \cdot u \: \mu_g = 0.
        \end{equation}
    \item \label{item:phi-in-kernel-of-propagator} $\varphi$ is in the
        kernel of $F_M$, i.e.
        \begin{equation}
            \label{eq:phi-in-kernel-of-propagator}
            F_M \varphi = 0.
        \end{equation}
    \end{lemmalist}
\end{lemma}
\begin{proof}
    We show \refitem{item:phi-is-casimir} $\Rightarrow$
    \refitem{item:phi-in-kernel-of-propagator} $\Rightarrow$
    \refitem{item:phi-vanishes-on-solutions} $\Rightarrow$
    \refitem{item:phi-is-casimir}. Assume $\PB{\varphi, \argument} =
    0$, then $0 = \PB{\varphi,\psi} = \int_M h^{-1}(F_M\varphi, \psi)
    \mu_g$ for all $\psi \in \Secinfty_0(E^*)$ which implies $F_M
    \varphi = 0$ since the pairing is non-degenerate. Now, if $F_M
    \varphi = 0$ then by Theorem~\ref{theorem:propagator-complex},
    \refitem{item:propagator-complex-exact-on-globally-hyperbolic}
    applied to $D^\Trans$ we know $\varphi = D^\Trans \chi$ for some
    $\chi \in \Secinfty_0(E^*)$. Thus
    \[
    \int_M \varphi \cdot u \: \mu_g
    = \int_M D^\Trans \chi \cdot u \: \mu_g
    = \int_M \varphi \cdot Du \: \mu_g = 0
    \]
    for any solution $u \in \Secinfty_{\mathrm{sc}}(E)$ of the wave
    equation. Finally, assume that
    \refitem{item:phi-vanishes-on-solutions} holds and let $\psi \in
    \Secinfty_0(E^*)$ be arbitrary. Then $(F_M \psi)^\# = G_M \psi^\#$
    solves the homogeneous wave equation. Thus
    \[
    0 = \int_M \varphi \cdot (F_M \psi)^\# \mu_g
    = \int_M h^{-1}( \varphi, F_M \psi) \mu_g
    = - \PB{\varphi, \psi}
    \]
    for all $\psi \in \Secinfty_0(E^*)$. By the Leibniz rule this
    implies $\PB{\varphi, \argument} = 0$ in general, since these
    $\psi$ generate the whole algebra.
\end{proof}

We can rephrase the result of the lemma as follows: the kernel of
$F_M$ is a subspace $\ker F_M \subseteq \Secinfty_0(E^*)$ which
generates an ideal inside $\Sym^\bullet \Secinfty_0(E^*)$. The
generators of this ideal are Casimir elements whence the ideal is in
fact even a \emIndex{Poisson ideal}. Thus the quotient algebra of
$\Sym^\bullet \Secinfty_0(E^*)$ by this ideal becomes a Poisson
algebra itself. Now we want to relate this quotient to the canonical
Poisson algebra defined on a Cauchy hypersurface $\Sigma$ as
constructed in the previous subsection. We want to establish a Poisson
isomorphism which is compatible with the evaluation on solutions of
the wave equation. To make these things more precise we again consider
the result from
Theorem~\ref{theorem:global-integral-formula-for-solution}. If $u \in
\Secinfty_{\mathrm{sc}}(E)$ is the unique solution of $Du = 0$ with
initial conditions $u_0, \dot{u}_0$ on $\Sigma$ then the evaluation of
$\varphi \in \Secinfty_0(E^*)$ on $u$ can be expressed by
\begin{equation}
    \label{eq:evalution-and-the-famous-integral-formula}
    \int_M \varphi \cdot u \: \mu_g
    = \int_\Sigma
    \left(
        (i^\# \nabla^{E^*}_{\mathfrak{n}} F_M \varphi) \cdot u_0
        - (i^\# F_M \varphi) \cdot \dot{u}_0
    \right) \mu_\Sigma,
\end{equation}
according to
Theorem~\ref{theorem:global-integral-formula-for-solution}. Comparing
this with the evaluation of a section $(\varphi_0, \dot{\varphi}_0)
\in \Secinfty_0(\iota^\# (E^* \oplus E^*))$ on initial conditions
according to \eqref{eq:poisson-bracket-of-linear-polynomials}, i.e.
\begin{equation}
    \label{eq:evaluation-of-section-on-initial-conditions}
    (\varphi_0, \dot{\varphi}_0) \at{(u_0,\dot{u}_0)}
    = \int_\Sigma
    (\varphi_0 \cdot u_0 + \dot{\varphi}_0 \cdot \dot{u}_0)
    \: \mu_\Sigma,
\end{equation}
suggests to map $\varphi \in \Secinfty_0(E^*)$ to the section
$(\varphi_0, \dot{\varphi}_0) \in \Secinfty_0(\iota^\# (E^* \oplus
E^*))$ given by
\begin{equation}
    \label{eq:map-for-identification-of-poisson-algebras}
    \varphi_0 = \iota^\# \nabla^{E^*}_{\mathfrak{n}} F_M \varphi
    \quad \textrm{and} \quad
    \dot{\varphi}_0 = - \iota^\# F_M \varphi.
\end{equation}
We denote this ``restriction map'' by
\begin{equation}
    \label{eq:restriction-map-to-construct-poisson-iso}
    \varrho_\Sigma: \Secinfty_0(E^*) \ni \varphi
    \; \mapsto \;
    \left(
        \iota^\# \nabla^{E^*}_{\mathfrak{n}} F_M \varphi,
        - \iota^\# F_M \varphi
    \right)
    \in \Secinfty_0(\iota^\#(E^* \oplus E^*)).
\end{equation}
Since $\Sym^\bullet \Secinfty_0(E^*)$ is freely generated by
$\Secinfty_0(E^*)$ we can extend $\varrho_\Sigma$ in a unique way to a
unital algebra homomorphism to $\Sym^\bullet \Secinfty_0(\iota^\#(E^*
\oplus E^*))$ which we still denote by
\begin{equation}
    \label{eq:building-the-algebra-homomorphism}
    \varrho_\Sigma: \Sym^\bullet \Secinfty_0(E^*)
    \longrightarrow
    \Sym^\bullet \Secinfty_0(\iota^\# (E^* \oplus E^*)).
\end{equation}
Then the above discussion results in the following lemma:
\begin{lemma}
    \label{lemma:the-algebra-homomorhpism}
    Let $u \in \Secinfty_{\mathrm{sc}}(E)$ be a solution of the
    homogeneous wave equation with initial conditions $u_0, \dot{u}_0
    \in \Secinfty_0(\iota^\#E)$ on $\Sigma$. Then for every $\Phi \in
    \Sym^\bullet \Secinfty_0(E^*)$ we have
    \begin{equation}
        \label{eq:the-algebra-homomorphism}
        \Phi(u) = \varrho_\Sigma(\Phi)(u_0, \dot{u}_0).
    \end{equation}
\end{lemma}
\begin{proof}
    We know \eqref{eq:the-algebra-homomorphism} for $\Phi = \varphi
    \in \Secinfty_0(E^*)$ by construction. For the constants we have
    by definition $\varrho_\Sigma(\mathbb{1}) = \mathbb{1}$ whence
    \eqref{eq:the-algebra-homomorphism} is also true here. For higher
    symmetric tensors $\Phi \in \Sym^\bullet \Secinfty_0(E^*)$ the
    evaluation on $u$ was defined to be compatible with the
    $\vee$-product, i.e.
    \[
    (\varphi_1 \vee \cdots \vee \varphi_k)(u)
    = \varphi_1(u) \cdots \varphi_k(u).
    \]
    Since we used the same sort of evaluation also for the symmetric
    tensors in $\Sym^\bullet \Secinfty_0(\iota^\# (E^* \oplus E^*))$
    the statement follows from the algebra homomorphism property of
    $\varrho_\Sigma$.
\end{proof}

Since the initial conditions determine the solution uniquely and vice
versa it is tempting to use the algebra homomorphism $\varrho_\Sigma$
to relate the Poisson algebras on $M$ and on $\Sigma$. Indeed, we have
the following result:
\begin{lemma}
    \label{lemma:alg-hom-is-poisson-hom}%
    \index{Poisson algebra!homomorphism}%
    The algebra homomorphism $\varrho_\Sigma$ is a homomorphism of
    Poisson algebras
    \begin{equation}
        \label{eq:alg-hom-is-poisson-hom}
        \varrho_\Sigma:
        (\Sym^\bullet \Secinfty_0(E^*), \PB{\argument, \argument})
        \longrightarrow
        (\Sym^\bullet \Secinfty_0(\iota^\# (E^* \oplus E^*)),
        \PB{\argument, \argument}_h).
    \end{equation}
\end{lemma}
\begin{proof}
    Since the Poisson brackets satisfy a Leibniz rule by definition
    and since $\varrho_\Sigma$ is a unital algebra homomorphism it
    suffices to check the claim on generators. Thus let $\varphi, \psi
    \in \Secinfty_0(E^*)$ be given and let $(\varphi_0,
    \dot{\varphi}_0) = \varrho_\Sigma(\varphi)$ and
    $(\psi_0,\dot{\psi}_0) = \varrho_\Sigma(\psi)$ be the
    corresponding sections in $\Secinfty_0(\iota^\#(E^* \oplus
    E^*))$. Moreover, both Poisson brackets
    $\{(\varphi_0,\dot{\varphi}_0),(\psi_0, \dot{\psi}_0)\}_h$ and
    $\PB{\varphi, \psi}$ are constants, i.e. multiples of the unit
    elements, respectively. Thus we only have to compute these number
    as $\varrho_\Sigma(\mathbb{1}) = \mathbb{1}$ by definition. We
    have
    \begin{align*}
        \PB{(\varphi_0, \dot{\varphi}_0), (\psi_0, \dot{\psi}_0)}
        &= \int_\Sigma
        \left(
            h_\Sigma^{-1}(\varphi_0,\dot{\psi}_0)
            - h_\Sigma^{-1}(\dot{\varphi}_0,\psi_0)
        \right) \mu_\Sigma \\
        &= - \int_\Sigma
        \left(
            (\iota^\# \nabla^{E^*}_{\mathfrak{n}} F_M \varphi) \cdot
            (\iota^\# F_M \psi)^\#
            - (\iota^\# F_M \varphi) \cdot
            (\iota^\# \nabla^{E^*}_{\mathfrak{n}} F_M \psi)^\#
        \right) \mu_\Sigma.
        \tag{$*$}
    \end{align*}
    Now $F_M \psi$ is a solution of the wave equation, $D^\Trans F_M
    \psi = 0$. Since $D$ is symmetric, $u = (F_M \psi)^\# = G_M
    \psi^\#$ is a solution of $Du = 0$. The initial conditions for $u$
    on $\Sigma$ are given by
    \[
    u_0 = \iota^\# u = \iota^\# (F_M \psi)^\#
    \quad
    \textrm{and}
    \quad
    \dot{u}_0 = \iota^\# \nabla^E_{\mathfrak{n}} u
    = \iota^\# (\nabla^{E^*}_{\mathfrak{n}} F_M \psi)^\#
    \]
    since the connections $\nabla^E$ and $\nabla^{E^*}$ are compatible
    with the musical isomorphisms as $\nabla^E$ is assumed to be
    metric with respect to $h$. By
    Theorem~\ref{theorem:global-integral-formula-for-solution} we
    conclude
    \begin{align*}
        \PB{(\varphi_0, \dot{\varphi}_0), (\psi_0, \dot{\psi}_0)}
        & \stackrel{(*)}{=} - \int_M \varphi \cdot u \: \mu_g \\
        &= - \int_M \varphi \cdot (F_M \psi)^\# \mu_g \\
        &= - \int_M h^{-1}(\varphi, F_M \psi) \: \mu_g \\
        &= \int_M h^{-1}(F_M \varphi, \psi) \: \mu_g \\
        &= \PB{\varphi, \psi}.
    \end{align*}
    This shows that the constants coincide and thus the claim follows.
\end{proof}

\begin{lemma}
    \label{lemma:kernel-of-poisson-hom}
    The Poisson homomorphism $\varrho_\Sigma$ is surjective and its
    kernel coincides with the ideal generated by the Casimir elements
    in $\Secinfty_0(E^*)$, which coincides with all those $\Phi \in
    \Sym^\bullet \Secinfty_0(E^*)$ which vanish on all solutions $u
    \in \Secinfty_{\mathrm{sc}}(E)$ of the homogeneous wave equation
    $Du = 0$.
\end{lemma}
\begin{proof}
    By definition we have $\varrho_\Sigma(\mathbb{1}) =
    \mathbb{1}$. Now let $(\varphi_0, \dot{\varphi}_0) \in
    \Secinfty_0(\iota^\#(E^* \oplus E^*))$ be given. Then there is a
    unique solution $\Phi \in \Secinfty_{\mathrm{sc}}(E^*)$ of the
    homogeneous wave equation $D^\Trans \Phi = 0$ with initial
    conditions
    \[
    \iota^\# \Phi = - \dot{\varphi}_0
    \quad \textrm{and} \quad
    \iota^\# \nabla^{E^*}_{\mathfrak{n}} \Phi = \varphi_0
    \tag{$*$}
    \]
    by
    Theorem~\ref{theorem:uniqueness-of-hom-equation-with-initial-values},
    \refitem{item:smooth-global-solutions} applied to $D^\Trans$. By
    Theorem~\ref{theorem:propagator-complex},
    \refitem{item:propagator-complex-exact-on-globally-hyperbolic} we
    know that $\Phi = F_M \varphi$ for some $\varphi \in
    \Secinfty_0(E^*)$. But then $\varrho_\Sigma(\varphi) = (\varphi_0,
    \dot{\varphi}_0)$ follows directly from ($*$). Since the sections
    $(\varphi_0, \dot{\varphi}_0) \in \Secinfty_0(\iota^\#(E^* \oplus
    E^*))$ generate the whole symmetric algebra and $\varrho_\Sigma$
    is an algebra homomorphism, the surjectivity follows. Now let
    $\Phi \in \Sym^\bullet \Secinfty_0(E^*)$. Then
    $\varrho_\Sigma(\Phi) = 0$ iff for all $(u_0, \dot{u}_0)$ we have
    $\varrho_\Sigma(\Phi)(u_0, \dot{u}_0) = 0$. But this is equivalent
    to $\Phi(u) = 0$ for all solutions $u \in
    \Secinfty_{\mathrm{sc}}(E)$ of the homogeneous wave equation $Du =
    0$ by Lemma~\ref{lemma:the-algebra-homomorhpism}. Thus the kernel
    of $\varrho_\Sigma$ consists precisely of those $\Phi \in
    \Sym^\bullet \Secinfty_0(E^*)$ which vanish on solutions. Since
    the kernel is clearly a (Poisson) ideal as $\varrho_\Sigma$ is a
    (Poisson) algebra homomorphism and since the Casimir elements
    $\varphi \in \Secinfty_0(E^*)$ vanish on solutions by
    Lemma~\ref{lemma:degeneration-space-of-poisson-bracket} it follows
    that the ideal generated by the Casimir elements is part of the
    kernel. Now in symmetric degree one the converse is true: $\varphi
    \in \Secinfty_0(E^*)$ is a Casimir element iff it is in the
    kernel. Thus we see that the induced map
    \[
    \varrho_\Sigma:
    \Secinfty_0(E^*) \big/
    \left\{
        \varphi \in \Secinfty_0(E^*)
        \; \big| \;
        \PB{\varphi, \argument} = 0
    \right\}
    \longrightarrow \Secinfty_0(\iota^\#(E^* \oplus E^*))
    \]
    is already a linear isomorphism. Thus we have an algebra
    isomorphism
    \[
    \varrho_\Sigma:
    \Sym^\bullet
    \left(
        \Secinfty_0(E^*) \big/
        \left\{
            \varphi \in \Secinfty_0(E^*)
            \; \big| \;
            \PB{\varphi, \argument} = 0
        \right\}
    \right)
    \longrightarrow
    \Sym^\bullet \Secinfty_0(\iota^\#(E^* \oplus E^*).
    \]
    By a general argument, one has canonically $\Sym^\bullet (V \big/
    W) = \Sym^\bullet V \big/ \mathcal{J}(W)$ for every linear
    subspace $W \subseteq V$, where $\mathcal{J}(W) \subseteq
    \Sym^\bullet V$ is the ideal generated by the elements in
    $W$. Hence we can conclude that $\varrho_\Sigma$ is already
    injective on $\Sym^\bullet \Secinfty_0(E^*)$ modulo the ideal
    generated by the Casimir elements in $\Secinfty_0(E^*)$. Thus the
    two ideals coincide.
\end{proof}

The covariant Poisson bracket gives us automatically the correct
quotient procedure: the \emph{vanishing ideal}\index{Vanishing ideal}
of the subspace of solutions to the wave equation is a Poisson ideal,
which can now be characterized in many equivalent ways: it is the
ideal generated by the Casimir elements (and hence easily seen to be a
Poisson ideal), or, equivalently, the ideal generated by the kernel of
$F_M$, or, equivalently, the kernel of any of the Poisson
homomorphisms $\varrho_\Sigma$ for any Cauchy hypersurface
$\varrho_\Sigma$. However, the physically important interpretation is
the first: two $\Phi, \Psi \in \Sym^\bullet \Secinfty_0(E)$ should be
considered to be the same observables if they yield the same
``expectation values''
\begin{equation}
    \label{eq:same-expectation-values}
    \Phi(u) = \Psi(u)
\end{equation}
for all physically relevant $u \in \Secinfty_{\mathrm{sc}}(E)$,
i.e. for all \emph{solutions} of the wave equation. Note that a priori
it is not clear whether this vanishing ideal of the subspace of
solutions is a Poisson ideal at all. We can now summarize the results
so far.
\begin{theorem}[Covariant Poisson algebra]
    \label{theorem:the-poisson-isomorphisms}%
    \index{Covariant Poisson algebra}%
    Let $(M,g)$ be a globally hyperbolic spacetime and let $E$ be a
    real valued vector bundle with fibre metric and metric connections
    $\nabla^E$. Let $D = \dAlembert^\nabla + B$ be a symmetric,
    normally hyperbolic differential operator on $E$. Moreover, let
    $\PB{\argument, \argument}$ be the covariant Poisson bracket for
    $\Sym^\bullet \Secinfty_0(E^*)$ and let $\iota: \Sigma
    \hookrightarrow M$ be a smooth spacelike Cauchy hypersurface.
    \begin{theoremlist}
    \item \label{item:several-char-of-degen-space} The following
        subspaces of $\Sym^\bullet \Secinfty_0(E^*)$ coincide:
        \begin{compactitem}
        \item The vanishing ideal of the solutions of the wave
            equation $Du = 0$, i.e.
            \begin{equation}
                \label{eq:vanishing-ideal-of-solutions}
                \left\{
                    \Phi \in \Sym^\bullet \Secinfty_0
                    \; \big| \;
                    \Phi(u) = 0
                    \; \textrm{for all} \;
                    u \in \Secinfty_{\mathrm{sc}}(E)
                    \; \textrm{with} \; Du = 0
                \right\}.
            \end{equation}
        \item The ideal generated by the Casimir elements $\varphi \in
            \Secinfty_0(E^*)$.
        \item The ideal generated by the kernel of $F_M:
            \Secinfty_0(E^*) \longrightarrow
            \Secinfty_{\mathrm{sc}}(E^*)$.
        \item The kernel of the Poisson homomorphism
            \begin{equation}
                \label{eq:the-poisson-homomorphism}
                \varrho_\Sigma: \Sym^\bullet \Secinfty_0(E^*)
                \longrightarrow
                \Sym^\bullet \Secinfty_0(\iota^\# (E^* \oplus E^*)).
            \end{equation}
        \end{compactitem}
    \item \label{item:degen-space-is-poisson-ideal} The subspace in
        \refitem{item:several-char-of-degen-space} is a Poisson ideal.
    \item \label{item:the-poisson-isomorphism} The quotient Poisson
        algebra $\Sym^\bullet \Secinfty_0(E^*) \big/ \ker
        \varrho_\Sigma$ is canonically isomorphic to the Poisson
        algebra $\Sym^\bullet (\Secinfty_0(E^*) \big/ \ker F_M)$
        endowed with the induced bracket coming from
        \eqref{eq:poisson-bracket-on-generators} and
        \begin{equation}
            \label{eq:the-poisson-isomorphism}
            \varrho_\Sigma:
            \Sym^\bullet \Secinfty_0(E^*) \big/ \ker \varrho_\Sigma
            \longrightarrow
            \Sym^\bullet \Secinfty_0(\iota^\#(E^* \oplus E^*))
        \end{equation}
        is an isomorphism of Poisson algebra. It is compatible with
        evaluation on solutions and initial data, respectively, in the
        sense of \eqref{eq:the-algebra-homomorphism}.
    \end{theoremlist}
\end{theorem}
\begin{proof}
    All the statements are clear from the preceding lemmas.
\end{proof}

\begin{remark}[Covariant Poisson bracket]
    \label{remark:covariant-poisson-bracket}%
    \index{Dynamical Poisson bracket}%
    The remarkable feature of the Poisson bracket $\PB{\argument,
      \argument}$ on $\Sym^\bullet \Secinfty_0(E^*)$ as well as on the
    quotient $\Sym^\bullet (\Secinfty_0(E^*) \big/ \ker F_M)$ is that
    it does not refer to a splitting $\mathbb{R} \times \Sigma$ of
    $M$. Instead it is ``\emph{fully covariant}'', i.e. defined in
    global and canonical terms only. Nevertheless, via
    $\varrho_\Sigma$ it is isomorphic to the Poisson algebra on the
    Cauchy hypersurface $\Sigma$. The price is that for the
    construction of $\PB{\argument, \argument}$ we have to use the
    dynamics already. This is a new feature as in geometrical
    mechanics the Poisson structure is understood as a \emph{purely
      kinematical} ingredient of the theory. The dynamics comes only
    after specifying a Hamiltonian as an element of the a priori given
    Poisson algebra. Thus the above ``covariant'' Poisson bracket may
    also deserve the name ``\emph{dynamical Poisson bracket}''.
\end{remark}
\begin{remark}[Time evolution]
    \label{remark:time-evolution}%
    \index{Time evolution}%
    Using the Poisson isomorphisms $\varrho_\Sigma$ for different
    Cauchy hypersurfaces we get a time evolution from one Cauchy
    hypersurface to another one. For smooth Cauchy hypersurfaces
    $\Sigma, \Sigma'$ we have
    \begin{equation}
        \label{eq:time-evolution}
        \varrho_{\Sigma'} \circ \varrho_\Sigma^{-1}:
        \Sym^\bullet \Secinfty_0(\iota^\#(E^* \oplus E^*))
        \longrightarrow
        \Sym^\bullet \Secinfty_0(\iota'^\# (E^* \oplus E^*))
    \end{equation}
    with $\varrho_\Sigma , \varrho_{\Sigma'}$ as in
    \eqref{eq:the-poisson-isomorphism}. This is an isomorphism of
    Poisson algebras. In this sense, the time evolution of the wave
    equation is ``symplectic''.
\end{remark}

The next observation would be indeed very complicated and almost
impossible to detect inside the canonical Poisson algebras of
polynomials on the initial data. Here the global point of view indeed
tuns out to be superior: Since we interpret the $\Phi \in \Sym^\bullet
\Secinfty_0(E^*)$ as polynomial observables we can speak of a support
of them. Indeed, we define for
\begin{equation}
    \label{eq:pol-observable}%
    \index{Observable!support}%
    \Phi
    = \sum_{k=0}^N \varphi_1^{(k)} \vee \cdots \vee \varphi_k^{(k)}
\end{equation}
the \emph{support} of $\Phi$ to be the (finite) union of the supports
of the $\varphi_1^{(k)}, \ldots, \varphi_k^{(k)}$. In this sense we
can speak of an observable being located in a certain region of the
spacetime. The physical interpretation is that $\Phi$ corresponds to
an observation (measurement) performed on the solution $u$ in the
spacetime region determined by $\supp \Phi$. Since we consider only
those $\Phi$ coming from compactly supported $\varphi \in
\Secinfty_0(E^*)$ the support of $\Phi$ is also compact. Causality now
means that two measurements $\Phi$ and $\Phi'$ should not influence
each other in any way if they are performed in spacelike regions of
$M$. The next proposition says that this is indeed the case:
\begin{proposition}[Locality]
    \label{proposition:locality}%
    \index{Locality}%
    Let $U, U' \subseteq M$ be open subsets such that $U$ is spacelike
    to $U'$. Then for all $\Phi, \Phi' \in \Sym^\bullet
    \Secinfty_0(E^*)$ with $\supp \Phi \subseteq U$ and $\supp \Phi'
    \subseteq U'$ we have
    \begin{equation}
        \label{eq:locality}
        \PB{\Phi, \Phi'} = 0.
    \end{equation}
\end{proposition}
\begin{proof}
    By the Leibniz rule it is again sufficient to consider $\varphi,
    \varphi' \in \Secinfty_0(E^*)$ with $\supp \varphi \subseteq U$
    and $\supp \varphi' \subseteq U'$ only. But here
    \eqref{eq:locality} is obvious since $\supp \varphi \subseteq U$
    and $\supp F_M \varphi' \subset J_M(\supp \varphi') \subseteq
    J_M(U')$ have no overlap. Thus the integral
    \eqref{eq:poisson-bracket-on-generators} vanishes, see also
    Figure~\ref{fig:locality}.
    \begin{figure}
        \centering
        \input{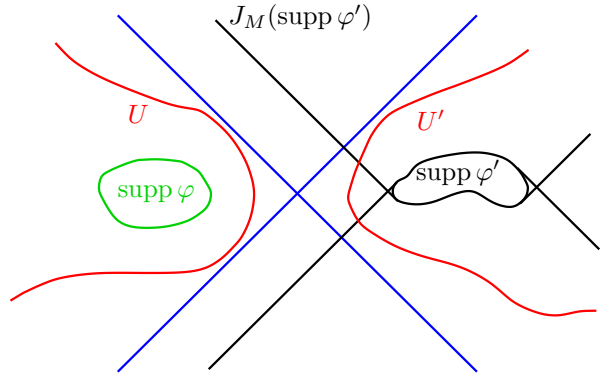}
        \caption{\label{fig:locality} Illustration of the locality
          concept.}
    \end{figure}
\end{proof}

For a later yet to be found transition to a quantum field theory,
i.e. a quantization of the classical observable algebra, it is useful
to consider also the complexification of $\Sym^\bullet
\Secinfty_0(E^*)$. This is the ultimate definition of the classical
observable algebra.
\begin{definition}[Classical observable algebra]
    \label{definition:classicle-observables}%
    \index{Observable}%
    \index{Classical Observable algebra}%
    \index{Poisson salgebra@Poisson $^*$-algebra}
    The classical observable algebra of the classical field theory
    determined by the wave equation is the unital Poisson $^*$-algebra
    \begin{equation}
        \label{eq:classical-observables}
        \mathcal{A}(M)
        = \Sym^\bullet
        \left(
            \Secinfty_0(E^*) \big/ \ker F_M
        \right) \tensor \mathbb{C}
    \end{equation}
    endowed with the complex conjugation as $^*$-involution, the
    symmetric tensor product as associative and commutative product,
    and the covariant Poisson bracket induced from $\Sym^\bullet
    \Secinfty_0(E^*)$.
\end{definition}

Here a Poisson $^*$-algebra means that the $^*$-involution is
compatible with the Poisson bracket in the sense that $\PB{\argument,
  \argument}$ is real, i.e. for $\Phi, \Psi \in \Sym^\bullet
(\Secinfty_0(E^*) \big/ \ker F_M) \tensor \mathbb{C}$ we have
\begin{equation}
    \label{eq:poisson-bracket-is-real}
    \cc{\PB{\Phi,\Psi}} = \PB{\cc{\Phi}, \cc{\Psi}},
\end{equation}
which is obvious as we complexified a Poisson algebra over
$\mathbb{R}$. This is the ultimate reason that we insisted on a real
vector bundle from the beginning. Alternatively, we can write the
complexification as
\begin{equation}
    \label{eq:alternative-complexification}
    \Sym^\bullet \Secinfty_0(E^*) \tensor \mathbb{C}
    = \Sym^\bullet_{\mathbb{C}} \Secinfty_0(E^* \tensor \mathbb{C}),
\end{equation}
where we take the symmetric algebra over the complex numbers of the
section of the complexified bundle. Again, this is compatible with the
quotient procedure since $F_M$ behaves well under complexification,
according to
Proposition~\ref{proposition:symmetry-of-green-operators},
\refitem{item:greenops-and-complexification}.

The locality property clearly passes to the quotient in the following
sense: for an open subset $U \subseteq M$ we define analogously to
\eqref{eq:classical-observables}
\begin{equation}
    \label{eq:compex-local-observables}
    \mathcal{A}_M(U)
    = \Sym^\bullet
    \left(
        \Secinfty_0(E\at{U}) \big/ \ker F_M\at{\Secinfty_0(E|_{U})}
    \right) \tensor \mathbb{C},
\end{equation}
and call this the subalgebra of observables located in $U$. Clearly,
we have natural embeddings
\begin{equation}
    \label{eq:embeddings-of-local-observables-in-each-other}
    \mathcal{A}_M(U) \hookrightarrow
    \mathcal{A}_M(U') \hookrightarrow
    \mathcal{A}_M(M)
\end{equation}
for all $U \subseteq U' \subseteq M$ and each $\mathcal{A}_M(U)$ is a
Poisson $^*$-algebra itself. In this sense, $\mathcal{A}_M(M)$ becomes
the inductive limit (int the category of Poisson $^*$-algebras) of the
collection of the $\mathcal{A}_M(U)$. The important consequence of
Proposition~\ref{proposition:locality} says that we have a local net
of observable algebras:
\begin{theorem}[Local net of observables]
    \label{theorem:local-net-of-observables}%
    \index{Observable!Local net}%
    \index{Causality}%
    The collection of Poisson $^*$-algebras $\left\{ \mathcal{A}_M(U)
        \; \big| \; U \subseteq M \; \textrm{is open}\right\}$ forms a
    net of local observables with inductive limit $\mathcal{A}(M)$,
    satisfying the causality condition
    \begin{equation}
        \label{eq:causality-condition-of-net-of-observables}
        \PB{ \mathcal{A}_M(U), \mathcal{A}_M(U') } = 0
    \end{equation}
    for $U, U' \subseteq M$ spacelike to each other.
\end{theorem}
\begin{remark}
    \label{remark:local-net-of-observables}%
    \index{Haag-Kastler axioms}%
    This property is the classical analogy of one of the Haag-Kastler
    axioms for an (algebraic or axiomatic) quantum field theory:
    observables in spacelike regions should commute. We refer to
    \cite{haag:1993a} for further information on algebraic quantum
    field theory. Note that it would be extremely complicated to
    encode this net structure in the canonical Poisson algebra over
    $\Sigma$: here the covariant approach turns out to be the better
    choice.
\end{remark}
We also have the following version of the time slice axiom:
\begin{theorem}[Time slice axiom]
    \label{theorem:time-slice-axiom}%
    \index{Time slice axiom}%
    Let $\iota: \Sigma \hookrightarrow M$ be a smooth spacelike Cauchy
    hypersurface coming from a splitting $M \simeq \mathbb{R} \times
    \Sigma$. Let $\epsilon > 0$, then we have
    \begin{equation}
        \label{eq:time-slice-axiom}
        \mathcal{A}_M( (-\epsilon,\epsilon) \times \Sigma)
        = \mathcal{A}_M(M).
    \end{equation}
\end{theorem}
\begin{proof}
    This equality is of course not true on the level of the polynomial
    algebra $\Sym^\bullet \Secinfty_0(E)$ itself since there are
    certainly elements with support \emph{outside}
    $(-\epsilon,\epsilon) \times \Sigma$. The point is that they are
    equivalent to elements in $\Sym^\bullet
    \Secinfty_0\left(E^*\at{(-\epsilon,\epsilon) \times
          \Sigma}\right)$ modulo the kernel of $F_M$. First we note
    that $(-\epsilon,\epsilon) \times \Sigma$ is again a globally
    hyperbolic spacetime by its own. Moreover, the embedding of
    $(-\epsilon,\epsilon) \times \Sigma$ into $M = \mathbb{R} \times
    \Sigma$ is causally compatible. We can now apply our theory of
    Green operators to $D$ and $D^\Trans$ restricted to
    $(-\epsilon,\epsilon) \times \Sigma$ and obtain unique Green
    operators $F^\pm_{(-\epsilon,\epsilon) \times \Sigma}$ for
    $D^\Trans$ on $(-\epsilon,\epsilon) \times \Sigma$, too. Since
    $(-\epsilon,\epsilon) \times \Sigma$ is causally compatible in
    $M$, the support properties of $F^\pm_{(-\epsilon,\epsilon) \times
      \Sigma}$ match those of $F^\pm_M$ ``restricted'' to
    $(-\epsilon,\epsilon) \times \Sigma$. Thus by uniqueness we
    conclude that for $\varphi \in
    \Secinfty_0\left(E^*\at{(-\epsilon,\epsilon) \times
          \Sigma}\right)$ we have
    \[
    F^\pm_{(-\epsilon,\epsilon) \times \Sigma} \varphi
    = F^\pm_M \varphi \at{(-\epsilon,\epsilon) \times \Sigma}.
    \tag{$*$}
    \]
    This implies that on $\Sym^\bullet
    \Secinfty_0\left(E^*\at{(-\epsilon,\epsilon) \times
          \Sigma}\right)$ the covariant Poisson bracket coming from
    $F_{(-\epsilon,\epsilon) \times \Sigma}$ coincides with the
    restriction of the covariant Poisson bracket coming from
    $F_M$. Moreover, if $\varphi \in \Secinfty_0(E^*)$ with support in
    $(-\epsilon,\epsilon) \times \Sigma$ vanishes on $u \in
    \Secinfty_{\mathrm{sc}}(E)$ satisfying $Du=0$ on $M$ it also
    vanishes on $u \in \Secinfty_{\mathrm{sc}}(E
    \at{(-\epsilon,\epsilon) \times \Sigma})$ satisfying $Du = 0$ on
    $(-\epsilon, \epsilon) \times \Sigma$. Indeed, in the condition
    $\varphi(u)=0$ only $u\at{(-\epsilon,\epsilon) \times \Sigma}$
    enters. This shows that
    \[
    \ker F_{(-\epsilon,\epsilon) \times \Sigma}
    = \ker \left(
        F_M \at{\Secinfty_0\left(
              E\at{(-\epsilon,\epsilon) \times \Sigma}
          \right)
        }
    \right).
    \]
    Therefore the Poisson $^*$-algebra
    $\mathcal{A}_M((-\epsilon,\epsilon) \times \Sigma)$ built using
    $F_M$ and the Poisson $^*$-algebra
    $\mathcal{A}_{(-\epsilon,\epsilon) \times
      \Sigma}((-\epsilon,\epsilon) \times \Sigma)$ coincide. Now we
    have the Poisson $^*$-isomorphisms
    \[
    \varrho_\Sigma:
    \mathcal{A}_{(-\epsilon,\epsilon) \times \Sigma}
    \longrightarrow
    \Sym^\bullet (\Secinfty_0(\iota^\#(E^* \oplus E^*)))
    \tensor \mathbb{C},
    \]
    according to Theorem~\ref{theorem:the-poisson-isomorphisms}
    applied to the spacetime $(-\epsilon,\epsilon) \times \Sigma$ as
    well as
    \[
    \varrho_\Sigma^{-1}:
    \Sym^\bullet( \Secinfty_0(\iota^\#(E^* \oplus E^*))) \tensor
    \mathbb{C}
    \longrightarrow \mathcal{A}_M(M),
    \]
    also using Theorem~\ref{theorem:the-poisson-isomorphisms}, now for
    the spacetime $M$. But this shows the equality
    \eqref{eq:time-slice-axiom}.
\end{proof}

\begin{remark}[Time slice axiom]
    \label{remark:time-slice-axiom}%
    \index{Time slice axiom}%
    We can rephrase this statement by saying that for every $\varphi
    \in \Secinfty_0(E^*)$ with arbitrary compact support there is also
    a $\psi \in \Secinfty_0\left(E^* \at{(-\epsilon,\epsilon) \times
          \Sigma}\right)$ having compact support very close to
    $\Sigma$ such that their images under $\varrho_\Sigma$ in
    $\Secinfty_0(\iota^\#(E^* \oplus E^*))$ coincide, see also
    Figure~\ref{fig:time-slice-axiom}.
    \begin{figure}
        \centering
        \input{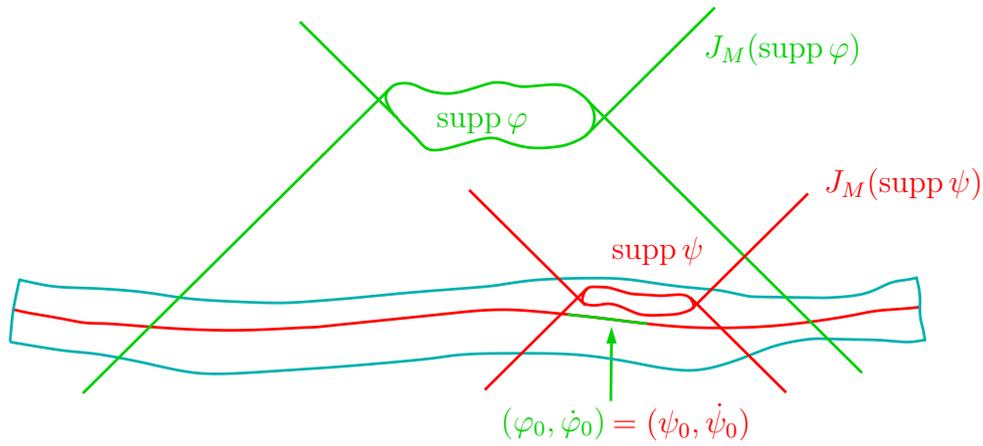}
        \caption{\label{fig:time-slice-axiom}%
          The time slice axiom
        }
    \end{figure}
    Since we know that $\varrho_\Sigma$ is injective up to elements in
    $\ker F_M$ which is the image of $D^\Trans$ by
    \eqref{eq:propagator-compex} according to
    Theorem~\ref{theorem:propagator-complex}, we see that for every
    $\varphi \in \Secinfty_0(E^*)$ there is a $\psi \in
    \Secinfty_0(E^*)$ with support in $(-\epsilon,\epsilon) \times
    \Sigma$ such that
    \begin{equation}
        \label{eq:time-slice-consequence}
        \varphi - \psi \in \ker F_M = \image D^\Trans,
    \end{equation}
    see also \cite[Lem.~4.5.6]{baer.ginoux.pfaeffle:2007a} for another
    approach to this question. Physically speaking, the time slice
    feature \eqref{eq:time-slice-axiom} says that on the level of
    observables a Cauchy hypersurface already determines
    everything. In view of our previous results this is of course not
    very surprising.
\end{remark}


%% file: appendixA.tex
%
%

\chapter{Parallel Transport, Jacobi Vector Fields, and all that}
\label{cha:parall-transp-jacobi}

In this appendix we collect some facts on parallel transports, Taylor
expansions and Jacobi vector fields needed in the computation of the
derivatives of densities.

%
%

\section{Taylor Expansion of Parallel Transports}
\label{sec:taylor-parallel-transport}

\input{taylor}

%
%

\section{Jacobi Vector Fields and the Tangent Map of $\exp_p$}
\label{sec:jacobi-and-tangent-of-exp}

\input{jacobi}

%
%

\section{Jacobi Determinants of the Exponential Map}
\label{sec:jacobi-determinants-exp}

\input{jacobiexp}


%% file: taylor.tex
%
%

Let $\nabla$ be a torsion-free covariant derivative for $M$ and let
$\nabla^E$ be a covariant derivative for a vector bundle $E
\longrightarrow M$. The aim is to compute the Taylor expansion of the
parallel transport with respect to $E$ along curves in $M$. Of
particular interest will be the geodesics with respect to $\nabla$.

Out of $\nabla$ and $\nabla^E$ we can build covariant derivatives for
all kind of bundles constructed from $TM$ and $E$ via dualizing and
taking tensor products. We will denote them all by $\nabla$ or
$\nabla^E$ if $E$ is involved. If $\gamma: I \subseteq \mathbb{R}
\longrightarrow M$ is a smooth curve defined on some open interval
then the pull-back connection of $\nabla$ or $\nabla^E$ will be
denoted by $\nabla^\#$. The canonical vector field on $\mathbb{R}$ is
$\frac{\partial}{\partial t}$.

\begin{lemma}
    \label{lemma:time-derivative-parallel-transport}%
    \index{Parallel transport}%
    Let $\gamma: I \subseteq \mathbb{R} \longrightarrow M$ be a smooth
    curve in $M$ and let $s \in \Secinfty(\gamma^\# E)$ be a section
    of $E$ along $\gamma$. For $t, t_0 \in I$ and all $k \in
    \mathbb{N}_0$ we have
    \begin{equation}
        \label{eq:time-derivative-of-parallel-transport}
        \frac{\D^k}{\D\!t^k}
        \left( P_{\gamma, t_0 \rightarrow t} \right)^{-1} s(t)
        = \left( P_{\gamma, t_0 \rightarrow t} \right)^{-1}
        \left(
            \nabla^\#_{\frac{\partial}{\partial t}} \cdots 
            \nabla^\#_{\frac{\partial}{\partial t}}
            s(t)    
        \right),
    \end{equation}
    where $P_{\gamma, t_0 \rightarrow t}: E_{\gamma(t_0)} \rightarrow
    E_{\gamma(t)}$ denotes the parallel transport along $\gamma$ with
    respect to $\nabla^E$.
\end{lemma}
\begin{proof}
    We choose a vector space basis $e_\alpha(t_0) \in E_{\gamma(t_0)}$
    and define smooth sections $e_\alpha \in \Secinfty(\gamma^\# E)$
    along $\gamma$ by
    \[
    e_\alpha(t) = P_{\gamma, t_0 \rightarrow t} e_\alpha(t_0),
    \]
    i.e. by parallel transporting $e_\alpha(t_0)$ to every point
    $\gamma(t)$ for $t \in I$. Since the parallel transport is a
    linear isomorphism, for every $t$ the $e_\alpha(t)$ still form a
    basis of $E_{\gamma(t)}$.  By the very definition, the
    $e_\alpha(t)$ solve the differential equation
    \[
    \nabla^\#_{\frac{\partial}{\partial t}} e_\alpha(t) = 0
    \tag{$*$}
    \]
    with initial conditions $e_\alpha(t_0) \in \E_{\gamma(t_0)}$. Thus
    they are covariantly constant along $\gamma$. Now let $s \in
    \Secinfty(\gamma^\# E)$ be arbitrary. Then there are unique smooth
    functions $s^\alpha \in \Cinfty(I)$ with
    \[
    s(t) = s^\alpha(t) e_\alpha(t).
    \]
    By linearity of $P_{\gamma, t_0 \rightarrow t}$ we have
    \begin{align*}
        \left( P_{\gamma, t_0 \rightarrow t} \right)^{-1} (s(t))
        = s^\alpha(t)
        \left(P_{\gamma, t_0 \rightarrow t} \right)^{-1}
        (e_\alpha(t))
        = s^\alpha(t)
        \left( P_{\gamma, t_0 \rightarrow t} \right)^{-1}
        P_{\gamma, t_0 \rightarrow t} (e_\alpha(t_0))
        = s^\alpha(t) e_\alpha(t_0).
    \end{align*}
    This shows that we can express the left hand side of
    \eqref{eq:time-derivative-of-parallel-transport}, being a curve in
    the vector space $E_{\gamma(t_0)}$, with respect to the
    \emph{fixed} basis $e_\alpha(t_0)$. Thus the $t$-derivatives are
    easily computed giving
    \begin{align*}
        \frac{\D}{\D\! t}
        \left( P_{\gamma, t_0 \rightarrow t} \right)^{-1}
        s(t)
        &= \frac{\D}{\D\! t} (s^\alpha(t) e_\alpha(t_0))
        = \dot{s}^\alpha(t) e_\alpha(t_0)
        = \dot{s}^\alpha(t)
        \left( P_{\gamma, t_0 \rightarrow t} \right)^{-1}
        \left(
            P_{\gamma, t_0 \rightarrow t} (e_\alpha(t_0))
        \right) \\
        &= \left( P_{\gamma, t_0 \rightarrow t} \right)^{-1}
        (\dot{s}^\alpha(t) e_\alpha(t))
        \stackrel{\mathclap{(*)}}{=}
        \left( P_{\gamma, t_0 \rightarrow t} \right)^{-1}
        \left(
            \nabla^\#_{\frac{\partial}{\partial t}}
            (s^\alpha(t) e_\alpha(t))
        \right) \\
        &= \left( P_{\gamma, t_0 \rightarrow t} \right)^{-1}
        \left(
            \nabla^\#_{\frac{\partial}{\partial t}} s(t)
        \right),
    \end{align*}
    by the covariant constancy of the $e_\alpha(t)$. This shows
    \eqref{eq:time-derivative-of-parallel-transport} for $k=1$ and
    from here we can proceed by induction.
\end{proof}

The next lemma will be useful to compute the Taylor coefficients of a
function of several variables in an efficient way. The proof is a
simple computation.
\begin{lemma}
    \label{lemma:derivative-formula}
    Let $F \in \Cinfty(\mathbb{R}^n, \mathbb{R}^m)$ and $k \in
    \mathbb{N}_0$. Then one has
    \begin{equation}
        \label{eq:derivative-formula}
        \frac{\partial^k F}{\partial v^{i_1} \cdots \partial v^{i_k}}
        \At{v = 0}
        = \frac{1}{k!}
        \frac{\partial^k}{\partial v^{i_1} \cdots \partial v^{i_k}}
        \left(
            \frac{\D^k}{\D t^k} F(tv) \At{t=0}
        \right).
    \end{equation}
\end{lemma}
The following technical lemma will allow us to compute iterated
covariant derivatives in terms of the symmetrized covariant
derivative.
\begin{lemma}
    \label{lemma:iterated-covariant-derivatives}
    For $s \in \Secinfty(E)$ one inductively defines
    \begin{equation}
        \label{eq:covariant-0}
        \begin{split}
            \nabla^0 s &= s, \\            
            \left( \nabla^1 s \right) (X) &= \nabla^E_X s, \\
            \left( \nabla^k s \right) (X_1, \ldots, X_k)
            &= \left( \nabla^E_{X_1} \nabla^{k-1} s \right)
            (X_2, \ldots, X_k)
        \end{split}
    \end{equation}
    for $X_1, \ldots, X_k \in \Secinfty(TM)$. Then $\nabla^k s \in
    \Secinfty( \tensor^k T^*M \tensor E)$ is a well-defined tensor
    field and we have
    \begin{equation}
        \label{eq:symmetrized-covariant-k-is-symD}
        \sum_{\sigma \in \Sym_k} \left(\nabla^k s\right)
        (X_{\sigma(1)}, \ldots, X_{\sigma(k)})
        = \left(
            \left(\SymD^E\right)^k s
        \right)
        (X_1, \ldots, X_k)
    \end{equation}
    for the totally symmetric part of $\nabla^k s$.
\end{lemma}
\begin{proof}
    By induction it is clear that $\nabla^k s$ is $\Cinfty(M)$-linear
    in each argument. Thus it defines a tensor field of the above
    type. To prove \eqref{eq:symmetrized-covariant-k-is-symD} we first
    note that for $k = 0,1$ we have $\nabla^0 s = s = (\SymD^E)^0 s$
    and $\nabla^1 s = \SymD^E s$ as wanted. We proceed by induction
    and have
    \begin{align*}
        \sum_{\sigma \in \Sym_k} \left( \nabla^k s \right)
        (X_{\sigma(1)}, \ldots, X_{\sigma(k)})
        &= \sum_{\sigma \in \Sym_k}
        \left( \nabla^E_{X_{\sigma(1)}} \nabla^{k-1} s \right)
        (X_{\sigma(2)}, \ldots, X_{\sigma(k)})\\
        &= \sum_{\ell=1}^k
        \sum_{\substack{\sigma \in \Sym_k \\ \sigma(1) = \ell}}
        \left( \nabla^E_{X_\ell} \nabla^{k-1} s \right)
        (X_{\sigma(2)}, \ldots, X_{\sigma(k)}) \\
        &= \sum_{\ell=1}^k
        \left(
            \nabla^E_{X_\ell} (\SymD^E)^{k-1} s
        \right)
        (X_2, \ldots, \stackrel{\ell}{\wedge}, \ldots, X_k) \\
        &= \left(
            (D^E)^k s
        \right)
        (X_1, \ldots, X_k),
    \end{align*}
    since the permutations $\sigma \in \Sym_k$ with $\sigma(1) = \ell$
    are precisely the permutations of the remaining $1, \ldots,
    \stackrel{\ell}{\wedge}, \ldots, k$ entries.
\end{proof}

Since covariant derivatives are extended to tensor bundles and dual
bundles in such a way that we have Leibniz rules with respect to
tensor products and natural pairings, the parallel transport enjoys
homomorphism properties in the following sense:
\begin{lemma}
    \label{lemma:parallel-transport-tensor-product-pairing}
    Let $\gamma: I \subseteq \mathbb{R} \longrightarrow M$ be a smooth
    curve in $M$ and $f \in \Cinfty(M)$, $s_{t_0}, \widetilde{s}_{t_0}
    \in E_{\gamma(t_0)}$ and $\alpha_{t_0} \in E^*_{\gamma(t_0)}$.
    \begin{lemmalist}
    \item \label{item:parallel-functions} Viewing $f$ as a sections of
        $\tensor^0 E$, $\gamma^\# f = \gamma^* f$ is parallel if and
        only if $f$ is constant along $\gamma$.
    \item \label{item:parallel-tansport-of-tensor-product} For all $t
        \in I$ we have
        \begin{equation}
            \label{eq:parallel-transport-of-tensor-product}
            P_{\gamma, t_0 \rightarrow t}
            (s_{t_0} \tensor \widetilde{s}_{t_0})
            =
            P_{\gamma, t_0 \rightarrow t} (s_{t_0}) \tensor
            P_{\gamma, t_0 \rightarrow t} (\widetilde{s}_{t_0}).
        \end{equation}
    \item \label{item:parallel-transport-and-pairing} For all $t \in I$ we have
        \begin{equation}
            \label{eq:parallel-transport-and-pairing}
            \alpha_{t_0} (s_{t_0})
            = P_{\gamma, t_0 \rightarrow t} (\alpha_{t_0} (s_{t_0}))
            = P_{\gamma, t_0 \rightarrow t}(\alpha_{t_0})
            \left(
                P_{\gamma, t_0 \rightarrow t} (s_{t_0})
            \right).
        \end{equation}
    \end{lemmalist}
\end{lemma}
\begin{proof}
    For the first part we observe that by definition $\nabla^E_X f =
    X(f)$ when viewing a function as a tensor field. Moreover,
    $\gamma^\# f = \gamma^* f$ and hence
    $\nabla^\#_{\frac{\partial}{\partial t}} \gamma^\#f =
    \frac{\partial}{\partial t} \gamma^* f \at{t} = \dot{\gamma}(t) f
    \at{\gamma(t)}$ which is zero iff $f \circ \gamma$ is constant. It
    follows that for a \emph{number} $z \in \tensor^0 E_{\gamma(t_0)}
    = \mathbb{C}$ we simply have $P_{\gamma, t_0 \rightarrow t} (z) =
    z$ for all times. This shows the first part. For the second part
    we note that the left hand side is the unique solution of
    \[
    \nabla^\#_{\frac{\partial}{\partial t}}
    P_{\gamma, t_0 \rightarrow t}
    (s_{t_0} \tensor \widetilde{s}_{t_0})
    = 0
    \]
    with initial condition $s_{t_0} \tensor \widetilde{s}_{t_0}$ for
    $t=t_0$. For the right hand side we compute
    \begin{align*}
        &\nabla^\#_{\frac{\partial}{\partial t}}
        \left(
            P_{\gamma, t_0 \rightarrow t} (s_{t_0})
            \tensor P_{\gamma, t_0 \rightarrow t}
            (\widetilde{s}_{t_0})
        \right) \\
        &\quad=
        \nabla^\#_{\frac{\partial}{\partial t}}
        \left(
            P_{\gamma, t_0 \rightarrow t}(s_{t_0})
        \right)
        \tensor
        P_{\gamma, t_0 \rightarrow t} (\widetilde{s}_{t_0}) 
         + P_{\gamma, t_0 \rightarrow t}(s_{t_0}) \tensor
        \nabla^\#_{\frac{\partial}{\partial t}}
        P_{\gamma, t_0 \rightarrow t} (\widetilde{s}_{t_0})
        = 0,
    \end{align*}
    by the Leibniz rule of $\nabla^\#$ for sections
    $\Secinfty(\gamma^\# (E \tensor E))$ with respect to
    $\tensor$. Since the right hand side of
    \eqref{eq:parallel-transport-of-tensor-product} is $s_{t_0}
    \tensor \widetilde{s}_{t_0}$ for $t = t_0$ we have
    \eqref{eq:parallel-transport-of-tensor-product} by
    uniqueness. Analogously, one shows
    \refitem{item:parallel-transport-and-pairing}.
\end{proof}

By combination of \refitem{item:parallel-tansport-of-tensor-product}
and \refitem{item:parallel-transport-and-pairing} we obtain the
compatibility of parallel transport with the usual tensor product
constructions and multilinear pairings. We shall use this frequently
in the following.

Since a covariant derivative $\nabla^E$ also induces a covariant
derivative for the density bundles we consider the compatibility of
the parallel transport with the evaluation of a density on a basis. To
this end we first recall the definition of the covariant derivative of
a density. If $A_\alpha^\beta \in \Secinfty(T^*U)$ denote the local
connection one-forms of $\nabla^E$ with respect to a local frame
$e_\alpha \in \Secinfty\left(E\at{U}\right)$ then the covariant
derivative of a $z$-density $\mu \in \Secinfty(\Dichten^z E^*)$ is
defined locally by
\begin{equation}
    \label{eq:covariant-derivative-of-density}%
    \index{Covariant derivative!density}%
    (\nabla_X \mu) (e_1, \ldots, e_N)
    = X( \mu(e_1, \ldots, e_N) )
    - z \sum_{\alpha = 1}^N A_\alpha^\alpha(X) \mu(e_1, \ldots, e_N),
\end{equation}
where $N = \rank E$ and $z \in \mathbb{C}$, see e.g.
\cite[Sect.~2.2]{waldmann:2007a} for this approach and the proof that
\eqref{eq:covariant-derivative-of-density} indeed gives a globally
defined $\nabla_X \mu \in \Secinfty(\Dichten^z E^*)$. We shall
interpret \eqref{eq:covariant-derivative-of-density} in a more global
way. Since $\mu$ is \emph{not} multilinear in the arguments $e_1,
\ldots, e_N$ we can \emph{not} expect a simple Leibniz rule (and hence
an alternative global definition of $\nabla_X \mu$) for the covariant
derivative of a $z$-density. Instead we shall solve the differential
equation
\begin{equation}
    \label{eq:diffequation-for-covariant-differentiated-density}
    \nabla^\#_{\frac{\partial}{\partial t}} \mu = 0
\end{equation}
for a $\mu \in \Secinfty(\gamma^\# \Dichten^z E^*)$ explicitly. To this
end, we note that $\gamma^\#(\Dichten^z E^*) = \Dichten^z (\gamma^\#
E)^*$. Thus we can evaluate
\eqref{eq:diffequation-for-covariant-differentiated-density} in a
local frame of $\gamma^\#E$ giving the equivalent local condition
\begin{equation}
    \label{eq:local-diffeq-for-cov-diff-density}
    0 = (\nabla^\#_{\frac{\partial}{\partial t}} \mu)
    (e_1, \ldots, e_N)
    = \frac{\partial}{\partial t} (\mu(e_1, \ldots, e_N))
    - z \sum_{\alpha=1}^N A^\#{}^\alpha_\alpha
    \left(\frac{\partial}{\partial t}\right)
    \mu(e_1, \ldots, e_N)
\end{equation}
where $A^\#{}^\alpha_\beta \in \Secinfty(T^*I)$ are the local
connection one-forms of $\nabla^\#$ with respect to the local frame
$e_\alpha \in \Secinfty(\gamma^\# E)$. Note that
\eqref{eq:local-diffeq-for-cov-diff-density} is valid not only for
frames of the form $\gamma^\# e_\alpha$ but for all frames. In
particular, we can choose a covariantly constant frame as in the proof
of Lemma~\ref{lemma:time-derivative-parallel-transport}. This simply
means that $A^\#{}^\beta_\alpha = 0$ for such a frame. Thus we arrive
at the statement that for a covariantly constant frame we have
\eqref{eq:local-diffeq-for-cov-diff-density} iff
\begin{equation}
    \label{eq:local-diffeq-cov-diff-dens-cov-constant-frame}
    \frac{\partial}{\partial t} (\mu(e_1, \ldots, e_N)) = 0.
\end{equation}
This means that for a covariantly constant frame the function
$\mu(e_1, \ldots, e_N)$ is constant. Conversely, if $\mu(e_1, \ldots,
e_N)$ is constant for a covariantly constant frame then
$A^\#{}^\beta_\alpha = 0$. Hence by
\eqref{eq:covariant-derivative-of-density} we conclude that $\mu$ is
covariantly constant. From this we obtain the following statement:
\begin{lemma}[Parallel transport of densities]
    \label{lemma:parallel-transport-of-densities}%
    \index{Parallel transport!density}%
    Let $z \in \mathbb{C}$ and $\gamma: I \subseteq \mathbb{R}
    \longrightarrow M$ a smooth curve. For a $z$-density $\mu \in
    \Dichten^z E_{\gamma(t_0)}^*$ and a basis $e_1, \ldots, e_N \in
    E_{\gamma(t_0)}$ we have
    \begin{equation}
        \label{eq:parallel-transport-of-densities}
        \mu(e_1, \ldots, e_N)
        = \left(
            P_{\gamma, t_0 \rightarrow t} (\mu)
        \right)
        \left(
            P_{\gamma, t_0 \rightarrow t} (e_1), \ldots,
            P_{\gamma, t_0 \rightarrow t} (e_N)
        \right).
    \end{equation}
\end{lemma}
\begin{proof}
    Let $\mu(t) = P_{\gamma, t_0 \rightarrow t}(\mu) \in \Dichten^z
    E_{\gamma(t)}^*$ and let $e_\alpha \in \Secinfty(\gamma^\# E)$ be a
    covariantly constant frame, i.e. $e_\alpha(t) = P_{\gamma, t_0
      \rightarrow t} (e_\alpha(t_0))$. Then we know that
    \[
    \mu( e_1(t_0), \ldots, e_N(t_0))
    = \mu(t) (e_1(t), \ldots, e_N(t))
    \]
    by our previous considerations. But this is
    \eqref{eq:parallel-transport-of-densities}.
\end{proof}

Thus also here the parallel transport has ``homomorphism
properties''. Note however that the covariant derivative does not obey
a simple Leibniz rule with respect to the ``pairing'' of a $z$-density
and a frame.

Now we consider geodesics $\gamma(t) = \exp_p(tv)$ with respect to
$\nabla$ instead of arbitrary curves. Since in this case
$\dot{\gamma}$ is covariantly constant along $\gamma$ we obtain the
following lemma:
\begin{lemma}
    \label{lemma:derivatives-of-cov-pullbacked-section}
    Let $s \in \Secinfty(E)$ and let $\gamma: I \subseteq \mathbb{R}
    \longrightarrow M$ by a geodesic. Then we have for all $k \in
    \mathbb{N}_0$ and $t \in I$
    \begin{equation}
        \label{eq:derivatives-of-cov-pullback-section}
        \Big(
            \underbrace{
              \nabla^\#_{\frac{\partial}{\partial t}} \cdots
              \nabla^\#_{\frac{\partial}{\partial t}}}_ {k \; \textrm{times}}
            \gamma^\# s
        \Big) (t)
        = \left(
            \nabla^k s \at{\gamma(t)}
        \right)
        (\dot{\gamma}(t), \ldots, \dot{\gamma}(t)).
    \end{equation}
\end{lemma}       
\begin{proof}
    For $k = 0$ the statement is clearly correct. For $k = 1$ we have
    \begin{align*}
        \left(
            \nabla^\#_{\frac{\partial}{\partial t}} \gamma^\#s
        \right)(t)
        =
        \left( \nabla_{\dot{\gamma}(t)} s \right) (\gamma(t))
        =
        \nabla s \at{\gamma(t)} (\dot{\gamma}(t))
    \end{align*}
    by definition of $\nabla^\#$. Thus
    \eqref{eq:derivatives-of-cov-pullback-section} holds for $k=1$ as
    well. The general case follows by induction since
    \begin{align*}
        \nabla^k s \at{\gamma(t)}
        (\dot{\gamma}(t), \ldots, \dot{\gamma}(t))
        & = \left(
            \nabla^\#_{\frac{\partial}{\partial t}} \nabla^{k-1} s
        \right) \At{\gamma(t)}
        (\dot{\gamma}(t), \ldots, \dot{\gamma}(t)) \\
        & = \left(
            \nabla^\#_{\frac{\partial}{\partial t}}
            \left(
                \gamma^\# \nabla^{k-1}s
            \right)
        \right) \At{t}
        (\dot{\gamma}(t), \ldots, \dot{\gamma}(t)) \\
        & = \nabla^\#_{\frac{\partial}{\partial t}} 
        \left(
            (\gamma^\# \nabla^{k-1} s)
        \right) \At{t}
        (\dot{\gamma}(t), \ldots, \dot{\gamma}(t)) \\
        & \quad - \sum_{\ell = 1}^{k-1} \gamma^\# \nabla^{k-1} s \at{t}
        (\dot{\gamma}(t), \ldots, 
        \nabla^\#_{\frac{\partial}{\partial t}} \dot{\gamma} \at{t},
        \ldots,  \dot{\gamma}(t)) \\
        & = \nabla^\#_{\frac{\partial}{\partial t}} \cdots
        \nabla^\#_{\frac{\partial}{\partial t}} s \at{t}
        - 0,
    \end{align*}
    using that $\dot{\gamma}$ is covariantly constant.
\end{proof}

Using Lemma~\ref{lemma:iterated-covariant-derivatives} we can rephrase
the statement \eqref{eq:derivatives-of-cov-pullback-section} using the
symmetrized covariant derivative since we only evaluate $\nabla^k
s\at{\gamma(t)}$ on $k$ times the same vector $\dot{\gamma}(t)$. Thus
we have
\begin{equation}
    \label{eq:k-derivative-of-section-and-SymD}
    \nabla^\#_{\frac{\partial}{\partial t}} \cdots
    \nabla^\#_{\frac{\partial}{\partial t}} s
    = \frac{1}{k!}
    \left(
        \gamma^\# \left(\SymD^E\right)^k s
    \right)
    (\dot{\gamma}, \ldots, \dot{\gamma}),
\end{equation}
taking into account the correct combinatorics. We can use this now to
compute the Taylor coefficients of the parallel transport along
geodesics in general:
\begin{proposition}[Taylor coefficients of the parallel transport]
    \label{proposition:taylor-coefficients-of-partrans}%
    \index{Parallel transport!Taylor coefficients}%
    Let $k \in \mathbb{N}_0$ and $s \in \Secinfty(E)$ be given. Denote
    by $\gamma_v(t) = \exp_p(tv)$ the geodesic starting at $p$ with
    velocity $v \in T_pM$. Then the Taylor coefficients of the
    parallel transport in radial directions are given by
    \begin{equation}
        \label{eq:taylor-coeffients-of-partrans}
        \frac{\partial^k}{\partial v^{i_1} \cdots \partial v^{i_k}}
        \left(
            P_{\gamma_v, 0 \rightarrow 1}
        \right)^{-1}
        s(\gamma_v(1))\At{v=0}
        =
        \inss(e_{i_1}) \cdots \inss(e_{i_k}) \frac{1}{k!}
        \left(\SymD^E\right)^k s\at{p},
    \end{equation}
    where $v^1, \ldots, v^k$ are the linear coordinates on $T_pM$ with
    respect to a vector space basis $e_1, \ldots, e_n \in T_pM$.
\end{proposition}
\begin{proof}
    First note that $s(\gamma_v(1)) \in E_{\gamma_v(1)}$ whence
    $\left( P_{\gamma_v, 0 \rightarrow 1} \right)^{-1} s(\gamma_v(1))
    \in E_{\gamma_v(0)} = E_p$ is indeed a vector in $E_p$ for all $v
    \in T_pM$. Thus the map
    \[
    v
    \; \mapsto \;
    \left( P_{\gamma_v, 0 \rightarrow 1} \right)^{-1} s(\gamma_v(1))
    \]
    is a smooth $E_p$-valued function on $T_pM$ defined on an open
    neighborhood of $0$. Thus we can apply
    Lemma~\ref{lemma:derivative-formula} to compute its Taylor
    coefficients. We obtain
    \begin{align*}
        \frac{\partial^k}{\partial v^{i_1} \cdots \partial v^{i_k}} &
        \left( P_{\gamma_v, 0 \rightarrow 1} \right)^{-1}
        s(\gamma_v(1))\At{v=0} \\
        &\stackrel{\mathclap{\textrm{Lem.~\ref{lemma:derivative-formula}}}}{=}
        \qquad
        \frac{\partial^k}{\partial v^{i_1} \cdots \partial v^{i_k}}
        \frac{1}{k!} \frac{\partial^k}{\partial t^k}\At{t=0}
        \left(
            P_{\gamma_{tv}, 0 \rightarrow 1}
        \right)^{-1}
        s(\gamma_{tv}(1)) \\
        &=
        \frac{\partial^k}{\partial v^{i_1} \cdots \partial v^{i_k}}
        \frac{1}{k!} \frac{\partial^k}{\partial t^k}\At{t=0}
        \left(
            P_{\gamma_{v}, 0 \rightarrow t}
        \right)^{-1}
        s(\gamma_{v}(t)) \\
        &\stackrel{\mathclap{\textrm{Lem.~\ref{lemma:time-derivative-parallel-transport}}}}{=}
        \qquad
        \frac{\partial^k}{\partial v^{i_1} \cdots \partial v^{i_k}}
        \frac{1}{k!} \frac{\partial^k}{\partial t^k}\At{t=0}
        \left(
            P_{\gamma_{v}, 0 \rightarrow t}
        \right)^{-1}
        \left(
            \nabla^\#_{\frac{\partial}{\partial t}} \cdots
            \nabla^\#_{\frac{\partial}{\partial t}} s
            (\gamma_v(t))
        \right) \\
        &\stackrel{\mathclap{\eqref{eq:k-derivative-of-section-and-SymD}}}{=}
        \quad
        \frac{\partial^k}{\partial v^{i_1} \cdots \partial v^{i_k}}
        \frac{1}{k!} \frac{\partial^k}{\partial t^k}\At{t=0}
        \left(
            P_{\gamma_{v}, 0 \rightarrow t}
        \right)^{-1}
        \left(
            \frac{1}{k!}
            \left(
                \gamma^\# \left(\SymD^E\right)^k s
            \right)
            (\dot{\gamma}, \ldots, \dot{\gamma})
        \right) \At{t=0} \\
        &=
        \frac{\partial^k}{\partial v^{i_1} \cdots \partial v^{i_k}}
        \frac{1}{k!} \frac{1}{k!}
        \left(\SymD^E\right)^k s \at{p} (v, \ldots, v) \\
        &=
        \frac{\partial^k}{\partial v^{i_1} \cdots \partial v^{i_k}}
        \frac{1}{k!} \frac{1}{k!} v^{j_1} \cdots v^{j_k}
        \left(\SymD^E\right)^k s \at{p}
        (e_{j_1}, \ldots, e_{j_k}) \\
        &=
        \frac{1}{k!} \left(\SymD^E\right)^k s \at{p}
        (e_{i_1}, \ldots, e_{i_k}),
    \end{align*}
    using $\gamma_v(0) = p$ and $\dot{\gamma}_v(0) = v$.  
\end{proof}

Of course, the Taylor expansion of $\left( P_{\gamma_v, 0 \rightarrow
      1} \right)^{-1} s(\gamma_v(1))$ around $0$ needs not to converge
at all. In fact, the Borel Lemma, see e.g.
\cite[Remark~5.3.34]{waldmann:2007a}, shows that \emph{all} possible
numerical values appear as Taylor coefficients of smooth functions.
Nevertheless, we can use this proposition to obtain the \emph{formal}
Taylor series in a very nice way:
\begin{corollary}
    \label{corollary:formal-taylor-series-of-partrans}
    The formal Taylor series of the function $T_pM \ni v \mapsto
    \left( P_{\gamma_v, 0 \rightarrow 1} \right)^{-1} s(\gamma_v(1))
    \in E_p$ is given by
    \begin{equation}
        \label{eq:formal-taylor-series-of-partrans}
        \left( P_{\gamma_v, 0 \rightarrow 1} \right)^{-1}
        s(\gamma_v(1))
        \sim \mathcal{J}\left(\E^{\SymD^E} s\right)(v),
    \end{equation}
    where $\mathcal{J}: \bigoplus_{k=0}^\infty \Sym^k T_p^*M \tensor E_p
    \longrightarrow \Pol^\bullet(T_pM) \tensor E_p$ is the canonical
    isomorphism, extended to formal series in the symmetric and
    polynomial degree, respectively.
\end{corollary}
\begin{proof}
    This is now just a matter of computation. By
    Proposition~\ref{proposition:taylor-coefficients-of-partrans} we
    have in the sense of a formal series in $v$
    \begin{align*}
        \sum_{k=0}^\infty \frac{1}{k!} 
        \frac{\partial^k}{\partial v^{i_1} \cdots \partial v^{i_k}}
        \left(
            P_{\gamma_v, 0 \rightarrow 1}
        \right)^{-1}
        s(\gamma_v(1))\At{v=0} v^{i_1} \cdots v^{i_k}
        &=
        \sum_{k=0}^\infty \frac{1}{k!} \frac{1}{k!}
        \left(\SymD^E\right)^k s \At{p} (e_{i_1}, \ldots, e_{i_k})
        v^{i_1} \cdots v^{i_k} \\
        & =
        \sum_{k=0}^\infty \frac{1}{k!} \frac{1}{k!}
        \left(\SymD^E\right)^k s (v, \ldots, v) \\
        & =
        \sum_{k=0}^\infty \frac{1}{k!} 
        \mathcal{J}\left(\left(\SymD^E\right)^k s\right)(v) \\
        & =
        \mathcal{J}
        \left(
            \sum_{k=0}^\infty \frac{1}{k!} \left(\SymD^E\right)^k s
        \right) \\
        & =
        \mathcal{J} \left(\E^{\SymD^E} s \right)(v).
    \end{align*}
\end{proof}

In a more informal way one can say that the Taylor expansion of the
parallel transport along geodesics around initial velocity $0$ is
given by the exponential of the symmetrized covariant derivative.

We can specialize this statement to functions instead of general
sections. Here we simply have for $f \in \Cinfty(M)$
\begin{align*}
    \left(
        P_{\gamma_v, 0 \rightarrow 1}
    \right)^{-1}
    f(\gamma_v(1))
    =
    f(\gamma_v(1))
    =
    f( \exp_p(v) )
    = (\exp_p^*f) (v),
\end{align*}
since by Lemma~\ref{lemma:parallel-transport-tensor-product-pairing}
the parallel transport of numbers is trivial. Thus we obtain the
Taylor expansion of $\exp_p^*$ around $0$:
\begin{corollary}[Taylor expansion of $\exp_p^*$]
    \label{corollary:taylor-of-exp}%
    \index{Exponential map!Taylor expansion}%
    Let $V \subseteq T_pM$ be an open neighborhood of $0$ such that
    $\exp_p \at{V}$ is a diffeomorphism onto $U = \exp_p(V) \subseteq
    M$. Moreover, let $f \in \Cinfty(U)$. Then the formal Taylor
    series of $\exp_p^*f \in \Cinfty(V)$ around $0$ is given by
    \begin{equation}
        \label{eq:taylor-of-exp}
        \exp_p^* f \sim \mathcal{J} \left( \E^{\SymD} f \right).
    \end{equation}
    With other words, the Taylor expansion in normal coordinates
    around $p$ coincides with the Taylor expansion using $\SymD$.
\end{corollary}


%% file: jacobi.tex
%
%

In this section we consider not a single curve $\gamma$ in $M$ but
families of curves which are smoothly parametrized by an additional
variable. With other words, we consider smooth \emph{surfaces}
\begin{equation}
    \label{eq:smooth-surfaces}
    \sigma: \Sigma \longrightarrow M
\end{equation}
in $M$ where $\Sigma \subseteq \mathbb{R}^2$ is open. For convenience,
we mainly restrict to $\Sigma = I \times I'$ where $I, I' \subseteq
\mathbb{R}$ are open intervals. Hence $\Sigma$ is an open
rectangle. The two variables will be denoted by $(t, s) \in
\Sigma$. The canonical vector fields $\frac{\partial}{\partial t}$ and
$\frac{\partial}{\partial s}$ on $\Sigma$ give now rise to vector
fields
\begin{equation}
    \label{eq:canonical-vector-fields-on-surfaces}
    \dot{\sigma} = T \sigma \left( \frac{\partial}{\partial t} \right)
    \quad
    \textrm{and}
    \quad
    \sigma' = T \sigma \left( \frac{\partial}{\partial s} \right),
\end{equation}
which we can view as vector fields along $\sigma$, i.e. sections
\begin{equation}
    \label{eq:can-vecs-as-sections-along-sigma}
    \dot{\sigma}, \sigma' \in \Secinfty(\sigma^\# TM)
\end{equation}
of the pulled back tangent bundle. The first lemma gives a geometric
interpretation of the torsion of a covariant derivative. Note that for
$\nabla^\#$ there is no intrinsic definition of torsion possible.
\begin{lemma}
    \label{lemma:torsion-and-surfaces}%
    \index{Torsion}%
    Let $\nabla$ be a covariant derivative for $M$ and $\sigma: \Sigma
    \longrightarrow M$ a smooth surface. Then
    \begin{equation}
        \label{eq:torsion-and-canonical-vector-fields}
        \nabla^\#_{\frac{\partial}{\partial t}} \sigma'
        - \nabla^\#_{\frac{\partial}{\partial s}} \dot{\sigma}
        = \sigma^\# \Tor (\dot{\sigma}, \sigma').
    \end{equation}
    In particular, if $\nabla$ is torsion-free we have
    \begin{equation}
        \label{eq:canonical-vecs-for-torsion-free}
        \nabla^\#_{\frac{\partial}{\partial t}} \sigma'
        - \nabla^\#_{\frac{\partial}{\partial s}} \dot{\sigma}
        = 0.
    \end{equation}
\end{lemma}
\begin{proof}
    This is just a simple consequence of the definition of the
    pull-back connection $\nabla^\#$. If $(U, x)$ is a local chart we
    have
    \[
    \dot{\sigma}(t,s)
    = \frac{\partial \sigma^i}{\partial t} (t,s)
    \frac{\partial}{\partial x^i} \At{\sigma(t,s)}
    \quad
    \textrm{and}
    \quad
    \sigma'(t,s)
    = \frac{\partial \sigma^i}{\partial s} (t,s)
    \frac{\partial}{\partial x^i} \At{\sigma(t,s)},
    \]
    where $\sigma^i = x^i \circ \sigma$. Then
    \begin{align*}
        \nabla^\#_{\frac{\partial}{\partial t}} \sigma' \At{t,s}
        = \frac{\partial}{\partial t}
        \frac{\partial \sigma^i}{\partial s}(t,s)
        \frac{\partial}{\partial x^i}\At{\sigma(t,s)}
        + \frac{\partial \sigma^i}{\partial s}(t,s)
        \Gamma_{ji}^k (\sigma(t,s))
        \frac{\partial \sigma^j}{\partial s}(t,s)
        \frac{\partial}{\partial x^k} \At{\sigma(t,s)},
    \end{align*}
    and analogously for $\nabla^\#_{\frac{\partial}{\partial s}}
    \dot{\sigma}$.  From this the claim
    \eqref{eq:torsion-and-canonical-vector-fields} follows since
    $\Tor_{ij}^k = \Gamma_{ij}^k - \Gamma_{ji}^k$. But then
    \eqref{eq:canonical-vecs-for-torsion-free} is clear.
\end{proof}

\begin{lemma}
    \label{lemma:curvature-and-commutation-of-partial-derivatives}
    Let $\nabla^E$ be a covariant derivative for $E \longrightarrow M$
    and $\sigma: \Sigma \longrightarrow M$ a smooth surface in
    $M$. Then for $e \in \Secinfty(\sigma^\# E)$ we have
    \begin{equation}
        \label{eq:curvature-and-commutation-of-partial-derivatives}
        \nabla^\#_{\frac{\partial}{\partial t}}
        \nabla^\#_{\frac{\partial}{\partial s}} e
        -
        \nabla^\#_{\frac{\partial}{\partial s}}
        \nabla^\#_{\frac{\partial}{\partial t}} e
        =
        R^E \at{\sigma} (\dot{\sigma}, \sigma') e.
    \end{equation}
\end{lemma}
\begin{proof}
    This is just a particular case of the statement that the local
    curvature two-forms of $\nabla^E$ are the pull-backs of the local
    curvature two-forms of $\nabla$ together with $\left[
        \frac{\partial}{\partial t},\frac{\partial}{\partial s}
    \right] = 0$.
\end{proof}

We can now turn to Jacobi vector fields: they will turn out to be the
infinitesimal version of a family of geodesics. One defines for a yet
arbitrary curve a Jacobi vector field as follows:
\begin{definition}[Jacobi vector field]
    \label{definition:jacobi-vector-fields}%
    \index{Jacobi vector field}%
    Let $\gamma: I \subseteq \mathbb{R} \longrightarrow M$ be a smooth
    curve in $M$. Then a vector field $J \in \Secinfty(\gamma^\# TM)$
    is called Jacobi vector field along $\gamma$ if it satisfies the
    differential equation
    \begin{equation}
        \label{eq:jacobi-field-equation}
        \nabla^\#_{\frac{\partial}{\partial t}}
        \nabla^\#_{\frac{\partial}{\partial t}}
        J(t)
        =
        R_{\gamma(t)} ( \dot{\gamma}(t), J(t) ) \dot{\gamma}(t)
    \end{equation}
    for all $t \in I$.
\end{definition}
Up to now it is not necessary for $\gamma$ to be a geodesic, though
later on in most applications $\gamma$ will be a geodesic. We
investigate \eqref{eq:jacobi-field-equation} in a local chart $(U,
x)$. As usual we set $\gamma^i = x^i \circ \gamma$. Then we have for
\begin{equation}
    \label{eq:jacobi-field-in-local-coordinates}
    J(t) = J^i(t) \frac{\partial}{\partial x^i} \At{\gamma(t)}
\end{equation}
the first covariant derivative
\begin{equation}
    \label{eq:covariant-derivitive-of-jacobi}
    \nabla^\#_{\frac{\partial}{\partial t}} J
    = \frac{\D\! J^i}{\D\! t} \frac{\partial}{\partial x^i}
    + \Gamma_{ij}^k J^i \dot{\gamma}^j \frac{\partial}{\partial x^k}.
\end{equation}
Analogously, one computes the second covariant derivative
\begin{align}
    \nabla^\#_{\frac{\partial}{\partial t}}
    \nabla^\#_{\frac{\partial}{\partial t}}
    J
    &=
    \frac{\D^2\! J^i}{\D\! t^2} \frac{\partial}{\partial x^i}
    +
    \Gamma_{ij}^k \frac{\D\! J^i}{\D\! t} \dot{\gamma}^j
    \frac{\partial}{\partial x^k}
    +
    \frac{\D}{\D\! t} \left( \Gamma_{ij}^k \dot{\gamma}^j \right)
    J^i \frac{\partial}{\partial x^k} \nonumber \\
    &\quad+
    \Gamma_{ij}^k \frac{\D\! J^i}{\D\! t} \dot{\gamma}^j
    \frac{\partial}{\partial x^k}
    +
    \Gamma_{ij}^k J^i \dot{\gamma}^j \Gamma_{kl}^m
    \frac{\partial}{\partial x^m},
    \label{eq:jacobi-field-equation-local}
\end{align}
where always the data on $M$ has to be evaluated at $\gamma(t)$. On
the other hand we have for the right hand side of
\eqref{eq:jacobi-field-equation}
\begin{equation}
    \label{eq:jacobi-curvature-term-local}
    R(\dot{\gamma}, J) \dot{\gamma}
    = R^\ell{}_k{}_{ij} \dot{\gamma}^i J^j \dot{\gamma}^k
    \frac{\partial}{\partial x^\ell}.
\end{equation}
It follows that \eqref{eq:jacobi-field-equation} is locally a system
of linear second order differential equations for the coefficient
functions $J^i$ on $I \subseteq \mathbb{R}$ having the identity as
leading symbol and time-dependent coefficients for the first and
zeroth order terms. Thus we can apply the well-known theorems on
existence and uniqueness of solutions for such ordinary differential
equations:
\begin{proposition}
    \label{proposition:jacobi-vector-fields}
    Let $\gamma: I \subseteq \mathbb{R} \longrightarrow M$ be a smooth
    curve and $a \in I$. Then for every $v, w \in T_{\gamma(t)} M$
    there exists a unique Jacobi vector field $J_{v,w}$ along $\gamma$
    with
    \begin{equation}
        \label{eq:jacobi-vector-fields-inital-conditions}
        J_{v,w}(a) = v
        \quad
        \textrm{and}
        \quad
        \nabla^\#_{\frac{\partial}{\partial t}} J_{v,w} (a)
        = w.
    \end{equation}
    Moreover, the map
    \begin{equation}
        \label{eq:jadobi-depends-linear-on-initial-values}
        T_{\gamma(t)} M \oplus T_{\gamma(t)} M \ni (v,w)
        \; \mapsto \;
        J_{v,w} \in \Secinfty(\gamma^\# TM)
    \end{equation}
    is a linear injection.
\end{proposition}
\begin{proof}
    We cover the image of $\gamma$ by local charts. Then locally we
    have existence and uniqueness by the local form of
    \eqref{eq:jacobi-field-equation}. The uniqueness then guarantees
    that the local solutions patch together nicely on the overlaps of
    the charts. Then the linearity of
    \eqref{eq:jadobi-depends-linear-on-initial-values} is a
    consequence of the linearity of \eqref{eq:jacobi-field-equation}.
\end{proof}

Now we consider the particular case of a geodesic $\gamma(t) =
\exp_p(tv)$. In this case we can describe the Jacobi vector fields
with initial values $J(a) = 0$ explicitly as follows:
\begin{theorem}
    \label{theorem:jacobi-for-geodesic}
    Let $v,w \in T_pM$ and let $I \times I' \subseteq \mathbb{R}^2$ be
    a small enough open rectangle around $(0,0)$ such that
    \begin{equation}
        \label{eq:geodesic-surface}
        \sigma: I \times I' \ni (t,s)
        \; \mapsto \;
        \sigma(t,s) = \exp_p( t(v+sw) )
    \end{equation}
    is well-defined. Moreover, let $\gamma(t) = \sigma(t,0)$ be the
    geodesic with initial velocity $v$ at $p \in M$. Then
    \begin{equation}
        \label{eq:jacobi-for-geodesic}
        J(t) = \sigma'(t,0) \in \Secinfty(\gamma^\# TM)
    \end{equation}
    is the Jacobi vector field along $\gamma$ with initial values
    \begin{equation}
        \label{eq:initial-values-for-jadobi-geodesic}
        J(0) = 0
        \quad
        \textrm{and}
        \quad
        \nabla^\#_{\frac{\partial}{\partial t}} J (0)
        = w.
    \end{equation}
\end{theorem}
\begin{proof}
    First we notice that for small enough $I, I'$ around $0$ the map
    $\sigma$ is well-defined and hence $J$ is a smooth vector field
    along the geodesic $\gamma$. We compute by the chain rule
    \begin{align*}
        \sigma'(t,s)
        & = \frac{\partial}{\partial s} \exp_p( t(v+sw))
        =
        \left(
            T_{t(v+sw)} \exp_p
        \right)
        \left(
            \frac{\D}{\D\! s'} \At{s'=0}
            (s' \mapsto t(v+(s+s')w))
        \right) \\
        & =
        \left(
            T_{t(v+sw)} \exp_p
        \right)
        (tw)
        =
        t \left(
            T_{t(v+sw)} \exp_p
        \right)
        (w),
    \end{align*}
    where we have used the linearity of the tangent map and the
    canonical identification $T_{t(v+sw)} T_pM \simeq T_pM$ as
    usual. It follows that $J(0) = \sigma'(0,0) = 0$ is satisfied
    indeed. Moreover, we compute
    \[
    \nabla^\#_{\frac{\partial}{\partial t}} \sigma'(t,s)
    =
    \nabla^\#_{\frac{\partial}{\partial t}}
    \left(
        t (T_{t(v+sw)} \exp_p) (w)
    \right)
    =
    \left(T_{t(v+sw)} \exp_p\right) (w)
    +
    t \nabla^\#_{\frac{\partial}{\partial t}}
    \left(
        (T_{t(v+sw)} \exp_p) (w)
    \right),
    \]
    by the Leibniz rule for a covariant derivative. It follows that
    \[
    \nabla^\#_{\frac{\partial}{\partial t}} J (0)
    =
    \nabla^\#_{\frac{\partial}{\partial t}} \sigma'(t,0) \At{t=0}
    =
    (T_{t(v+sw)} \exp_p) (w) \At{t=s=0} + 0
    =
    T_0 \exp_p (w)
    =
    w,
    \]
    since $T_0 \exp_p = \id$. This shows that $J$ has the correct
    initial conditions
    \eqref{eq:initial-values-for-jadobi-geodesic}. Finally we compute
    \begin{align*}
       \nabla^\#_{\frac{\partial}{\partial t}}
       \nabla^\#_{\frac{\partial}{\partial t}}
       J(t)
       &=
       \nabla^\#_{\frac{\partial}{\partial t}}
       \nabla^\#_{\frac{\partial}{\partial t}}
       \sigma'(t,s) \At{s=0}
       =
       \nabla^\#_{\frac{\partial}{\partial t}}
       \nabla^\#_{\frac{\partial}{\partial s}}
       \dot{\sigma}(t,s) \At{s=0} \\
       &=
       \nabla^\#_{\frac{\partial}{\partial t}}
       \nabla^\#_{\frac{\partial}{\partial t}}
       \dot{\sigma}(t,s) \At{s=0}
       +
       R \At{\sigma(t,0)} (\dot{\sigma}(t,0), \sigma'(t,0)
       \dot{\sigma}(t,0),
   \end{align*}
   by Lemma~\ref{lemma:torsion-and-surfaces} and the torsion-freeness
   of $\nabla$ as well as by
   Lemma~\ref{lemma:curvature-and-commutation-of-partial-derivatives}.
   Now for all $s$ the curve $t \mapsto \sigma(t,s) = \exp_p(t(v+sw))$
   is a geodesic whence $\nabla^\#_{\frac{\partial}{\partial t}}
   \dot{\sigma}(t,s) = 0$ identically in $s$. This finally shows that
   the Jacobi equation, i.e. \eqref{eq:jacobi-field-equation}, is
   satisfied.
\end{proof}

\begin{corollary}
    \label{corollary:jacobi-formula-for-initial-conditions}
    Let $v,w \in T_pM$. Then
    \begin{equation}
        \label{eq:jacobi-formula-for-initial-conditions}
        J(t) = t \left( T_{tv} \exp_p \right) (w)
    \end{equation}
    is the unique Jacobi vector field along $\gamma(t) = \exp_p(tv)$
    with $J(0) = 0$ and $\nabla^\#_{\frac{\partial}{\partial t}} J(0)
    = w$.
\end{corollary}

By covariant differentiation of the Jacobi differential equation we
obtain the covariant derivatives of the Jacobi vector field up to all
orders, at least recursively. To this end, we first notice that the
right hand side of \eqref{eq:jacobi-field-equation} can be viewed as a
natural pairing of $\gamma^\# R \in \Secinfty(\gamma^\# \End{TM}
\tensor \Anti^2 TM)$ with $\dot{\gamma}, J \in \Secinfty(\gamma^\#
TM)$. Thus using the covariant derivative $\nabla^\#$ on all the
involved bundles gives
\begin{equation}
    \label{eq:differentiation-of-jacobi-equation-left-side}
    \nabla^\#_{\frac{\partial}{\partial t}}
    \left(
        \nabla^\#_{\frac{\partial}{\partial t}}
        \nabla^\#_{\frac{\partial}{\partial t}}
        J
    \right)
    =
    \nabla^\#_{\frac{\partial}{\partial t}}
    \left(
        \gamma^\# R (\dot{\gamma}, J) \dot{\gamma}
    \right)
    =
    \left(
        \nabla^\#_{\frac{\partial}{\partial t}}
        \gamma^\# R
    \right)
    (\dot{\gamma}, J) \dot{\gamma}
    +
    (\gamma^\# R) (\dot{\gamma},
    \nabla^\#_{\frac{\partial}{\partial t}} J )
    \dot{\gamma},
\end{equation}
since $\nabla^\#_{\frac{\partial}{\partial t}} \dot{\gamma} = 0$ for a
geodesic. Moreover,
\begin{equation}
    \label{eq:differentiation-of-jacobi-equation-right-side}
    \nabla^\#_{\frac{\partial}{\partial t}} \gamma^\# R \At{t}
    =
    \left(
        \nabla_{\dot{\gamma}(t)} R
    \right)
    \At{\gamma(t)}
\end{equation}
allows to compute the covariant derivatives of $\gamma^\# R$ in terms
of the covariant derivatives of $R$ on $M$. By iteration, the
successive use of the Leibniz rule of $\nabla^\#$ with respect to
natural pairings yields the following statement:
\begin{lemma}
    \label{lemma:exponential-cov-diffs-of-J}
    Let $J \in \Secinfty(\gamma^\# TM)$ be a Jacobi vector field along
    a geodesic $\gamma$. Then
    \begin{equation}
        \label{eq:14}
        \exp\left(
            \lambda \nabla^\#_{\frac{\partial}{\partial t}}
        \right)
        \nabla^\#_{\frac{\partial}{\partial t}}
        \nabla^\#_{\frac{\partial}{\partial t}}
        J
        =
        \left(
            \exp\left(
                \lambda \nabla^\#_{\frac{\partial}{\partial t}}
            \right)
            \gamma^\# R
        \right)
        (\dot{\gamma},J) \dot{\gamma}
        +
        (\gamma^\# R)
        \left(
            \dot{\gamma},
            \exp\left(
                \lambda \nabla^\#_{\frac{\partial}{\partial t}}
            \right)
        \right)
        \dot{\gamma}
    \end{equation}
    in the sense of a formal power series in the formal parameter
    $\lambda$.
\end{lemma}
\begin{proof}
    Either this is shown by differentiating both sides with respect to
    $\lambda$ and observing that the resulting differential equations
    coincide thanks to
    \eqref{eq:differentiation-of-jacobi-equation-left-side}, or by
    induction in the summation parameter of the exponential series.
\end{proof}

\begin{remark}
    \label{remark:recursive-computation-of-cov-diffs-of-J}
    The lemma can be used to efficiently compute $(
    \nabla^\#_{\frac{\partial}{\partial t}} )^k J$ at $t=0$ for $k \in
    \mathbb{N}_0$. Indeed, it provides a recursion scheme giving
    \begin{equation}
        \label{eq:second-cov-diff-of-J}
        \left( \nabla^\#_{\frac{\partial}{\partial t}} \right)^2 J(0)
        = R(v, J(0) ) = 0,
    \end{equation}
    \begin{equation}
        \label{eq:third-cov-diff-of-J}
        \left( \nabla^\#_{\frac{\partial}{\partial t}} \right)^3 J(0)
        =
        \left(
            \nabla^\#_{\frac{\partial}{\partial t}} \gamma^\# R
        \right)
        (v, J(0) ) v
        +
        R \left(
            v, \nabla^\#_{\frac{\partial}{\partial t}} J(0)
        \right) v
        =
        0 + R(v,w) v,
    \end{equation}
    since $\dot{\gamma}(0) = v$ and $J(0)=0$ as well as
    $\nabla^\#_{\frac{\partial}{\partial t}} J(0) = w$. The next terms
    are
    \begin{align}
        \left( \nabla^\#_{\frac{\partial}{\partial t}} \right)^4
        J(0)
        &=
        \left(
            \left(
                \nabla^\#_{\frac{\partial}{\partial t}}
            \right)^2
            \gamma^\# R
        \right) (v,0) v
        +
        2 \left(
            \nabla^\#_{\frac{\partial}{\partial t}} \gamma^\# R
        \right) (v,w) v
        +
        R \left(
            v, \left(
                \nabla^\#_{\frac{\partial}{\partial t}}
            \right)^2 J(0)
        \right) v \nonumber \\
        &=
        2 \left(
            \nabla^\#_{\frac{\partial}{\partial t}}
            \gamma^\# R
        \right) (v,w) v          
        \label{eq:fourth-cov-diff-of-J}
    \end{align}
    and
    \begin{align}
        \left( \nabla^\#_{\frac{\partial}{\partial t}} \right)^5
        J(0)
        & =
        \left(
            \left( \nabla^\#_{\frac{\partial}{\partial t}} \right)^3 
            \gamma^\# R
        \right) (v,0) v
        +
        3 \left(
            \left( \nabla^\#_{\frac{\partial}{\partial t}} \right)^2 
            \gamma^\# R
        \right)
        \left(
            v, \nabla^\#_{\frac{\partial}{\partial t}} J(0)
        \right) v \nonumber \\
        &\quad+
        3 \left(
            \nabla^\#_{\frac{\partial}{\partial t}} \gamma^\# R
        \right)
        \left(
            v, \left(
                \nabla^\#_{\frac{\partial}{\partial t}} \right)^2
            J(0) 
        \right) v
        +
        R \left(
            v, \left(
                \nabla^\#_{\frac{\partial}{\partial t}}
            \right)^3 J(0)
        \right) v \nonumber \\
        &=
        3 \left(
            \left( \nabla^\#_{\frac{\partial}{\partial t}} \right)^2 
            \gamma^\# R
        \right)
        (v,w) v
        +
        3 R(v, R(v,w) v) v,
        \label{eq:fifth-cov-diff-of-J}
    \end{align}
    using successively those computations done in lower orders.
    Moreover, an easy induction shows that $(
    \nabla^\#_{\frac{\partial}{\partial t}} )^k J_w (0)$ is a
    homogeneous polynomial in $v$ of order $k-1$ and linear in $w$.
    Here one uses that $\nabla^\#_{\frac{\partial}{\partial t}}$ of an
    arbitrary tensor field $\gamma^\# T$ is linear in $v$.
\end{remark}

We can use this to compute the Taylor expansion of the tangent map of
the exponential map $\exp_p$. For any $v \in T_pM$ the tangent map
$T_v \exp_p$ is a linear map $T_v \exp_p: T_pM \longrightarrow
T_{\exp_p(v)} M$. In order to compute its Taylor expansion around
$v=0$ we first have to identify $T_{\exp_p(v)} M$ with $T_pM$ again by
using the parallel transport $P_{\gamma_v, 0 \rightarrow 1}: T_pM
\longrightarrow T_{\exp_p(v)}M$ along the geodesic $t \mapsto
\gamma_v(t) = \exp_p(tv)$. This way we obtain a linear map
\begin{equation}
    \label{eq:Texp-with-partrans-identification}
    \left(
        P_{\gamma_v, 0 \rightarrow 1}
    \right)^{-1}
    \circ T_v \exp_p :
    T_pM \longrightarrow T_pM
\end{equation}
for every $v \in T_pM$ small enough. We want to compute now the Taylor
coefficients of
\begin{equation}
    \label{eq:17}
    T_pM \ni v
    \; \mapsto \;
    \left( P_{\gamma_v, 0 \rightarrow 1} \right)^{-1} \circ T_v \exp_p
    \in \End(T_pM)
\end{equation}
around $v=0$. To do so we evaluate the endomorphism on a fixed vector
$w \in T_pM$ and consider the map
\begin{equation}
    \label{eq:evaulation-of-Texp-with-partrans-identification}
    v
    \; \mapsto \;
    \left( P_{\gamma_v, 0 \rightarrow 1} \right)^{-1} \circ T_v \exp_p
    (w).
\end{equation}
In order to compute the partial derivatives of
\eqref{eq:evaulation-of-Texp-with-partrans-identification} in the
$v$-variable it suffices to consider the derivatives of the map
\begin{equation}
    \label{eq:backtracing-to-t-derivatives}
    t
    \; \mapsto \;
    \left( P_{\gamma_v, 0 \rightarrow 1} \right)^{-1}
    \circ T_{tv} \exp_p (w)
\end{equation}
around $t=0$ instead and use Lemma~\ref{lemma:derivative-formula}
afterwards. Since $T_{tv} \exp_p (w) = \frac{1}{t} J_w(t)$ is a
multiple of the unique Jacobi vector field $J_w \in
\Secinfty(\gamma_v^\# TM)$ along $\gamma_v$ with $J_w(0) = 0$ and
$\nabla^\#_{\frac{\partial}{\partial t}} J_w(0) = w$ we can compute
its covariant derivatives by means of
Lemma~\ref{lemma:exponential-cov-diffs-of-J} and
Remark~\ref{remark:recursive-computation-of-cov-diffs-of-J}
recursively. Finally, we note that
\begin{equation}
    \label{eq:partrans-and-rescaling}
    P_{\gamma_{tv}, 0 \rightarrow 1} = P_{\gamma_v, 0 \rightarrow t},
\end{equation}
whence we have to consider the map
\begin{equation}
    \label{eq:what-we-finally-want-to-differentiate}
    t \mapsto
    \left( P_{\gamma_v, 0 \rightarrow t} \right)^{-1}
    \left( \frac{1}{t} J_w(t) \right),
\end{equation}
of which we want to compute the Taylor coefficients around
$t=0$. Collecting things we obtain the following result:
\begin{theorem}[Taylor coefficients of $T\exp_p$]
    \label{theorem:taylor-expansion-of-exp}%
    Let $p \in M$ and $v,w \in T_pM$. Then for all $k \in
    \mathbb{N}_0$ we have
    \begin{equation}
        \label{eq:partial-derivatives-of-exp}
        \frac{\partial^k}{\partial v^{i_1} \cdots \partial v^{i_k}}
        \At{v=0}
        \left( P_{\gamma_v, 0 \rightarrow 1} \right)^{-1}
        \circ T_v \exp_p (w)
        =
        \frac{1}{(k+1)!}
        \frac{\partial^k}{\partial v^{i_1} \cdots \partial v^{i_k}}
        \underbrace{
          \nabla^\#_{\frac{\partial}{\partial t}}
          \cdots
          \nabla^\#_{\frac{\partial}{\partial t}}
        }_{k+1 \; \textrm{times}}
        J_w(t) \At{t=0}.
    \end{equation}
    The first terms of the (formal) Taylor expansion around $v=0$ are
    therefore given by
    \begin{equation}
        \label{eq:first-terms-of-taylor-of-exp}
        \left( P_{\gamma_v, 0 \rightarrow 1} \right)^{-1}
        \circ T_v \exp_p (w)
        =
        w + \frac{1}{6} R_p(v,w)v
        + \frac{1}{12} ( \nabla_v R )_p (v,w) v
        + \cdots.
    \end{equation}
\end{theorem}
\begin{proof}
    By Corollary~\ref{corollary:jacobi-formula-for-initial-conditions}
    we have $t T_{tv}\exp_p(w) = J_w(t)$ whence we can compute the
    $\nabla^\#_{\frac{\partial}{\partial t}}$-derivatives of the
    vector field $t \mapsto T_{tv}\exp_p(w)$ at $t=0$ as follows. By
    the Leibniz rule we have
    \begin{align*}
        \underbrace{
          \nabla^\#_{\frac{\partial}{\partial t}}
          \cdots
          \nabla^\#_{\frac{\partial}{\partial t}}
        }_{k \; \textrm{times}}
        J_w(t) \At{t=0}
        &=
        \underbrace{
          \nabla^\#_{\frac{\partial}{\partial t}}
          \cdots \nabla^\#_{\frac{\partial}{\partial t}}
        }_{k \; \textrm{times}}
        (t T_{tv} \exp_p(w)) \At{t=0}
        \\
        &=
        t T_{tv} \exp_p(w) \At{t=0}
        + k 
        \underbrace{
          \nabla^\#_{\frac{\partial}{\partial t}}
          \cdots 
          \nabla^\#_{\frac{\partial}{\partial t}}
        }_{k-1 \; \textrm{times}}
        T_{tv}\exp_p(w) \At{t=0}
        + 0,
    \end{align*}
    whence for $k \geq 1 $ we get
    \[
    \underbrace{
      \nabla^\#_{\frac{\partial}{\partial t}}
      \cdots \nabla^\#_{\frac{\partial}{\partial t}}
    }_{k-1 \; \textrm{times}}
    (t T_{tv} \exp_p(w)) \At{t=0}
    =
    \frac{1}{k}
    \underbrace{
      \nabla^\#_{\frac{\partial}{\partial t}}
      \cdots \nabla^\#_{\frac{\partial}{\partial t}}
    }_{k \; \textrm{times}}
    J_w(t) \At{t=0}.
    \tag{$*$}
    \]
    Now the right hand side is recursively computable by
    Lemma~\ref{lemma:exponential-cov-diffs-of-J}, see
    Remark~\ref{remark:recursive-computation-of-cov-diffs-of-J} for
    the first terms. We can collect the results and obtain
    \begin{align*}
        \frac{\partial^k}{\partial v^{i_1} \cdots \partial v^{i_k}}
        \At{v=0} 
        &
        \left(
            \left( P_{\gamma_v, 0 \rightarrow 1} \right)^{-1}
            \circ T_v \exp_p
        \right)(w) \\        
        &\stackrel{\mathclap{\eqref{eq:derivative-formula}}}{=}
        \quad
        \frac{1}{k!}
        \frac{\partial^k}{\partial v^{i_1} \cdots \partial v^{i_k}}
        \frac{\D^k}{\D\! t^k} \At{t=0}
        \left(
            \left( P_{\gamma_{tv}, 0 \rightarrow 1} \right)^{-1}
            \circ T_{tv} \exp_p
        \right)(w) \\
        &\stackrel{\mathclap{\eqref{eq:partrans-and-rescaling}}}{=}
        \quad
        \frac{1}{k!}
        \frac{\partial^k}{\partial v^{i_1} \cdots \partial v^{i_k}}
        \frac{\D^k}{\D\! t^k} \At{t=0}
        \left( P_{\gamma_{v}, 0 \rightarrow t} \right)^{-1}
        \left(
            T_{tv} \exp_p(w)
        \right) \\
        &\stackrel{\mathclap{\eqref{eq:time-derivative-of-parallel-transport}}}{=}
        \quad
        \frac{1}{k!}
        \frac{\partial^k}{\partial v^{i_1} \cdots \partial v^{i_k}}
        \left( P_{\gamma_{v}, 0 \rightarrow t} \right)^{-1}
        \underbrace{
          \nabla^\#_{\frac{\partial}{\partial t}}
          \cdots \nabla^\#_{\frac{\partial}{\partial t}}
        }_{k \; \textrm{times}}
        T_{tv} \exp_p (w) \At{t=0} \\
        &\stackrel{\mathclap{(*)}}{=}
        \frac{1}{k!}
        \frac{\partial^k}{\partial v^{i_1} \cdots \partial v^{i_k}}
        \frac{1}{k+1}
        \underbrace{
          \nabla^\#_{\frac{\partial}{\partial t}}
          \cdots \nabla^\#_{\frac{\partial}{\partial t}}
        }_{k+1 \; \textrm{times}}
        J_w(t) \At{t=0},
    \end{align*}
    which shows \eqref{eq:partial-derivatives-of-exp}. As we know from
    Remark~\ref{remark:recursive-computation-of-cov-diffs-of-J}, the
    $(k+1)$-st covariant derivative of $J_w$ at $0$ is a homogeneous
    polynomial in $v$ of order $k$. This is also clear from the proof
    of Lemma~\ref{lemma:derivative-formula}. Now we compute the first
    orders of the Taylor expansion explicitly. Since we already know
    $T_0 \exp_p = \id$ the zeroth order is given as in
    \eqref{eq:first-terms-of-taylor-of-exp}. In fact, this was used to
    show $\nabla^\#_{\frac{\partial}{\partial t}} J_w(0) = w$. For the
    first order $k=1$ we get
    \[
    \frac{\partial}{\partial v^i}
    \left( P_{\gamma_{v}, 0 \rightarrow 1} \right)^{-1}
    T_v \exp_p(w) \At{v=0}
    =
    \frac{\partial}{\partial v^i}
    \frac{1}{2}
    \nabla^\#_{\frac{\partial}{\partial t}}
    \nabla^\#_{\frac{\partial}{\partial t}}
    J_w(0)
    =
    0
    \]
    by \eqref{eq:second-cov-diff-of-J}. The next order gives
    \begin{align*}
        \frac{\partial^2}{\partial v^i \partial v^j}
        \left( P_{\gamma_{v}, 0 \rightarrow 1} \right)^{-1}
        \circ T_v \exp_p(w) \At{v=0}
        &=
        \frac{\partial^2}{\partial v^i \partial v^j}
        \frac{1}{6}
        \left(
            \nabla^\#_{\frac{\partial}{\partial t}}
        \right)^3 J_w(0) \\
        &=
        \frac{1}{6} \frac{\partial^2}{\partial v^i \partial v^j}
        R_p(v,w)v \\
        &=
        \frac{1}{6}
        \left(
            R_p( \frac{\partial}{\partial x^j},w )
            \frac{\partial}{\partial x^i}
            +
            R_p( \frac{\partial}{\partial x^i},w )
            \frac{\partial}{\partial x^j}
        \right).
    \end{align*}
    Thus
    \[
    \frac{1}{2!}
    \frac{\partial^2}{\partial v^i \partial v^j}
    \left( P_{\gamma_{v}, 0 \rightarrow 1} \right)^{-1}
    \circ T_v \exp_p(w) v^i v^j
    =
    \frac{1}{6} R(v,w)v,
    \]
    explaining the quadratic term in
    \eqref{eq:first-terms-of-taylor-of-exp}. The cubic term is
    obtained from \eqref{eq:fourth-cov-diff-of-J}
    \begin{align*}
        \frac{\partial^3}{\partial v^i \partial v^j \partial v^k}
        \left( P_{\gamma_{v}, 0 \rightarrow 1} \right)^{-1}
        \circ T_v \exp_p(w) \At{v=0}
        & =
        \frac{1}{4!}
        \frac{\partial^3}{\partial v^i \partial v^j \partial v^k}
        2 \left(
            \nabla^\#_{\frac{\partial}{\partial t}} \gamma^\#R
        \right) (v,w)v \\
        & =
        \frac{1}{4!}
        \frac{\partial^3}{\partial v^i \partial v^j \partial v^k}
        2 \left( \nabla_v R \right)_p (v,w)v,
    \end{align*}
    from which we get
    \begin{align*}
        \frac{1}{3!}
        \frac{\partial^3}{\partial v^i \partial v^j \partial v^k}
        \left( P_{\gamma_{v}, 0 \rightarrow 1} \right)^{-1}
        \circ T_v \exp_p(w) \At{v=0} v^i v^j v^k
        =
        \frac{1}{12} \left( \nabla_v R \right)_p (v,w)v,
    \end{align*}
    as claimed in \eqref{eq:first-terms-of-taylor-of-exp}.
\end{proof}

\begin{remark}
    \label{remark:taylor-expansion-of-exp}
    More symbolically we can write the (formal) Taylor expansion of
    the tangent map of $\exp_p$ as
    \begin{align*}
        \left( P_{\gamma_{v}, 0 \rightarrow 1} \right)^{-1}
        \circ T_v \exp_p(w)
        & \sim_{v \rightarrow 0}
        \sum_{k=0}^\infty \frac{1}{k!}
        \frac{\partial^k}{\partial v^{i_1} \cdots \partial v^{i_k}}
        \left(
            \left( P_{\gamma_v, 0 \rightarrow 1} \right)^{-1}
            \circ T_v \exp_p
        \right)(w) \At{v=0} v^{i_1} \cdots v^{i_k} \\
        & =
        \sum_{k=0}^\infty \frac{1}{k!}
        \frac{\partial^k}{\partial v^{i_1} \cdots \partial v^{i_k}}
        \left(
            \frac{1}{(k+1)!}
            \underbrace{
              \nabla^\#_{\frac{\partial}{\partial t}}
              \cdots \nabla^\#_{\frac{\partial}{\partial t}}
            }_{k+1 \; \textrm{times}}
            J_w(t) \At{t=0}
        \right) v^{i_1} \cdots v^{i_k} \\
        & = 
        \sum_{k=0}^\infty \frac{1}{(k+1)!}
        \underbrace{
          \nabla^\#_{\frac{\partial}{\partial t}}
          \cdots \nabla^\#_{\frac{\partial}{\partial t}}
        }_{k+1 \; \textrm{times}}
        J_w(t) \At{t=0} \\
        & =
        \exp\left(\nabla^\#_{\frac{\partial}{\partial t}}\right)
        J_w(t) \At{t=0},
    \end{align*}
    since on one hand $ \nabla^\#_{\frac{\partial}{\partial t}} \cdots
    \nabla^\#_{\frac{\partial}{\partial t}} J_w(t)$ is a homogeneous
    polynomial in $v$ of degree $k$ for $k+1$ derivatives and since
    the zeroth term of the exponential series does not contribute due
    to $J_w(0) = 0$. Of course the formula
    \begin{equation}
        \label{eq:FormalTaylor-of-exp}
        \left( P_{\gamma_{v}, 0 \rightarrow 1} \right)^{-1}
        \circ T_v \exp_p(w)
        \sim_{v \rightarrow 0}
        \exp\left(\nabla^\#_{\frac{\partial}{\partial t}}\right)
        J_w(t) \At{t=0} 
    \end{equation}
    is only the formal Taylor expansion: in general, the right hand
    side will not converge in any reasonable sense. Note however that
    the combinatorics to compute the covariant derivatives of $J_w$ at
    $t=0$ is fairly simple and given by universal polynomials in the
    curvature and its covariant derivatives at $p$.
\end{remark}

As a last application of our investigations of Jacobi vector fields we
specialize to the case of a semi-Riemannian manifold $(M,g)$ and the
Levi-Civita connection $\nabla$. Then one has the following result,
known as the \emph{Gauss Lemma}:
\begin{proposition}[Gauss Lemma]
    \label{proposition:gauss-lemma}%
    \index{Gauss Lemma}%
    Let $(M,g)$ be a semi-Riemannian manifold and $p \in M$. Then for
    $v,w \in T_pM$ we have
    \begin{equation}
        \label{eq:gauss-lemma}
        g_{\exp_p(v)}
        \left(
            T_v \exp_p(v), T_v \exp_p(w)
        \right)
        =
        g_p(v,w),
    \end{equation}
    whenever $v$ is still in the domain of $\exp_p$.
\end{proposition}
\begin{proof}
    We consider the surface $\sigma(t,s) = \exp_p(t(v+sw))$ which is
    defined for $t \in [0,1]$ and $s$ small enough. Then we have
    \[
    \dot{\sigma}(t,s) = T_{\exp_p(t(v+sw))} \exp_p(v+sw)
    \quad
    \textrm{and}
    \quad
    \sigma'(t,s) = T_{\exp_p(t(v+sw))} \exp_p(tw)
    \]
    by the chain rule as we computed already in the proof of
    Theorem~\ref{theorem:jacobi-for-geodesic}. Thus we have to compute
    $g_{\exp_p(v)} \left( \dot{\sigma}(1,0), \sigma'(1,0) \right)$.
    We consider the geodesic $t \mapsto \exp_p(t(v+sw)) = \gamma_s(t)$
    with initial velocity vector $v+sw$. First we note by
    $\nabla^\#_{\frac{\partial}{\partial t}} \dot{\gamma}_s = 0$ that
    \[
    \frac{\partial}{\partial t}
    g_{\gamma_s (t)} \left(
        \dot{\gamma}_s(t), \dot{\gamma}_s(t) 
    \right)
    =
    2 g_{\gamma_s (t)}
    \left(
        \nabla^\#_{\frac{\partial}{\partial t}} \dot{\gamma}_s(t),
        \dot{\gamma}_s(t)
    \right) 
    =
    0,
    \tag{$*$}
    \]
    by the fact that $g$ is covariantly constant. It follows that
    $g_{\gamma_s (t)} (\dot{\gamma}_s(t), \dot{\gamma}_s(t) ) =
    g_p(v+sw,v+sw)$. In fact, this is the Gauss Lemma for $w=v$. To
    proceed we compute
    \begin{align*}
        \frac{\partial}{\partial t}
        g_{\sigma(t,s)} ( \dot{\sigma}(t,s), \sigma'(t,s) )
        & =
        g_{\sigma(t,s)}
        \left(
            \nabla^\#_{\frac{\partial}{\partial t}}
            \dot{\sigma}(t,s), 
            \sigma'(t,s)
        \right)
        +
        g_{\sigma(t,s)}
        \left(
            \dot{\sigma}(t,s), 
            \nabla^\#_{\frac{\partial}{\partial t}} \sigma'(t,s)
        \right) \\
        & =
        0
        +
        g_{\sigma(t,s)}
        \left(
            \dot{\sigma}(t,s), 
            \nabla^\#_{\frac{\partial}{\partial t}}
            \dot{\sigma}(t,s)
        \right) \\
        & =
        \frac{1}{2} \frac{\partial}{\partial s}
        \left(
            g_{\sigma(t,s)}
            \left(
                \dot{\sigma}(t,s), \dot{\sigma}(t,s)
            \right)
        \right),
    \end{align*}
    using the fact that $\nabla$ is torsion-free, see
    Lemma~\ref{lemma:torsion-and-surfaces}. Since all curves $t
    \mapsto \sigma(t,s) = \gamma_s(t)$ are geodesics we know that
    \[
    g_{\sigma(t,s)} ( \dot{\sigma}(t,s), \dot{\sigma}(t,s) )
    =
    g_{\sigma(0,s)} ( \dot{\sigma}(0,s), \dot{\sigma}(0,s) )
    +
    g_p( v+sw, v+sw ),
    \]
    whence
    \[
    \frac{1}{2} \frac{\partial}{\partial s}
    \left(
        g_{\sigma(t,s)}
        \left(
            \dot{\sigma}(t,s), \dot{\sigma}(t,s)
        \right)
    \right)
    =
    g_p(v,w) + s g_p(w,w).
    \]
    Putting things together we have for $s=0$ and all $t$
    \[
    \frac{\partial}{\partial t}
    g_{\sigma(t,0)} ( \dot{\sigma}(t,0), \sigma'(t,0) )
    =
    g_p(w,w)
    \]
    independent of $t$. Hence we conclude $g_{\sigma(t,0)} (
    \dot{\sigma}(t,0), \sigma'(t,0) ) = t g_p(v,w)$ and setting $t=1$
    gives the desired result \eqref{eq:gauss-lemma}.
\end{proof}

\begin{remark}[Gauss Lemma]
    \label{remark:gauss-lemma}
    The geometric interpretation of the Gauss Lemma is two-fold. For
    $v=w$ we see that the length-square of the tangent vector of a
    geodesic is constant. In the Riemannian setting this simply means
    that the length itself stays constant whence geodesics are curves
    with ``constant velocity''. In the Hamiltonian picture, this part
    of the Gauss Lemma can be interpreted as energy conservation under
    the Hamiltonian time evolution, see e.g. \cite[Aufgabe~3.10,
    vii.)]{waldmann:2007a} for this point of view. The case with
    arbitrary $w$ means that along a geodesic at least the ``angles''
    with respect to the tangent vector of the geodesic are preserved.
\end{remark}


%% file: jacobiexp.tex
%
%

Now we will use the formal Taylor expansion of $T_v \exp_p$ around
$v=0$ to consider the following problem. Given a positive density $\mu
\in \Secinfty(\Dichten T^*M)$ on $M$, i.e. $\mu > 0$ everywhere, we
can compare the constant density $\mu_p$ on $T_pM$ with $\mu$ via the
exponential map $\exp_p$ of $\nabla$. More precisely, we consider an
open neighborhood of the zero section such that
\begin{equation}
    \label{eq:projection-times-exp}
    \pi \times \exp : V \subseteq TM \longrightarrow U \subseteq M
    \times M
\end{equation}
is a diffeomorphism onto its image, denoted by $U$. In fact, $U$ is an
open neighborhood of the diagonal since $(\pi \times \exp) (0_p) =
(p,p)$ for $0_p \in T_pM$.
\begin{definition}
    \label{definition:density-comparison-function}
    Let $\mu \in \Secinfty(\Dichten T^*M)$ be a positive density on
    $M$. Then the function $\rho: U \longrightarrow \mathbb{R}$ is
    defined by
    \begin{equation}
        \label{eq:density-comparison-function}
        \rho(p,q) (\exp_p{}_* \mu_p)_q = \mu _q
    \end{equation}
    for $(p,q) \in U$.
\end{definition}
\begin{lemma}
    \label{lemma:density-comparison-function}
    Let $\mu \in \Secinfty(\Dichten T^*M)$ by a positive density. Then
    $\rho \in \Cinfty(U)$ and $\rho > 0$.
\end{lemma}
\begin{proof}
    We have
    \[
    \rho(p,q) = \frac{\mu_q}{(\exp_p{}_* \mu_p)\at{q}} > 0.
    \]
    Moreover, the map $(p,q) \mapsto \exp_{p*} \mu_p \at{q}$ is a
    smooth map on $U$ with values in $\Dichten T^*_qM$.  Since at
    every point $(p,q)$ the value is a \emph{positive} density the
    quotient is well-defined and smooth.
\end{proof}

\begin{remark}
    \label{remark:density-comparison-function}
    Geometrically speaking, the function $\rho$ measures how much the
    density $\mu$ at $q$ differs from the density $\mu$ at $p$ when
    the latter is moved to $q$ by means of the exponential map. Thus
    $\rho$ encodes the change of volume as one moves around in $M$.
    Note that $\rho$ is \emph{not} symmetric.
\end{remark}

Sometimes we fix a reference point $p \in M$ and consider the function
$\rho_p : U_p \longrightarrow \mathbb{R}$ defined by
\begin{equation}
    \label{eq:densitiy-comparison-function-p-fixed}
    \rho_p(q) = \rho(p, q)
\end{equation}
for $q \in U_p \subseteq M$ where $U_p$ is an open neighborhood on
which we have normal coordinates, i.e. $U_p = \exp_p(V_p)$ with $V_p =
V \cap T_pM$. Thus we have
\begin{equation}
    \label{eq:fixed-comparison-function-times-pushed-fixed-density}
    \rho_p \exp_{p*} \mu_p = \mu
\end{equation}
on $U_p$. Moreover, it will also be convenient to compare the
densities on the tangent space of $p$ and not on $M$. Thus one defines
the function $\widetilde{\rho}: V \longrightarrow \mathbb{R}$ by
\begin{equation}
    \label{eq:density-comparison-function-on-TpM}
    \widetilde{\rho}(v_p) \mu_p = (\exp_p^* \mu) (v_p).
\end{equation}
Thus $\widetilde{\rho}$ is the prefactor of the \emph{constant
  density} $\mu_p$ on $T_pM$ such that we obtain the pull-back of
$\mu$. Clearly, we have
\begin{equation}
    \label{eq:fixed-comparison-and-comparison}
    \widetilde{\rho}(v_p) = \rho (p, \exp_p(v_p)),
\end{equation}
whence also $\widetilde{\rho} = \rho \circ (\pi \times \exp) \in
\Cinfty(V)$ is smooth and positive. Again, we write
$\widetilde{\rho}_p \in \Cinfty(V_p)$ for the restriction of
$\widetilde{\rho}$ to a particular tangent space of a fixed reference
point $p \in M$.

The aim is now to compute the (formal) Taylor expansion of
$\widetilde{\rho}_p$ around $v=0$ which is equivalent to the (formal)
Taylor expansion of $\rho_p$ in normal coordinates around $p$. To this
end, we first give another interpretation of $\rho$ and
$\widetilde{\rho}$. In fact, we have \emph{two} aspects of comparing
the volumes. On one hand, the density $\mu$ is not ``constant'' along
$M$ since there is simply no intrinsic way to formulate such a
statement. On the other hand, the exponential map needs not to be
volume preserving. We try to separate these two effects as follows:
Using the unique geodesic $t \mapsto \exp_p(tv)$ from $p$ to $q =
\exp_p(v)$ we can parallel transport $\mu_q$ back to $p$ using the
parallel transport induced by $\nabla$ on the density bundle. This
gives a \emph{constant} density $\left( P_{\gamma_v, 0 \rightarrow 1}
\right)^{-1} \mu_{\exp_p(v)} \in \Dichten T^*_pM$ on $T_pM$ for every
$v \in V_p$. Thus this will be a \emph{constant} multiple of $\mu_p$
depending parametrically on $v$. This $v$-dependence measures how much
$\mu$ is \emph{not parallel} with respect to $\nabla$. Secondly, we
consider the tangent map
\begin{equation}
    \label{eq:tangent-map-of-exp}
    T_v \exp_p : T_pM \longrightarrow T_{\exp_p(v)} M,
\end{equation}
and want to determine its change of volume features. Since source and
target are different vector spaces there is no way to define a
``determinant'' of this linear map, we first have to take care that we
get a map from a tangent space into the \emph{same} tangent space.
Thus we consider
\begin{equation}
    \label{eq:par-circ-exp}
    \left( P_{\gamma_v, 0 \rightarrow 1} \right)^{-1}
    \circ T_v \exp_p : T_pM \longrightarrow T_pM
\end{equation}
instead.

Combining the effects we use the density $\left( P_{\gamma_v, 0
      \rightarrow 1} \right)^{-1} \mu_{\exp_p(v)}$ and evaluate on a
basis $e_1, \ldots, e_n \in T_pM$ after applying $\left( P_{\gamma_v,
      0 \rightarrow 1} \right)^{-1} \circ T_v \exp_p$ to it. In order
to get a result which is independent of the chosen basis we normalize
it by $\mu_p(e_1, \ldots, e_n)$, i.e. we consider the quantity
\begin{equation}
    \label{eq:normalized-partrans-density}
    \frac{1}{\mu_p(e_1, \ldots, e_n)}
    \left(
        \left( P_{\gamma_v, 0 \rightarrow 1} \right)^{-1}
        (\mu_{\exp_p(v)})
    \right)
    \left(
        P_{\gamma_v, 0 \rightarrow 1} ^{-1} \circ T_v\exp_p (e_1),
        \ldots,
        P_{\gamma_v, 0 \rightarrow 1}^{-1} \circ T_v\exp_p (e_n)
    \right)
\end{equation}
for $v \in V_p \subseteq T_pM$. Then we have the following statement:
\begin{lemma}
    \label{lemma:comp-func-on-TpM-is-norm-partrans-density}
    For $p \in M$ and $v \in V_p \subseteq T_pM$ we have
    \begin{equation}
        \label{eq:comp-func-on-TpM-is-norm-partrans-density}
        \widetilde{\rho}(p) =
        \frac{1}{\mu_p(e_1, \ldots, e_n)}
        \left(
            \left( P_{\gamma_v, 0 \rightarrow 1} \right)^{-1}
            (\mu_{\exp_p(v)})
        \right)
        \left(
            P_{\gamma_v, 0 \rightarrow 1} ^{-1} \circ T_v\exp_p (e_1), 
            \ldots,
            P_{\gamma_v, 0 \rightarrow 1}^{-1} \circ T_v\exp_p (e_n)
        \right),
    \end{equation}
    where $e_1, \ldots, e_n \in T_pM$ is a basis.
\end{lemma}
\begin{proof}
    Using Lemma~\ref{lemma:parallel-transport-of-densities} we compute
    \begin{align*}
        &\left(
            P_{\gamma_v, 0 \rightarrow 1}^{-1} (\mu_{\exp_p(v)})
        \right)
        \left(
            P_{\gamma_v, 0 \rightarrow 1} ^{-1} \circ T_v\exp_p e_1,
            \ldots,
            P_{\gamma_v, 0 \rightarrow 1}^{-1} \circ T_v\exp_p e_n
        \right) \\
        &\quad=
        \mu_{\exp_p(v)} (T_v \exp_p (e_1), \ldots, T_v \exp_p (e_n))
        \\ 
        &\quad=
        (\exp_p^*\mu) \at{v} (e_1, \ldots, e_n)
    \end{align*}
    by the definition of the pull-back of a density. But then the
    right hand side of
    \eqref{eq:comp-func-on-TpM-is-norm-partrans-density} is
    \[
    \frac{\exp_p^* \mu \at{v} (e_1, \ldots, e_n)}
    {\mu_p(e_1, \ldots, e_n)}
    = \widetilde{\rho}_p(v)
    \]
    by \eqref{eq:density-comparison-function-on-TpM} proving the
    lemma.
\end{proof}

Since we have a good understanding of the Taylor expansion of $
P_{\gamma_v, 0 \rightarrow 1} ^{-1} \circ T_v\exp_p$ as well as of the
parallel transport $P_{\gamma_v, 0 \rightarrow 1}$ itself, we can use
these results to obtain the complete Taylor expansion of the function
$\widetilde{\rho}_p$ around $v=0$, at least up to the usual recursive
computation of the derivatives of the Jacobi vector fields.
\begin{theorem}
    \label{theorem:taylor-expansion-of-dens-comp-func}
    Let $\mu \in \Secinfty(\Dichten T^*M)$ be a positive density on
    $M$ and let $p \in M$. Then the function $\widetilde{\rho}_p$ from
    \eqref{eq:comp-func-on-TpM-is-norm-partrans-density} has the
    following formal Taylor expansion around $v=0$
    \begin{equation}
        \label{eq:taylor-expansion-of-dens-comp-func}
        \widetilde{\rho}_p(v) \sim_{v \rightarrow 0}
        \mathcal{J}(\E^{\SymD} \mu)(v) \cdot
        \det \left(
            P_{\gamma_v, 0 \rightarrow 1} ^{-1} \circ T_v\exp_p
        \right),
    \end{equation}
    the first orders of which are explicitly given by
    \begin{equation}
        \label{eq:taylo-exp-of-dens-com-func-first-terms}
        \widetilde{\rho}_p(v)
        =
        1 + \alpha_p(v) + \frac{1}{2} (\nabla_v \alpha \at{p})(v)
        + \frac{1}{2} \alpha(v)^2 - \frac{1}{6} \Ric_p(v,v)
        + \cdots
    \end{equation}
    up to terms of order higher than $2$. Here $\alpha \in
    \Secinfty(T^*M)$ is the one-form with $\nabla_X \mu = \alpha(X)
    \mu$.
\end{theorem}
\begin{proof}
    We fix a basis $e_1, \ldots, e_n \in T_pM$, then we first have
    \[
    \left(
        (P_{\gamma_v, 0 \rightarrow 1})^{-1} \mu_{\exp_p(v)}
    \right)
    (A e_1, \ldots, A e_n)
    =
    |\det A| \left(
        (P_{\gamma_v, 0 \rightarrow 1})^{-1} \mu_{\exp_p(v)}
    \right)
    (e_1, \ldots, e_n)       
    \]
    for any linear map $A: T_pM \longrightarrow T_pM$. Since in our
    case $A = P_{\gamma_v, 0 \rightarrow 1} ^{-1} \circ T_v\exp_p$ is
    continuously connected to $\id_{T_pM}$ via $v \longrightarrow 0$,
    we see that the determinant is always positive. Thus we can
    evaluate the determinant of $P_{\gamma_v, 0 \rightarrow 1} ^{-1}
    \circ T_v\exp_p$ in the usual multilinear way. Since in general $A
    e_1 \wedge \cdots \wedge A e_n = \det(A) e_1 \wedge \cdots \wedge
    e_n$, we have to compute
    \begin{align*}
        &P_{\gamma_v, 0 \rightarrow 1}^{-1}
        \left(T_v\exp_p (e_1)\right)
        \wedge \cdots \wedge
        P_{\gamma_v, 0 \rightarrow 1}^{-1}
        \left(T_v\exp_p (e_n)\right)
        \\
        &\qquad \sim_{v \rightarrow 0}
        \left(
            \exp\left(\nabla^\#_{\frac{\partial}{\partial t}}\right)
            J_{e_1}(t)
            \At{t=0}
        \right)
        \wedge \cdots \wedge
        \left(
            \exp\left(\nabla^\#_{\frac{\partial}{\partial t}}\right)
            J_{e_n}(t)
            \At{t=0}
        \right),
    \end{align*}
    and compare it to $e_1 \wedge \cdots \wedge e_n$. From here we get
    the first orders explicitly by
    \eqref{eq:first-terms-of-taylor-of-exp}.
    \begin{align*}
        &
        \left(
            \exp\left(\nabla^\#_{\frac{\partial}{\partial t}}\right)
            J_{e_1}(t)
            \At{t=0}
        \right)
        \wedge \cdots \wedge
        \left(
            \exp\left(\nabla^\#_{\frac{\partial}{\partial t}}\right)
            J_{e_n}(t)
            \At{t=0}
        \right)
        \\
        &\qquad=
        \left(
            e_1 + \frac{1}{6} R_p(v,e_1) v + \cdots
        \right)
        \wedge \cdots \wedge
        \left(
            e_n + \frac{1}{6} R_p(v,e_n) v + \cdots
        \right) \\
        &\qquad=
        e_1 \wedge \cdots \wedge e_n
        + \frac{1}{6} \sum_{\ell=1}^n e_1 \wedge \cdots \wedge
        R_p(v,e_\ell) v \wedge \cdots \wedge e_n
        + \cdots \\
        &\qquad=
        e_1 \wedge \cdots e_n
        + \frac{1}{6} \sum_{\ell=1}^n e_1 \wedge \cdots \wedge
        e^k (R_p(v, e_\ell) v) e_k \wedge \cdots \wedge e_n
        + \cdots \\
        &\qquad=
        e_1 \wedge \cdots e_n
        + \frac{1}{6} e^\ell (R_p(v,e_\ell)v)
        e_1 \wedge \cdots e_n
        + \cdots \\
        &\qquad=
        \left(1 - \frac{1}{6} \Ric_p(v,v) + \cdots \right)
        e_1 \wedge \cdots \wedge e_n,
    \end{align*}
    whence up to second order we get
    \[
    \det \left(
        P_{\gamma_v, 0 \rightarrow 1} ^{-1} \circ T_v\exp_p
    \right)
    = 1 - \frac{1}{6} \Ric_p(v,v) + \cdots.
    \]
    The second step consists in Taylor expanding the parallel
    transport of $\mu$. Here we have by
    Corollary~\ref{corollary:formal-taylor-series-of-partrans} the
    formal Taylor expansion
    \[
     P_{\gamma_v, 0 \rightarrow 1} ^{-1} \mu_{\exp_p(v)}
     \sim_{v \rightarrow 0}
     \mathcal{J}(\E^{\SymD} \mu) (v),
     \tag{$*$}
     \]
     where $\SymD: \Secinfty(\Sym^\bullet T^*M \tensor \Dichten T^*M)
     \longrightarrow \Secinfty(\Sym^{\bullet+1} T^*M \tensor \Dichten
     T^*M)$ is the symmetrized covariant derivative on the density
     bundle. This shows
     \eqref{eq:taylor-expansion-of-dens-comp-func}. The first orders
     of $(*)$ are given by
     \[
     \mathcal{J}(\E^{\SymD} \mu)(v)
     = \mu + (\SymD \mu)(v) 
     + \frac{1}{2} \frac{1}{2} (\SymD^2 \mu)(v,v) + \cdots.
     \]
     Using $\nabla_X \mu = \alpha(X) \mu$ we get
     \begin{align*}
         \SymD \mu &= \alpha \tensor \mu, \\
         \SymD^2 \mu
         &= \SymD \alpha \tensor \mu + \alpha \vee \SymD \mu \\
         &= \SymD \alpha \tensor \mu + \alpha \vee \alpha \tensor \mu, \\
         \SymD^3 \mu 
         &= \SymD^2 \alpha \tensor \mu
         + (\SymD \alpha \vee \alpha + \alpha \vee \SymD \alpha) \tensor
         \mu
         + \alpha \vee \alpha \vee \alpha \tensor \mu \\
         &= \left(
             \SymD^2 \alpha + 2 \SymD \alpha \vee \alpha
             + \alpha \vee \alpha \vee \alpha
         \right) \tensor \mu, \\
         \SymD^4 \mu
         &=
         \left(
             \SymD^3 \alpha + 2 \SymD^2 \alpha \vee \alpha 
             + 2 \SymD \alpha \vee \SymD \alpha 
             + 3 \SymD \alpha \vee \alpha \vee \alpha
         \right) \tensor \mu \\
         &\quad+ \left(
             \SymD^2 \alpha + 2\SymD \alpha \vee \alpha
             + \alpha \vee \alpha \vee \alpha
         \right) \vee \alpha \tensor \mu \\
         &=
         \left(
             \SymD^3 \alpha + 3 \SymD^2 \alpha \vee \alpha
             + 2 \SymD \alpha \vee \SymD \alpha
             + 5 \SymD \alpha \vee \alpha \vee \alpha 
             + \alpha \vee \alpha \vee \alpha \vee \alpha
         \right) \tensor \mu,
     \end{align*}
     and so on by the Leibniz rule. Thus in particular
     \[
     (\SymD \mu) (v) = \alpha(v) \mu
     \quad
     \textrm{and}
     \quad
     (\SymD^2 \mu) (v,v)
     = (\SymD \alpha)(v,v) \mu + 2 \alpha(v)\alpha(v) \mu.
     \]
     Note that $(\SymD \alpha)(v,v) = 2(\nabla_v \alpha)(v)$ by the
     definition of $\SymD$ acting on a one-form. Collecting all terms
     gives
     \begin{align*}
         \mathcal{J} (\E^{\SymD} \mu) (v)
         &= \mu + \alpha(v) \mu 
         + \frac{1}{4} (2 (\nabla_v \alpha)(v) + 2\alpha(v)\alpha(v))
         \mu + \cdots \\
         &=
         \left(
             1 + \alpha(v) 
             + \frac{1}{2} \left((\nabla_v \alpha)(v) +\alpha(v)^2\right)
             + \cdots
         \right) \mu.
     \end{align*}
     Putting things together we have up to second order in $v$
     \begin{align*}
         \widetilde{\rho}_p(v)
         &=
         \left(
             1 + \alpha(v)+\frac{1}{2}(\alpha_v(\alpha))(v)
             +\frac{1}{2} \alpha(v)^2 + \cdots
         \right)
         \cdot
         \left(
             1 - \frac{1}{6} \Ric_p(v,v) + \cdots
         \right) \\
         &=
         1 + \alpha_p(v)
         +\frac{1}{2} 
         \left(
             (\nabla_v\alpha) \at{p} (v) + \alpha_p(v)^2
         \right)
         - \frac{1}{6} \Ric_p(v,v) + \cdots.
     \end{align*}
 \end{proof}

 \begin{remark}
     \label{remark:tayloo-expansion-of-dens-comp-func}
     Again, we note that the evaluation of arbitrarily high orders of
     the Taylor expansion of $\widetilde{\rho}_p$ is reduced to the
     fairly easy computation of $\SymD ^k \mu$ for arbitrary $k$ as
     well as to the slightly more involved Taylor expansion of the
     determinant of the tangent map of $\exp$. However, for the
     tangent map itself we have a fairly easy and completely algebraic
     procedure via the Jacobi fields. Since also $\SymD^k \mu$ can be
     computed in terms of covariant derivatives of the one-form
     $\alpha$ by a simple recursion, we can consider the problem of
     finding higher orders in $\widetilde{\rho}_p$ to be algebraic and
     simple.
 \end{remark}


%% file: appendixB.tex
%
%

\chapter{A Brief Reminder on Stokes Theorem}
\label{satz:stokes-theorem}

%
%

\numberwithin{equation}{chapter}
\renewcommand{\theequation}{\thechapter.\arabic{equation}}

In this appendix we collect a few basic facts on Stokes' Theorem and
its applications in semi-Riemannian geometry.

We start with the following situation: let $U \subseteq M$ be an open
subset and assume that its topological boundary $\iota: \partial U
\hookrightarrow M$ is an embedded submanifold of codimension one. In
this situation we say that $U$ has a \emph{smooth
  boundary}\index{Smooth boundary}.
\begin{Blemma}[Transverse vector field]
    \label{lemma:transversal-vector-field}%
    \index{Transverse vector field}%
    Let $U \subseteq M$ be a non-empty open subset with smooth
    boundary. Then there exists a transverse vector field
    $\mathfrak{n} \in \Secinfty(\iota^\# TM)$ on $\partial U$,
    i.e. for all $p \in \partial U$ the vector $\mathfrak{n}(p)$ is
    transverse to $T_p (\partial U) \subseteq T_pM$.
\end{Blemma}
\begin{proof}
    By assumption we have an atlas of submanifold charts. Since the
    codimension is one, we can label the coordinates $(x^1, \ldots,
    x^n)$ in such a chart in a way that $x^n = 0$ corresponds to the
    boundary $\partial U$ and $x^n > 0$ yields points inside $U$, see
    Figure~\ref{fig:chart-for-boundary}.
    \begin{figure}
        \centering
        \input{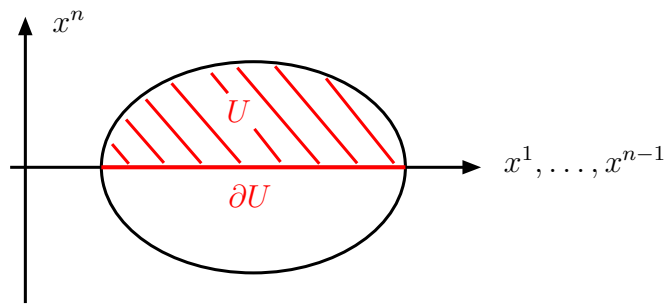}
        \caption{\label{fig:chart-for-boundary}%
          A chart for the boundary $\partial U$.}
    \end{figure}
    Clearly, we can find an atlas with this feature. Now $-
    \frac{\partial}{\partial x^n}$ is pointing outwards of $U$ in such
    a chart. It is now easy to check that the property of pointing
    outwards is \emph{convex}, i.e. convex combinations of (locally
    defined) vector fields which point outwards point outwards
    again. Thus a partition of unity argument gives a smooth vector
    field $\mathfrak{n}$ on $\partial U$ which points outwards at
    every point. In particular $\mathfrak{n}(p)$ is transverse to
    $T_p (\partial U)$ at every $p \in \partial U$.
\end{proof}

On a connected component of $\partial U$ a transverse vector field is
either pointing outwards or pointing inwards. We can use transverse
vector field to induce orientations:
\begin{Blemma}
    \label{lemma:induced-orientation-by-transversal-vectorfield}%
    \index{Orientable}%
    Assume $M$ is orientable and $\omega \in
    \Secinfty(\Anti^{\mathrm{top}} T^*M)$ is a nowhere vanishing
    $n$-form. If $\mathfrak{n} \in \Secinfty(\iota^\# TM)$ is a
    transverse vector field to $\partial U$ then $\ins_\mathfrak{n}
    \omega \in \Secinfty(\Anti^{\mathrm{top}} T^* \partial U)$ is a
    nowhere vanishing $(n-1)$-form on $\partial U$.
\end{Blemma}
\begin{proof}
    Of course, here we view $\ins_\mathfrak{n} \omega$ as a
    $(n-1)$-form defined on $\partial U$ only. If $e_2, \ldots, e_n$
    form a basis in $T_p \partial U$ then $\mathfrak{n}(p), e_2,
    \ldots, e_n \in T_pM$ form a basis by transversality. Thus,
    $\omega$ evaluated on this basis is non-zero, hence
    $\ins_\mathfrak{n} \omega$ is nowhere vanishing.
\end{proof}

\begin{Bdefinition}[Induced orientation]
    \label{definition:induced-orientation}%
    \index{Orientation!induced}%
    Let $U \subseteq M$ be open with smooth boundary $\partial U$. If
    $M$ is oriented then the induced orientation of $\partial U$ is
    defined by the $(n-1)$-form $\ins_\mathfrak{n} \omega$ where
    $\omega \in \Secinfty(\Anti^{\mathrm{top}} T^*M)$ is a positively
    oriented $n$-form and $\mathfrak{n} \in \Secinfty(\iota^\# TM)$ is
    a transverse vector field pointing outwards.
\end{Bdefinition}
\begin{Bremark}
    \label{remark:induced-orientation}
    It is an easy check that this is indeed well-defined, i.e. the
    induced orientation of $\partial U$ only depends on the
    orientation of $M$ but not on the choices of $\omega$ and
    $\mathfrak{n}$.
\end{Bremark}
With respect to these orientations we can integrate top degree
forms. The fundamental feature of such integrations is then formulated
in Stokes' Theorem:
\begin{Btheorem}[Stokes]
    \label{theorem:stokes}%
    \index{Stokes' Theorem}%
    Let $M$ be oriented and let $U \subseteq M$ be a non-empty open
    subset with smooth boundary $\iota: \partial U \hookrightarrow M$,
    equipped with the induced orientation. Then for all $\omega \in
    \Secinfty_0(\Anti^{n-1} T^*M)$ we have
    \begin{equation}
        \label{eq:stokes}
        \int_U \D \omega = \int_{\partial U} \iota^* \omega.
    \end{equation}
\end{Btheorem}
For a proof of this well-known theorem one may consult any textbook on
differential geometry, see. e.g. \cite[Thm.~8.11]{michor:2001a:script}
or \cite[Thm.~14.9]{lee:2003a}.
\begin{Bremark}
    \label{remark:stokes}
    There are many generalizations of \eqref{eq:stokes} for forms and
    boundaries of less regularity than $\Cinfty$: this is reasonable
    to expect since ultimately \eqref{eq:stokes} is an equation
    between integrals whence only measure-theoretic properties should
    be relevant. In particular, the theorem still holds for boundaries
    with corners, see \cite[Thm.~14.20]{lee:2003a}.
\end{Bremark}
We shall now use this theorem to obtain similar results for the
non-oriented situation: this is still plausible to be possible as
changing the orientation from $\omega$ to $-\omega$ should
produce the same sign on both sides of \eqref{eq:stokes}. We shall now
see how this can be made precise.
\begin{Blemma}
    \label{lemma:induces-density}
    Let $U \subseteq M$ be open with smooth boundary and let
    $\mathfrak{n} \in \Secinfty(\iota^\# TM)$ be a transverse vector
    field.
    \begin{lemmalist}
    \item \label{item:induced-density} For $\mu \in \Secinfty(
        \Dichten T^*M)$, $p \in \partial U$ and $e_2, \ldots e_n \in
        T_p (\partial U)$ the definition
        \begin{equation}
            \label{eq:induced-density}
            (\ins_\mathfrak{n} \mu) \at{p} (e_2, \ldots, e_n)
            = \mu_p(\mathfrak{n}(p), e_2, \ldots, e_n)
        \end{equation}
        defines a smooth density $\ins_\mathfrak{n} \mu \in \Secinfty(
        \Dichten T^*(\partial U))$.
    \item \label{item:inducing-densites-is-continuous} The map
        \begin{equation}
            \label{eq:inducing-densities-is-continuous}
            \Secinfty( \Dichten T^*M) \ni \mu
            \; \mapsto \;
            \ins_\mathfrak{n} \mu \in
            \Secinfty(\Dichten T^*(\partial U))
        \end{equation}
        is continuous and $\Cinfty(M)$-linear in the sense that for $f
        \in \Cinfty(M)$ we have
        \begin{equation}
            \label{eq:densitiy-inducint-function-linear}
            \ins_\mathfrak{n}(f\mu) = \iota^*f \ins_\mathfrak{n} \mu.
        \end{equation}
    \item \label{item:positive-induces-positive} For a positive
        density $\mu$ also $\ins_\mathfrak{n} \mu$ is positive.
    \end{lemmalist}
\end{Blemma}
\begin{proof}
    We choose a submanifold chart $(V,x)$ of $M$ such that $x^n = 0$
    corresponds to $\partial U$ in this chart. Then any transverse
    vector field $\mathfrak{n}$ has a nontrivial
    $\frac{\partial}{\partial x^n}$-component along $x^n = 0$,
    i.e. writing
    \[
    \mathfrak{n} \at{V} = \mathfrak{n}^i \frac{\partial}{\partial x^i}
    \]
    with $\mathfrak{n}^i \in \Secinfty( \partial U \cap V)$ we have
    $\mathfrak{n}^n(x^1, \ldots, x^{n-1}) \neq 0$. If $e_2, \ldots
    e_n$ are a frame at $p \in \partial U$ the it is easy to check
    that $\ins_\mathfrak{n} \mu$ transforms correctly under the change
    of frames. Thus \eqref{eq:induced-density} defines a density
    indeed. Moreover, if $\mu \at{V} = \mu_V | \D\!x^1 \wedge \cdots
    \wedge x^n |$ where $\mu_V \in \Cinfty(V)$ is the local form of
    $\mu$ in this chart then $\ins_\mathfrak{n} \mu \at{\partial U
      \cap V} = \iota^*(\mu_V) |\D\!x^1 \wedge \cdots \wedge
    \D\!x^{n-1}| |\mathfrak{n}^n|$ whence the local function
    representing $\ins_\mathfrak{n} \mu$ is $\iota^* \mu_V
    |\mathfrak{n}^n|$. Since $\mathfrak{n}^n$ is everywhere different
    from zero, this is smooth again, showing that $\ins_\mathfrak{n}
    \mu$ is indeed smooth. Moreover, if $\mu$ is positive we see that
    $\ins_\mathfrak{n} \mu$ is positive as well. The continuity is
    again a consequence of the above local expression as we can use
    these submanifolds charts to characterize the Fr\'echet topologies
    of $\Secinfty(\Dichten T^*M)$ and $\Secinfty(\Dichten T^*
    (\partial U))$, respectively. Finally,
    \eqref{eq:densitiy-inducint-function-linear} is clear from the
    definition.
\end{proof}

Thus having specified a transverse vector field $\mathfrak{n}$ of
$\partial U$ we can speak of the \emph{induced density}
$\ins_\mathfrak{n} \mu$ coming from a density $\mu$ on $M$. From the
above definition it is clear that
\begin{equation}
    \label{eq:inducing-of-scaled-trans-field}
    \ins_{f \mathfrak{n}} \mu = |f| \ins_\mathfrak{n} \mu
\end{equation}
for any nowhere vanishing function $f \in \Cinfty(\partial U)$.

We now specialize to the following situation: assume that $M$ is in
addition a semi-Riemannian manifold with metric $g$. Moreover, we
assume that $\partial U$ allows for a transverse vector field
$\mathfrak{n}$ which is nowhere lightlike, where we shall use the
notions of timelike, spacelike and lightlike vectors as in the
Lorentzian situation.  Then on each connected component
$g(\mathfrak{n},\mathfrak{n})$ is either positive or negative whence
$\mathfrak{n}$ is either timelike or spacelike everywhere on this
connected component. We can now achieve two things: first we can
arrange $\mathfrak{n}$ in such a way that $\mathfrak{n}(p)$ is not
only transverse to $T_p (\partial U)$ but orthogonal.  Moreover, we
can normalize $\mathfrak{n}(p)$ at every $p \in \partial U$. Finally,
we choose $\mathfrak{n}(p)$ to point outwards: his determines
$\mathfrak{n}(p)$ uniquely. Indeed, since $T_p (\partial U) \subseteq
T_pM$ has codimension one the annihilator space $T_p (\partial U)^\ann
\subseteq T_p^*M$ of one-forms annihilating $T_p(\partial U)$ is
one-dimensional. Then $\mathfrak{n}(p) \in (T_p(\partial U)^\ann)^\#$
is orthogonal to all of $T_p \partial U$ and uniquely determined as
$\mathfrak{n}(p)^\flat \in T_p (\partial U)^\ann$ by definition. Then
normalizing and orienting it gives a unique vector.
\begin{Bdefinition}[Normal vector field]
    \label{definition:normal-vector-field}%
    \index{Normal vector field}%
    Let $(M,g)$ be semi-Riemannian and let $U \subseteq M$ be open
    with smooth boundary. Assume that the annihilator spaces $T_p
    (\partial U)^\ann \subseteq T_p^*M$ of $\partial U$ are never
    lightlike (with respect to $g^{-1}$). Then the unique normalized
    transverse vector field $\mathfrak{n} \in \Secinfty(\iota^\# TM)$
    which is orthogonal to $\partial U$ and pointing outward is called
    the normal vector field of $\partial U$.
\end{Bdefinition}
This allows us to obtain a uniquely determined metric and density on
$\partial U$ as follows:
\begin{Bdefinition}
    \label{definition:induced-metric}
    Let $M$ be semi-Riemannian and let $U \subseteq M$ be open with
    connected smooth boundary such that $T_p (\partial U)^\ann
    \subseteq T_p^*M$ is never lightlike. Then the induced metric on
    $\partial U$ is $\iota^*g \in \Secinfty(\Sym^2 T^* \partial U)$.
\end{Bdefinition}
\begin{Blemma}
    \label{lemma:induced-metric}
    Under the above assumptions, $\iota^*g$ is a semi-Riemannian
    metric on $\partial U$. Moreover,
    \begin{equation}
        \label{eq:induced-metric}
        \mu_{\iota^*g} = \ins_\mathfrak{n} \mu_g,
    \end{equation}
    where $\mathfrak{n} \in \Secinfty(\iota^\# TM)$ is the normal
    vector field of $\partial U$.
\end{Blemma}
\begin{proof}
    Let $p \in \partial U$. Then we have to show that $\iota^*g
    \at{p}$ is indeed non-degenerate (the Riemannian case is
    trivial). We find in $T_pM$ a semi-Riemannian frame $e_1, \ldots,
    e_n$ such that $e_1 = \mathfrak{n}(p)$. Then $e_2, \ldots e_n$ are
    a basis of $T_p \partial U$ with
    \[
    \iota^*g \at{p} (e_i,e_j) = g_p(e_i,e_j) = \pm \delta_{ij}
    \]
    for $i,j = 2, \ldots, n$. Thus $\iota^*g$ is non-degenerate. Its
    signature can be obtained from knowing whether $\mathfrak{n}$ is
    time- or spacelike and from the signature of $g$. In particular,
    $e_2, \ldots e_n$ is a semi-Riemannian frame for $\iota^*g$. From
    this we see that by definition $\mu_{\iota^*g} \at{p}(e_2, \ldots,
    e_n) = 1$. On the other hand $(\ins_\mathfrak{n} \mu_g)(e_2,
    \ldots, e_n) = \mu_g \at{p} (\mathfrak{n}(p), e_2, \ldots, e_n) =
    1$ as $\mathfrak{n}(p), e_2, \ldots, e_n$ is a semi-Riemannian
    frame for $g$. Thus the two densities coincide as they coincide on
    one frame.
\end{proof}

We can now use the normal vector field $\mathfrak{n}$ to formulate
Gauss' Theorem as a consequence of Stokes' Theorem:
\begin{Btheorem}[Gauss]
    \label{theorem:gauss-theorem}%
    \index{Gauss' Theorem}%
    Let $(M,g)$ be a semi-Riemannian manifold and $U \subseteq M$ open
    with smooth connected boundary $\iota: \partial U \hookrightarrow
    M$. Assume that $T(\partial U)^\ann$ is never lightlike. Then for
    all vector fields $X \in \Secinfty_0(TM)$ we have
    \begin{equation}
        \label{eq:gauss}
        \int_U \divergenz (X) \mu_g
        = \epsilon \int_{\partial U} g(\iota^\#X, \mathfrak{n})
        \mu_{\iota^*g},
    \end{equation}
    where $\epsilon = \SP{\mathfrak{n},\mathfrak{n}} \in \{1, -1\}$.
\end{Btheorem}
\begin{proof}
    First we consider the oriented case. Thus the left hand side is
    \[
    \int_U \divergenz(X) \mu_g = \int_U \divergenz(X) \Omega_g,
    \]
    with the positively oriented volume form $\Omega_g$ yielding
    $\mu_g$ under the canonical map from forms to densities, see
    \cite[Prop.~2.2.42]{waldmann:2007a}. Note that $\divergenz(X)$ can
    alternatively be computed via $\divergenz(X) \Omega_g = \Lie_X
    \Omega_g = \D (\ins_X \Omega_g)$. Thus we can apply Stokes'
    Theorem and get
    \[
    \int_U \divergenz(X) \Omega_g
    = \int_U \D (\ins_X \Omega_g)
    = \int_{\partial U} \iota^* (\ins_X \Omega_g).
    \tag{$*$}
    \]
    Now along $\partial U$ we can decompose $X$ into its
    $\mathfrak{n}$-component and parallel components. We have
    \[
    X(p) = \epsilon g_p(X(p),\mathfrak{n}(p)) \mathfrak{n}(p)
    + X_{\parallel}(p)
    \]
    where $X_\parallel(p)$ is orthogonal to $\mathfrak{n}(p)$ and
    hence in $T_p \partial U$. Note that we need the constant
    $\epsilon$ here since $g_p(\mathfrak{n}(p), \mathfrak{n}(p)) =
    \epsilon$ may be $-1$ instead of $1$. However, $\epsilon$ is
    constant on $\partial U$. Now we note that
    \[
    \iota^* \ins_{X_\parallel(p)} \Omega_g \at{p} = 0,
    \]
    since evaluating $\ins_{X_\parallel(p)} \Omega_g \at{p}$ on $n-1$
    tangent vectors in $T_p \partial U$ means evaluating
    $\Omega_g\at{p}$ on $n$ tangent vectors in $T_p \partial U$. Thus
    they are necessarily linear dependent. This shows that $\iota^*
    \ins_{X(p)} \Omega_g \at{p} = \epsilon g_p(X(p), \mathfrak{n}(p))
    \ins_{\mathfrak{n}(p)} \Omega_g \at{p}$. Finally, it is easy to
    see that $\ins_{\mathfrak{n}(p)} \Omega_g \at{p}$ is the (by
    definition positively oriented) semi-Riemannian volume form of
    $\iota^*g$. This is clear be the same argument as for
    $\mu_{\iota^*g}$ in Lemma~\ref{lemma:induced-metric}. This finally
    shows
    \[
    \int_U \divergenz(X) \mu_g
    = \int_{\partial U} \divergenz(X) \Omega_g
    = \epsilon \int_{\partial U} g(i^\#X, \mathfrak{n})
    \Omega_{i^*g}
    = \epsilon \int_{\partial U} g(i^\#X, \mathfrak{n})
    \mu_{i^*g},
    \tag{$*$*}
    \]
    and hence \eqref{eq:gauss}. If we change the orientation from
    $\Omega_g$ to $- \Omega_g$ then the induced orientation
    $\Omega_{\iota^*g}$ changes to $- \Omega_{\iota^*g}$ since the
    normal vector field $\mathfrak{n}$ remains \emph{unchanged}:
    ``pointing outwards'' does not depend on any choice of
    orientation. Thus we see that the left and right side of ($**$)
    both change their sign. From this we conclude that
    \eqref{eq:gauss} also holds in the non-oriented case: indeed, by a
    partition of unity argument we can chop down $X$ into small pieces
    having support in a chart. There we can choose an orientation and
    use ($**$). Summing up again is allowed as the validity of ($**$)
    does not depend on the local choices.
\end{proof}

A particular case of interest is the following. Assume $(M,g)$ is a
Lorentzian manifold and the boundary $\partial U$ is spacelike. Then
the normal vector field $\mathfrak{n}$ is timelike and we have
\begin{equation}
    \label{eq:gauss-for-lorentzian-case}
    \int_U \divergenz(X) \mu_g
    = \int_{\partial U} g(\iota^\#X, \mathfrak{n}) \mu_{\iota^*g}
\end{equation}
for all $X \in \Secinfty_0(TM)$.
